%% file: main.tex
\documentclass[11pt]{article}

\usepackage[numbers]{natbib}
\usepackage{amsmath}
\usepackage{amssymb}
\usepackage{amsthm}
\usepackage{xcolor}
\usepackage{array}
\usepackage{booktabs}
\usepackage{tabularx}
\usepackage{longtable}
\usepackage{float}
\usepackage{hyperref}
\usepackage{enumitem}
\usepackage[margin=1in]{geometry}

\hypersetup{colorlinks=true,citecolor=blue,linkcolor=blue,hypertexnames=false}
\expandafter\def\expandafter\UrlBreaks\expandafter{\UrlBreaks\do\-}
\allowdisplaybreaks

\input{_0_vars}
\input{_1_maths}

\theoremstyle{plain}
\newtheorem{theorem}{Theorem}[section]
\newtheorem{lemma}[theorem]{Lemma}

\theoremstyle{remark}

\begin{document}

\date{}
\title{\paperTitle}
\author{\paperAuthor}
\maketitle

\begin{abstract}
\input{00_abstract}
\end{abstract}

\clearpage
\begingroup
\small
\tableofcontents
\endgroup
\clearpage

\input{_2_body}

\section*{AI Disclosure}
\addcontentsline{toc}{section}{AI Disclosure}
We used OpenAI Codex with \texttt{gpt-6.0-astra}, \texttt{gpt-6.2-astra} and \texttt{gpt-6.2-sol}
to develop and check bounds and revise the exposition. Among Anthropic's
models, \texttt{claude-fable-5-2} assisted with proof development,
\texttt{claude-opus-5-5} with writing out and checking the proofs of
Section~\ref{sec:joint_accuracy}, proof auditing, and manuscript revision,
and \texttt{claude-haiku-5-5} with citation checks.

\bibliographystyle{alpha}
\begingroup
\raggedright
\bibliography{ref}
\endgroup

\clearpage
\appendix
\input{_3_app}

\end{document}

%% file: _0_vars.tex
\newcommand{\paperTitle}{Optimal $\ell_p$ Row Sparsification}
\newcommand{\paperAuthor}{Zhao Song\thanks{\texttt{magic.linuxkde@gmail.com}.} 
\and
Lichen Zhang\thanks{\texttt{lichenz@mit.edu}. Massachusetts Institute of Technology.}
}

%% file: _1_maths.tex
\DeclareMathOperator{\supp}{supp}

\DeclareMathOperator{\rank}{rank}
\newcommand{\E}{\mathbb E}
\newcommand{\R}{\mathbb R}
\DeclareMathOperator{\tr}{tr}
\renewcommand{\d}{\mathrm{d}}

%% file: 00_abstract.tex
Let $S_p(r,\epsilon)$ be the worst-case number of reweighted original
rows needed to preserve $\|Ax\|_p^p$ within a factor $1\pm\epsilon$
for every query $x$ on a rank-$r$ matrix $A$. Write $L:=\log(2/\epsilon)$
and take $0<\epsilon\le1/2$. We obtain optimal $S_p(r, \epsilon)$ up to constant, in several regimes.

For each fixed $0<p<2$, we prove the tight bound
$S_p(r,\epsilon)=\Theta_{p,a}(r/\epsilon^2)$ throughout polynomial
accuracy $r^{-a}\le\epsilon\le1/2$ for every fixed $a>0$.

For even $p$ the bounds match over the full accuracy range.
At $p=2$ the classical count is $\Theta(\min\{r/\epsilon^2,r^2\})$.
For each fixed even integer $p\ge4$, we prove
\[
 S_p(r,\epsilon)=\Theta_p(\min \{ \epsilon^{-2}
   r^{p/2}, \epsilon^{-1} r^{(p+1)/2} ,r^p
  \}).
\]

For each fixed $p>0$ outside the even integers and sufficiently large
$r$, the tight bound at very high accuracy $L\ge C_p r\log(2r)$ is
\[
 S_p(r,\epsilon)=\Theta_p(
   (\epsilon^{-2}r^{1/2-p}\log(2r))^{(r-1)/(r+2p)}).
\]

We further provide some improved upper bound results on row count. For $1\leq p<2$ we prove an $O(r/\epsilon^2)$ upper bound, recovering that of Reis--Rothvoss at $p=1$ and removing a factor of $\log(2r)$ for $1<p<2$. For $0<p<1$ we show $O(\epsilon^{-2}r\min\{1/p,\log(2r)\})$ rows suffice. For any non-even $p>2$ and $L\leq c_pr$, we show that the count is $r^{p/2}/\epsilon^2$ up to polylogarithmic factors in $r, 1/\epsilon$.

%% file: _2_body.tex
\input{01_intro}
\input{05_overview}
\input{10_preli}
\input{30_lp_p_below_2}
\input{35_small_p}
\input{40_lp_p_above_2}
\input{60_lower}
\input{80_proof}
\input{81_proof}
\input{82_proof}

%% file: 01_intro.tex
\section{Introduction}

Given a matrix $A\in\R^{m\times n}$ with rows $a_i^\top$ and an
exponent $p>0$, the $\ell_p$ row sparsification problem asks for a
subset $S\subseteq\{1,\ldots,m\}$ and positive weights $w_i$ such that
\begin{equation}\label{eq:row_sparsification}
 (1-\epsilon)\|Ax\|_p^p
 \le\sum_{i\in S}w_i|a_i^\top x|^p
 \le(1+\epsilon)\|Ax\|_p^p
 \qquad\text{for every }x\in\R^n,
\end{equation}
where $\|Ax\|_p^p=\sum_i|a_i^\top x|^p$ and $0<\epsilon<1$. We call
such a pair $(S,w)$ a \emph{sparsifier} of $A$. With $r:=\rank(A)$,
we study the worst-case support
\begin{equation}\label{eq:support_function}
 S_p(r,\epsilon):=\sup_{\rank(A)=r}\ \min\{|S|:(S,w)\text{ satisfies
 Eq.~\eqref{eq:row_sparsification}}\},
\end{equation}
and assume $r\ge2$, since one row suffices for rank one. Row
sparsification is a weighted discretization of a subspace of
$\ell_p^m$, studied by Lewis~\cite{lew78}, Schechtman~\cite{sch87},
and Bourgain, Lindenstrauss, and Milman~\cite{blm89}. Cohen and
Peng~\cite{cp15} used Lewis weights for efficient reweighted row
sampling, and such sparsifiers are the basic primitive behind $\ell_p$
subspace embeddings and $\ell_p$ regression with few rows.

For $p=2$, spectral sparsification~\cite{bss12} gives
$O(r/\epsilon^2)$ rows. For $p=1$, Talagrand~\cite{tal90} proved
$O(r\log(2r)/\epsilon^2)$, and Reis and Rothvoss~\cite{rr26} improved
it to $O(r/\epsilon^2)$. For $1<p<2$, the best previous bound is
$O(r\log(2r)/\epsilon^2)$, due to Dai and
Kosov~\cite[Corollary~1.4]{dk26}. For fixed $p>2$, some rank-$r$
matrix needs $\Omega_p(r^{p/2})$ rows~\cite{bdgjn77}, the best previous
upper bound is $O_p(r^{p/2}\epsilon^{-2}\log(2r/\epsilon))$~\cite[Corollary~4.5 and Remark~4.6]{kos21},
and $O_p(r^{p/2}/\epsilon^2)$ holds for even $p$~\cite{sch11}. Cut
sparsification gives the lower bound $\Omega(r/\epsilon^2)$ for every
$p$ when $r\ge C/\epsilon^2$~\cite{akw14}.

The cut lower bound leaves open the rate at finer polynomial
accuracies below two. Above two, the prior bounds leave open the joint
dependence on rank and accuracy. In this paper, we remove the
logarithmic factor from the upper bound for $1<p<2$ and determine the
worst-case count at every polynomial accuracy for fixed $0<p<2$. For
each fixed even integer $p\ge4$, we determine the count throughout
$0<\epsilon\le1/2$. For fixed \emph{non-even} $p>2$, meaning $p$ lies
outside the even integers, we give matching bounds up to polylogarithmic
factors in $r$ and $1/\epsilon$ when
$r\ge C_p\log(2/\epsilon)$.\footnote{Unsubscripted constants $c,C$ are
absolute, and constants with subscripts, such as $C_p$, depend only on
the indicated parameters. Their values may change between
occurrences.} For every fixed non-even $p>0$, we also give a tight bound
when $\log(2/\epsilon)\ge C_p r\log(2r)$ and $r$ is sufficiently large.

Our upper bounds use the partial coloring of Reis and
Rothvoss~\cite{rr26}. Its key input is a lower bound on the Gaussian
measure of the cone of signed weight changes that decrease no energy
value. We prove it for $0<p\le2$ by a one-sided estimate for a Gaussian
trigonometric field in Lewis coordinates, which loses only a constant
per dimension. Above two, the coloring also cancels the leading
directions of a polynomial approximation of the energy. Our lower
bounds test the energy on structured queries, built from mutually
unbiased bases, binary quadratic forms, tensor products of hard
instances, and mixed coefficients (Section~\ref{sec:technical_overview}).

\subsection{Our Results}

Our results below two are as follows.

\begin{theorem}[Exponents up to two, informal version of Theorems~\ref{thm:general_lp_field} and~\ref{thm:small_p_uniform}]\label{thm:lp_one_to_two_informal}
Let $A\in\R^{m\times n}$ have rank $r\ge1$. For $1\le p\le2$ and
$0<\epsilon<1$, there is a sparsifier satisfying
Eq.~\eqref{eq:row_sparsification} with $O(r/\epsilon^2)$ rows. For
$0<p<1$ and $0<\epsilon\le1/2$, there is one with
$O(r\epsilon^{-2}\min\{1/p,\log(2r)\})$ rows. Both implicit constants
are absolute.
\end{theorem}

\begin{theorem}[Optimality below two, informal version of Theorem~\ref{thm:below_two_main}(a)]\label{thm:joint_polynomial_informal}
Fix $0<p<2$ and $a>0$. Uniformly for $r^{-a}\le\epsilon\le1/2$,
\[
 S_p(r,\epsilon)=\Theta_{p,a}(r/\epsilon^2).
\]
\end{theorem}

In Theorem~\ref{thm:lp_one_to_two_informal}, the bound for $1\le p\le2$
recovers Reis and Rothvoss at $p=1$ and removes the factor $\log(2r)$
of Dai and Kosov for $1<p<2$. The bound for $0<p<1$ does not depend on
$p$ once $p\le1/\log(2r)$. In the lower bound, random block flips of a partite
sign dictionary isolate the top coefficient of the energy, on which a
small support leaves a large discrepancy. Above two, our upper bound
is $r^{p/2}/\epsilon^2$ up to a factor $O_p(\log(2r))$, and our lower
bounds determine the count for even $p$, and for non-even $p$ up to
logarithmic factors when $r\ge C_p\log(2/\epsilon)$.

\begin{theorem}[Above two, informal version of Theorem~\ref{thm:lp_above_two}]\label{thm:lp_above_two_informal}
Let $p>2$, let $A\in\R^{m\times n}$ have rank $r\ge2$, and let
$0<\epsilon\le1/2$. With $t:=\log(2\sqrt{\log(2r)}/\epsilon)$, there is
a sparsifier satisfying Eq.~\eqref{eq:row_sparsification} with rows
\[
 O_p(r^{p/2}\epsilon^{-2}(1+\log(2r)/t^{p-2})).
\]
\end{theorem}

\begin{theorem}[Even exponents, informal version of Theorems~\ref{thm:lp_above_two} and~\ref{thm:lp_lower_main}(a)]\label{thm:lp_even_lower_informal}
For every fixed even integer $p\ge4$, every $r\ge2$, and every
$0<\epsilon\le1/2$,
\[
 S_p(r,\epsilon)=\Theta_p(\min\{r^{p/2}/\epsilon^2,\,
                 r^{(p+1)/2}/\epsilon,\,r^p\}).
\]
\end{theorem}

\begin{theorem}[Non-even exponents above two, informal version of Theorem~\ref{thm:lp_lower_main}(b)]\label{thm:lp_noneven_amplified_informal}
Fix a non-even $p>2$. If $r\ge C_p\log(2/\epsilon)$, there is a
rank-$r$ matrix for which every sparsifier satisfying
Eq.~\eqref{eq:row_sparsification} has
$\widetilde{\Omega}_p(r^{p/2}/\epsilon^2)$ rows, where the tilde hides
a polylogarithmic factor in $1/\epsilon$.
\end{theorem}

The extra factor in Theorem~\ref{thm:lp_above_two_informal} is at most
two once $t^{p-2}\ge\log(2r)$. For even $p$ we also get
$O_p(r^{(p+1)/2}/\epsilon)$ rows and an exact representation with
$\binom{r+p-1}{p}$ rows, and the three terms of the minimum are
attained for $\epsilon\ge r^{-1/2}$, for
$r^{-(p-1)/2}\le\epsilon\le r^{-1/2}$, and for
$\epsilon\le r^{-(p-1)/2}$. For non-even $p$, the previous lower bounds
with accuracy dependence $\epsilon^{-2}$ grew only linearly in
$r$~\cite{lww21,li26}, and
Theorem~\ref{thm:lp_noneven_amplified_informal} reaches $r^{p/2}$. Its
proof tensors a low-rank hard instance with nearly orthogonal vectors, so that product queries force a sparsifier to
approximate every copy separately.

At high accuracy, write $L:=\log(2/\epsilon)$ and
$\eta:=(r-1)/(r+2p)$. For fixed rank and non-even $p\ge1$,
Li~\cite[Theorem~1.2]{li26fixed} gives $O_{r,p}(\epsilon^{-2\eta})$
reweighted original rows, matching the lower bound recorded by Li, Lin,
and Woodruff~\cite[Eq.~(19)]{llw23}. Theorem~\ref{thm:high_accuracy_main}(c)
proves both bounds for every non-even $p>0$, extending the rate to
$0<p<1$.
Our main high-accuracy result is uniform in the rank.

\begin{theorem}[High accuracy, informal version of Theorem~\ref{thm:high_accuracy_main}(a) and~(b)]\label{thm:high_accuracy_informal}
Fix a non-even $p>0$ and $a>0$. For all sufficiently large $r$ and
$L\ge ar$,
\[
 (\epsilon^{-2}r^{1/2-p}\log(2r))^\eta
 \lesssim_{p,a}S_p(r,\epsilon)\lesssim_{p,a}
 (\epsilon^{-2}r^{1-p})^\eta .
\]
Moreover, $S_p(r,\epsilon)=\Theta_p((\epsilon^{-2}r^{1/2-p}\log(2r))^\eta)$
for $L\ge C_pr\log(2r)$.
\end{theorem}

The ratio of the upper bound to the lower bound is
$O_{p,a}((\sqrt r/\log(2r))^\eta)$ uniformly as $L/r$ grows. The optimal
rank prefactor at $L\asymp r$ remains open.
Section~\ref{sec:high_accuracy} proves these results.
Table~\ref{tab:tight_regimes} collects the regimes in which our bounds
are tight, and Table~\ref{tab:joint_regimes} locates our bounds at
higher accuracy by regime.

\begin{table}[!ht]
\centering
\small
\begin{tabularx}{\linewidth}{@{}>{\raggedright\arraybackslash}p{0.17\linewidth} >{\raggedright\arraybackslash}X l >{\raggedright\arraybackslash}p{0.15\linewidth} l@{}}
\toprule
Exponent & Regime & $S_p(r,\epsilon)$ & Upper bound & Lower bound \\
\midrule
Fixed $0<p<1$ & Fixed $r$, $\epsilon\to0$
 & $\Theta_{r,p}(\epsilon^{-2\eta})$
 & Thm.~\ref{thm:high_accuracy_main}(c) & Thm.~\ref{thm:high_accuracy_main}(c) \\
Fixed $0<p<2$ & $r^{-a}\le\epsilon\le1/2$, fixed $a>0$
 & $\Theta_{p,a}(r/\epsilon^2)$
 & Thm.~\ref{thm:general_lp_field} & Thm.~\ref{thm:below_two_main}(a) \\
$1\le p\le2$ & $r\ge C/\epsilon^2$, absolute constants
 & $\Theta(r/\epsilon^2)$
 & Thm.~\ref{thm:general_lp_field} (new for $1<p<2$) & \cite{akw14} \\
\addlinespace
Fixed non-even $p>0$ & $L\ge C_pr\log(2r)$, $r$ large
 & $\Theta_p((\epsilon^{-2}r^{1/2-p}\log(2r))^\eta)$
 & Thm.~\ref{thm:high_accuracy_main}(b) & Thm.~\ref{thm:high_accuracy_main}(a) \\
Fixed non-even $p>2$ & $r\ge C_pL$
 & $r^{p/2}/\epsilon^2$ up to polylog$(r/\epsilon)$
 & Thm.~\ref{thm:lp_above_two} & Thm.~\ref{thm:lp_lower_main}(b) \\
\addlinespace
Fixed even $p\ge4$ & $r^{-1/2}\le\epsilon\le1/2$
 & $\Theta_p(r^{p/2}/\epsilon^2)$
 & \cite{sch11} & Thm.~\ref{thm:lp_lower_main}(a) \\
 & $r^{-(p-1)/2}\le\epsilon\le r^{-1/2}$
 & $\Theta_p(r^{(p+1)/2}/\epsilon)$
 & Thm.~\ref{thm:lp_above_two} & Thm.~\ref{thm:lp_lower_main}(a) \\
 & $0<\epsilon\le r^{-(p-1)/2}$
 & $\Theta_p(r^p)$
 & Conic Carath\'eodory & Thm.~\ref{thm:lp_lower_main}(a) \\
\bottomrule
\end{tabularx}
\caption{Regimes in which the worst-case support is determined, with
the source of each bound. Here $r\ge2$, $0<\epsilon\le1/2$,
$L:=\log(2/\epsilon)$ and $\eta:=(r-1)/(r+2p)$.
The subscripts on $\Theta$ specify the allowed dependence of the
constants. The constant $C$ is absolute, and $C_p$ depends only on $p$. The upper
bounds for $p=1$ and $p=2$ in the third row are due to~\cite{rr26}
and~\cite{bss12}. The matching bounds in the fourth row are proved in
Section~\ref{sec:high_accuracy}.}
\label{tab:tight_regimes}
\end{table}

\begin{table}[!ht]
\centering
\small
\begin{tabularx}{\linewidth}{@{}>{\raggedright\arraybackslash}p{0.21\linewidth} >{\raggedright\arraybackslash}X >{\raggedright\arraybackslash}X@{}}
\toprule
Accuracy regime & Exponents $0<p<2$ & Non-even exponents $p>2$ \\
\midrule
Polynomial, $\epsilon\ge r^{-a}$
 & Tight: $\Theta_{p,a}(r/\epsilon^2)$
   (Thms.~\ref{thm:general_lp_field}, \ref{thm:below_two_main}(a))
 & Tight up to polylog when $r\ge C_pL$
   (Thms.~\ref{thm:lp_above_two}, \ref{thm:lp_lower_main}(b)) \\
Subexponential, $C_p\Lambda\le L\le r$
 & Upper bounds in Thm.~\ref{thm:below_two_main}(b), Lem.~\ref{lem:joint-positive-caps}.
   Lower bounds in Thm.~\ref{thm:below_two_main}(b), Lem.~\ref{lem:joint-adaptive-shifts}.
   Tight for regular components (Lem.~\ref{lem:joint-regular-components}).
 & Upper bounds in Thms.~\ref{thm:lp_above_two}, \ref{thm:above_two_main}(a).
   Lower bounds in Thm.~\ref{thm:above_two_main}(b), (c). \\
Exponential, $L\asymp r$
 & \multicolumn{2}{>{\raggedright\arraybackslash}p{0.74\linewidth}@{}}{Every
   non-even $p>0$: bounds within a factor $O((\sqrt r/\Lambda)^\eta)$
   (Lem.~\ref{lem:joint-adaptive-shifts}, Thm.~\ref{thm:high_accuracy_main}(a)).
   Tight for regular rows (Lem.~\ref{lem:regular-high-accuracy}).} \\
Beyond exponential, $L\ge C_pr\Lambda$
 & \multicolumn{2}{>{\raggedright\arraybackslash}p{0.74\linewidth}@{}}{Every
   non-even $p>0$: tight, $\Theta_p((\epsilon^{-2}r^{1/2-p}\Lambda)^\eta)$
   (Thm.~\ref{thm:high_accuracy_main}(b))} \\
Fixed $r$, $\epsilon\to0$
 & \multicolumn{2}{>{\raggedright\arraybackslash}p{0.74\linewidth}@{}}{Every
   non-even $p>0$: tight, $\Theta_{r,p}(\epsilon^{-2\eta})$
   (Thm.~\ref{thm:high_accuracy_main}(c))} \\
\bottomrule
\end{tabularx}
\caption{Results of Sections~\ref{sec:joint_accuracy},
\ref{sec:high_accuracy} and~\ref{sec:joint_above_two} by accuracy
regime, with $\Lambda:=\log(2r)$. For even $p\ge4$ the count is determined throughout
(Table~\ref{tab:tight_regimes}).}
\label{tab:joint_regimes}
\end{table}

\subsection{Further Results}

For every real $p>2$, real mutually unbiased bases give, for each rank, one input
that requires $\Omega_p(\min\{r^{p/2}/\epsilon,r^2\})$ rows at every
accuracy (Lemma~\ref{lem:lp_real_lower}). For non-even $p>2$ and
fixed $s\ge2$, rows built from binary quadratic forms require
$\Omega_{p,s}(\min\{r^{p/2}/\epsilon,r^s\})$ rows
(Lemma~\ref{lem:lp_noneven_lower}). Hsieh, Lee, Mohanty, Putterman,
and Zhang~\cite[Section~3]{hlmpz25} isolate individual constraints
using a Boolean AND implemented by noncommutative group operations,
whereas we isolate a row with a signed combination of real queries.
Below two, polynomial moments
save a power of the rank over $r/\epsilon^2$ for every fixed $0<p<2$
when $L\asymp r^{2/(p+2)}$ (Lemma~\ref{lem:joint-positive-caps}), and
adaptive shifts gain a factor $\log r$ in the lower bounds near
exponential accuracy (Lemma~\ref{lem:joint-adaptive-shifts}). Above
two, shared query coordinates and overlapping transverse groups
strengthen the large-rank
lower bound for non-even $p$ (Theorem~\ref{thm:above_two_main}(b)
and~(c)). We also give a short proof
of the rank-two rate $O_p(\epsilon^{-1/(p+1)})$ for $p>2$ by local
moment matching (Lemma~\ref{lem:lp_fixed_rank}).
Siegel~\cite[Proposition~1]{siegel25} used related local sampling
estimates for approximation by integer powers of the ReLU.

\subsection{Related Work}

Tables~\ref{tab:lp_below_one_history}, \ref{tab:lp_row_count_history},
and~\ref{tab:lp_above_two_history} list the previous bounds for
$0<p<1$, $1<p<2$, and $p>2$ with precise references and conditions,
writing $\Lambda:=\log(2r)$. For $0<p<1$, Zvavitch~\cite{zva00} and
Schechtman and Zvavitch~\cite{sz01} gave the first bounds linear in $r$
up to logarithmic factors. Yasuda~\cite[p.~202]{yasuda24thesis} records
an $\Omega(r^{p/2}/\epsilon)$ lower bound for sampling-based
embeddings, attributing it to~\cite{lww21}. Our
Lemma~\ref{lem:lp_noneven_lower} holds for every reweighting of
original rows and specifies the dependence on $r$ and $\epsilon$ for
each fixed $p$ and $s$.

\begin{table}[!ht]
\centering
\small
\begin{tabular}{@{}r l l@{}}
\toprule
Year & Reference & Row count \\
\midrule
\multicolumn{3}{@{}l}{Upper bounds} \\
2000 & \cite[Theorem~1]{zva00} & $O_p(r\Lambda H^2/\epsilon^2)$ \\
2001 & \cite[Theorem~1.1]{sz01} & $O_{p,\epsilon}(r\Lambda^3)$ \\
2022 & \cite[Theorem~2.11]{mmwy22} & $O_p(r\Lambda^2\log(2r/\epsilon)/\epsilon^2)$ \\
2023 & \cite[Theorem~A.2]{wy23} & $O_p(r\Lambda^3/\epsilon^2)$ \\
& Theorem~\ref{thm:lp_one_to_two_informal} & $O(r\epsilon^{-2}\min\{1/p,\Lambda\})$ \\
& Theorem~\ref{thm:joint_polynomial_informal} & $\Theta_{p,a}(r/\epsilon^2)$ for $\epsilon\ge r^{-a}$ \\
\midrule
\multicolumn{3}{@{}l}{Lower bound} \\
2014 & \cite[Theorem~2.2]{akw14} & $\Omega(r/\epsilon^2)$ \\
\bottomrule
\end{tabular}
\caption{Bounds for $0<p<1$ and $0<\epsilon\le1/2$, with
$\Lambda:=\log(2r)$ and $H:=\log\log(2r/\epsilon^2)+\log(2/\epsilon)$.
The 2001 bound is for fixed $\epsilon$. The 2022 bound uses
Remark~2.12 of~\cite{mmwy22}, and the 2022 and 2023 bounds are for
constant success probability. The cut-sparsification lower bound
uses incidence matrices and assumes $r\ge C/\epsilon^2$.}
\label{tab:lp_below_one_history}
\end{table}

\begin{table}[!ht]
\centering
\small
\begin{tabular}{@{}r l l@{}}
\toprule
Year & Reference & Row count \\
\midrule
\multicolumn{3}{@{}l}{Upper bounds} \\
1995 & \cite[Theorem~1.1]{tal95} & $O_{p,\epsilon}(r\Lambda\log^2(2\Lambda))$ \\
2015 & \cite[Theorem~7.1]{cp15} & $O_p(r\epsilon^{-2}\log(2r/\epsilon)\log^2(2\Lambda/\epsilon))$ \\
2026 & \cite[Corollary~1.4]{dk26} & $O(r\Lambda/\epsilon^2)$ \\
& Theorem~\ref{thm:lp_one_to_two_informal} & $O(r/\epsilon^2)$ \\
\midrule
\multicolumn{3}{@{}l}{Lower bounds} \\
2014 & \cite[Theorem~2.2]{akw14} & $\Omega(r/\epsilon^2)$ \\
& Theorem~\ref{thm:joint_polynomial_informal} & $\Omega_{p,a}(r/\epsilon^2)$ for $\epsilon\ge r^{-a}$ \\
\bottomrule
\end{tabular}
\caption{Bounds for $1<p<2$ and $0<\epsilon\le1/2$, with
$\Lambda:=\log(2r)$. The 1995 bound is for fixed $\epsilon$.
The cut-sparsification lower bound assumes $r\ge C/\epsilon^2$.}
\label{tab:lp_row_count_history}
\end{table}

\begin{table}[!ht]
\centering
\small
\begin{tabular}{@{}r l l@{}}
\toprule
Year & Reference & Row count \\
\midrule
\multicolumn{3}{@{}l}{Upper bounds} \\
1989 & \cite[Theorem~7.3]{blm89} & $O_{p,\epsilon}(r^{p/2}\Lambda)$ \\
2011 & \cite{sch11} & $O_p(r^{p/2}/\epsilon^2)$, even $p$ \\
2015 & \cite[Theorem~7.1]{cp15} & $O_p(r^{p/2}\Lambda\log(2/\epsilon)/\epsilon^5)$ \\
2021 & \cite[Corollary~4.5]{kos21} & $O_p(r^{p/2}\epsilon^{-2}\log(2r/\epsilon))$ \\
2026 & \cite[Theorem~1.2]{li26fixed} & $O_{r,p}(\epsilon^{-2\eta})$, fixed $r,p$ \\
& Theorem~\ref{thm:lp_above_two_informal} & $O_p(r^{p/2}\epsilon^{-2}(1+\Lambda/t^{p-2}))$ \\
& Theorem~\ref{thm:lp_above_two} & $O_p(\min\{r^{p/2}/\epsilon^2,r^{(p+1)/2}/\epsilon,r^p\})$, even $p$ \\
\midrule
\multicolumn{3}{@{}l}{Lower bounds} \\
1977 & \cite[Theorem~1.1]{bdgjn77} & $\Omega_p(r^{p/2})$ \\
2014 & \cite[Theorem~2.2]{akw14} & $\Omega(\min\{r/\epsilon^2,r^2\})$ \\
2021 & \cite[Theorem~6.1]{lww21} & $\widetilde{\Omega}_p(r/\epsilon^2)$, non-even $p$ \\
2023 & \cite[Eq.~(19)]{llw23} & $\Omega_{r,p}(\epsilon^{-2\eta})$, non-even $p\ge1$ \\
2026 & \cite[Theorem~1.2]{li26} & $\widetilde{\Omega}_p(r/\epsilon^2)$, non-even $p$ \\
& Lemma~\ref{lem:lp_noneven_lower} & $\Omega_{p,s}(\min\{r^{p/2}/\epsilon,r^s\})$, non-even $p$ \\
& Theorem~\ref{thm:lp_noneven_amplified_informal} & $\widetilde{\Omega}_p(r^{p/2}/\epsilon^2)$, non-even $p$ \\
& Theorem~\ref{thm:lp_even_lower_informal} & $\Omega_p(\min\{r^{p/2}/\epsilon^2,r^{(p+1)/2}/\epsilon,r^p\})$, even $p$ \\
\bottomrule
\end{tabular}
\caption{Bounds for $p>2$, $r\ge2$, and $0<\epsilon\le1/2$, with
$\Lambda:=\log(2r)$, $t:=\log(2\sqrt\Lambda/\epsilon)$, and
$\eta:=(r-1)/(r+2p)$. The 1989 bound is for fixed $\epsilon$.
The 2021 lower bound assumes $r\ge C\log(1/\epsilon)$ and the tilde
hides a polylogarithmic factor in $1/\epsilon$. The 2023 lower bound
is for fixed $r,p$ as $\epsilon$ tends to zero. Li's 2026 lower bound
assumes $r\ge C_p\log(1/\epsilon)$ and the tilde hides a
polylogarithmic factor in $r/\epsilon$. These two lower bounds also
apply to arbitrary linear embeddings. Our non-even lower bound
with $\epsilon^{-2}$ dependence assumes $r\ge C_p\log(2/\epsilon)$
and hides a polylogarithmic factor in $1/\epsilon$. Our other
lower bounds hold throughout the displayed accuracy range, with
$s\ge2$ an arbitrary fixed integer.}
\label{tab:lp_above_two_history}
\end{table}

\subsection{Roadmap}

Section~\ref{sec:technical_overview} outlines the proofs and
Section~\ref{sec:preliminaries} collects the tools.
Section~\ref{sec:general_lp_field} proves the upper bound for
$0<p\le2$, Section~\ref{sec:small_p_uniform} removes the factor $1/p$
for small $p$, and Section~\ref{sec:lp_above_two} proves the upper
bounds for $p>2$. Section~\ref{sec:lp_even_lower} proves the lower
bounds for $p>2$. Sections~\ref{sec:joint_accuracy}, \ref{sec:high_accuracy}
and~\ref{sec:joint_above_two} treat the joint
dependence on rank and accuracy, including the optimality below two
and the high-accuracy bounds, and Appendix~\ref{app:common_tools}
proves their common tools.

%% file: 05_overview.tex
\section{Technical Overview}\label{sec:technical_overview}

\subsection{\texorpdfstring{A Gaussian Cone for $0<p\le2$ (Section~\ref{sec:general_lp_field})}{A Gaussian Cone for Exponents at Most Two}}

Reis and Rothvoss~\cite{rr26} obtain linear-size sparsifiers at $p=1$
by partial coloring: at each step they find a signed change in the row
weights that sets a constant fraction of the weights to zero and
changes every energy value by a small relative amount. The geometric
input is a lower bound on the Gaussian measure of the cone on which
the signed energy is nonnegative for every query. For $1<p<2$ the previous
bound of Dai and Kosov~\cite{dk26} carries a factor $\log(2r)$.
Independent random sampling of rows needs such a factor even at
$p=2$. We prove the cone bound for every $0<p\le2$ in Lewis coordinates by
constructing a Gaussian trigonometric field. At a global maximum of a
confined and tilted field, integrating the centered second difference
against $t^{-1-p}$ gives an element of a shifted cone. We vary the tilt
and use the area formula to bound the probability of the cone from
below. This gives a round with error of order
$\sqrt{r/(N\min\{p,1\})}$ on $N$ current rows, and the errors of the
rounds form a geometric series.

\subsection{\texorpdfstring{Uniformity as $p$ Tends to Zero (Section~\ref{sec:small_p_uniform})}{Uniformity as p Tends to Zero}}

The preceding round pays a factor $p^{-1/2}$ in its error. For
$p\lesssim 1/\Lambda$, we instead normalize every energy vector to have
total mass one and split its coordinates into dyadic ramps. We retain
high-sensitivity rows and protect sets of coordinates that are
fat-shattered by the ramp classes. To bound the total protected set,
we turn the fat-shattering witnesses back into linear-form witnesses.
Taking $p$th roots makes the interference between most pairs of blocks
small. If too many coordinates were protected, then adding the
witnesses would realize more sign patterns than a rank-$r$ hyperplane
arrangement allows. On the remaining rows we use random signs and
control the discrepancy uniformly over the ramp classes.

\subsection{\texorpdfstring{High Accuracy for $p>2$ (Section~\ref{sec:lp_above_two})}{High Accuracy for p>2}}

For even $p$, tensor lifting reduces the problem to spectral
sparsification~\cite{sch11}. For other $p>2$ the previous bound carries
a factor $\log(2r/\epsilon)$~\cite{kos21}. We write the
row energy as a polynomial in the normalized inner product plus a
remainder. In the partial-coloring step, we impose linear constraints
that cancel the leading directions of the polynomial feature
covariance. The remainder controls the variance, and a local chaining
estimate controls the oscillation within each member of a cover. We
choose the degree so that both parts fit within the error budget as
the support decreases. When $p$ is even, the polynomial representation
is exact. The projection then controls the pointwise variance of the
whole Gaussian process, and polarization lets us estimate its supremum
using a constant-radius net of the $\ell_p$ unit ball. This gives a
per-round error proportional to $r^{(p+1)/2}/N$. For fixed rank, we
also use moment matching directly on small cells of the sphere,
preserving the local Taylor terms exactly and averaging the remaining
errors across cells.

\subsection{\texorpdfstring{Lower Bounds for $p>2$ (Section~\ref{sec:lp_even_lower})}{Lower Bounds for p>2}}

All our row lower bounds test the energy at a few structured queries
and show that a sparsifier with small support violates one of them. We
use four constructions.
\begin{itemize}[nosep]
\item \emph{A row-kernel criterion} (Lemmas~\ref{lem:lp_real_lower}
 and~\ref{lem:lp_even_lower}). We test the energy on the rows of the
 input matrix. If the Gram matrix of the kernel $|a_i^\top a_j|^p$ is
 bounded below after its constant component is removed, these tests
 limit how much mass a sparsifier can put on a small support. We
 obtain this condition from real mutually unbiased bases and, for even
 $p$, from random sign rows whose tensor powers are nearly orthogonal.
 In the even case, polarization also turns a uniform error bound into
 a bound on the distance between the input and output weight vectors.
\item \emph{Isolating one row} (Lemma~\ref{lem:lp_noneven_lower}).
 For non-even $p$, we construct rows from binary quadratic forms, so
 that the off-diagonal inner products take values in a fixed finite
 set after normalization. We shift queries in a coordinate shared by
 all rows. A signed combination of these queries cancels the analytic
 off-diagonal terms and retains the contribution of one row.
\item \emph{Tensoring hard instances} (Theorem~\ref{thm:lp_lower_main}(b)).
 We tensor a low-rank hard instance with nearly orthogonal directions.
 Product queries give a diagonally dominant system in the errors of the
 individual copies. Hence approximation of the full instance forces
 approximation within each copy, and we can add their lower bounds on
 the support size.
\item \emph{Mixed coefficients for even $p$} (Theorem~\ref{thm:lp_lower_main}(a)).
 We reuse coordinates across groups and extract a mixed coefficient
 with one query variable per block. We recover the error tensor of each
 group in Frobenius norm using a lower bound on the Gram matrix of the
 first block. Independent orthogonal frames keep the original energy
 bounded on these queries, and their product basis makes each omitted
 row contribute one to the squared error. Choosing unequal block
 dimensions gives the intermediate term $r^{(p+1)/2}/\epsilon$.
\end{itemize}

\subsection{\texorpdfstring{Joint Rank and Accuracy (Sections~\ref{sec:joint_accuracy}--\ref{sec:joint_above_two})}{Joint Rank and Accuracy}}

We outline the main ideas behind each group of results.
\begin{itemize}[nosep]
\item \emph{Optimality below two} (Theorem~\ref{thm:below_two_main}(a)).
 For a partite dictionary of sign rows, Boolean block flips extract the
 top coefficient of the query energy, and a small support leaves a
 large discrepancy across the slices of this coefficient tensor.
\item \emph{Lower bounds at higher accuracy}
 (Theorem~\ref{thm:below_two_main}(b)
 and Lemma~\ref{lem:joint-adaptive-shifts}). We pack positive input weights
 using the distance induced by adaptive query differences. Sharing the
 first query block lets several transverse dictionaries force retained
 rows with a common rank budget.
\item \emph{High-accuracy upper bounds}
 (Theorem~\ref{thm:high_accuracy_main}(a)
 and~(b)). Polynomial moment
 constraints leave a remainder that can be integrated against threshold
 indicators. The query sphere has dimension $r-1$, which controls both
 the number of crossings in a paired row update and the metric entropy
 of the integrated discrepancy. Below two, a cover of low-energy
 queries lets us apply the same argument to localized approximants.
\item \emph{High-accuracy lower bounds}
 (Theorem~\ref{thm:high_accuracy_main}(a)). The hard input has one tag
 coordinate shared by $\Xi=e^{\Theta(r)}$ Gaussian rows $(1,u_t)$, tested
 by shifted scalar functions along random query directions. Two
 independent tests have squared correlation of order $1/\Xi$. Hence a
 test thresholded at Gaussian height of order $\sqrt{r\log r}$
 separates a Gaussian direction in the space of masses with constant
 probability, by the Paley--Zygmund inequality, and Gaussian packing
 bounds the number of supports. The shift spacing $\Xi^{-1/(r-1)}$
 produces the accuracy exponent $\eta$.
\end{itemize}

%% file: 10_preli.tex
\section{Preliminaries}\label{sec:preliminaries}

\subsection{Notation}

We write $B_2^d$ for the Euclidean unit ball in $\R^d$ and
$[m]:=\{1,\ldots,m\}$. Logarithms are natural unless a base is
specified. Constants with subscripts may depend on the indicated
parameters. Unsubscripted constants are absolute.
We write $S^{d-1}$ for the unit sphere of $\R^d$.
For a real random variable $X$ define
$\|X\|_{\psi_2}:=\inf\{t>0:\E[\exp(X^2/t^2)]\le2\}$. For a metric space
$(T,\rho)$ and $u>0$, $N(T,\rho,u)$ is the least number of balls of
radius $u$ that cover $T$.

In Sections~\ref{sec:joint_accuracy}, \ref{sec:high_accuracy}
and~\ref{sec:joint_above_two} and in Appendix~\ref{app:common_tools},
$0<\epsilon\le1/2$, $L:=\log(2/\epsilon)$,
$\Lambda:=\log(2r)$, and $\eta:=(r-1)/(r+2p)$. At subexponential accuracy
we also use
\[
 K:=\frac{L}{\log(er/L)} .
\]
Bounds with a subscript $p$ treat $p$ as fixed.

\subsection{Probability Inequalities}

\begin{lemma}[Markov and Chebyshev inequalities]\label{lem:pre_markov}
Let $X\ge0$ be a random variable and let $t>0$. Then
$\Pr[X\ge t]\le\E[X]/t$. Consequently, for a real random variable $Y$
and $t>0$, $\Pr[|Y-\E[Y]|\ge t]\le\operatorname{Var}(Y)/t^2$, and
$\Pr[Y\ge t]\le e^{-\lambda t}\E[e^{\lambda Y}]$ for every $\lambda>0$.
\end{lemma}
\begin{proof}
The first bound follows by taking expectations in
$t\mathbf1_{\{X\ge t\}}\le X$. The other two are the first bound applied
to $(Y-\E[Y])^2$ and to $e^{\lambda Y}$.
\end{proof}

\begin{lemma}[Chernoff bound for a Gaussian variable]\label{lem:pre_gaussian_chernoff}
Let $G$ be a standard Gaussian variable and let $t\ge0$. Then
$\Pr[G>t]\le e^{-t^2/2}$ and $\Pr[|G|>t]\le2e^{-t^2/2}$.
\end{lemma}
\begin{proof}
For $t>0$, Lemma~\ref{lem:pre_markov} with $\lambda:=t$ and
$\E[e^{tG}]=e^{t^2/2}$ give $\Pr[G>t]\le e^{-t^2}e^{t^2/2}$. The case
$t=0$ is trivial, and the second bound follows by symmetry.
\end{proof}

\begin{lemma}[Gaussian tails {\cite[Proposition~2.1.2]{ver18}}]\label{lem:pre_gaussian_tail}
Let $G$ be a standard Gaussian variable with density $\phi$. For every
$t>0$,
\[
 (1/t-1/t^3)\phi(t)\le\Pr[G>t]\le\phi(t)/t .
\]
\end{lemma}

\begin{lemma}[Hoeffding's inequality {\cite[Theorem~2.2.6]{ver18}}]\label{lem:pre_hoeffding}
Let $X_1,\ldots,X_M$ be independent random variables with
$X_l\in[a_l,b_l]$. For every $t>0$,
\[
 \Pr[\sum_{l=1}^M(X_l-\E[X_l])\ge t]
 \le\exp(-2t^2/\sum_{l=1}^M(b_l-a_l)^2).
\]
\end{lemma}

\begin{lemma}[Bernstein's inequality {\cite[Theorem~2.8.4]{ver18}}]\label{lem:pre_bernstein}
Let $X_1,\ldots,X_M$ be independent mean-zero random variables with
$|X_l|\le K$, and let $\sigma^2:=\sum_{l=1}^M\E[X_l^2]$. For every $t>0$,
\[
 \Pr[|\sum_{l=1}^MX_l|\ge t]\le2\exp(-\frac{t^2/2}{\sigma^2+Kt/3}).
\]
\end{lemma}

\begin{lemma}[Paley--Zygmund inequality {\cite{pz32}}]\label{lem:pre_paley_zygmund}
Let $Z\ge0$ be a random variable with $0<\E[Z^2]<\infty$ and let
$0\le\theta\le1$. Then
$\Pr[Z\ge\theta\E[Z]]\ge(1-\theta)^2(\E[Z])^2/\E[Z^2]$.
\end{lemma}

\begin{lemma}[\v{S}id\'ak's inequality {\cite{sid67}}]\label{lem:pre_sidak}
Let $X_1,\ldots,X_n$ be centered jointly Gaussian random variables and
let $c_1,\ldots,c_n\ge0$. Then
\[
 \Pr[|X_i|\le c_i\text{ for all }i]\ge\prod_{i=1}^n\Pr[|X_i|\le c_i].
\]
\end{lemma}

\begin{lemma}[Gaussian concentration {\cite[Theorem~5.2.2]{ver18}}]\label{lem:pre_gaussian_concentration}
There is an absolute constant $c_G>0$ such that, for every function $f$ on $\R^n$ that is $R$-Lipschitz in the Euclidean
norm, a standard Gaussian vector $g'$ in $\R^n$ and
every $t>0$,
\[
 \Pr[|f(g')-\E[f(g')]|\ge t]\le2\exp(-c_Gt^2/R^2).
\]
\end{lemma}

\begin{lemma}[Sub-gaussian variables and processes {\cite[Example~2.5.8, Propositions~2.5.2 and~2.6.1, and Theorem~8.1.3]{ver18}}]\label{lem:pre_subgaussian}
There is an absolute constant $C$ with the following properties.
\begin{enumerate}[label=(\alph*),nosep]
\item If $X$ is a centered random variable with $\|X\|_{\psi_2}\le K$,
then $\E[e^{\lambda X}]\le e^{C\lambda^2K^2}$ for every real $\lambda$.
\item If $\sigma_1,\ldots,\sigma_M$ are independent uniform signs and
$a\in\R^M$, then $\|\sum_{l=1}^Ma_l\sigma_l\|_{\psi_2}\le C\|a\|_2$.
\item If $(X_t)_{t\in T}$ is a centered process on a finite metric space
$(T,\rho)$ with $\|X_t-X_s\|_{\psi_2}\le K\rho(t,s)$ for all $s,t\in T$,
then
\[
 \E[\sup_{t\in T}X_t]\le CK\int_0^\infty\sqrt{\log N(T,\rho,u)}\,\d u .
\]
\end{enumerate}
\end{lemma}

\begin{lemma}[Symmetrization {\cite[Lemma~6.4.2]{ver18}}]\label{lem:pre_symmetrization}
Let $X_1,\ldots,X_M$ be independent mean-zero random vectors in a
finite-dimensional normed space, and let $\sigma_1,\ldots,\sigma_M$ be
independent uniform signs that are independent of them. Then
\[
 \E[\|\sum_{l=1}^MX_l\|]\le2\E[\|\sum_{l=1}^M\sigma_lX_l\|].
\]
\end{lemma}

\subsection{Covering and Chaining}

\begin{lemma}[Nets and separated sets {\cite[Corollary~4.2.13]{ver18}}]\label{lem:pre_nets}
Let $\|\cdot\|$ be a norm on a real vector space of dimension $d\ge1$
with unit ball $B$, and let $\rho>0$ and $R>0$. Every subset of $RB$
whose points are pairwise at distance at least $\rho$ has at most
$(1+2R/\rho)^d$ points, and every point of $RB$ is within $\rho$ of a
maximal such subset. The same holds with $RB$ replaced by the sphere
$\{x:\|x\|=R\}$. For $R=1$ and $\rho\le1$ the bound is at most
$(3/\rho)^d$.
\end{lemma}

\begin{lemma}[Maurey's lemma]\label{lem:pre_maurey}
Let $y=\sum_j\lambda_jy_j$ be a finite combination of points of a
Euclidean space with $\lambda_j\ge0$, $\sum_j\lambda_j\le1$ and
$\|y_j\|_2\le1$, and let $k\ge1$ be an integer. Then some mean of $k$
points of $\{0\}\cup\{y_j\}$ is within $k^{-1/2}$ of $y$.
\end{lemma}
\begin{proof}
Let $Y$ take the value $y_j$ with probability $\lambda_j$ and the value
$0$ otherwise. The mean of $k$ independent copies of $Y$ has expected
squared distance $\E[\|Y-\E[Y]\|_2^2]/k\le1/k$ from $y$.
\end{proof}

For $U\subseteq\mathbb R^I$ and $\Delta>0$, a set $Q\subseteq I$ is
$\Delta$-fat-shattered if there are thresholds $\gamma_j$, $j\in Q$,
such that every $z\in\{0,1\}^Q$ has a witness $u\in U$ satisfying
$u_j\ge\gamma_j+\Delta$ when $z_j=1$ and
$u_j\le\gamma_j-\Delta$ when $z_j=0$.
Write $\operatorname{fat}_\Delta(U)$ for the maximum cardinality of such
a set. A proper Euclidean $h$-net of $U$ is a finite subset
$V\subseteq U$ such that every $u\in U$ has $\|u-v\|_2\le h$
for some $v\in V$.

\begin{lemma}[Euclidean covering {\cite[Theorem~1]{mv03}}]\label{lem:small_p_cover}
There are absolute constants $c_0\in(0,1)$ and $C_0\ge2$ with the
following property. For nonempty $U\subseteq[0,a]^m$ and
$0<h<a\sqrt m$, there is a proper Euclidean $h$-net $V$ satisfying
$\log|V|\le C_0(1+\operatorname{fat}_{c_0h/\sqrt m}(U))
\log(C_0a\sqrt m/h)$.
\end{lemma}

For a pseudometric $\rho$ on $T$, write $e_k(T,\rho)$ for the infimum of
the radii of covers using centers in $T$. We allow at most $2^{2^k}$
balls for $k\ge1$ and one ball for $k=0$. For $\alpha>0$ and $r\ge1$,
we use the chaining functional
\begin{equation}\label{eq:chaining_functional_definition}
 \gamma_{\alpha,r}(T,\rho)
 :=\inf_{(\mathcal A_k)}\sup_{t\in T}
 (\sum_{k\ge0}2^{kr/\alpha}
 \operatorname{diam}_{\rho}(\mathcal A_k(t))^r)^{1/r},
\end{equation}
where the infimum ranges over nested partitions with $|\mathcal A_0|=1$ and $|\mathcal A_k|\le2^{2^k}$, and $\mathcal A_k(t)$ is the cell containing $t$.

\begin{lemma}[Gaussian chaining {\cite[Theorems~2.3 and~2.6, Lemma~2.7]{kos21}}]\label{lem:kosov_gaussian_chaining}
Let $p>2$, let $\nu$ be a probability measure on a finite set, and let $U$ be a centered Euclidean ball in a finite-dimensional subspace of $L_{2p}(\nu)$. Set
\[
 \rho(f,g):=\|(f-g)(|f|^{p-1}+|g|^{p-1})\|_{L_2(\nu)},
 \qquad \|v\|_h:=\|v|h|^{p-1}\|_{L_2(\nu)}.
\]
For $\eta>0$, write $\|v\|_{h,\eta}:=(\|v\|_h^2+\eta\|v\|_{L_2(\nu)}^2)^{1/2}$.
For every nonempty compact $T\subseteq U$, every continuous centered Gaussian process $(Z_f)_{f\in T}$ satisfying $(\E[|Z_f-Z_g|^2])^{1/2}\le\rho(f,g)$, and every $f_0\in T$,
\[
 \E[\sup_{f\in T}|Z_f-Z_{f_0}|]
 \le C_p\sup_{h\in T}(\sum_{k\ge0}2^k e_k(U,\|\cdot\|_{h,\eta})^2)^{1/2}
 +C_p\gamma_{2p,p}(T,L_{2p}(\nu))^p.
\]
Here $C_p$ depends only on $p$.
\end{lemma}

To apply Kosov's theorem, use the Euclidean ball as the
$2$-convex body and $L_{2p}(\nu)$ as the auxiliary metric.
Lemma~2.7 of~\cite{kos21} gives the hypotheses on $\rho$ with
the weighted seminorms $\|\cdot\|_h$. The norms $\|\cdot\|_{h,\eta}$
majorize those seminorms, so they satisfy the same upper bound on
increments with constants independent of $\eta$.

\begin{lemma}[Euclidean entropy bound {\cite[proof of Lemma~3.4(1)]{kos21}}]\label{lem:kosov_euclidean_entropy}
Let $Q$ be a positive semidefinite $d\times d$ matrix, let $R>0$, and set $\|x\|_Q:=(x^\top Qx)^{1/2}$. Then
\[
 (\sum_{k\ge0}2^k e_k(RB_2^d,\|\cdot\|_Q)^2)^{1/2}
 \le CR\sqrt{\tr[Q]},
\]
where $C$ is absolute.
\end{lemma}

\subsection{Analysis and Linear Algebra}

\begin{lemma}[Minkowski's inequality {\cite[Chapter~6]{fol99}}]\label{lem:pre_minkowski}
Let $1\le p<\infty$ and let $\mu$ and $\nu$ be $\sigma$-finite measures.
\begin{enumerate}[label=(\alph*),nosep]
\item For measurable functions $f$ and $g$,
$\|f+g\|_{L_p(\mu)}\le\|f\|_{L_p(\mu)}+\|g\|_{L_p(\mu)}$.
\item For a measurable function $F$ on the product space such that the
right side is finite,
\[
 \|\int F(\cdot,y)\,\d\nu(y)\|_{L_p(\mu)}\le
 \int\|F(\cdot,y)\|_{L_p(\mu)}\,\d\nu(y).
\]
\end{enumerate}
\end{lemma}

\begin{lemma}[Bessel's inequality and Parseval's identity {\cite[Chapter~5]{fol99}}]\label{lem:pre_parseval}
Let $(e_n)$ be an orthonormal family in a Hilbert space. Then
$\sum_n|\langle f,e_n\rangle|^2\le\|f\|^2$ for every $f$, with equality
when $(e_n)$ is an orthonormal basis.
\end{lemma}

\begin{lemma}[Stone--Weierstrass theorem {\cite[Chapter~4]{fol99}}]\label{lem:pre_stone_weierstrass}
Let $X$ be a compact Hausdorff space and let $\mathcal A$ be a subalgebra
of the real continuous functions on $X$ that contains the constants and
separates the points of $X$. Then $\mathcal A$ is dense in the continuous
functions on $X$ in the uniform norm.
\end{lemma}

\begin{lemma}[Ratio of gamma functions {\cite{te51}}]\label{lem:pre_gamma_ratio}
For fixed real $a$ and $b$,
$\Gamma(x+a)/\Gamma(x+b)=x^{a-b}(1+O(1/x))$ as $x\to\infty$.
\end{lemma}

\begin{lemma}[Area formula {\cite[Theorem~3.2.3]{fed96}}]\label{lem:federer_area_formula}
Let $F:\mathbb R^d\to\mathbb R^d$ be continuously differentiable and let $E\subseteq\mathbb R^d$ be Borel. Then
\[
 \int_E|\det\d F(x)|\d x
 =\int_{\mathbb R^d}\#\{x\in E:F(x)=y\}\d y.
\]
In particular, $\operatorname{vol}(F(E))\le\int_E|\det\d F(x)|\d x$.
\end{lemma}

\begin{lemma}[Gershgorin's theorem {\cite[Theorem~6.1.1]{hj13}}]\label{lem:pre_gershgorin}
Every eigenvalue of a complex $n\times n$ matrix $A$ lies in one of the
discs $\{z:|z-A_{ii}|\le\sum_{j\ne i}|A_{ij}|\}$, $i\in[n]$.
\end{lemma}

\begin{lemma}[Partition of an isotropic set {\cite[Corollary~1.5]{mss15}}]\label{lem:pre_mss}
Let $k\ge2$ be an integer and let $u_1,\ldots,u_m\in\R^d$ satisfy
$\sum_{i=1}^mu_iu_i^\top=I_d$ and $\|u_i\|_2^2\le\delta$ for all $i$. Then
$[m]$ can be partitioned into sets $S_1,\ldots,S_k$ with
$\|\sum_{i\in S_j}u_iu_i^\top\|_{\rm op}\le(1/\sqrt k+\sqrt\delta)^2$ for
every $j$.
\end{lemma}

\begin{lemma}[Johnson--Lindenstrauss lemma {\cite[Theorem~2.1]{dg03}}]\label{lem:pre_jl}
Let $0<\varepsilon<1$ and let $X$ be a set of $n\ge2$ points of a Euclidean
space. For every integer $k\ge24\varepsilon^{-2}\log n$ there is a linear
map $f$ into $\R^k$ with
$(1-\varepsilon)\|u-v\|_2^2\le\|f(u)-f(v)\|_2^2\le(1+\varepsilon)\|u-v\|_2^2$
for all $u,v\in X$.
\end{lemma}

\begin{lemma}[Random sign tensors {\cite[Theorem~1.2(2a)]{knh24}}]\label{lem:random_sign_tensor_gap}
Fix an integer $k\ge1$. Define $M:=\lfloor r^k/(4k!)\rfloor$, and let $\xi_1,\ldots,\xi_M$ be independent uniform vectors in $\{-1,1\}^r$. Let $X$ be the matrix whose columns are indexed by the $k$-element subsets $S\subseteq\{1,\ldots,r\}$ and whose entries are
\[
 X_{i,S}:=\prod_{j\in S}\xi_{ij}.
\]
Then
\[
 \Pr[\frac{k!}{r^k}XX^\top\succeq\tfrac18 I_M]\longrightarrow1
 \qquad\text{as }r\longrightarrow\infty.
\]
\end{lemma}

\subsection{Convex Geometry and Spherical Harmonics}

\begin{lemma}[Brunn--Minkowski inequality {\cite[Lecture~5]{bal97}}]\label{lem:pre_brunn_minkowski}
Let $A$ and $B$ be nonempty compact subsets of $\R^n$ and let
$0\le\lambda\le1$. Then
\[
 \operatorname{vol}(\lambda A+(1-\lambda)B)^{1/n}
 \ge\lambda\operatorname{vol}(A)^{1/n}+(1-\lambda)\operatorname{vol}(B)^{1/n}.
\]
\end{lemma}

\begin{lemma}[Geometric Brascamp--Lieb inequality {\cite[Theorem~6.4]{bal97}}]\label{lem:pre_brascamp_lieb}
Let $b_1,\ldots,b_m$ be unit vectors of $\R^s$ and let $c_i>0$ satisfy
$\sum_ic_ib_ib_i^\top=I_s$. For nonnegative measurable functions
$f_1,\ldots,f_m$ on $\R$,
\[
 \int_{\R^s}\prod_{i=1}^mf_i(\langle b_i,y\rangle)^{c_i}\,\d y\le
 \prod_{i=1}^m(\int_\R f_i)^{c_i}.
\]
\end{lemma}

\begin{lemma}[Log-concave marginals {\cite[Lecture~5]{bal97}}]\label{lem:pre_logconcave_marginal}
Let $U$ be uniformly distributed on a convex body in $\R^n$ and let
$v\in\R^n$ be nonzero. Then $\langle v,U\rangle$ has a log-concave
density.
\end{lemma}

\begin{lemma}[Busemann's intersection inequality {\cite{bus53}, \cite[Section~9.4]{gar06}}]\label{lem:pre_busemann}
Let $r\ge2$, let $\omega$ be the normalized surface measure on
$S^{r-1}$, and let $\mathsf K:=\{\varrho v:v\in S^{r-1},\ 0\le\varrho\le
\varrho_0(v)\}$ for a continuous positive function $\varrho_0$ on
$S^{r-1}$. Then
\[
 \int_{S^{r-1}}\operatorname{vol}_{r-1}(\mathsf K\cap x^\perp)^r\,\d\omega(x)
 \le\frac{\operatorname{vol}_{r-1}(B_2^{r-1})^r}
 {\operatorname{vol}_r(B_2^r)^{r-1}}\operatorname{vol}_r(\mathsf K)^{r-1}.
\]
\end{lemma}

\begin{lemma}[Spherical harmonics {\cite[Theorems~1.1.2, 1.1.3 and~1.2.9 and Corollary~1.1.5]{dx13}}]\label{lem:pre_harmonics}
Let $r\ge2$, let $\omega$ be the normalized surface measure on
$S^{r-1}$, and let $\mathcal H_n$ be the space of restrictions to
$S^{r-1}$ of the homogeneous harmonic polynomials of degree $n$ on
$\R^r$.
\begin{enumerate}[label=(\alph*),nosep]
\item The spaces $\mathcal H_n$ are pairwise orthogonal in $L_2(\omega)$,
and the restriction to $S^{r-1}$ of every polynomial of degree at most
$n$ is a sum of elements of $\mathcal H_0,\ldots,\mathcal H_n$.
\item For every bounded Borel function $\psi_0$ on $[-1,1]$ and every
$n\ge0$ there is a scalar $\lambda_n(\psi_0)$ such that
$\int\psi_0(\langle u,v\rangle)Y(v)\,\d\omega(v)=\lambda_n(\psi_0)Y(u)$ for
all $Y\in\mathcal H_n$ and $u\in S^{r-1}$.
\end{enumerate}
\end{lemma}

\begin{lemma}[Multipliers of the $p$-cosine transform {\cite[Equations~(2.5) and~(2.6)]{opr13}}]\label{lem:pre_cosine_multipliers}
Let $r\ge2$, let $p>0$ be non-even, and let $\lambda_n$ be the scalar
$\lambda_n(\psi_0)$ of Lemma~\ref{lem:pre_harmonics}(b) for
$\psi_0(t):=|t|^p$. There is a number $c_{r,p}\ne0$ such that, for every
integer $j\ge0$,
\[
 \lambda_{2j}=c_{r,p}(-1)^j\frac{\Gamma(j-p/2)}{\Gamma(j+(r+p)/2)}.
\]
\end{lemma}

\subsection{Combinatorial Bounds}

\begin{lemma}[Hyperplane region count {\cite[Proposition~2.4]{sta07}}]\label{lem:small_p_regions}
An arrangement of $q$ affine hyperplanes in $\mathbb R^r$ has at most
$\sum_{j=0}^r\binom qj$ open regions. If $1\le r\le q$, this number is
at most $(eq/r)^r$.
\end{lemma}

\begin{lemma}[Warren's theorem {\cite[Theorem~3]{war68}}]\label{lem:pre_warren}
Let $G_1,\ldots,G_m$ be real polynomials of degree at most $k\ge1$ in $d$
variables, with $m\ge d$. Then $\R^d\setminus\bigcup_l\{G_l=0\}$ has at
most $(4ekm/d)^d$ connected components.
\end{lemma}

\begin{lemma}[Real mutually unbiased bases {\cite[Lemma~3.3, Theorem~3.4, and Example~9.2]{ccks97}}]\label{lem:pre_mub}
Let $q\ge4$ be a power of four. There are $q/2+1$ orthonormal bases of
$\R^q$ such that any two vectors from distinct bases have inner product
of magnitude $q^{-1/2}$.
\end{lemma}

\subsection{Partial Coloring}

\begin{lemma}[Partial coloring from volume {\cite[Theorem~7]{rr26}}]\label{lem:volume_partial_coloring}
For every $c>0$ there is a constant $C>0$ depending only on $c$ such that any centrally symmetric convex body $K\subseteq[-1,1]^N$ with $\operatorname{vol}(K)\ge c^N$ admits a point $z\in CK\cap[-1,1]^N$ satisfying
\[
 |\{i\in\{1,\ldots,N\}:|z_i|=1\}|\ge N/2.
\]
\end{lemma}

\begin{lemma}[Partial coloring in a subspace {\cite[Theorem~8]{rot16}}]\label{lem:subspace_partial_coloring}
There are absolute constants $c_0,c_1>0$ and $N_0$ with the following property. If $N\ge N_0$, $H\subseteq\mathbb R^N$ is a subspace of codimension at most $c_0N$, and $K\subseteq H$ is a centrally symmetric convex body with standard Gaussian measure $\gamma_H(K)\ge1/2$ in $H$, then some $z\in K\cap[-1,1]^N$ satisfies
\[
 |\{i\in\{1,\ldots,N\}:|z_i|=1\}|\ge c_1N.
\]
\end{lemma}

\begin{lemma}[Partial coloring in a subspace at small measure {\cite[Theorem~8]{rot16}}]\label{lem:pre_rothvoss}
Let $0<\varepsilon\le1/60000$ and
$\delta:=\frac32\varepsilon\log_2(1/\varepsilon)$. There is an integer
$n_{\rm c}$ such that the following holds for $n\ge n_{\rm c}$. If
$H\subseteq\R^n$ is a subspace defined by at most $2\delta n$ linear
equations and $\mathcal K'\subseteq H$ is a symmetric convex body in $H$
with $\gamma_H(\mathcal K')\ge e^{-\delta n}$, then some point of
$\mathcal K'\cap[-1,1]^n$ has at least $\varepsilon n$ coordinates in
$\{-1,1\}$.
\end{lemma}

\subsection{Sparsification}

\begin{lemma}[Even-power lifting {\cite{sch11}}]\label{lem:even_power_lifting}
Let $a_1,\ldots,a_m\in\mathbb R^n$ span a space of dimension $d\ge1$, and let $k\ge1$ be an integer. Equip $(\mathbb R^n)^{\otimes k}$ with its Euclidean tensor inner product and define $b_i:=a_i^{\otimes k}$. Then
\begin{align*}
 \dim\operatorname{span}\{b_1,\ldots,b_m\}
 \le & ~ \binom{d+k-1}{k}\le d^k,\\
 \langle b_i,x^{\otimes k}\rangle^2
 = & ~ |a_i^\top x|^{2k}\qquad(x\in\R^n).
\end{align*}
\end{lemma}
\begin{proof}
The tensors $b_i$ are symmetric tensors over the $d$-dimensional
row span. The symmetric power has dimension $\binom{d+k-1}{k}$ and
lies in the full tensor power of dimension $d^k$. The identity follows
from the tensor inner product.
\end{proof}

\begin{lemma}[Spectral sparsification {\cite[Theorem~3.1]{bss12}}]\label{lem:bss_spectral_sparsification}
Let $b_1,\ldots,b_m\in\mathbb R^N$, set $M:=\sum_{i=1}^m b_i b_i^\top$, and write $r:=\operatorname{rank}M$. For every $0<\epsilon<1$, there are a set $S\subseteq\{1,\ldots,m\}$ with $|S|\le Cr/\epsilon^2$ and positive weights $s_i$, $i\in S$, such that
\[
 (1-\epsilon)M\preceq\sum_{i\in S}s_i b_i b_i^\top\preceq(1+\epsilon)M.
\]
Here $C$ is absolute, and $X\preceq Y$ means that $Y-X$ is positive semidefinite.
\end{lemma}

\begin{lemma}[Row sensitivities {\cite[Lemma~B.5]{wy24}}]\label{lem:small_p_sensitivity}
Let $0<p<1$, let $A\in\mathbb R^{m\times n}$ have rank $r\ge1$, and
write $F_A(x):=\sum_i|a_i^\top x|^p$. For
$\sigma_i:=\sup_{x:F_A(x)>0}|a_i^\top x|^p/F_A(x)$, one has
$\sum_i\sigma_i\le r$.
\end{lemma}

\begin{lemma}[Sensitivity sampling {\cite[Theorem~1]{bflsz21}}]\label{lem:small_p_sampling}
Let $f_1,\ldots,f_m:X\to\mathbb R_{\ge0}$, define $F:=\sum_i f_i$, and
suppose $F(x)>0$ for some $x\in X$. Let
$\sigma_i:=\sup_{x:F(x)>0}f_i(x)/F(x)$, discard indices with
$\sigma_i=0$, and define $T:=\sum_i\sigma_i$.
Suppose the range family
$\{\{i:f_i(x)/\sigma_i\ge t\}:x\in X,\ t\in\mathbb R\}$
has VC dimension at most $v$.
For every $0<\delta<1/2$, there exist nonnegative weights $w_i$ with
at most $CT\delta^{-2}(v+1)\log(2T)$ nonzero entries such that
$(1-\delta)F(x)\le\sum_iw_i f_i(x)\le(1+\delta)F(x)$ for every
$x\in X$. The constant $C$ is absolute.
\end{lemma}

\begin{lemma}[Row lower bound at logarithmic rank {\cite[Theorem~6.1]{lww21}}]\label{lem:pre_lww}
Fix a non-even $p>2$. There are constants $c'_p,C'_p,b_p,\epsilon_p>0$
such that, for every $0<\epsilon\le\epsilon_p$ and
$k:=\lceil C'_p\log(2/\epsilon)\rceil$, some rank-$k$ matrix admits no
$(1\pm3\epsilon)$ approximation of its $\ell_p$ energy by fewer than
$c'_pk/(\epsilon^2\log^{b_p}(2/\epsilon))$ reweighted original rows.
\end{lemma}

\begin{lemma}[Sparsification at fixed dimension {\cite[Theorem~1.2]{li26fixed}}]\label{lem:pre_li_sparsification}
Fix an integer $d\ge2$ and $p\ge1$. For every probability measure $\mu$
with finitely many atoms on $S^{d-1}$ and every $\epsilon'\in(0,1/2)$,
there is a probability measure $\mu'$ with positive masses on at most
$C_{d,p}\epsilon'^{-2(d-1)/(d+2p)}$ of the atoms of $\mu$ such that, for
every unit vector $x$,
\[
 |\int|\langle x,u\rangle|^p\,\d\mu'(u)-\int|\langle x,u\rangle|^p\,\d\mu(u)|
 \le\epsilon' .
\]
\end{lemma}

\begin{lemma}[Band discrepancy {\cite[Theorem~2.1]{li26fixed}}]\label{lem:pre_li_band}
Fix an integer $n\ge1$. There are constants $c_n,C_n>0$ with the
following property. Let
$N\ge2$ be an integer, let $\mu^0$ be a probability measure with finitely
many atoms on $S^n$, let $V\subseteq C(S^n)$ be a subspace containing the
constants with $\dim V\le c_nN$, and let $\xi$ be a finite signed Borel
measure on $[0,1]$. Then there is a probability measure $\mu^1$ with
positive masses on at most $N$ atoms of $\mu^0$ such that $\int f\,\d\mu^1=
\int f\,\d\mu^0$ for $f\in V$ and
\[
 \sup_{x\in S^n}|\int_{[0,1]}(\mu^1(B_{x,s})-\mu^0(B_{x,s}))
 \,\d\xi(s)|\le C_nN^{-1/2-1/(2n)}\|\xi\|,\qquad
 B_{x,s}:=\{u\in S^n:|\langle x,u\rangle|\le s\}.
\]
\end{lemma}

We use the model of quantum computation of \cite[Section~2]{ko23}. A
procedure acts on a finite-dimensional complex Hilbert space with a fixed
orthonormal basis, by unitary operations and by measurements in that basis,
and it may choose later operations as functions of earlier outcomes.
Classical random bits come from measuring qubits prepared in a uniform
superposition. Code for a real random variable on a finite probability space
$(\Omega,\lambda)$ consists of a synthesizer, a unitary that maps a fixed
start vector to $\sum_\omega\lambda_\omega^{1/2}|\omega\rangle|{\rm
garbage}_\omega\rangle$, and a unitary that writes the value at $\omega$ into
an output register of finitely many bits. A procedure may apply both, their
inverses and their controlled versions \cite[Definitions~2.2, 2.5
and~2.10]{ko23}.

\begin{lemma}[Quantum mean estimation {\cite[Theorem~1.1 and footnote~2]{ko23}}]\label{lem:pre_ko}
There is an absolute constant $C_{\rm KO}$ such that, for each integer
$n\ge1$, there is a procedure that, given code for a real random variable
$Y$ on a finite probability space, uses the code at most $C_{\rm KO}n$
times in every run and outputs $\widehat\mu$ with
$\Pr[|\widehat\mu-\E[Y]|>\operatorname{sd}(Y)/n]\le1/3$.
\end{lemma}

Four classical facts are proved where they are used: Harris's
inequality~\cite{har60}, Mehler's formula~\cite{meh66}, the bound of
Haussler, Littlestone and Warmuth on the edges of a one-inclusion
graph~\cite{hau95}, and Gaussian concentration by the argument of Maurey
and Pisier~\cite{pis86}.

%% file: 30_lp_p_below_2.tex
\section{\texorpdfstring{Row Sparsification for $0<p\le2$}{Row sparsification for 0 < p <= 2}}
\label{sec:general_lp_field}

\subsection{Main Result}

We prove the following bound by iterated partial coloring. Each round
removes at least a quarter of the current rows and changes every
energy value by a relative error of order $\sqrt{r/(N\min\{p,1\})}$
on $N$ rows. The Gaussian estimate behind a round is
Lemma~\ref{lem:general_lp_cone}, and the volume form of partial
coloring is Lemma~\ref{lem:volume_partial_coloring}. Since the
support shrinks geometrically, the errors of the rounds form a
convergent geometric series. Section~\ref{sec:small_p_uniform} removes
the factor $1/\min\{p,1\}$ for small $p$ at the cost of a factor
$\log(2r)$.

\begin{theorem}[Linear row bound for fixed $p$]\label{thm:general_lp_field}
Let $0<p\le2$ and let $A\in\R^{m\times n}$ have rows $a_i^\top$ and
rank $r$. For every $0<\epsilon<1$, there are a set $S$ of at most
$Cr/(\min\{p,1\}\epsilon^2)$ original rows and positive weights $w_i$
such that for every $x\in\R^n$,
\[
 (1-\epsilon)\|Ax\|_p^p
 \le \sum_{i\in S} w_i |a_i^\top x|^p
 \le (1+\epsilon)\|Ax\|_p^p .
\]
The constant $C$ is absolute. In particular, $Cr/\epsilon^2$ rows suffice for every $1\le p\le2$.
\end{theorem}
\begin{proof}
If $r=0$, the empty set preserves the energy. We discard zero rows
and assume $r\ge1$. Consider a current system
of $N$ nonzero rows of rank at most $r$, with the current weights
absorbed into the rows. Write $f_i(x):=|a_i^\top x|^p$ and
$F:=\sum_i f_i$. We use the convex cone
\[
 \mathcal D:=\{z\in\R^N:\sum_i z_i f_i(x)\le0\text{ for all }x\}.
\]
Let $G$ be standard Gaussian in $\R^N$ and define $s:=80r/\min\{p,1\}$.
For every positive diagonal matrix $D$, we can absorb $D_{ii}$ into
row $i$ by rescaling the row by $D_{ii}^{1/p}$.
Lemma~\ref{lem:general_lp_cone} and Gaussian symmetry then give
$\Pr[DG\in\mathcal D]\ge e^{-s}$. Assume $N\ge4s$ and define
$\eta:=2\sqrt{s/N}$ and $J:=\mathcal D\cap3\sqrt N B_2^N$.

\emph{Every halfspace containing the cube captures volume.}
Fix $a\in\R^N\setminus\{0\}$ and let
$J_a:=J\cap\{z:a^\top z\le\eta\|a\|_1\}$. We show that
$\operatorname{vol}(J_a)\ge e^{-2N}$. Take
$D_{ii}:=\min(1,\|a\|_1/(N|a_i|))$ when $a_i\ne0$ and set $D_{ii}:=1$
otherwise. Since $\|Da\|_\infty\le\|a\|_1/N$ and
$\|Da\|_1\le\|a\|_1$, we have $\|Da\|_2\le\|a\|_1/\sqrt N$.
Using $\log t\le t$ on the entries where $D_{ii}<1$ gives
$-\log\det D\le\sum_i N|a_i|/\|a\|_1=N$.
The Gaussian tail for $a^\top DG$ bounds the probability of leaving
the halfspace by $e^{-2s}$. Since $D_{ii}\le1$, the Gaussian
norm tail gives $\Pr[\|DG\|_2>3\sqrt N]\le2^{N/2}e^{-9N/4}\le e^{-N}$.
Thus $\Pr[DG\in J_a]\ge e^{-s}-e^{-2s}-e^{-N}\ge e^{-s}/2$.
Dividing by the maximum density $(2\pi)^{-N/2}/\det D$ proves the
volume bound.

\emph{Symmetrization and partial coloring.}
For a convex body $K$ define $v_K(b):=\operatorname{vol}(K\cap(2b-K))$.
Integrating in $b$ gives $\int v_K(b)\d b=2^{-N}\operatorname{vol}(K)^2$.
Since $v_K$ vanishes outside $K$, averaging with $K=J_a$ gives a point
$b\in J_a$ with $v_J(b)\ge v_{J_a}(b)\ge2^{-N}e^{-2N}$.
By continuity and Lemma~\ref{lem:pre_brunn_minkowski}, this superlevel set of $v_J$
is compact and convex. It meets every halfspace
$\{z:a^\top z\le\eta\|a\|_1\}$, so by separation it intersects
$\eta[-1,1]^N$. Choose $b$ in the intersection and take
$T:=(J-b)\cap(b-J)$. This body is symmetric and has volume at least
$2^{-N}e^{-2N}$. For $z\in T$, we have $b\pm z\in J$, so
$\|z\|_2\le3\sqrt N$ and
$|\sum_i z_i f_i(x)|\le-\sum_i b_i f_i(x)\le\eta F(x)$ for all $x$.

Define $Q:=[-1,1]^N$. At most $C_0^N$ tiles $2k+Q$, $k\in\mathbb Z^N$,
meet $T$, since they lie in $5\sqrt N B_2^N$ and have disjoint interiors.
Lemma~\ref{lem:pre_brunn_minkowski} applied to $T\cap(2k+Q)$ and its reflection shows that
each tile contributes at most $\operatorname{vol}(T\cap Q)$. Thus
$\operatorname{vol}(T\cap Q)\ge c_0^N$ for an absolute $c_0>0$.
By Lemma~\ref{lem:volume_partial_coloring} with $K:=T\cap Q$, some
$z\in C_1(T\cap Q)\cap Q$ has at least $N/2$ entries of magnitude one.
Choose its sign so at least $N/4$ entries equal $-1$. Multiplying row
weights by $1+z_i$ leaves at most $3N/4$ positive weights and gives
simultaneous relative error at most $C_1\eta$ in this round.

\emph{Iteration.}
We stop when the support is at most $R:=Cr/(\min\{p,1\}\epsilon^2)$.
Before the last update, the support exceeds $R$. Reading the updates
backwards, each earlier support size is at least $4/3$ times the next
one. For $C$ sufficiently large, every update is eligible and the
sum of its relative errors is bounded by
\[
 2C_1\sqrt{s/R}\sum_{j\ge0}(\sqrt{3/4})^j\le\epsilon/2 .
\]
The lower factors multiply to at least $1-\epsilon/2$ and the upper factors
to at most $e^{\epsilon/2}\le1+\epsilon$. The product of the multipliers
on each surviving row gives its positive weight.
\end{proof}

\subsection{A One-Sided Gaussian Estimate}
\begin{lemma}[Gaussian measure of the positive cone]\label{lem:general_lp_cone}
Let $A$ and $p$ be as in Theorem~\ref{thm:general_lp_field} and let
$g_1,\ldots,g_m$ be independent standard Gaussians. Then
\[
 \Pr[\sum_{i=1}^m g_i|a_i^\top x|^p\ge0
             \text{ for every }x\in\mathbb R^n]
 \ge \exp(-80r/\min\{p,1\}).
\]
\end{lemma}
\begin{proof}
The proof has three steps. First, at a global maximum of a confined
and tilted Gaussian field, a fractional second difference produces a
point of a shifted copy of the cone. Second, the area formula converts
the set of admissible tilts into an integral of a Jacobian determinant
over the field's domain. Third, a bound on that integrand at each
point, which involves the measure $q$ of the cone, yields the claimed
lower bound on $q$.

We discard zero rows and restrict to the $r$-dimensional row span.
If $r=0$ the event has probability one. We first assume $0<p<2$ and
define $\rho:=\min\{p,1\}$ and $\lambda:=\sqrt\rho$.
We use Lewis coordinates with exponent $2p$. To obtain them, minimize
$\sum_i(a_i^\top Q a_i)^p-p\log\det Q$ over positive definite $Q$.
Since the rows span, the first term grows at least as a positive multiple
of $\|Q\|_{\rm op}^p$ when $\|Q\|_{\rm op}$ tends to infinity. On bounded
sets the logarithm gives divergence at the singular boundary. A minimizer
therefore exists. Its first-order condition gives, on putting
$b_i:=Q^{1/2}a_i$, $u_i:=b_i/\|b_i\|_2$ and $c_i:=\|b_i\|_2^{2p}$,
\begin{equation}\label{eq:general_lp_frame}
 \sum_i c_i u_i u_i^\top=I_r,\qquad
 \sum_i c_i=r,\qquad |b_i^\top x|^p=\sqrt{c_i}|u_i^\top x|^p.
\end{equation}
The change of coordinates is invertible. Write
$\mathcal C:=\{z:\sum_i\sqrt{c_i}z_i|u_i^\top x|^p\ge0\text{ for every }x\}$
and let $q$ be its standard Gaussian measure.

\paragraph{The field and its maxima.}
Define $w(s):=\log(1+s^2)$, $\phi(s):=\arctan s$ and $s_i:=u_i^\top y$. Define
\[
 W(y):=\sum_i c_i w(s_i),\qquad
 \Psi(y):=\sum_i c_i\phi(s_i)u_i,\qquad
 H(y):=\d\Psi(y)=\sum_i\frac{c_i}{1+s_i^2}u_i u_i^\top .
\]
The matrix $H(y)$ is positive definite. For the zonotope $Z:=\sum_i c_i[-u_i,u_i]$,
Cauchy--Schwarz and Eq.~\eqref{eq:general_lp_frame} give
$\sum_i c_i|u_i^\top x|\le\sqrt r\|x\|_2$, so $Z\subseteq\sqrt r B_2^r$.
Since $\Psi(y)\in(\pi/2)Z$, every $\tau\in\lambda Z/4$ satisfies
$|\tau^\top\Psi(y)|\le\lambda c_0r$ with $c_0:=\pi/8$.
For $y=R\theta$, $\|\theta\|_2=1$, concavity gives
$\log(1+R^2 t)\ge t\log(1+R^2)$ on $0\le t\le1$.
Taking $t:=(u_i^\top\theta)^2$ and summing gives $W(y)\ge\log(1+R^2)$.

We use independent standard Gaussians $g_i,h_i$ to define
\[
 X_i(y):=g_i\cos s_i+h_i\sin s_i,\qquad f(y):=\sum_i\sqrt{c_i}X_i(y).
\]
For each realization, $f$ is bounded. Since $W$ diverges at infinity
and $\tau^\top\Psi$ is bounded, $f-\lambda W+\tau^\top\Psi$ attains
a global maximum for every $\tau\in\lambda Z/4$. At a maximizer $y$,
stationarity gives $T(y)=\tau$, where
$T(y):=H(y)^{-1}(\lambda\nabla W(y)-\nabla f(y))$.
On the event $f(0)\ge0$, comparison with the origin gives
\begin{equation}\label{eq:general_lp_height}
 f(y)\ge\lambda(W(y)-c_0r).
\end{equation}

At the maximizer, the centered difference is nonnegative in every
direction $x$ and at every step $t>0$. We integrate against
$t^{-1-p}\d t/k_p$ with $k_p$ from
Lemma~\ref{lem:general_lp_fractional_difference}. The cosine addition
formula and the substitution $t\mapsto t|u_i^\top x|$ give
$X(y)+b(y,\tau)\in\mathcal C$, where
\[
 b_i(y,\tau):=\sqrt{c_i}\Gamma(p)
 (2\lambda\operatorname{Re} z_i+(u_i^\top\tau)\operatorname{Im} z_i),
 \qquad z_i:=(1-\mathrm i s_i)^{-p}.
\]
Indeed each scalar fractional difference acquires the factor $|u_i^\top x|^p$,
which is also zero when $u_i^\top x=0$. Since $|z_i|\le1$,
\begin{equation}\label{eq:lp_gaussian_shift_size}
 \|b(y,\tau)\|_2^2\le\Gamma(p)^2(4\lambda^2r+\|\tau\|_2^2)\le B_p^2r,
 \qquad B_p^2:=65\lambda^2\Gamma(p)^2/16 .
\end{equation}

\paragraph{The Jacobian and its integral.}
Let $a(s):=(1+s^2)^{-1}$ and define
$D(y,\tau):=\sum_i c_i(\lambda w''(s_i)-(u_i^\top\tau)a'(s_i))u_i u_i^\top$.
Differentiating $H(y)T(y)=\lambda\nabla W(y)-\nabla f(y)$ gives
\begin{equation}\label{eq:lp_field_jacobian}
 \d T(y)=H(y)^{-1}(D(y,T(y))+\sum_i\sqrt{c_i}X_i(y)u_i u_i^\top).
\end{equation}
Write $U$ for the matrix with columns $\sqrt{c_i}u_i$, so $UU^\top=I_r$.
For real $v_i$, the squared Frobenius norm of $U\operatorname{diag}(v_i)U^\top$
is at most that of $U\operatorname{diag}(v_i)$, namely $\sum_i c_i v_i^2$.
Thus $|w''|\le2$ and $|a'|\le1$ imply
$\tr[D(y,\tau)^2]\le8\lambda^2r+2\|\tau\|_2^2\le65\lambda^2r/8$ for $\tau\in\lambda Z/4$.
Lemma~\ref{lem:general_lp_determinant} and the arithmetic--geometric mean
inequality bound the determinant second moment for this $D(y,\tau)$ by $(1+65\lambda^2/8)^r\le(73/8)^r$, since $\lambda\le1$.

We draw a basis $I$ of $r$ indices with probability $c_I\det(U_I)^2$, where
$c_I:=\prod_{i\in I}c_i$ and $U_I$ has columns $u_i$ for $i\in I$.
Cauchy--Binet gives total mass one, and differentiating
$\det(\sum_i t_i c_i u_i u_i^\top)$ at $t_i=1$ gives inclusion probabilities $c_i$.
Another application of Cauchy--Binet and Jensen's inequality gives
\begin{equation}\label{eq:lp_confinement_determinant}
 \det H(y)=\E_I[\prod_{i\in I}a(s_i)]\ge e^{-W(y)}.
\end{equation}
Integrating Eq.~\eqref{eq:lp_confinement_determinant} and changing variables in each basis yields
\begin{equation}\label{eq:lp_confinement_integral}
 \int_{\mathbb R^r}e^{-W(y)}\d y
 \le \pi^r\sum_{I}c_I|\det U_I|
 =(\pi/2)^r\operatorname{vol}(Z),
\end{equation}
where the last equality is the zonotope volume formula.

\paragraph{The probability at a fixed point.}
Let $E$ be the closed set of points $y$ for which $T(y)\in\lambda Z/4$,
$X(y)+b(y,T(y))\in\mathcal C$ and Eq.~\eqref{eq:general_lp_height} holds.
On $f(0)\ge0$ every tilt in $\lambda Z/4$ has a maximizer in $E$.
Applying Lemma~\ref{lem:federer_area_formula} for each realization and then taking expectations gives
\begin{equation}\label{eq:general_lp_area}
 \tfrac12(\lambda/4)^r\operatorname{vol}(Z)
 \le\int_{\mathbb R^r}\E[|\det\d T(y)|{\bf1}_{\{y\in E\}}]\d y .
\end{equation}

At fixed $y$, the variables $Y_i(y):=-g_i\sin s_i+h_i\cos s_i$ form a standard
Gaussian vector independent of $X(y)$. Since $\nabla f(y)=\sum_i\sqrt{c_i}Y_i(y)u_i$,
conditioning on $T(y)=\tau$ leaves $X(y)$ standard Gaussian and $f(y)$ centered
Gaussian of variance $r$. For $\tau\in\lambda Z/4$ write $b:=b(y,\tau)$. The Gaussian
density ratio and H\"older's inequality with exponents $5/4$ and $5$ give, for a
standard Gaussian $G$,
\begin{equation}\label{eq:lp_shifted_cone_probability}
 \Pr[G+b\in\mathcal C]
 =\E[{\bf1}_{\{G\in\mathcal C\}}e^{b^\top G-\|b\|_2^2/2}]
 \le q^{4/5}e^{2\|b\|_2^2}.
\end{equation}
For $T(y)=\tau\in\lambda Z/4$, apply H\"older's inequality with exponents $2,4,4$ to
the determinant of $H(y)\d T(y)$ from Eq.~\eqref{eq:lp_field_jacobian}, the cone indicator controlled by Eq.~\eqref{eq:lp_shifted_cone_probability}, and the height indicator.
Using Eq.~\eqref{eq:lp_confinement_determinant}, Eq.~\eqref{eq:lp_gaussian_shift_size}, and the Gaussian height tail gives a uniform
conditional bound. Averaging over $T(y)$ bounds the integrand in
Eq.~\eqref{eq:general_lp_area} as follows:
\begin{equation}\label{eq:lp_area_integrand_bound}
\E[|\det\d T(y)|\mathbf1_{\{y\in E\}}]
\le e^{W(y)}(73/8)^{r/2}q^{1/5}e^{B_p^2r/2}
 \exp(-\lambda^2(W(y)-c_0r)_+^2/(8r)).
\end{equation}
For $v\ge0$, we have $2v-\lambda^2(v-c_0r)_+^2/(8r)\le(2c_0+8/\lambda^2)r$ by completing the square.
Thus Eq.~\eqref{eq:lp_confinement_integral} lets us integrate Eq.~\eqref{eq:lp_area_integrand_bound} using $e^{-W(y)}$ and cancel $\operatorname{vol}(Z)$.
Since $r\ge1$, we obtain
\begin{equation}\label{eq:lp_scaled_cone_exponent}
 -\log q\le 5r(\log(4\pi/\lambda)+\tfrac12\log(73/8)+\pi/4+8/\lambda^2+B_p^2/2),
\end{equation}
where the bound follows by rearranging the integrated area estimate.
Log-convexity of $\Gamma$ and $\Gamma(1)=\Gamma(2)=1$ give $\Gamma(p)\le1$ on
$[1,2]$ and $\Gamma(p)=\Gamma(1+p)/p\le1/p$ on $(0,1)$. Thus $\Gamma(p)\le1/\rho$ and $\rho B_p^2\le65/16$.
Using $\log t\le t-1$ with $t=1/\rho$ and $\log(4\pi)\ge1/2$ gives
\begin{equation}\label{eq:lp_scale_logarithm}
 \rho\log(4\pi/\sqrt\rho)=\rho\log(4\pi)+\tfrac\rho2\log(1/\rho)\le\rho\log(4\pi)+\tfrac{1-\rho}{2}\le\log(4\pi),
\end{equation}
Multiplying Eq.~\eqref{eq:lp_scaled_cone_exponent} by $\rho$, the logarithm term is bounded by Eq.~\eqref{eq:lp_scale_logarithm}, the shift term by $\rho B_p^2\le65/16$, and the height term satisfies $8\rho/\lambda^2=8$. The remaining positive terms can only decrease because $\rho\le1$.
Since $5(\log(4\pi)+\tfrac12\log(73/8)+\pi/4+8+65/32)<80$, we obtain $-\log q\le80r/\rho$, proving the claimed bound for $0<p<2$.

Finally, take $p_j\uparrow2$ with $1\le p_j<2$, using the original rows and
the same Gaussian coefficients for all $j$. If the inequalities hold for
infinitely many $j$, continuity in $p$ gives the inequalities at $p=2$ for
every $x$. The probability of this upper limit of events is at least the
upper limit of their probabilities, which is at least $e^{-80r}$.
\end{proof}

\subsection{Fractional Differences of the Confinement and the Tilt}
\begin{lemma}[Fractional difference identities]\label{lem:general_lp_fractional_difference}
For $0<p<2$ define $k_p:=2\int_0^\infty(1-\cos t)t^{-1-p}\d t$ and define
\[
 \mathcal D_p v(s):=\frac1{k_p}\int_0^\infty
 (2v(s)-v(s+t)-v(s-t))t^{-1-p}\d t .
\]
For $w(s):=\log(1+s^2)$ and $\phi(s):=\arctan s$,
\[
 \mathcal D_p w(s)=-2\Gamma(p)\operatorname{Re}(1-\mathrm i s)^{-p},
 \qquad
 \mathcal D_p\phi(s)=\Gamma(p)\operatorname{Im}(1-\mathrm i s)^{-p}.
\]
\end{lemma}
\begin{proof}
The bounds $0\le1-\cos t\le\min\{t^2/2,2\}$ give $0<k_p<\infty$.
Differentiation in $s$ and evaluation at zero give the representations
\[
 w(s)=2\int_0^\infty e^{-\xi}(1-\cos(\xi s))\frac{\d \xi}{\xi},
 \qquad \phi(s)=\int_0^\infty e^{-\xi}\sin(\xi s)\frac{\d \xi}{\xi}.
\]
The centered differences of $1-\cos(\xi s)$ and $\sin(\xi s)$ are
$-2\cos(\xi s)(1-\cos(\xi t))$ and $2\sin(\xi s)(1-\cos(\xi t))$, respectively.
Since $\int_0^\infty(1-\cos(\xi t))t^{-1-p}\d t=k_p \xi^p/2$ for $\xi>0$,
absolute integration is bounded by a constant times
$k_p\int_0^\infty e^{-\xi}\xi^{p-1}\d \xi$. Fubini's theorem and
$\int_0^\infty e^{-(1-\mathrm i s)\xi}\xi^{p-1}\d \xi=\Gamma(p)(1-\mathrm i s)^{-p}$
give both identities.
\end{proof}

\subsection{A Gaussian Determinant Moment}
\begin{lemma}[Second moment of a rank-one Gaussian perturbation]\label{lem:general_lp_determinant}
Let $u_1,\ldots,u_m$ be unit vectors in $\mathbb R^r$ and $c_i>0$ satisfy
$\sum_i c_i u_i u_i^\top=I_r$. For every real symmetric $D$ and independent
standard Gaussians $G_1,\ldots,G_m$,
\[
 \E[\det(D+\sum_i\sqrt{c_i}G_i u_i u_i^\top)^2]\le\det(I_r+D^2).
\]
\end{lemma}
\begin{proof}
Use an orthonormal eigenbasis of $D$, with eigenvalues $\lambda_1,\ldots,\lambda_r$,
and let $V$ have columns $u_i$ in this basis. Write $V_I$ for its columns in $I$,
and $V_{J,I}$ for its submatrix
with rows $J$ and columns $I$, and $c_I:=\prod_{i\in I}c_i$.
Expanding the determinant along the diagonal of $D$ and applying Cauchy--Binet
to each complementary minor, the coefficient of $\prod_{i\in I}G_i$ is
\[
 \sqrt{c_I}\sum_{|J|=|I|}\det(V_{J,I})^2\prod_{j\notin J}\lambda_j .
\]
By Cauchy--Binet and Hadamard's inequality,
$\sum_{|J|=|I|}\det(V_{J,I})^2=\det(V_I^\top V_I)\le1$.
Hence the squared coefficient is at most
$c_I\sum_{|J|=|I|}\det(V_{J,I})^2\prod_{j\notin J}\lambda_j^2$.
For fixed $J$, Cauchy--Binet and isotropy give
$\sum_{|I|=|J|}c_I\det(V_{J,I})^2=\det((V\operatorname{diag}(c_i)V^\top)_{J,J})=1$.
Distinct square-free Gaussian monomials are orthogonal and have squared norm one.
Summing the squared coefficients gives at most
$\sum_J\prod_{j\notin J}\lambda_j^2=\prod_j(1+\lambda_j^2)$.
\end{proof}

%% file: 35_small_p.tex
\section{\texorpdfstring{Row Sparsification Uniformly in $0<p<1$}{Row sparsification uniformly in 0 < p < 1}}
\label{sec:small_p_uniform}

\subsection{Main Result}

For $p$ above $c/\log(2r)$, Theorem~\ref{thm:general_lp_field} already
gives the bound below. For smaller $p$ we replace its Gaussian round
by the round of Lemma~\ref{lem:small_p_round}, which removes an eighth
of the rows at relative error $\theta$ once there are
$Cr\theta^{-2}\log(2r)$ rows, and we iterate it as in
Section~\ref{sec:general_lp_field}.

\begin{theorem}[Uniform bound in the quasi-norm range]\label{thm:small_p_uniform}
Fix $0<p<1$. Let $A\in\mathbb R^{m\times n}$ have rank $r\ge1$.
For every $0<\epsilon\le1/2$, there exist a set
$S\subseteq\{1,\ldots,m\}$ with
$|S|\le Cr\epsilon^{-2}\log(2r)$ and positive weights $w_i$, $i\in S$,
such that
\[
 (1-\epsilon)\|Ax\|_p^p
 \le\sum_{i\in S}w_i|a_i^\top x|^p
 \le(1+\epsilon)\|Ax\|_p^p
 \qquad\text{for every }x\in\mathbb R^n.
\]
The constant $C$ is absolute and, in particular, independent of $p$.
\end{theorem}
\begin{proof}
We write $c_*,C_*$ for the constants in Lemma~\ref{lem:small_p_round}.
If $p>c_*/\log(2r)$, Theorem~\ref{thm:general_lp_field} gives at most
$C_3r/(p\epsilon^2)\le(C_3/c_*)r\epsilon^{-2}\log(2r)$ rows.
We use the round from Lemma~\ref{lem:small_p_round} for
$0<p\le c_*/\log(2r)$.

Define $\beta:=7/8$, choose an absolute $a>0$ with
$a/(1-\sqrt\beta)\le1/4$, and set
$R:=C'r\epsilon^{-2}\log(2r)$, where $C'\ge C_*/a^2$.
Whenever the current support $m$ exceeds $R$, we apply the round with
\[
 \theta(m):=a\epsilon\sqrt{R/m}.
\]
Its size hypothesis holds because
$m\theta(m)^2=a^2\epsilon^2R\ge C_*r\log(2r)$.
We absorb each current weight into its row by multiplying by the
$p$th root. Each step reduces the support to at most $\beta m$.
For the counts $m_0,\ldots,m_{q-1}>R$ before deletion, we have
$m_j\ge R\beta^{-(q-1-j)}$. Consequently,
\[
 \sum_{j=0}^{q-1}\theta(m_j)
 \le a\epsilon\sum_{h\ge0}\beta^{h/2}\le\epsilon/4.
\]
The lower multiplicative factors have product at least
$1-\sum_j\theta(m_j)\ge1-\epsilon/4$. The upper factors have product
at most $\exp(\sum_j\theta(m_j))\le e^{\epsilon/4}\le1+\epsilon$,
where the last step uses $e^t\le1+2t$ for $0\le t\le1/2$.
The final support is at most $R$, and all surviving weights are positive.
\end{proof}

\subsection{Counting Patterns with Exceptional Coordinates}

\begin{lemma}[A bound on the total size of shattered blocks]\label{lem:small_p_blocks}
Let $I$ be partitioned into nonempty blocks $Q_1,\ldots,Q_k$, and let
$\ell_j:\R^d\to\R$, $j\in I$, have joint rank at most $r\ge1$.
Suppose every binary pattern on $Q_i$ has a witness whose own high
coordinates satisfy $|\ell_j|\ge1$ and whose own low coordinates
satisfy $|\ell_j|\le1-g$. Suppose each such witness has at most
$D|Q_i|$ coordinates outside $Q_i$ where $|\ell_j|>\eta$.
Then $\sum_i|Q_i|\le C_b rD$ under the rank-dependent condition
\[
 0<g<1,\quad D\ge1,\qquad 0\le\eta\le\frac{g}{40960r}.
\]
One may take $C_b:=1280\cdot1024$.
\end{lemma}
\begin{proof}
Define $n_i:=|Q_i|$ and $N_0:=\sum_i n_i$. Within each block, thresholding
at $1-g/2$ realizes all $2^{n_i}$ patterns. If $n_i\ge5r$,
Lemma~\ref{lem:small_p_regions} would give
$2^{n_i}\le(2en_i/r)^r$, which is impossible. Thus $n_i<5r$.

Suppose $N_0\ge1280\cdot1024rD$. We select each block independently
with probability $\rho:=1280r/N_0$, so that
$\rho D\le1/1024$. For the random selected set $T$, define
$N_T:=\sum_{i\in T}n_i$ and $\mu:=\E[N_T]=1280r$. Then
\[
 \operatorname{Var}(N_T)\le5r\mu,
 \qquad \Pr[N_T<\mu/2]\le20r/\mu=1/64.
\]
Fix one witness for each pattern in each block. We sum the witnesses
for the restrictions of a pattern to the selected blocks. A coordinate
is bad if a witness from another selected block has absolute value
greater than $\eta$ there. Write $\beta(T)$ for the mean number of
bad coordinates over uniform patterns on the selected coordinates.
Each potentially bad incidence requires two blocks to be selected, so
\[
 \E_T[\beta(T)]\le\rho^2DN_0=\rho D\mu\le\mu/1024.
\]
Consequently $\Pr[\beta(T)>\mu/64]\le1/16$. Also
$\E[|T|]=\rho k\le\mu$, whence $\Pr[|T|>8\mu]\le1/8$.
These three failure probabilities sum to less than one. Fix $T$ with
\[
 N:=N_T\ge\mu/2=640r,\qquad
 \beta(T)\le N/32,\qquad |T|\le10240r.
\]
For at least half the $2^N$ input patterns, there are at most $N/16$
bad coordinates. On every other coordinate, the sum of interference
has absolute value at most $|T|\eta\le g/4$. Thus thresholding its absolute linear-form values
at $1-g/2$ recovers the input pattern on all its good coordinates.
One threshold pattern accounts for at most
$\sum_{j\le N/16}\binom Nj\le2^{N/2}$ input patterns. Hence at least
$2^{N/2-1}\ge2^{N/3}$ threshold patterns occur. A small scalar increase
makes weak high comparisons strict without changing any low comparison.
Lemma~\ref{lem:small_p_regions} implies
\[
 2^{N/3}\le(2eN/r)^r,
\]
which is false for $N/r\ge32$, in particular for $N\ge640r$.
\end{proof}

\subsection{Deleting a Constant Fraction of the Rows}

\begin{lemma}[One round for small $p$]\label{lem:small_p_round}
There are absolute constants $c,C>0$ such that the following holds.
Let $A\in\R^{m\times n}$ have rank at most $r\ge1$, let $0<\theta<1/2$,
and suppose
\[
 0<p\le \frac{c}{\log(2r)},\qquad
 m\ge C r\theta^{-2}\log(2r).
\]
There are weights $w_i\in\{0,1,2\}$ with
$|\supp(w)|\le7m/8$ such that, simultaneously for every $x\in\R^n$,
\[
 |\sum_iw_i|a_i^\top x|^p-\sum_i|a_i^\top x|^p|
 \le\theta\sum_i|a_i^\top x|^p.
\]

\end{lemma}
\begin{proof}
Write $F(x):=\sum_i|a_i^\top x|^p$. We take all weights zero if $A=0$
and choose $c$ small enough that $p<1$ throughout the stated range.
For each row define
\[
 \sigma_i:=\sup_{F(x)>0}\frac{|a_i^\top x|^p}{F(x)}.
\]
Lemma~\ref{lem:small_p_sensitivity} gives $\sum_i\sigma_i\le r$.

We retain at weight one every row in
\[
 S_0:=\{i:\sigma_i>8r/m\},\qquad I_0:=[m]\setminus S_0,
 \qquad s:=8r/m.
\]
Then $|S_0|<m/8$. On $I_0$ consider the compact class
\[
 \mathcal C:=\{(|a_i^\top x|^p/F(x))_{i\in I_0}:F(x)>0\}\cup\{0\}.
\]
Every member satisfies $0\le c_i\le s$ and $\sum_{i\in I_0}c_i\le1$.
Every nonzero member has the form $c_i=|a_i^\top y|^p$ with
$F(y)=1$. We use $y=0$ for the zero vector. On the row span, the
level set $F(y)=1$ is compact, so $\mathcal C$ is compact as well.

\paragraph{Parameters.}

We fix the constants in the order they are used. Take $c_0,C_0$ from
Lemma~\ref{lem:small_p_cover}, decreasing $c_0$ to at most $1/4$,
and take $C_v$ from Lemma~\ref{lem:small_p_increments}. Take an absolute $K_0$ large enough for the
remainder estimate below and set
\[
 J:=\lceil\log_2(K_0r)\rceil,\qquad
 B:=\lceil\log_2(128J)\rceil,\qquad b_0:=2^{-B}.
\]
In particular $(256J)^{-1}<b_0\le(128J)^{-1}$. Choose an absolute
$\delta_0:=1/(20000C_v)$. Define
\begin{equation}\label{eq:small_p_parameters}
 \begin{split}
 H := & ~ K_H\log(2r),\qquad
 h:=\frac{\delta_0\theta b_0}{\sqrt{rH}},\qquad
 \Delta:=\frac{c_0h}{\sqrt m},\\
 Q := & ~ c_1\theta^2/H,\qquad
 \mathcal A:=\{s/2,s/4,\ldots,s/2^J\},\qquad
 \mathcal B:=\{1,1/2,\ldots,b_0\},
 \end{split}
\end{equation}
where $c_1:=1/(8192C_0)$. Choose $K_H$ sufficiently large, and
then $K\ge\max\{1,64C_b/c_1\}$. We prove the result
assuming $m\theta^2\ge KrH$, which has the required form.
Writing $a_*:=s/2=4r/m$, this assumption gives
\begin{equation}\label{eq:small_p_ratio}
 \frac{a_*\sqrt m}{h}
 =\frac{4r\sqrt{rH}}{\delta_0\theta b_0\sqrt m}
 \le\frac{4r}{\delta_0b_0}\le\frac{1024rJ}{\delta_0}.
\end{equation}
Since $J=O(\log(2r))$, the choice of $K_H$ ensures, for every $r\ge1$,
\begin{equation}\label{eq:small_p_entropy_budget}
 \log(2+C_0a_*\sqrt m/h)\le H,
 \qquad H\ge\log(200J(B+1)),\qquad H\ge1.
\end{equation}
These requirements only involve absolute constants and the rank.
Indeed, define $J_0:=\lceil\log_2K_0\rceil+1$. Since $J\le J_0r$
and $B+1\le10J$, it suffices to take
\[
 K_H:=2+\frac{\max\{\log(2+1024C_0J_0/\delta_0),
                         \log(2000J_0^2),1\}}{\log2}.
\]

For $a>0$ define $\psi_a(t):=\min\{\max\{t-a,0\},a\}$ and
\[
 u_a(c):=(\psi_a(c_i))_{i\in I_0},\qquad
 b_a(c):=\sum_{i\in I_0}\psi_{a/4}(c_i).
\]
For all $c\in\mathcal C$,
\begin{equation}\label{eq:small_p_mass_bounds}
 \sum_{a\in\mathcal A}b_a(c)\le1,\quad
 \|u_a(c)\|_1\le4b_a(c),\quad
 |\{i\in I_0:c_i\ge a/2\}|\le4b_a(c)/a.
\end{equation}
Also, for $0\le t\le s$,
\begin{equation}\label{eq:small_p_ramps}
 t=\sum_{a\in\mathcal A}\psi_a(t)+\min\{t,s2^{-J}\}.
\end{equation}
The dyadic decomposition telescopes to give Eq.~\eqref{eq:small_p_ramps}.
For Eq.~\eqref{eq:small_p_mass_bounds}, we use three facts. The ramps at
scales $a/4$ sum pointwise to at most $t$, we have
$\psi_a(t)\le4\psi_{a/4}(t)$, and $\psi_{a/4}(t)=a/4$ for $t\ge a/2$.

\paragraph{Protecting coordinates.}

Call a scale active if $a>h/\sqrt m$. Only active scales need nets.
For an active $a$ and $b\in\mathcal B$ define
\[
 \mathcal U_{a,b}:=\{u_a(c):c\in\mathcal C,\ b_a(c)\le b\},
 \qquad \tau_{a,b}:=\max\{1,\lfloor Qb/a\rfloor\}.
\]
We start with $S_1=\varnothing$. Whenever a class restricted to
$I_0\setminus S_1$ has fat dimension at least $\tau_{a,b}$ at scale
$\Delta$, we add a shattered set of exactly that size to $S_1$.
This process terminates because every step adds a coordinate.
Write $\alpha_i,\beta_i$ for the parameters of the $i$th selected block
and $n_i:=\tau_{\alpha_i,\beta_i}$ for its size.

We normalize each block so that the high threshold is one.
For a fat-shattering threshold $\gamma_j$ in block $Q_i$, define
$t_j:=\alpha_i+\gamma_j+\Delta$, so that
$\alpha_i\le t_j\le2\alpha_i$. Define the linear forms
\[
 \ell_j(z):=\frac{a_j^\top z}{(t_j/\alpha_i)^{1/p}}
 \quad(j\in Q_i).
\]
For a witness represented by $c_j=|a_j^\top y|^p$, use
$z:=y/\alpha_i^{1/p}$. Its own high coordinates have $|\ell_j(z)|\ge1$,
and its own low coordinates satisfy
\[
 |\ell_j(z)|\le(1-2\Delta/t_j)^{1/p}
 \le1-\Delta/\alpha_i\le1-g,
 \qquad g:=\Delta/a_*.
\]
If any active scale exists, $0<g<c_0\le1/4$. At a coordinate of
another block $Q_d$, whenever $c_j<\alpha_i/2$,
\[
 |\ell_j(z)|^p=\frac{c_j}{\alpha_i(t_j/\alpha_d)}<1/2.
\]
Thus the background bound is $\eta:=2^{-1/p}$. By Eq.~\eqref{eq:small_p_mass_bounds},
there are at most $4\beta_i/\alpha_i$ exceptions. Since
$n_i\ge Q\beta_i/(2\alpha_i)$, the number of exceptions is at most
$Dn_i$ with $D:=8/Q$.

Eq.~\eqref{eq:small_p_parameters} and $m\theta^2\ge rH$ give
\begin{equation}\label{eq:small_p_gap}
 g=\frac{c_0\delta_0\theta b_0\sqrt m}{4r\sqrt{rH}}
 \ge\frac{c_0\delta_0b_0}{4r}
 >\frac{c_0\delta_0}{1024rJ}.
\end{equation}
Consequently the hypothesis of Lemma~\ref{lem:small_p_blocks} follows if
\begin{equation}\label{eq:small_p_exponent}
 2^{-1/p}\le\frac{c_0\delta_0}{1024\cdot40960\cdot r^2J}.
\end{equation}
There is an absolute $c>0$ for which $p\le c/\log(2r)$ implies
Eq.~\eqref{eq:small_p_exponent} for all $r\ge1$. Explicitly, define
\[
 \kappa:=\frac{c_0\delta_0}{1024\cdot40960},\qquad
 d_0:=3+\log_2(J_0/\kappa),\qquad c:=\frac{\log2}{d_0}.
\]
Then $p\le c/\log(2r)$ gives
$2^{-1/p}\le(2r)^{-d_0}\le\kappa/(8J_0r^3)
\le\kappa/(r^2J)$, and $p\le1$.

Lemma~\ref{lem:small_p_blocks} now gives
\[
 |S_1|\le 8C_b r/Q=8C_b rH/(c_1\theta^2)\le m/8.
\]
If there are no active scales, take $S_1=\varnothing$. In either case
$S:=S_0\cup S_1$ has size less than $m/4$.

\paragraph{Concentration at coarse scales and uniform fine increments.}

We draw independent uniform signs $\xi_i$ for $i\notin S$ and define
$\xi_i:=0$ on $S$. We can equivalently draw signs on all coordinates
and pad the classes by zero on $S$. Let
$W_a:=\{u_a(c)|_{[m]\setminus S}:c\in\mathcal C\}$, padded by zero.
The dimension bound in Lemma~\ref{lem:small_p_increments} applies to every $W_a$. Lemma~\ref{lem:small_p_increments} and Eq.~\eqref{eq:small_p_entropy_budget} give
\[
 \E[\max_{a\in\mathcal A}
   \sup_{\substack{u,v\in W_a\\\|u-v\|_2\le h}}
   |\langle\xi,u-v\rangle|]
 \le \sqrt2C_v h\sqrt{rH}\le\frac{\theta b_0}{10000}.
\]
By Markov's inequality (Lemma~\ref{lem:pre_markov}), with probability at least $99/100$,
\begin{equation}\label{eq:small_p_fine}
 |\langle\xi,u-v\rangle|\le\theta b_0/64
 \quad\text{for all these pairs simultaneously}.
\end{equation}

For inactive scales, every $u\in W_a$ satisfies $\|u\|_2\le a\sqrt m
\le h$. Since zero belongs to $W_a$, Eq.~\eqref{eq:small_p_fine} controls every
such $u$. For active scales with $\tau_{a,b}=1$, the stopping rule
implies every coordinate of every remaining member is less than
$2\Delta$: a value at least $2\Delta$, together with zero, would
fat-shatter its coordinate at threshold $\Delta$. Thus its norm is at
most $2\Delta\sqrt m\le h$, so Eq.~\eqref{eq:small_p_fine} controls this case too.

For an active pair with $\tau_{a,b}\ge2$, let $U_{a,b}$ be the
restricted and padded class. By Lemma~\ref{lem:small_p_cover}, it has a proper
$h$-net $V_{a,b}$ satisfying, with $R_{a,b}:=\theta^2b/a$,
\[
 \log|V_{a,b}|\le C_0\tau_{a,b}H
 \le C_0c_1R_{a,b}\le R_{a,b}/8192.
\]
Every center has squared norm at most $4ab$. Therefore
\[
 \Pr[|\langle\xi,v\rangle|>\theta b/16]
 \le2\exp(-R_{a,b}/2048).
\]
The failure probability over all the centers is at most
$2\exp(-3R_{a,b}/8192)$. As $\tau_{a,b}\ge2$ implies
$R_{a,b}\ge2H/c_1$, this is at most $2e^{-12H}$, and hence at most
$1/(100J(B+1))$ by Eq.~\eqref{eq:small_p_entropy_budget}. A union bound over pairs shows
that with probability at least $99/100$ all centers satisfy their
bounds. On this event and Eq.~\eqref{eq:small_p_fine}, every remaining
$u\in\mathcal U_{a,b}$ satisfies
\[
 |\langle\xi,u\rangle|\le\theta b/16+\theta b_0/64<\theta b/8.
\]
The inactive and $\tau=1$ cases satisfy the same conclusion.

\paragraph{The truncated remainder and total error.}

The class $\{(\min\{c_i,s2^{-J}\})_{i\notin S}:c\in\mathcal C\}$,
padded by zero, also has fat dimension at most $5r$, contains zero,
and is bounded coordinatewise by $s2^{-J}$. Its expected Rademacher
supremum is at most
\[
 C_v s2^{-J}\sqrt{mr}
 =\frac{8C_v r\sqrt r}{2^J\sqrt m}
 \le\frac{8C_v}{K_0\sqrt{KH}}\theta.
\]
Choosing $K_0\ge6400C_v$ makes this at most $\theta/800$.
Thus with probability at least $99/100$ the remainder is bounded by
$\theta/8$, simultaneously for all $c$. All three accuracy events
hold with probability at least $97/100$.

On the intersection of the three events, fix any $c\in\mathcal C$.
If $b_a(c)>b_0$, choose $b\in\mathcal B$ with
$b_a(c)\le b<2b_a(c)$, giving
$|\langle\xi,u_a(c)\rangle|\le\theta b_a(c)/4$. If
$b_a(c)\le b_0$, use $b=b_0$ and obtain $\theta b_0/8$. Hence
\[
 \sum_{a\in\mathcal A}|\langle\xi,u_a(c)\rangle|
 \le\frac\theta4\sum_a b_a(c)+\frac{J\theta b_0}{8}
 \le\theta/4+\theta/1024.
\]
Adding the remainder in Eq.~\eqref{eq:small_p_ramps} gives
\[
 |\langle\xi,c\rangle|\le\theta/4+\theta/1024+\theta/8<\theta.
\]

We use weight one on $S$ and weight $1+\xi_i$ elsewhere. The expected
number of retained rows outside $S$ is $(m-|S|)/2$. By Lemma~\ref{lem:pre_markov}, with
probability at least $1/3$ this number is at most $3(m-|S|)/4$.
This support event intersects the accuracy event, since
$97/100+1/3>1$. On the intersection the support is at most
$3m/4+|S|/4<13m/16<7m/8$. Multiplying the normalized error bound by
$F(x)$ proves the claim for $F(x)>0$. Queries with $F(x)=0$ are
preserved exactly.
\end{proof}

\subsection{Uniform Rademacher Increments}

\begin{lemma}[Uniform Rademacher increments]\label{lem:small_p_increments}
Let $r,J\ge1$, and let $U_1,\ldots,U_J$ be subsets of $[0,a_*]^m$.
Suppose each class contains zero and has fat-shattering dimension at most
$5r$ at every positive scale. Let $\xi_1,\ldots,\xi_m$ be independent uniform signs. For $h>0$,
\begin{align}
 &\E[\max_{j\le J}\sup_{\substack{u,v\in U_j\\
                                  \|u-v\|_2\le h}}
                 |\langle\xi,u-v\rangle|]\notag\\
 &\hspace{10mm}\le C_v h
       \sqrt{r\log(2+a_*\sqrt m/h)+\log(2J)}.\label{eq:small_p_increment}
 \intertext{For a single such class $U\subseteq[0,a]^m$,}
 &\E[\sup_{u\in U}|\langle\xi,u\rangle|]
       \le C_v a\sqrt{mr}.\label{eq:small_p_wholeclass}
\end{align}
Here $C_v\ge1$ is absolute. The dimension hypothesis holds for every class
of the form
\[
 \{(q(|a_i^\top y|^p))_{i\in I}:y\in Y\},
\]
where $p>0$, the forms have joint rank at most $r$, and $q$ is a fixed
nondecreasing ramp or truncation. Subclasses, coordinate restrictions,
and padding by zero preserve this dimension bound.
\end{lemma}
\begin{proof}
Write $N(U,t)$ for the minimum cardinality of a proper Euclidean $t$-net of $U$.
Applying Lemma~\ref{lem:small_p_cover} with the dimension bound gives
\[
 \log N(U_j,t)\le C r\log(Ca_*\sqrt m/t)
 \qquad(0<t<a_*\sqrt m).
\]
For $D_{j,h}:=\{u-v:u,v\in U_j,\ \|u-v\|_2\le h\}$, covering each
factor at radius $t/8$ and replacing nonempty covering balls by points
of $D_{j,h}$ gives
\[
 \log N(\bigcup_{j\le J}(D_{j,h}\cup-D_{j,h}),t)
 \le \log(2J)+Cr\log(2+Ca_*\sqrt m/t).
\]
This union contains zero and lies in $hB_2^m$. Moreover,
\[
 \E[\exp(s\langle\xi,z-z'\rangle)]
 =\prod_i\cosh(s(z_i-z'_i))
 \le\exp(s^2\|z-z'\|_2^2/2).
\]
We take proper nets at radii $2^{-k}h$, starting with the center zero.
For the edges between two successive nets, the expectation of the
largest absolute increment is at most its maximum Euclidean length
times the square root of twice the logarithm of twice the number of
edges. Summing these bounds yields at most
\[
 C\sum_{k\ge0}2^{-k}h
   \sqrt{\log(2J)+Cr\log(2+Ca_*\sqrt m/h)+Cr(k+1)}.
\]
The geometric sums are finite and give Eq.~\eqref{eq:small_p_increment}. The
remaining error tends to zero since $|\langle\xi,z-z'\rangle|
\le\sqrt m\|z-z'\|_2$. If $h>a_*\sqrt m$, the whole-class estimate
with a union over $J$ gives the same upper bound. Starting the same
argument at radius $a\sqrt m$ for a single class proves
Eq.~\eqref{eq:small_p_wholeclass}.

Finally, fix a nondecreasing ramp or truncation $q$ and consider
\[
 \{(q(|a_i^\top y|^p))_{i\in I}:y\in Y\}.
\]
Shattering $k\ge r$ coordinates would yield $2^k$ distinct strict patterns of
$2k$ affine hyperplanes in $\R^r$, while Lemma~\ref{lem:small_p_regions} bounds their number of regions by $(2ek/r)^r$. For $k\ge5r$, $2^{k/r}>2e(k/r)$ gives a contradiction.
The gap in fat shattering permits choosing a strict threshold inside
each gap before inverting $q$ and the positive power. Constant
coordinates cannot be shattered. Restrictions to a fixed coordinate
set, padding by zero, and subclasses preserve this bound.
\end{proof}

%% file: 40_lp_p_above_2.tex
\section{\texorpdfstring{Row Sparsification for $p>2$}{Row sparsification for p > 2}}
\label{sec:lp_above_two}

\subsection{Main Result}

The general bound of Theorem~\ref{thm:lp_above_two} comes from iterated partial coloring in Lewis
coordinates (Lemma~\ref{lem:lp_above_two_round}). In each round we
write the row energy as a polynomial of degree $q$ in the normalized
inner product plus a small remainder
(Lemma~\ref{lem:lp_above_two_remainder}), restrict the coloring to a
subspace that cancels the leading directions of the polynomial
features, and control the oscillation of the resulting Gaussian
process by chaining (Lemma~\ref{lem:lp_above_two_oscillation}). The
degree $q$ is chosen of order $\log(2\sqrt\Lambda/\epsilon)$. For even
$p$ the polynomial is exact, which gives the bound
$C_pr^{(p+1)/2}/\epsilon$.

\begin{theorem}[Row sparsification above two]\label{thm:lp_above_two}
Fix $p>2$. Let $A\in\mathbb R^{m\times n}$ have rows $a_i^\top$ and rank $r\ge2$. Define $D:=r^{p/2}$ and $\Lambda:=\log(2r)$. For every $0<\epsilon\le1/2$, we can choose a set $S$ of original rows and positive weights $w_i$, $i\in S$, such that
\[
 |S|\le \frac{C_p D}{\epsilon^2}(1+\frac{\Lambda}{(\log(2\sqrt \Lambda/\epsilon))^{p-2}})
\]
and, for every $x\in\mathbb R^n$,
\[
 (1-\epsilon)\|Ax\|_p^p
 \le \sum_{i\in S}w_i|a_i^\top x|^p
 \le (1+\epsilon)\|Ax\|_p^p .
\]
In particular, $C_p r^{p/2}/\epsilon^2$ rows suffice whenever $\log(2\sqrt \Lambda/\epsilon)\ge \Lambda^{1/(p-2)}$. If $p$ is an even integer, we can use at most
\[
 \min\{\frac{C_p r^{p/2}}{\epsilon^2},
             \frac{C_p r^{(p+1)/2}}{\epsilon},\binom{r+p-1}{p}\}
\]
rows, where the last bound gives exact equality for every $x$.
In addition, if
\[
 0<\epsilon\le1/16,
 \qquad 2<p\le2+\frac{\epsilon}{4\log r},
\]
then the same approximation guarantee holds with
\[
 |S|\le C\min\{r/\epsilon^2,r^2\},
\]
where $C$ is an absolute constant.
\end{theorem}
\begin{proof}
\emph{Even exponents.} Suppose $p$ is an even integer and define $k:=p/2$. By Lemma~\ref{lem:even_power_lifting}, the vectors $b_i:=a_i^{\otimes k}$ span a space of dimension at most $r^k$ and satisfy $\langle b_i,x^{\otimes k}\rangle^2=|a_i^\top x|^p$. We apply Lemma~\ref{lem:bss_spectral_sparsification} to $\sum_i b_i b_i^\top$ and evaluate the resulting quadratic form at $x^{\otimes k}$. This gives the bound $C r^k/\epsilon^2$. For an exact representation, we restrict to the row span and view each $(a_i^\top x)^p$ as a homogeneous polynomial in $r$ variables. These polynomials lie in a space of dimension $\binom{r+p-1}{p}$. By conic Carath\'eodory, their sum is a nonnegative combination of at most that many original polynomials.

For the intermediate bound, apply Lemma~\ref{lem:lp_above_two_round}
with $q=p$ until the support is at most
$M:=C'_p r^{(p+1)/2}/\epsilon$. Absorb the current row weights before
each round. The even-exponent estimate gives relative error at most
$C_p r^{(p+1)/2}/N$ in a round of size $N$, so each round with $N>M$
preserves rank for sufficiently large $C'_p$. The support shrinks by a
constant factor per round, and the sum of the relative errors is at
most $C''_p r^{(p+1)/2}/M$. Choose $C'_p$ so that this sum is at most
$\epsilon/3$. Multiplying the relative factors gives the claimed
accuracy, as in the argument below.

\emph{The rank-dependent bound.} Discard zero rows and define $t:=\log(2\sqrt \Lambda/\epsilon)$, $q:=\max\{2,\lfloor c_p t\rfloor\}$, and $H:=1+\Lambda/q^{p-2}$. We choose $c_p>0$ so that $C_pc_p\le1/2$, where $C_p$ is the constant in the exponential in Lemma~\ref{lem:lp_above_two_round}. Since $q\le2+c_pt$ and $\epsilon=2\sqrt \Lambda e^{-t}$, we have
\begin{equation}\label{eq:lp_high_accuracy_parameter}
 \epsilon q^{p-1}\exp(C_pq)\le 2e^{2C_p}\sqrt \Lambda(2+c_pt)^{p-1}e^{-t/2}\le C'_p\sqrt \Lambda .
\end{equation}
The choice of $q$ also gives $q\asymp_p t$.

We apply Lemma~\ref{lem:lp_above_two_round} until the support size $N$ is at most $M:=C''_p D H/\epsilon^2$. Before each round we absorb the current weights into the rows by multiplying each row by the $p$th root of its weight. By Eq.~\eqref{eq:lp_high_accuracy_parameter}, we have $\epsilon q\sqrt \Lambda\exp(C_pq)\le C'_p \Lambda/q^{p-2}\le C'_p H$. For $C''_p$ sufficiently large, substituting $N>M$ into the round bound gives $\delta_N<1$. Thus the first round preserves rank $r$, and the same holds inductively for all later rounds. Each round deletes a fixed fraction of the rows. Summing backwards from the last round gives $\sum N^{-1/2}\le C/\sqrt M$ and $\sum N^{-1}\le C/M$. The sum of the relative errors is at most
\[
 C_p\sqrt{D H/M}+(C_pD/M)q\sqrt \Lambda\exp(C_pq)\le\epsilon/3
\]
for $C''_p$ sufficiently large. The product of the upper factors is at most $e^{\epsilon/3}\le1+\epsilon$, and the product of the lower factors is at least $1-\epsilon/3$. We take $w_i$ to be the product of the multipliers on row $i$. Substituting $q$ into $M$ gives the support bound. When $t\ge \Lambda^{1/(p-2)}$, we have $\Lambda/t^{p-2}\le1$.

\emph{Exponents near two.}
Suppose $0<\epsilon\le1/16$ and define $\delta:=p-2$, so
$0<\delta\log r\le\epsilon/4$. Discard zero rows and apply
Lemma~\ref{lem:lp_above_two_lewis}. In the resulting coordinates,
the original energy is $rF_p^\nu(x)$, where
\[
 F_s^\lambda(x):=\sum_i\lambda_i|u_i^\top x|^s,
 \qquad \sum_i\nu_i u_i u_i^\top=I_r/r,
 \qquad \sum_i\nu_i=1.
\]
Here $\nu_i>0$ and $\|u_i\|_2=1$, and the indices still label the
original rows. For every probability vector $\lambda$ on these
indices and every unit vector $x$, Jensen's inequality and
$|u_i^\top x|\le1$ give
\[
 (F_2^\lambda(x))^{p/2}\le F_p^\lambda(x)\le F_2^\lambda(x).
\]
In particular, $r^{-p/2}\le F_p^\nu(x)\le r^{-1}$.

Set $t:=\epsilon/8$. Apply
Lemma~\ref{lem:bss_spectral_sparsification} to the vectors
$\sqrt{\nu_i}u_i$. It gives nonnegative masses $\mu_i$, supported
on at most $C_0r/t^2$ original indices, such that
\[
 (1-t)I_r/r\preceq\sum_i\mu_i u_i u_i^\top\preceq(1+t)I_r/r.
\]
Taking traces shows that $M:=\sum_i\mu_i$ belongs to
$[1-t,1+t]$. Thus $\widehat\mu_i:=\mu_i/M$ is a probability
vector with the same support. Define $h:=\log((1+t)/(1-t))$.
For every unit vector $x$ we have
\[
 e^{-h}/r\le F_2^{\widehat\mu}(x)\le e^h/r.
\]
Applying the moment comparison to $\widehat\mu$ and $\nu$ gives
\[
 e^{-(p/2)h}r^{-\delta/2}
 \le\frac{F_p^{\widehat\mu}(x)}{F_p^\nu(x)}
 \le e^h r^{\delta/2}.
\]
Since $r\ge2$ and $\epsilon\le1/16$, the assumed window implies
$p/2\le2$. Also $h\le2t/(1-t)$, so both logarithmic errors are
at most
\[
 \rho:=\frac p2h+\frac\delta2\log r
 \le\frac{4t}{1-t}+\frac\epsilon8
 \le(\frac{64}{127}+\frac18)\epsilon
 <\frac34\epsilon.
\]
We have $e^{-\rho}\ge1-\rho\ge1-\epsilon$ and
$e^\rho\le1/(1-\rho)\le1+\epsilon$; the last inequality follows
from $\rho\le3\epsilon/4\le\epsilon/(1+\epsilon)$.
Consequently $F_p^{\widehat\mu}$ approximates $F_p^\nu$ within
$1\pm\epsilon$ on unit vectors, and homogeneity extends this to
all $x$. On the retained indices set $w_i:=\widehat\mu_i/\nu_i$.
These weights are positive. Undoing the coordinate change gives
the required approximation on original rows with at most
$64C_0r/\epsilon^2$ rows.

For the second support bound, we represent $I_r/r$ exactly as
$\sum_i\mu_i u_i u_i^\top$ using at most $d:=r(r+1)/2$ positive
coefficients. To see this directly, start with $\mu=\nu$.
If more than $d$ coefficients are positive, the corresponding
symmetric matrices are linearly dependent. Choose a nonzero
real vector $c$ on this support with
$\sum_i c_i u_i u_i^\top=0$. Taking traces gives $\sum_i c_i=0$,
so $c$ has both signs. Subtract $\alpha c$ from $\mu$, where
$\alpha:=\min_{c_i>0}\mu_i/c_i$. This preserves the represented
matrix, keeps all coefficients nonnegative, and deletes at least
one positive coefficient. Repeating gives the claimed support.
Taking traces now gives $\sum_i\mu_i=1$. The same moment
comparison, with $h=0$, yields
\[
 r^{-\delta/2}\le\frac{F_p^\mu(x)}{F_p^\nu(x)}\le r^{\delta/2}
\]
for every unit vector $x$. Since $\delta\log r/2\le\epsilon/8$,
these factors lie in $[1-\epsilon,1+\epsilon]$ by the preceding
exponential estimates. Homogeneity and the weights
$w_i:=\mu_i/\nu_i$ again give the required approximation on
original rows, now using at most $d\le r^2$ rows. Choosing the
smaller of the two constructions proves the additional bound
with an absolute constant $C$.

\end{proof}

\subsection{Lewis Coordinates}

\begin{lemma}[Isotropic probability weights]\label{lem:lp_above_two_lewis}
Fix $p>0$ and let $A$ have nonzero rows and rank $r\ge1$. After restricting queries to its row span and changing coordinates, we can write its energy, up to a common factor $r$, as $\sum_i\nu_i|u_i^\top x|^p$. Here $\nu$ is a probability vector, the $u_i$ are unit vectors in $\mathbb R^r$, and $\sum_i\nu_i u_i u_i^\top=I_r/r$.
\end{lemma}
\begin{proof}
In coordinates on the row span, we minimize $\sum_i(a_i^\top Q a_i)^{p/2}-(p/2)\log\det Q$ over positive definite $Q$. Since the rows span, the first term grows at least as a positive multiple of $\|Q\|_{\mathrm{op}}^{p/2}$ at infinity, whereas $\log\det Q\le r\log\|Q\|_{\mathrm{op}}$. On bounded sets the logarithmic barrier gives divergence at the singular boundary, so a minimizer exists. Its first-order condition is
\begin{equation}\label{eq:lp_lewis_stationarity}
 \sum_i(a_i^\top Q a_i)^{p/2-1}a_i a_i^\top=Q^{-1} .
\end{equation}
We set $c_i:=\|Q^{1/2}a_i\|_2^p$ and $u_i:=Q^{1/2}a_i/\|Q^{1/2}a_i\|_2$.
Conjugate Eq.~\eqref{eq:lp_lewis_stationarity} by $Q^{1/2}$ to obtain
$\sum_i c_i u_i u_i^\top=I_r$, and take the trace to get $\sum_i c_i=r$.
The coordinate change maps the row vectors to $c_i^{1/p}u_i$,
so the assertion follows with $\nu_i:=c_i/r$.
The argument uses only $p>0$: existence of the minimizer needs only
that the first term grows like a positive power of $\|Q\|_{\mathrm{op}}$.
\end{proof}

\subsection{Deleting a Constant Fraction of the Rows}

\begin{lemma}[One round of row deletion]\label{lem:lp_above_two_round}
Fix $p>2$ and an integer $q\ge2$. Let $A$ have $N$ nonzero rows and rank $r\ge2$. Define $D:=r^{p/2}$ and $\Lambda:=\log(2r)$. For all sufficiently large $N$ there are multipliers $t_i\in[0,2]$, at least $cN$ of them zero, such that for every $x$,
\[
 \lvert\sum_i t_i|a_i^\top x|^p-\|Ax\|_p^p\rvert\le\delta_N\|Ax\|_p^p,
\]
where
\[
 \delta_N:=C_p\sqrt{D/N}(1+\sqrt{\Lambda/q^{p-2}})+(D/N)q\sqrt \Lambda\exp(C_pq).
\]
Here $c>0$ is an absolute constant.
If $p$ is an even integer and $q=p$, we can use
$\delta_N:=C_p r^{(p+1)/2}/N$.
\end{lemma}
\begin{proof}
\emph{The coloring subspace.} We use the coordinates from Lemma~\ref{lem:lp_above_two_lewis}. Write $f_x(i):=u_i^\top x$, $B:=\{x:\|f_x\|_{L_p(\nu)}\le1\}$, and $R:=\sqrt r$, so $D=R^p$. Isotropy gives $\|x\|_2/R=\|f_x\|_{L_2(\nu)}\le\|f_x\|_{L_p(\nu)}$, hence $B\subseteq R B_2^r$. We set $J:=\{i:\nu_i\le2/N\}$, so $\ell:=|J|\ge N/2$.

We take $P(t):=\sum_{k=0}^q b_k t^k$ from Lemma~\ref{lem:lp_above_two_remainder} and define $A_q:=\max\{1,\sum_k|b_k|\}\le\exp(C_pq)$. In the orthogonal direct sum of the Euclidean tensor powers through degree $q$, we define
\[
 h(u):=(\sqrt{|b_k|/A_q} u^{\otimes k})_{k=0}^q,
 \qquad b(y):=(\operatorname{sgn}(b_k)\sqrt{|b_k|/A_q} y^{\otimes k})_{k=0}^q .
\]
For $u,y$ in the Euclidean unit ball, both vectors have norm at
most one. They satisfy $P(u^\top y)=A_q\langle h(u),b(y)\rangle$.
The positive operator $\Sigma:=\sum_i\nu_i h(u_i)h(u_i)^\top$ satisfies
$\tr[\Sigma]\le1$. We define $k_0:=\lfloor c_0\ell\rfloor$ for a sufficiently
small absolute constant $c_0>0$ and let $V$ be the span of the $k_0$
leading eigenvectors of $\Sigma$. If its rank is smaller, we take $V$
to be the range. Write $\Pi$ for the orthogonal projection onto
$V^\perp$. Since $\tr[\Sigma]\le1$, we have
$\|\Pi\Sigma\Pi\|_{\mathrm{op}}\le1/(k_0+1)$. We impose orthogonality
to $V$ by restricting the coloring to the subspace
\begin{equation}\label{eq:lp_coloring_subspace}
 E:=\{z\in\mathbb R^J:\sum_{i\in J}z_i\nu_i h(u_i)\perp V\}
\end{equation}
This subspace has codimension at most $c_0\ell$.

\emph{The variance and oscillation.}
Let $G$ be standard Gaussian in $E$ and define
\[
 Z_x:=\sqrt N\sum_{i\in J}G_i\nu_i|f_x(i)|^p.
\]
We show that the standard deviation of $Z_x$ on $B$ is at most
\begin{equation}\label{eq:lp_process_deviation}
 \sigma:=C_p\sqrt D/q^{p/2}+D\exp(C_pq)/\sqrt N .
\end{equation}
We define $\operatorname{rem}_i(x):=|f_x(i)|^p-DP(f_x(i)/R)$. Lemma~\ref{lem:lp_above_two_remainder} gives $|\operatorname{rem}_i(x)|\le C_p q^{-p/2}\sqrt D |f_x(i)|^{p/2}$. Orthogonal projection onto $E$ contracts coefficient norms, and $N\nu_i^2\le2\nu_i$ on $J$. Thus the variance of the remainder sum is at most $2\sum_i\nu_i \operatorname{rem}_i(x)^2\le C_p D/q^p$. By Eq.~\eqref{eq:lp_coloring_subspace}, we can replace $h(u_i)$ by $\Pi h(u_i)$ in the polynomial sum. Its variance is at most
\[
 2D^2 A_q^2\langle b(x/R),\Pi\Sigma\Pi b(x/R)\rangle\le 2D^2 A_q^2/(k_0+1) .
\]
Since $k_0+1\ge c_0\ell\ge c_0N/2$, the triangle inequality for standard deviations gives Eq.~\eqref{eq:lp_process_deviation}.

\emph{Oscillation for even $p$.}
When $p$ is even and $q=p$, we have $P(t)=t^p$, $A_q=1$, and
$\operatorname{rem}_i(x)=0$. Thus the polynomial variance bound gives
$\E[Z_x^2]\le C D^2/N$ for every $x\in B$.
Use the norm with unit ball $B$ and choose an $\eta_p$-net
$\mathcal N\subseteq B$ of size at most $(1+2/\eta_p)^r$.
For a fixed realization, $Z$ is a homogeneous polynomial of degree
$p$. Define $M_Z:=\sup_{x\in B}|Z_x|$. By polarization its symmetric
$p$-linear form has norm at most $(p^p/p!)M_Z$ on $B^p$.
Expanding the difference between the evaluations at $x\in B$ and
its net point gives an error of at most
$p(p^p/p!)M_Z\eta_p$. Choose $\eta_p:=p!/(2p^{p+1})$.
We obtain $M_Z\le2\max_{x\in\mathcal N}|Z_x|$, so the
Gaussian maximum bound gives $\E[M_Z]\le C_p D\sqrt{r/N}$.

\emph{Oscillation for general $p$.}
For the general bound we use chaining for every $p$ and $q$.
For the full process, the scalar mean value bound, projection onto $E$, and $N\nu_i^2\le2\nu_i$ on $J$ imply that the increment standard deviation is at most $C_p\rho(x,y)$, where
\[
 \rho(x,y):=\|(f_x-f_y)(|f_x|^{p-1}+|f_y|^{p-1})\|_{L_2(\nu)} .
\]
We use the cover from Lemma~\ref{lem:lp_above_two_oscillation} with $K:=\Lambda q^2$ and write $U_q:=\sqrt D(1+\sqrt{\Lambda/q^{p-2}})$. For each member $T$ of the cover, we fix $x_0\in T$ and add $\E[|Z_{x_0}|]\le\sigma$ to the oscillation bound to get $\E[\sup_{x\in T}|Z_x|]\le C_p U_q+\sigma$. The supremum is a $\sigma$-Lipschitz function of $G$ because every coefficient vector in $E$ has norm at most $\sigma$. Gaussian concentration and the maximum bound over at most $\exp(C_pK)$ sets give $\E[\sup_{x\in B}|Z_x|]\le C_p(U_q+\sigma\sqrt K)$.

\emph{The partial coloring.} Write $\mu:=\E[\sup_{x\in B}|Z_x|]$.
By Markov's inequality (Lemma~\ref{lem:pre_markov}), the symmetric convex set of $z\in E$ satisfying
\begin{equation}\label{eq:lp_coloring_discrepancy}
 \sup_{x\in B}\lvert\sqrt N\sum_{i\in J}z_i\nu_i|f_x(i)|^p\rvert\le4\mu
\end{equation}
has Gaussian measure at least $3/4$ in $E$. We intersect it with a sufficiently large Euclidean ball so its measure is still at least $1/2$. Since the codimension is at most $c_0\ell$, Lemma~\ref{lem:subspace_partial_coloring} gives a point $z$ in this set and in $[-1,1]^J$ with at least $c_1\ell$ coordinates of absolute value one. This holds for $c_0$ sufficiently small and $\ell$ sufficiently large. By symmetry, we can choose the sign of $z$ so that at least $c_1\ell/2$ coordinates equal $-1$. We set $t_i:=1+z_i$ on $J$ and $t_i:=1$ elsewhere.
Dividing Eq.~\eqref{eq:lp_coloring_discrepancy} by $\sqrt N$ gives
error at most $4\mu/\sqrt N$ on $B$. In the general bound we use
$\sqrt K=q\sqrt \Lambda$ and $\sqrt K/q^{p/2}=\sqrt{\Lambda/q^{p-2}}$.
In the even case the net bound gives error at most
$C_p D\sqrt r/N$. Homogeneity gives the relative error for all
queries, and we undo the Lewis change of coordinates.
\end{proof}

\subsection{Polynomial Approximation of the Row Energy}

\begin{lemma}[Polynomial approximation with a vanishing remainder]\label{lem:lp_above_two_remainder}
Fix $p>2$ and an integer $q\ge2$. There is an even polynomial $P(t):=\sum_{j=0}^q b_jt^j$ with $\sum_{j=0}^q|b_j|\le\exp(C_pq)$ such that $\operatorname{rem}(t):=|t|^p-P(t)$ is $C_p$-Lipschitz on $[-1,1]$ and
\[
 |\operatorname{rem}(t)|\le C_p q^{-p/2}|t|^{p/2}\qquad(-1\le t\le1).
\]
\end{lemma}
\begin{proof}
For even $p$ we take $P(t):=t^p$ when $q\ge p$, and take $P:=0$ and enlarge $C_p$ when $q<p$. For the remaining exponents, we expand $|t|^p$ in Chebyshev polynomials $T_j(t):=\cos(j\arccos t)$. The odd coefficients vanish. For positive even $j$, the coefficient $a_j$ is a fixed multiple of $I_j:=\int_0^{\pi/2}\cos^p(\theta)\cos(j\theta)\d\theta$. We integrate the derivative of $\cos^{p+1}(\theta)\sin((j+1)\theta)$ to obtain
\begin{equation}\label{eq:lp_chebyshev_recurrence}
 (j+p+2)I_{j+2}=(p-j)I_j .
\end{equation}
For even $j>p$, Eq.~\eqref{eq:lp_chebyshev_recurrence} gives the absolute ratio $1-(2p+2)/(j+p+2)$. We take logarithms and sum over even indices to get $|a_j|\le C_p(j+1)^{-p-1}$. The series converges absolutely, and uniqueness of Fourier coefficients identifies its sum with $|t|^p$.

\emph{Matching derivatives at zero.} We define $h:=\lceil p/4\rceil$ and $m:=2h$, so $p/2\le m<p$, and let $n:=\lfloor q/(2h)\rfloor$. If $n=0$, then $q<2h$ and $P:=0$ suffices after increasing $C_p$. Otherwise $n$ is comparable to $q$ up to constants depending on $p$. We truncate the Chebyshev series through degree $2n$ to obtain $P_0$ and define $U(\theta):=|\sin\theta|^p-P_0(\sin\theta)$. Since $T_{2j}(\sin\theta)=(-1)^j\cos(2j\theta)$ for $|\theta|\le\pi/2$, we can differentiate the tail term by term and use the coefficient bound to get
\begin{equation}\label{eq:lp_truncation_derivatives}
 \sup_{|\theta|\le\pi/2}|U^{(\ell)}(\theta)|\le C_p n^{\ell-p}\qquad(0\le\ell\le m).
\end{equation}
We choose $V(\theta):=\sum_{b=0}^{h-1}v_b\cos(2bn\theta)$ so that $V^{(2j)}(0)=U^{(2j)}(0)$ for $0\le j<h$. Dividing the $j$th equation by $(-1)^j(2n)^{2j}$ gives the Vandermonde matrix on $0^2,1^2,\ldots,(h-1)^2$. Its inverse depends only on $p$ and every normalized right-hand side is at most $C_p n^{-p}$ in absolute value. Hence $|v_b|\le C_p n^{-p}$ for each $b$.

We set $W:=U-V$. Its even derivatives of orders less than $m$ vanish at zero by construction, and its odd derivatives there vanish by symmetry. By Eq.~\eqref{eq:lp_truncation_derivatives} and the bounds on $v_b$, we have $\|W\|_\infty\le C_p n^{-p}$ and $\|W^{(m)}\|_\infty\le C_p n^{m-p}$ on $[-\pi/2,\pi/2]$. Taylor's formula therefore gives
\begin{equation}\label{eq:lp_polynomial_remainder}
 |W(\theta)|\le C_p n^{-p}\min\{1,(n|\theta|)^m\}
 \le C_p n^{-p/2}|\theta|^{p/2} .
\end{equation}
We define $P(t):=P_0(t)+\sum_{b=0}^{h-1}(-1)^{bn}v_bT_{2bn}(t)$. It is even and has degree at most $2hn\le q$. For $t=\sin\theta$ with $|\theta|\le\pi/2$, we have $\operatorname{rem}(t)=W(\theta)$ and $|\theta|\le(\pi/2)|t|$. Substituting into Eq.~\eqref{eq:lp_polynomial_remainder} gives the error bound.

\emph{Bounds on the polynomial.} Since $\|T_j'\|_\infty\le j^2$ on $[-1,1]$ and $p>2$, we have $\|P_0'\|_\infty\le\sum_{j\le2n}|a_j|j^2\le C_p$. The correction contributes at most $\sum_{b<h}|v_b|(2bn)^2\le C_p n^{2-p}$ to $\|P'\|_\infty$. The function $|t|^p$ is $p$-Lipschitz on $[-1,1]$, so $\operatorname{rem}$ is $C_p$-Lipschitz there. The recurrence $T_{j+1}(t)=2tT_j(t)-T_{j-1}(t)$ bounds the sum of the absolute monomial coefficients of $T_j$ by $3^j$. Since $\sum_j|a_j|+\sum_b|v_b|\le C_p$, the corresponding sum for $P$ is at most $\exp(C_pq)$.
\end{proof}

\subsection{Gaussian Oscillation on the Unit Ball}

\begin{lemma}[Oscillation on an entropy cover]\label{lem:lp_above_two_oscillation}
Fix $p>2$, $r\ge2$, and an integer $q\ge2$. Let $u_i\in\mathbb R^r$ be unit vectors and let $\nu$ be a probability vector satisfying $\sum_i\nu_i u_i u_i^\top=I_r/r$. Define $f_x(i):=u_i^\top x$, $B:=\{x:\|f_x\|_{L_p(\nu)}\le1\}$, $R:=\sqrt r$, $\Lambda:=\log(2r)$, and $K:=\Lambda q^2$. Define
\[
 \rho(x,y):=\|(f_x-f_y)(|f_x|^{p-1}+|f_y|^{p-1})\|_{L_2(\nu)} .
\]
There is a cover of $B$ by at most $\exp(C_p K)$ sets $T$ such that, for every continuous centered Gaussian process $(Z_x)_{x\in B}$ with $(\E[|Z_x-Z_y|^2])^{1/2}\le\rho(x,y)$ and every $x_0\in T$,
\[
 \E[\sup_{x\in T}|Z_x-Z_{x_0}|]\le C_p R^{p/2}(1+\sqrt{\Lambda/q^{p-2}}) .
\]
\end{lemma}
\begin{proof}
\emph{Separating the two chaining terms.} We identify $x$ with $f_x$ and define $U:=R B_2^r$. Isotropy gives $B\subseteq U$. For each nonempty compact $T\subseteq B$ and $\eta>0$, Lemma~\ref{lem:kosov_gaussian_chaining} gives
\begin{equation}\label{eq:lp_oscillation_chaining}
\begin{aligned}
 \E[\sup_{x\in T}|Z_x-Z_{x_0}|]\le & ~ C_p\sup_{h\in T}(\sum_{k\ge0}2^k e_k(U,\|\cdot\|_{h,\eta})^2)^{1/2}
     +C_p\gamma_{2p,p}(T,L_{2p}(\nu))^p .
\end{aligned}
\end{equation}
The squared norm $\|x\|_{h,\eta}^2$ is the quadratic form of
$Q_h+\eta I_r/r$, where $Q_h:=\sum_i\nu_i |f_h(i)|^{2p-2}u_i u_i^\top$.
Lemma~\ref{lem:kosov_euclidean_entropy} bounds the first term in
Eq.~\eqref{eq:lp_oscillation_chaining} by
$C_p R\sup_{h\in T}\sqrt{\tr[Q_h]+\eta}$.
Taking $\eta\downarrow0$ in this upper bound gives
\begin{equation}\label{eq:lp_ellipsoid_entropy}
 C_p R\sup_{h\in B}(\sum_i\nu_i\|u_i\|_2^2|f_h(i)|^{2p-2})^{1/2}\le C_p R^{p/2} .
\end{equation}
Here we used $|f_h(i)|^{2p-2}\le R^{p-2}|f_h(i)|^p$ and $\sum_i\nu_i|f_h(i)|^p\le1$. To bound the second term in Eq.~\eqref{eq:lp_oscillation_chaining}, we use the common refinements of the partitions induced by successive $L_{2p}(\nu)$ covers in Eq.~\eqref{eq:chaining_functional_definition}. Their product cardinality through level $k$ is at most $2^{2^{k+1}}$. Shifting the levels by one gives
\begin{equation}\label{eq:lp_chaining_entropy_sum}
 \gamma_{2p,p}(T,L_{2p}(\nu))^p\le C_p\sum_{k\ge0}2^{k/2}e_k(T,L_{2p}(\nu))^p .
\end{equation}

\emph{Covering the whole unit ball.} We set $s:=4p\log(2r)$, $a:=p/2$, $\theta:=a/(1-p/s)$, and $\|x\|_*:=\|f_x\|_{L_s(\nu)}$. For $z\in B-B$, interpolation gives
\begin{equation}\label{eq:lp_norm_interpolation}
 \|f_z\|_{L_{2p}(\nu)}^p\le \|f_z\|_{L_p(\nu)}^{p-\theta}\|f_z\|_{L_s(\nu)}^\theta\le 2^{p-\theta}\|z\|_*^\theta .
\end{equation}
Also $\|f_z\|_{L_{2p}(\nu)}^p\le C_pR^a$ by $\|f_z\|_\infty\le2R$ and $\|f_z\|_{L_p(\nu)}\le2$. The choice of $s$ gives $R^\theta\le C_pR^a$.

For $g\sim N(0,I_r)$, Jensen's inequality gives $\E[\|g\|_*]\le (\E[|g_1|^s])^{1/s}\le C\sqrt s$. The centered $\|\cdot\|_*$-ball of radius $v:=2C\sqrt s$ therefore has Gaussian measure at least $1/2$. Translating this symmetric ball by $w$ leaves measure at least $e^{-\|w\|_2^2/2}/2$, by averaging the Gaussian density ratios at $x$ and $-x$. We take a maximal $t$-separated subset of $B$ in $\|\cdot\|_*$. Dilating its points by $2v/t$ gives centers for disjoint radius-$v$ balls up to their boundaries. Their Euclidean norms are at most $2vR/t$. Thus the cardinality is at most $2\exp(CR^2s/t^2)$, and maximality gives $e_k(B,\|\cdot\|_*)\le CR\sqrt{s/2^k}$ for $k\ge1$. Applying Eq.~\eqref{eq:lp_norm_interpolation} to each point and its center gives, with $E_k:=e_k(B,L_{2p}(\nu))$,
\begin{equation}\label{eq:lp_global_entropy}
 E_k^p\le C_pR^a\min\{1,(s/2^k)^{a/2}\} .
\end{equation}
For $2^k\ge s$ this uses $\theta\ge a$.

\emph{Bounding oscillation on each set.} We choose $k_*$ so that $K\le2^{k_*}<2K$ and take a cover at this level by at most $\exp(CK)$ closed $L_{2p}(\nu)$-balls of radius $2E_{k_*}$. Each nonempty intersection $T$ of $B$ with one of these balls has diameter at most $4E_{k_*}$. We also have $e_k(T,L_{2p}(\nu))\le2E_k$, by recentering each ball of a cover of $B$ that intersects $T$ at a point in $T$. We use the diameter bound for $k<k_*$ and Eq.~\eqref{eq:lp_global_entropy} for $k\ge k_*$. Since $a>1$, both geometric sums give
\begin{equation}\label{eq:lp_local_entropy_sum}
 \sum_{k\ge0}2^{k/2}e_k(T,L_{2p}(\nu))^p
 \le C_pR^a s^{a/2}K^{(1-a)/2}=C'_pR^a\sqrt{\Lambda/q^{p-2}} .
\end{equation}
The last equality uses $s=4p\Lambda$, $K=\Lambda q^2$, and $a=p/2$. Substituting Eq.~\eqref{eq:lp_local_entropy_sum} into Eq.~\eqref{eq:lp_chaining_entropy_sum}, and then using Eq.~\eqref{eq:lp_ellipsoid_entropy} in Eq.~\eqref{eq:lp_oscillation_chaining}, proves the claim.
\end{proof}

\subsection{Sparsification at Fixed Rank}

Lemma~\ref{lem:lp_fixed_rank_li} gives the sharp rate at fixed rank.
Lemma~\ref{lem:lp_fixed_rank} uses local polynomial moments to control
the error near the equator of each query, which gives an elementary
second proof of the rank-two rate and weaker bounds in higher ranks.

\begin{lemma}[Fixed rank]\label{lem:lp_fixed_rank_li}
Let $p>2$, let $A$ have rank $r\ge2$, and let $0<\epsilon\le1/2$. There
is a sparsifier satisfying Eq.~\eqref{eq:row_sparsification} with at most
$C_{r,p}\epsilon^{-2(r-1)/(r+2p)}$ original rows.
\end{lemma}
\begin{proof}
Discard zero rows and use Lemma~\ref{lem:lp_above_two_lewis} to write
the energy in Lewis coordinates as $rF_\nu(x)$, where
$F_\nu(x):=\sum_i\nu_i|u_i^\top x|^p$. For every unit vector $x$,
\[
 F_\nu(x)\ge(\sum_i\nu_i|u_i^\top x|^2)^{p/2}=r^{-p/2},
\]
by Jensen's inequality and isotropy. Apply Lemma~\ref{lem:pre_li_sparsification}
to the probability measure $\sum_i\nu_i\delta_{u_i}$ with additive
error $\epsilon r^{-p/2}$. The resulting measure is supported on
original directions. We divide each retained mass by $\nu_i$ to
obtain its row weight. If several input rows have the same direction,
we assign its entire output mass to one of them. The additive error
bound gives the required relative error for unit $x$, and homogeneity
extends it to all $x$. Absorbing the factor $r^{-p/2}$ in the constant
gives the bound.
\end{proof}

\begin{lemma}[Local moment matching]\label{lem:lp_fixed_rank}
Let $A$ have rank $r\ge2$ and fix $p>2$ and $0<\epsilon\le1/2$.
The original rows admit a $(1\pm\epsilon)$ approximation of their
$\ell_p$ energy with at most
\begin{equation}\label{eq:lp_fixed_rank_local_bound}
 C_{r,p}(\epsilon^{-2}\log(2/\epsilon))^{(r-1)/(r-1+2p)}
\end{equation}
rows. When $r<2p+2$, the denominator in the exponent improves to
$r+2p$. When $r=2$, $C_p\epsilon^{-1/(p+1)}$ rows suffice.
\end{lemma}
\begin{proof}
We discard zero rows and use Lemma~\ref{lem:lp_above_two_lewis}. Write $\nu$ for the probability measure on the unit row directions $u_i\in S^{r-1}$ and define $F_\nu(x):=\int |u^\top x|^p\d\nu(u)$. We keep the original row indices as labels, including when two directions coincide. On the Euclidean unit sphere, isotropy gives $F_\nu(x)\ge r^{-p/2}$. It suffices to find a probability measure $\widehat\nu$ on a subset of the original directions with $|F_{\widehat\nu}(x)-F_\nu(x)|\le\epsilon r^{-p/2}$ there. We then extend by homogeneity and assign each retained row the multiplier $\widehat\nu_i/\nu_i$.

\emph{Matching moments on each piece.} We write $\nu$ as a sum of small measures $\nu_C$ supported in spherical cells. Write $w_C$ for the mass of $\nu_C$ and fix $m:=\lceil p\rceil+1$. Consider the polytope of nonnegative weights on the support of $\nu_C$ that agree with $\nu_C$ on all polynomial moments through degree $m$. Every vertex has at most $\binom{r+m}{m}$ nonzero entries. Indeed, if the columns of the moment matrix on its support are linearly dependent, we can perturb its positive weights in both directions along a null vector, contradicting extremality. The constant moment fixes the total mass at $w_C$, so the polytope is bounded. We express $\nu_C$ as a convex combination of its vertices and draw a vertex $\widehat\nu_C$ from this combination, independently for each piece. Thus $\E[\widehat\nu_C]=\nu_C$. We define $\widehat\nu:=\sum_C\widehat\nu_C$ and merge contributions to the same original row.

For any unit $x$ and any polynomial $P$ of degree at most $m$, moment matching gives the following bound on the centered error from piece $C$:
\begin{equation}\label{eq:lp_fixed_rank_piece_error}
 |\int |u^\top x|^p\d(\widehat\nu_C-\nu_C)(u)|
 \le 2w_C\sup_{u\in C}||u^\top x|^p-P(u)| .
\end{equation}
We will choose the cells and the approximating polynomials to bound the sum of the squares of these errors.

\emph{Cells of a common size.} Let $s:=r-1$ and let $h$ be a negative power of two. We partition the sphere into $O_r(h^{-s})$ cells of diameter at most $C_r h$. We use the radial projections of dyadic cubes of side $h$ in the faces of $[-1,1]^r$. At a boundary ray we select the first incident face in a fixed ordering, and within each face we use nested half-open cubes, including the last endpoint. These conventions make each finer cell a subset of its parent. In each base cell we divide its measure into pieces of mass at most $h^s$, allowing an atom to contribute to more than one piece. There are at most $h^{-s}+O_r(h^{-s})$ pieces, since a cell of mass $w$ requires at most $1+w h^{-s}$ pieces.

The derivative of order $\lceil p\rceil-1$ of $g(t):=|t|^p$ is H\"older on $[-1,1]$ of exponent $p-\lceil p\rceil+1$. The Taylor polynomial of this order about $u_C^\top x$, composed with $u^\top x$, therefore approximates $|u^\top x|^p$ on a cell $C$ to within $C_{r,p} h^p$. Here $u_C$ is any fixed point in the cell. By Eq.~\eqref{eq:lp_fixed_rank_piece_error}, the sum of the squared error bounds is at most $C_{r,p} h^{2p}\sum_C w_C^2\le C_{r,p} h^{2p+s}$.

\emph{Refining cells of large mass.} Suppose now that $r<2p+2$. We start at the same base scale $h$ and subdivide each cell of side $h_j:=h2^{-j}$ whose mass exceeds $h^s2^{bj}$, where $0<b<p-(s-1)/2$. We use the measure on each nonempty leaf as a piece. This procedure terminates because the threshold eventually exceeds one. At level $j$, at most $h^{-s}2^{-bj}$ cells exceed the threshold. Each has $2^s$ children, so the total number of pieces is $O_{r,p}(h^{-s})$.

We estimate the error at each scale using the location of the equator $\{u:u^\top x=0\}$. For a level-$j$ cell $C$ choose a point $u_C$ in that cell. If $|u_C^\top x|\le 2C_r h_j$, use the zero polynomial in Eq.~\eqref{eq:lp_fixed_rank_piece_error}. Otherwise use the Taylor polynomial of degree $m$ of $g$ about $u_C^\top x$. On every segment between $u_C^\top x$ and $u^\top x$ for $u\in C$, the absolute value of the argument is at least $|u_C^\top x|/2$ in the latter case. The derivative of order $m+1$ is bounded there in absolute value by $C_p|u_C^\top x|^{p-m-1}$. In both cases, the approximation error on $C$ is at most
\[
 C_{r,p} h_j^p (1+|u_C^\top x|/h_j)^{p-m-1} .
\]
Each level-$j$ spherical cell has surface measure at least $c_r h_j^s$. This follows because the radial projection of each cube face has area Jacobian bounded below by a positive constant depending only on $r$. The band $|u^\top x|\le a$ has surface measure at most $C_r a$ for $a>0$. Since the diameter of each level-$j$ cell is at most $C_r h_j$, there are at most $C_r 2^k h_j^{1-s}$ level-$j$ cells with $|u_C^\top x|\le 2^k h_j$, for any integer $k\ge0$.

We sum the squared approximation bounds over the band $|u_C^\top x|\le h_j$ and over successive bands $2^{k-1}h_j<|u_C^\top x|\le2^k h_j$. The band counts and $2(p-m-1)\le-2$ bound the resulting sum by $C_{r,p}h_j^{2p+1-s}$. Multiplying by the square of the mass threshold and using Eq.~\eqref{eq:lp_fixed_rank_piece_error}, the sum of squared errors over the level-$j$ leaves is at most
\[
 C_{r,p} h^{2p+s+1}2^{-j(2p+1-s-2b)} .
\]
Our choice of $b$ makes the sum over levels at most $C_{r,p}h^{2p+s+1}$.

When $r=2$, we sum the absolute errors directly. Now $s=1$, so each band contains at most $C2^k$ cells, and $p-m-1\le-2$. Summing the approximation bounds over the bands therefore gives at most $C_p h_j^p$ at level $j$. After multiplying by the mass bound $h2^{bj}$, the total error over all levels is at most $C_p h^{p+1}$ because $b<p$. We may take any moment-matching vertex on each piece. Taking $h$ dyadic and of order $c_p\epsilon^{1/(p+1)}$ gives the rank-two bound.

\emph{Uniform approximation.} In the common-size construction set $\alpha:=p+s/2$. In the refined construction set $\alpha:=p+(s+1)/2$. In either case, the independent piece errors have mean zero, and the sum of their squared bounds is at most $C_{r,p}h^{2\alpha}$. By Hoeffding's inequality (Lemma~\ref{lem:pre_hoeffding}), the probability of error larger than $\epsilon r^{-p/2}/2$ at any fixed unit $x$ is at most $2\exp(-c_{r,p}\epsilon^2/h^{2\alpha})$. Both $F_{\widehat\nu}$ and $F_\nu$ are $p$-Lipschitz in Euclidean distance on the unit sphere because their coefficients have total mass one and $g$ is $p$-Lipschitz on $[-1,1]$. We thus take a Euclidean net of radius $\epsilon r^{-p/2}/(4p)$ with at most $(C_{r,p}/\epsilon)^r$ points. The union bound gives a choice of $\widehat\nu$ with the desired accuracy as soon as
\[
 h^{2\alpha}\le c'_{r,p}\epsilon^2/\log(2/\epsilon) .
\]
We take the largest dyadic $h\le1/2$ satisfying this condition. The support bound is $C_{r,p}h^{-s}$, which gives Eq.~\eqref{eq:lp_fixed_rank_local_bound} and its refinement.
\end{proof}

%% file: 60_lower.tex
\section{\texorpdfstring{Row Lower Bounds for $p>2$}{Row lower bounds for p > 2}}
\label{sec:lp_even_lower}

\subsection{Main Result}

Theorem~\ref{thm:lp_lower_main} gives the row lower bounds for $p>2$.
For even $p$ it matches the upper bounds of
Theorem~\ref{thm:lp_above_two}, and its proof extracts a mixed
coefficient of the query error. For non-even $p$ the proof amplifies
a low-rank hard instance by tensor products
(Lemma~\ref{lem:lp_lower_tensor_product}).
Lemma~\ref{lem:lower_kernel_criterion} turns a lower bound on the row
kernel $|a_i^\top a_j|^p$ into a bound on the support of any
sparsifier. Lemma~\ref{lem:lp_real_lower} applies it to real mutually
unbiased bases for every real $p>2$. For non-even $p$,
Lemma~\ref{lem:lp_noneven_lower} isolates single rows using binary
quadratic forms. For even $p$, Lemma~\ref{lem:lp_even_lower} uses
random sign rows and polarization.

\begin{theorem}[Row lower bounds above two]\label{thm:lp_lower_main}
Fix $p>2$ and let $0<\epsilon\le1/2$. There are constants
$c_p,C_p,a_p>0$ depending only on $p$ such that the following hold.
\begin{enumerate}[label=(\alph*),nosep]
\item If $p$ is an even integer, then for every $r\ge2$ some rank-$r$
matrix requires at least
\[
 c_p\min\{r^{p/2}\epsilon^{-2},
          r^{(p+1)/2}\epsilon^{-1},r^p\}
\]
reweighted original rows.
\item If $p$ is not an even integer and $r\ge C_p\log(2/\epsilon)$,
then some rank-$r$ matrix requires at least
\[
 \frac{c_p r^{p/2}}{\epsilon^2\log^{a_p}(2/\epsilon)}
\]
reweighted original rows.
\end{enumerate}
\end{theorem}
\begin{proof}
\emph{Part (a).}
We use several query blocks with different dimensions. Extracting a
mixed coefficient gives a system on the error tensors in these blocks.
The system remains well conditioned when we restrict each query to
the subspaces assigned to its group.
Write $b:=p/2$ and $d:=\lfloor r/(b+1)\rfloor$. We first assume
$d$ sufficiently large and $\epsilon$ sufficiently small in terms
of $p$. The identity matrix handles bounded $r$, since the displayed
minimum is at most $r^p$. For $\epsilon$ bounded away from zero in
terms of $p$, Lemma~\ref{lem:lp_even_lower} gives
$\Omega_p(r^{p/2})$, which suffices after changing $c_p$.

\emph{The families of vectors.}
Let $g:=\lfloor d^b/(4b!)\rfloor$. Choose $g$ unit vectors $u_i$ in
$\R^d$ such that the matrix $K_{\ell i}:=(u_\ell^\top u_i)^b$
satisfies $K\succeq I_g/8$ and
$\sup_{\|x\|_2=1}\sum_i|u_i^\top x|^p\le C_p$. To obtain both
properties, take independent normalized sign vectors. The
distinct-index contribution to the Gram matrix of their $b$th
tensor powers is bounded below by Lemma~\ref{lem:random_sign_tensor_gap},
and the energy bound follows from Lemma~\ref{lem:sign_energy_upper}.
Both events have probability tending to one.

For integers $1\le k_h\le d$ with $h\in[b]$, define
$N:=\prod_{h=1}^b k_h$. For each $h$ take $g$ orthonormal frames
$(v_{i,h,a})_{a\in[k_h]}$ in $\R^d$, indexed by $i\in[g]$.
We choose these so that
\begin{equation}\label{eq:lp_mixed_frame}
 \sup_{\|y\|_2=1}\sum_{i=1}^g\sum_{a=1}^{k_h}|v_{i,h,a}^\top y|^p
 \le C_p k_h\qquad(h\in[b]).
\end{equation}
By Lemma~\ref{lem:rotated_frame_energy}, independent Haar frames
satisfy these conditions with probability tending to one, since
$g=O_p(d^{p/2})$ and $b$ is fixed. Write $U_{i,h}$ for the span of
the $i$th frame in block $h$.

\emph{The mixed coefficient of the query error.}
Use the coordinates
$\R^r=(\R^d)^{b+1}\oplus\R^{r-(b+1)d}$. For each $i\in[g]$ and
$\mathbf a\in[k_1]\times\cdots\times[k_b]$, include the row whose
vector is
\[
 (u_i,v_{i,1,a_1},\ldots,v_{i,b,a_b},0).
\]
Include the rows of $I_r$ to ensure rank $r$. Let $E(x)$ be the
difference between the reweighted and original energies, and
write $w_{i,\mathbf a}$ for the weights on the constructed rows,
with zero weight on omitted rows.
For a query
\[
 x=(u_\ell,t_1 y_1,\ldots,t_b y_b,0),\qquad
 |t_h|\le1,\qquad y_h\in U_{\ell,h},\quad \|y_h\|_2=1,
\]
the original energy is at most
\[
 (b+1)^{p-1}(
 N\sum_i|u_i^\top u_\ell|^p+
 \sum_{h=1}^b\frac{N}{k_h}\sum_{i,a}|v_{i,h,a}^\top y_h|^p
 )+b+1\le C'_p N.
\]
Here the first bound separates the coordinate blocks, and the
second uses Eq.~\eqref{eq:lp_mixed_frame} and the bound on the $u_i$.
A valid sparsifier therefore satisfies $|E(x)|\le C'_p\epsilon N$
throughout the parameter cube. Since $E(x)$ is a polynomial in the
$t_h$ of degree at most $p$ in each variable, interpolation on a
fixed $(p+1)$-point grid in each coordinate bounds its coefficient
of $\prod_h t_h$ by $C''_p\epsilon N$.

Define tensors in the common Euclidean space $H:=(\R^d)^{\otimes b}$ by
\[
 T_i:=\sum_{\mathbf a}(w_{i,\mathbf a}-1)
         v_{i,1,a_1}\otimes\cdots\otimes v_{i,b,a_b}.
\]
Let $W_i:=\bigotimes_{h=1}^b U_{i,h}$ and let $P_i$ be the
orthogonal projection onto $W_i$. The mixed coefficient is
\[
 \frac{p!}{b!}
 \langle\widetilde T_\ell,y_1\otimes\cdots\otimes y_b\rangle,
 \qquad \widetilde T_\ell:=P_\ell\sum_i K_{\ell i}T_i .
\]
The identity rows contribute zero to this coefficient. If
$\iota:\bigoplus_i W_i\to H^g$ is the isometric inclusion, the map
$T\mapsto\widetilde T$ is $\iota^*(K\otimes I_H)\iota\succeq I/8$.
Thus $\sum_i\|T_i\|_F^2\le64\sum_\ell\|\widetilde T_\ell\|_F^2$.

\emph{Counting omissions.}
Take $k_1$ to be the largest block dimension. For each choice of
basis vectors in blocks $2,\ldots,b$, the vector obtained by
contracting $\widetilde T_\ell$ against those basis vectors has
norm at most $C_p\epsilon N$. This follows from the coefficient
bound by maximizing over the unit vector in block one. Summing
squares over these choices gives
$\|\widetilde T_\ell\|_F^2\le C_p(N/k_1)\epsilon^2 N^2$. Since each
$T_i$ is expressed in an orthonormal product basis, we conclude that
\[
 \sum_{i,\mathbf a}(w_{i,\mathbf a}-1)^2
 \le C_p gN\,\epsilon^2 N^2/k_1.
\]
Each omitted row contributes one to the left side. In particular,
if $\epsilon^2 N^2/k_1\le a_p$ for a sufficiently small $a_p>0$,
at least $gN/2$ rows survive.

To meet this inequality, define
$k_1:=\lfloor(a_p/4)\min\{d,\epsilon^{-2}\}\rfloor$ and
$B:=\min\{k_1^{b-1},\sqrt{a_p}/(\epsilon\sqrt{k_1})\}\ge1$.
Indeed, $k_1\ge1$ in the range under consideration and
$\epsilon^2k_1\le a_p/4$.
Choose the remaining $k_h\le k_1$ with product in $[B/2,B]$.
Indeed, use as many factors $k_1$ as possible without exceeding $B$,
then round down the remaining ratio for one more factor and fill
the remaining slots with ones. If $k_1=1$, take every factor one.
This choice gives $\epsilon^2 N^2/k_1\le a_p$ and
$N\asymp_p\min\{k_1^b,\sqrt{k_1}/\epsilon\}$. If $\epsilon^{-2}\le d$,
then $N\asymp_p\epsilon^{-2}$. Otherwise
$N\asymp_p\min\{d^b,\sqrt d/\epsilon\}$. Using $g\asymp_p d^{p/2}$
and $d\asymp_p r$ proves the result.

\emph{Part (b).}
Take $\epsilon$ sufficiently small in terms of $p$ and
$k:=\lceil C'_p\log(2/\epsilon)\rceil$. Lemma~\ref{lem:pre_lww} gives a rank-$k$
instance with $t\ge c'_p k/(\epsilon^2\log^{b_p}(2/\epsilon))$
at accuracy $3\epsilon$, for constants $c'_p,b_p>0$. We choose $C_p$
large enough that $r\ge d_p k$ and apply
Lemma~\ref{lem:lp_lower_tensor_product}. The factor $k^{1-p/2}$
in the resulting bound is absorbed into $\log^{a_p}(2/\epsilon)$.
For the remaining accuracies, Lemma~\ref{lem:lp_noneven_lower}
with $s\ge p/2$ gives $\Omega_p(r^{p/2})$, which suffices after
decreasing $c_p$.
\end{proof}

\subsection{A Support Bound from the Row Kernel}

We start with the argument used for the two Gram-matrix constructions.
We write $a_i^\top$ for row $i$ of $A$ and use weights $w_i$ on the
terms $|a_i^\top x|^p$. The corresponding diagonal row multiplier
is $D_{ii}:=w_i^{1/p}$.

\begin{lemma}[Support bound from the row kernel]\label{lem:lower_kernel_criterion}
Let $p>2$ and let $A\in\R^{N\times r}$ have unit Euclidean row norms. Set $K_{ij}:=|a_i^\top a_j|^p$. Suppose that for some $B\ge1$, $\kappa>0$, and $0\le\mu\le\min\{1,Br^{-p/2}\}$,
\[
 K\succeq\mu\mathbf1\mathbf1^\top+\kappa I_N,
 \qquad \|h_+\|_2^2\le BN,
\]
where $h:=K\mathbf1-(\mu N+1-\mu)\mathbf1$ and $(h_+)_i:=\max\{h_i,0\}$. For every $0<\epsilon\le1/2$, any nonnegative weights preserving $\|Ax\|_p^p$ within $1\pm\epsilon$ for every $x$ have support size at least
\[
 c_{B,\kappa}\frac{N}{1+\epsilon N/r^{p/2}},
\]
where $c_{B,\kappa}>0$ depends only on $B$ and $\kappa$.
\end{lemma}
\begin{proof}
Write $W:=\sum_iw_i$ and $Z:=\|w\|_2^2$. We average the approximation over the Euclidean unit sphere. Each row contributes the same positive average because its Euclidean norm is one. Thus
\begin{equation}\label{eq:lower_kernel_total_weight}
 (1-\epsilon)N\le W\le(1+\epsilon)N.
\end{equation}
We apply the upper approximation at $x=a_j$ to obtain $(Kw)_j\le(1+\epsilon)(K\mathbf1)_j$. Multiplying by $w_j\ge0$ and summing, we get
\begin{equation}\label{eq:lower_kernel_weighted_bound}
 w^\top Kw\le(1+\epsilon)w^\top K\mathbf1.
\end{equation}
The lower bound on $K$ lets us subtract the contribution of its constant component. Using Eq.~\eqref{eq:lower_kernel_weighted_bound} and the definition of $h$, we obtain
\begin{align}
 \kappa Z \le & ~ (1+\epsilon)w^\top K\mathbf1-\mu W^2\notag\\
 = & ~ \mu W((1+\epsilon)N-W)
   +(1+\epsilon)(1-\mu)W+(1+\epsilon)w^\top h.
\label{eq:lower_weight_second_moment_chain}
\end{align}
By Eq.~\eqref{eq:lower_kernel_total_weight} and the bound on $\mu$, the first term is at most $3B\epsilon N^2/r^{p/2}$. The second is at most $9N/4$ because $\mu\ge0$. We bound the last term using $w\ge0$, Cauchy--Schwarz, and Young's inequality.
\[
 (1+\epsilon)w^\top h\le\tfrac32\sqrt{BNZ}\le\tfrac\kappa2 Z+\frac{9BN}{8\kappa}.
\]
Substituting these bounds into Eq.~\eqref{eq:lower_weight_second_moment_chain} gives $Z\le C_{B,\kappa}(N+\epsilon N^2/r^{p/2})$. Finally, Cauchy--Schwarz gives $|\supp(w)|\ge W^2/Z$, and Eq.~\eqref{eq:lower_kernel_total_weight} gives $W\ge N/2$.
\end{proof}

\subsection{Mutually Unbiased Bases}

\begin{lemma}[Lower bound for real exponents]\label{lem:lp_real_lower}
Fix a real $p>2$. For every integer $r\ge2$, there exists a rank-$r$ matrix $A\in\R^{N\times r}$ such that, for every $0<\epsilon\le1/2$, any nonnegative weights preserving $\|Ax\|_p^p$ within $1\pm\epsilon$ for every $x\in\R^r$ have support size at least
\[
 c_p\min\{r^{p/2}/\epsilon,r^2\},
\]
where $c_p>0$ depends only on $p$.
\end{lemma}
\begin{proof}
Let $q$ be the largest power of four at most $r$. We first assume $q\ge4$ and $q^{1-p/2}\le1/2$. Lemma~\ref{lem:pre_mub} gives $q/2+1$ real orthonormal bases in dimension $q$ such that vectors from distinct bases have inner products of magnitude $q^{-1/2}$. We form $A_0\in\R^{N_0\times q}$ from all $N_0:=q(q/2+1)$ basis vectors.

Set $\mu:=q^{-p/2}$ and let $J_q$ denote the $q\times q$ all-ones matrix. Rows from different bases contribute $\mu$ to the kernel. Within a basis, the kernel is the identity. Hence
\[
 K=\mu\mathbf1\mathbf1^\top+I_{N_0}-\mu\operatorname{diag}(J_q,\ldots,J_q)
 \succeq\mu\mathbf1\mathbf1^\top+\tfrac12 I_{N_0},
\]
where the inequality uses $\mu q\le1/2$. We have $K\mathbf1=(1+\mu(N_0-q))\mathbf1$, so $h=-\mu(q-1)\mathbf1$ in Lemma~\ref{lem:lower_kernel_criterion}. Since $h_+=0$, the lemma gives
\[
 |\supp(w)|\ge c\frac{N_0}{1+\epsilon N_0/q^{p/2}}
 \ge c'\min\{q^{p/2}/\epsilon,q^2\},
\]
where we used $N_0\ge q^2/2$. If $r=q$, we take $A:=A_0$. If $r>q$, we take $A:=\operatorname{diag}(A_0,I_{r-q})$ and restrict the approximation to vectors supported in the first $q$ coordinates. The same lower bound holds for $A$, and $q>r/4$ gives the claimed bound in terms of $r$.

If $q<4$ or $q^{1-p/2}>1/2$, then $r$ is bounded by a constant depending only on $p$. We take $A:=I_r$ in these cases. Each coordinate query forces one nonzero weight, and decreasing $c_p$ gives the asserted bound since the minimum is at most $r^2$.
\end{proof}

\subsection{Isolating Rows for Non-Even Exponents}

We isolate a row with a signed combination of real queries. We use a family of rows that agree in one coordinate and have finitely many inner-product magnitudes. Shifting the queries in the shared coordinate makes every other row analytic on the query interval.

\begin{lemma}[Lower bound for non-even exponents]\label{lem:lp_noneven_lower}
Fix a non-even $p>2$ and an integer $s\ge2$. For every integer $r\ge2$ and every $0<\epsilon\le1/2$, there exists a rank-$r$ matrix $A\in\R^{N\times r}$ such that any nonnegative weights preserving $\|Ax\|_p^p$ within $1\pm\epsilon$ for every $x\in\R^r$ have support size at least
\[
 c_{p,s}\min\{r^{p/2}/\epsilon,r^s\},
\]
where $c_{p,s}>0$ depends only on $p$ and $s$.
\end{lemma}
\begin{proof}
Let $q$ be the largest power of two with an odd exponent and $q\le r$. We first work in dimension $q$, assuming it is sufficiently large in terms of $p$ and $s$. We use the rows in Lemma~\ref{lem:finite_inner_products}, including an orthonormal basis, to form an $N\times q$ matrix $A_0$. Write $D_q:=q^{p/2}$, let $e_0$ be the coordinate vector from the lemma, and set $B:=2^s+1$.

\emph{A signed combination of queries.} Define $f(t):=|t|^p$ on $[-1,1]$. Let $V$ be the span of the functions $t\mapsto(B+ut)^p$ for $u\in U_s\cup(-U_s)$, where $U_s$ is the finite set in Lemma~\ref{lem:finite_inner_products}. All functions in $V$ are real analytic on a neighborhood of $[-1,1]$ since $B>|u|$. The singularity of $f$ at zero for non-even $p$ gives $f\notin V$. Point evaluations span the dual of $V+\R f$ because a proper span would give a nonzero function in $V+\R f$ vanishing at every point of $[-1,1]$. We can therefore choose finitely many $t_\ell\in[-1,1]$ and coefficients $b_\ell\in\R$, depending only on $p,s$, such that the functional $L(h):=\sum_\ell b_\ell h(t_\ell)$ vanishes on $V$, $L(f)=\lambda>0$, and $\sum_\ell|b_\ell|=1$.

For row $a_j^\top$, we use queries of the form $x=t a_j+B e_0$. Every row has $a_i^\top e_0=q^{-1/2}$, so
\[
 |a_i^\top x|^p=q^{-p/2}|B+\sqrt q~a_i^\top a_j~t|^p.
\]
For $i\ne j$ this function of $t$ belongs to $V$ because $\sqrt q~|a_i^\top a_j|\in U_s$ and $B>2^s$. For $i=j$, its image under $L$ is
\[
 \lambda_q:=L(t\mapsto|t+B/\sqrt q|^p).
\]
Since $L$ uses finitely many points, $\lambda_q\to\lambda$. Thus $\lambda_q\ge\lambda/2$ for sufficiently large $q$.

\emph{The cost of a query.} The same magnitude bound and $|t|\le1$ give, for all $j$,
\[
 \|A_0(t a_j+B e_0)\|_p^p
 \le C_{p,s}(1+N/D_q).
\]
Indeed, the term $i=j$ is bounded by $(1+B/\sqrt q)^p$. Each other term is at most $(B+2^s)^p/D_q$. Applying the assumed relative error at the resulting queries and taking the signed sum with coefficients $b_\ell$ yields
\begin{equation}\label{eq:lp_noneven_weight_error}
 |w_j-1|\lambda_q\le C_{p,s}\epsilon(1+N/D_q).
\end{equation}
Here the left-hand side contains only the term $i=j$ because $L$ vanishes on $V$, and the right-hand side uses $\sum_\ell|b_\ell|=1$.

Choose $0<\delta\le1/2$ such that $2C_{p,s}\delta<\lambda/2$. For $\epsilon\le\delta$, take $N:=\min\{q^s,\lfloor\delta D_q/\epsilon\rfloor\}$. We have $q\le N$ for sufficiently large $q$, because $p>2$. Eq.~\eqref{eq:lp_noneven_weight_error} gives $|w_j-1|<1$ for every $j$. Thus all $N$ rows must be kept.

For $\delta<\epsilon\le1/2$, we take $N:=\lfloor\min\{q^s,D_q\}\rfloor$. Since $N\le D_q$, evaluating at $x=a_j$ gives $w_j\le(1+\epsilon)\|A_0 a_j\|_p^p\le C'_{p,s}$ using the same inner-product bound. Sphere averaging as in Eq.~\eqref{eq:lower_kernel_total_weight} gives $\sum_j w_j\ge N/2$. We obtain $|\supp(w)|\ge N/(2C'_{p,s})$, which is the desired rate because $\epsilon>\delta$.

The matrix $A_0$ has rank $q$ since it includes an orthonormal basis. When $r>q$, we take $A:=\operatorname{diag}(A_0,I_{r-q})$ and restrict to queries on the first $q$ coordinates. Since $q>r/4$, this proves the theorem for sufficiently large $r$. For the remaining ranks, the identity matrix forces $r$ nonzero weights and we decrease $c_{p,s}$ using $\min\{r^{p/2}/\epsilon,r^s\}\le r^s$.
\end{proof}

\begin{lemma}[A finite set of inner-product magnitudes]\label{lem:finite_inner_products}
Let $s\ge2$ and $q=2^k$ with odd $k>2s$. There is a family of $q^s$ Euclidean unit vectors in $\R^q$ that contains an orthonormal basis. There is a coordinate vector $e_0$ for which every member $a$ has $a^\top e_0=q^{-1/2}$. For any distinct members $a,b$,
\[
 \sqrt q~|a^\top b|\in U_s:=\{0\}\cup\{2^{d/2}:0\le d\le2s,\ d\in\mathbb Z\}.
\]
\end{lemma}
\begin{proof}
Let $\mathbb F_q$ be the field of $q$ elements and let $\operatorname{Tr}:\mathbb F_q\to\mathbb F_2$ be its trace. For $v=(v_1,\ldots,v_s)\in\mathbb F_q^s$, define
\[
 Q_v(x):=\operatorname{Tr}(\sum_{j=1}^s v_j x^{2^j+1}),\qquad
 a_v:=q^{-1/2}((-1)^{Q_v(x)})_{x\in\mathbb F_q}.
\]
We show that every nonzero $v$ gives a quadratic form whose polar form has a small radical. Using invariance of the trace under squaring, we can write the polar form as
\[
 Q_v(x+y)+Q_v(x)+Q_v(y)=\operatorname{Tr}(yP_v(x)),\qquad
 P_v(x):=\sum_{j=1}^s(v_j x^{2^j}+(v_jx)^{2^{k-j}}).
\]
The trace pairing is nondegenerate, so the radical is $R_v:=\{x:P_v(x)=0\}$. On $\mathbb F_q$, the function $P_v(x)^{2^s}$ equals the polynomial
\[
 \sum_{j=1}^s(v_j^{2^s}x^{2^{s+j}}+v_j^{2^{s-j}}x^{2^{s-j}}).
\]
Its exponents are distinct, and for $v\ne0$ the polynomial is nonzero and has degree at most $2^{2s}$. Hence $R_v$ is an $\mathbb F_2$-linear space of dimension $d_v\le2s$.

Now set $S_v:=\sum_x(-1)^{Q_v(x)}$. In the square of this sum we write the two arguments as $x$ and $x+z$. Summing over $x$ gives zero when $z\notin R_v$ and gives $q(-1)^{Q_v(z)}$ otherwise. Thus
\[
 S_v^2=q\sum_{z\in R_v}(-1)^{Q_v(z)}\in\{0,q2^{d_v}\},
\]
where the last step uses that $Q_v$ is linear on $R_v$. For distinct $u,v$, we have $a_u^\top a_v=S_{u+v}/q$. The bound $d_{u+v}\le2s<k$ ensures the vectors are distinct and gives the claimed magnitudes.

Let $e_0$ be the coordinate vector indexed by the zero element of $\mathbb F_q$. All rows have this coordinate equal to $q^{-1/2}$ because $Q_v(0)=0$. To see that the family contains an orthonormal basis, take the $q$ vectors with $v_j=0$ for $j>1$. The map $x\mapsto x^3$ permutes $\mathbb F_q$ because $k$ is odd and $3$ is coprime to $2^k-1$. Thus the inner product of distinct rows in this collection is a sum of a nontrivial additive character over $\mathbb F_q$, divided by $q$, and equals zero.
\end{proof}

\subsection{Random Sign Rows for Even Exponents}

\begin{lemma}[Lower bound for even exponents from random sign rows]\label{lem:lp_even_lower}
Fix an even integer $p\ge4$. For every integer $r\ge2$, there exists a rank-$r$ matrix $A\in\R^{N\times r}$ such that, for every $0<\epsilon\le1/2$, any nonnegative weights preserving $\|Ax\|_p^p$ within $1\pm\epsilon$ for every $x\in\R^r$ have support size at least
\[
 c_p\min\{r^p,\max\{r^{p/2}/\epsilon,r/\epsilon^2\}\},
\]
where $c_p>0$ depends only on $p$.
\end{lemma}
\begin{proof}
Suppose first that $r$ is sufficiently large in terms of $p$. We take $N:=\lfloor r^p/(4p!)\rfloor$ and the rank-$r$ sign matrix in Lemma~\ref{lem:even_lower_hard_matrix}. Its row kernel satisfies Lemma~\ref{lem:lower_kernel_criterion} with $\kappa:=1/8$ and $B$ depending only on $p$. Hence
\[
 |\supp(w)|\ge c_{1,p} \frac{N}{1+\epsilon N/r^{p/2}}
 \ge c_{2,p}\min\{r^{p/2}/\epsilon,r^p\}.
\]
For the second bound, we use the same matrix and form the symmetric tensor
\[
 T:=\sum_{i=1}^N(w_i-1)a_i^{\otimes p}.
\]
Its diagonal polynomial $P(x):=\langle T,x^{\otimes p}\rangle$ is the difference between the reweighted and original sums, since $p$ is even. We use the bound on $\|Ax\|_p^p$ from Lemma~\ref{lem:even_lower_hard_matrix} to obtain
\[
 \sup_{\|x\|_2=1}|P(x)|\le C_p\epsilon N/r^{p/2}.
\]
Write $T(x_1,\ldots,x_p)$ for the corresponding multilinear form. For Euclidean unit vectors $x_1,\ldots,x_p$, the real polarization identity gives
\[
 T(x_1,\ldots,x_p)=\frac{1}{2^p p!}
 \sum_{\sigma\in\{-1,1\}^p}(\prod_{j=1}^p\sigma_j)
 P(\sum_{j=1}^p\sigma_jx_j).
\]
Write $M:=\sup_{\|x\|_2=1}|P(x)|$. By $p$-homogeneity and the
triangle inequality, every summand has magnitude at most $p^p M$.
It follows that each vector $T(e_{j_1},\ldots,e_{j_{p-1}},\cdot)$ has
Euclidean norm at most $p^p M/p!$. Summing over the $r^{p-1}$
coordinate tuples gives
\[
 \|T\|_F^2\le C'_p r^{p-1}(\epsilon N/r^{p/2})^2
 = C'_p\epsilon^2N^2/r.
\]
The row kernel is the Gram matrix of $a_i^{\otimes p}$. Its eigenvalue bound therefore gives $\|T\|_F^2\ge\|w-\mathbf1\|_2^2/8$. Using $\|w\|_2^2\le2N+2\|w-\mathbf1\|_2^2$ and the mass bound in Eq.~\eqref{eq:lower_kernel_total_weight}, we obtain
\[
 |\supp(w)|\ge\frac{(\sum_i w_i)^2}{\|w\|_2^2}
 \ge c_{3,p}\frac{N}{1+\epsilon^2N/r}
 \ge c_{4,p}\min\{r/\epsilon^2,r^p\}.
\]
Taking the larger of the two lower bounds gives the result. For the remaining ranks, we take $A:=I_r$ and decrease $c_p$, using that the expression is at most $r^p$ and that the rank threshold depends only on $p$.
\end{proof}

\begin{lemma}[Random sign instance for even exponents]\label{lem:even_lower_hard_matrix}
Fix an even integer $p\ge4$ and define $N:=\lfloor r^p/(4p!)\rfloor$. For all sufficiently large $r$ in terms of $p$, there exists a rank-$r$ matrix $A\in\{-r^{-1/2},r^{-1/2}\}^{N\times r}$ with the following properties. Its kernel $K_{ij}:=|a_i^\top a_j|^p$ satisfies
\[
 K\succeq\mu\mathbf1\mathbf1^\top+\tfrac18 I_N,
 \qquad \|K\mathbf1-(\mu N+1-\mu)\mathbf1\|_2^2\le C_pN,
\]
for some $0<\mu\le\min\{1,C_pr^{-p/2}\}$. Moreover,
\[
 \sup_{\|x\|_2=1}\|Ax\|_p^p\le C_p N/r^{p/2}.
\]
Here $C_p$ depends only on $p$.
\end{lemma}
\begin{proof}
We draw $\xi_1,\ldots,\xi_N$ independently from $\{-1,1\}^r$ and set $a_i:=\xi_i/\sqrt r$. For independent uniform signs $\eta_1,\ldots,\eta_r$, set $S_r:=\sum_{j=1}^r\eta_j$ and $\mu:=r^{-p}\E[|S_r|^p]$. In the expansion of an even moment of $S_r$, a term has nonzero expectation only if every index occurs an even number of times. Thus a moment of order $2k$ has at most $C_kr^k$ nonzero terms. We obtain
\begin{equation}\label{eq:lower_sign_moments}
 r^{p/2}\le\E[|S_r|^p]\le C_pr^{p/2},
 \qquad \E[|S_r|^{2p}]\le C_pr^p,
\end{equation}
where the first inequality follows from Jensen's inequality and $\E[S_r^2]=r$. This gives $0<\mu\le C_pr^{-p/2}$, and $|S_r|\le r$ gives $\mu\le1$.

\emph{Kernel lower bound.} For $S\subseteq[r]$, define $\chi_S(\xi):=\prod_{j\in S}\xi_j$ and $v_S:=(\chi_S(\xi_i))_{i=1}^N$. In the expansion of $K_{ij}=r^{-p}(\sum_{\ell=1}^r\xi_{i\ell}\xi_{j\ell})^p$, we group the ordered tuples by the set of indices of odd multiplicity. We obtain
\begin{equation}\label{eq:lower_walsh_expansion}
 K=\sum_{S\subseteq[r]}b_Sv_Sv_S^\top,
 \qquad b_S\ge0,
 \qquad b_\varnothing=\mu,
 \qquad b_S=p!/r^p\text{ for }|S|=p.
\end{equation}
Averaging over independent signs selects the empty-set coefficient, and a set of size $p$ arises from its $p!$ orderings. Let $X$ have columns $v_S$ for $|S|=p$. Lemma~\ref{lem:random_sign_tensor_gap} gives
\begin{equation}\label{eq:lower_tensor_spectral_gap}
 (p!/r^p)XX^\top\succeq\tfrac18 I_N
\end{equation}
with probability tending to one as $r\to\infty$. Keeping the empty-set term and the terms of degree $p$ in Eq.~\eqref{eq:lower_walsh_expansion} gives the desired lower bound on $K$.

\emph{Kernel row sums.} Since $K_{ii}=1$, we have $h_i=\sum_{j\ne i}(K_{ij}-\mu)$. Given $\xi_i$, these summands are independent and centered because $\xi_i$ multiplied coordinatewise by $\xi_j$ has the uniform sign distribution. By Eq.~\eqref{eq:lower_sign_moments}, the conditional variance of each summand is at most $C_pr^{-p}$. Thus
\begin{equation}\label{eq:lower_residual_second_moment}
 \E[\|h\|_2^2]\le C_p N(N-1)/r^p\le C_pN.
\end{equation}
By Markov's inequality (Lemma~\ref{lem:pre_markov}), $\|h\|_2^2\le4 C_pN$ with probability at least $3/4$.

\emph{Uniform bound on the sum.} Lemma~\ref{lem:sign_energy_upper}
gives the uniform energy bound with probability tending to one,
since $N=\Theta_p(r^p)\ge r^{p/2}$ for sufficiently large $r$.

\emph{Full rank.} A proper subspace of $\R^r$ contains at most half of the sign cube. Indeed, it lies in a hyperplane with a nonzero coordinate coefficient, and fixing all other signs leaves at most one choice in that hyperplane. Until we reach full rank, each new row therefore increases the dimension with conditional probability at least $1/2$. The expected number of rows needed for full rank is at most $2r$, so the probability that $N$ rows fail to span is at most $2r/N$. For sufficiently large $r$, the union bound gives a realization satisfying all four properties.
\end{proof}

\begin{lemma}[Uniform energy of random sign rows]\label{lem:sign_energy_upper}
Fix $p>2$. Let $a_1,\ldots,a_N$ be independent and uniform in
$\{-d^{-1/2},d^{-1/2}\}^d$. With probability at least
$1-5^d N e^{-4d}$,
\[
 \sup_{\|x\|_2=1}\sum_{i=1}^N |a_i^\top x|^p
 \le C_p(1+N/d^{p/2}).
\]
\end{lemma}
\begin{proof}
Take a $1/2$-net $\mathcal N$ of the Euclidean unit sphere of size
at most $5^d$. For fixed $x\in\mathcal N$, each variable
$|\sqrt d~ a_i^\top x|$ has tail probability at most $2e^{-t^2/2}$
at threshold $t$. Thus the probability that its $j$th largest value
over $i\in[N]$ exceeds $t$ is at most $(2eN/j)^j e^{-jt^2/2}$.
Taking $t_j:=\sqrt{2\log(2eN/j)+8d/j}$ and a union bound over $x,j$,
all order statistics satisfy these bounds except with probability
$5^dNe^{-4d}$. For every $x\in\mathcal N$ on this event,
\[
 \sum_i|\sqrt d~ a_i^\top x|^p
 \le\sum_j t_j^p\le C_p(N+d^{p/2}).
\]
The last inequality follows by bounding the logarithm sum by
$N\int_0^1(\log(2e/u))^{p/2}\d u$ and the other sum by
$C_p d^{p/2}\sum_{j\ge1}j^{-p/2}$. For the matrix $A$ of these
rows, the triangle inequality gives
$\sup_{\|x\|_2=1}\|Ax\|_p\le2\max_{x\in\mathcal N}\|Ax\|_p$.
\end{proof}

\subsection{Amplifying Lower Bounds by Tensor Products}

We can amplify a low-rank hard instance by tensoring it with a family
of nearly orthogonal vectors. On product queries the contributions of
the different copies are mixed by a diagonally dominant matrix.
Inverting this matrix gives an approximation guarantee within each copy.

\begin{lemma}[Tensor-product amplification]\label{lem:lp_lower_tensor_product}
Fix $p>2$. There are constants $c_p>0$ and $d_p\ge2$ with the
following property. Let $0<\epsilon\le1/6$, let $B$ be a matrix of
rank $k$, and let $r\ge d_p k$. Suppose any $(1\pm3\epsilon)$
approximation of the energy of $B$ by reweighted original rows
requires at least $t$ rows. Then some rank-$r$ matrix requires at least
$c_p(r/k)^{p/2}t$ rows at accuracy $\epsilon$.
\end{lemma}
\begin{proof}
We work in coordinates on the row span of $B$ and write its rows as
$b_j^\top$, $j\in[m]$. Define $s:=\max\{2,\lceil p/2\rceil\}$. Choose
$d=2^h$ with $h$ odd and $d$ as large as possible subject to
$d\le r/k$. Then $r/(4k)<d\le r/k$. Choose $d_p$ large enough for
Lemma~\ref{lem:finite_inner_products} to apply in dimension $d$. Its
family contains a basis and at least $d^{p/2}$ unit vectors, with
$|u_i^\top u_j|\le2^s/\sqrt d$ for distinct vectors. We keep
$g:=\lfloor 2^{-sp-1}d^{p/2}\rfloor$ vectors, including a basis.
We have $g\ge d$ by increasing $d_p$, since $p/2>1$.

Use the rows $(u_i\otimes b_j)^\top$ for $i\in[g]$ and $j\in[m]$.
They span $\R^{dk}$. If $dk<r$, pad the instance by a direct sum
with $I_{r-dk}$ and restrict all queries below to the first $dk$
coordinates. Define the $g\times g$ matrix $K$ by
$K_{\ell i}:=|u_\ell^\top u_i|^p$. Its diagonal entries are one,
and every off-diagonal row sum is at most
$(g-1)2^{sp}d^{-p/2}\le1/2$. Hence the operator norm of $K-I_g$
on $\ell_\infty^g$ is at most $1/2$.

For any nonnegative sparsifier weights $w_{ij}$, write
\[
 F(y):=\|By\|_p^p,\qquad
 e_i(y):=\sum_j w_{ij}|b_j^\top y|^p-F(y).
\]
At a query $u_\ell\otimes y$, the input energy is
$(K\mathbf1)_\ell F(y)$ and the signed error is $(Ke(y))_\ell$.
The relative approximation at all $g$ such queries therefore gives
$\|Ke(y)\|_\infty\le(3/2)\epsilon F(y)$. Since
$\|Ke(y)\|_\infty\ge\|e(y)\|_\infty/2$, each of the $g$ copies
approximates $B$ within $1\pm3\epsilon$ for every $y$.
Each needs at least $t$ rows, so the support is at least
$gt\ge c_p(r/k)^{p/2}t$.
\end{proof}

\subsection{Random Orthonormal Frames}

\begin{lemma}[Uniform energy of random orthonormal frames]\label{lem:rotated_frame_energy}
Fix $p>2$. For $1\le k\le d$, let $V_1,\ldots,V_g$ be independent
$d\times k$ matrices, each consisting of the first $k$ columns of a
Haar orthogonal matrix. With probability at least $1-5^d g e^{-4d}$,
\[
 \sup_{\|y\|_2=1}\sum_{i=1}^g\|V_i^\top y\|_p^p
 \le C_p(1+gk/d^{p/2}).
\]
\end{lemma}
\begin{proof}
Fix a unit vector $y$. Each $V_i^\top y$ has the distribution of the
first $k$ coordinates of a uniform point on the unit sphere. For a
standard Gaussian $G$ in $\R^d$, its direction is independent of its
radius, and $\Pr[\|G\|_2\ge\sqrt d/2]$ is bounded below by an absolute
constant by Lemma~\ref{lem:pre_paley_zygmund}. The function $G\mapsto\|G_{[k]}\|_p$ is
$1$-Lipschitz and has mean at most $C_p k^{1/p}$. Gaussian concentration
and independence of the radius therefore give
\[
 \Pr[\sqrt d~\|V_i^\top y\|_p>2C_p k^{1/p}+t]\le C e^{-ct^2}
 \qquad(t\ge0)
\]
for absolute $C,c>0$.

Let $Z_i:=(\sqrt d~\|V_i^\top y\|_p-2C_p k^{1/p})_+$.
For each $1\le h\le g$, a union bound over the $h$-subsets of $[g]$
bounds the probability that the $h$th largest $Z_i$ exceeds $t$ by
$(Ce g/h)^h e^{-cht^2}$. Use
$t_h^2:=c^{-1}(\log(Ce g/h)+4d/h)$ and take a union bound over $h$
and a $1/2$-net of the unit sphere of size at most $5^d$. On the
resulting event, $\sum_h t_h^p\le C_p(g+d^{p/2})$, by the same
logarithm and power summations as in Lemma~\ref{lem:sign_energy_upper}.
The claimed upper bound holds on the net because
\[
 \sum_i\|V_i^\top y\|_p^p
 \le 2^{p-1}d^{-p/2}(g(2C_p)^p k+\sum_i Z_i^p).
\]
Taking the $p$th root and using the triangle inequality extends the
bound to the sphere, with a factor at most $2^p$ after taking the
$p$th power again.
\end{proof}

%% file: 80_proof.tex
\section{Joint Rank and Accuracy Below Two}
\label{sec:joint_accuracy}

\subsection{Main Result}

Theorem~\ref{thm:below_two_main} determines the worst-case support at
polynomial accuracy and bounds it at subexponential accuracy. In
part~(a) the upper bound is Theorem~\ref{thm:general_lp_field}, and the
lower bound flips the blocks of a partite sign dictionary to extract
the top coefficient of the query energy.
In part~(b), the upper bound rounds the Lewis masses by iterated partial coloring. On
each cap of low-energy queries, linear constraints preserve the polynomial
part of every row energy up to a sketching error, and the Gaussian width
of the remainder controls the rest. Sketching the polynomial constraints
keeps their number small, so the degree can grow to order $L$. The lower
bound shares the first-block query coordinates across several groups of
transverse coordinates. On a query for one group, the errors of the other
groups cancel, but their energy still enters the error tolerance.

\begin{theorem}[Polynomial and subexponential accuracy below two]
\label{thm:below_two_main}
Fix $0<p<2$.
\begin{enumerate}[label=(\alph*),nosep]
\item For every fixed $a>0$, uniformly for $r^{-a}\le\epsilon\le1/2$,
\[
 S_p(r,\epsilon)=\Theta_{p,a}(r/\epsilon^2).
\]
\item For $C_p\log(2r)\le L\le r$,
\[
 S_p(r,\epsilon)\lesssim_p
 \epsilon^{-2}\min\{r,\frac{r^3}{K^{2-p}L^{2p}}\}.
\]
For $C_p\log r\le L\le c_\circ r$, define $b:=1+p/2$. Then
\[
 S_p(r,\epsilon)\gtrsim_p \frac{r}{\epsilon^2}
 \begin{cases}
 (\log(e r/L^b)/L)^{(p+1)/2},& L^b\le r,\\
 r/L^{p+3/2},& L^b\ge r.
 \end{cases}
\]
Here $c_\circ>0$ is an absolute constant.
\end{enumerate}
\end{theorem}
\begin{proof}
\emph{Part (a).}
The upper bound is Theorem~\ref{thm:general_lp_field}. For the lower
bound choose a fixed even integer $m>2a+1$ and define
$n:=\lfloor r/m\rfloor$. Use $m$ blocks of $n$ coordinates and a row
$a_{\mathbf i}=\sum_{j=1}^m e_{j,i_j}$ for each tuple
$\mathbf i\in[n]^m$. In each slice indexed by $i_1$, give $d$ rows
mass one and distribute at most another $d$ units of mass positively
among the other rows. Include identity rows of mass $1/m$ to ensure
rank $r$. As usual, a row with mass $w$ means a row scaled by $w^{1/p}$.

Take a Boolean query $x=(x^{(1)},\ldots,x^{(m)})$, with remaining
coordinates zero. Flip its blocks by independent signs
$\sigma_1,\ldots,\sigma_m$ and average the signed error against
$\prod_j\sigma_j$. This removes the identity-row terms and extracts
\[
 \lambda_{p,m}\sum_{\mathbf i}\Delta_{\mathbf i}
                               \prod_{j=1}^m x^{(j)}_{i_j},
 \qquad
 \lambda_{p,m}=\E_\sigma[
        (\prod_j\sigma_j)|\sum_j\sigma_j|^p],
\]
where $\Delta$ is the output mass minus the input mass. The integral
formula for $|t|^p$ in terms of $1-\cos(st)$ expresses
$|\lambda_{p,m}|$ as a positive constant times
$\int_0^\infty\sin^m(s)s^{-p-1}\,\mathrm ds$. This is finite and
positive because $m>p$ is even. The input energy on these probes is
$O_m(nd)$ for $d\ge1$, so relative accuracy implies that the extracted
Boolean tensor has supremum $O_{p,m}(\epsilon nd)$.

Conversely, choose the first block of signs to maximize the
contribution from all slices, and average over the other blocks.
For each slice its decoupled Rademacher sum has second moment equal
to the squared Euclidean norm of the slice and fourth moment at
most $3^{m-1}$ times the square of that second moment, by successive
conditioning. This gives a supremum at least
$c_m\sum_{i_1}\|\Delta_{i_1,\cdot}\|_2$. If fewer than $nd/2$
mass-one rows survive, the sum is at least $n\sqrt d/2$. Indeed,
if $d_{i_1}$ such rows are missing from a slice, then the slice norm
is at least $\sqrt{d_{i_1}}\ge d_{i_1}/\sqrt d$.
Choose $d:=\lfloor c_{p,m}/\epsilon^2\rfloor$ with $c_{p,m}$ small.
For sufficiently large $r$, the assumption on $\epsilon$ gives
$d\le n^{m-1}$. The upper and lower estimates for the extracted
supremum are incompatible. If $d$ is bounded or zero, the rank
lower bound suffices after adjusting the constant.

\emph{Part (b), upper bound.}
Theorem~\ref{thm:general_lp_field} gives $S_p(r,\epsilon)\le
Cr/(\min\{p,1\}\epsilon^2)$, so it suffices to find sparsifiers with
$O_p(\epsilon^{-2}r^3K^{p-2}L^{-2p})$ rows. Zero rows do not matter, so
we may use the Lewis setup of the common tools, Eq.~\eqref{lem:t-eq-lewis},
with rows $i\in I$, unit directions $u_i$ and masses $\nu_i$. By
homogeneity and the last sentence of the subsection on the Lewis setup, it
suffices to find $\mu^*\in[0,\infty)^I$ with $|\supp\mu^*|\le N_*$, for
the number $N_*$ chosen in Step~3, such that
\begin{equation}\label{lem:x-subexp-eq-goal}
 |\sum_{i\in I}\mu^*_i|\langle u_i,x\rangle|^p-F(x)|\le\epsilon
 F(x)\qquad(x\in S^{r-1}).
\end{equation}
The argument is the same for all $0<p<2$. Below one, the caps come from the
ensemble cover of Lemma~\ref{lem:t-cap-cover}(c), and
Lemma~\ref{lem:t-shifted-tv}(c) holds for every $p>0$.

\emph{Parameters.} Let $c_1\in(0,1]$ depend only on $p$ and $c_2\in(0,1]$
be absolute, both fixed in Step~4, and define
\[
 k:=\lfloor c_1K\rfloor,\qquad q:=\lfloor c_2L\rfloor,\qquad
 \tau:=(2r)^{-3}.
\]
Step~4 shows $1\le k\le r$ and $2\le q\le r$. Let $\mathcal C$ be the cap
family of Lemma~\ref{lem:t-cap-family} for this $k$. For $c\in\mathcal C$
and $i\in I$ let $Q_{c,i}$, $R_{c,i}$ and $s_{c,i}$ be as in
Lemma~\ref{lem:t-cap-width}, and let $Y\subseteq B_2^r$ be the norming set
of Lemma~\ref{lem:t-norming-set}(b) for degree $q$.

\emph{Step 1: sketches.} For $c\in\mathcal C$ we have
$|\langle u_i,v_c\rangle|\le\|v_c\|_2\le2$ and $\varrho_c\le1$, so
Lemma~\ref{lem:t-shifted-tv}(c) gives $\sup_{|s|\le1}|Q_{c,i}(s)|\le
m_p:=3^p+C''_p$, with $C''_p$ the constant of that lemma. For each cap
$c$ we apply Lemma~\ref{lem:t-sketch} with $a=1$, $r'=r$, the rows $i\in I$,
$w_{i,1}=u_i$, the polynomials $P_i=Q_{c,i}/m_p$ of degree at most $q$,
the set $Y$, $\mathsf A=3^{q+1}$ (allowed since
$\sup_{[-1,1]}|Q_{c,i}/m_p|\le1$) and $\tau=(2r)^{-3}$. Since
$\log(e|Y|)\le1+(r+1)\log(8(q+1)^2)\le15r\log(2r)$, we have
$\log\log(e|Y|)\le4\log(2r)$. This gives a matrix $\widetilde W^c$
indexed by $I\times Y$ of rank at most $d_{\rm sk}$, which is $1+(q+1)d_0$
in the notation of that proof, such that
\[
 |Q_{c,i}(\langle u_i,y\rangle)-m_p\widetilde W^c_{i,y}|\le m_p\tau
 \quad(i\in I,\ y\in Y),\qquad
 \log(2d_{\rm sk})\le C(q+4\log(2r)).
\]

\emph{Step 2: one round.} Let $\mu\in[0,\infty)^I$ have total mass one and
$N:=|\supp\mu|$ positive entries. Let $J:=\{i:0<\mu_i\le2/N\}$ and
$\mu_J:=(\mu_i)_{i\in J}$. By Lemma~\ref{lem:t-light-rows},
$|J|\ge N/2$ and $\|\mu_J\|_2^2\le2/N$. For $c\in\mathcal C$ let
$\mathcal B_c\subseteq\R^J$ be a basis of the span of the columns
$(\widetilde W^c_{i,y})_{i\in J}$, $y\in Y$, so that $|\mathcal B_c|\le
d_{\rm sk}$. Let $H$ be the subspace of $z\in\R^J$ with
\[
 \sum_{i\in J}\mu_iz_i=0\quad\text{and}\quad\sum_{i\in J}\mu_iz_i\xi_i=0
 \ \text{ for all }\xi\in\textstyle\bigcup_{c\in\mathcal C}\mathcal B_c .
\]
Its codimension is at most $1+|\mathcal C|d_{\rm sk}$. For $z\in\R^J$ define
\[
 \Psi(z):=\sup_{x\in S^{r-1}}\frac1{F(x)}|\sum_{i\in J}\mu_iz_i
 |\langle u_i,x\rangle|^p| .
\]
This is a seminorm on $\R^J$, finite because $F\ge1/r$ on the sphere,
hence continuous. We claim that a standard Gaussian vector $G$ in $H$
satisfies
\begin{equation}\label{lem:x-subexp-eq-round}
 \E[\Psi(G)]\le C_{\rm rd}\,r^{3/2}k^{p/2-1}q^{-p}N^{-1/2},
\end{equation}
where $C_{\rm rd}\ge1$ depends only on $p$.

For $c\in\mathcal C$ consider the random polynomial
$\mathcal Q_c(y):=\sum_{i\in J}\mu_iG_iQ_{c,i}(\langle u_i,y\rangle)$ on
$\R^r$, of degree at most $q$. Fix $y\in Y$. The column
$(\widetilde W^c_{i,y})_{i\in J}$ lies in the span of $\mathcal B_c$, so
$\sum_{i\in J}\mu_iG_i\widetilde W^c_{i,y}=0$ because $G\in H$. Hence
$\mathcal Q_c(y)=\langle G,\omega\rangle$ with
$\omega_i:=\mu_i(Q_{c,i}(\langle u_i,y\rangle)-m_p\widetilde W^c_{i,y})$.
This is a centered Gaussian variable with variance at most
$\|\omega\|_2^2\le m_p^2\tau^2\|\mu_J\|_2^2$. By
Lemma~\ref{lem:t-norming-set}(b), applied to each realization, and
Lemma~\ref{lem:t-sidak-chaining}(a),
\[
 \E[\max_{c\in\mathcal C}\sup_{y\in B_2^r}|\mathcal Q_c(y)|]
 \le2\,\E[\max_{c\in\mathcal C,\,y\in Y}|\mathcal Q_c(y)|]
 \le2m_p\tau\|\mu_J\|_2\sqrt{2\log(2|\mathcal C||Y|)} .
\]
Since $x\mapsto x\log(er/x)$ increases on $(0,r]$, we have
$k\log(er/k)\le r$, and Lemma~\ref{lem:t-cap-family}(i) gives
$\log(2|\mathcal C|)\le C_pr$. Also
$\log|Y|\le(r+1)\log(8(q+1)^2)\le14r\log(2r)$. Since $\log(2r)\le r$ and
$\tau=(2r)^{-3}$, the right side above is at most $C_pr^{-2}N^{-1/2}$.

Now let $c\in\mathcal C$ and $x\in X_c$, and define
$y:=(x-v_c)/\varrho_c\in B_2^r$. Then $s_{c,i}(x)=\langle u_i,y\rangle$,
and therefore
\[
 \sum_{i\in J}\mu_iG_i|\langle u_i,x\rangle|^p=\mathcal Q_c(y)
 +\sum_{i\in J}\mu_iG_iR_{c,i}(s_{c,i}(x)).
\]
Every $x\in S^{r-1}$ lies in some $X_c$, and $F(x)\ge1/r$. Hence the
bound on $\mathcal Q_c$ and Eq.~\eqref{lem:t-eq-round} give
\[
 \E[\Psi(G)]\le C_pr^{-1}N^{-1/2}+C_p\,rk^{p/2-1}q^{-p}\sqrt{2/N}\,
 (\sqrt r+\sqrt{C_pr}).
\]
Since $k,q\le r$, we have $r^{3/2}k^{p/2-1}q^{-p}\ge r^{(1-p)/2}\ge
r^{-1}$, and Eq.~\eqref{lem:x-subexp-eq-round} follows.

\emph{Step 3: iteration.} Let $\delta_0,\delta_1,n_{\rm c}$ be the constants of
Lemma~\ref{lem:t-subspace-coloring}, and define
\[
 a_*:=2C_{\rm rd}\,r^{3/2}k^{p/2-1}q^{-p},\qquad
 C_*:=(\frac{2C_{\rm rd}}{1-(1-\delta_1/10)^{1/2}})^2,\qquad
 N_*:=\lceil\frac{C_*r^3k^{p-2}q^{-2p}}{\epsilon^2}\rceil .
\]
Step~4 shows that
\begin{equation}\label{lem:x-subexp-eq-budget}
 N_*\ge\max\{5n_{\rm c},\ 20,\ 2/\delta_0,\ 2(1+|\mathcal C|d_{\rm sk})
 /\delta_0\}.
\end{equation}
We apply Lemma~\ref{lem:t-iterated-coloring} with the index set $I$,
$\mathcal X:=S^{r-1}$, $\Phi_x(\mu):=\sum_i\mu_i|\langle u_i,x\rangle|^p$,
$\omega:=F$, $\mathcal S:=\{\mu\in[0,\infty)^I:\sum_i\mu_i=1\}$,
$b_*:=1/2$, and $a_*,N_*$ as above. Let $\mu\in\mathcal S$ have $N>N_*$
positive entries. The units are the singletons of $J$, so there are
$|J|\ge N/2$ of them. Take $H$ from Step~2 and
$\mathcal K:=\{z\in H:\Psi(z)\le2\E[\Psi(G)]\}$, which is closed, symmetric
and convex. By Eq.~\eqref{lem:x-subexp-eq-budget} and $|J|>N_*/2$, the
codimension of $H$ is at most $\delta_0N_*/2<\delta_0|J|$, and $|J|\ge
1/\delta_0$. Lemma~\ref{lem:t-light-rows} then gives
$\gamma_H(\mathcal K)\ge1/2\ge e^{-\delta_0|J|}$. For
$z\in\mathcal K\cap[-1,1]^J$ the masses $\mu^z_i=\mu_i(1+z_i)$ are
nonnegative and have total mass one, because $\sum_{i\in J}\mu_iz_i=0$.
Thus $\mu^z\in\mathcal S$, and Eq.~\eqref{lem:x-subexp-eq-round} gives
\[
 |\Phi_x(\mu^z)-\Phi_x(\mu)|=|\sum_{i\in J}\mu_iz_i|\langle u_i,x
 \rangle|^p|\le\Psi(z)F(x)\le a_*N^{-1/2}F(x).
\]
Lemma~\ref{lem:t-iterated-coloring} with $\mu^0:=\nu$, for which
$\Phi_x(\nu)=F(x)$, gives $\mu^*\in\mathcal S$ with
$|\supp\mu^*|\le N_*$ and
\[
 |\Phi_x(\mu^*)-F(x)|\le\frac{a_*N_*^{-1/2}}{1-(1-\delta_1/10)^{1/2}}F(x)
 \le\epsilon F(x)
\]
by the choice of $C_*$. This is Eq.~\eqref{lem:x-subexp-eq-goal}.

\emph{Step 4: the parameters.} For $r\ge2$ we have $\log(er)\le2\log(2r)$,
so $K\ge L/\log(er)\ge C_p/2$, and also $L\ge C_p$. For $C_p\ge4/(c_1c_2)$
this gives $c_1K\ge2$ and $c_2L\ge4$, hence $k\ge1$, $q\ge2$,
$k\ge c_1K/2$ and $q\ge c_2L/2$. Moreover $k\le
K\le L\le r$ and $q\le r$. To check Eq.~\eqref{lem:x-subexp-eq-budget}, note
that $er/K=(er/L)\log(er/L)$, so $\log(er/K)\le2\log(er/L)$ and
$K\log(er/K)\le2L$. Monotonicity of $x\mapsto x\log(er/x)$ gives
\[
 k\log(er/k)\le c_1K(\log(er/K)+\log(1/c_1))
 \le c_1(2+\log(1/c_1))L .
\]
By Lemma~\ref{lem:t-cap-family}(i), with its constant $C'_p$, and Step~1,
\[
 \log(2(1+|\mathcal C|d_{\rm sk})/\delta_0)\le
 C'_pc_1(2+\log(1/c_1))L+Cc_2L+C_p\log(2r).
\]
Choose $c_1$ with $C'_pc_1(2+\log(1/c_1))\le1/4$ and $c_2$ with
$Cc_2\le1/4$. Then the right side is at most $L/2+C_p\log(2r)$. On the
other hand $\epsilon^{-2}\ge e^{2L}/4$, $k^{p-2}\ge r^{-2}$,
$q^{-2p}\ge r^{-4}$ and $C_*\ge4$, so $\log N_*\ge2L-3\log r$. For
$L\ge C_p\log(2r)$ with $C_p$ large, this
gives Eq.~\eqref{lem:x-subexp-eq-budget}. Finally, since
$\epsilon^{-2}r^3K^{p-2}L^{-2p}\ge e^{2L}r^{-3}/4\ge1$,
\[
 N_*\le1+C_*\epsilon^{-2}r^3(c_1K/2)^{p-2}(c_2L/2)^{-2p}
 \le C_p\,\epsilon^{-2}r^3K^{p-2}L^{-2p}.
\]

\emph{Part (b), lower bound.}
\emph{Constants.} Define $C_{\rm k}:=64$ and fix an absolute
$\kappa\in(0,1/2]$ with $2C_{\rm k}\kappa(2+\log(e/\kappa))\le1/2$. Let
$k_1$ and $C_{p,\kappa}$ be given by
Lemma~\ref{lem:x-subexp-covariance}, define
$a_0:=\max\{1,1/c_p,C_{p,\kappa}\}$, which is at least one as
Lemma~\ref{lem:x-subexp-hard-array} requires, with $c_p$ from
Lemma~\ref{lem:x-subexp-walsh}, and let $\theta_{\rm h},\delta_{\rm h},C_{\rm h}$ be given
by Lemma~\ref{lem:x-subexp-hard-array} for this $a_0$. Define
$c_0:=2^{-25}$, so that every $H\ge\log(e/c_0)$ satisfies
$e^H\ge2^{19}eC_{\rm k}$, and let $c_\circ:=c_0/4$. All of these depend only
on $p$, except $C_{\rm k},\kappa,c_0,c_\circ$, which are absolute. The range
$C_p\log r\le L\le c_\circ r$ is empty unless $r$ is large, and we assume
$C_p\ge4$. Then $r\ge3$, and $\log(er)\le2\log r\le2L/C_p$.

\emph{A construction with $M$ groups.} Let $M\ge1$ be an integer with
$ML\le c_0r$ and define $H:=\log(er/(ML))\ge\log(e/c_0)$. Let $k$ be the odd
integer with $C_{\rm k} L/H\le k<C_{\rm k} L/H+2$. Define $N:=\lfloor
r/2\rfloor$, $n:=\lfloor(r-N)/(Mk)\rfloor$ and
\[
 D:=\lfloor(\frac{\delta_{\rm h}}3)^2\epsilon^{-2}k^{-1/2}
 (k^{p/2}+M)^{-2}\rfloor .
\]
We check the hypotheses of the lemmas, for $C_p$ large.
\begin{enumerate}[label=(\roman*),nosep,leftmargin=*]
\item Since $H\le\log(er)\le L$, we have $C_{\rm k} L/H\ge C_{\rm k}$, so
$k\le2C_{\rm k} L/H$, and $k\ge C_{\rm k} L/\log(er)\ge C_{\rm k} C_p/2\ge
\max\{k_1,100\}$.
\item Since $r-N\ge r/2$ and $H\ge1$, we have $(r-N)/(Mk)\ge rH/(4C_{\rm k}
ML)=He^{H-1}/(4C_{\rm k})\ge e^H/(4eC_{\rm k})$. As $e^H\ge2^{19}eC_{\rm k}$
and $2^{16}\ge8eC_{\rm k}$, this gives $n\ge e^H/(8eC_{\rm k})\ge2^{16}$ and
$\log n\ge H-\log(8eC_{\rm k})\ge H/2$. Also $n\le rH/(C_{\rm k} ML)\le e^{2H}$, so
$\log n\le2H$. Since $n\ge1$, we have $Mk\le r$ and $N+Mkn\le r$.
\item Since $\epsilon^{-2}=e^{2L}/4$, $k\le r$ and $k^{p/2}+M\le2r$, we
get $2L\ge\log D\ge2L-3\log(2r)-\log(144/\delta_{\rm h}^2)\ge3L/2$. Hence
$D\ge\sqrt k$ and $ND\ge C_{\rm h}$.
\item By (i) and (ii), $(k/16)\log n\ge(C_{\rm k} L/(16H))(H/2)=2L\ge\log
D$, so $D\le n^{k/16}$.
\item Let $d:=\lfloor\kappa k\rfloor\ge1$. Since
$\sum_{j\le d}\binom kj\le(ek/d)^d\le(e/\kappa)^{\kappa k}$,
\[
 \log\sum_{j\le\kappa k}\binom kjn^j\le\kappa k(\log n+\log(e/\kappa)
 )\le2C_{\rm k}\kappa(2+\log(e/\kappa))L\le L/2\le\log(D/8).
\]
\end{enumerate}
By (ii) and (iv), Lemma~\ref{lem:x-subexp-tuples} gives a set
$\mathcal D\subseteq[n]^k$ of $D$ tuples, any two agreeing in at most
$k/4$ positions. By (i), (iii) and (v),
Lemma~\ref{lem:x-subexp-covariance} applies. With
$\varsigma:=k^{-1/4}$, its conclusions are the hypotheses of
Lemma~\ref{lem:x-subexp-hard-array} for $a_0$, and (iii) gives $ND\ge
C_{\rm h}$. Let $w\in[1/2,3/2]^{N\times\mathcal D}$ be the array given
there, and let $(S,\lambda)$ be a sparsifier of $A_w$ with accuracy
$\epsilon$. By Lemma~\ref{lem:x-subexp-sign-test} and the choice of $D$,
\[
 \|\tilde w^g-w\|_*\le3\epsilon ND(k^{p/2}+M)\le\delta_{\rm h}N\sqrt D\,k^{-1/4}
 \qquad(g\in[M]).
\]
Hence each $\tilde w^g$ has at least $\theta_{\rm h}ND$ nonzero entries, that is,
$S$ contains at least $\theta_{\rm h}ND$ rows of group $g$. The groups use
disjoint rows, so $|S|\ge\theta_{\rm h}MND$. Since $N\ge r/3$ and
$D\ge(\delta_{\rm h}/3)^2\epsilon^{-2}k^{-1/2}(k^{p/2}+M)^{-2}/2$ by (iii),
\begin{equation}\label{lem:x-subexp-eq-groups}
 S_p(r,\epsilon)\ge c'_p\,\frac r{\epsilon^2}\cdot
 \frac M{k^{1/2}(k^{p/2}+M)^2}.
\end{equation}

\emph{The case $L^b\le r$.} Define $H_0:=\log(er/L^b)\ge1$ and
$M:=\lfloor c_0(L/H_0)^{p/2}\rfloor$. Since $H_0\le\log(er)\le2L/C_p$, we
have $c_0(L/H_0)^{p/2}\ge2$ for $C_p$ large, so $M\ge
c_0(L/H_0)^{p/2}/2\ge1$. Also $ML\le c_0L^bH_0^{-p/2}$, which is at most
$c_0r$ and at most $L^b$. Hence $H\ge H_0$, and
\[
 H\le\log(\frac{2er}{c_0L^bH_0^{-p/2}})=H_0+\frac p2\log H_0
 +\log\frac2{c_0}\le C_1H_0,\qquad C_1:=2+\log(2/c_0).
\]
By (i), $k\le2C_{\rm k} L/H_0$ and $k\ge(C_{\rm k}/C_1)L/H_0$. The latter
gives $M\le c_0(L/H_0)^{p/2}\le C_2k^{p/2}$ with
$C_2:=c_0(C_1/C_{\rm k})^{p/2}$. Then Eq.~\eqref{lem:x-subexp-eq-groups} gives
\[
 S_p(r,\epsilon)\ge\frac{c'_p}{(1+C_2)^2}\cdot\frac r{\epsilon^2}\cdot
 \frac{c_0}2(\frac L{H_0})^{p/2}(\frac{2C_{\rm k} L}{H_0}
 )^{-p-1/2}\ge c''_p\,\frac r{\epsilon^2}(\frac{H_0}L
 )^{(p+1)/2}.
\]

\emph{The case $L^b\ge r$.} Define $M:=\lfloor c_0r/(2L)\rfloor$. Since
$L\le c_\circ r=c_0r/4$, we have $c_0r/(2L)\ge2$ and $M\ge c_0r/(4L)$. Also
$ML\le c_0r/2$, so the construction applies, with $H\ge1$. By (i),
$k\le2C_{\rm k} L$. Moreover $M\le c_0r/(2L)\le c_0L^{b-1}=c_0L^{p/2}$, so
$k^{p/2}+M\le((2C_{\rm k})^{p/2}+c_0)L^{p/2}$.
Then Eq.~\eqref{lem:x-subexp-eq-groups} gives
\[
 S_p(r,\epsilon)\ge\frac{c'_p}{((2C_{\rm k})^{p/2}+c_0)^2}\cdot
 \frac r{\epsilon^2}\cdot\frac{c_0r}{4L}(2C_{\rm k} L)^{-1/2}L^{-p}
 \ge c''_p\,\frac r{\epsilon^2}\cdot\frac r{L^{p+3/2}}.\qedhere
\]
\end{proof}

\subsection{The Quadratic Endpoint}

\begin{lemma}[The quadratic endpoint]\label{lem:quadratic_endpoint}
For $p=2$, every $r\ge2$ and every $0<\epsilon\le1/2$,
\[
 S_2(r,\epsilon)=\Theta(\min\{r/\epsilon^2,r^2\}).
\]
\end{lemma}
\begin{proof}
The two upper bounds follow from
Lemma~\ref{lem:bss_spectral_sparsification} and
conic Carathéodory in the space of symmetric $r$ by $r$ matrices.
For the lower bound take the incidence matrix of the complete graph
on $n:=r+1$ vertices. A reweighted Laplacian differs from the original
by Frobenius norm at most $\epsilon n\sqrt r$. Thus, writing $w_{ij}$
for the output weights, including zeros,
$2\sum_{i<j}(w_{ij}-1)^2\le\epsilon^2 n^2r$. Its total weight
$W=\sum_{i<j}w_{ij}$ is $\Theta(nr)$ by the trace comparison.
If $s$ edges survive, $\sum w_{ij}^2\ge W^2/s$.
Expanding the Frobenius estimate gives
$W^2/s\le O(n^2(1+\epsilon^2r))$, which proves the lower bound.
\end{proof}

\subsection{Tools for the Subexponential Lower Bound}

Throughout the lower bound, $0<p<2$ and
\[
 h(s):=\tfrac12(|s+1|^p-|s-1|^p)\qquad(s\in\R).
\]
For a finite set $\mathcal T$, a \emph{uniform sign vector} is a uniform
random point $x$ of $\{-1,1\}^{\mathcal T}$, and
$\chi_T(x):=\prod_{v\in T}x_v$ for $T\subseteq\mathcal T$. The
characters $\chi_T$ are orthonormal in $L_2$ of the uniform measure. We
use the following bound on point masses. If $x$ is a uniform sign vector
in $\{-1,1\}^a$, $a\ge0$, then
\begin{equation}\label{lem:x-subexp-eq-atom}
 \max_{s\in\R}\Pr[\sum_{l=1}^ax_l=s]\le(a+1)^{-1/2}.
\end{equation}
Indeed, for $a=2l$ the maximum is $v_l:=\binom{2l}l4^{-l}$, and
$v_{l+1}^2(2l+3)=v_l^2(2l+1)\cdot(2l+1)(2l+3)/(2l+2)^2\le v_l^2(2l+1)$,
so $v_l^2(2l+1)\le v_0^2=1$. For $a=2l+1$ the maximum is
$\binom{2l+1}l2^{-2l-1}=v_l(2l+1)/(2l+2)\le(2l+2)^{-1/2}$.

\begin{lemma}[Separated tuples]\label{lem:x-subexp-tuples}
Let $n\ge2^{16}$ and $k\ge1$ be integers and $1\le D\le n^{k/16}$. There is
a set $\mathcal D\subseteq[n]^k$ of $D$ tuples such that any two distinct
tuples of $\mathcal D$ agree in at most $k/4$ positions.
\end{lemma}
\begin{proof}
Let $t,t'$ be independent uniform tuples. For a set $T\subseteq[k]$ of
positions, $t$ and $t'$ agree on $T$ with probability $n^{-|T|}$. A union
bound over the at most $2^k$ sets with $|T|>k/4$ shows that they agree in
more than $k/4$ positions with probability at most $2^kn^{-k/4}\le
n^{-3k/16}$, since $2\le n^{1/16}$. For $D$ independent uniform tuples the
expected number of such pairs is at most $D^2n^{-3k/16}\le n^{-k/16}<1$.
Fix a realization with no such pair. Its tuples are distinct, since equal
tuples agree in $k$ positions.
\end{proof}

\begin{lemma}[The high Walsh band of one tuple]\label{lem:x-subexp-walsh}
Let $k\ge100$ be odd, let $x$ be a uniform sign vector in $\{-1,1\}^k$ and
$s:=x_1+\cdots+x_k$. For $T\subseteq[k]$ define
$\hat h(T):=\E[h(s)\chi_T(x)]$. Then $\hat h(T)=0$ when $|T|$ is even, and
\[
 \sum_{T\subseteq[k],\ |T|\ge k-k/100}\hat h(T)^2\ge c_pk^{-1/2},
\]
where $c_p>0$ depends only on $p$.
\end{lemma}
\begin{proof}
For $0<p<2$ and $y\in\R$ the substitution $u=t|y|$ gives
\begin{equation}\label{lem:x-subexp-eq-cosine}
 |y|^p=\varpi_p\int_0^\infty(1-\cos(ty))t^{-1-p}\,\d t,
 \qquad\varpi_p^{-1}:=\int_0^\infty(1-\cos u)\,u^{-1-p}\,\d u ,
\end{equation}
and $\varpi_p^{-1}\in(0,\infty)$ because $0\le1-\cos u\le\min\{2,u^2/2\}$.
Apply Eq.~\eqref{lem:x-subexp-eq-cosine} to $y=s\pm1$ and use
$\cos(t(s-1))-\cos(t(s+1))=2\sin(t)\sin(ts)$. This gives
\[
 h(s)=\varpi_p\int_0^\infty\sin(t)\sin(ts)\,t^{-1-p}\,\d t,
\]
where the integral converges absolutely since
$|\sin(t)\sin(ts)|\le\min\{1,|s|t^2\}$. By independence,
$\E[\chi_T(x)e^{\mathrm its}]=(\mathrm i\sin t)^{|T|}(\cos t)^{k-|T|}$. Its
imaginary part $\E[\chi_T(x)\sin(ts)]$ vanishes when $|T|$ is even and
equals $\pm\sin^{|T|}(t)\cos^{k-|T|}(t)$ when $|T|$ is odd. Exchanging the
finite expectation with the integral, we get $\hat h(T)=0$ for even $|T|$,
and for odd $|T|$
\[
 |\hat h(T)|=\varpi_p|\int_0^\infty\sin^{|T|+1}(t)\cos^{k-|T|}(t)
 \,t^{-1-p}\,\d t| .
\]
Now let $|T|=k-j$ with $j$ even, so that $|T|$ is odd. Both exponents of
$f_j(t):=\sin^{k+1-j}(t)\cos^j(t)$ are even, so $f_j$ is nonnegative and
$\pi$-periodic, and $f_j(\pi-t)=f_j(t)$. Keeping only the interval
$[\pi,2\pi]$ of the integral, we get
\[
 |\hat h(T)|\ge\varpi_p(2\pi)^{-1-p}\int_\pi^{2\pi}f_j
 =2\varpi_p(2\pi)^{-1-p}I_j,\qquad
 I_j:=\int_0^{\pi/2}f_j=\frac{\Gamma(\alpha)\Gamma(\beta)}{2\Gamma(\alpha+\beta)},
\]
with $\alpha:=(k+2-j)/2$ and $\beta:=(j+1)/2$, by the Beta integral.

We bound $\binom kjI_j^2$ from below. Define $\gamma:=\alpha+\beta=(k+3)/2$
and $\rho(y):=\Gamma(y+\frac12)/\Gamma(y)$ for $y>0$. The duplication
formula $\Gamma(2y)=2^{2y-1}\pi^{-1/2}\Gamma(y)\Gamma(y+\frac12)$ applied
to $j!=\Gamma(2\beta)$, $(k-j)!=\Gamma(2\alpha)/(2\alpha-1)$ and
$k!=\Gamma(2\gamma)/((2\gamma-1)(2\gamma-2))$ gives
\[
 \binom kjI_j^2=\frac{\sqrt\pi}2\cdot\frac{k+1-j}{(k+1)(k+2)}\cdot
 \frac{\rho(\gamma)}{\rho(\alpha)\rho(\beta)} .
\]
Log-convexity of $\Gamma$ gives $\Gamma(y+\frac12)^2\le\Gamma(y)\Gamma(y+1)$
and $\Gamma(y+1)^2\le\Gamma(y+\frac12)\Gamma(y+\frac32)$, hence
$y/\sqrt{y+1/2}\le\rho(y)\le\sqrt y$. Since $\gamma\ge\alpha$ and
$\gamma\ge1/2$, we get $\rho(\gamma)/\rho(\alpha)\ge\sqrt{\gamma/2}/\sqrt\alpha
\ge1/\sqrt2$, and $\rho(\beta)\le\sqrt{(j+1)/2}$. For $k\ge100$ and $j\le
k/100$ we have $(k+1-j)/((k+1)(k+2))\ge0.96/k$, and therefore
$\binom kjI_j^2\ge1/(2k\sqrt{j+1})$.

There are $\binom kj$ sets $T$ with $|T|=k-j$. The even integers
$j\in[0,k/100]$ number at least $k/200$, and each has $j+1\le k/50$. Hence
\[
 \sum_{|T|\ge k-k/100}\hat h(T)^2\ge4\varpi_p^2(2\pi)^{-2-2p}
 \cdot\frac k{200}\cdot\frac{\sqrt{50}}{2k\sqrt k}
 \ge\frac{\varpi_p^2(2\pi)^{-2-2p}}{15\sqrt k}.\qedhere
\]
\end{proof}

\begin{lemma}[Localized Chebyshev truncation]\label{lem:x-subexp-localized}
For every integer $d\ge1$ there is an even polynomial $Q_d$ of degree at
most $d$ such that
\[
 ||t|^p-Q_d(t)|\le\frac{C_pd^{-p}}{1+d|t|}\qquad(|t|\le1).
\]
\end{lemma}
\begin{proof}
The function $|\cos\theta|^p$ is even, continuous and $\pi$-periodic. Let
$A_0:=\frac1\pi\int_{-\pi/2}^{\pi/2}\cos^p\theta\,\d\theta$ and
$A_j:=\frac2\pi\int_{-\pi/2}^{\pi/2}\cos^p\theta\cos(2j\theta)\,\d\theta$
for $j\ge1$ be its cosine coefficients. For $j\ge0$ the function
$\cos^{p+1}\theta\sin((2j+1)\theta)$ vanishes at $\pm\pi/2$, so the integral
of its derivative over $[-\pi/2,\pi/2]$ vanishes:
\[
 \int_{-\pi/2}^{\pi/2}((2j+1)\cos^{p+1}\theta\cos((2j+1)\theta)
 -(p+1)\cos^p\theta\sin\theta\sin((2j+1)\theta))\d\theta=0 .
\]
Insert $2\cos\theta\cos((2j+1)\theta)=\cos(2j\theta)+\cos((2j+2)\theta)$
and $2\sin\theta\sin((2j+1)\theta)=\cos(2j\theta)-\cos((2j+2)\theta)$. For
$j=0$ this gives $(p+2)\int\cos^p\theta\cos(2\theta)\,\d\theta
=p\int\cos^p\theta\,\d\theta$, so $0<A_1\le2$. For $j\ge1$ it gives the
recurrence
\[
 (p+2j+2)A_{j+1}=(p-2j)A_j .
\]
Since $p<2\le2j$, the ratio $A_{j+1}/A_j=-(2j-p)/(2j+2+p)$ is negative
with absolute value less than one, for every $0<p<2$ and $j\ge1$. Hence
$A_j=(-1)^{j+1}|A_j|$ for $j\ge1$, and $|A_j|$ decreases in $j\ge1$.
Moreover $1-(2j-p)/(2j+2+p)=(2p+2)/(2j+2+p)\ge(p+1)/(j+2)$, so
\[
 |A_j|\le A_1\exp(-(p+1)\sum_{l=1}^{j-1}\frac1{l+2})
 \le2(\frac3{j+2})^{p+1}\le54\,j^{-p-1}.
\]
In particular $\sum_j|A_j|<\infty$. The cosine series therefore converges
uniformly to a continuous function with the same coefficients, which is
$|\cos\theta|^p$. With the Chebyshev polynomials,
$T_m(\cos\theta)=\cos(m\theta)$, this reads
$|t|^p=\sum_{j\ge0}A_jT_{2j}(t)$ uniformly on $[-1,1]$.

Let $J:=\lfloor d/2\rfloor+1$, so that $d/2<J\le d$, and
$Q_d:=\sum_{j<J}A_jT_{2j}$, an even polynomial of degree $2J-2\le d$.
For $t=\cos\theta$ the error is
$\mathcal E(t):=|t|^p-Q_d(t)=\sum_{j\ge J}A_j\cos(2j\theta)$. First,
$|\mathcal E(t)|\le54\sum_{j\ge J}j^{-p-1}\le54(1+1/p)J^{-p}$. Second, define
$\phi:=\theta-\pi/2$, so that $(-1)^j\cos(2j\theta)=\cos(2j\phi)$ and
$|\sin\phi|=|t|$. Then $\mathcal E(t)=-\sum_{j\ge J}|A_j|\cos(2j\phi)$. For
$t\ne0$ the partial sums $S_m:=\sum_{j=J}^m\cos(2j\phi)$ satisfy
$|S_m|\le|\sin((m-J+1)\phi)/\sin\phi|\le1/|t|$. Abel summation,
\[
 \sum_{j=J}^m|A_j|\cos(2j\phi)=\sum_{j=J}^{m-1}(|A_j|-|A_{j+1}|)S_j
 +|A_m|S_m ,
\]
and the monotonicity of $|A_j|$ give $|\mathcal E(t)|\le|A_J|/|t|\le
54J^{-p-1}/|t|$. Since $J>d/2$, the two bounds give
$|\mathcal E(t)|\le C_pd^{-p}\min\{1,1/(d|t|)\}\le2C_pd^{-p}/(1+d|t|)$.
\end{proof}

\begin{lemma}[A large subspace with controlled covariance]
\label{lem:x-subexp-covariance}
Let $0<\kappa\le1/2$. There are $k_1$ and $C_{p,\kappa}$, depending only
on $p$ and $\kappa$, with the following property. Let $k\ge k_1$ be odd,
let $n\ge1$, and let $\mathcal D\subseteq[n]^k$ be a set of $D\ge\sqrt k$
tuples, any two of which agree in at most $k/4$ positions. Suppose
$\sum_{j\le\kappa k}\binom kjn^j\le D/8$. Let $x$ be a uniform sign vector
in $\{-1,1\}^{[k]\times[n]}$, and define
\[
 s_t(x):=\sum_{l=1}^kx_{l,t_l},\qquad
 \mathbf h(x):=(h(s_t(x)))_{t\in\mathcal D},\qquad
 \Sigma:=\E[\mathbf h(x)\mathbf h(x)^\top] .
\]
Then $\Sigma\succeq c_pk^{-1/2}I_{\mathcal D}$, with $c_p$ from
Lemma~\ref{lem:x-subexp-walsh}, and there is an orthogonal projection
$\Pi$ of $\R^{\mathcal D}$ of rank at least $D/2$ such that
\[
 \|\Pi\Sigma\Pi\|\le C_{p,\kappa}k^{-1/2},\qquad
 \E[\|\Pi\mathbf h(x)\|_2^4]\le C_{p,\kappa}D^2/k .
\]
\end{lemma}
\begin{proof}
We take $k_1:=\max\{100,\lceil2/\kappa\rceil\}$. Write
$\mathcal I_t:=\{(l,t_l):l\in[k]\}$ for the coordinates used by $t$. Two
distinct tuples share at most $k/4$ coordinates, since coordinates in
different blocks are distinct.

\emph{The lower bound.} The function $h(s_t(x))$ depends only on the
coordinates in $\mathcal I_t$. Its Walsh coefficient at $T\subseteq[k]\times[n]$
vanishes unless $T\subseteq\mathcal I_t$, and for $T\subseteq\mathcal I_t$ it equals
$\hat h(\{l:(l,t_l)\in T\})$ in the notation of
Lemma~\ref{lem:x-subexp-walsh}. Call $T$ \emph{high} if $|T|\ge k-k/100$.
A high $T$ is contained in $\mathcal I_t$ for at most one $t\in\mathcal D$,
because $k-k/100>k/4$. For $a\in\R^{\mathcal D}$, Parseval's identity (Lemma~\ref{lem:pre_parseval}) and
Lemma~\ref{lem:x-subexp-walsh} give
\[
 a^\top\Sigma a=\sum_T(\sum_ta_t\E[h(s_t)\chi_T])^2
 \ge\sum_{t\in\mathcal D}a_t^2\sum_{T\subseteq\mathcal I_t\text{ high}}
 \E[h(s_t)\chi_T]^2\ge c_pk^{-1/2}\|a\|_2^2 .
\]

\emph{A low-degree part.} Let $d:=\lfloor\kappa k\rfloor$. Since
$\kappa k\ge2$, we have $d\ge\kappa k/2\ge\kappa(k+1)/4$. Let $Q_d$ be the
polynomial of Lemma~\ref{lem:x-subexp-localized} and
\[
 \tilde h(s):=\frac{(k+1)^p}2(Q_d(\frac{s+1}{k+1})
 -Q_d(\frac{s-1}{k+1})),
\]
a polynomial of degree at most $d$. Since
$h(s)=\frac12(k+1)^p(|\frac{s+1}{k+1}|^p-|\frac{s-1}{k+1}|^p)$, for odd $s$
with $|s|\le k$ Lemma~\ref{lem:x-subexp-localized} gives
\[
 |h(s)-\tilde h(s)|\le\frac{(k+1)^pC_pd^{-p}}2\sum_\pm\frac1{1+d|s\pm1|/(k+1)}
 \le\frac{C'_{p,\kappa}}{1+|s|}.
\]
Here we used $(k+1)/d\le4/\kappa$ and, since $|s|\ge1$,
$1+d|s\pm1|/(k+1)\ge(\kappa/4)(1+|s\pm1|)\ge(\kappa/8)(1+|s|)$. Every
$s_t(x)$ is odd with $|s_t(x)|\le k$. Expanding powers of $s_t(x)$ and
using $x_v^2=1$, we can write $\tilde h(s_t(x))=\sum_S\gamma_{t,S}
\chi_S(x)$, where $\gamma_{t,S}=0$ unless $S\subseteq\mathcal I_t$ and $|S|\le
d$. Let $\mathcal S$ be the family of sets of at most $d$ coordinates with
at most one coordinate in each block, so that $|\mathcal S|=\sum_{j\le
d}\binom kjn^j\le D/8$. Then $(\tilde h(s_t(x)))_{t\in\mathcal D}=
\sum_{S\in\mathcal S}\chi_S(x)(\gamma_{t,S})_{t\in\mathcal D}$ lies, for
every $x$, in the span $V$ of the vectors $(\gamma_{t,S})_{t\in\mathcal
D}$, $S\in\mathcal S$, and $\dim V\le D/8$.

\emph{The residual.} Define $U_t:=h(s_t)-\tilde h(s_t)$ and
$U:=(U_t)_{t\in\mathcal D}$, so that $|U_t|\le C'_{p,\kappa}/(1+|s_t|)$.
By Eq.~\eqref{lem:x-subexp-eq-atom} and $\sum_{s\in\mathbb Z}(1+|s|)^{-2}\le3$,
\[
 \E[U_t^2]\le3C'^2_{p,\kappa}(k+1)^{-1/2},\qquad
 \E[U_t^4]\le3C'^4_{p,\kappa}(k+1)^{-1/2}.
\]
For distinct $t,t'\in\mathcal D$, condition on the coordinates in
$W:=\mathcal I_t\cap\mathcal I_{t'}$, of which there are at most $k/4$. Given them,
$s_t$ and $s_{t'}$ are shifts of sums over the disjoint sets
$\mathcal I_t\setminus W$ and $\mathcal I_{t'}\setminus W$, each of at least $3k/4$
independent signs. The conditional point masses obey the same bound with
$3k/4$ in place of $k$, so
$\E[U_t^2U_{t'}^2]\le9C'^4_{p,\kappa}(3k/4)^{-1}=12C'^4_{p,\kappa}/k$.
With $D\ge\sqrt k$ we get
\[
 \E[\|U\|_2^2]\le3C'^2_{p,\kappa}Dk^{-1/2},\qquad
 \E[\|U\|_2^4]=\sum_t\E[U_t^4]+\sum_{t\ne t'}\E[U_t^2U_{t'}^2]
 \le15C'^4_{p,\kappa}D^2/k .
\]

\emph{The projection.} Let $\Pi_0$ be the orthogonal projection onto
$V^\perp$. Then $\Pi_0\mathbf h=\Pi_0U$, and
$\Sigma_0:=\Pi_0\Sigma\Pi_0=\E[\Pi_0UU^\top\Pi_0]$ has trace at most
$\E[\|U\|_2^2]\le3C'^2_{p,\kappa}Dk^{-1/2}$. Let $E_+$ be the span of the
eigenvectors of $\Sigma_0$ with eigenvalues greater than
$24C'^2_{p,\kappa}k^{-1/2}$. Then $\dim E_+\le D/8$, and $E_+\subseteq
V^\perp$ since it lies in the range of $\Sigma_0$. Let $\Pi$ be the
orthogonal projection onto $V^\perp\cap E_+^\perp$, of rank at least
$D-D/8-D/8\ge D/2$. This subspace is invariant under $\Sigma_0$, with
eigenvalues at most $24C'^2_{p,\kappa}k^{-1/2}$. Since $\Pi=\Pi\Pi_0$, we
have $\Pi\mathbf h=\Pi U$ and $\Pi\Sigma\Pi=\Pi\Sigma_0\Pi$. Hence
$\|\Pi\Sigma\Pi\|\le24C'^2_{p,\kappa}k^{-1/2}$ and
$\E[\|\Pi\mathbf h\|_2^4]\le\E[\|U\|_2^4]\le15C'^4_{p,\kappa}D^2/k$.
\end{proof}

\begin{lemma}[Gaussian concentration]\label{lem:x-subexp-gauss-conc}
Let $Z$ be a standard Gaussian vector in $\R^m$ and let $f:\R^m\to\R$
satisfy $|f(z)-f(z')|\le\ell\|z-z'\|_2$. Then for every $u>0$,
\[
 \Pr[f(Z)-\E[f(Z)]>u]\le\exp(-\frac{2u^2}{\pi^2\ell^2}),
\]
and the same bound holds for $\Pr[f(Z)-\E[f(Z)]<-u]$.
\end{lemma}
\begin{proof}
The second claim follows from the first applied to $-f$. Since
$|f(z)|\le|f(0)|+\ell\|z\|_2$, all exponential moments of $f(Z)$ are
finite. Suppose first that $f$ is continuously differentiable, so that
$\|\nabla f\|_2\le\ell$. Let $Z'$ be an independent copy of $Z$, and for
$0\le\vartheta\le\pi/2$ define $Z_\vartheta:=Z\sin\vartheta+Z'\cos\vartheta$
and $\dot Z_\vartheta:=Z\cos\vartheta-Z'\sin\vartheta$. Then
\[
 f(Z)-f(Z')=\int_0^{\pi/2}\langle\nabla f(Z_\vartheta),\dot Z_\vartheta
 \rangle\,\d\vartheta .
\]
For $\lambda>0$, Jensen's inequality for the uniform probability on
$[0,\pi/2]$ gives
\[
 e^{\lambda(f(Z)-f(Z'))}\le\frac2\pi\int_0^{\pi/2}\exp(\frac{\pi\lambda}2
 \langle\nabla f(Z_\vartheta),\dot Z_\vartheta\rangle)\d\vartheta .
\]
For each $\vartheta$ the pair $(Z_\vartheta,\dot Z_\vartheta)$ is a
rotation of $(Z,Z')$, so it has the law of $(Z,Z')$. Conditioning on
$Z_\vartheta$ gives
\[
 \E[\exp(\frac{\pi\lambda}2\langle\nabla f(Z_\vartheta),\dot
 Z_\vartheta\rangle)]=\E[\exp(\frac{\pi^2\lambda^2}8\|\nabla
 f(Z)\|_2^2)]\le e^{\pi^2\lambda^2\ell^2/8}.
\]
Jensen's inequality in $Z'$ gives
$\E[e^{\lambda(f(Z)-\E[f(Z)])}]\le\E[e^{\lambda(f(Z)-f(Z'))}]\le
e^{\pi^2\lambda^2\ell^2/8}$. Markov's inequality (Lemma~\ref{lem:pre_markov}) with
$\lambda:=4u/(\pi^2\ell^2)$ gives the claim. For general $f$ and
$s>0$, the function $f_s(z):=\E[f(z+s Z')]$ is smooth and
$\ell$-Lipschitz, and $|f_s-f|\le\ell s\,\E[\|Z'\|_2]\le\ell s\sqrt
m$. Hence $\Pr[f(Z)-\E[f(Z)]>u]\le\exp(-2(u-2\ell s\sqrt m)^2/(\pi^2
\ell^2))$, and we let $s\to0$.
\end{proof}

\begin{lemma}[A mass array far from all sparse arrays]
\label{lem:x-subexp-hard-array}
Let $a_0\ge1$. There are $\theta_{\rm h},\delta_{\rm h},C_{\rm h}>0$, depending only on
$a_0$, with the following property. Let $N,D\ge1$ be integers with $ND\ge
C_{\rm h}$, let $\mathbf h$ be a random vector in $\R^D$ on a finite
probability space, $\Sigma:=\E[\mathbf h\mathbf h^\top]$ and
$\varsigma>0$. Suppose that $\Sigma\succeq a_0^{-1}\varsigma^2I_D$, and that
some orthogonal projection $\Pi$ of rank at least $D/2$ satisfies
$\|\Pi\Sigma\Pi\|\le a_0\varsigma^2$ and $\E[\|\Pi\mathbf h\|_2^4]\le
a_0D^2\varsigma^4$. For arrays $a\in\R^{N\times D}$ with rows $a_i$ define
\[
 \|a\|_*:=\sum_{i=1}^N\E[|\langle a_i,\mathbf h\rangle|] .
\]
Then $\|\cdot\|_*$ is a norm, and there is $w\in[1/2,3/2]^{N\times D}$
such that $\|y-w\|_*>\delta_{\rm h} N\sqrt D\varsigma$ for every
$y\in\R^{N\times D}$ with fewer than $\theta_{\rm h} ND$ nonzero entries.
\end{lemma}
\begin{proof}
If $\|a\|_*=0$, then $\langle a_i,\mathbf h\rangle=0$ almost surely, so
$a_i^\top\Sigma a_i=0$ and $a_i=0$ for every $i$. Thus $\|\cdot\|_*$ is a
norm.

\emph{Gaussian arrays.} Let $Z=(Z_i)_{i\le N}$ consist of independent
standard Gaussian vectors in $\R^D$, define $X_i:=\Pi Z_i$, and let
$f(Z):=\|X\|_*=\sum_i\E_{\mathbf h}[|\langle Z_i,\Pi\mathbf h\rangle|]$.
Write $Y:=\|\Pi\mathbf h\|_2$. For fixed $\mathbf h$ the variable
$\langle Z_i,\Pi\mathbf h\rangle$ is centered Gaussian with standard
deviation $Y$, so $\E[f(Z)]=\sqrt{2/\pi}\,N\,\E[Y]$. We have
$\E[Y^2]=\tr[\Pi\Sigma\Pi]$, which lies between $a_0^{-1}\varsigma^2D/2$
and $a_0\varsigma^2D$, since $\Pi\Sigma\Pi\succeq a_0^{-1}\varsigma^2\Pi$.
H\"older's inequality $\E[Y^2]\le(\E[Y])^{2/3}(\E[Y^4])^{1/3}$ gives
$\E[Y]\ge(\E[Y^2])^{3/2}(\E[Y^4])^{-1/2}$. Hence
\[
 c_4N\sqrt D\varsigma\le\E[f(Z)]\le\sqrt{a_0}\,N\sqrt D\varsigma,
 \qquad c_4:=\sqrt{2/\pi}\,(2a_0)^{-3/2}a_0^{-1/2}.
\]
The function $f$ is Lipschitz with constant $\ell:=\sqrt{a_0N}\varsigma$,
since
\[
 |f(Z)-f(Z')|\le\sum_i(\E[\langle Z_i-Z'_i,\Pi\mathbf h\rangle^2])^{1/2}
 \le\|\Pi\Sigma\Pi\|^{1/2}\sum_i\|Z_i-Z'_i\|_2
 \le\sqrt{a_0N}\,\varsigma\,\|Z-Z'\|_2 .
\]
Let $c_5:=2(1-2^{-1/2})^2c_4^2/(\pi^2a_0)$, so that $c_5\le2/\pi^2$.
Lemma~\ref{lem:x-subexp-gauss-conc} with
$u=\sqrt{a_0}N\sqrt D\varsigma$ and with $u=(1-2^{-1/2})c_4N\sqrt
D\varsigma$ gives
\[
 \Pr[f(Z)>2\sqrt{a_0}N\sqrt D\varsigma]\le e^{-c_5ND},\qquad
 \Pr[f(Z)<2^{-1/2}c_4N\sqrt D\varsigma]\le e^{-c_5ND}.
\]
If $Z'$ is an independent copy of $Z$ and $X'_i:=\Pi Z'_i$, then $Z-Z'$
has the law of $\sqrt2Z$ and $f$ is positively homogeneous, so the second
bound gives $\Pr[\|X-X'\|_*<c_4N\sqrt D\varsigma]\le e^{-c_5ND}$.

\emph{Bounded entries and sampling.} The array $X$ is a standard Gaussian
vector in the subspace $(\operatorname{ran}\Pi)^N$ of $\R^{N\times D}$, and
each entry $X_{it}$ has variance $\Pi_{tt}\le1$. Fix $B\ge1$ with
$e^{-B^2/2}\le c_5/8$. For a standard Gaussian $g$ we have
$\Pr[|g|>B]\le e^{-B^2/2}\le1/2$. Lemma~\ref{lem:t-sidak-chaining}(d),
applied to the $ND$ strips $\{|X_{it}|\le B\}$, gives
$\Pr[\max_{i,t}|X_{it}|\le B]\ge\exp(-2NDe^{-B^2/2})\ge e^{-c_5ND/4}$. Let
$\Omega_1$ be the event that $\max_{i,t}|X_{it}|\le B$ and
$\|X\|_*\le2\sqrt{a_0}N\sqrt D\varsigma$. For $ND\ge2/c_5$,
\[
 \Pr[\Omega_1]\ge e^{-c_5ND/4}-e^{-c_5ND}\ge\tfrac12e^{-c_5ND/4}.
\]
Let $Q:=\lceil e^{c_5ND/8}\rceil$ and let $X^{(1)},\ldots,X^{(Q)}$ be
independent, each with the law of $X$ conditioned on $\Omega_1$. For
$l\ne l'$,
\[
 \Pr[\|X^{(l)}-X^{(l')}\|_*<c_4N\sqrt D\varsigma]
 \le\frac{\Pr[\|X-X'\|_*<c_4N\sqrt D\varsigma]}{\Pr[\Omega_1]^2}
 \le4e^{-c_5ND/2}.
\]
There are at most $Q^2/2\le2e^{c_5ND/4}$ pairs, so for $ND\ge12/c_5$ some
realization has all pairwise distances at least $c_4N\sqrt D\varsigma$.
Define $w^{(l)}:=\mathbf1+X^{(l)}/(2B)\in[1/2,3/2]^{N\times D}$, where
$\mathbf1$ is the all-ones array, and $\Delta:=c_4N\sqrt D\varsigma/(2B)$.
Then $\|w^{(l)}-w^{(l')}\|_*\ge\Delta$ for $l\ne l'$, and
$\|w^{(l)}-\mathbf1\|_*\le\sqrt{a_0}N\sqrt D\varsigma/B=(2\sqrt{a_0}/c_4)
\Delta$.

\emph{Counting sparse arrays.} Let $\delta_{\rm h}:=c_4/(8B)$, so that
$\delta_{\rm h}N\sqrt D\varsigma=\Delta/4$, let $C_3:=3+16\sqrt{a_0}/c_4$, and fix
$\theta_{\rm h}\in(0,1/2]$ with $\theta_{\rm h}\log(eC_3/\theta_{\rm h})\le c_5/16$. Suppose that
every $w^{(l)}$ has an array $y^{(l)}$ with fewer than $\theta_{\rm h} ND$ nonzero
entries and $\|y^{(l)}-w^{(l)}\|_*\le\Delta/4$. For $l\ne l'$ we get
$\|y^{(l)}-y^{(l')}\|_*\ge\Delta/2$, so the $y^{(l)}$ are distinct, and
$\|y^{(l)}-y^{(l')}\|_*\le2(2\sqrt{a_0}/c_4+1/4)\Delta$. Fix a set $S$ of
entries and consider the $y^{(l)}$ supported in $S$. In the coordinate
space $\R^S$, on which $\|\cdot\|_*$ is a norm, the balls of radius
$\Delta/4$ about them are disjoint and lie in the ball of radius
$(4\sqrt{a_0}/c_4+3/4)\Delta$ about any one of them. Comparing volumes,
there are at most $C_3^{|S|}$ of them. Summing over the supports,
\[
 Q\le\sum_{j<\theta_{\rm h} ND}\binom{ND}jC_3^j .
\]
If $\theta_{\rm h} ND<1$ the right side is $1<Q$. Otherwise the bound
$\sum_{j\le m'}\binom{ND}j\le(eND/m')^{m'}$ for $1\le m'\le ND$, with
$m':=\lfloor\theta_{\rm h} ND\rfloor$, bounds the right side by
$\exp(\theta_{\rm h} ND\log(eC_3/\theta_{\rm h}))\le e^{c_5ND/16}<Q$. Both cases are
contradictions, so some $w:=w^{(l)}$ has the stated property. We may take
$C_{\rm h}:=12/c_5$.
\end{proof}

\subsection{Instances for the Subexponential Lower Bound}

Let $N,M,k,n\ge1$ be integers with $k$ odd and $N+Mkn\le r$, and let
$\mathcal D\subseteq[n]^k$ be a set of $D$ tuples. Let $e_1,\ldots,e_r$ be
the standard basis of $\R^r$. The first block consists of
$e_1,\ldots,e_N$. The next $Mkn$ basis vectors are relabeled
$e_{g,l,\iota}$, $g\in[M]$, $l\in[k]$, $\iota\in[n]$, where $\iota$ runs
over block $l$ of the transverse group $g$. The remaining basis vectors are
padding. For $g\in[M]$ and $t\in\mathcal D$ define
$u_{g,t}:=\sum_{l=1}^ke_{g,l,t_l}$. For a mass array
$w\in(0,\infty)^{N\times\mathcal D}$ let $A_w$ be the matrix with the rows
\[
 w_{it}^{1/p}(e_i+u_{g,t}),\quad(g,i,t)\in[M]\times[N]\times\mathcal D,
 \qquad\text{and}\qquad r^{-1/p}e_j,\quad j\in[r].
\]
The second kind are the identity rows, one for every basis vector,
including those of the first block and the padding. Hence $A_w$ has rank
$r$. Let
$\mathbf h$ be as in Lemma~\ref{lem:x-subexp-covariance} for this
$\mathcal D$, and let $\|\cdot\|_*$ be the norm of
Lemma~\ref{lem:x-subexp-hard-array} on $\R^{N\times\mathcal D}$ for this
$\mathbf h$.

\begin{lemma}[The sign test]\label{lem:x-subexp-sign-test}
Let $w\in[1/2,3/2]^{N\times\mathcal D}$ and let $(S,\lambda)$ be a
sparsifier of $A_w$ with accuracy $\epsilon$. We write its weights as
$\lambda$ because $w$ denotes the mass array. For $g\in[M]$ let
$\tilde w^g_{it}:=\lambda_{g,i,t}w_{it}$ if the row indexed by $(g,i,t)$
lies in $S$, and $\tilde w^g_{it}:=0$ otherwise. Then
$\|\tilde w^g-w\|_*\le3\epsilon ND(k^{p/2}+M)$ for every $g\in[M]$.
\end{lemma}
\begin{proof}
Fix $g$ and define $a:=\tilde w^g-w$. For a sign vector
$x\in\{-1,1\}^{[k]\times[n]}$ and $\sigma\in\{-1,1\}^N$ let
$y=y(\sigma,x)\in\R^r$ have $y_i:=\sigma_i$ for $i\in[N]$,
$y_{g,l,\iota}:=x_{l,\iota}$, and all other entries zero. Then
$\langle e_i+u_{g,t},y\rangle=\sigma_i+s_t(x)$, while
$\langle e_i+u_{g',t},y\rangle=\sigma_i$ for $g'\ne g$, and
$\sum_j|y_j|^p=N+kn\le r$. Let $\mathrm{err}(y)$ be the output energy of
the sparsifier at $y$ minus $\|A_wy\|_p^p$. The contributions of the
groups $g'\ne g$ and of the identity rows to $\mathrm{err}(y)$ do not
change when $\sigma$ is replaced by $-\sigma$, since they depend on
$\sigma$ only through $|\sigma_i|^p=1$. Since
$|\sigma_i+s|^p-|-\sigma_i+s|^p=2\sigma_ih(s)$ for $\sigma_i=\pm1$,
\[
 \tfrac12(\mathrm{err}(y(\sigma,x))-\mathrm{err}(y(-\sigma,x)))
 =\sum_{i=1}^N\sigma_i\langle a_i,\mathbf h(x)\rangle .
\]
Take $\sigma_i:=\operatorname{sgn}\langle a_i,\mathbf h(x)\rangle$, with
$\operatorname{sgn}0:=1$. The accuracy of the sparsifier at both queries
gives
\[
 \sum_{i=1}^N|\langle a_i,\mathbf h(x)\rangle|\le\frac\epsilon2(
 \|A_wy(\sigma,x)\|_p^p+\|A_wy(-\sigma,x)\|_p^p),
\]
and the right side equals
\[
 \epsilon(\sum_{i,t}w_{it}\frac{|s_t(x)+1|^p+|s_t(x)-1|^p}2
 +(M-1)\sum_{i,t}w_{it}+\frac{N+kn}r).
\]
Since $p\le2$, Jensen's inequality gives $\E[|s_t(x)\pm1|^p]\le
(k+1)^{p/2}\le2k^{p/2}$. With $w_{it}\le3/2$ and $ND\ge1$, taking
expectations over a uniform $x$ gives
$\|a\|_*\le\epsilon(3NDk^{p/2}+\frac32(M-1)ND+1)\le3\epsilon ND(k^{p/2}+M)$.
\end{proof}

\subsection{Improved Upper Bound at the Subexponential Crossover}

For fixed $0<\alpha<1$ and $L=r^\alpha$, the general bounds of
Theorem~\ref{thm:below_two_main}(b) and Lemma~\ref{lem:joint-positive-caps} below place $\epsilon^2S_p$ between $r^{D_p(\alpha)}$ and
$r^{U_p(\alpha)}$ up to logarithmic factors in $r$, where
\[
\begin{aligned}
 D_p(\alpha)
   &=\begin{cases}
       1-(p+1)\alpha/2,&\alpha\le 2/(p+2),\\
       2-(p+3/2)\alpha,&\alpha\ge 2/(p+2),
     \end{cases}\\
 U_p(\alpha)
   &=\begin{cases}
       \min\{1,\,3-p-2\alpha\},&0<p\le1,\\
       \min\{1,\,p+1-2p\alpha\},&1<p<2.
     \end{cases}
\end{aligned}
\]
The lower powers follow from Theorem~\ref{thm:below_two_main}(b)
and the upper powers from Lemma~\ref{lem:joint-positive-caps}.
The largest gap in these powers for fixed $p$ occurs where
$U_p$ first drops below one: at $\alpha=1-p/2$ for $p\le1$
and at $\alpha=1/2$ for $p>1$.

\begin{lemma}[Improved upper bound at the subexponential crossover]
\label{lem:joint-positive-caps}
Fix $0<p<2$ and suppose $C_p\log(2r)\le L\le r$.
\begin{enumerate}[label=(\alph*),nosep,leftmargin=*]
\item We have
\[
 S_p(r,\epsilon)\lesssim_p \epsilon^{-2}\log(2r)
 (\frac{A_p(r,K)}{L^{2p}}+\frac{r^3K^{p-2}}{L^{4p}}),
 \qquad
 A_p(r,K)=
 \begin{cases}
 r^{3-p/2}K^{3p/2-2}, & 0<p<1,\\
 r^{2+p/2}K^{p/2-1}, & 1\le p<2.
 \end{cases}
\]
\item For $L\ge C_p\log^3(2r)$, set $Q:=L/\Lambda^3$. Then
\[
 S_p(r,\epsilon)\lesssim_p \epsilon^{-2}\Lambda^{1+4\max\{1,p\}}
 (\frac{B_p(r,K)}{Q^{2p}}+\frac{r^3K^{p-2}}{Q^{4p}}),
 \quad
 B_p(r,K)=
 \begin{cases}
 r^{3-p}K^{2p-2}, & 0<p<1,\\
 r^{p+1}, & 1\le p<2.
 \end{cases}
\]
For $L/r^{2/(p+2)}$ in a fixed compact interval $I\subset(0,\infty)$,
this gives
\[
 S_p(r,\epsilon)\lesssim_{p,I}\epsilon^{-2}r^{1-\delta_p}\log^{O_p(1)}(2r),
 \qquad\delta_p=\frac{p\min\{p,2-p\}}{p+2}>0.
\]
\item For $1\le p<2$ define
$T:=L/\sqrt r$ if $p>1$ and $T:=K/\sqrt r$ if $p=1$. Then
\[
 S_p(r,\epsilon)\lesssim_p \epsilon^{-2}r
      \min\{1,T^{-2p}\log(2+T)\}.
\]
In particular $T\to\infty$ suffices for $S_p=o_p(r/\epsilon^2)$.
\item For $0<p<1$ define
$T_0:=K/r^{1-p/2}$. Then
\[
 S_p(r,\epsilon)\lesssim_p\epsilon^{-2}r
       \min\{1,T_0^{-2}\log(2+T_0)\}.
\]
For $L/r^{2/(p+2)}$ in a fixed compact interval $I\subset(0,\infty)$ this is
$O_{p,I}(\epsilon^{-2}r^{1-p^2/(p+2)}\log^3(2r))$.
\end{enumerate}
\end{lemma}

All four parts follow one scheme. A frame from the ensemble of
Lemma~\ref{lem:t-frame-ensemble} gives each low-energy query coordinates
whose rows and columns have bounded energy, and rounding the largest
coordinates gives a center with a small residual. For each center the row
energies are approximated by polynomials in nonnegative features of the
residual. In part~(a) the features are the positive and negative
contributions of the residual to each row. In part~(b) and in part~(c) for
$1<p<2$ they are the contributions from pairs of magnitude bands. In
part~(d) and in part~(c) for $p=1$ they are the square roots of the
normalized absolute residual coordinates, and the polynomial comes from a
quantum mean estimator. The approximants for all centers lie in a space of
dimension at most $e^{L/8}$, and calibrated sampling
(Lemma~\ref{lem:s-cal}) turns their error into a support bound. Parts~(c)
and~(d) use a refinement of calibration, Lemma~\ref{lem:x-caps-cal}, whose
logarithmic factor depends on the approximation error instead of the rank.

\subsection{Setting for the Crossover Bound}

Throughout the proofs below, $0<p<2$, $\theta_p:=\min\{p,1\}$, and, as in the
rest of this section, $L:=\log(2/\epsilon)$. Zero rows do not
matter, so we assume that $A$ has nonzero rows $a_i^\top$, $i\in I$, and rank
$r\ge1$, and we use the coordinates Eq.~\eqref{lem:s-eq-lewis}:
$a_i^\top x=c_i^{1/p}u_i^\top y$ with unit vectors $u_i\in\R^r$,
$\sum_ic_iu_iu_i^\top=I_r$ and $\sum_ic_i=r$. All quantities below depend on a
query $x$ only through $y$, so we write queries as $y\in\R^r$. Thus
\[
 F(y):=\|Ax\|_p^p=\sum_{i\in I}c_i|u_i^\top y|^p,\qquad
 \mathbb B:=\{y:F(y)\le1\},
\]
and we use the base law $\pi_i:=c_i/r$ of Lemma~\ref{lem:s-cal} with
$\beta=1$, for which $\kappa_y(i)=r|u_i^\top y|^p$. With $\nu_i:=c_i/r$ these
coordinates form the Lewis setup Eq.~\eqref{lem:t-eq-lewis}, whose energy
$\sum_i\nu_i|\langle u_i,y\rangle|^p$ equals $F(y)/r$. On $\mathbb B$ we have
$\|y\|_2\le1$, as shown at the start of the proof of Lemma~\ref{lem:s-cal},
so $|u_i^\top y|\le1$ for every $i$. By Theorem~\ref{thm:general_lp_field},
$S_p(r,\epsilon)\le Cr/(\theta_p\epsilon^2)$ in general. For $1\le L\le r$
we use the elementary bounds
\begin{equation}\label{lem:x-caps-eq-params}
 \frac L{2\Lambda}\le K\le L,\qquad
 k\log\frac{er}k\le c(2+\log\frac1c)L
 \quad\text{for }0<c\le1\text{ and integers }1\le k\le cK .
\end{equation}
Indeed, $1\le\log(er/L)\le\log(er)\le2\Lambda$. Next, the identity
$er/K=(er/L)\log(er/L)$ gives $er/K\le(er/L)^2$, hence
$K\log(er/K)\le2L$. Finally, $x\mapsto x\log(er/x)$ increases on $(0,r]$,
so $k\log(er/k)\le cK(\log(er/K)+\log(1/c))$.

\subsection{Frame Coordinates}

\begin{lemma}[Frame coordinates]\label{lem:x-caps-frame}
Let $1\le k\le r$ be an integer. Let $C_0$ and $C'$ be the absolute
constants of Lemma~\ref{lem:t-frame-ensemble} and define
$R:=\max\{1,C'\}^{1/p}$. There is a finite set $\mathcal P$ with
$\log|\mathcal P|\le C_p(\Lambda+k\log(er/k))$, and for each
$\varpi\in\mathcal P$ an integer $N_\varpi\in[r,2C_0r]$ and reals
$b_{il}=b^\varpi_{il}$ and $\alpha_i=\alpha^\varpi_i$, $i\in I$,
$l\in[N_\varpi]$, with
\begin{equation}\label{lem:x-caps-eq-frame}
 |b_{il}|\le1,\qquad\sum_lb_{il}^2\le1,\qquad\sum_ic_ib_{il}^2\le1,
 \qquad|\alpha_i|\le R,
\end{equation}
such that the following holds. Every $y\in\mathbb B$ has some
$\varpi\in\mathcal P$ and some $e\in\R^{N_\varpi}$ with
$u_i^\top y=\alpha_i+\sum_lb_{il}e_l$ for all $i\in I$ and
\[
 \|e\|_\infty\le Rk^{-1/p},\quad\|e\|_p\le R,\quad
 \|e\|_2\le Rk^{1/2-1/p},\quad
 \|e\|_1\le\begin{cases}Rk^{1-1/p},&p\le1,\\
 2C_0Rr^{1-1/p},&p\ge1.\end{cases}
\]
We call $\varpi$ a \emph{pair} for $y$ and $e$ its \emph{residual}.
\end{lemma}

\begin{proof}
Let $\mathcal L$ be the family of lists of Lemma~\ref{lem:t-frame-ensemble}.
For $y\in\mathbb B\setminus\{0\}$ that lemma, applied to $y/\|y\|_2$ in the
Lewis setup described above, gives a list $\mathbf u=(u_{i_1},\ldots,
u_{i_N})$ with $\sum_l|\langle u_{i_l},y\rangle|^p\le C'r\cdot F(y)/r\le
R^p$, by homogeneity. For $y=0$ take any list. Let $U$ be the matrix with
columns $u_{i_l}$ and $S:=UU^\top$, so that $(C_0/2)I_r\preceq S\preceq
2C_0I_r$ and $r\le N\le2C_0r$. Define $z:=U^\top y$ and
$b_{il}:=\langle u_i,S^{-1}u_{i_l}\rangle$. Then $y=S^{-1}Uz$, so
$u_i^\top y=\sum_lb_{il}z_l$, and $\|z\|_p\le R$. Since $C_0\ge2$ we have
$\|S^{-1}\|\le1$. Hence $|b_{il}|\le1$,
$\sum_lb_{il}^2=u_i^\top S^{-1}UU^\top S^{-1}u_i=u_i^\top S^{-1}u_i\le1$,
and $\sum_ic_ib_{il}^2=\|S^{-1}u_{i_l}\|_2^2\le1$ because
$\sum_ic_iu_iu_i^\top=I_r$.

Apply Lemma~\ref{lem:t-pball-net} with $n:=N\ge r\ge k$, $k$ and $R$. It
gives a center $z^0\in\mathcal N_{N,k}(R)$ of $z$, and $e:=z-z^0$ satisfies
the four displayed bounds: those for $p\le1$ directly, and for $p\ge1$ the
last one follows from H\"older's inequality,
$\|e\|_1\le N^{1-1/p}\|e\|_p\le(2C_0r)^{1-1/p}R$. Define
$\alpha_i:=\sum_lb_{il}z^0_l$, so $u_i^\top y=\alpha_i+\sum_lb_{il}e_l$.
Every center in $\mathcal N_{N,k}(R)$ is the center of some vector with
$\ell_p$ norm at most $R$, and its entries are dominated by those of that
vector. Hence $\|z^0\|_2\le\|z^0\|_p\le R$ and
$|\alpha_i|\le\|b_{i\cdot}\|_2\|z^0\|_2\le R$. Let $\mathcal P$ be the set of
pairs $(\mathbf u,z^0)$ with $\mathbf u\in\mathcal L$ and
$z^0\in\mathcal N_{|\mathbf u|,k}(R)$. Then
$\log|\mathcal P|\le\log(C_pr\log(2r))+k\log(2eC_0r/k)+C_pk$, which gives
the bound on $|\mathcal P|$.
\end{proof}

\subsection{Polynomials on Boxes}

For $\theta\in\R$ define $\mathfrak s(\theta):=(1-\cos\theta)/2=\sin^2(\theta/2)$.
It maps $[0,\pi]$ onto $[0,1]$. We write $\mathsf T_j$ for the Chebyshev polynomials,
so that $\cos(j\theta)=\mathsf T_j(\cos\theta)=\mathsf T_j(1-2\mathfrak s(\theta))$, and
$\mathbb T:=[-\pi,\pi]$. For an integrable function $J$ on $\mathbb T$ let
$\mathfrak d(J):=\int_{\mathbb T}J(t)(1-\cos t)\,\d t$ be its cosine deficit.

\begin{lemma}[Coefficients from the supremum]\label{lem:x-caps-coef}
Let $D\ge1$ and $q_1,\ldots,q_D\ge0$ be integers, and let $P$ be a real
polynomial in $D$ variables of degree at most $q_h$ in the $h$th variable.
Then the absolute values of the monomial coefficients of $P$ sum to at most
$2^D\prod_{h\le D}((q_h+1)8^{q_h})\sup_{[0,1]^D}|P|$.
\end{lemma}

\begin{proof}
Define $\Gamma(\theta):=P(\mathfrak s(\theta_1),\ldots,\mathfrak s(\theta_D))$.
Each power $\mathfrak s(\theta)^m$ is a polynomial of degree $m$ in
$\cos\theta$, and each power $\cos^m\theta$ is a combination of the
$\cos(j\theta)$, $j\le m$. Hence
$\Gamma(\theta)=\sum_{\mathbf j}\gamma_{\mathbf j}\prod_h\cos(j_h\theta_h)$,
where $\mathbf j$ ranges over $\prod_h\{0,\ldots,q_h\}$. The functions
$\cos(j\theta)$ are orthogonal on $[0,\pi]$, with squared norm $\pi$ for
$j=0$ and $\pi/2$ for $j\ge1$. Writing $w_0:=1$ and $w_j:=2$ for $j\ge1$, we
get
\[
 |\gamma_{\mathbf j}|=|\pi^{-D}\prod_hw_{j_h}\int_{[0,\pi]^D}\Gamma(\theta)
 \prod_h\cos(j_h\theta_h)\,\d\theta|\le2^D\sup_{[0,1]^D}|P| .
\]
The polynomials $P(s)$ and $\sum_{\mathbf j}\gamma_{\mathbf j}\prod_h
\mathsf T_{j_h}(1-2s_h)$ agree on $[0,1]^D$, so they coincide. By the recurrence
$\mathsf T_{j+1}(u)=2u\mathsf T_j(u)-\mathsf T_{j-1}(u)$ with $\mathsf T_0=1$ and $\mathsf T_1(u)=u$, the absolute
coefficient sum of $\mathsf T_j$ is at most $(1+\sqrt2)^j$, since $x=1+\sqrt2$
satisfies $x^2=2x+1$. As $(1-2s)^m$ has absolute coefficient sum $3^m$, the
polynomial $\mathsf T_j(1-2s)$ has absolute coefficient sum at most
$(3+3\sqrt2)^j\le8^j$. Absolute coefficient sums multiply for products of
polynomials in disjoint variables, and there are $\prod_h(q_h+1)$
multi-indices.
\end{proof}

\begin{lemma}[A pair of Jackson kernels]\label{lem:x-caps-kernel}
There are absolute constants $q_0\ge8$ and $C_J\ge1$ such that for every
integer $q\ge q_0$ there are nonnegative even trigonometric polynomials
$\mathcal J_q,\mathcal J'_q$ on $\mathbb T$ of degree at most $q$ with
$\int_{\mathbb T}\mathcal J_q=\int_{\mathbb T}\mathcal J'_q=1$, $\mathfrak d(\mathcal J'_q)=2\mathfrak
d(\mathcal J_q)$, and $\int_{\mathbb T}(\mathcal J_q+\mathcal J'_q)(t)|t|^v\,\d t\le C_Jq^{-v}$ for
$0\le v\le4$.
\end{lemma}

\begin{proof}
For $n\ge1$ let $J_n(s):=\lambda_n^{-1}(\sin(ns/2)/\sin(s/2))^8$ with
$\int_{\mathbb T}J_n=1$. These are the kernels of Step~2 of the proof of
Lemma~\ref{lem:t-tv-approx} with $\ell=4$: nonnegative even trigonometric
polynomials of degree at most $4(n-1)$. The derivation
of Eq.~\eqref{lem:t-eq-jackson} uses only $\ell\ge1$ for the pointwise bound and
$2\ell\ge v+2$ for the moment bound. With $\ell=4$ it gives
$\int_{\mathbb T}J_n(s)|s|^v\,\d s\le2C_4n^{-v}$ for $0\le v\le6$, where
$C_4:=(\pi(1+\pi)/2)^8/2$. Hence $\mathfrak d(J_n)\le\frac12\int J_ns^2\le
C_4n^{-2}$, using $1-\cos s\le s^2/2$.

For a lower bound on the deficit, recall from that derivation that the base
$|\sin(ns/2)/\sin(s/2)|$ is at most $\min\{n,\pi/|s|\}$ on $\mathbb T$ and at
least $2n/\pi$ for $|s|\le1/n$. The first bound gives
$\lambda_n\le2\pi n^7+2\pi^8\int_{\pi/n}^\infty s^{-8}\,\d s\le3\pi n^7$, so
$J_n(s)\ge(2n/\pi)^8/(3\pi n^7)$ for $|s|\le1/n$. Since
$1-\cos s=2\sin^2(s/2)\ge2s^2/\pi^2$ for $|s|\le\pi$,
\[
 \mathfrak d(J_n)\ge\frac{2^8n}{3\pi^9}\int_{-1/n}^{1/n}\frac{2s^2}{\pi^2}\,\d s
 =c_\delta n^{-2},\qquad c_\delta:=\frac{2^{10}}{9\pi^{11}} .
\]
Let $D_0:=\lceil(2C_4/c_\delta)^{1/2}\rceil$ and $q_0:=8D_0$. For $q\ge q_0$
define $n:=\lfloor q/4\rfloor$ and $n':=\lfloor n/D_0\rfloor$. Then $n\ge2D_0$,
$n\ge q/8$, $n'\ge n/(2D_0)\ge q/(16D_0)$, and both $J_n$ and $J_{n'}$ have
degree at most $q$. Moreover
$\mathfrak d(J_{n'})\ge c_\delta n'^{-2}\ge c_\delta D_0^2n^{-2}\ge
2C_4n^{-2}\ge2\mathfrak d(J_n)$. Define
$\vartheta:=(\mathfrak d(J_{n'})-2\mathfrak d(J_n))/(\mathfrak
d(J_{n'})-\mathfrak d(J_n))\in[0,1)$, $\mathcal J_q:=J_n$ and
$\mathcal J'_q:=\vartheta J_n+(1-\vartheta)J_{n'}$. Then $\mathfrak
d(\mathcal J'_q)=\mathfrak d(J_{n'})-\vartheta(\mathfrak d(J_{n'})-\mathfrak
d(J_n))=2\mathfrak d(J_n)$. Finally $\mathcal J_q+\mathcal J'_q\le2J_n+J_{n'}$, so the
moment bound holds with $C_J:=6C_4(16D_0)^4$.
\end{proof}

\begin{lemma}[Approximation on a box]\label{lem:x-caps-box}
There is an absolute constant $C_{\rm box}$ with the following property.
Let $D\ge1$, $\alpha\in\R$, $\rho_1,\ldots,\rho_D>0$,
$\sigma_1,\ldots,\sigma_D\in\{-1,1\}$, and let $q_1,\ldots,q_D\ge q_0$ be
integers. There is a polynomial $P$ in $D$ variables, of degree at most
$q_h$ in the $h$th variable, with
$\sup_{[0,1]^D}|P|\le3(|\alpha|+\sum_h\rho_h)^p$ and
\begin{equation}\label{lem:x-caps-eq-box}
 ||\alpha+\sum_{h\le D}\sigma_ht_h|^p
 -P(\frac{t_1}{\rho_1},\ldots,\frac{t_D}{\rho_D})|
 \le C_{\rm box}(\sum_{h\le D}(\frac{(\rho_ht_h)^{1/2}}{q_h}
 +\frac{\rho_h}{q_h^2}))^p
 \qquad(t\in\prod_h[0,\rho_h]).
\end{equation}
\end{lemma}

\begin{proof}
Let $\mathcal J_h:=\mathcal J_{q_h}$ and $\mathcal J'_h:=\mathcal J'_{q_h}$ be the kernels of
Lemma~\ref{lem:x-caps-kernel} and
\[
 \mathcal U(\phi):=2\prod_{h\le D}\mathcal J_h(\phi_h)-\prod_{h\le D}\mathcal J'_h(\phi_h)
 \qquad(\phi\in\mathbb T^D).
\]
It is a trigonometric polynomial of degree at most $q_h$ in $\phi_h$, even in
each variable, with $\int\mathcal U=1$ and $\int|\mathcal U|\le3$. Integrating
out the other variables gives $\int\mathcal U(\phi)(1-\cos\phi_h)\,\d\phi=
2\mathfrak d(\mathcal J_h)-\mathfrak d(\mathcal J'_h)=0$, and $\int\mathcal U(\phi)\sin\phi_h\,
\d\phi=0$ by evenness in $\phi_h$.

Define $\psi(\theta):=|\alpha+\sum_h\sigma_h\rho_h\mathfrak s(\theta_h)|^p$ for
$\theta\in\R^D$. It is continuous, $2\pi$-periodic and even in each
variable. Let
\[
 \widehat P(\theta):=\int_{\mathbb T^D}\mathcal U(\phi)\psi(\theta-\phi)\,
 \d\phi=\int_{\mathbb T^D}\psi(\vartheta)\,\mathcal U(\theta-\vartheta)\,
 \d\vartheta .
\]
The second form shows that $\widehat P$ is a real trigonometric polynomial of
degree at most $q_h$ in $\theta_h$. The substitution
$\vartheta_h\mapsto-\vartheta_h$ shows that it is even in each $\theta_h$.
Averaging over the $2^D$ sign changes of the variables, $\widehat P$ is a
combination of the products $\prod_h\cos(j_h\theta_h)$ with $j_h\le q_h$.
Since $\cos(j\theta)=\mathsf T_j(1-2\mathfrak s(\theta))$, there is a polynomial $P$
of degree at most $q_h$ in the $h$th variable with
$\widehat P(\theta)=P(\mathfrak s(\theta_1),\ldots,\mathfrak s(\theta_D))$.
As $\mathfrak s([0,\pi])=[0,1]$, we get $\sup_{[0,1]^D}|P|\le\int|\mathcal
U|\sup|\psi|\le3(|\alpha|+\sum_h\rho_h)^p$.

Fix $t\in\prod_h[0,\rho_h]$ and choose $\theta_h\in[0,\pi]$ with
$\rho_h\mathfrak s(\theta_h)=t_h$, so that $P(t_1/\rho_1,\ldots,t_D/\rho_D)=
\widehat P(\theta)$. Define $\chi:=\alpha+\sum_h\sigma_ht_h$ and
$\Delta(\phi):=\sum_h\sigma_h\rho_h(\mathfrak s(\theta_h-\phi_h)-\mathfrak
s(\theta_h))$. Then $\psi(\theta-\phi)=|\chi+\Delta(\phi)|^p$. Since
\[
 \mathfrak s(\theta_h-\phi_h)-\mathfrak s(\theta_h)
 =\tfrac12\cos\theta_h\,(1-\cos\phi_h)-\tfrac12\sin\theta_h\sin\phi_h ,
\]
the two moment identities give $\int\mathcal U\Delta=0$. Moreover
$|\sin\theta_h|=2|\sin(\theta_h/2)\cos(\theta_h/2)|\le2\mathfrak
s(\theta_h)^{1/2}$ and $1-\cos\phi\le\phi^2/2$, hence
\[
 |\Delta(\phi)|\le\sum_h\Delta_h(\phi_h),\qquad
 \Delta_h(\phi_h):=(\rho_ht_h)^{1/2}|\phi_h|+\rho_h\phi_h^2/4 .
\]
For reals $u,w$ we have $||u+w|^p-|u|^p|\le|w|^p$ if $p\le1$, by
subadditivity of $s\mapsto s^p$ on $[0,\infty)$. If $1<p<2$, define
$\gamma:=p-1\in(0,1)$. The map $u\mapsto\operatorname{sgn}(u)|u|^\gamma$
satisfies $|\operatorname{sgn}(u)|u|^\gamma-\operatorname{sgn}(u')|u'|^\gamma|
\le2^{1-\gamma}|u-u'|^\gamma$, by subadditivity for equal signs and by
concavity for opposite signs. Hence
$||u+w|^p-|u|^p-p\operatorname{sgn}(u)|u|^{p-1}w|\le\int_0^1p2^{1-\gamma}
|\lambda w|^\gamma|w|\,\d\lambda=2^{2-p}|w|^p$. Since $\int\mathcal U=1$ and
$\int\mathcal U\Delta=0$,
\[
 \widehat P(\theta)-|\chi|^p=\int\mathcal U(\phi)(|\chi+\Delta(\phi)|^p-|\chi|^p
 -\mathbf 1_{\{p>1\}}\,p\operatorname{sgn}(\chi)|\chi|^{p-1}\Delta(\phi))\,\d\phi ,
\]
and therefore $|\widehat P(\theta)-|\chi|^p|\le2\int(2\Pi+\Pi')(\sum_h
\Delta_h)^p$, where $\Pi:=\prod_h\mathcal J_h$ and $\Pi':=\prod_h\mathcal J'_h$. Let
$\Pi''\in\{\Pi,\Pi'\}$ have factors $\pi_h$, and define
$E:=\sum_h((\rho_ht_h)^{1/2}q_h^{-1}+\rho_hq_h^{-2})$. If $p\ge1$,
Minkowski's inequality (Lemma~\ref{lem:pre_minkowski}(a)) in $L_p(\Pi'')$, then in $L_p(\pi_h)$, and the moment
bounds of Lemma~\ref{lem:x-caps-kernel} at the orders $p$ and $2p\le4$ give
\[
 (\int\Pi''(\sum_h\Delta_h)^p)^{1/p}
 \le\sum_h((\rho_ht_h)^{1/2}(C_Jq_h^{-p})^{1/p}
 +\frac{\rho_h}4(C_Jq_h^{-2p})^{1/p})\le C_J^{1/p}E .
\]
If $p<1$, Jensen's inequality and the moments of orders one and two give
$\int\Pi''(\sum_h\Delta_h)^p\le(\sum_h\int\pi_h\Delta_h)^p\le(C_JE)^p$. As
$C_J\ge1$, in both cases the integral is at most $C_JE^p$, so the error is at
most $6C_JE^p$. Take $C_{\rm box}:=6C_J$.
\end{proof}

\begin{lemma}[A sketched common space]\label{lem:x-caps-space}
There is an absolute constant $C_{\rm sk}$ with the following property. Let
$D,q,r'\ge1$ be integers, $\mathsf A\ge1$ and $0<\tau<1/2$. For each $i$ in
a finite set $I$ let $P_i$ be a polynomial in $D$ variables, of degree at
most $q$ in each variable, whose monomial coefficients have absolute values
summing to at most $\mathsf A$, and let $w_{i,1},\ldots,w_{i,D}\in
B_2^{r'}$. There is a space $\mathcal W$ of functions on $I$ with
\[
 \log(2\dim\mathcal W)\le D\log(q+1)+C_{\rm sk}(\log(\mathsf A/\tau)
 +\log(2Dr'q))
\]
such that every $z\in B_2^{r'}$ has some $f\in\mathcal W$ with
$|f(i)-P_i(\langle w_{i,1},z\rangle,\ldots,\langle w_{i,D},z\rangle)|\le2\tau$
for all $i\in I$.
\end{lemma}

\begin{proof}
Define $\eta:=\tau/(\mathsf ADq)$, and let $Y\subseteq B_2^{r'}$ be a set such
that every point of $B_2^{r'}$ is within $\eta$ of $Y$, with
$|Y|\le(1+2/\eta)^{r'}$ (Lemma~\ref{lem:pre_nets}).
Lemma~\ref{lem:t-sketch} with $a:=D$ gives a matrix $\widetilde M$ indexed by
$I\times Y$, of rank at most $d_{\rm sk}$, with
$|\widetilde M_{i,y}-P_i(\langle w_{i,1},y\rangle,\ldots,\langle
w_{i,D},y\rangle)|\le\tau$. Let $\mathcal W$ be its column space. For
$z\in B_2^{r'}$ choose $y\in Y$ with $\|z-y\|_2\le\eta$ and let $f$ be the
column of $y$. All inner products $\langle w_{i,h},z\rangle$ and $\langle
w_{i,h},y\rangle$ lie in $[-1,1]$, where every partial derivative of $P_i$
is at most $q\mathsf A$ in absolute value, and they differ by at most
$\eta$. Hence the values of $P_i$ at the two points differ by at most
$Dq\mathsf A\eta=\tau$, which proves the approximation.

The constant $C_a$ of Lemma~\ref{lem:t-sketch} hides the dependence on $a$,
so we use the explicit form stated there. With $a=D$ it gives
$d_{\rm sk}:=1+(q+1)^Dd_0$ with
$\log d_0\le C+2\log(8\mathsf A/\tau)+\log(2Dr')+\log\log(e+16Dq\mathsf
A/\tau)+\log\log(e|Y|)$, where $C$ is absolute. Here
$\log(e|Y|)\le1+r'\log(3/\eta)\le2r'\log(3\mathsf ADq/\tau)$, so
$\log\log(e|Y|)\le\log(2r')+\log(3\mathsf ADq/\tau)$. Since
$\log(2d_{\rm sk})\le\log4+D\log(q+1)+\log d_0$ and $\log(\mathsf A/\tau)
\ge\log2$, the claim follows.
\end{proof}

\subsection{Calibration with the Entropy of the Energy Ball}

For a nonnegative function $\xi$ on $I$ define
$d_\xi(y,y'):=(\sum_i\xi(i)(\kappa_y(i)-\kappa_{y'}(i))^2)^{1/2}$, a
pseudometric on $\R^r$. As in Lemma~\ref{lem:s-cal}, for a pseudometric $d$
and $t>0$ let $N(\mathbb B,d,t)$ be the least number of $d$-balls of
radius~$t$ with centers in $\mathbb B$ that cover $\mathbb B$.

\begin{lemma}[Entropy of the energy ball]\label{lem:x-caps-pack}
Let $a\ge1$, and let $\xi\ge0$ on $I$ satisfy $\sum_i\xi(i)\le a$ and
$\sup_{y\in\mathbb B}\sum_i\xi(i)\kappa_y(i)\le a$. Then
$N(\mathbb B,d_\xi,t)\le(C_{p,a}/t)^{r/\theta_p}$ for $0<t\le1$.
\end{lemma}

\begin{proof}
\emph{Volume.} With $F_\nu:=F/r$, the energy of the Lewis setup, the set
$\mathbb B=\{rF_\nu\le1\}$ is the set of Lemma~\ref{lem:t-inverse-energy}. It
is star-shaped with radial function $F(\vartheta)^{-1/p}$ on
the sphere, so integration in polar coordinates and that lemma give
\[
 \operatorname{vol}(\mathbb B)=\operatorname{vol}(B_2^r)\int_{S^{r-1}}
 F^{-r/p}\,\d\omega=\operatorname{vol}(B_2^r)\,r^{-r/p}
 \|F_\nu^{-1}\|_{L_{r/p}(\omega)}^{r/p}\le(C_{\rm v}b_r)^r
 \operatorname{vol}(B_2^r),
\]
where $b_r:=r^{1/2-1/p}$ and $C_{\rm v}$ is the $1/p$th power of the constant
of Lemma~\ref{lem:t-inverse-energy}.

\emph{Perturbations.} Let $Z$ be uniform on $\varrho B_2^r$, $\varrho>0$,
and let $u$ be a unit vector. Write $Z=\|Z\|_2\vartheta$ with $\vartheta$
uniform on $S^{r-1}$ and independent of $\|Z\|_2$, whose density on
$[0,\varrho]$ is $rs^{r-1}/\varrho^r$. Then $\E[\|Z\|_2^{2m}]=r\varrho^{2m}/(r+2m)$.
The same decomposition of a standard Gaussian vector gives
$\E[\vartheta_1^2]=1/r$ and $\E[\vartheta_1^4]=3/(r(r+2))$. By rotation
invariance, $\E[\langle u,Z\rangle^2]=\varrho^2/(r+2)\le\varrho^2/r$ and
$\E[\langle u,Z\rangle^4]=3\varrho^4/((r+2)(r+4))\le3\varrho^4/r^2$. Hence, by
Jensen's inequality, $\E[|\langle u,Z\rangle|^v]\le(\varrho^2/r)^{v/2}$ for
$0<v\le2$ and $\E[|\langle u,Z\rangle|^v]\le3(\varrho^2/r)^{v/2}$ for
$2\le v\le4$.

Fix $y\in\mathbb B$ and $0<t\le1$. Suppose first $1\le p<2$ and take
$\varrho:=tb_r$, so $\varrho^2/r=t^2r^{-2/p}$. Then
$\E[F(Z)]\le\sum_ic_i(\varrho^2/r)^{p/2}=t^p$. For reals $w,h$ the mean value
theorem and subadditivity of $s\mapsto s^{p-1}$ give
$||w+h|^p-|w|^p|\le p|h|(|w|^{p-1}+|h|^{p-1})$, with $0^0:=1$. Hence
\[
 \E[d_\xi(y,y+Z)^2]\le2p^2r^2\sum_i\xi(i)(\E[\langle u_i,Z\rangle^2
 \,|u_i^\top y|^{2p-2}]+\E[|\langle u_i,Z\rangle|^{2p}]).
\]
Define $\gamma:=(2p-2)/p\in[0,1)$. By H\"older's inequality and the hypotheses,
$\sum_i\xi(i)|u_i^\top y|^{2p-2}\le(\sum_i\xi(i))^{1-\gamma}(\sum_i\xi(i)
|u_i^\top y|^p)^\gamma\le a^{1-\gamma}(a/r)^\gamma=ar^{-\gamma}$, because
$\sum_i\xi(i)|u_i^\top y|^p=r^{-1}\sum_i\xi(i)\kappa_y(i)$. As
$r^2t^2r^{-2/p}r^{-\gamma}=t^2$ and $2p\le4$, we get
$\E[d_\xi(y,y+Z)^2]\le2p^2a(t^2+3t^{2p})\le32at^2$. By Markov's inequality (Lemma~\ref{lem:pre_markov}),
with probability at least $1/2$ both $F(Z)\le4t^p$ and
$d_\xi(y,y+Z)\le C_at$ hold, where $C_a:=(128a)^{1/2}$. On this event the
triangle inequality for $F^{1/p}$ gives $F(y+Z)^{1/p}\le1+4t$, so
$y+Z\in5\mathbb B$. If $0<p<1$, take instead $\varrho:=t^{1/p}b_r$, so
$(\varrho^2/r)^{p/2}=t/r$. Subadditivity gives $|\kappa_{y+Z}(i)-\kappa_y(i)|
\le r|\langle u_i,Z\rangle|^p$ and $F(y+Z)\le F(y)+F(Z)$. Hence
$\E[d_\xi(y,y+Z)^2]\le r^2a(\varrho^2/r)^p=at^2$ and $\E[F(Z)]\le t$, and with
probability at least $1/2$ we have $d_\xi(y,y+Z)\le C_at$ and
$y+Z\in5^{1/p}\mathbb B$.

\emph{Packing.} In both cases let $G_y$ be the set of points $y+z$ with
$z\in\varrho B_2^r$, $y+z\in5^{1/\theta_p}\mathbb B$ and $d_\xi(y,y+z)\le
C_at$. It has volume at least $\frac12\operatorname{vol}(\varrho B_2^r)$. If
$\mathcal S\subseteq\mathbb B$ has pairwise $d_\xi$-distances greater than
$2C_at$, the sets $G_y$, $y\in\mathcal S$, are disjoint by the triangle
inequality. They lie in $5^{1/\theta_p}\mathbb B$, so
\[
 |\mathcal S|\cdot\tfrac12\varrho^r\operatorname{vol}(B_2^r)\le
 5^{r/\theta_p}\operatorname{vol}(\mathbb B)\le5^{r/\theta_p}(C_{\rm v}b_r)^r
 \operatorname{vol}(B_2^r).
\]
As $\varrho=t^{1/\theta_p}b_r$, this gives $|\mathcal S|\le2(5C_{\rm
v}^{\theta_p}/t)^{r/\theta_p}$. A maximal such $\mathcal S$ is a
$2C_at$-cover of $\mathbb B$ with centers in $\mathbb B$. Since $C_a\ge1$,
replacing $t$ by $t/(2C_a)$ proves the lemma with
$C_{p,a}:=20C_aC_{\rm v}^{\theta_p}$.
\end{proof}

\begin{lemma}[Calibration with energy entropy]\label{lem:x-caps-cal}
Fix $B_0>0$. There is $C\ge1$ depending only on $p$ and $B_0$ with the
following property. Let $\mathcal V$ be a space of functions on $I$ with
$1\le\dim\mathcal V\le e^{L/8}$, and let $r^{-B_0}\le\delta\le1$ satisfy
$\operatorname{dist}_{L_2(\pi)}(\kappa_y,\mathcal V)\le\delta$ for every
$y\in\mathbb B$. If $L\ge C\log(2r)$, there are a set $S\subseteq I$ with
$|S|\le C\epsilon^{-2}r\delta^2\log(2+1/\delta)$ and positive weights
$\widehat w_i$, $i\in S$, with
$|\sum_{i\in S}\widehat w_i|a_i^\top x|^p-\|Ax\|_p^p|\le\epsilon\|Ax\|_p^p$
for every $x$.
\end{lemma}

\begin{proof}
We verify the hypothesis of Lemma~\ref{lem:s-cal}(b) with $\beta=1$, using
the facts listed in Lemma~\ref{lem:s-cal-entropy}. Let $\xi$ be the measure
on $I$ defined there. It has total mass at most $8$, and
$\widehat d\le\sqrt2\,d_\xi$. Outside an event of probability at most
$e^{-2L}$, which requires the lower limit for $L$ of that lemma, we have
$\sup_{y\in\mathbb B}\sum_i\xi(i)\kappa_y(i)\le5$. On the complement,
Lemma~\ref{lem:x-caps-pack} with $a=8$ applies. A $d_\xi$-ball of radius
$t/\sqrt2$ lies in the $\widehat d$-ball of radius $t$ with the same center,
so $N(\mathbb B,\widehat d,t)\le(C_1/t)^{r/\theta_p}$ for $0<t\le1$, with
$C_1\ge e$ depending only on $p$. For $t>1$ we use
$N(\mathbb B,\widehat d,t)\le N(\mathbb B,\widehat d,1)$.

Define $a_1:=\min\{8\delta,1\}$ and $t_1:=eC_1$. For $0<a'\le t_1/e$ we have
$\int_0^{a'}\sqrt{\log(t_1/t)}\,\d t\le2a'\sqrt{\log(t_1/a')}$, as in the
proof of Lemma~\ref{lem:s-cal}(a). Since $1/a_1\le1/\delta$ and $\delta\le1$,
\[
 \int_0^{8\delta}\sqrt{\log N(\mathbb B,\widehat d,t)}\,\d t
 \le\sqrt{\frac r{\theta_p}}(2a_1\sqrt{\log(eC_1/a_1)}
 +8\delta\sqrt{\log C_1})\le C_2\,\delta\sqrt{r\log(2+1/\delta)} .
\]
Let $M:=\lceil C_3\epsilon^{-2}r\delta^2\log(2+1/\delta)\rceil$ with
$C_3:=(C_2^2+1)/c_{\rm cal}^2$. Then Eq.~\eqref{lem:s-eq-entropy} holds outside
an event of probability at most $e^{-2L}\le e^{-3L/2}$, and
$\delta\le c_{\rm cal}\epsilon\sqrt M$. Moreover $M\ge C_3\epsilon^{-2}r
\delta^2\ge e^{2L}r^{-2B_0}/4\ge e^{1.8L}$ once $C$ is large. Part~(b) of
Lemma~\ref{lem:s-cal} gives the claim with $|S|\le M\le2C_3\epsilon^{-2}r
\delta^2\log(2+1/\delta)$.
\end{proof}

\subsection{Magnitude Bands}

For $t\in\R$ write $t_+:=\max\{t,0\}$ and $t_-:=\max\{-t,0\}$.

\begin{lemma}[Magnitude bands]\label{lem:x-caps-bins}
Let $1\le k\le r$ be an integer, let $\varpi\in\mathcal P$ be as in
Lemma~\ref{lem:x-caps-frame}, with $N:=N_\varpi$ and coefficients $b_{il}$
and $\alpha_i$, and define $\rho:=Rk^{1/2-1/p}$. There are finite sets
$\mathcal H$ and $\mathcal G$ with $|\mathcal H|\le C_p\Lambda^2$ and
$|\mathcal G|\le C_p\Lambda$, maps $f_g:\R^N\to[0,\infty)^N$ for
$g\in\mathcal G$, and for each band $h\in\mathcal H$ an integer
$|h|\in[0,C_p\Lambda]$, an index $g(h)\in\mathcal G$, a sign $\sigma_h$,
numbers $0<\rho_h\le\rho$ and $\bar\beta_h,\bar f_h,R_h>0$, and vectors
$\beta_{i,h}\in[0,\infty)^N$, $i\in I$, with the following properties. At
most $16(j+1)$ bands have $|h|=j$. Moreover $\bar\beta_h\bar f_h=R_h\rho_h$,
$\log(1+R_h)\le C(1+|h|)$ and $\|\beta_{i,h}\|_2\le\bar\beta_h$. Let
$y\in\mathbb B$ have the pair $\varpi$ and residual $e$, and define
$t_{ih}:=\langle\beta_{i,h},f_{g(h)}(e)\rangle$. Then
$\|f_{g(h)}(e)\|_2\le\bar f_h$, $0\le t_{ih}\le\rho_h$,
\begin{equation}\label{lem:x-caps-eq-trunc}
 ||u_i^\top y|^p-|\alpha_i+\sum_{h\in\mathcal H}\sigma_ht_{ih}
 |^p|\le C_pr^{-10}\qquad(i\in I),
\end{equation}
and $W_h:=\sum_ic_i(\rho_ht_{ih})^p$ satisfies
\begin{equation}\label{lem:x-caps-eq-overlap}
 W_h\le C_pr^{p-1}2^{-\gamma_p|h|}\ \ \text{if }1\le p<2,\
 \gamma_p:=\frac{p(p-1)}2;\qquad
 W_h\le C_pr^{1-p}k^{2p-2}\ \ \text{if }0<p<1 .
\end{equation}
\end{lemma}

\begin{proof}
Define $C_{\rm t}:=11+10/p$, $\varsigma_a:=2^ar^{-1/2}$ and
$v_b:=R2^br^{-1/p}$ for integers $a,b$. Let $\mathcal A$ be the set of
integers $a$ for which $[\varsigma_a,2\varsigma_a)$ meets
$[r^{-C_{\rm t}},1]$, and $\mathcal B$ the set of integers $b$ for which
$[v_b,2v_b)$ meets $[r^{-C_{\rm t}},Rk^{-1/p}]$. Then $\varsigma_a\le1$ and
$2\varsigma_a>r^{-C_{\rm t}}$ give $|a|\le C_{\rm t}\log_2(2r)+1$, and
$v_b\le Rk^{-1/p}\le R$ and $2v_b>r^{-C_{\rm t}}$ give
$|b|\le(C_{\rm t}+1/p)\log_2(2r)+|\log_2R|+1$. Let
$\mathcal G:=\mathcal B\times\{\pm1\}$ and
$\mathcal H:=\mathcal A\times\{\pm1\}\times\mathcal B\times\{\pm1\}$. For
$h=(a,\sigma_1,b,\sigma_2)$ define $\sigma_h:=\sigma_1\sigma_2$,
$|h|:=|a|+|b|$ and $g(h):=(b,\sigma_2)$. There are at most $4j+1$ integer
pairs $(a,b)$ with $|a|+|b|=j$, hence at most $16(j+1)$ bands with $|h|=j$.
Define
\begin{align*}
 f_{(b,\sigma_2)}(e)_l&:=|e_l|\ \text{ if }\operatorname{sgn}e_l=\sigma_2,\
 |e_l|\ge r^{-C_{\rm t}}\text{ and }|e_l|\in[v_b,2v_b),\ \text{ and }0
 \text{ otherwise},\\
 (\beta_{i,h})_l&:=|b_{il}|\ \text{ if }\operatorname{sgn}b_{il}=\sigma_1,\
 |b_{il}|\ge r^{-C_{\rm t}}\text{ and }|b_{il}|\in[\varsigma_a,
 2\varsigma_a),\ \text{ and }0\text{ otherwise}.
\end{align*}

\emph{Truncation.} Let $y,e$ be as in the statement. An index $l$ with
$|b_{il}|\ge r^{-C_{\rm t}}$ and $|e_l|\ge r^{-C_{\rm t}}$ lies in the
supports of $\beta_{i,h}$ and $f_{g(h)}(e)$ for exactly one band $h$, because
$|b_{il}|\le1$ and $|e_l|\le\|e\|_\infty\le Rk^{-1/p}$, and its term
$b_{il}e_l$ equals $\sigma_h|b_{il}||e_l|$. Hence $s'_i:=\sum_h\sigma_h
t_{ih}$ differs from $s_i:=\sum_lb_{il}e_l$ by at most
$r^{-C_{\rm t}}(\|e\|_1+\|b_{i\cdot}\|_1)\le r^{-C_{\rm t}}\sqrt N(\|e\|_2
+1)\le C_{\rm R}r^{1/2-C_{\rm t}}$, where $C_{\rm R}:=(R+1)(2C_0)^{1/2}$ and
we used $\|e\|_2\le R$. Now $|\alpha_i+s_i|=|u_i^\top y|\le1$. For $p\ge1$
the mean value theorem bounds the left side of Eq.~\eqref{lem:x-caps-eq-trunc}
by $p(1+C_{\rm R})^{p-1}C_{\rm R}r^{1/2-C_{\rm t}}$, and for $p<1$
subadditivity bounds it by $C_{\rm R}^pr^{-p(C_{\rm t}-1/2)}$. Both are at
most $C_pr^{-10}$.

\emph{Counts.} Fix $h=(a,\sigma_1,b,\sigma_2)$ and define
\[
 d_a:=\min\{N,\varsigma_a^{-2}\},\qquad
 n_b:=\min\{N,R^pv_b^{-p}\}=\min\{N,r2^{-pb}\},\qquad
 \ell_h:=\min\{d_a,n_b\}.
\]
As $\sum_lb_{il}^2\le1$,
at most $d_a$ indices have $|b_{il}|\ge\varsigma_a$. As $\|e\|_p\le R$, at
most $n_b$ indices have $|e_l|\ge v_b$. Hence the number $g_{ih}$ of indices
in both supports is at most $\ell_h$, each of them contributes less than
$4\varsigma_av_b$ to $t_{ih}$, and $t_{ih}\le\|b_{i\cdot}\|_2\|e\|_2\le\rho$
by Cauchy--Schwarz. Thus $t_{ih}\le\rho_h:=\min\{4\varsigma_av_b\ell_h,
\rho\}$. As $\sum_ic_ib_{il}^2\le1$ and $\sum_ic_i=r\le N$,
\[
 \sum_ic_ig_{ih}\le\sum_{l\in\supp f_{g(h)}(e)}\ \sum_ic_i
 \mathbf 1_{\{|b_{il}|\ge\varsigma_a\}}\le n_b\min\{r,\varsigma_a^{-2}\}
 \le d_an_b .
\]
The same counts give $\|\beta_{i,h}\|_2\le\bar\beta_h:=\min\{2\varsigma_a
\sqrt{d_a},1\}$ and $\|f_{g(h)}(e)\|_2\le\bar f_h:=\min\{2v_b\sqrt{n_b},
\rho\}$. Define $R_h:=\bar\beta_h\bar f_h/\rho_h$. If $\rho_h=\rho$, then
$R_h\le1$. Otherwise $R_h\le\sqrt{d_an_b}/\ell_h$. Since
$r\le N\le2C_0r$, with $a_+:=\max\{a,0\}$ and $b_+:=\max\{b,0\}$,
\begin{equation}\label{lem:x-caps-eq-counts}
 r2^{-2a_+}\le d_a\le2C_0r2^{-2a_+},\qquad
 r2^{-pb_+}\le n_b\le2C_0r2^{-pb_+} .
\end{equation}
Hence $\sqrt{d_an_b}/\ell_h\le(2C_0)^{1/2}2^{a_++b_+}$ and
$\log(1+R_h)\le C(1+|h|)$.

\emph{The overlap sums for $1\le p<2$.} Since $\rho_ht_{ih}\le
16\varsigma_a^2v_b^2\ell_hg_{ih}$ and $g_{ih}\le\ell_h$,
\[
 W_h\le16^p\varsigma_a^{2p}v_b^{2p}\ell_h^{2p-1}\sum_ic_ig_{ih}
 \le16^p\varsigma_a^{2p}v_b^{2p}\ell_h^{2p}\max\{d_a,n_b\} .
\]
Define $\mathsf x:=2^{-2a_+}$, $\mathsf z:=2^{-pb_+}$, $a_-:=\max\{-a,0\}$ and
$b_-:=\max\{-b,0\}$. Then Eq.~\eqref{lem:x-caps-eq-counts} gives
\[
 \ell_h\le2C_0r\min\{\mathsf x,\mathsf z\},\qquad\max\{d_a,n_b\}\le2C_0r\max\{\mathsf x,\mathsf z\},
\]
while $\varsigma_a^{2p}=r^{-p}\mathsf x^{-p}2^{-2pa_-}$ and
$v_b^{2p}=R^{2p}r^{-2}\mathsf z^{-2}2^{-2pb_-}$. Hence
$W_h\le C_pr^{p-1}2^{-2p(a_-+b_-)}\Phi$ with
$\Phi:=\mathsf x^{-p}\mathsf z^{-2}\min\{\mathsf x,\mathsf z\}^{2p}\max\{\mathsf x,\mathsf z\}$. If $\mathsf x\le\mathsf z$, then
$\Phi=\mathsf x^p\mathsf z^{-1}$, which is at most $\mathsf x^{p-1}$ and at most
$(\mathsf x/\mathsf z)^p\mathsf z^{p-1}\le\mathsf z^{p-1}$. If $\mathsf x\ge\mathsf z$, then $\Phi=(\mathsf z^2/\mathsf x)^{p-1}$ and
$\mathsf z^2/\mathsf x\le\min\{\mathsf x,\mathsf z\}$. In both cases $\Phi\le(\mathsf x\mathsf z)^{(p-1)/2}=
2^{-(p-1)a_+-p(p-1)b_+/2}$. As $|a|=a_++a_-$, $|b|=b_++b_-$ and
$\gamma_p\le\min\{p-1,2p\}$, the first bound in Eq.~\eqref{lem:x-caps-eq-overlap}
follows.

\emph{The overlap sums for $0<p<1$.} By concavity,
$\sum_ic_ig_{ih}^p\le r^{1-p}(\sum_ic_ig_{ih})^p\le r^{1-p}(d_an_b)^p$. The
definitions of $d_a$ and $n_b$ give $\varsigma_a^2\le d_a^{-1}$ and
$v_b^p\le R^pn_b^{-1}$. Hence
\[
 W_h\le16^p\varsigma_a^{2p}v_b^{2p}\ell_h^p\,r^{1-p}(d_an_b)^p
 \le16^pR^{2p}r^{1-p}\ell_h^pn_b^{p-2}\le16^pR^{2p}r^{1-p}n_b^{2p-2} .
\]
Finally $v_b\le Rk^{-1/p}$ gives $R^pv_b^{-p}\ge k$, and $N\ge r\ge k$, so
$n_b\ge k$ and $n_b^{2p-2}\le k^{2p-2}$.
\end{proof}

\subsection{Adaptive Degrees Above One}

\begin{lemma}[Adaptive degrees]\label{lem:x-caps-adaptive}
Fix $1<p<2$. There are $r_p$ and $C_p$ such that for $r\ge r_p$ and
$\sqrt r\le L\le r$ some space $\mathcal V$ of functions on $I$ with
$1\le\dim\mathcal V\le e^{L/8}$ satisfies
$\operatorname{dist}_{L_2(\pi)}(\kappa_y,\mathcal V)\le C_p(r/L^2)^{p/2}$ for
every $y\in\mathbb B$. The constants $r_p$ and $C_p$ given by the proof
grow without bound as $p\downarrow1$.
\end{lemma}

\begin{proof}
Constants $C_1,C_2,\ldots$ depend only on $p$. Define
$\zeta_p:=\gamma_p/(4p)=(p-1)/8$, with $\gamma_p$ as
in Eq.~\eqref{lem:x-caps-eq-overlap}, and let $c_1,c_3\in(0,1]$ depend only on
$p$, to be fixed at the end. Define $k:=\lfloor c_1K\rfloor$, and let $r$ be
large in terms of $p$, $c_1$ and $c_3$. Since $L\ge\sqrt r$,
Eq.~\eqref{lem:x-caps-eq-params} gives $K\ge\sqrt r/(2\Lambda)$, hence
$k\ge c_1K/2\ge c_1L/(4\Lambda)$, and $k\le r$. Let $\mathcal P$ be the set
of Lemma~\ref{lem:x-caps-frame} for this $k$. Fix $\varpi\in\mathcal P$,
take the bands of Lemma~\ref{lem:x-caps-bins}, define $D:=|\mathcal H|\le
C_1\Lambda^2$, and for $h\in\mathcal H$ let
\[
 q_h:=\lceil L^{3/4}+c_3L2^{-\zeta_p|h|}\rceil,\qquad
 \lambda_h:=2^{-3-|h|/2}.
\]
Then $q_0\le q_h\le3L$, and $\sum_h\lambda_h^2\le\sum_j16(j+1)2^{-6-j}=1$.

\emph{Approximation.} For $i\in I$ let $P_i$ be the polynomial of
Lemma~\ref{lem:x-caps-box} with this $D$, $\alpha:=\alpha_i$, the signs
$\sigma_h$, the sides $\rho_h$ and the degrees $q_h$. Let $y\in\mathbb B$
have the pair $\varpi$ and residual $e$, and define
$E_i:=\sum_h((\rho_ht_{ih})^{1/2}q_h^{-1}+\rho_hq_h^{-2})$. By
Lemmas~\ref{lem:x-caps-box} and~\ref{lem:x-caps-bins},
\[
 ||u_i^\top y|^p-P_i((t_{ih}/\rho_h)_h)|
 \le C_2r^{-10}+C_{\rm box}E_i^p\qquad(i\in I).
\]
Minkowski's inequality (Lemma~\ref{lem:pre_minkowski}(a)) in $L_{2p}$ of the weights $c_i$, whose total mass is
$r$, gives
\[
 (\sum_ic_iE_i^{2p})^{1/(2p)}\le\sum_h\frac{W_h^{1/(2p)}}{q_h}
 +r^{1/(2p)}\sum_h\frac{\rho_h}{q_h^2} .
\]
Let $C_3$ be the constant of Eq.~\eqref{lem:x-caps-eq-overlap}. By that bound
and $q_h\ge c_3L2^{-\zeta_p|h|}$, the first sum is at most
\[
 \frac{C_3^{1/(2p)}r^{(p-1)/(2p)}}{c_3L}\sum_h2^{-\zeta_p|h|}\le
 \frac{C_4r^{(p-1)/(2p)}}{c_3L},
\]
because $\gamma_p/(2p)-\zeta_p=\zeta_p$ and
$\sum_h2^{-\zeta_p|h|}\le16\sum_j(j+1)2^{-\zeta_pj}$. By $q_h\ge L^{3/4}$,
$\rho_h\le\rho=Rk^{-(2-p)/(2p)}$ and $k\ge c_1L/(4\Lambda)$, the second term
is at most
\[
 C_1\Lambda^2r^{1/(2p)}R(\frac{4\Lambda}{c_1L})^{(2-p)/(2p)}
 L^{-3/2} .
\]
Its ratio to $r^{(p-1)/(2p)}/L$ is at most
$C_5\Lambda^3r^{(2-p)/(2p)}L^{-1/2-(2-p)/(2p)}\le C_5\Lambda^3
r^{-(p-1)/(2p)}$, using $L\ge\sqrt r$. This is at most one for $r$ large,
since $p>1$. Hence $\sum_ic_iE_i^{2p}\le C_6r^{p-1}L^{-2p}$.

\emph{A common space.} Let $w_{i,h}\in\R^{DN}$ have the block
$\beta_{i,h}/\bar\beta_h$ in position $h$ and zeros elsewhere, and let
$z(e)\in\R^{DN}$ have the block $\lambda_hf_{g(h)}(e)/\bar f_h$ in position
$h$. Then $\|w_{i,h}\|_2\le1$, $\|z(e)\|_2^2\le\sum_h\lambda_h^2\le1$, and
$t_{ih}/\rho_h=(R_h/\lambda_h)\langle w_{i,h},z(e)\rangle$. Define
$\widetilde P_i(s):=P_i((R_hs_h/\lambda_h)_h)$, so that
$P_i((t_{ih}/\rho_h)_h)=\widetilde P_i((\langle w_{i,h},z(e)\rangle)_h)$. By
Lemma~\ref{lem:x-caps-coef} and $\sup_{[0,1]^D}|P_i|\le3(R+D\rho)^p$, the
monomial coefficients of $\widetilde P_i$ have absolute sum at most
$\mathsf A$, where
\[
 \log\mathsf A\le D\log2+\sum_h(\log(q_h+1)+q_h\log8)
 +\log(3(R+DR)^2)+\sum_hq_h\log\max\{1,R_h/\lambda_h\}.
\]
As $\log(1+R_h)\le C(1+|h|)$, $\log(1/\lambda_h)\le3+|h|$ and $q_h\le3L$, this
gives $\log\mathsf A\le C_7\Lambda^3+C_8\sum_hq_h(1+|h|)$ with $C_8$
absolute. Moreover, as $|h|\le C_1\Lambda$,
\[
 \sum_hq_h(1+|h|)\le(L^{3/4}+1)D(1+C_1\Lambda)+c_3L\cdot16\sum_j(j+1)^2
 2^{-\zeta_pj}\le C_9\Lambda^3L^{3/4}+C_{10}c_3L ,
\]
where $C_{10}$ depends only on $p$. Lemma~\ref{lem:x-caps-space} with $D$,
$q:=\lfloor3L\rfloor$, $r':=DN$, $\mathsf A$ and $\tau:=(2r)^{-10}$ gives a space
$\mathcal W_\varpi$ with
$\log(2\dim\mathcal W_\varpi)\le C_{11}\Lambda^3(1+L^{3/4})
+C_{\rm sk}C_8C_{10}c_3L$. Applied at $z(e)$, it gives $f\in\mathcal
W_\varpi$ with $|f(i)-P_i((t_{ih}/\rho_h)_h)|\le2\tau$ for all $i$.

\emph{Conclusion.} Let $\mathcal V$ be the sum of the constants and the
spaces $\mathcal W_\varpi$, $\varpi\in\mathcal P$. For $y\in\mathbb B$ with
pair $\varpi$, the function $rf$ above lies in $\mathcal V$, and
\[
 \|\kappa_y-rf\|_{L_2(\pi)}^2=r\sum_ic_i(|u_i^\top y|^p-f(i))^2
 \le3r\sum_ic_i(C_2^2r^{-20}+C_{\rm box}^2E_i^{2p}+4\tau^2)
 \le C_{12}(\frac{r^p}{L^{2p}}+r^{-18}),
\]
which is at most $2C_{12}(r/L^2)^p$ since $L\le r$. By
Lemma~\ref{lem:x-caps-frame} and Eq.~\eqref{lem:x-caps-eq-params},
\[
 \log\dim\mathcal V\le\log2+C_{13}(\Lambda+c_1(2+\log\frac1{c_1}
 )L)+C_{11}\Lambda^3(1+L^{3/4})+C_{\rm sk}C_8C_{10}c_3L .
\]
Choose $c_1$ with $C_{13}c_1(2+\log(1/c_1))\le1/32$ and then $c_3$ with
$C_{\rm sk}C_8C_{10}c_3\le1/32$. As $L\ge\sqrt r$, the remaining terms are at
most $L/32$ for $r$ large, so $\dim\mathcal V\le e^{L/8}$. As $p\downarrow1$
we have $\zeta_p\to0$. Then $C_4$ and $C_{10}$ grow, $c_3$ must shrink, and
the factor $r^{-(p-1)/(2p)}$ that absorbs the second term degenerates, so
the constants $C_p$ and $r_p$ given by this proof blow up.
\end{proof}

\subsection{A Variance-Sensitive Polynomial}

We use quantum procedures and code for a random variable as defined
before Lemma~\ref{lem:pre_ko}, and $C_{\rm KO}$ is the constant of that
lemma.

\begin{lemma}[A variance-sensitive polynomial]\label{lem:x-caps-poly}
There is an absolute constant $C_{\rm po}$ with the following property. Let
$\Omega$ be a finite set, $b\ge1$ an integer, $d:\Omega\to[-1,1]\cap
2^{-b}\mathbb Z$, $f:[-1,1]\to\R$, and $q\ge2$ an integer. For a probability
vector $\lambda$ on $\Omega$ define $\mu(\lambda):=\sum_\omega\lambda_\omega
d(\omega)$ and $s(\lambda)^2:=\sum_\omega\lambda_\omega d(\omega)^2$. For
every $\lambda$ there is a random variable $Z_\lambda\in[-1,1]$ with
\[
 \E[|Z_\lambda-\mu(\lambda)|^2]\le C_{\rm po}(s(\lambda)^2q^{-2}
 +q^{-4}),
\]
and a single polynomial $\Pi$ in variables $(v_\omega)_{\omega\in\Omega}$, of
total degree at most $C_{\rm po}q$, satisfies
$\E[f(Z_\lambda)]=\Pi((\lambda_\omega^{1/2})_{\omega\in\Omega})$ for every
probability vector $\lambda$.
\end{lemma}

\begin{proof}
\emph{The code.} Let $e_\star$ and $e_\omega$, $\omega\in\Omega$, be
orthonormal basis vectors of a register of qubits, with $e_\star$ the all-zero
basis vector, on whose remaining basis
vectors the operators below act as the identity. For a probability vector
$\lambda$ define $v:=\sum_\omega\lambda_\omega^{1/2}e_\omega$ and
$P_\lambda:=I-(e_\star-v)(e_\star-v)^\top$. As $\|e_\star-v\|_2^2=2$, the
operator $P_\lambda$ is the reflection in the hyperplane orthogonal to
$e_\star-v$. It is real, orthogonal and symmetric, so
$P_\lambda^{-1}=P_\lambda^*=P_\lambda$, and $P_\lambda e_\star=v$. Thus
$P_\lambda$ is a synthesizer for $\lambda$ with start vector $e_\star$ and no
garbage. On the product of this register with an output register of $2b+2$
bits, let $Y_1$ and $Y_2$ be the permutation unitaries
$e_\omega\otimes|z\rangle\mapsto e_\omega\otimes|z\oplus\operatorname{enc}
(d(\omega))\rangle$ and $e_\omega\otimes|z\rangle\mapsto e_\omega\otimes
|z\oplus\operatorname{enc}(d(\omega)^2)\rangle$, where $\operatorname{enc}$
is a fixed-point binary encoding, $\oplus$ is bitwise addition modulo two,
and the value $0$ is used at $e_\star$ and at the remaining basis vectors.
They do not depend on $\lambda$. With
$P_\lambda$ they form code for $Y:=d(\omega)$ and for $Y^2$, where $\omega$
has law $\lambda$. Every entry of $P_\lambda$ and of its controlled version
$\operatorname{diag}(I,P_\lambda)$ is a polynomial of degree at most two in
$v$.

\emph{The estimator.} Write $\mu:=\mu(\lambda)$ and $s:=s(\lambda)$. Let
$H:=\lfloor\log_2q\rfloor\ge1$ and $n_j:=2^{H-j}$ for $0\le j\le H$, so
$n_0>q/2$. For each $j<H$, run the procedure of Lemma~\ref{lem:pre_ko}
with $n:=n_j$, independently $120(j+1)$ times on the code for $Y$ and
$120(j+1)$ times on the code for $Y^2$. Let $m_j$ and $m'_j$ be the medians
of the outputs, clipped to $[-1,1]$ and $[0,1]$, respectively, and define
\[
 \mathcal I_j:=[m_j-\xi_j,\,m_j+\xi_j]\cap[-1,1],\qquad
 \xi_j:=(m_j'^{1/2}+1/n_j)/n_j ,
\]
and $\mathcal I_H:=[-1,1]$. Let $j_*$ be the least $j$ such that
$\bigcap_{l=j}^H\mathcal I_l\ne\emptyset$, and let $Z_\lambda$ be the
midpoint of the closed interval $\bigcap_{l=j_*}^H\mathcal I_l$.

\emph{Error.} Call a level $j<H$ good if $|m_j-\mu|\le s/n_j$ and
$|m'_j-s^2|\le s/n_j$, and call the level $H$ good. As $|Y|\le1$, the standard
deviations of $Y$ and $Y^2$ are at most $(\E[Y^2])^{1/2}=s$ and
$(\E[Y^4])^{1/2}\le s$. A median lies within $s/n_j$ of the mean unless at
most half of the runs do. By Hoeffding's
inequality (Lemma~\ref{lem:pre_hoeffding}), this has probability at most
$e^{-120(j+1)/18}\le2^{-8(j+1)-1}$ for each of the two medians. Clipping does
not increase the errors, as $\mu\in[-1,1]$ and $s^2\in[0,1]$. Hence level $j$
is bad with probability at most $2^{-8(j+1)}$. On a good level $j<H$ we have
$|m_j'^{1/2}-s|\le1/n_j$. For $s>0$ divide $|m'_j-s^2|\le s/n_j$ by
$m_j'^{1/2}+s\ge s$, and for $s=0$ note that $m'_j=0$. Hence $\mu\in\mathcal
I_j$ and $|\mathcal I_j|\le2\xi_j\le2s/n_j+4/n_j^2$. Both also hold for
$j=H$, where $n_H=1$. If all levels $l\ge j$ are good, then
$\mu\in\bigcap_{l\ge j}\mathcal I_l$, so $j_*\le j$, $Z_\lambda\in\mathcal
I_j$ and $|Z_\lambda-\mu|\le2s/n_j+4/n_j^2$. Let $j_{\rm b}:=0$ if all levels
are good, and otherwise one plus the largest bad index. Then
$\Pr[j_{\rm b}=j]\le2^{-8j}$, and with $n_j=2^{-j}n_0$,
\[
 \E[|Z_\lambda-\mu|^2]\le\sum_{j=0}^H2^{-8j}(\frac{2s2^j}{n_0}+
 \frac{4\cdot4^j}{n_0^2})^2\le\frac{9s^2}{n_0^2}+\frac{35}{n_0^4}
 \le\frac{36s^2}{q^2}+\frac{560}{q^4} .
\]

\emph{Polynomiality.} The estimator is a procedure in the model above: the
runs are independent, and $Z_\lambda$ is a function of their outcomes. In
every execution it uses $P_\lambda$, its inverse or a controlled version at
most $n_{\rm use}:=\sum_{j<H}240(j+1)C_{\rm KO}n_j\le960C_{\rm KO}q$ times, because
$\sum_j(j+1)2^{-j}\le4$. For a complete record $\varsigma$ of measurement
outcomes, let $\psi_\varsigma$ be the vector obtained from the initial state
by applying, in order, the operations prescribed along $\varsigma$, with each
measurement replaced by the orthogonal projection onto its observed outcome.
The probability of $\varsigma$ is $\|\psi_\varsigma\|_2^2$, the output is a
function $Z(\varsigma)$, and there are finitely many records. The operations
other than the uses of $P_\lambda$ do not depend on $\lambda$, and each use
multiplies by a matrix whose entries are polynomials of degree at most two in
$v$. Hence the coordinates of $\psi_\varsigma$ are polynomials of degree at
most $2n_{\rm use}$ in $v$ with complex coefficients, and $\|\psi_\varsigma\|_2^2$ is a
real polynomial of degree at most $4n_{\rm use}$ in $v$. Therefore
$\E[f(Z_\lambda)]=\sum_\varsigma f(Z(\varsigma))\|\psi_\varsigma\|_2^2$ is a
polynomial of degree at most $4n_{\rm use}$ in $(\lambda_\omega^{1/2})_\omega$, the
same for every $\lambda$. Take $C_{\rm po}:=\max\{560,3840C_{\rm KO}\}$.
\end{proof}

\begin{lemma}[Square-root features below one]\label{lem:x-caps-var}
Fix $0<p\le1$. There is $C_p$ such that for $C_p\log(2r)\le L\le r$ some
space $\mathcal V$ of functions on $I$ with $1\le\dim\mathcal V\le e^{L/8}$
satisfies, for every $y\in\mathbb B$,
\[
 \operatorname{dist}_{L_2(\pi)}(\kappa_y,\mathcal V)^2\le C_p(r^{2-p}
 K^{-2}+r^2K^{-2p-2}+r^{-18}).
\]
\end{lemma}

\begin{proof}
Constants $C_1,C_2,\ldots$ depend only on $p$. Let $c_1\in(0,1]$ depend only
on $p$, to be fixed below, define $k:=q:=\lfloor c_1K\rfloor$, and take
$C_p\ge8/c_1$, so that $k\ge\max\{2,c_1K/2\}$
by Eq.~\eqref{lem:x-caps-eq-params}. Let $\mathcal P$ be the set of
Lemma~\ref{lem:x-caps-frame} for this $k$, and define $R_1:=Rk^{1-1/p}$, so
that $\|e\|_1\le R_1$ for every residual. Let $b_0\ge1$ be an integer with
$2^{-b_0}\le r^{-10/p}$. For $\varpi\in\mathcal P$ let
$\tilde b_{il}:=\operatorname{sgn}(b_{il})2^{-b_0}\lfloor2^{b_0}|b_{il}|
\rfloor$, so that $|\tilde b_{il}|\le|b_{il}|$ and
$|\tilde b_{il}-b_{il}|\le r^{-10/p}$.

\emph{The polynomials.} Fix $\varpi\in\mathcal P$ with $N:=N_\varpi$, and let
$\Omega$ consist of the labels $(l,+)$ and $(l,-)$, $l\in[N]$, and a null
label $\varnothing$. For a residual $e$ let $\lambda(e)$ be the probability
vector with $\lambda_{(l,\pm)}:=(e_l)_\pm/R_1$ and
$\lambda_\varnothing:=1-\|e\|_1/R_1$. For $i\in I$ let
$d_i(l,+):=\tilde b_{il}$, $d_i(l,-):=-\tilde b_{il}$ and
$d_i(\varnothing):=0$, which lie in $[-1,1]\cap2^{-b_0}\mathbb Z$, and let
$f_i(z):=|\alpha_i+R_1z|^p$. Lemma~\ref{lem:x-caps-poly} for
$(\Omega,d_i,f_i,q)$ gives random variables $Z_{i,\lambda}$ and a polynomial
$\Pi_i$ of total degree at most $C_{\rm po}q$ with
$\E[f_i(Z_{i,\lambda})]=\Pi_i(\lambda^{1/2})$, where
$\lambda^{1/2}:=(\lambda_\omega^{1/2})_\omega$.

\emph{Error.} Let $y\in\mathbb B$ have the pair $\varpi$ and residual $e$,
define $\lambda:=\lambda(e)$, and let $\mu_i$ and $s_i$ be the numbers
$\mu(\lambda)$ and $s(\lambda)$ for the table $d_i$. Then
$R_1\mu_i=\sum_l\tilde b_{il}e_l$ and $s_i^2=R_1^{-1}\sum_l|e_l|\tilde
b_{il}^2$. Since $u_i^\top y=\alpha_i+\sum_lb_{il}e_l$, $|R_1\mu_i-\sum_l
b_{il}e_l|\le r^{-10/p}\|e\|_1\le R_1r^{-10/p}$ and $||u+w|^p-|u|^p|\le|w|^p$
for $p\le1$,
\[
 |\Pi_i(\lambda^{1/2})-|u_i^\top y|^p|
 \le R_1^p\,\E[|Z_{i,\lambda}-\mu_i|^p]+R_1^pr^{-10}
 \le C_{\rm po}^{p/2}R_1^p(\frac{s_i}q+\frac1{q^2})^p+R_1^pr^{-10},
\]
by Jensen's inequality and Lemma~\ref{lem:x-caps-poly}. By
$\sum_ic_i\tilde b_{il}^2\le\sum_ic_ib_{il}^2\le1$ we have
$\sum_ic_is_i^2\le R_1^{-1}\|e\|_1\le1$, so concavity gives
$\sum_ic_is_i^{2p}\le r^{1-p}$. Squaring and summing,
\[
 r\sum_ic_i(\Pi_i(\lambda^{1/2})-|u_i^\top y|^p)^2
 \le C_1R_1^{2p}(r^{2-p}q^{-2p}+r^2q^{-4p}+r^{-18}).
\]
As $R_1^{2p}=R^{2p}k^{2p-2}$, $k=q\ge c_1K/2$ and $2p-2\le0$, the right side
is at most $C_2(r^{2-p}K^{-2}+r^2K^{-2p-2}+r^{-18})$.

\emph{The space.} For fixed $\varpi$, the vector $(\Pi_i(v))_{i\in I}$ is a
combination of the coefficient vectors $(\text{coefficient of }\mathfrak m
\text{ in }\Pi_i)_{i\in I}$, one for each monomial $\mathfrak m$ of total
degree at most $C_{\rm po}q$ in the $2N+1$ variables. Let $\mathcal
W_\varpi$ be their span, and $\mathcal V$ the sum of the constants and the
spaces $\mathcal W_\varpi$. For $y$ as above, the function $i\mapsto
r\Pi_i(\lambda(e)^{1/2})$ lies in $\mathcal V$, which proves the distance
bound. Using $\log\binom{n'+m'}{m'}\le m'\log(e(n'+m')/m')$, $N\le2C_0r$ and
$C_{\rm po}\ge1$, we get
$\log\dim\mathcal W_\varpi\le C_{\rm po}q(\log(er/q)+\log(6C_0))$. By
Lemma~\ref{lem:x-caps-frame} and Eq.~\eqref{lem:x-caps-eq-params}, with $q\le
c_1K\le c_1L$,
\[
 \log\dim\mathcal V\le\log2+C_3(\Lambda+c_1(2+\log\frac1{c_1})
 L)+C_{\rm po}c_1(2+\log\frac1{c_1}+\log(6C_0))L .
\]
Choose $c_1$ so that both coefficients of $L$ are at most $1/24$, and then
$C_p$ so large that $\log2+C_3\Lambda\le L/24$. Then
$\dim\mathcal V\le e^{L/8}$.
\end{proof}

\subsection{Proof of the Crossover Bound}

In each part, $C_1,C_2,\ldots$ denote constants depending only on $p$, and
the constant $C_p$ in the hypothesis $L\ge C_p\log(2r)$ is taken large in
terms of them and of the constants of the lemmas used.

\begin{proof}[Proof of Lemma~\ref{lem:joint-positive-caps}(a)]
Let $c_1,c_2\in(0,1]$ depend only on $p$, to be fixed below, and define
$k:=\lfloor c_1K\rfloor$, $q:=\lfloor c_2L\rfloor$ and $\rho:=Rk^{1/2-1/p}$.
By Eq.~\eqref{lem:x-caps-eq-params}, $k\ge c_1K/2\ge1$ and $q\ge c_2L/2\ge q_0$,
and $k,q\le r$. Let $\mathcal P$ be the set of
Lemma~\ref{lem:x-caps-frame} for this $k$, and fix $\varpi\in\mathcal P$
with $N:=N_\varpi$. For $i\in I$ let $P_i$ be the polynomial of
Lemma~\ref{lem:x-caps-box} with $D:=2$, $\alpha:=\alpha_i$, signs $(1,-1)$,
sides $\rho_1=\rho_2:=\rho$ and degrees $q_1=q_2:=q$. By
Lemma~\ref{lem:x-caps-coef}, $R+2\rho\le3R$ and $R\ge1$, its monomial
coefficients have absolute sum at most $\mathsf A:=4(q+1)^264^q\cdot27R^2$.
Let $w_{i,1}\in\R^{2N}$ consist of the coordinates $(b_{il})_+$, $l\in[N]$,
followed by the coordinates $(b_{il})_-$, and let $w_{i,2}$ consist of the
$(b_{il})_-$ followed by the $(b_{il})_+$. Then $\|w_{i,h}\|_2=
\|b_{i\cdot}\|_2\le1$. Lemma~\ref{lem:x-caps-space} with $D=2$, $q$,
$r':=2N$, $\mathsf A$ and $\tau:=(2r)^{-4}$ gives a space $\mathcal
W_\varpi$ with $\log(2\dim\mathcal W_\varpi)\le C_1(q+\Lambda)$.

Let $y\in\mathbb B$ have the pair $\varpi$ and residual $e$, and let
$z\in\R^{2N}$ consist of the coordinates $(e_l)_+/\rho$ followed by the
$(e_l)_-/\rho$. Then $\|z\|_2=\|e\|_2/\rho\le1$. Moreover
$a_i^+:=\rho\langle w_{i,1},z\rangle=\sum_l(b_{il}e_l)_+$ and
$a_i^-:=\rho\langle w_{i,2},z\rangle=\sum_l(b_{il}e_l)_-$, so
$u_i^\top y=\alpha_i+a_i^+-a_i^-$, and
$t_i:=a_i^++a_i^-=\sum_l|b_{il}||e_l|\le\|b_{i\cdot}\|_2\|e\|_2\le\rho$.
Lemma~\ref{lem:x-caps-box} at $(a_i^+,a_i^-)\in[0,\rho]^2$, with
$(\rho a_i^+)^{1/2}+(\rho a_i^-)^{1/2}\le2(\rho t_i)^{1/2}$, and
Lemma~\ref{lem:x-caps-space} at $z$ give $f\in\mathcal W_\varpi$ with
\[
 (|u_i^\top y|^p-f(i))^2\le C_2((\rho t_i)^pq^{-2p}
 +\rho^{2p}q^{-4p}+\tau^2)\qquad(i\in I).
\]
By H\"older's inequality with exponents $2/p$ and $2/(2-p)$, and then
Lemma~\ref{lem:pre_minkowski}(a) in $L_2$ of the weights $c_i$
together with Eq.~\eqref{lem:x-caps-eq-frame},
\[
 \sum_ic_i(\rho t_i)^p\le\rho^pr^{1-p/2}(\sum_ic_it_i^2)^{p/2}
 \le\rho^pr^{1-p/2}(\sum_l|e_l|(\sum_ic_ib_{il}^2)^{1/2})^p
 \le\rho^pr^{1-p/2}\|e\|_1^p .
\]
By Lemma~\ref{lem:x-caps-frame}, $\rho^p\|e\|_1^p\le C_3k^{p/2-1}r^{p-1}$
if $p\ge1$ and $\rho^p\|e\|_1^p\le C_3k^{3p/2-2}$ if $p<1$. Let $\mathcal V$
be the sum of the constants and the spaces $\mathcal W_\varpi$. Then
$rf\in\mathcal V$ and
\[
 \|\kappa_y-rf\|_{L_2(\pi)}^2=r\sum_ic_i(|u_i^\top y|^p-f(i))^2
 \le C_4(r^{2-p/2}\rho^p\|e\|_1^pq^{-2p}+r^2\rho^{2p}q^{-4p}+r^2\tau^2
 ).
\]
Here $r^2\tau^2\le r^{-6}\le r^{-3p}\le r^2\rho^{2p}q^{-4p}$, since
$\rho^{2p}=R^{2p}k^{p-2}\ge r^{p-2}$ and $q\le r$. Hence the hypothesis
of Lemma~\ref{lem:s-cal} holds with $\beta=1$, $B_0=3$ and
\[
 \delta^2:=2C_4(\Xi_p+r^2k^{p-2}q^{-4p})\ge r^{-6},\qquad
 \Xi_p:=\begin{cases}r^{2-p/2}k^{3p/2-2}q^{-2p},&0<p<1,\\
 r^{1+p/2}k^{p/2-1}q^{-2p},&1\le p<2,\end{cases}
\]
where we absorbed $C_3R^{2p}$ into $C_4$. By
Lemma~\ref{lem:x-caps-frame} and Eq.~\eqref{lem:x-caps-eq-params},
\[
 \log\dim\mathcal V\le\log2+C_5(\Lambda+c_1(2+\log\frac1{c_1})
 L)+C_1(c_2L+\Lambda),
\]
which is at most $L/8$ once $c_1,c_2$ are small and $C_p$ is large.
Lemma~\ref{lem:s-cal}(a) gives a reweighting of at most
$C_6\epsilon^{-2}r\delta^2\Lambda$ original rows with relative error
$\epsilon$ for every query. All exponents of $k$ and $q$ in $\delta^2$ are
negative, so $k\ge c_1K/2$ and $q\ge c_2L/2$ give
$r\delta^2\le C_7(A_p(r,K)L^{-2p}+r^3K^{p-2}L^{-4p})$.
\end{proof}

\begin{proof}[Proof of Lemma~\ref{lem:joint-positive-caps}(b)]
Let $c_1,c_2\in(0,1]$ depend only on $p$, to be fixed below, and define
$k:=\lfloor c_1K\rfloor$ and $q:=\lfloor c_2Q\rfloor$, where
$Q:=L/\Lambda^3$. As $L\ge C_p\Lambda^3$, we have $q\ge c_2Q/2\ge q_0$ and
$k\ge c_1K/2\ge1$. Let $\mathcal P$ be the set of
Lemma~\ref{lem:x-caps-frame} for this $k$. Fix $\varpi\in\mathcal P$, take
the bands of Lemma~\ref{lem:x-caps-bins}, and define $D:=|\mathcal H|\le
C_1\Lambda^2$. For $i\in I$ let $P_i$ be the polynomial of
Lemma~\ref{lem:x-caps-box} with this $D$, $\alpha:=\alpha_i$, the signs
$\sigma_h$, the sides $\rho_h$ and all degrees equal to $q$.

\emph{A common space by counting monomials.} Let $\Phi(e)\in\R^{|\mathcal
G|N}$ be the concatenation of the vectors $f_g(e)$, $g\in\mathcal G$. Each
$t_{ih}$ is a linear form in $\Phi(e)$ whose coefficients depend only on $i$
and $h$. Hence $P_i((t_{ih}/\rho_h)_h)$ is a polynomial of total degree at
most $Dq$ in $\Phi(e)$, with coefficients depending only on $i$. Let
$\mathcal W_\varpi$ be the span of the coefficient vectors, indexed by
$i\in I$, of the monomials of total degree at most $Dq$ in $|\mathcal G|N$
variables. As $|\mathcal G|N\le C_2\Lambda r$,
\[
 \log\dim\mathcal W_\varpi\le\log\binom{|\mathcal G|N+Dq}{Dq}\le
 Dq\log(e(1+|\mathcal G|N))\le C_3q\Lambda^3\le C_3c_2L .
\]

\emph{Error.} Let $y\in\mathbb B$ have the pair $\varpi$ and residual $e$,
and let $f(i):=P_i((t_{ih}/\rho_h)_h)$, a function in $\mathcal W_\varpi$.
By Lemmas~\ref{lem:x-caps-box} and~\ref{lem:x-caps-bins},
$||u_i^\top y|^p-f(i)|\le C_4r^{-10}+C_{\rm box}E_i^p$ with
$E_i:=\sum_h((\rho_ht_{ih})^{1/2}q^{-1}+\rho_hq^{-2})$. As $E_i$ is a sum of
$2D$ terms, the power-mean inequality gives
$E_i^{2p}\le(2D)^{(2p-1)_+}\sum_h((\rho_ht_{ih})^pq^{-2p}+\rho_h^{2p}
q^{-4p})$. Summing with the weights $c_i$ and using
Eq.~\eqref{lem:x-caps-eq-overlap}, $\rho_h\le\rho$ and $\sum_ic_i=r$,
\[
 \sum_ic_iE_i^{2p}\le(2D)^{1+(2p-1)_+}(W_*q^{-2p}+r\rho^{2p}q^{-4p}),
 \qquad W_*:=\begin{cases}C_5r^{1-p}k^{2p-2},&0<p<1,\\
 C_5r^{p-1},&1\le p<2.\end{cases}
\]
Let $\mathcal V$ be the sum of the constants and the spaces $\mathcal
W_\varpi$. Then $rf\in\mathcal V$ and, as $1+(2p-1)_+\le2\max\{1,p\}$ and
$\rho^{2p}=R^{2p}k^{p-2}$,
\[
 \|\kappa_y-rf\|_{L_2(\pi)}^2\le3r\sum_ic_i(C_4^2r^{-20}+C_{\rm box}^2
 E_i^{2p})\le C_6D^{2\max\{1,p\}}(rW_*q^{-2p}+r^2k^{p-2}q^{-4p}
 )+C_6r^{-18}.
\]
As in part~(a), $r^{-18}\le r^2k^{p-2}q^{-4p}$, so the hypothesis of
Lemma~\ref{lem:s-cal} holds with $\beta=1$, $B_0=3$ and
$\delta^2:=2C_6D^{2\max\{1,p\}}(rW_*q^{-2p}+r^2k^{p-2}q^{-4p})\ge r^{-6}$.
By Lemma~\ref{lem:x-caps-frame} and Eq.~\eqref{lem:x-caps-eq-params},
$\log\dim\mathcal V\le\log2+C_7(\Lambda+c_1(2+\log(1/c_1))L)+C_3c_2L\le
L/8$ once $c_1,c_2$ are small and $C_p$ is large. Lemma~\ref{lem:s-cal}(a)
and $D^{2\max\{1,p\}}\le C_1^4\Lambda^{4\max\{1,p\}}$ give a reweighting with
at most
\[
 C_8\epsilon^{-2}\Lambda^{1+4\max\{1,p\}}(r^2W_*q^{-2p}+r^3k^{p-2}
 q^{-4p})
\]
original rows. Now $q\ge c_2Q/2$ and $k\ge c_1K/2$, all exponents of $k$
and $q$ are negative, $r^2W_*=C_5B_p(r,K)$ for $p\ge1$, and
$r^2W_*\le C_9r^{3-p}K^{2p-2}=C_9B_p(r,K)$ for $p<1$. This proves the
first claim of~(b).

\emph{The crossover.} Let $I\subseteq[\iota_1,\iota_2]$ with
$0<\iota_1\le\iota_2$, and define $\varkappa:=2/(p+2)$, so that
$\iota_1r^\varkappa\le L\le\iota_2r^\varkappa$. If $L<C_p\Lambda^3$ or
$L>r$, then $r$ is bounded in terms of $p$ and $I$, and
Theorem~\ref{thm:general_lp_field} suffices. Otherwise the first claim
applies, with $K\ge L/(2\Lambda)\ge\iota_1r^\varkappa/(2\Lambda)$, $K\le L$
and $Q\ge\iota_1r^\varkappa/\Lambda^3$. Up to factors
$C_{p,I}\Lambda^{O_p(1)}$ the first term is $r^{p+1-2p\varkappa}$ for
$p\ge1$ and $r^{3-p+(2p-2)\varkappa-2p\varkappa}=r^{3-p-2\varkappa}$ for
$p<1$, and the second term is $r^{3+(p-2)\varkappa-4p\varkappa}=
r^{(2-3p)/(p+2)}$. Now
\[
 p+1-\frac{4p}{p+2}=\frac{p^2-p+2}{p+2}=1-\frac{p(2-p)}{p+2},\qquad
 3-p-\frac4{p+2}=\frac{2+p-p^2}{p+2}=1-\frac{p^2}{p+2},
\]
which is $1-\delta_p$ in each range. Finally
$(2-3p)/(p+2)\le1-\delta_p$, because $2-3p\le p^2-p+2$ and
$2-3p\le2+p-p^2$ for $0<p<2$.
\end{proof}

\begin{proof}[Proof of Lemma~\ref{lem:joint-positive-caps}(c)]
Let $1\le p<2$. Since $L\le r$, we have $T\le\sqrt r$. Let $T_1\ge1$ depend
only on $p$, fixed below. If $T\le T_1$, then
$\min\{1,T^{-2p}\log(2+T)\}\ge T_1^{-2p}/2$, and
Theorem~\ref{thm:general_lp_field} gives the claim. Let $T>T_1$.

\emph{The case $1<p<2$.} Here $L=T\sqrt r\ge\sqrt r$. Let $r_p$ be as in
Lemma~\ref{lem:x-caps-adaptive}, enlarged so that $\sqrt r\ge C\log(2r)$ for
$r\ge r_p$, with $C$ the constant of Lemma~\ref{lem:x-caps-cal} for $B_0=1$.
If $r<r_p$, then $T^{-2p}\log(2+T)\ge r_p^{-p}\log2$, and
Theorem~\ref{thm:general_lp_field} suffices. If $r\ge r_p$,
Lemma~\ref{lem:x-caps-adaptive} gives a space $\mathcal V$ with
$\dim\mathcal V\le e^{L/8}$ and distance at most $C_1(r/L^2)^{p/2}=
C_1T^{-p}$ on $\mathbb B$. Define $\delta:=\max\{C_1,1\}T^{-p}$ and choose
$T_1$ with $\max\{C_1,1\}T_1^{-p}\le1$. Then $\delta\le1$, and
$\delta\ge r^{-p/2}\ge r^{-1}$ because $T\le\sqrt r$.
Lemma~\ref{lem:x-caps-cal} with $B_0=1$ gives a reweighting with at most
$C_2\epsilon^{-2}r\delta^2\log(2+1/\delta)\le C_3\epsilon^{-2}rT^{-2p}
\log(2+T)$ original rows, since $2+T^p\le(2+T)^2$. The constants of
Lemma~\ref{lem:x-caps-adaptive} blow up as $p\downarrow1$, so this argument
gives no bound at $p=1$ that is uniform in $p$. For each fixed $p>1$ it is
complete.

\emph{The case $p=1$.} Here $T:=K/\sqrt r$. Lemma~\ref{lem:x-caps-var} with
$p=1$ gives a space $\mathcal V$ with $\dim\mathcal V\le e^{L/8}$ and
squared distance at most $C_4(rK^{-2}+r^2K^{-4}+r^{-18})=C_4(T^{-2}+T^{-4}
+r^{-18})\le3C_4T^{-2}$ on $\mathbb B$, using $T\ge1$ and $T^{-2}\ge r^{-1}$.
Define $\delta:=(3C_4)^{1/2}T^{-1}$ and choose $T_1\ge(3C_4)^{1/2}$. Then
$r^{-1}\le\delta\le1$, and Lemma~\ref{lem:x-caps-cal} with $B_0=1$ gives at
most $C_5\epsilon^{-2}rT^{-2}\log(2+T)$ rows.

In both cases $T^{-2p}\log(2+T)\to0$ as $T\to\infty$, which gives the last
assertion of~(c).
\end{proof}

\begin{proof}[Proof of Lemma~\ref{lem:joint-positive-caps}(d)]
Let $0<p<1$. As $K\le L\le r$, we have $T_0\le r^{p/2}$. Let $T_1\ge1$ depend
only on $p$. If $T_0\le T_1$, Theorem~\ref{thm:general_lp_field} gives the
claim as in part~(c). Let $T_0>T_1$. Then $K\ge T_0r^{1-p/2}\ge\sqrt r$, so
$r^2K^{-2p-2}=T_0^{-2}r^pK^{-2p}\le T_0^{-2}$, and $r^{-18}\le r^{-p}\le
T_0^{-2}$. Lemma~\ref{lem:x-caps-var} gives a space $\mathcal V$ with
$\dim\mathcal V\le e^{L/8}$ and squared distance at most
$C_1(r^{2-p}K^{-2}+r^2K^{-2p-2}+r^{-18})\le3C_1T_0^{-2}$ on $\mathbb B$.
Define $\delta:=(3C_1)^{1/2}T_0^{-1}$ and choose $T_1\ge(3C_1)^{1/2}$. Then
$r^{-1}\le\delta\le1$, and Lemma~\ref{lem:x-caps-cal} with $B_0=1$ gives at
most $C_2\epsilon^{-2}rT_0^{-2}\log(2+T_0)$ original rows. This proves the
bound in~(d).

Let $\iota_1r^{2/(p+2)}\le L\le\iota_2r^{2/(p+2)}$ with fixed
$0<\iota_1\le\iota_2$. Then $K\ge L/(2\Lambda)$, and since
$2/(p+2)-1+p/2=p^2/(2(p+2))$,
\[
 T_0\ge\frac{\iota_1}{2\Lambda}\,r^{p^2/(2(p+2))} .
\]
If $T_0\ge1$, then $T_0^{-2}\log(2+T_0)\le4\iota_1^{-2}\Lambda^2
r^{-p^2/(p+2)}\cdot2\Lambda$, because $T_0\le r$. If $T_0<1$, then $r$ is
bounded in terms of $p$ and $\iota_1$. In both cases
$S_p(r,\epsilon)\le C_{p,\iota_1,\iota_2}\epsilon^{-2}r^{1-p^2/(p+2)}
\log^3(2r)$.
\end{proof}

The shared-coordinate lower bound at this crossover has order
$\epsilon^{-2}r^{1/(p+2)}$. The upper and lower bounds still differ
in their rank dependence. Lemma~\ref{lem:joint-adaptive-shifts}
tracks the gain from adaptive shifts uniformly from the crossover
to exponential accuracy.

\subsection{Regular Affine Components}

\begin{lemma}[Matching bounds for regular affine components]
\label{lem:joint-regular-components}
Fix $1\le p<2$, $\delta>0$ and $0<c\le1/16$, and let $B\ge B_*$, where
$B_*$ is an absolute constant.
Consider rank-$r$ inputs with energy
\[
 F(x)=\sum_{g=1}^G c_g\,
       \E[\theta_g(Z_g)|\langle m_g+Z_g,x\rangle|^p],
 \qquad G\le e^{cL},\quad c_g>0,\quad B^{-1}\le\theta_g\le B .
\]
Each reference $Z_g$ is a finitely supported centered random
vector with covariance $Q_g$, and the rows of the input are the vectors
$(c_g\Pr[Z_g=\zeta]\theta_g(\zeta))^{1/p}(m_g+\zeta)$, one for each $g$
and each atom $\zeta$ of $Z_g$.
Assume
$\E[e^{\langle v,Z_g\rangle}]\le e^{B^2 v^\top Q_gv/2}$ for every $v$.
For every $v$ with $v^\top Q_gv>0$, $a\in\mathbb R$, and
$L^{-1/2}\le z\le1$, assume
\[
 \Pr[0<|a+\langle v,Z_g\rangle|
                  \le z\sqrt{v^\top Q_gv}]\le Bz .
\]
For $r^{2/(p+2)+\delta}\le L\le r$, the worst-case original-row
support in this class is
\[
 \Theta_{p,B,\delta,c}(
          \frac{\epsilon^{-2}r^2\log r}{L^{p+3/2}}).
\]
\end{lemma}

The lemma makes two assertions. Every input of the class has a
sparsifier in the sense of Eq.~\eqref{eq:row_sparsification} supported on
$O_{p,B,\delta,c}(\epsilon^{-2}r^2\log r/L^{p+3/2})$ of its rows. Some
input of the class admits no reweighting of its rows by real weights
with fewer nonzero weights that preserves every energy within relative
error $\epsilon$. The constant $B_*$ only ensures that the hard inputs
of the lower bound lie in the class. The proof of the upper bound works
for every $B\ge1$ and uses only $C_{p,B}\log^3(2r)\le L\le r$.

The upper bound applies calibrated sampling, Lemma~\ref{lem:s-cal},
with one approximation space for all queries at once. The space is a
direct sum over the components. On one component the row is random and
the query direction is fixed, and we approximate the cosine integral of
$|t|^p$ by polynomials in the eigenspaces of the covariance $Q_g$, with
degrees adapted to the variance of the query direction. The interval
condition controls the error of the frequency truncation near the zeros
of the score. The base law of the sampling mixes the Lewis law with a
reference law that weights each component by $\tr[Q_g^{p/2}]$, after a
change of coordinates that makes these traces sum to the transverse
dimension. The lower bound is the adaptive-shift input of
Lemma~\ref{lem:s-shift}, which lies in the class by
Lemma~\ref{lem:s-regular}.

\subsection{Four Auxiliary Lemmas}

The first lemma gives the change of coordinates.

\begin{lemma}[Normalized covariances]\label{lem:x-rc-normal}
Let $0<p<2$, let $W$ be a Euclidean space of dimension $d_W\ge1$, let
$c_1,\ldots,c_G>0$, and let $Q_1,\ldots,Q_G$ be positive semidefinite
operators on $W$ whose sum is positive definite. Then there is an
invertible linear map $T$ of $W$ with
\begin{equation}\label{lem:x-rc-eq-trace}
 \sum_{g=1}^Gc_g\,(TQ_gT^\top)^{p/2}=I_W,
 \qquad\text{hence}\qquad
 \sum_{g=1}^Gc_g\tr[(TQ_gT^\top)^{p/2}]=d_W .
\end{equation}
\end{lemma}

\begin{proof}
Fix an orthonormal basis of $W$, so that operators are matrices. Let
$\mathcal S$ be the closed set of matrices with determinant one and define
$\Phi(T):=\sum_gc_g\tr[(TQ_gT^\top)^{p/2}]$, a continuous function. For
positive semidefinite $X$ the subadditivity of $t\mapsto t^{p/2}$ on the
eigenvalues gives $\tr[X^{p/2}]\ge(\tr[X])^{p/2}$. Define
$c_{\min}:=\min_gc_g$ and let $\mu>0$ be the least eigenvalue of
$M:=\sum_gc_gQ_g$. Applying subadditivity once more,
\[
 \Phi(T)\ge\sum_gc_g^{1-p/2}(c_g\tr[TQ_gT^\top])^{p/2}
 \ge c_{\min}^{1-p/2}(\tr[TMT^\top])^{p/2}
 \ge c_{\min}^{1-p/2}(\mu\|T\|_{\rm HS}^2)^{p/2}.
\]
Hence the sublevel sets of $\Phi$ on $\mathcal S$ are compact, and
$\Phi$ attains its minimum on $\mathcal S$ at some $T$.

Let $H$ be symmetric with $\tr[H]=0$, so that $e^{tH}T\in\mathcal S$. If
$Q_g\ne0$, write $Q_g=A_gA_g^\top$ with $A_g$ of full column rank $k_g$.
The nonzero eigenvalues of $e^{tH}TQ_gT^\top e^{tH}=Y_tY_t^\top$,
$Y_t:=e^{tH}TA_g$, are those of the positive definite $k_g\times k_g$
matrix $X_t:=Y_t^\top Y_t=A_g^\top T^\top e^{2tH}TA_g$. Thus the
summand of $g$ in $\Phi(e^{tH}T)$ is $c_g\tr[X_t^{p/2}]$, a smooth
function of $t$. For a smooth path $X_t$ of positive definite matrices
and $f(\xi):=\xi^{p/2}$, the Cauchy integral
$f(X)=(2\pi\mathrm i)^{-1}\oint f(\xi)(\xi-X)^{-1}\,\d\xi$
over a contour in the right half plane around the spectrum gives
$\frac{\d}{\d t}\tr[f(X_t)]=\tr[f'(X_t)\dot X_t]$, after one integration by
parts in $\xi$. With $\dot X_0=2Y_0^\top HY_0$, where $Y_0:=TA_g$, and the identity
$Y(Y^\top Y)^{p/2-1}Y^\top=(YY^\top)^{p/2}$ for $Y$ of full column rank,
which follows from a singular value decomposition of $Y$, the derivative
of the summand of $g$ at $t=0$ is
$c_gp\tr[Y_0(Y_0^\top Y_0)^{p/2-1}Y_0^\top H]=c_gp\tr[H(TQ_gT^\top)^{p/2}]$.
Hence
\[
 \frac{\d}{\d t}\Phi(e^{tH}T)|_{t=0}=p\tr[HN],\qquad
 N:=\sum_gc_g(TQ_gT^\top)^{p/2}.
\]
This derivative vanishes at the minimizer. Taking
$H:=N-(\tr[N]/d_W)I_W$ gives $0=\tr[HN]=\tr[H^2]$, so
$N=\lambda I_W$ with $\lambda:=\tr[N]/d_W$. Here $\lambda>0$, since
otherwise every $TQ_gT^\top$ vanishes and $M=0$. Replacing $T$ by
$\lambda^{-1/p}T$ multiplies $N$ by $\lambda^{-1}$
and gives Eq.~\eqref{lem:x-rc-eq-trace}. The second identity is the trace of
the first.
\end{proof}

The second lemma collects the scalar facts used for one component score.

\begin{lemma}[Moments of a component score]\label{lem:x-rc-moments}
Let $B\ge1$ and let $Y$ be a finitely supported centered random
variable with variance $\sigma^2$ and $\E[e^{\lambda Y}]\le
e^{B^2\sigma^2\lambda^2/2}$ for all real $\lambda$.
\begin{enumerate}[label=(\alph*),nosep,leftmargin=*]
\item $\E[Y^4]\le16B^4\sigma^4$, and $\E[|Y|^m]\le768(B\sigma)^m$ for
$0<m\le8$.
\item For $1\le p<2$ and $a\in\R$,
$\E[|a+Y|^p]\ge(8\sqrt2B^2)^{-1}(a^2+\sigma^2)^{p/2}$.
\item Let $L\ge1$ and $a\in\R$, let $\sigma>0$, and suppose that
$\Pr[0<|a+Y|\le z\sigma]\le Bz$ for $L^{-1/2}\le z\le1$. Then for every
$T>0$,
\[
 \E[(1+T|a+Y|)^{-2}\mathbf 1_{a+Y\ne0}]
 \le B(L^{-1/2}+\min\{1,(T\sigma)^{-1}\}).
\]
\end{enumerate}
\end{lemma}

\begin{proof}
(a) If $\sigma=0$ then $Y=0$ almost surely. Otherwise define $s:=B\sigma$.
Markov's inequality (Lemma~\ref{lem:pre_markov}) for $e^{\lambda Y}$ with $\lambda=t/s^2$ gives
$\Pr[Y\ge t]\le e^{-t^2/(2s^2)}$, and similarly for $-Y$. Hence for
$m>0$
\[
 \E[|Y|^m]=\int_0^\infty mt^{m-1}\Pr[|Y|>t]\,\d t
 \le2m\int_0^\infty t^{m-1}e^{-t^2/(2s^2)}\,\d t
 =m2^{m/2}\Gamma(m/2)s^m .
\]
For $m=4$ and $m=8$ the factor is $16$ and $768$. For $0<m\le8$,
Lyapunov's inequality gives $\E[|Y|^m]\le(\E[Y^8])^{m/8}\le768s^m$.

(b) Define $U:=a+Y$, so $\E[U^2]=a^2+\sigma^2$, and by~(a)
$\E[U^4]\le8a^4+8\E[Y^4]\le128B^4(\E[U^2])^2$. With
$\vartheta:=2/(4-p)\in(0,1]$ we have $2=\vartheta p+4(1-\vartheta)$, and
H\"older's inequality gives
$\E[U^2]\le(\E[|U|^p])^\vartheta(\E[U^4])^{1-\vartheta}$. Therefore
\[
 \E[|U|^p]\ge(\E[U^2])^{(4-p)/2}(\E[U^4])^{-(2-p)/2}
 \ge(128B^4)^{-(2-p)/2}(\E[U^2])^{p/2},
\]
and $(128B^4)^{-(2-p)/2}\ge(128B^4)^{-1/2}=(8\sqrt2B^2)^{-1}$.

(c) Define $\mu(\varrho):=\Pr[0<|a+Y|\le\varrho]$. For
$\varrho\le L^{-1/2}\sigma$ monotonicity and the hypothesis at
$z=L^{-1/2}$ give $\mu(\varrho)\le BL^{-1/2}$. For
$L^{-1/2}\sigma<\varrho\le\sigma$ the hypothesis at $z=\varrho/\sigma$
gives $\mu(\varrho)\le B\varrho/\sigma$, and for $\varrho>\sigma$ we have
$\mu(\varrho)\le1\le B$. Hence
$\mu(\varrho)\le BL^{-1/2}+B\min\{1,\varrho/\sigma\}$ for all
$\varrho>0$. Since
$(1+T|y|)^{-2}=\int_{|y|}^\infty2T(1+T\varrho)^{-3}\,\d\varrho$ for
$y\ne0$, Tonelli's theorem gives
\[
 \E[(1+T|a+Y|)^{-2}\mathbf 1_{a+Y\ne0}]
 =\int_0^\infty\frac{2T\mu(\varrho)}{(1+T\varrho)^3}\,\d\varrho .
\]
Now $\int_0^\infty2T(1+T\varrho)^{-3}\,\d\varrho=1$ and
$\int_0^\infty2T\varrho(1+T\varrho)^{-3}\,\d\varrho=1/T$, so the right
side is at most $BL^{-1/2}+B\min\{1,1/(T\sigma)\}$.
\end{proof}

The third lemma truncates the cosine integral Eq.~\eqref{lem:x-subexp-eq-cosine}
at a frequency $T$, with a correction that keeps the value zero at the
origin.

\begin{lemma}[A truncated cosine integral]\label{lem:x-rc-cosine}
Let $1\le p<2$ and $T>0$, let $\varpi_p$ be the constant
of Eq.~\eqref{lem:x-subexp-eq-cosine}, and let $\kappa_T$ be the positive
measure on $(0,T]$ with density
$\varpi_p(t^{-1-p}+p^{-1}T^{-1-p})$. Define
$g_T(y):=\int(1-\cos(ty))\,\d\kappa_T(t)$ for $y\in\R$. Then $g_T(0)=0$,
\[
 ||y|^p-g_T(y)|\le\frac{6\varpi_pT^{-p}}{1+T|y|}
 \quad(y\in\R),\qquad
 \int(t/T)^3\,\d\kappa_T(t)\le2\varpi_pT^{-p}.
\]
\end{lemma}

\begin{proof}
The integral defining $g_T$ converges, as $1-\cos(ty)\le t^2y^2/2$ and
$p<2$. Clearly $g_T(0)=0$. Define $\operatorname{sinc}(w):=\sin(w)/w$ and
$\operatorname{sinc}(0):=1$, so that
$T^{-1}\int_0^T\cos(ty)\,\d t=\operatorname{sinc}(Ty)$. By Eq.~\eqref{lem:x-subexp-eq-cosine}
and $\int_T^\infty t^{-1-p}\,\d t=T^{-p}/p$,
\[
 |y|^p-g_T(y)=-\varpi_p\int_T^\infty\cos(ty)\,t^{-1-p}\,\d t
   +\frac{\varpi_p}{pT^p}\operatorname{sinc}(Ty).
\]
The integral is at most $T^{-p}/p\le T^{-p}$ in absolute value. For
$y\ne0$, one integration by parts writes it as
$-\sin(Ty)T^{-1-p}/y+((1+p)/y)\int_T^\infty\sin(ty)t^{-2-p}\,\d t$,
which is at most $2T^{-1-p}/|y|$ in absolute value. Also
$|\operatorname{sinc}(w)|\le\min\{1,1/|w|\}$. Hence
$||y|^p-g_T(y)|\le3\varpi_pT^{-p}\min\{1,1/(T|y|)\}$, and
$\min\{1,1/w\}\le2/(1+w)$. Finally
$\int(t/T)^3\d\kappa_T(t)=\varpi_pT^{-p}(1/(3-p)+1/(4p))\le2\varpi_p
T^{-p}$.
\end{proof}

The fourth lemma approximates a product of phases of block scores by
Taylor polynomials. For an integer $m\ge1$ define
$\varphi_m(w):=\sum_{l=0}^{m-1}(\mathrm iw)^l/l!$. For real $w$ the
integral form of the Taylor remainder gives
$|e^{\mathrm iw}-\varphi_m(w)|\le|w|^m/m!$.

\begin{lemma}[Products of phases]\label{lem:x-rc-phase}
Let $B\ge1$ and $r\ge2$, and define $C_B:=2e^{4+B^2/2}$. Let $Z$ be a
finitely supported centered random vector in a Euclidean space $W$ with
covariance $Q$ and $\E[e^{\langle v,Z\rangle}]\le e^{B^2v^\top Qv/2}$ for
$v\in W$. Let $E_j$, $j\in\mathcal J$, be mutually orthogonal subspaces
of $W$ invariant under $Q$, with $|\mathcal J|\le4\log(2r)$ and
orthogonal projections $P_j$. Let $u\in W$, and define
$Z_j:=\langle u,P_jZ\rangle$ and $\sigma_j^2:=\E[Z_j^2]$. Let $T>0$ and
$\rho_\circ\ge0$, let $\mathcal G$ be an event and $w_\circ$ a random
variable with $|w_\circ|\le\rho_\circ$ on $\mathcal G$. Let
$k_\circ\ge C_B(T\rho_\circ+\log(2r))$ and
$k_j\ge C_B(T^2\sigma_j^2+\log(2r))$, $j\in\mathcal J$, be integers.
Then for $0<t\le T$,
\[
 \|\mathbf 1_{\mathcal G}(e^{\mathrm itw_\circ}
   \prod_{j\in\mathcal J}e^{\mathrm itZ_j}
   -\varphi_{k_\circ}(tw_\circ)\prod_{j\in\mathcal J}
   \varphi_{k_j}(tZ_j))\|_{L_2}\le r^{-10}(t/T)^3 .
\]
\end{lemma}

\begin{proof}
If $\sigma_j=0$, then $Z_j=0$ almost surely and the two $j$-th factors
agree, so we may assume $\sigma_j>0$ for every $j$. Define
$Y_j:=Z_j/\sigma_j$. For $\lambda\in\R^{\mathcal J}$ let
$v:=\sum_j\lambda_jP_ju/\sigma_j$, so that $\langle v,Z\rangle=\sum_j
\lambda_jY_j$. Since $QP_ju\in E_j$ and the $E_j$ are orthogonal,
$v^\top Qv=\sum_j\lambda_j^2(P_ju)^\top QP_ju/\sigma_j^2=|\lambda|^2$,
as $(P_ju)^\top QP_ju=\E[\langle P_ju,Z\rangle^2]=\sigma_j^2$. Hence
$\E[e^{\sum_j\lambda_jY_j}]\le e^{B^2|\lambda|^2/2}$.

\emph{Mixed moments.} Let $\mathcal A\subseteq\mathcal J$. For $w\ge0$
and $k\ge1$ we have $w^{2k}\le(4k/e^2)^ke^{\sqrt kw}$, by maximizing
$w^{2k}e^{-\sqrt kw}$ at $w=2\sqrt k$. Also
$\prod_je^{\sqrt{k_j}|Y_j|}\le\sum_{\varsigma}\exp(\sum_j\varsigma_j
\sqrt{k_j}Y_j)$, the sum running over all signs
$\varsigma\in\{-1,1\}^{\mathcal A}$. Therefore
\begin{equation}\label{lem:x-rc-eq-mixed}
 \E[\prod_{j\in\mathcal A}Y_j^{2k_j}]
 \le\prod_{j\in\mathcal A}(\frac{4k_j}{e^2})^{k_j}
     \cdot2^{|\mathcal A|}e^{B^2\sum_{j\in\mathcal A}k_j/2}
 \le\prod_{j\in\mathcal A}(C_6k_j)^{k_j},\qquad C_6:=2e^{B^2/2},
\end{equation}
using $2\le2^{k_j}$ and $8/e^2\le2$. No independence of the blocks is
used.

\emph{Expansion.} Define $D_j:=e^{\mathrm itZ_j}-\varphi_{k_j}(tZ_j)$ and
$D_\circ:=e^{\mathrm itw_\circ}-\varphi_{k_\circ}(tw_\circ)$. Writing
each Taylor polynomial as the phase minus its error and expanding, the
difference in the lemma is a signed sum over the nonempty subsets
$\mathcal A\subseteq\mathcal J\cup\{\circ\}$ of terms whose absolute
values are at most $\prod_{j\in\mathcal A}|D_j|$, since the phases have
modulus one. On $\mathcal G$,
$|D_\circ|\le(t\rho_\circ)^{k_\circ}/k_\circ!\le x_\circ:=
(et\rho_\circ/k_\circ)^{k_\circ}$. By Eq.~\eqref{lem:x-rc-eq-mixed} and
$k!\ge(k/e)^k$, for $\mathcal A\subseteq\mathcal J$,
\[
 \|\prod_{j\in\mathcal A}D_j\|_{L_2}
 \le\prod_{j\in\mathcal A}\frac{(t\sigma_j)^{k_j}(C_6k_j)^{k_j/2}}{k_j!}
 \le\prod_{j\in\mathcal A}x_j,\qquad
 x_j:=(et\sigma_j(\frac{C_6}{k_j})^{1/2})^{k_j}.
\]
Hence the norm in the lemma is at most
$(1+x_\circ)\prod_j(1+x_j)-1$. As $C_B=e^4C_6\ge e^2$, the lower bounds
on the $k$'s give $et\sigma_j(C_6/k_j)^{1/2}\le(t/T)/e$ and
$et\rho_\circ/k_\circ\le(t/T)/e$. Since every $k$ is at least
$C_B\log(2r)\ge3$ and $t\le T$, each $x$ is at most
$(t/T)^3e^{-C_B\log(2r)}\le(t/T)^3(2r)^{-100}$. Their sum $y$ is at most
$(4\log(2r)+1)(2r)^{-100}(t/T)^3\le r^{-10}(t/T)^3/2\le1$, and the norm
is at most $e^y-1\le2y$.
\end{proof}

\subsection{One Approximation Space for a Regular Component}

The next lemma is the heart of the upper bound. The row is random and
the query direction $u$ and the shift $a$ are fixed, and the bound must
hold for all $u$ and $a$ with one space.

\begin{lemma}[A common space for one component]\label{lem:x-rc-budget}
Fix $1\le p<2$ and $B\ge1$. There is $C_1\ge1$, depending only on $p$
and $B$, with the following property. Let $r\ge2$, let $L$ be real and
$\ell$ an integer with $1\le\ell\le L\le r$, and let $Z$ be a finitely
supported centered random vector in a Euclidean space $W$ of dimension
at most $r$, with covariance $Q\ne0$. Suppose that
$\E[e^{\langle v,Z\rangle}]\le e^{B^2v^\top Qv/2}$ for $v\in W$ and that
\[
 \Pr[0<|a+\langle v,Z\rangle|\le z(v^\top Qv)^{1/2}]\le Bz
\]
for all $v\in W$ with $v^\top Qv>0$, all $a\in\R$ and all
$L^{-1/2}\le z\le1$. Define $\Sigma:=\tr[Q^{p/2}]$. Then there is a space
$\mathcal V$ of real functions on the support of $Z$, containing the
constants, with $\dim\mathcal V\le\exp(C_1(\ell+\log^3(2r)))$, such that
for every unit vector $u\in W$ and every $a\in\R$,
\begin{equation}\label{lem:x-rc-eq-budget}
 \min_{f\in\mathcal V}\E[(|a+\langle u,Z\rangle|^p-f(Z))^2]
 \le C_1\,\ell^{-p-3/2}\,\Sigma\,(u^\top Qu)^{p/2}.
\end{equation}
\end{lemma}

\begin{proof}
Define $s:=p+1/2$, let $C_B$ be the constant of Lemma~\ref{lem:x-rc-phase},
and let $\varpi_p$, $\kappa_T$ and $g_T$ be as in
Lemma~\ref{lem:x-rc-cosine}. Constants $C$ below depend only on $p$ and
$B$.

\emph{Spectral blocks.} Let $q_1\ge q_2\ge\cdots$ be the positive
eigenvalues of $Q$, repeated according to multiplicity, so that $\Sigma=\sum_iq_i^{p/2}$. For an interval
$J\subseteq(0,\infty)$ let $E(J)$ be the span of the eigenvectors of $Q$
with eigenvalues in $J$, a subspace invariant under $Q$. Each eigenvalue
above $q_*>0$ contributes more than $q_*^{p/2}$ to $\Sigma$, so
\begin{equation}\label{lem:x-rc-eq-count}
 \dim E((q_*,\infty))\le\Sigma q_*^{-p/2}\qquad(q_*>0).
\end{equation}
For $v\in\ker Q$ we have $\E[\langle v,Z\rangle^2]=0$, so the support of
$Z$ lies in the range of $Q$. For a subspace $E$ with orthogonal
projection $P_E$ and an integer $m\ge1$, a \emph{polynomial of order $m$
on $E$} is a function $\zeta\mapsto f(P_E\zeta)$ on the support of $Z$,
where $f$ is a real polynomial on $E$ of degree less than $m$. These
functions form a space of dimension at most
$\binom{\dim E+m-1}{m-1}$.

For $q_*>0$ and a real $n$ with $1\le n\le r$ we define blocks. Let
$j_*:=\min\{j\ge0:n2^{-j}\le1\}$, so $j_*\le1+\log_2n$. Define
$E_0:=E((q_*,\infty))$, $E_j:=E((q_*2^{-j},q_*2^{1-j}])$ for
$1\le j\le j_*$, and $E_{j_*+1}:=E((0,q_*2^{-j_*}])$. These $j_*+2\le
4\log(2r)$ subspaces are mutually orthogonal, invariant under $Q$, and
their sum is the range of $Q$. Give them the loads $n_0:=n$,
$n_j:=n2^{1-j}$ for $1\le j\le j_*$, and $n_{j_*+1}:=1$, so that
$n_j\ge1$ for every $j$, and the orders
$k_j:=\lceil C_B(n_j+\log(2r))\rceil$. Let $\mathcal P(q_*,n)$ be the
span of the products $\prod_{j=0}^{j_*+1}f_j$, where $f_j$ is a
polynomial of order $k_j$ on $E_j$. Define $T:=(n/q_*)^{1/2}$. For a unit
vector $u$ and $1\le j\le j_*+1$, the score $\langle u,P_jZ\rangle$ has
variance $(P_ju)^\top QP_ju\le q_*2^{1-j}$ if $j\le j_*$ and at most
$q_*2^{-j_*}$ if $j=j_*+1$, so in both cases
\begin{equation}\label{lem:x-rc-eq-load}
 T^2\,\E[\langle u,P_jZ\rangle^2]\le n_j\qquad(1\le j\le j_*+1).
\end{equation}

We claim that, with $h_0:=\Sigma q_*^{-p/2}$,
\begin{equation}\label{lem:x-rc-eq-dim}
 \log\dim\mathcal P(q_*,n)\le22C_Bn\log(e(1+h_0/n))+30C_B\log^3(2r).
\end{equation}
The dimension of $\mathcal P(q_*,n)$ is at most the product over $j$ of
the dimensions of the spaces of polynomials of order $k_j$ on $E_j$. Let
$m_j:=\dim E_j$, $a_{j,1}:=\lceil C_Bn_j\rceil$ and
$a_{j,2}:=\lceil C_B\log(2r)\rceil$, so $k_j-1\le a_{j,1}+a_{j,2}$.
Every monomial of degree at most $a_{j,1}+a_{j,2}$ is a product of
monomials of degrees at most $a_{j,1}$ and $a_{j,2}$. Hence the
$j$-th dimension is at most $\binom{m_j+a_{j,1}}{a_{j,1}}
\binom{m_j+a_{j,2}}{a_{j,2}}$. We use
$\log\binom{m+a}a\le a\log(e(1+m/a))$, which increases in $m$ and in
$a$. As $n_j\le a_{j,1}\le2C_Bn_j$, $a_{j,2}\le2C_B\log(2r)$, $m_j\le r$
and $\log(e(1+r))\le3\log(2r)$,
\[
 \log[\binom{m_j+a_{j,1}}{a_{j,1}}\binom{m_j+a_{j,2}}{a_{j,2}}]
 \le2C_Bn_j\log(e(1+m_j/n_j))+6C_B\log^2(2r).
\]
By Eq.~\eqref{lem:x-rc-eq-count}, $m_0\le h_0$ and $m_j\le h_02^{jp/2}$ for
$1\le j\le j_*$. With $\omega_0:=\log(e(1+h_0/n))\ge1$ and
$2^{jp/2}2^{j-1}\le2^{2j}$, the first terms are at most
$2C_Bn\omega_0$ for $j=0$, at most $2C_Bn2^{1-j}(\omega_0+2j)$ for
$1\le j\le j_*$, and at most $6C_B\log(2r)$ for $j=j_*+1$. Since
$\sum_{j\ge1}2^{1-j}(\omega_0+2j)=2\omega_0+8\le10\omega_0$, and since
there are at most $4\log(2r)$ blocks and $\log(2r)\ge\log4>1$,
the claim follows.

\emph{Three cases.} Fix a unit vector $u\in W$ and $a\in\R$. Define
$y:=a+\langle u,Z\rangle$ and $\sigma^2:=u^\top Qu$, and let $b_u$ be the
right side of Eq.~\eqref{lem:x-rc-eq-budget} without the factor $C_1$, that is,
$b_u:=\ell^{-s-1}\Sigma\sigma^p$. If $\sigma=0$, then
$\langle u,Z\rangle=0$ almost surely and $y=a$ is constant. Let
$\sigma>0$ and define $h_u:=\Sigma/\sigma^p$. As $\sigma^2\le q_1$, we have
$h_u\ge1$. The score $Y:=\langle u,Z\rangle$ is centered with variance
$\sigma^2$, satisfies the moment generating function bound of
Lemma~\ref{lem:x-rc-moments} (take $v=\lambda u$), and satisfies the
interval hypothesis of Lemma~\ref{lem:x-rc-moments}(c) (take $v=u$).
By Lemmas~\ref{lem:x-rc-cosine} and~\ref{lem:x-rc-moments}(c), for every
$T>0$,
\begin{equation}\label{lem:x-rc-eq-trunc}
 \E[(|y|^p-g_T(y))^2]
 \le36\varpi_p^2T^{-2p}\E[(1+T|y|)^{-2}\mathbf 1_{y\ne0}]
 \le36\varpi_p^2BT^{-2p}(L^{-1/2}+\min\{1,(T\sigma)^{-1}\}),
\end{equation}
because $|y|^p-g_T(y)$ vanishes where $y=0$.

\emph{Case 1: $h_u\ge\ell^{s+1}$.} Use the affine function
$f(Z):=|a|^p+p|a|^{p-1}\operatorname{sgn}(a)\langle u,Z\rangle$, with
$\operatorname{sgn}(0):=0$. For $1<p\le2$ the derivative of $|t|^p$ is
H\"older of order $p-1$ with constant $2^{2-p}p$, so the Taylor
remainder satisfies $||a+w|^p-|a|^p-p|a|^{p-1}\operatorname{sgn}(a)w|
\le2^{2-p}|w|^p\le2|w|^p$, and for $p=1$ the same bound with constant two
is checked directly. By Lemma~\ref{lem:x-rc-moments}(a),
$\E[(|y|^p-f(Z))^2]\le4\E[|Y|^{2p}]\le4(16B^4\sigma^4)^{p/2}
\le64B^4\sigma^{2p}$. Since $\sigma^p=\Sigma/h_u\le\Sigma\ell^{-s-1}$, this
is at most $64B^4b_u$. Let $\mathcal V_1$ be the span of the constants
and the linear functions $\zeta\mapsto\langle v,\zeta\rangle$, of
dimension at most $r+1$.

\emph{Case 2: $\ell\le h_u<\ell^{s+1}$.} Let $\mathcal H$ be the set of
numbers $h=\ell2^k$, $k\ge0$ an integer, with $h<\ell^{s+1}$, so
$|\mathcal H|\le1+s\log_2\ell\le5\log(2r)$. For $h\in\mathcal H$ define
\[
 q_*(h):=(\Sigma/h)^{2/p},\qquad n(h):=\ell(\ell/h)^{1/s},\qquad
 \mathcal V_h:=\mathcal P(q_*(h),n(h)).
\]
As $\ell\le h<\ell^{s+1}$, we have $1<n(h)\le\ell\le r$. Here
$h_0=\Sigma q_*(h)^{-p/2}=h$. With $x:=h/\ell\ge1$ we get
$h/n(h)=x^{1+1/s}$ and $n(h)\log(e(1+h/n(h)))\le\ell
x^{-1/s}(1+\log2+(1+1/s)\log x)\le3\ell$, since $x^{-1/s}\log x\le s/e$
and $s<5/2$. By Eq.~\eqref{lem:x-rc-eq-dim},
$\log\dim\mathcal V_h\le66C_B\ell+30C_B\log^3(2r)$.

Let $h\in\mathcal H$ satisfy $h\le h_u<2h$, and write $q_*$, $n$ for
$q_*(h)$, $n(h)$, with the blocks $E_j$, orders $k_j$ and frequency
$T:=(n/q_*)^{1/2}$ defined above. Then
$\sigma^p=\Sigma/h_u\in(q_*^{p/2}/2,q_*^{p/2}]$, so
$\sigma^2\le q_*$, and Eq.~\eqref{lem:x-rc-eq-load} extends to $j=0$ because
$\E[\langle u,P_0Z\rangle^2]\le\sigma^2$. Write $Z_j:=\langle u,P_jZ\rangle$,
so that $y=a+\sum_jZ_j$ almost surely, and for $\zeta$ in the support of
$Z$ and $0<t\le T$ define
\[
 \Phi_t(\zeta):=1-\operatorname{Re}[e^{\mathrm ita}
   \prod_{j=0}^{j_*+1}\varphi_{k_j}(t\langle u,P_j\zeta\rangle)],
 \qquad
 \Phi(\zeta):=\int\Phi_t(\zeta)\,\d\kappa_T(t).
\]
For each $t$, the real and imaginary parts of the product are in
$\mathcal V_h$, and so is the constant one, hence $\Phi_t\in\mathcal V_h$.
For fixed $\zeta$, $t\mapsto\Phi_t(\zeta)$ is smooth with
$\Phi_0(\zeta)=0$ and vanishing derivative at $t=0$, since every
$\varphi_{k_j}(w)=1+\mathrm iw+O(w^2)$. Thus $|\Phi_t(\zeta)|=O(t^2)$ and
the integral converges. Every linear functional on the finite-dimensional
space of functions on the support of $Z$ that vanishes on $\mathcal V_h$
is a finite combination of point evaluations, so it also vanishes on
$\Phi$. Hence $\Phi\in\mathcal V_h$. As
$\cos(ty)=\operatorname{Re}[e^{\mathrm ita}\prod_je^{\mathrm itZ_j}]$,
Minkowski's inequality (Lemma~\ref{lem:pre_minkowski}(b)), Lemma~\ref{lem:x-rc-phase} with all blocks
$E_0,\ldots,E_{j_*+1}$, $\mathcal G$ the whole space, $w_\circ:=0$,
$\rho_\circ:=0$ and $k_\circ:=\lceil C_B\log(2r)\rceil$, and
Lemma~\ref{lem:x-rc-cosine} give
\[
 \|g_T(y)-\Phi(Z)\|_{L_2}
 \le\int\|\prod_je^{\mathrm itZ_j}-\prod_j\varphi_{k_j}(tZ_j)
   \|_{L_2}\d\kappa_T(t)
 \le r^{-10}\int(t/T)^3\,\d\kappa_T(t)\le2\varpi_pr^{-10}T^{-p}.
\]
Now $T\sigma\ge(n/q_*)^{1/2}2^{-1/p}q_*^{1/2}\ge\sqrt n/2$ and
$L^{-1/2}\le n^{-1/2}$, since $n\le\ell\le L$. With $T^{-2p}=(q_*/n)^p$
and Eq.~\eqref{lem:x-rc-eq-trunc},
\[
 \E[(|y|^p-\Phi(Z))^2]
 \le2\cdot36\varpi_p^2B(\frac{q_*}n)^p\frac3{\sqrt n}
   +8\varpi_p^2r^{-20}(\frac{q_*}n)^p
 \le224\varpi_p^2B\,q_*^pn^{-s}.
\]
Finally $n^{-s}=\ell^{-s}(h/\ell)=h\ell^{-s-1}$, while
$b_u=\ell^{-s-1}\Sigma\sigma^p\ge\ell^{-s-1}hq_*^{p/2}\cdot q_*^{p/2}/2$.
The error is therefore at most $448\varpi_p^2Bb_u$.

\emph{Case 3: $h_u<\ell$.} Define $q_*:=(\Sigma/\ell)^{2/p}$, $n:=\ell$ and
$T:=(\ell/q_*)^{1/2}$, with the blocks $E_j$ and orders $k_j$ defined
above for these $q_*$ and $n$. Now $h_0=\ell$, so
$\dim E_0\le\ell$ by Eq.~\eqref{lem:x-rc-eq-count}, and
$\sigma^p=\Sigma/h_u>q_*^{p/2}$, so $\sigma^2>q_*$ and $E_0\ne0$. The
variance of $\langle u,P_0Z\rangle$ may be much larger than $q_*$, so
we partition the support of $Z$ according to the position of $P_0Z$.

Let $Q_0$ be the restriction of $Q$ to $E_0$, which is invertible, and
let $\mathcal G$ be the event $\|Q_0^{-1/2}P_0Z\|_2\le R_0:=6B\sqrt\ell$.
For a unit vector $w\in E_0$, the vector $v:=Q_0^{-1/2}w$ satisfies
$\langle w,Q_0^{-1/2}P_0Z\rangle=\langle v,Z\rangle$ and $v^\top Qv=1$,
so $\Pr[\langle w,Q_0^{-1/2}P_0Z\rangle\ge t]\le e^{-t^2/(2B^2)}$. Let
$\mathcal N_0$ be a $1/2$-net of the unit sphere of $E_0$ with at most
$6^{\dim E_0}\le6^\ell$ points (Lemma~\ref{lem:pre_nets}). For every
$x\in E_0$ there is $w\in\mathcal N_0$ with $\langle w,x\rangle\ge
\|x\|_2/2$, so $\|Q_0^{-1/2}P_0Z\|_2>R_0$ forces
$\langle w,Q_0^{-1/2}P_0Z\rangle>R_0/2$ for some $w\in\mathcal N_0$.
Hence
\[
 \Pr[\mathcal G^c]\le6^\ell e^{-R_0^2/(8B^2)}=e^{\ell(\log6-4.5)}
 \le e^{-2\ell}.
\]
The set $\mathcal E_0:=\{x\in E_0:\|Q_0^{-1/2}x\|_2\le R_0\}$ is an
ellipsoid whose semiaxes $\mathfrak a_i:=R_0q_i^{1/2}$ correspond to the
eigenvalues $q_i>q_*$. Define $\rho:=(\ell q_*)^{1/2}$, so $T\rho=\ell$,
and let $\mathfrak c_1,\ldots,\mathfrak c_{N_\circ}$ be a maximal subset of $\mathcal E_0$ with
pairwise distances greater than $\rho$. Then every point of
$\mathcal E_0$ is within $\rho$ of some $\mathfrak c_i$. The balls of radius
$\rho/2$ around the $\mathfrak c_i$ are disjoint and lie in
$\mathcal E_0+(\rho/2)B_2$, which is contained in twice the ellipsoid
with semiaxes $(\mathfrak a_i^2+\rho^2/4)^{1/2}$, since $(x_i+x'_i)^2\le
2x_i^2+2x_i'^2$ coordinatewise in the eigenbasis. Comparing volumes,
\[
 N_\circ\le\prod_i\frac{2(\mathfrak a_i^2+\rho^2/4)^{1/2}}{\rho/2}
 \le\prod_i(2+\frac{4\mathfrak a_i}\rho)
 =\prod_i(2+24B(\frac{q_i}{q_*})^{1/2}).
\]
As $q_i/q_*>1$ and $\log t\le(2/p)t^{p/2}$ for $t\ge1$,
$\log N_\circ\le\sum_i(\log(26B)+\frac1p(q_i/q_*)^{p/2})\le
\ell(\log(26B)+1)$, using $\dim E_0\le\ell$ and $\sum_iq_i^{p/2}\le
\Sigma=\ell q_*^{p/2}$.

We also write $\mathcal G$ for the set of atoms $\zeta$ of $Z$ with
$\|Q_0^{-1/2}P_0\zeta\|_2\le R_0$. For $\zeta\in\mathcal G$ we have
$P_0\zeta\in\mathcal E_0$. Let $i(\zeta)$ be the least $i$ with
$\|P_0\zeta-\mathfrak c_i\|_2\le\rho$, and let $\mathcal C_i$ be the set
of $\zeta\in\mathcal G$ with $i(\zeta)=i$. Let $\mathcal V_3$ be
the span of the functions $\mathbf 1_{\mathcal C_i}f$, $i\in[N_\circ]$,
$f\in\mathcal P(q_*,\ell)$, and of $\mathbf 1_{\mathcal G^c}$ and
$\mathbf 1_{\mathcal G^c}\langle v,\cdot\rangle$, $v\in W$. By the bound
on $N_\circ$ and Eq.~\eqref{lem:x-rc-eq-dim} with $h_0=n=\ell$,
$\log\dim\mathcal V_3\le\ell(\log(26B)+1)+38C_B\ell+30C_B\log^3(2r)
+\log(r+2)$.

Write $Z_j:=\langle u,P_jZ\rangle$ for $j\ge1$, define
$w_\circ:=\langle u,P_0Z-\mathfrak c_{i(Z)}\rangle$ on $\mathcal G$ and
$w_\circ:=0$ off $\mathcal G$, so that
$|w_\circ|\le\|P_0Z-\mathfrak c_{i(Z)}\|_2\le\rho$ and
$y=a+\langle u,\mathfrak c_{i(Z)}\rangle+w_\circ+\sum_{j\ge1}Z_j$ on $\mathcal G$.
For $i\in[N_\circ]$, an atom $\zeta$ of $Z$ and $0<t\le T$ define
\[
 \Phi_{i,t}(\zeta):=1-\operatorname{Re}[e^{\mathrm it(a+\langle
   u,\mathfrak c_i\rangle)}\varphi_{k_0}(t\langle u,P_0\zeta-\mathfrak c_i\rangle)
   \prod_{j=1}^{j_*+1}\varphi_{k_j}(t\langle u,P_j\zeta\rangle)
   ],\qquad\Phi_i:=\int\Phi_{i,t}\,\d\kappa_T(t),
\]
and use the approximant
$\Psi:=\sum_i\mathbf 1_{\mathcal C_i}\Phi_i+\mathbf 1_{\mathcal G^c}f$, with $f$
the affine function of Case~1. The factor $\varphi_{k_0}(t\langle
u,P_0\zeta-\mathfrak c_i\rangle)$ is a polynomial of order $k_0$ on $E_0$, so the
argument of Case~2 shows $\mathbf 1_{\mathcal C_i}\Phi_i\in\mathcal V_3$, and
$\Psi\in\mathcal V_3$. Lemma~\ref{lem:x-rc-phase}, applied to the blocks
$E_1,\ldots,E_{j_*+1}$ with Eq.~\eqref{lem:x-rc-eq-load}, to $\mathcal G$,
$w_\circ$, $\rho_\circ:=\rho$ and $k_\circ:=k_0\ge C_B(\ell+\log(2r))=
C_B(T\rho+\log(2r))$, gives as in Case~2
$\|\mathbf 1_{\mathcal G}(g_T(y)-\Psi(Z))\|_{L_2}\le2\varpi_pr^{-10}
T^{-p}$. Since $T\sigma>(\ell/q_*)^{1/2}q_*^{1/2}=\sqrt\ell$, $L^{-1/2}\le
\ell^{-1/2}$ and $T^{-2p}=(q_*/\ell)^p$, Eq.~\eqref{lem:x-rc-eq-trunc} gives
\[
 \E[(|y|^p-\Psi(Z))^2\mathbf 1_{\mathcal G}]
 \le2\cdot72\varpi_p^2B(\frac{q_*}\ell)^p\ell^{-1/2}
   +8\varpi_p^2r^{-20}(\frac{q_*}\ell)^p
 \le152\varpi_p^2Bq_*^p\ell^{-s}\le152\varpi_p^2Bb_u,
\]
because $b_u=\ell^{-s-1}\cdot\ell q_*^{p/2}\cdot\sigma^p\ge\ell^{-s}
q_*^p$. Off $\mathcal G$, the Taylor bound of Case~1, the Cauchy--Schwarz
inequality and Lemma~\ref{lem:x-rc-moments}(a) give
\[
 \E[(|y|^p-f(Z))^2\mathbf 1_{\mathcal G^c}]
 \le4(\E[|Y|^{4p}])^{1/2}\Pr[\mathcal G^c]^{1/2}
 \le112B^4\sigma^{2p}e^{-\ell}.
\]
As $\sigma^p\le\Sigma$, this is at most
$112B^4\ell^{s+1}e^{-\ell}b_u\le300B^4b_u$.

\emph{The space.} Let $\mathcal V:=\mathcal V_1+\sum_{h\in\mathcal H}
\mathcal V_h+\mathcal V_3$. It contains the constants, every case above
uses one of its subspaces, and it depends only on the law of $Z$ and on
$\ell$, $r$, $p$ and $B$. Its dimension is at most the sum of the dimensions,
and the bounds above give $\dim\mathcal V\le\exp(C_1(\ell+\log^3(2r)))$
for a suitable $C_1$ depending on $p$ and $B$. Taking $C_1$ also at least
$448\varpi_p^2B+300B^4$ proves Eq.~\eqref{lem:x-rc-eq-budget}.
\end{proof}

\subsection{Proof of the Bounds for Regular Components}

\begin{proof}[Proof of the upper bound in Lemma~\ref{lem:joint-regular-components}]
Let $A$ be an input of the class. Constants $C_2,C_3,\ldots$ depend only
on $p$ and $B$, and we show that the bound holds whenever
$C_2\log^3(2r)\le L\le r$. This covers the stated range once $r\ge r_0$,
with $r_0$ depending on $p$ and $B$, because $L\ge r^{2/(p+2)}\ge r^{1/2}$
there. For $r<r_0$, Theorem~\ref{thm:general_lp_field} gives a sparsifier
with at most $Cr/\epsilon^2$ rows, and
$r^2\log r/L^{p+3/2}\ge(\log2)r^{1/2-p}\ge(\log2)r_0^{-3/2}$, so its
support is at most $C'_{p,B}\epsilon^{-2}r^2\log r/L^{p+3/2}$.

\emph{Rows.} Let $I_0$ be the set of pairs $(g,\zeta)$ with $\zeta$ an
atom of $Z_g$, and define $\mu_{g\zeta}:=c_g\Pr[Z_g=\zeta]\theta_g(\zeta)>0$.
The row of $(g,\zeta)$ is $a_{g\zeta}=\mu_{g\zeta}^{1/p}(m_g+\zeta)$, so
$F(x)=\sum_{I_0}\mu_{g\zeta}|\langle m_g+\zeta,x\rangle|^p$. Let
$I\subseteq I_0$ be the set of pairs with $m_g+\zeta\ne0$. Zero rows do
not matter. Let $W$ be the sum of the ranges of the $Q_g$ in the query
space and $d_W:=\dim W$. Since $Q_g=\sum_\zeta\Pr[Z_g=\zeta]\zeta
\zeta^\top$, the range of $Q_g$ is spanned by the atoms of $Z_g$. As
$Z_g$ is centered, every atom is
$\zeta=\sum_{\zeta'}\Pr[Z_g=\zeta']((m_g+\zeta)-(m_g+\zeta'))$, and each
$m_g+\zeta$ is zero or a positive multiple of a row. Thus $W$ lies in the
row space and $d_W\le r$.

If $d_W=0$, every $Z_g$ is zero, every component has one row, and
keeping all $G$ rows with weight one is exact. As $G\le e^{cL}\le
e^{L/16}$ and
$\epsilon^{-2}r^2\log r/L^{p+3/2}\ge e^{2L}(\log2)r^{-3/2}/4\ge e^{L/16}$
for $L\ge C_2\log^3(2r)$, this case is done. Let $d_W\ge1$.

\emph{Normalization.} The operator $\sum_gc_gQ_g$ is positive definite
on $W$, since a vector of $W$ in its kernel is orthogonal to every range.
Let $T$ be given by Lemma~\ref{lem:x-rc-normal}, extend it by the
identity on $W^\perp$ to $\widetilde T$, and replace every row $a$ by
$\widetilde Ta$ and every query $x$ by $\widetilde T^{-\top}x$. Inner
products of rows with queries are unchanged, so the energies and the
sparsifiers are unchanged, and the rank stays $r$. The new references
$\widetilde TZ_g$ are centered with covariances $\widetilde T
Q_g\widetilde T^\top$, and the two hypotheses of the class hold for them,
because $\langle v,\widetilde TZ_g\rangle=\langle\widetilde T^\top
v,Z_g\rangle$. We keep the old symbols for the new data. Then $W$ is
unchanged and, by Eq.~\eqref{lem:x-rc-eq-trace},
\[
 \sum_gc_gQ_g^{p/2}=I_W,\qquad\sum_gc_g\Sigma_g=d_W,\qquad
 \Sigma_g:=\tr[Q_g^{p/2}].
\]

\emph{Energy.} Let $x$ be a query and $z:=P_Wx$. We claim
\begin{equation}\label{lem:x-rc-eq-energy}
 F(x)\ge\kappa_B\sum_gc_g(z^\top Q_gz)^{p/2},\qquad
 \sum_gc_g(z^\top Q_gz)^{p/2}\ge\|z\|_2^p,\qquad
 \kappa_B:=(8\sqrt2B^3)^{-1}.
\end{equation}
The score $Y_g:=\langle Z_g,x\rangle=\langle Z_g,z\rangle$ is centered
with variance $x^\top Q_gx=z^\top Q_gz$ and satisfies
$\E[e^{\lambda Y_g}]\le e^{B^2\lambda^2z^\top Q_gz/2}$. As $\theta_g\ge
B^{-1}$, Lemma~\ref{lem:x-rc-moments}(b) gives
\[
 F(x)\ge B^{-1}\sum_gc_g\E[|\langle m_g,x\rangle+Y_g|^p]
 \ge\kappa_B\sum_gc_g(z^\top Q_gz)^{p/2}.
\]
For the second inequality let $z\ne0$,
define $e:=z/\|z\|_2$, and expand $e$ in an orthonormal eigenbasis
$(f_i)$ of $Q_g$ on $W$ with eigenvalues $q_{g,i}$. The weights
$\langle f_i,e\rangle^2$ sum to one, so concavity of $t\mapsto t^{p/2}$
gives $(e^\top Q_ge)^{p/2}\ge\sum_i\langle f_i,e\rangle^2q_{g,i}^{p/2}
=e^\top Q_g^{p/2}e$. Summing with weights $c_g$ gives
$\sum_gc_g(e^\top Q_ge)^{p/2}\ge e^\top I_We=1$, and homogeneity
gives the claim.

\emph{The approximation space.} Let $C_1$ be the constant of
Lemma~\ref{lem:x-rc-budget}, define $c_\ell:=1/(32C_1)$ and
$\ell:=\lfloor c_\ell L\rfloor$, and take $C_2$ so large that $c_\ell L\ge1$ and
$C_1\log^3(2r)\le L/32$. Then $1\le\ell\le L\le r$, $\ell\ge c_\ell L/2$,
and $C_1(\ell+\log^3(2r))\le L/16$. For each $g$ with $Q_g\ne0$,
Lemma~\ref{lem:x-rc-budget} applies to $Z_g$ in $W$, whose dimension is
at most $r$, and gives a space $\mathcal V_g$ of functions on the atoms
of $Z_g$ of dimension at most $e^{L/16}$. If $Q_g=0$, then $Z_g=0$, and
$\mathcal V_g$ is the space of constants on the single atom.

Let $\gamma_i$, $i\in I$, be the Lewis weights of Eq.~\eqref{lem:s-eq-lewis}
for the rows of $I$ (written $c_i$ there), so $\sum_i\gamma_i=r$. Define
\[
 \Sigma_0:=\sum_{(g,\zeta)\in I}c_g\Sigma_g\Pr[Z_g=\zeta],\qquad
 \pi_{g\zeta}:=\frac12(\frac{c_g\Sigma_g\Pr[Z_g=\zeta]}{\Sigma_0}
   +\frac{\gamma_{g\zeta}}r)\qquad((g,\zeta)\in I).
\]
Some $Q_g$ is nonzero, its reference has at least two atoms, and at most
one of them has $m_g+\zeta=0$. Hence $0<\Sigma_0\le\sum_gc_g\Sigma_g
=d_W$. Thus $\pi$ is a probability vector on $I$ with
\[
 \pi_{g\zeta}\ge\frac{\gamma_{g\zeta}}{2r}\qquad\text{and}\qquad
 \pi_{g\zeta}\ge\frac{c_g\Sigma_g\Pr[Z_g=\zeta]}{2d_W},
\]
so Eq.~\eqref{lem:s-eq-importance} holds with $\beta=2$. Let $\mathcal V$ be
the space of functions $(g,\zeta)\mapsto\mu_{g\zeta}f_g(\zeta)/
\pi_{g\zeta}$ on $I$, where $f_g\in\mathcal V_g$ for every $g$. Then
$1\le\dim\mathcal V\le\sum_g\dim\mathcal V_g\le e^{cL+L/16}\le e^{L/8}$,
as $c\le1/16$.

\emph{The approximation error.} Let $F(x)\le1$, define $z:=P_Wx$, and
recall $\kappa_x(i)=|a_i^\top x|^p/\pi_i$ from Lemma~\ref{lem:s-cal}. For
an atom $\zeta$ of $Z_g$ we have $\langle m_g+\zeta,x\rangle=\langle
m_g,x\rangle+\langle\zeta,z\rangle$. Let $e_g(x)^2$ be the least value
of $\E[(|\langle m_g,x\rangle+\langle Z_g,z\rangle|^p-f(Z_g))^2]$ over
$f\in\mathcal V_g$, attained at $f_g$. If $Q_g=0$ or $z=0$, this function
of $Z_g$ is constant and $e_g(x)=0$. Otherwise apply
Lemma~\ref{lem:x-rc-budget} to $u:=z/\|z\|_2$ and
$a:=\langle m_g,x\rangle/\|z\|_2$, and multiply by $\|z\|_2^{2p}$:
\[
 e_g(x)^2\le C_1\ell^{-p-3/2}\Sigma_g(z^\top Q_gz)^{p/2}\|z\|_2^p .
\]
Let $\phi\in\mathcal V$ be built from these $f_g$. As
$\mu_{g\zeta}\le Bc_g\Pr[Z_g=\zeta]$, the second lower bound on $\pi$
and Eq.~\eqref{lem:x-rc-eq-energy} give
\[
 \|\kappa_x-\phi\|_{L_2(\pi)}^2
 =\sum_{(g,\zeta)\in I}\frac{\mu_{g\zeta}^2}{\pi_{g\zeta}}
   (|\langle m_g+\zeta,x\rangle|^p-f_g(\zeta))^2
 \le2B^2d_W\sum_{g:\,Q_g\ne0}\frac{c_g}{\Sigma_g}e_g(x)^2
 \le\frac{2B^2C_1d_W}{\kappa_B^2\ell^{p+3/2}}F(x)^2 .
\]
With $F(x)\le1$, $d_W\le r$ and $\ell\ge c_\ell L/2$, the
hypothesis Eq.~\eqref{lem:s-eq-cal-hyp} holds with
\[
 \delta^2:=C_3\frac r{L^{p+3/2}},\qquad
 C_3:=\max\{1,\frac{2B^2C_1}{\kappa_B^2}(\frac2{c_\ell})^{p+3/2}
 \}.
\]
As $L\le r$, $\delta^2\ge r^{-p-1/2}\ge r^{-4}$, so $\delta\ge r^{-B_0}$
with $B_0=2$.

\emph{Sampling.} Take $C_2$ so large that $L\ge C\log(2r)$ for the
constant $C$ of Lemma~\ref{lem:s-cal} with $\beta=2$ and $B_0=2$. Part~(a)
of that lemma gives positive weights on a set of at most
\[
 C\epsilon^{-2}r\delta^2\log(2r)\le2CC_3\,
 \frac{\epsilon^{-2}r^2\log r}{L^{p+3/2}}
\]
rows of $I$ that preserve every energy within relative error $\epsilon$.
These are original rows of the input, and undoing the change of
coordinates gives the same statement for the original data.
\end{proof}

\begin{proof}[Proof of the lower bound in Lemma~\ref{lem:joint-regular-components}]
We use Lemma~\ref{lem:s-shift} with the parameters
\[
 (\lambda_1,\lambda_2,\beta_0):=(1,2,2^{-(p+2)/(2p)})\quad\text{if }
 L\le r/4,\qquad
 (\lambda_1,\lambda_2,\beta_0):=(1/4,2,1)\quad\text{if }L>r/4.
\]
For each of the two triples let $c'_S$ and $d_0$ be the constants of
Lemma~\ref{lem:s-shift}, and let $B_1\ge3$ and $c_V$ be those of
Lemma~\ref{lem:s-regular} for $(\lambda_1,\lambda_2)$. Let $d_1$,
depending only on $p$, be so large that for both triples every integer
$d\ge d_1$ satisfies $d\ge d_0$ and $(2d)^{-1/2}\ge e^{-c_Vd}$, and that
for $d\ge d_1$ the event of Lemma~\ref{lem:s-halfspace} with
$a'=\lambda_1$ has probability at least $9/10$ and
Lemma~\ref{lem:s-regular} applies. Let $r_1\ge2$, depending only on
$p$, $\delta$ and $c$, be so large that $\lfloor\sqrt{r_1}\rfloor\ge
d_1$ and every $r\ge r_1$ satisfies $2r^2\le e^{c\sqrt r}$ and
$r^{\delta(p+2)/(4p)}\ge2^{(p+2)/(2p)}$. Define
$B_*:=\max\{B_1,4\}$, the maximum taken over both triples. It is an
absolute constant, because the two triples are fixed. Let
$B\ge B_*$ and $r^{2/(p+2)+\delta}\le L\le r$.

\emph{Small ranks.} Let $r<r_1$. As $\epsilon\le1/2$ we have
$L\ge\log4>1$, so $r\ge2$. Define $K:=2\lceil r/2\rceil$ and
$t_k:=(2k-1-K)/K$ for $k\in[K]$. Take one component, $G=1$, with $c_1=1$,
$\theta_1=1$, $m_1=0$ and $Z_1$ uniform on $\{t_1,\ldots,t_K\}^r$. No
$t_k$ vanishes, as $K$ is even, so the rows $K^{-r/p}\zeta$ are nonzero.
They have rank $r$, because the support contains $(t_1,\ldots,t_1)$ and
the $r$ points that differ from it in one coordinate. The coordinates of
$Z_1$ are independent and symmetric with values in $[-1,1]$ and variance
$s^2:=(K^2-1)/(3K^2)\in[1/4,1/3]$. Hence $Q_1=s^2I_r$ and
\[
 \E[e^{\langle v,Z_1\rangle}]=\prod_k\E[\cosh(v_kZ_{1,k})]\le\prod_k\cosh v_k
 \le e^{\|v\|_2^2/2}\le e^{4v^\top Q_1v/2},
\]
which is the moment condition for every $B\ge2$.

For the interval condition let $U$ be uniform on $[-1,1]^r$, and let
$Z_k$ be the midpoint of the interval $[-1+2(j-1)/K,-1+2j/K)$,
$j\in[K]$, that contains $U_k$. Then $Z$ has the law of $Z_1$ and
$|Z_k-U_k|\le1/K\le1/r$, so
$|\langle v,Z\rangle-\langle v,U\rangle|\le\|v\|_1/r\le\|v\|_2/\sqrt r$.
Fix $v\ne0$ and let $f$ be the density of $\langle v,U\rangle$. It is
even, and it is log-concave by Lemma~\ref{lem:pre_logconcave_marginal},
being a marginal of the uniform law on a convex body. Hence $M:=\sup f=f(0)$. Let
$g(t):=Me^{-2Mt}$ for $t\ge0$. Both $f$ and $g$ have integral $1/2$ over
$[0,\infty)$. The function $\log(f/g)$ is concave on $[0,\infty)$, with
values in $[-\infty,\infty)$, and vanishes at $0$, so
$\{t\ge0:f(t)\ge g(t)\}$ is an interval with left endpoint $0$. If it is
all of $[0,\infty)$, then $f=g$ almost everywhere. Otherwise let $t_0$ be
its right endpoint. Then $(t^2-t_0^2)(f(t)-g(t))\le0$ for almost every
$t\ge0$. In both cases
\[
 \int_0^\infty t^2f(t)\,\d t\le\int_0^\infty t^2g(t)\,\d t
 +t_0^2\int_0^\infty(f-g)(t)\,\d t=\frac1{4M^2},
\]
where the middle term is absent in the first case. As
$\operatorname{Var}\langle v,U\rangle=\|v\|_2^2/3$, this gives
$\|v\|_2^2/3\le1/(2M^2)$, that is, $M\le\sqrt{3/2}/\|v\|_2$. Now let
$a\in\R$, $\sigma:=\sqrt{v^\top Q_1v}=s\|v\|_2\le\|v\|_2/\sqrt3$ and
$L^{-1/2}\le z\le1$. As $L\le r$, we have $\|v\|_2/\sqrt r\le z\|v\|_2$, so
\begin{align*}
 \Pr[0<|a+\langle v,Z_1\rangle|\le z\sigma]
 &\le\Pr[|a+\langle v,U\rangle|\le(1/\sqrt3+1)z\|v\|_2]\\
 &\le2M(1/\sqrt3+1)z\|v\|_2\le(\sqrt2+\sqrt6)z\le4z .
\end{align*}
Thus this input lies in the class for every $B\ge4$, and every
sparsifier has at least one row. Since
$\epsilon^{-2}=e^{2L}/4\le e^{2r_1}$, the right side of the lemma is
bounded in terms of $p$, $\delta$ and $c$, which proves the lower bound
for $r<r_1$. This is the only place where $c$ enters the constant of the
lower bound.

\emph{Large ranks: parameters.} Let $r\ge r_1$. If $L\le r/4$, define
\[
 N:=\lfloor r/2\rfloor,\qquad d:=\lfloor L\rfloor,\qquad
 H:=\lfloor\frac{r-N}d\rfloor,\qquad b:=H^{-1/p}.
\]
Then $d\ge\lfloor\sqrt r\rfloor$ and $d\le L<d+1\le2d$, so
$\lambda_1d\le L\le\lambda_2d$. As $r-N\ge r/2$ and $d\le r/4$, we have
$(r-N)/d\ge2$, hence $H\ge(r-N)/(2d)\ge r/(4d)$, and $Hd\le r-N$. Also
$Hb^p=1$ and $b\le1$. Since $H\le r/d$ and $d\ge L/2$,
\[
 b\sqrt d\ge(\frac{d^{1+p/2}}r)^{1/p}
 \ge2^{-(p+2)/(2p)}(\frac{L^{1+p/2}}r)^{1/p}
 \ge2^{-(p+2)/(2p)}r^{\delta(p+2)/(2p)},
\]
because $L^{1+p/2}\ge r^{(2/(p+2)+\delta)(p+2)/2}=r^{1+\delta(p+2)/2}$. In
particular $b\sqrt d\ge\beta_0$, and by the choice of $r_1$,
$b\sqrt d\ge r^{\delta(p+2)/(4p)}$. If $L>r/4$, define
$N:=\lfloor r/2\rfloor$, $H:=1$, $d:=r-N\ge r/2$ and $b:=1$. Then
$d/4\le r/4<L\le r\le2d$, $Hb^p=1$, $b\sqrt d\ge1=\beta_0$, and
$r=N+Hd$. In both cases $d\ge\lfloor\sqrt r\rfloor\ge d_1$, $N\ge r/3$
and $r\ge N+Hd$.

\emph{The input.} Let $\mathcal E$ be the event of
Lemma~\ref{lem:s-halfspace} for the vectors $\sqrt du_t$ with
$a'=\lambda_1$. It has probability at least $9/10$, so
Lemma~\ref{lem:s-shift} applies with it. Let $u_1,\ldots,u_\Xi$, the law
$\mathcal Q$ and the masses $\alpha\in[1/2,3/2]^{N\times\Xi}$ be as given
there, and let $A$ be the rank-$r$ input of Lemma~\ref{lem:s-cost} built
from these data. By Lemma~\ref{lem:s-regular}, its energy has the form
of the lemma, with one component for each pair of a tag $i\in[N]$ and
a transverse copy, and one for each coordinate row. For the pair of $i$
and a copy with inclusion $\iota$, the reference is uniform on the
$\Xi$ distinct points $\zeta_t:=\iota(u_t-\bar u)$, the weight is $\Xi$,
the mean is $(e_i,\iota\bar u,0)$ and $\theta(\zeta_t)=\alpha_{it}$. The
vector $(\Xi\cdot\Xi^{-1}\alpha_{it})^{1/p}(m+\zeta_t)$ is
$\alpha_{it}^{1/p}(e_i,\iota u_t,0)$, which is the row of $A$ for this
tag, copy and $t$. A coordinate row is its own component. So the rows of
$A$ are exactly the rows prescribed by the lemma. We check the
remaining conditions of the class.
\begin{enumerate}[label=(\roman*),nosep,leftmargin=*]
\item The number of components is $NH+(r-N-Hd)\le r^2+r\le2r^2$, as
$N\le r$ and $H\le r$. Since $L\ge\sqrt r$ and $r\ge r_1$, this is at
most $e^{cL}$.
\item The tilts lie in $[1/2,3/2]\subseteq[B^{-1},B]$, as $B\ge3$, and
the references are finitely supported and centered.
\item By Lemma~\ref{lem:s-regular}, $\E[e^{\langle v,Z_g\rangle}]\le
e^{B_1^2v^\top Q_gv/2}\le e^{B^2v^\top Q_gv/2}$ for every $v$.
\item By Lemma~\ref{lem:s-regular}, the interval probability is at most
$3z\le Bz$ for $e^{-c_Vd}\le z\le1$. As $L\le\lambda_2d=2d$ and
$(2d)^{-1/2}\ge e^{-c_Vd}$, this range contains $[L^{-1/2},1]$.
\end{enumerate}
Hence $A$ lies in the class.

\emph{Counting.} By Lemma~\ref{lem:s-shift}, every reweighting of the
rows of $A$ by real weights with relative error at most $\epsilon$, and
in particular every sparsifier, gives nonzero weight to more than
$c'_S\epsilon^{-2}NHd^{-p-1/2}(1+\log_+(b\sqrt d))$ rows. If $L\le r/4$,
then $NH\ge(r/3)(r/(4d))$, $d\le L$, and
$\log_+(b\sqrt d)\ge(\delta(p+2)/(4p))\log r\ge(\delta/2)\log r$, as
$p\le2$. The count is therefore at least
\[
 \frac{c'_S\delta}{24}\cdot\frac{\epsilon^{-2}r^2\log r}{L^{p+3/2}} .
\]
If $L>r/4$, then $H=1$, $N\ge r/3$, $d\le r$, and $1+\log\sqrt d\ge
1+\frac12\log(r/2)\ge\frac12\log r$. The count is at least
$(c'_S/6)\epsilon^{-2}r^{1/2-p}\log r$, and $r^{1/2-p}=r^2/r^{p+3/2}\ge
4^{-p-3/2}r^2/L^{p+3/2}$ because $r<4L$. This proves the lower bound for
$r\ge r_1$, with a constant depending only on $p$ and $\delta$.
\end{proof}

%% file: 81_proof.tex
\section{High Accuracy for Non-Even Exponents}
\label{sec:high_accuracy}

\subsection{Main Result}

Theorem~\ref{thm:high_accuracy_main} treats accuracy $L\ge ar$ and
fixed rank for every non-even $p>0$.
In part~(a), the lower bound is the one-tag construction of Lemma~\ref{lem:s-onetag}.
The upper bound rounds the Lewis masses by iterated partial coloring on
pairs of light rows. A low-crossing matching makes the error of each query
a sum over few crossing pairs, and both this crossing count and the
entropy of the integrated error vectors have exponent $r-1$. Below two the
error is measured on caps of queries with comparable energy, as in
Lemma~\ref{lem:t-cap-family}. Above two a single residual on the whole
sphere suffices.

In part~(b), the lower bound is the one-tag construction of
Lemma~\ref{lem:s-onetag}(a). For the upper bound we round the Lewis
masses by iterated partial coloring until
$N\asymp(\epsilon^{-2}r^{1/2-p}\Lambda)^\eta$ rows remain. With
$q\asymp r N^{1/(r-1)}$, a round on $n$ rows must keep the relative
error at every query of order
\[
 (r\Lambda)^{1/2}(r^{p+1/2}q^{-2p-1}/n)^{1/2}.
\]
We do not know a single space of row functions of dimension $o(N)$ that
achieves this for all queries. Instead, for each probability law on the
queries we build a space adapted to it (Lemma~\ref{lem:x-sc-criterion}).
On a pair of small cells of rows and queries its polynomial degree grows
with the local densities of both, and Busemann's intersection inequality
(Lemma~\ref{lem:pre_busemann}) bounds both its dimension and its mean squared error through the exponent
$r/(r-1)$. A minimax argument gives one projection that serves every
query, and chaining over energy tiers gives the round bound
(Lemma~\ref{lem:x-sc-round}).

In part~(c), the upper bound uses the sparsification and
band-discrepancy theorems of Li (Lemmas~\ref{lem:pre_li_sparsification}
and~\ref{lem:pre_li_band}), and the lower bound
tests the spherical harmonics of one degree.

\begin{theorem}[High accuracy for non-even exponents]
\label{thm:high_accuracy_main}
Fix a non-even $p>0$ and let $\eta:=(r-1)/(r+2p)$.
\begin{enumerate}[label=(\alph*),nosep]
\item Fix $a>0$. For all sufficiently large $r$ and all $L\ge ar$,
\[
 (\epsilon^{-2}r^{1/2-p}\log(2r))^\eta
 \lesssim_{p,a}S_p(r,\epsilon)\lesssim_{p,a}
 (\epsilon^{-2}r^{1-p})^\eta ,
\]
with constants that are uniform as $L/r$ grows.
\item For sufficiently large $r$ and uniformly for
$L\ge C_p r\log(2r)$,
\[
 S_p(r,\epsilon)=\Theta_p(
   (\epsilon^{-2}r^{1/2-p}\log(2r))^\eta).
\]
\item For fixed rank $r\ge2$,
\[
 S_p(r,\epsilon)=\Theta_{r,p}(\epsilon^{-2(r-1)/(r+2p)})
\]
as $\epsilon\to0$.
\end{enumerate}
\end{theorem}
\begin{proof}
\emph{Part (a), upper bound.}
Constants $c,C$ with or without indices depend only on $p$ and $a$, and
$r\ge r_1(p,a)$ is large. Define $d:=r-1$. Zero rows do not matter, so we use
the Lewis setup Eq.~\eqref{lem:t-eq-lewis}, with unit vectors $u_1,\ldots,u_n
\in S^d$, a probability vector $\nu$ and the energy $F$. Write
$\Phi_x(\mu):=\sum_i\mu_i|\langle u_i,x\rangle|^p$. It suffices to find a
probability vector $\mu^*$ with $\supp\mu^*\subseteq\supp\nu$, at most $M$
nonzero entries, where $M$ is defined below, and
$|\Phi_x(\mu^*)-F(x)|\le\epsilon F(x)$ for unit $x$. Homogeneity then gives
all $x$, and the weights $\mu^*_i/\nu_i$ on the original rows form a
sparsifier.

\emph{Degree and stopping size.} Let $c_a$ be as in
Lemma~\ref{lem:x-ha-binom} for this $a$, define $c_q:=c_a/2$, and let
\[
 M:=\lceil C_M(\epsilon^{-2}r^{1-p})^\eta\rceil,\qquad
 q:=\lfloor c_qrM^{1/d}\rfloor,
\]
where $C_M\ge1$ is fixed at the end in terms of $p$ and $a$. The following
hold for every $C_M\ge1$ once $r$ is large.
\begin{enumerate}[label=(\alph*),nosep,leftmargin=*]
\item $\log M\ge ad$. Indeed $\epsilon^{-2}=e^{2L}/4\ge e^{2ar}/4$, so
$\log M\ge\eta(2ar-\log4-|p-1|\log r)$. Here $\eta\ge3/4$ for $r\ge6p+4$,
and $\log4+|p-1|\log r\le2ar/3$ for large $r$, so $\log M\ge ar$.
\item $q\ge c_0r$ with $c_0:=c_q/2$, since $c_qrM^{1/d}\ge c_qr\ge2$.
\item $10e^{ar/8}\binom{q+d}d\le\delta_0M$. Indeed
$q\le2c_qdM^{1/d}=c_adM^{1/d}$ as $r\le2d$, so
Lemma~\ref{lem:x-ha-binom} and (a) give $\binom{q+d}d\le e^{-ad/2}M$.
Moreover $ar/8-ad/2\le-ad/4$ and $10e^{-ad/4}\le\delta_0$.
\item $M\ge\max\{20,5n_{\rm c},5e^{ar/8}/\delta_0\}$, by (a).
\item $d\log(eN)\le N^{1-1/d}$ for every $N\ge M$. The function
$N\mapsto N^{1-1/d}/\log(eN)$ increases for $N\ge3$ when $d\ge2$, and at
$N=M$ it is at least $e^{y/2}/(1+y)\ge y/a\ge d$, where $y:=\log M\ge ad$
is large.
\end{enumerate}

\emph{The residual system.} We now define $\mathcal C$, $X_c$, $F_c$,
$G_c$, $P_{c,g}$ and $\beta_{c,g}$ satisfying Eq.~\eqref{lem:x-ha-eq-split}
with $V_0:=C_Vr^{p/2}q^{-p}$, such that $\log|\mathcal C|\le ar/8$, every
unit $x$ lies in some $X_c$, and $F\ge F_c$ on $X_c$.

Let $0<p<2$. Apply Lemma~\ref{lem:t-cap-family}(ii) with $a_0:=a/8$. For
large $r$ it gives a family $\mathcal C$ of caps with $\log|\mathcal C|\le
ar/8$ whose query sets $X_c$ cover $S^d$, with $\|v_c\|_2\le2$,
$r^{1/2-1/p}\le\varrho_c\le1$ and $r\varrho_c^p/t_c\le C_{p,a_0}r^{p/2}$.
Define $F_c:=t_c/(2r)$, so $F\ge F_c$ on $X_c$. For each cap take the groups
$G_c$, the polynomials $P_{c,g}$ and the measures $\beta_{c,g}$ of
Lemma~\ref{lem:t-shifted-tv}(d) for $v=v_c$ and $\varrho=\varrho_c$. Since
$|\langle u,v_c\rangle|\le2$ for unit $u$, the groups partition $S^d$, and
each is the intersection of $S^d$ with at most two affine tests. The
identity in Eq.~\eqref{lem:x-ha-eq-split} holds for $u\in g$ and $x\in X_c$,
because $\|x-v_c\|_2\le\varrho_c$. By (b) and the last assertion
of Lemma~\ref{lem:t-shifted-tv}(d) with $b_0:=1/p-1/2$,
\[
 \sum_g\|\beta_{c,g}\|\le(3^pC_p+1)\varrho_c^pq^{-p}
 =(3^pC_p+1)\frac{r\varrho_c^p}{t_c}\,2F_c\,q^{-p}\le C_Vr^{p/2}q^{-p}F_c .
\]

Let $p>2$. Take a single member $c$ with $X_c:=S^d$, $F_c:=r^{-p/2}$ and
$G_c:=\{S^d\}$, the intersection of $S^d$ with no test. The Lewis setup
gives $F\ge F_c$ on $S^d$. Let $P$ be the polynomial of
Lemma~\ref{lem:t-tv-approx} for $(v,\sigma)=(p,0)$ and degree $q$, define
$R:=|t|^p-P$, $P_{c,S^d}:=P+R(-1)$, and let $\beta_{c,S^d}$ be the measure
$R'(s)\,\d s$ on $(-1,1]$. Lemma~\ref{lem:t-threshold-form} gives
$|t|^p=P_{c,S^d}(t)+\int\mathbf1_{\{t\ge s\}}\,\d\beta_{c,S^d}(s)$ for
$|t|\le1$, which is Eq.~\eqref{lem:x-ha-eq-split} since $|\langle u,x\rangle|
\le1$, and $\|\beta_{c,S^d}\|=\operatorname{TV}_{[-1,1]}(R)\le C_pq^{-p}
=C_pr^{p/2}q^{-p}F_c$.

\emph{Rounding.} By (c), (d) and (e), every probability vector $\mu$ with
$N:=|\supp\mu|>M$ satisfies the hypotheses of Lemma~\ref{lem:x-ha-round}
for this residual system, since $|\mathcal C|\le e^{ar/8}$ and $N\ge M$.
Apply Lemma~\ref{lem:t-iterated-coloring} with $I:=[n]$, $\mathcal S$ the
set of probability vectors on $[n]$, $\mathcal X:=S^d$, the functions
$\Phi_x$, $\omega(x):=F_{c(x)}$ for some $c(x)$ with $x\in X_{c(x)}$,
$N_*:=M$, $b_*:=\frac12+\frac1{2d}$ and
$a_*:=C_{\rm R}e^{a/8}\sqrt d\,V_0$. Lemma~\ref{lem:x-ha-round} supplies its
hypothesis, because $|\mathcal C|^{1/(2d)}\le e^{ar/(16d)}\le e^{a/8}$.
Starting from $\nu$ we obtain $\mu^*$ with at most $M$ nonzero entries,
$\supp\mu^*\subseteq\supp\nu$ and, as $\omega\le F$ and $b_*\ge1/2$,
\[
 |\Phi_x(\mu^*)-F(x)|\le C_5\sqrt d\,V_0M^{-1/2-1/(2d)}F(x),\qquad
 C_5:=\frac{C_{\rm R}e^{a/8}}{1-(1-\delta_1/10)^{1/2}} .
\]
By (b), $q\ge c_qrM^{1/d}/2$, so $V_0\le C_V(2/c_q)^pr^{-p/2}M^{-p/d}$.
Since $\frac12+\frac{p+1/2}d=\frac{r+2p}{2(r-1)}=\frac1{2\eta}$, the
relative error is at most
\[
 C_6\,r^{(1-p)/2}M^{-1/(2\eta)},\qquad C_6:=C_5C_V(2/c_q)^p .
\]
As $M\ge C_M(\epsilon^{-2}r^{1-p})^\eta$ and $1/(2\eta)\ge1/2$, we have
$M^{1/(2\eta)}\ge C_M^{1/2}\epsilon^{-1}r^{(1-p)/2}$. With $C_M:=\max\{1,C_6^2\}$ the
relative error is at most $\epsilon$. Hence $S_p(r,\epsilon)\le M\le
2C_M(\epsilon^{-2}r^{1-p})^\eta$, and no constant depends on $L\ge ar$.

\emph{Part (a), lower bound.}
Lemma~\ref{lem:s-onetag}(a) applies to the same $p$ and $a$. Its quantity
$R_{\rm 1t}=(\epsilon^{-2}r^{1/2-p}\log(2r))^\eta$ has the same exponent
$\eta=(r-1)/(r+2p)$ as the theorem. For $r\ge r_0$ and every $L\ge ar$ it
gives a rank-$r$ input for which every reweighting with real weights and
relative error at most $\epsilon$ has more than $c_{\rm 1t}R_{\rm 1t}$
nonzero weights. Sparsifiers have positive weights, so
$S_p(r,\epsilon)>c_{\rm 1t}(\epsilon^{-2}r^{1/2-p}\log(2r))^\eta$, with
$c_{\rm 1t}$ and $r_0$ depending only on $p$ and $a$.

\emph{Part (b).}
\emph{Lower bound.} Lemma~\ref{lem:s-onetag}(a), applied with its
parameter $a$ equal to $1$, gives $c_{\rm 1t}>0$ and $r_0$, depending only
on $p$, such that $S_p(r,\epsilon)>c_{\rm 1t}R_{\rm 1t}$ for every non-even
$p>0$, $r\ge r_0$ and $L\ge r$. Since $C_p\log(2r)\ge1$, the range
$L\ge C_pr\log(2r)$ lies in $L\ge r$.

\emph{Upper bound.} The upper bound holds for every
$p\in(0,\infty)\setminus\{2\}$. Let $C_{\rm rd}$ be the constant of
Lemma~\ref{lem:x-sc-round} and define
\[
 C_{\rm sp}:=\max\{1,(\frac{20}{\delta_1})^2C_{\rm rd}^2
 (\frac2{c_q})^{2p+1}\},
\]
which depends only on $p$. Let $C_{\rm end}$ and $r_0$ be as in the
standing assumptions for this $C_{\rm sp}$, and let $r\ge r_0$ and
$L\ge C_{\rm end}r\log(2r)$. Thus the constant $C_p$ of the theorem is
$C_{\rm end}$. Zero rows do not affect the energy, so we may use the
Lewis setup Eq.~\eqref{lem:t-eq-lewis} with rows $i\in I$, directions $u_i$ and
masses $\nu_i$. If $|I|\le N$ there is nothing to prove. Otherwise apply
Lemma~\ref{lem:t-iterated-coloring} with the set $\mathcal S$ of
probability vectors on $I$, $\mathcal X:=\R^r\setminus\{0\}$,
$\Phi_x(\mu):=\sum_i\mu_i|\langle u_i,x\rangle|^p$, the weight function
$F$, $N_*:=N$, $b_*:=1/2$ and $a_*:=C_{\rm rd}(r^{p+1}\Lambda q^{-2p}
\gamma)^{1/2}$. For $\mu\in\mathcal S$ with $n:=|\supp\mu|>N$, take as
units the singletons of $J:=\{i:0<\mu_i\le2/n\}$, so $|J|\ge n/2\ge n/5$,
and take $H$ and $\mathcal K$ from Lemma~\ref{lem:x-sc-round} with $w:=\mu$.
For $z\in\mathcal K\cap[-1,1]^J$ the vector $\mu^z$, with
$\mu^z_i=\mu_i(1+z_i)$ on $J$, is nonnegative and has total mass one, and
$\Phi_x(\mu^z)-\Phi_x(\mu)=\sum_{i\in J}\mu_iz_i|\langle u_i,x\rangle|^p$
is bounded by Eq.~\eqref{lem:x-sc-eq-round}. Also $N\ge\max\{5n_{\rm c},20\}$.
The lemma gives a probability vector $\mu^*$ with
$\supp\mu^*\subseteq\supp\nu$, $|\supp\mu^*|\le N$ and
\[
 |\Phi_x(\mu^*)-F(x)|\le\frac{a_*N^{-1/2}}{1-(1-\delta_1/10)^{1/2}}F(x)
 \le\frac{20a_*}{\delta_1\sqrt N}F(x)\qquad(x\ne0).
\]
By Lemma~\ref{lem:x-sc-params}(b), $q\ge c_qrY/2$, so with
$1/\eta=1+(2p+1)/d$,
\[
 \frac{a_*^2}N=C_{\rm rd}^2\frac{r^{p+3/2}\Lambda}{q^{2p+1}N}\le
 C_{\rm rd}^2(\frac2{c_q})^{2p+1}r^{1/2-p}\Lambda N^{-1/\eta}
 \le C_{\rm rd}^2(\frac2{c_q})^{2p+1}C_{\rm sp}^{-1}\epsilon^2,
\]
using $N\ge C_{\rm sp}R_{\rm 1t}$, $C_{\rm sp}^{-1/\eta}\le C_{\rm sp}^{-1}$
and $R_{\rm 1t}^{-1/\eta}=\epsilon^2r^{p-1/2}\Lambda^{-1}$. By the choice of
$C_{\rm sp}$ the relative error is at most $\epsilon$. The weights
$\mu^*_i/\nu_i$ on the original rows preserve $\|Ax\|_p^p$ within the same
factors, by the last sentence of the subsection on the Lewis setup. Hence
$S_p(r,\epsilon)\le N\le2C_{\rm sp}R_{\rm 1t}$.

\emph{Part (c).}
Constants in this part may depend on $r$ and $p$, and $\eta=(r-1)/(r+2p)$.

\emph{Upper bound.} Zero rows do not matter, so we may use the Lewis setup
Eq.~\eqref{lem:t-eq-lewis}. It suffices to find a probability vector $\mu$ with
$\supp\mu\subseteq\supp\nu$, at most $C_{r,p}\epsilon^{-2\eta}$ nonzero
entries and $|\sum_i\mu_i|\langle u_i,x\rangle|^p-F(x)|\le\epsilon F(x)$ for
unit $x$. By homogeneity this extends to all $x$, and the weights
$\mu_i/\nu_i$ on the original rows form a sparsifier. For unit $x$ the Lewis
setup gives $F(x)\ge\varphi_r:=r^{-\max\{1,p/2\}}$, so an additive error
$\epsilon\varphi_r$ on the unit sphere suffices.

For $p\ge1$ we apply Lemma~\ref{lem:pre_li_sparsification} with
$d=r$, $\epsilon'=\epsilon\varphi_r$ and the law of $\nu$, after merging
coinciding directions. This gives $O_{r,p}(\epsilon^{-2\eta})$ rows. For
$r\ge r_0(p)$ the same bound also follows from
Theorem~\ref{thm:high_accuracy_main}(a) with $a=1$, whose proof does not
use Lemma~\ref{lem:pre_li_sparsification}, because $L\ge r$ for small $\epsilon$.

For $0<p<1$ we combine Lemma~\ref{lem:t-vc-variation} with
Lemma~\ref{lem:pre_li_band}, whose constants are $c_n,C_n$ and whose
bands are $B_{x,s}$.
Let $n:=r-1$, let $\mu^0$ be the law of $\nu$ on $S^{r-1}$, and for an
integer $q\ge1$ let $V$ be the space of restrictions to $S^{r-1}$ of
polynomials of degree at most $q$ on $\R^r$. A monomial of degree $k\le q$
agrees on the sphere with its products by powers of $\|u\|_2^2$, so $V$ is
spanned by restrictions of homogeneous polynomials of degrees $q$ and
$q-1$. Hence $\dim V\le2\binom{q+r-1}{r-1}\le2(q+1)^{r-1}$. Fix $\gamma>0$
with $2(2\gamma)^{r-1}\le c_n$. For an integer $N\ge\gamma^{-(r-1)}$ define
$q:=\lfloor\gamma N^{1/(r-1)}\rfloor$, so that $q\ge1$,
$q\ge\gamma N^{1/(r-1)}/2$ and $\dim V\le2(2q)^{r-1}\le c_nN$.
Lemma~\ref{lem:t-vc-variation} gives an even polynomial $P$ of degree at
most $q$ and a signed measure $\beta$ on $(0,1]$ with $\|\beta\|\le
Cq^{-p}$ and $|t|^p=P(t)+\int\mathbf1_{\{|t|>s\}}\,\d\beta(s)$ for $|t|\le1$.
Take $\xi:=\beta$ and let $\mu^1$ be as above. For unit $x$ the function
$u\mapsto P(\langle x,u\rangle)$ lies in $V$. Since $\mu^1$ and $\mu^0$
both have mass one, Fubini's theorem gives
\[
 \int|\langle x,u\rangle|^p\,\d(\mu^1-\mu^0)(u)
 =\int(\mu^1-\mu^0)(\{|\langle x,u\rangle|>s\})\,\d\beta(s)
 =-\int(\mu^1(B_{x,s})-\mu^0(B_{x,s}))\,\d\beta(s).
\]
Its absolute value is at most
\[
 C_nCN^{-1/2-1/(2(r-1))}q^{-p}\le C'_{r,p}N^{-1/2-1/(2(r-1))-p/(r-1)}
 =C'_{r,p}N^{-(r+2p)/(2(r-1))},
\]
which is at most $\epsilon\varphi_r$ for $N:=\lceil C''_{r,p}
\epsilon^{-2\eta}\rceil$. This proves the upper bound for $0<p<1$.

\emph{Lower bound.} Let $\sigma$ be the rotation-invariant probability
measure on $S^{r-1}$. For a finite signed Borel measure $\rho$ on
$S^{r-1}$ define $(K\rho)(x):=\int|\langle x,u\rangle|^p\,\d\rho(u)$, and let
$\beta:=(K\sigma)(x)>0$, which does not depend on the unit vector $x$. By
Lemma~\ref{lem:pre_harmonics}(b) the operator $g\mapsto K(g\,\sigma)$ acts on
spherical harmonics of degree $k$ as multiplication by a number
$\lambda_k$, with $\lambda_0=\beta$. By Lemma~\ref{lem:pre_cosine_multipliers}, the multiplier of an even
degree $2j$ is a multiple of $(-1)^j\Gamma(j-p/2)/\Gamma(j+(r+p)/2)$ that
does not depend on $j$, so
\[
 \frac{\lambda_{2j}}\beta=(-1)^j\frac{(-p/2)_j}{((r+p)/2)_j},\qquad
 (y)_j:=y(y+1)\cdots(y+j-1).
\]
No factor of $(-p/2)_j$ vanishes, so $\lambda_{2j}\ne0$. Moreover
$|\lambda_{2j+2}/\lambda_{2j}|=|j-p/2|/(j+(r+p)/2)<1$, so $|\lambda_{2j}|$
decreases in $j$. By Lemma~\ref{lem:pre_gamma_ratio},
$\Gamma(j-p/2)/\Gamma(j+(r+p)/2)$ is asymptotic to $j^{-(r+2p)/2}$, hence
there is $c_1=c_1(r,p)>0$ with $|\lambda_{2j}|/\beta\ge c_1j^{-(r+2p)/2}$
for every $j\ge1$.

\emph{The input.} Partition $S^{r-1}$ into finitely many Borel cells of
positive measure and diameter less than $\delta<1$, choose a point $u_c$ in
each cell, and let $\nu:=\sum_c\sigma(\text{cell }c)\,\delta_{u_c}$. As
$|t|^p$ is uniformly continuous on $[-1,1]$, for small $\delta$ we have
$|K\nu(x)-\beta|\le\epsilon\beta$ for all unit $x$. The points $u_c$ span
$\R^r$, since a unit normal of a hyperplane containing them would be at
distance $\sqrt2>\delta$ from all of them. So the rows
$\sigma(\text{cell }c)^{1/p}u_c^\top$ form a rank-$r$ input with energy
$K\nu$. Let $(S,w)$ be a sparsifier of it and let $\mu$ be the positive
measure giving $u_c$ the mass $w_c\sigma(\text{cell }c)$ for $c\in S$. Then
$|K\mu-K\nu|\le\epsilon K\nu\le\epsilon(1+\epsilon)\beta$, and so
$|K(\mu-\sigma)|\le(2+\epsilon)\epsilon\beta\le3\epsilon\beta$ on the unit
sphere.

\emph{A vanishing polynomial.} Let $m\ge1$, let $\mathcal P_m$ be the space
of homogeneous polynomials of degree $m$ on $\R^r$, and define
$D_m:=\dim\mathcal P_m=\binom{r+m-1}{m}$. Restriction to $S^{r-1}$ is
injective on $\mathcal P_m$ by homogeneity. Suppose $|S|<D_m/2$. The
subspace $W\subseteq\mathcal P_m$ of polynomials vanishing on
$\supp\mu$ has dimension $k\ge D_m-|S|>D_m/2$. Take an
$L_2(\sigma)$-orthonormal basis $f_1,\ldots,f_k$ of $W$ and extend it to
an orthonormal basis $f_1,\ldots,f_{D_m}$ of $\mathcal P_m$. The function
$\sum_{i\le D_m}f_i(x)^2$ does not depend on the basis, so it is
invariant under rotations because $\mathcal P_m$ is. Hence it is
constant, and integrating shows that it equals $D_m$. In particular
$\sum_{i\le k}f_i(x)^2\le D_m$. For $\zeta$ uniform on $S^{k-1}$ and fixed
$y\in\R^k$ we have $\E[\langle\zeta,y\rangle^4]=3\|y\|_2^4/(k(k+2))$. Hence
$Q_\zeta:=\sum_{i\le k}\zeta_if_i$ satisfies
\[
 \E[\int Q_\zeta^4\,\d\sigma]\le\frac{3D_m^2}{k(k+2)}<12 .
\]
Fix $Q:=Q_\zeta$ with $\int Q^4\,\d\sigma\le12$. Then $\int Q^2\,\d\sigma=1$
and $\int Q^2\,\d\mu=0$.

\emph{The contradiction.} The polynomial $Q^2$ lies in $\mathcal P_{2m}$.
By the classical decomposition of homogeneous polynomials, its
restriction to the sphere is $\sum_{j=0}^mY_{2j}$ with $Y_{2j}$ a spherical
harmonic of degree $2j$, and these are orthogonal in $L_2(\sigma)$. Define
$g:=\sum_{j=0}^m\lambda_{2j}^{-1}Y_{2j}$, so that $K(g\,\sigma)=Q^2$ on the
sphere and $\|g\|_{L_2(\sigma)}\le|\lambda_{2m}|^{-1}\|Q^2\|_{L_2(\sigma)}
\le\sqrt{12}/|\lambda_{2m}|$ by the monotonicity of $|\lambda_{2j}|$. By
Fubini's theorem and the symmetry of the kernel,
\[
 1=|\int Q^2\,\d(\mu-\sigma)|=|\int g(u)\,K(\mu-\sigma)(u)\,
 \d\sigma(u)|\le\|g\|_{L_2(\sigma)}\cdot3\epsilon\beta
 \le\frac{3\sqrt{12}\,\epsilon\beta}{|\lambda_{2m}|} .
\]
Thus $|\lambda_{2m}|/\beta<11\epsilon$. For
$m:=\lfloor(c_1/(11\epsilon))^{2/(r+2p)}\rfloor$, which is at least one
for small $\epsilon$, this contradicts $|\lambda_{2m}|/\beta\ge
c_1m^{-(r+2p)/2}\ge11\epsilon$. Hence every sparsifier has
$|S|\ge D_m/2\ge m^{r-1}/(2(r-1)!)\ge c_{r,p}\epsilon^{-2\eta}$.
\end{proof}

\subsection{One Round on Paired Light Rows}

\begin{lemma}[Binomial counts at high accuracy]\label{lem:x-ha-binom}
For every $a>0$ there is $c_a\in(0,1]$ such that for all integers $d\ge1$
and $q\ge0$ and every real $M\ge e^{ad}$ with $q\le c_adM^{1/d}$,
\[
 \binom{q+d}d\le e^{-ad/2}M .
\]
\end{lemma}

\begin{proof}
For integers $0\le k\le n$ the term $(k/n)^k((n-k)/n)^{n-k}\binom nk$ of the
binomial expansion of $1=(k/n+(n-k)/n)^n$ is at most one. With $n=q+d$ and
$k=q$ this gives $\binom{q+d}d\le e^{dh(q/d)}$ for
$h(\theta):=(1+\theta)\log(1+\theta)-\theta\log\theta$, with $h(0)=0$. The
function $h$ is continuous and increasing on $[0,\infty)$, since
$h'(\theta)=\log(1+1/\theta)>0$. Choose $c_a\in(0,1]$ with
$h(c_ae^a)\le a/2$. Define $\ell_M:=(\log M)/d\ge a$ and
$g(s):=s-h(c_ae^s)$. Then $g'(s)=1-\vartheta\log(1+1/\vartheta)$ with
$\vartheta:=c_ae^s$, which is positive because $\log(1+y)<y$. Hence
$g(\ell_M)\ge g(a)\ge a/2$. Since $q/d\le c_ae^{\ell_M}$ and $h$ increases,
$\binom{q+d}d\le e^{dh(c_ae^{\ell_M})}=e^{d\ell_M-dg(\ell_M)}\le Me^{-ad/2}$.
\end{proof}

The next lemma is one round of the rounding. Its residual systems carry,
for each member $c$ of a finite family, a set $X_c$ of unit queries and a
number $F_c>0$, which in the applications is a lower bound for the
energy on $X_c$. Here $\delta_0,\delta_1,n_{\rm c}$ are the constants of
Lemma~\ref{lem:t-subspace-coloring}, units and $\mu^z$ are as defined
before Lemma~\ref{lem:t-iterated-coloring}, and light rows are as in
Lemma~\ref{lem:t-light-rows}.

\begin{lemma}[One round on paired light rows]\label{lem:x-ha-round}
Let $p>0$, $r\ge3$, $d:=r-1$, let $u_1,\ldots,u_n\in S^d$, let $q\ge1$
be an integer and $V_0>0$, and let $\mathcal C$ be a nonempty finite set.
For each $c\in\mathcal C$ let $X_c\subseteq S^d$, let $F_c>0$, and let
$G_c$ be a finite partition of $S^d$ into sets that are intersections of
$S^d$ with at most two affine tests on $\R^r$. For $g\in G_c$ let
$P_{c,g}$ be a polynomial of degree at most $q$ and $\beta_{c,g}$ a finite
signed Borel measure on $\R$ such that
\begin{equation}\label{lem:x-ha-eq-split}
 |\langle u,x\rangle|^p=P_{c,g}(\langle u,x\rangle)+\int\mathbf1_{\{\langle
 u,x\rangle\ge s\}}\,\d\beta_{c,g}(s)\quad(u\in g,\ x\in X_c),\qquad
 \sum_{g\in G_c}\|\beta_{c,g}\|\le V_0F_c .
\end{equation}
Let $\mu$ be a probability vector on $[n]$ with $N:=|\supp\mu|$ such that
\[
 N\ge\max\{20,5|\mathcal C|/\delta_0\},\qquad d\log(eN)\le N^{1-1/d},\qquad
 10|\mathcal C|\binom{q+d}d\le\delta_0N .
\]
Then there are a set $\mathcal P$ of pairwise disjoint pair units of $\mu$
with $|\mathcal P|\ge N/5$, a subspace $H\subseteq\R^{\mathcal P}$ of
codimension at most $\delta_0|\mathcal P|$ and a closed symmetric convex
set $\mathcal K\subseteq H$ with $\gamma_H(\mathcal K)\ge
e^{-\delta_0|\mathcal P|}$, such that for every $z\in\mathcal K\cap
[-1,1]^{\mathcal P}$ the vector $\mu^z$ is a probability vector and
\[
 |\sum_i(\mu^z_i-\mu_i)|\langle u_i,x\rangle|^p|\le C_{\rm R}
 |\mathcal C|^{1/(2d)}\sqrt d\,V_0N^{-1/2-1/(2d)}F_c\qquad(c\in\mathcal C,\
 x\in X_c),
\]
where $C_{\rm R}$ is absolute.
\end{lemma}

\begin{proof}
\emph{Pairs.} Let $J:=\{i:0<\mu_i\le2/N\}$, so $|J|\ge N/2$ by
Lemma~\ref{lem:t-light-rows}. Apply Lemma~\ref{lem:t-low-crossing} to the
points $u_i$, $i\in J$, with constant $C_{\rm lc}$, and orient each of its
pairs arbitrarily. This gives a set $\mathcal P$ of
$\lfloor|J|/2\rfloor\ge N/5$ disjoint pair units such that every
intersection of at most three affine tests crosses at most
\[
 3C_{\rm lc}(|J|^{1-1/d}+d\log(e|J|))\le6C_{\rm lc}N^{1-1/d}=:\chi
\]
of the pairs, by the hypothesis $d\log(eN)\le N^{1-1/d}$. For
$e=(i,j)\in\mathcal P$ we have $\tau_e=\min\{\mu_i,\mu_j\}\le2/N$ and
$\sum_e\tau_e\le\frac12\sum_{i\in J}\mu_i\le\frac12$, so
$\|\tau\|_2^2\le(2/N)\sum_e\tau_e\le1/N$. Pair units preserve the total
mass, so $\mu^z$ is a probability vector for every $z\in[-1,1]^{\mathcal P}$.

\emph{Polynomial constraints.} Let $\mathcal P_q(S^d)$ be the space of
restrictions to $S^d$ of polynomials of degree at most $q$ on $\R^r$. A
monomial of degree $k\le q$ agrees on $S^d$ with its products by powers of
$\|x\|_2^2$, so $\mathcal P_q(S^d)$ is spanned by restrictions of
homogeneous polynomials of degrees $q$ and $q-1$, and $\dim\mathcal
P_q(S^d)\le2\binom{q+d}d$. For $u\in S^d$ let $g_c(u)$ be the group of $G_c$
containing $u$. For $z\in\R^{\mathcal P}$ let $L_c(z)$ be the restriction
to $S^d$ of
\[
 x\mapsto\sum_{e=(i,j)\in\mathcal P}\tau_ez_e(P_{c,g_c(u_i)}(\langle
 u_i,x\rangle)-P_{c,g_c(u_j)}(\langle u_j,x\rangle)).
\]
Then $L_c$ is a linear map into $\mathcal P_q(S^d)$. Let $H$ be the
intersection of the kernels of the $L_c$, $c\in\mathcal C$. Its
codimension is at most $2|\mathcal C|\binom{q+d}d\le\delta_0N/5\le
\delta_0|\mathcal P|$.

\emph{The error as an integrated vector.} Fix $c\in\mathcal C$. On
$\mathcal Z_c:=G_c\times\R$, with the discrete $\sigma$-algebra on $G_c$,
let $\beta_c$ be the signed measure with $\beta_c(\{g\}\times B):=
\beta_{c,g}(B)$, so that $\|\beta_c\|=\sum_g\|\beta_{c,g}\|\le V_0F_c$. For
$e=(i,j)\in\mathcal P$, $(g,s)\in\mathcal Z_c$ and $x\in\R^r$ define
\[
 \phi_c(e,(g,s),x):=\mathbf1_{\{u_i\in g,\ \langle u_i,x\rangle\ge s\}}
 -\mathbf1_{\{u_j\in g,\ \langle u_j,x\rangle\ge s\}}\in\{-1,0,1\}.
\]
For fixed $(e,(g,s))$ this is a function of the indicators of the two
affine tests $\{x:\langle u_i,x\rangle\ge s\}$ and
$\{x:\langle u_j,x\rangle\ge s\}$. Since the groups of $G_c$ partition
$S^d$, Eq.~\eqref{lem:x-ha-eq-split} says that for $u\in S^d$ and $x\in X_c$,
\[
 |\langle u,x\rangle|^p=P_{c,g_c(u)}(\langle u,x\rangle)+\int\mathbf1_{\{u\in
 g,\ \langle u,x\rangle\ge s\}}\,\d\beta_c(g,s).
\]
Hence, for $z\in\R^{\mathcal P}$ and $x\in X_c$,
\[
 \sum_i(\mu^z_i-\mu_i)|\langle u_i,x\rangle|^p
 =\sum_{e=(i,j)}\tau_ez_e(|\langle u_i,x\rangle|^p-|\langle u_j,x
 \rangle|^p)=L_c(z)(x)+\sum_e\tau_ez_e\int\phi_c(e,\cdot,x)\,\d\beta_c .
\]
For $z\in H$ the first term vanishes, and the second equals
$F_c\langle z,w^c_x\rangle$, where
\[
 w^c_x:=\Pi_H(\tau_e\int\phi_c(e,\cdot,x)\,\d(\beta_c/F_c))_{e\in
 \mathcal P}.
\]
We apply Lemma~\ref{lem:t-integrated-vectors} to the measures
$\beta_c/F_c$, of total variation at most $V_0$, with $\ell=2$,
$X=X_c\subseteq S^d$ and $\Pi=\Pi_H$. Its quantity $D$ is at most
$V_0\|\tau\|_2\le V_0N^{-1/2}$. For fixed $(g,s)$ and $x$, the value
$\phi_c(e,(g,s),x)$ is nonzero exactly when the set
$g\cap\{u:\langle u,x\rangle\ge s\}$, an intersection of $S^d$ with at most
three affine tests, crosses the pair $\{u_i,u_j\}$. So at most $\chi$
indices $e$ give a nonzero value, and part~(c) of the lemma gives
\[
 \|w^c_x\|_2\le V_0\cdot\frac2N\cdot\sqrt\chi=C_2V_0N^{-1/2-1/(2d)},
 \qquad C_2:=2(6C_{\rm lc})^{1/2}.
\]
Part~(b) gives $N(\{w^c_x:x\in X_c\},u)\le\max\{1,(C_3V_0N^{-1/2}/u)^{2d}\}$
for every $u>0$, where $C_3$ is its constant for $\ell=2$.

\emph{Chaining.} Let $\mathcal V\subseteq H$ be the union over
$c\in\mathcal C$ of the sets $\{w^c_x:x\in X_c\}$. Then $N(\mathcal V,u)\le
|\mathcal C|\max\{1,(C_3V_0N^{-1/2}/u)^{2d}\}$. Define $\kappa_0:=\delta_0/5$,
$C_{\rm ch}:=\max\{C_2,C_3\kappa_0^{-1/2}\}$ and
\[
 \rho:=C_{\rm ch}|\mathcal C|^{1/(2d)}V_0N^{-1/2-1/(2d)} .
\]
Then $\sup_{v\in\mathcal V}\|v\|_2\le\rho$, and for $j\ge1$
\[
 N(\mathcal V,2^{-j/2}\rho)\le\max\{|\mathcal C|,(C_3/C_{\rm ch})^{2d}
 N2^{jd}\}\le\kappa_0N2^{dj},
\]
because $|\mathcal C|\le\kappa_0N$ by hypothesis and
$(C_3/C_{\rm ch})^{2d}\le(C_3/C_{\rm ch})^2\le\kappa_0$.
Lemma~\ref{lem:t-sidak-chaining}(c) with $n_0:=N$ gives countably many
strips $\{y\in H:|\langle y,w_l\rangle|\le c_l\}$, each of Gaussian measure
at least $1-e^{-2}$ and with total failure probability at most
$\kappa_0N/2$, such that every point $y$ of their intersection
$\mathcal K$ satisfies $\sup_{v\in\mathcal V}|\langle y,v\rangle|\le
17\sqrt d\,\rho$. The bound holds at every point of $\mathcal K$, since the
proof of that lemma bounds every realization. The set $\mathcal K\subseteq
H$ is closed, symmetric and convex, and part~(d) gives
$\gamma_H(\mathcal K)\ge e^{-\kappa_0N}=e^{-\delta_0N/5}\ge
e^{-\delta_0|\mathcal P|}$. For $z\in\mathcal K\cap[-1,1]^{\mathcal P}$,
$c\in\mathcal C$ and $x\in X_c$ the error is
$F_c|\langle z,w^c_x\rangle|\le17\sqrt d\,\rho F_c$. This proves the
lemma with $C_{\rm R}:=17C_{\rm ch}$.
\end{proof}

No step of this proof uses a triangle inequality for the $p$-energy, so
it applies for $0<p<1$ as well.

\subsection{Adaptive Shifts}

\begin{lemma}[Adaptive shifts at exponential and subexponential accuracy]
\label{lem:joint-adaptive-shifts}
Fix a non-even $p>0$. On any compact interval
$I\subset(0,\infty)$, for $L/r\in I$ and sufficiently large $r$,
\[
 S_p(r,\epsilon)\gtrsim_{p,I}\epsilon^{-2}r^{1/2-p}\log r.
\]
Uniformly for $r^{2/(p+2)}\le L\le r/4$,
\[
 S_p(r,\epsilon)\gtrsim_p
 \epsilon^{-2}\frac{r^2}{L^{p+3/2}}
 (1+\log\frac{L^{1+p/2}}{r}).
\]
\end{lemma}

\begin{proof}
Both bounds are instances of Lemma~\ref{lem:s-shift}, whose conclusion
holds for real weights and hence for sparsifiers. In this proof $N$, $H$
and $d$ are the numbers of tags, of transverse copies and the transverse
dimension of that lemma, and $c'_S$ and $d_0$ are its constants.

\emph{First bound.} Write $I:=[\alpha_1,\alpha_2]$ with
$0<\alpha_1\le\alpha_2$. Define $N:=\lfloor r/2\rfloor$, $H:=1$,
$d:=r-N=\lceil r/2\rceil$ and $b:=1$, and apply Lemma~\ref{lem:s-shift}
with $\lambda_1:=\alpha_1<\lambda_2:=2\alpha_2$ and $\beta_0:=1$. Since
$d\le r\le2d$, we get $\lambda_1d\le\alpha_1r\le L\le\alpha_2r\le
\lambda_2d$. Also $Hb^p=1$, $b\sqrt d\ge1=\beta_0$, $d\ge d_0$ once
$r\ge2d_0$, and $r=N+Hd$. Hence
\[
 S_p(r,\epsilon)>c'_S\epsilon^{-2}Nd^{-p-1/2}(1+\tfrac12\log d)
 \ge c'_S\epsilon^{-2}\cdot\frac r3\cdot r^{-p-1/2}\cdot\tfrac12\log r ,
\]
using $N\ge r/3$, $d\le r$ and $1+\frac12\log(r/2)\ge\frac12\log r$.

\emph{Second bound.} Let $r^{2/(p+2)}\le L\le r/4$, define
$\beta_0:=2^{-(p+2)/(2p)}$, and let $d_0$ be the constant of
Lemma~\ref{lem:s-shift} for $\lambda_1:=1$, $\lambda_2:=2$ and this
$\beta_0$. We may take $d_0$ to be an integer with $d_0\ge2$. If $L<d_0$, then $r\le L^{(p+2)/2}<
d_0^{(p+2)/2}$ and $\epsilon^{-2}=e^{2L}/4<e^{2d_0}$, so the right side is
bounded in terms of $p$ while $S_p(r,\epsilon)\ge1$. Let $L\ge d_0$ and
define
\[
 N:=\lfloor r/2\rfloor,\qquad d:=\lfloor L\rfloor,\qquad
 H:=\lfloor\frac{r-N}d\rfloor,\qquad b:=H^{-1/p}.
\]
Then $d\ge d_0$ and $d\le L<d+1\le2d$, so $\lambda_1d\le L\le\lambda_2d$.
Since $r-N\ge r/2$ and $L\le r/4$, we have $(r-N)/d\ge2$, so $H\ge2$ and
$H\ge(r-N)/(2d)\ge r/(4d)$. Also $Hb^p=1$, $b\le1$ and $r\ge N+Hd$. As
$H\le(r-N)/d\le r/d$ and $d\ge L/2$,
\[
 b\sqrt d\ge(\frac{d^{1+p/2}}r)^{1/p}\ge
 2^{-(p+2)/(2p)}(\frac{L^{1+p/2}}r)^{1/p}\ge\beta_0 ,
\]
because $L^{1+p/2}\ge r$. Lemma~\ref{lem:s-shift} gives, with
$y:=L^{1+p/2}/r\ge1$,
\[
 S_p(r,\epsilon)>c'_S\epsilon^{-2}NHd^{-p-1/2}(1+\log_+(b\sqrt d))
 \ge\frac{c'_S}{12}\,\epsilon^{-2}\frac{r^2}{L^{p+3/2}}
 (1+\log_+(\beta_0y^{1/p})),
\]
where we used $NH\ge r^2/(12d)$ and $d\le L$. Finally
$1+\log_+(\beta_0y^{1/p})\ge c_p(1+\log y)$ with $c_p>0$ depending only on
$p$. Indeed, define $\kappa_p:=-\log\beta_0$. If $\log y\le2p\kappa_p$, the left
side is at least $1\ge(1+\log y)/(1+2p\kappa_p)$. Otherwise it is at least
$1+(\log y)/(2p)\ge\min\{1,1/(2p)\}(1+\log y)$.
\end{proof}

Combining the upper bound of Theorem~\ref{thm:high_accuracy_main}(a) with Lemma~\ref{lem:joint-adaptive-shifts}
gives, on any compact positive range $I$ of $L/r$,
\[
 c_{p,I}\epsilon^{-2}r^{1/2-p}\log r\le S_p(r,\epsilon)
     \le C_{p,I}\epsilon^{-2}r^{1-p}.
\]
Indeed $\epsilon^{-2}r^{1-p}\ge1$ there for large $r$, so its power $\eta<1$
is at most itself. Uniformly for $L\ge ar$, the ratio of the bounds in
Theorem~\ref{thm:high_accuracy_main}(a) is
$O_{p,a}((\sqrt r/\log(2r))^\eta)$.

\subsection{Regular Rows}

\begin{lemma}[Matching bounds for regular rows at high accuracy]
\label{lem:regular-high-accuracy}
Under the hypotheses of Theorem~\ref{thm:high_accuracy_main}(a), define
\[
 R:=(\epsilon^{-2}r^{1/2-p}\log(2r))^\eta,\qquad
 M:=\lceil C R\rceil,\qquad q:=\lfloor c(r-1)M^{1/(r-1)}\rfloor .
\]
Fix a sufficiently small $c_*>0$ and a sufficiently large $A_*$.
Then choose $c\in(0,1]$ small in terms of $p$ and $a$, and then $C\ge1$
large in terms of $c$, $c_*$, $A_*$, $p$, and $a$. Consider rank-$r$ inputs that in some coordinates
have nonzero rows $\ell_iu_i^\top$, with $\ell_i>0$ and unit vectors
$u_i\in\R^r$, such that the probability law $\nu$ giving row $i$ the mass
$\ell_i^p/\sum_k\ell_k^p$ satisfies, for every unit query $x$,
\[
 \E_\nu[|u^\top x|^p]\ge c_*r^{-p/2},\qquad
 \E_\nu[(1+q|u^\top x|)^{-2}]\le A_*\sqrt r/q .
\]
Here ``in some coordinates'' means after restricting the queries to the
row span and an invertible linear change of coordinates on it. The worst-case
original-row support in this class is $\Theta_{p,a,c_*,A_*}(R)$.
\end{lemma}

Equivalently, up to these operations, the class consists of the finite probability laws $\nu$ on
unit rows $u_i\in\R^r$ satisfying the two conditions, each standing for the
input with rows $\nu_i^{1/p}u_i^\top$. Indeed an invertible change of
coordinates $T$ acts on queries by $x\mapsto T^{-\top}x$ and preserves every
energy, a common factor on all rows multiplies every energy by the same
number, and neither changes which pairs $(S,w)$ are sparsifiers. With
$Z:=\sum_k\ell_k^p$, weights $w_i$ on the rows $\ell_iu_i^\top$ give the
energy $Z\sum_iw_i\nu_i|\langle u_i,x\rangle|^p$. So $(S,w)$ is a
sparsifier if and only if the measure $\mu_i:=w_i\nu_i$ on $S$ satisfies
$|\sum_i\mu_i|\langle u_i,x\rangle|^p-F(x)|\le\epsilon F(x)$ for every $x$,
where $F(x):=\E_\nu[|u^\top x|^p]$.

The upper bound replaces the variation residuals by Chebyshev tails of
$|t|^p$, which are small away from the cusp. Exact moments of degree $q$
preserve the second condition through all rounds, so no low-crossing
matching is needed. We use three facts.

\begin{lemma}[Chebyshev tails of $|t|^p$]\label{lem:x-ha-cheb}
Let $p>0$ be non-even and $\kappa:=\min\{p,1\}$. Let
$a_k:=\frac2\pi\int_0^\pi|\cos\phi|^p\cos(k\phi)\,\d\phi$ for integers
$k\ge0$. For an integer $q\ge p+2$ let $P_q$ be the even polynomial of
degree at most $q$ with $P_q(\cos\phi)=\frac{a_0}2+\sum_{1\le j\le q/2}
a_{2j}\cos(2j\phi)$, and define $b_q(t):=|t|^p-P_q(t)$. There is $C_{\rm b}$,
depending only on $p$, such that for $|t|\le1$ and all real $\phi,\phi'$
with $|\phi-\phi'|\le1/q$,
\[
 |b_q(t)|\le\frac{C_{\rm b}q^{-p}}{1+q|t|},\qquad
 |b_q(\cos\phi)-b_q(\cos\phi')|\le\frac{C_{\rm b}q^{-p}(q|\phi-\phi'|)^\kappa}
 {1+q|\cos\phi|} .
\]
\end{lemma}

\begin{proof}
The polynomial $P_q$ exists because $\cos(2j\phi)=T_{2j}(\cos\phi)$ for the
Chebyshev polynomials $T_{2j}$. In this proof $C$ depends only on $p$.

\emph{Coefficients.} Define $I_k:=\int_{-\pi/2}^{\pi/2}\cos^p\phi\cos(k\phi)
\,\d\phi$ and $I'_k:=\int_{-\pi/2}^{\pi/2}\cos^{p-1}\phi\cos(k\phi)\,\d\phi$,
which converge since $p>0$. For even $k$ the integrand of $a_k$ has period
$\pi$, so $a_k=\frac2\pi I_k$, and $a_k=0$ for odd $k$ by the
substitution $\phi\mapsto\pi-\phi$. For $k\ge1$ the
identity $\cos\phi\cos(k\phi)=\frac12(\cos((k-1)\phi)+\cos((k+1)\phi))$
gives $I_k=\frac12(I'_{k-1}+I'_{k+1})$. Integrating by parts, with
boundary terms vanishing as $\cos^p(\pm\pi/2)=0$, and using
$\sin\phi\sin(k\phi)=\frac12(\cos((k-1)\phi)-\cos((k+1)\phi))$, gives
$I_k=\frac p{2k}(I'_{k-1}-I'_{k+1})$. Hence $I'_{k-1}=(1+k/p)I_k$ and
$I'_{k+1}=(1-k/p)I_k$ for $k\ge1$. The first identity at $k+2$ and the
second at $k$ give $(p+k+2)I_{k+2}=(p-k)I_k$ for $k\ge1$, and the first at
$k=2$ with $I'_1=I_0$ gives the same relation for $k=0$. Thus
\[
 \frac{a_{2j+2}}{a_{2j}}=\frac{p-2j}{p+2j+2}\qquad(j\ge0),\qquad a_0>0 .
\]
As $p$ is not even, every $a_{2j}$ is nonzero. Let $j_0:=\lceil p/2\rceil$,
so $j_0>p/2$. For $j\ge j_0$ the ratio is negative, so
$(-1)^ja_{2j}=\varsigma c_j$ with $c_j:=|a_{2j}|$ and a fixed sign
$\varsigma$. For $j\ge j_0$ we have $c_{j+1}/c_j=(2j-p)/(2j+p+2)<1$ and
$(j+1)c_{j+1}/(jc_j)<1$, the latter because $(j+1)(2j-p)<j(2j+p+2)$. So
$c_j$ and $jc_j$ do not increase for $j\ge j_0$. With $\omega_p:=p/2+1$,
$1-x\le e^{-x}$ and $\sum_{i=j_0}^{j-1}(i+\omega_p)^{-1}\ge
\log((j+\omega_p)/(j_0+\omega_p))$,
\[
 c_j=c_{j_0}\prod_{i=j_0}^{j-1}(1-\frac{p+1}{i+\omega_p})
 \le c_{j_0}(\frac{j_0+\omega_p}{j+\omega_p})^{p+1}\qquad(j\ge j_0).
\]
Since also $|a_k|\le2$ for all $k$, we get $c_j\le Cj^{-p-1}$ for $j\ge1$.

\emph{The tail series.} The cosine series of the continuous even function
$|\cos\phi|^p$ has the coefficients $a_k$, which are absolutely summable,
so it converges uniformly to $|\cos\phi|^p$. Let $k_q:=\lfloor q/2\rfloor$.
Then $k_q\ge q/4$ and $k_q+1>q/2\ge p/2+1>j_0$, and for real $\phi$
\[
 b_q(\cos\phi)=\sum_{j>k_q}a_{2j}\cos(2j\phi)=\varsigma\sum_{j>k_q}c_j
 \cos(2j\psi),\qquad\psi:=\phi-\pi/2,
\]
since $\cos(2j\phi)=(-1)^j\cos(2j\psi)$. Note that $|\sin\psi|=|\cos\phi|$.
We use Abel summation in the following form. If $e_j\ge0$ does not
increase for $j>k$ and the partial sums $\sum_{j=k+1}^ny_j$ are bounded by
$B$ in absolute value for all $n$, then $|\sum_{j=k+1}^ne_jy_j|\le e_{k+1}B$
for all $n$. Moreover $|\sum_{j=k+1}^n\cos(2j\psi)|$ and
$|\sum_{j=k+1}^n\sin(2j\psi)|$ are at most $1/|\sin\psi|$ when
$\sin\psi\ne0$, as the sum of the $e^{2ij\psi}$ has modulus
$|\sin((n-k)\psi)/\sin\psi|$.

\emph{The size bound.} Absolute summation gives $|b_q(t)|\le\sum_{j>k_q}c_j
\le Ck_q^{-p}\le Cq^{-p}$. Abel summation gives $|b_q(\cos\phi)|\le
c_{k_q+1}/|\cos\phi|\le Cq^{-p-1}/|\cos\phi|$. Together,
$|b_q(t)|\le Cq^{-p}\min\{1,1/(q|t|)\}\le2Cq^{-p}/(1+q|t|)$.

\emph{Derivatives.} On the set where $\cos\phi\ne0$ the differentiated
series $-\varsigma\sum_{j>k_q}2jc_j\sin(2j\psi)$ converges locally
uniformly, by Abel summation and $jc_j\to0$. Hence it is the derivative of
$\phi\mapsto b_q(\cos\phi)$ there, and
\[
 |\frac{\d}{\d\phi}b_q(\cos\phi)|\le\frac{2(k_q+1)c_{k_q+1}}
 {|\cos\phi|}\le\frac{Cq^{-p}}{|\cos\phi|}\qquad(\cos\phi\ne0).
\]
If $p>1$, then $\sum_j2jc_j<\infty$, the differentiated series converges
uniformly on $\R$, and $|\frac{\d}{\d\phi}b_q(\cos\phi)|\le\sum_{j>k_q}2jc_j
\le Cq^{1-p}$ everywhere.

\emph{Increments.} Let $\theta:=|\phi-\phi'|\le1/q$ and $t:=\cos\phi$.
Suppose first $q|t|\ge2$. On the segment between $\phi$ and $\phi'$ we have
$|\cos|\ge|t|-\theta\ge|t|/2>0$, so the derivative bound gives
$|b_q(\cos\phi)-b_q(\cos\phi')|\le2Cq^{-p}\theta/|t|$, which is at most
$3Cq^{-p}(q\theta)/(1+q|t|)\le3Cq^{-p}(q\theta)^\kappa/(1+q|t|)$ because
$1+q|t|\le\frac32q|t|$ and $q\theta\le1$. Suppose now $q|t|<2$. Then
$1+q|t|<3$, and it suffices to bound the increment by
$Cq^{-p}(q\theta)^\kappa$. For $p>1$ the global derivative bound gives
$Cq^{1-p}\theta$. For $p\le1$ write $b_q:=|t|^p-P_q$. We have
$||\cos\phi|^p-|\cos\phi'|^p|\le|\cos\phi-\cos\phi'|^p\le\theta^p$. The
derivative of $\phi\mapsto P_q(\cos\phi)$ is
$-\sum_{1\le j\le k_q}2ja_{2j}\sin(2j\phi)$. For $p<1$ its absolute value
is at most $\sum_{j\le k_q}2jc_j\le C\sum_{j\le k_q}j^{-p}\le Cq^{1-p}$. For
$p=1$ we have $j_0=1$, and with $e_j:=2jc_j\le C/j$ the derivative equals
$-\varsigma\sum_{j\le k_q}e_j\sin(2j\psi)$, where we may take
$|\psi|\le\pi/2$ by periodicity, and it vanishes if $\psi=0$. For $\psi\ne0$
define $j_1:=\lfloor1/|\psi|\rfloor$. The terms with $j\le j_1$ contribute at
most $\sum_{j\le j_1}(C/j)\,2j|\psi|\le2C$, and Abel summation bounds the
remaining terms by $e_{j_1+1}/|\sin\psi|\le C|\psi|\cdot\pi/(2|\psi|)$,
using $|\sin\psi|\ge2|\psi|/\pi$. So for $p\le1$ the derivative of
$P_q(\cos\phi)$ is at most $Cq^{1-p}$, and
$|b_q(\cos\phi)-b_q(\cos\phi')|\le\theta^p+Cq^{1-p}\theta\le
(1+C)q^{-p}(q\theta)^p$.
\end{proof}

\begin{lemma}[A polynomial envelope]\label{lem:x-ha-envelope}
For every integer $q\ge64$ there is an even polynomial $W_q$ of degree at
most $q$ such that
\[
 \frac{c_{\rm W}}{(1+q|t|)^2}\le W_q(t)\le\frac{C_{\rm W}}{(1+q|t|)^2}
 \qquad(|t|\le1),\qquad c_{\rm W}:=\frac1{4400},\quad C_{\rm W}:=1024 .
\]
\end{lemma}

\begin{proof}
Let $m:=\lfloor q/8\rfloor$, so $m\ge8$ and $q/16\le m\le q/8$. For $|t|\le1$
define $\psi:=\arcsin t\in[-\pi/2,\pi/2]$ and, for $m\le l<2m$,
\[
 K_l(t):=(\frac{\sin(l\psi)}{l\sin\psi})^2
 =\frac1{l^2}\sum_{|k|<l}(l-|k|)\cos(2k\psi)
 =\frac1{l^2}\sum_{|k|<l}(l-|k|)\,T_{|k|}(1-2t^2),
\]
with $K_l(0):=1$. The second expression is $l^{-2}|\sum_{j<l}e^{2ij\psi}|^2$,
and the third uses $\cos(2k\psi)=T_{|k|}(\cos2\psi)$ and
$\cos2\psi=1-2t^2$. So $K_l$ is an even polynomial of degree $2l-2$, and
$W_q:=\frac1m\sum_{l=m}^{2m-1}K_l$ is an even polynomial of degree at most
$4m-4<q$.

\emph{Upper bound.} We have $|\sin(l\psi)|\le l|\sin\psi|$ and
$|\sin(l\psi)|\le1$, with $\sin\psi=t$. So $K_l(t)\le\min\{1,(mt)^{-2}\}\le
4(1+m|t|)^{-2}$, and $1+q|t|\le16(1+m|t|)$ gives $W_q(t)\le1024(1+q|t|)^{-2}$.

\emph{Lower bound.} If $m|\psi|\le1/2$, then $l|\psi|<1$ for $l<2m$, and
$|\sin\psi|\le|\psi|$ gives $K_l(t)\ge(\sin(l|\psi|)/(l|\psi|))^2\ge
\sin^21\ge0.7$, so $W_q(t)\ge0.7$. Let $m|\psi|>1/2$, so $t\ne0$. Then
$K_l(t)=\sin^2(l\psi)/(l^2t^2)\ge\sin^2(l\psi)/(4m^2t^2)$, and it suffices to
show $\frac1m\sum_{l=m}^{2m-1}\sin^2(l\psi)\ge1/1100$, since then
$W_q(t)\ge(4400m^2t^2)^{-1}\ge(4400(1+q|t|)^2)^{-1}$. If $m|\psi|\ge\pi$,
then
\[
 \frac1m\sum_{l=m}^{2m-1}\sin^2(l\psi)=\frac12-\frac1{2m}\sum_{l=m}^{2m-1}
 \cos(2l\psi)\ge\frac12-\frac1{2m|\sin\psi|}\ge\frac12-\frac\pi{4m|\psi|}
 \ge\frac14 .
\]
If $1/2<m|\psi|<\pi$, the points $l|\psi|$, $m\le l<2m$, lie in $(1/2,2\pi)$
with spacing $|\psi|>1/(2m)$. A point of $(1/2,2\pi)$ at distance at least
$1/16$ from $\pi\mathbb Z$ has $\sin^2\ge\sin^2(1/16)$. The other points lie
in two intervals of length at most $1/8$, which contain fewer than
$2(m/4+1)$ of the $l|\psi|$. As $m\ge8$, at least $m/4$ indices remain,
and the average is at least $\frac14\sin^2(1/16)\ge1/1100$.
\end{proof}

\begin{lemma}[Nets on the sphere]\label{lem:x-ha-net}
Let $d\ge1$ and $0<\theta\le\pi$. There is a set $Y\subseteq S^d$ with at
most $(2\pi/\theta)^d$ points such that every point of $S^d$ is at geodesic
distance at most $\theta$ from $Y$.
\end{lemma}

\begin{proof}
In polar coordinates about its center, a geodesic ball of radius
$\varrho\le\pi/2$ in $S^d$ has normalized measure
\[
 \frac{\int_0^\varrho\sin^{d-1}s\,\d s}{\int_0^\pi\sin^{d-1}s\,\d s}
 \ge\frac{(2/\pi)^{d-1}\varrho^d/d}\pi\ge(\frac\varrho\pi)^d,
\]
by $\sin s\ge2s/\pi$ on $[0,\pi/2]$ and $2^{d-1}\ge d$. If the points of a
set are pairwise at geodesic distance more than $\theta$, the balls of
radius $\theta/2$ about them are disjoint, so there are at most
$(2\pi/\theta)^d$ points. A maximal such set is the required $Y$.
\end{proof}

\begin{proof}[Proof of the upper bound in Lemma~\ref{lem:regular-high-accuracy}]
Let $d:=r-1$. By the equivalence above we may take unit vectors
$u_1,\ldots,u_n\in S^d$ and a probability vector $\nu$ with
\[
 F(x):=\sum_i\nu_i|\langle u_i,x\rangle|^p\ge c_*r^{-p/2},\qquad
 \sum_i\nu_i(1+q|\langle u_i,x\rangle|)^{-2}\le A_*\sqrt r/q
\]
for unit $x$.
We find a probability vector $\mu^*$ with $\supp\mu^*\subseteq\supp\nu$, at
most $M$ nonzero entries and $|\Phi_x(\mu^*)-F(x)|\le\epsilon F(x)$ for unit
$x$, where $\Phi_x(\mu):=\sum_i\mu_i|\langle u_i,x\rangle|^p$.

\emph{Constants.} Let $\kappa:=\min\{p,1\}$, let $C_{\rm b}$, $c_{\rm W}$ and
$C_{\rm W}$ be as in Lemmas~\ref{lem:x-ha-cheb} and~\ref{lem:x-ha-envelope},
and let $c_a$ be as in Lemma~\ref{lem:x-ha-binom}. We take
$c\le\min\{c_a,1/8\}$, then $C$ large as below, and finally $r$ large in
terms of all constants. As in~(a) of the proof of
Theorem~\ref{thm:high_accuracy_main}(a), $R\ge1$ and $\log M\ge\log R\ge ad$.
Since $cdM^{1/d}\ge cd\ge2$, we have $q\ge cdM^{1/d}/2\ge\max\{p+2,64\}$.
Lemma~\ref{lem:x-ha-binom} gives $\binom{q+d}d\le e^{-ad/2}M$, so
$4\binom{q+d}d\le\delta_0M$, and $M\ge\max\{5n_{\rm c},20\}$.

\emph{Preserved envelope.} Let $\mathcal S$ be the set of $\mu\in[0,\infty)^n$
with $\sum_i\mu_if(u_i)=\sum_i\nu_if(u_i)$ for every polynomial $f$ of degree
at most $q$ on $\R^r$. Every $\mu\in\mathcal S$ is a probability vector. For
unit $x$ the function $u\mapsto W_q(\langle u,x\rangle)$ has degree at most
$q$, so Lemma~\ref{lem:x-ha-envelope} gives, for $\mu\in\mathcal S$,
\begin{equation}\label{lem:x-ha-eq-kept}
 \sum_i\mu_i(1+q|\langle u_i,x\rangle|)^{-2}\le\frac1{c_{\rm W}}\sum_i\mu_i
 W_q(\langle u_i,x\rangle)=\frac1{c_{\rm W}}\sum_i\nu_iW_q(\langle u_i,x
 \rangle)\le\frac{C_{\rm W}A_*\sqrt r}{c_{\rm W}q}.
\end{equation}

\emph{One round.} Let $\mu\in\mathcal S$ with $N:=|\supp\mu|>M$, let $J$ be
its light rows, and use the singletons of $J$ as units, so
$|\mathcal P|=|J|\ge N/2$ by Lemma~\ref{lem:t-light-rows}. Let
$H\subseteq\R^J$ consist of the $z$ with $\sum_{i\in J}\mu_iz_if(u_i)=0$ for
every polynomial $f$ of degree at most $q$. As in the proof of
Lemma~\ref{lem:x-ha-round}, its codimension is at most $2\binom{q+d}d\le
\delta_0N/2\le\delta_0|\mathcal P|$, and $\mu^z\in\mathcal S$ for $z\in
H\cap[-1,1]^J$. For $z\in H$ and unit $x$, the polynomial part
$u\mapsto P_q(\langle u,x\rangle)$ of Lemma~\ref{lem:x-ha-cheb} cancels, so
\[
 \Phi_x(\mu^z)-\Phi_x(\mu)=\sum_{i\in J}\mu_iz_ib_q(\langle u_i,x\rangle)
 =\langle z,w_x\rangle,\qquad w_x:=\Pi_H(\mu_ib_q(\langle u_i,x\rangle)
 )_{i\in J}.
\]
By Lemma~\ref{lem:x-ha-cheb}, $\mu_i\le2/N$ on $J$ and Eq.~\eqref{lem:x-ha-eq-kept},
\[
 \|w_x\|_2^2\le\frac2N\sum_{i\in J}\mu_ib_q(\langle u_i,x\rangle)^2\le
 \frac{2C_{\rm b}^2q^{-2p}}N\sum_i\mu_i(1+q|\langle u_i,x\rangle|)^{-2}
 \le\rho_N^2,
\]
where $\rho_N:=C_9q^{-p}(A_*\sqrt r/(qN))^{1/2}$ and
$C_9:=C_{\rm b}(2C_{\rm W}/c_{\rm W})^{1/2}$. Let $x,x'$ be unit vectors
at geodesic distance $\theta\le1/q$. For unit $u$ the angles
$\arccos\langle u,x\rangle$ and $\arccos\langle u,x'\rangle$ differ by at
most $\theta$, by the triangle inequality for angles. So the second bound
of Lemma~\ref{lem:x-ha-cheb}, applied with $\phi=\arccos\langle u_i,x
\rangle$, and the same computation give $\|w_x-w_{x'}\|_2\le(q\theta)^\kappa
\rho_N$.

\emph{Chaining.} Let $\mathcal V:=\{w_x:x\in S^d\}\subseteq H$, so
$\sup_{v\in\mathcal V}\|v\|_2\le\rho_N$. For $j\ge1$ define
$\theta_j:=2^{-j/(2\kappa)}/q\le1/q$ and let $Y_j$ be the net of
Lemma~\ref{lem:x-ha-net} for $\theta_j$. Every $w_x$ is within
$(q\theta_j)^\kappa\rho_N=2^{-j/2}\rho_N$ of some $w_y$, $y\in Y_j$. Since
$q^d\le(cd)^dM<(cd)^dN$,
\[
 N_j:=N(\mathcal V,2^{-j/2}\rho_N)\le(2\pi q)^d2^{jd/(2\kappa)}
 \le(2\pi cd)^d2^{jd/(2\kappa)}N .
\]
Define $\lambda_j:=(2d\log(2d)+2dj/\kappa)^{1/2}\ge1$. Since $\sqrt2/e<e^{-1/2}$
and $d/\kappa\ge1$,
\[
 \sum_{j\ge1}N_je^{-\lambda_j^2/2}\le(\pi c)^dN\sum_{j\ge1}(\sqrt2/e
 )^{jd/\kappa}\le2(\pi c)^dN\le\frac{\delta_0}4N
\]
for large $r$, as $\pi c<1/2$. Lemma~\ref{lem:t-sidak-chaining}(c) gives
countably many strips in $H$, each failing with probability at most
$e^{-\lambda_1^2/2}\le1/2$ and with total failure probability at most
$\delta_0N/4$, such that every point $y$ of their intersection $\mathcal K$
satisfies
\[
 \sup_{x\in S^d}|\langle y,w_x\rangle|\le2\rho_N\sum_{j\ge1}\lambda_j2^{-j/2}
 \le C_\kappa\rho_N\sqrt{d\log(2d)},\qquad C_\kappa:=7+12\kappa^{-1/2},
\]
by $\sum_j2^{-j/2}\le2.42$ and $\sum_j\sqrt j\,2^{-j/2}\le4.15$. By part~(d),
$\gamma_H(\mathcal K)\ge e^{-\delta_0N/2}\ge e^{-\delta_0|\mathcal P|}$.

\emph{Iteration.} Thus Lemma~\ref{lem:t-iterated-coloring} applies with
$I:=[n]$, the set $\mathcal S$, $\mathcal X:=S^d$, $\omega:=1$,
$N_*:=M$, $b_*:=1/2$ and
$a_*:=C_\kappa C_9\sqrt{d\log(2d)}\,q^{-p}(A_*\sqrt r/q)^{1/2}$. Starting
from $\nu$ it gives $\mu^*\in\mathcal S$ with at most $M$ nonzero entries,
$\supp\mu^*\subseteq\supp\nu$ and
$|\Phi_x(\mu^*)-F(x)|\le C_{10}a_*M^{-1/2}$ for unit $x$, where
$C_{10}:=1/(1-(1-\delta_1/10)^{1/2})$. Divide by $F(x)\ge c_*r^{-p/2}$ and use
$q\ge cdM^{1/d}/2\ge crM^{1/d}/4$ and $\frac12+\frac{p+1/2}d=\frac1{2\eta}$.
The relative error is at most
\[
 C_{11}\sqrt{\log(2r)}\,r^{1/4-p/2}M^{-1/(2\eta)},\qquad
 C_{11}:=C_{10}C_\kappa C_9(\frac4c)^{p+1/2}\frac{A_*^{1/2}}{c_*}.
\]
As $M\ge CR$ and $1/(2\eta)\ge1/2$, we have $M^{1/(2\eta)}\ge C^{1/2}
\epsilon^{-1}r^{1/4-p/2}\sqrt{\log(2r)}$. With $C\ge C_{11}^2$ the relative
error is at most $\epsilon$. So every input in the class has a sparsifier
with at most $M\le2CR$ rows, uniformly in $L\ge ar$.
\end{proof}

\begin{proof}[Proof of the lower bound in Lemma~\ref{lem:regular-high-accuracy}]
Let $c_{\rm 1t}$, $c_{\rm fl}$ and $A_{\rm sl}$ be the constants of
Lemma~\ref{lem:s-onetag}, which depend only on $p$ and $a$. Require
$c_*\le c_{\rm fl}$, $A_*\ge A_{\rm sl}$ and $c\le1\le C$, and let $r$ be
so large that $r-1\ge\max\{\log_2(2C),1/c\}$ and Lemma~\ref{lem:s-onetag}
applies with $Q_0:=2$. Here $R=R_{\rm 1t}\ge1$. Then $q\ge\lfloor c(r-1)
\rfloor\ge1$ and, as $M\le CR+1\le2CR$,
\[
 q\le c(r-1)(2CR)^{1/(r-1)}\le(r-1)(2C)^{1/(r-1)}R^{1/(r-1)}
 \le2rR^{1/(r-1)} .
\]
So Lemma~\ref{lem:s-onetag}(b) applies to the corollary's $q$. Take the
input of part~(a) with the property of part~(b). It has rows
$\alpha_t^{1/p}(1,u_t)^\top$ with masses $\alpha_t\in[1/2,3/2]$. After the
change of coordinates $(y_0,y)\mapsto(y_0,\sqrt d\,y)$ of the rows these
are $\alpha_t^{1/p}\ell_t\tilde u_t^\top$, and the law $\nu$ giving
$\tilde u_t$ mass proportional to $\alpha_t\ell_t^p$ satisfies both
conditions of the corollary by part~(b), since $c_{\rm fl}\ge c_*$ and
$A_{\rm sl}\le A_*$. So this input lies in the class. By part~(a) each of
its sparsifiers has more than $c_{\rm 1t}R$ rows.
\end{proof}

\subsection{Notation and Caps Beyond Exponential Accuracy}

Write $d:=r-1$, $\mathbb S:=S^{r-1}$, and let $\omega$ be the normalized
surface measure on $\mathbb S$. For functions on $\mathbb S$ define
$\langle f,g\rangle:=\int fg\,\d\omega$ and let $\|\cdot\|_s$ be the norm
of $L_s(\omega)$. For $v\in\mathbb S$ and $\varrho>0$ let
$\operatorname{Cap}(v,\varrho):=\{u\in\mathbb S:\|u-v\|_2\le\varrho\}$ and
$\omega_r(\varrho):=\omega(\operatorname{Cap}(v,\varrho))$. For $m\ge3$
let $\phi_m(t):=c_m(1-t^2)^{(m-3)/2}$ on $(-1,1)$ be the density of one
coordinate of a uniform point of $S^{m-1}$. Let $\mathcal P_m$ be the
space of real polynomials of degree at most $m$ on $\R^r$. Restricted to
$\mathbb S$, a polynomial of degree at most $m$ agrees with the sum of a
homogeneous polynomial of degree $m$ and one of degree $m-1$, obtained by
multiplying each homogeneous part by a power of $\|u\|_2^2$. Hence the
restrictions of $\mathcal P_m$ to any subset of $\mathbb S$ form a space of
dimension at most
\begin{equation}\label{lem:x-sc-eq-dim}
 d_m:=2\binom{m+d}{d}\le2(\frac{e(m+d)}d)^d .
\end{equation}
Define $\theta_p:=\min\{p,1\}$ and $\bar p:=\max\{1,p/2\}$.

\begin{lemma}[Caps and coordinates]\label{lem:x-sc-caps}
Let $r\ge3$.
\begin{enumerate}[label=(\alph*),nosep,leftmargin=*]
\item For $0<\varrho\le1$, $(2\varrho/\pi)^d/(2d)\le\omega_r(\varrho)\le
(\pi\varrho/2)^d$.
\item For $0<\varrho\le1$, a subset of $\mathbb S$ with pairwise distances
at least $\varrho$ has at most $2d(\pi/\varrho)^d$ points.
\item For $0<\tau\le\varrho\le1/2$, every cap of radius $\varrho$ is
covered by at most $2d(\pi^2\varrho/\tau)^d$ caps of radius $\tau$.
\item For $m\ge6$, $\phi_m\le\sqrt m$ on $(-1,1)$ and
$\phi_m(t)\ge\sqrt m/8$ for $|t|\le1/\sqrt m$.
\end{enumerate}
\end{lemma}

\begin{proof}
(a) The cap of Euclidean radius $\varrho$ has angular radius
$\vartheta=2\arcsin(\varrho/2)\in[\varrho,\pi\varrho/2]$, and
$\vartheta\le\pi/2$. In polar coordinates about its center,
$\omega_r(\varrho)=\int_0^\vartheta\sin^{d-1}\varphi\,\d\varphi/
\int_0^\pi\sin^{d-1}\varphi\,\d\varphi$. Since $\sin\varphi\ge2\varphi/\pi$
on $[0,\pi/2]$, the numerator is at least $(2/\pi)^{d-1}\vartheta^d/d$,
and the denominator is at most $\pi$. This gives the lower bound. The
numerator is at most $\vartheta^d/d$. On
$|\varphi-\pi/2|\le d^{-1/2}$ we have $\sin^{d-1}\varphi\ge
(1-1/(2d))^{d-1}\ge e^{-3/4}$, so the denominator is at least
$2e^{-3/4}d^{-1/2}$. Hence $\omega_r(\varrho)\le\vartheta^d e^{3/4}/(2\sqrt d)
\le\vartheta^d$ for $d\ge2$, and the case $d=1$ is direct.

(b) The caps of radius $\varrho/2$ about the points are disjoint, and (a)
bounds their number.

(c) Take a maximal subset of the cap with pairwise distances at least
$\tau$. The caps of radius $\tau$ about its points cover the cap, and the
caps of radius $\tau/2$ about them are disjoint and lie in a cap of radius
$\varrho+\tau/2\le1$. By (a) their number is at most
$(\pi(\varrho+\tau/2)/2)^d\cdot2d(\pi/\tau)^d\le2d(\pi^2\varrho/\tau)^d$.

(d) We have $c_m^{-1}=\int_{-1}^1(1-t^2)^{(m-3)/2}\d t\le\int
e^{-(m-3)t^2/2}\d t=(2\pi/(m-3))^{1/2}$, so $c_m\ge\sqrt m/4$ for
$m\ge5$. For $|t|\le m^{-1/2}$ the integrand is at least
$(1-1/m)^{m/2}\ge1/2$, so $c_m^{-1}\ge m^{-1/2}$. Both claims follow.
\end{proof}

\subsection{Section Operators}

For $v\in\mathbb S$ and $-1<t<1$ let $\sigma_{v,t}$ be the normalized
surface measure on the $(r-2)$-sphere
$\{x\in\mathbb S:\langle v,x\rangle=t\}$, and let $\sigma_{v,\pm1}$ be the
point mass at $\pm v$. For a bounded Borel function $g$ on $\mathbb S$ define
$(\mathcal T_tg)(v):=\int g\,\d\sigma_{v,t}$. Writing $\sigma_{v,t}$ as the
image of $\sigma_{e_1,t}$ under an orthogonal map that sends $e_1$ to $v$
and depends on $v$ in a Borel way, for instance a reflection, shows that
$\mathcal T_tg$ is bounded and Borel. Locally such maps can be chosen
continuous in $v$, so $\mathcal T_tg$ is continuous when $g$ is. For an orthogonal map $\Theta$, the map $\Theta$ carries
$\sigma_{v,t}$ to $\sigma_{\Theta v,t}$.

\begin{lemma}[Invariant pairs]\label{lem:x-sc-pairs}
Let $r\ge3$.
\begin{enumerate}[label=(\alph*),nosep,leftmargin=*]
\item Let $(u,y)$ be a random element of $\mathbb S\times\mathbb S$ whose
law is invariant under $(u,y)\mapsto(\Theta u,\Theta y)$ for every
orthogonal $\Theta$, and let $\psi$ be the law of $\langle u,y\rangle$.
For every bounded Borel function $\Phi$ on $\mathbb S\times\mathbb S$,
\[
 \E[\Phi(u,y)]=\int_{[-1,1]}\int_{\mathbb S}\int\Phi(v,x)\,
 \d\sigma_{v,t}(x)\,\d\omega(v)\,\d\psi(t).
\]
\item For bounded Borel $g$ and $-1\le t\le1$, $\int\mathcal T_tg\,\d\omega
=\int g\,\d\omega$ and $\|\mathcal T_tg\|_s\le\|g\|_s$ for $1\le s<\infty$.
\item For bounded Borel $f,g$ on $\mathbb S$ and $\psi_0$ on $[-1,1]$,
\[
 \iint\psi_0(\langle u,x\rangle)f(u)g(x)\,\d\omega(u)\,\d\omega(x)
 =\int_{-1}^1\psi_0(t)\phi_r(t)\langle f,\mathcal T_tg\rangle\,\d t .
\]
\end{enumerate}
\end{lemma}

\begin{proof}
(a) For $-1\le t\le1$ let $\varpi_t$ be the probability measure
$\int\Phi\,\d\varpi_t:=\int\int\Phi(v,x)\,\d\sigma_{v,t}(x)\,\d\omega(v)$
on $O_t:=\{(v,x)\in\mathbb S\times\mathbb S:\langle v,x\rangle=t\}$. It is
invariant under the diagonal action of the orthogonal group, since
$\Theta$ carries $\sigma_{v,t}$ to $\sigma_{\Theta v,t}$ and $\omega$ is
invariant. The group acts transitively on $O_t$. Let $\mathrm H$ be its
normalized Haar measure and fix $z_t\in O_t$. If $\nu$ is any invariant
probability on $O_t$ and $z=\Theta_zz_t$, then invariance, Fubini's
theorem and the invariance of $\mathrm H$ under right translation by
$\Theta_z$ give
$\int\Phi\,\d\nu=\int\int\Phi(\Theta z)\,\d\mathrm H(\Theta)\,\d\nu(z)
=\int\Phi(\Theta z_t)\,\d\mathrm H(\Theta)$. Hence $\nu=\varpi_t$, and
$t\mapsto\int\Phi\,\d\varpi_t$ is Borel, by this formula with
$z_t=(e_1,te_1+(1-t^2)^{1/2}e_2)$. Let $\Theta$ have law $\mathrm H$,
independent of $(u,y)$. Then $(\Theta u,\Theta y)$ has the law of $(u,y)$,
and conditionally on $(u,y)$ its law is an invariant probability on
$O_{\langle u,y\rangle}$, hence equal to $\varpi_{\langle u,y\rangle}$.
Thus $\E[\Phi(u,y)]=\E[\int\Phi\,\d\varpi_{\langle u,y\rangle}]$, which is the
claim.

(b) The second marginal of $\varpi_t$ is an invariant probability on
$\mathbb S$, hence $\omega$ by the same uniqueness argument. This gives
$\int\mathcal T_tg\,\d\omega=\int g\,\d\omega$. Jensen's inequality gives
$|\mathcal T_tg|^s\le\mathcal T_t|g|^s$, and we integrate.

(c) Apply (a) to independent $u,x$ with law $\omega$. Their joint law is
invariant, and $\langle u,x\rangle$ has density $\phi_r$. On $O_t$ the
function $\psi_0(\langle v,x\rangle)f(v)g(x)$ equals $\psi_0(t)f(v)g(x)$.
\end{proof}

\begin{lemma}[Section inequality]\label{lem:x-sc-section}
Let $r\ge8$ and let $f,g\ge0$ be bounded Borel functions on $\mathbb S$.
\begin{enumerate}[label=(\alph*),nosep,leftmargin=*]
\item $\|\mathcal T_0g\|_r\le\|g\|_{r/d}$.
\item For $0<h\le1/(4\sqrt r)$,
\begin{equation}\label{lem:x-sc-eq-section}
 \mathcal M_{f,g}(h):=\iint\mathbf1_{\{|\langle u,x\rangle|\le h\}}
 f(u)g(x)\,\d\omega(u)\,\d\omega(x)\le64h\sqrt r\,\|f\|_{r/d}\|g\|_{r/d} .
\end{equation}
\end{enumerate}
\end{lemma}

\begin{proof}
(a) First let $g$ be continuous and positive, and let $\mathsf K$ be the
star body $\{\varrho v:v\in\mathbb S,\ 0\le\varrho\le g(v)^{1/d}\}$.
Integration in polar coordinates gives
$\operatorname{vol}_r(\mathsf K)=\operatorname{vol}_r(B_2^r)\int g^{r/d}\,\d\omega$
and, for $x\in\mathbb S$, $\operatorname{vol}_{r-1}(\mathsf K\cap x^\perp)
=\operatorname{vol}_{r-1}(B_2^{r-1})\,\mathcal T_0g(x)$. Lemma~\ref{lem:pre_busemann}
states that
\[
 \int_{\mathbb S}\operatorname{vol}_{r-1}(\mathsf K\cap x^\perp)^r\,\d\omega(x)
 \le\frac{\operatorname{vol}_{r-1}(B_2^{r-1})^r}
 {\operatorname{vol}_r(B_2^r)^{r-1}}\operatorname{vol}_r(\mathsf K)^{r-1},
\]
with equality for the Euclidean ball. Substituting the two formulas gives
$\int(\mathcal T_0g)^r\,\d\omega\le(\int g^{r/d}\,\d\omega)^d$, which is
the claim. For general $g$ choose continuous positive $g_j$ with
$g_j\to g$ in $L_{r/d}(\omega)$ and in $L_1(\omega)$. By
Lemma~\ref{lem:x-sc-pairs}(b), $\mathcal T_0g_j\to\mathcal T_0g$ in
$L_1(\omega)$, hence almost everywhere along a subsequence, and Fatou's
lemma gives the claim.

(b) Let $\sin\alpha=4h\sqrt r\le1$ with $0<\alpha\le\pi/2$. Let $(u,x)$
have law $\varpi_0$, so that $x$ is a uniform unit vector orthogonal to
$u$, and given $(u,x)$ let $w$ have law $\sigma_{x,0}$. Define
$x':=\cos\alpha\,x+\sin\alpha\,w$. The map
$w\mapsto\cos\alpha\,x+\sin\alpha\,w$ carries $\sigma_{x,0}$ to
$\sigma_{x,\cos\alpha}$, so
$\E[f(u)g(x')]=\E[f(u)(\mathcal T_{\cos\alpha}g)(x)]=
\langle f,\mathcal T_0\mathcal T_{\cos\alpha}g\rangle$. The law of
$(u,x,w)$ is invariant under orthogonal maps, hence so is the law of
$(u,x')$. Moreover $\langle u,x'\rangle=\sin\alpha\langle u,w\rangle$, and
since $u$ is a unit vector of $x^\perp$, the variable $\langle u,w\rangle$
has density $\phi_{r-1}$. So $\langle u,x'\rangle$ has density
$\psi_\alpha(t):=\phi_{r-1}(t/\sin\alpha)/\sin\alpha$. For $|t|\le h$ we
have $|t|/\sin\alpha\le(4\sqrt r)^{-1}\le(r-1)^{-1/2}$, so
Lemma~\ref{lem:x-sc-caps}(d) gives
$\psi_\alpha(t)\ge\sqrt{r-1}/(32h\sqrt r)\ge1/(64h)$, while
$\phi_r(t)\le\sqrt r$. Since $\langle f,\mathcal T_tg\rangle\ge0$,
Lemma~\ref{lem:x-sc-pairs}(c), then Lemma~\ref{lem:x-sc-pairs}(a) for
$(u,x')$, give
\[
 \mathcal M_{f,g}(h)=\int_{-h}^h\phi_r\langle f,\mathcal T_tg\rangle\,\d t
 \le64h\sqrt r\int\psi_\alpha(t)\langle f,\mathcal T_tg\rangle\,\d t
 =64h\sqrt r\,\langle f,\mathcal T_0\mathcal T_{\cos\alpha}g\rangle .
\]
The exponents $r/d$ and $r$ are conjugate. By H\"older's inequality, part
(a) and Lemma~\ref{lem:x-sc-pairs}(b), the last pairing is at most
$\|f\|_{r/d}\|\mathcal T_{\cos\alpha}g\|_{r/d}\le\|f\|_{r/d}\|g\|_{r/d}$.
\end{proof}

\begin{lemma}[Pair densities]\label{lem:x-sc-pair-density}
Let $r\ge8$ and let $f,g$ be bounded Borel functions on $\mathbb S$. The
function $\varphi_{f,g}(t):=\langle f,\mathcal T_tg\rangle$ is continuous
on $(-1,1)$, so by Lemma~\ref{lem:x-sc-pairs}(c) the function
$\phi_r\varphi_{f,g}$ is the density of the image of
$f(u)g(x)\,\d\omega(u)\,\d\omega(x)$ under $(u,x)\mapsto\langle u,x\rangle$.
It is nonnegative if $f,g\ge0$. For every $k\le(r-4)/4$ it is $k$ times continuously differentiable on
$(-1/2,1/2)$, and
\[
 |\varphi_{f,g}^{(k)}(t)|\le2\cdot8^kk!\,\|f\|_2\|g\|_2\qquad
 (|t|<1/2,\ 1\le k\le(r-4)/4).
\]
\end{lemma}

\begin{proof}
Let $\mathcal H_n$ be the space of restrictions to $\mathbb S$ of the
homogeneous harmonic polynomials of degree $n$. By Lemma~\ref{lem:pre_harmonics}(a), these spaces are
orthogonal in $L_2(\omega)$, and the
restriction of every polynomial is a finite sum of their
elements. By the Stone--Weierstrass theorem (Lemma~\ref{lem:pre_stone_weierstrass})
restrictions of polynomials are dense in $L_2(\omega)$, so
$L_2(\omega)=\bigoplus_n\mathcal H_n$. The Funk--Hecke
formula (Lemma~\ref{lem:pre_harmonics}(b)) states that for bounded Borel $\psi_0$
on $[-1,1]$ and $Y\in\mathcal H_n$,
$\int\psi_0(\langle u,v\rangle)Y(v)\,\d\omega(v)=\lambda_n(\psi_0)Y(u)$, with
a scalar $\lambda_n(\psi_0)$ not depending on $Y$. By linearity this also
holds for complex combinations of elements of $\mathcal H_n$.

\emph{The section multipliers.} Let $Y\ne0$ be such a combination and fix
$t_0\in(-1,1)$. For every bounded Borel $f$,
Lemma~\ref{lem:x-sc-pairs}(c) and the Funk--Hecke formula give
$\int\psi_0(t)\phi_r(t)\langle f,\mathcal T_tY\rangle\,\d t=
\lambda_n(\psi_0)\langle f,Y\rangle$. The function
$t\mapsto\langle f,\mathcal T_tY\rangle$ is continuous on $(-1,1)$ by
dominated convergence, since $Y$ is continuous. Take for $\psi_0$
nonnegative continuous functions $\psi_j$ supported in
$[t_0-1/j,t_0+1/j]$ with $\int\psi_j\phi_r=1$. Then
$\langle f,\mathcal T_{t_0}Y\rangle=\lim_j\lambda_n(\psi_j)\langle f,Y
\rangle$. Choosing $f$ with $\langle f,Y\rangle\ne0$ shows that
$a_n(t_0):=\lim_j\lambda_n(\psi_j)$ exists, and then
$\mathcal T_{t_0}Y=a_n(t_0)Y$ almost everywhere, hence everywhere since
both sides are continuous in $u$. The scalar $a_n(t_0)$ does not depend on $Y$.
For $Y_0(x):=(x_1+\mathrm ix_2)^n$, which is harmonic, we have $Y_0(e_1)=1$.
A point of law $\sigma_{e_1,t}$ is $te_1+(1-t^2)^{1/2}w$ with $w$ uniform
on the unit sphere of $e_1^\perp$, and $V:=\langle w,e_2\rangle$ has
density $\phi_{r-1}$. Hence
\[
 a_n(t)=(\mathcal T_tY_0)(e_1)=\E[(t+\mathrm i(1-t^2)^{1/2}V)^n] .
\]
Since $V$ is symmetric, only even powers of $V$ contribute, so $a_n$ is a
polynomial. Also $|t+\mathrm i(1-t^2)^{1/2}V|^2=1-(1-t^2)(1-V^2)\le1$, so
$|a_n|\le1$.

\emph{The expansion.} Let $f_n,g_n$ be the components of $f,g$ in
$\mathcal H_n$, and $g_{\le M}:=\sum_{n\le M}g_n$. Since $g-g_{\le M}$ is
bounded, Lemma~\ref{lem:x-sc-pairs}(b) gives
$\|\mathcal T_t(g-g_{\le M})\|_2\le\|g-g_{\le M}\|_2\to0$. Hence
$\mathcal T_tg=\sum_na_n(t)g_n$ in $L_2(\omega)$ and
$\varphi_{f,g}(t)=\sum_na_n(t)\langle f_n,g_n\rangle$, where
$\sum_n|\langle f_n,g_n\rangle|\le\|f\|_2\|g\|_2$ by the Cauchy--Schwarz
inequality. The series converges uniformly, so $\varphi_{f,g}$ is
continuous. If $f,g\ge0$, then $\varphi_{f,g}\ge0$ because $\mathcal T_t$
is positive.

\emph{Derivatives.} Fix $V\in(-1,1)$ and let
$b(z):=z+\mathrm i(1-z^2)^{1/2}V$ with the principal branch, analytic for
$|z|<1$. Let $|t|\le1/2$ be real and $\varrho_V:=(1-V^2)/8$. For
$|z-t|\le\varrho_V$ we have $|z|\le5/8$, so the derivative of
$(1-z^2)^{1/2}$ has modulus at most $(5/8)(39/64)^{-1/2}<1$ and
$|b(z)-b(t)|\le2\varrho_V$. Also
$|b(t)|=(1-(1-t^2)(1-V^2))^{1/2}\le1-\frac38(1-V^2)$. Hence $|b(z)|\le1$ on
this disk, and Cauchy's estimate gives
$|\partial_t^kb(t)^n|\le k!\,8^k(1-V^2)^{-k}$. The variable $1-V^2$ has
the beta law with parameters $(r-2)/2$ and $1/2$. With $a:=(r-2)/2$,
\[
 \E[(1-V^2)^{-k}]=\frac{B(a-k,1/2)}{B(a,1/2)}
 =\prod_{j=1}^k(1+\frac1{2(a-j)})\le(1+\frac2r)^{r/4}<2
 \qquad(k\le\frac{r-4}4),
\]
since $2(a-j)\ge r/2$ for $j\le k$. Differentiation under the expectation
is justified by this bound, so $|a_n^{(k)}(t)|\le2\cdot8^kk!$ for all $n$
and $|t|\le1/2$, and the same holds for all orders up to $k$. The
differentiated series therefore converge uniformly, and summing gives the
bound.
\end{proof}

\subsection{Three Polynomial Facts}

\begin{lemma}[Growth, derivatives and central values]\label{lem:x-sc-poly}
Let $P$ be a real polynomial of degree at most $m\ge0$, and
$M:=\sup_{[-1,1]}|P|$.
\begin{enumerate}[label=(\alph*),nosep,leftmargin=*]
\item If $\zeta=a+\mathrm ib$ with $|a|\le5/6$ and $|b|\le\beta\le1/3$, then
$|P(\zeta)|\le2(m+1)Me^{3\beta m}$.
\item $|P'(y)|\le m(m+1)M(1-y^2)^{-1/2}$ for $|y|<1$.
\item If $P\ge0$ on $[-1,1]$, then $\sup_{|y|\le1/2}P(y)\le(m+3)
\int_{-1}^1P$.
\end{enumerate}
\end{lemma}

\begin{proof}
(a) Write $P=\sum_{j\le m}c_jT_j$ with the Chebyshev polynomials $T_j$.
The coefficients are
$c_j=\frac2\pi\int_0^\pi P(\cos\vartheta)\cos(j\vartheta)\,\d\vartheta$
for $j\ge1$, with $\frac1\pi$ in place of $\frac2\pi$ for $j=0$, so
$|c_j|\le2M$. For $w\ne0$ the
identity $T_j((w+w^{-1})/2)=(w^j+w^{-j})/2$ holds on the unit circle, hence
for all $w$. Define $\varrho:=1+3\beta$. The image of the circle
$|w|=\varrho'>1$ under $w\mapsto(w+w^{-1})/2$ is the ellipse with semi-axes
$(\varrho'\pm1/\varrho')/2$. These ellipses are confocal and increase with
$\varrho'$. Hence a point $\zeta$ in the closed region bounded by the ellipse
for $\varrho$ equals $(w+w^{-1})/2$ with $1\le|w|\le\varrho$, and there
$|T_j(\zeta)|\le\varrho^j$. Our $\zeta$ lies in this region. Indeed
the semi-axes satisfy
\[
 \frac{\varrho+1/\varrho}2\ge1,\qquad
 \frac{\varrho-1/\varrho}2=\frac{3\beta(2+3\beta)}{2+6\beta}\ge
 \frac{9\beta}4,
\]
and $(5/6)^2+16/81<1$. Thus $|P(\zeta)|\le2M\sum_{j\le m}\varrho^j\le
2(m+1)M\varrho^m$.

(b) The function $\vartheta\mapsto P(\cos\vartheta)$ is a trigonometric
polynomial $\sum_{|j|\le m}\gamma_je^{\mathrm ij\vartheta}$ with
$|\gamma_j|\le M$. Its derivative is bounded by $\sum_j|j||\gamma_j|\le
m(m+1)M$ and equals $-\sin\vartheta\,P'(\cos\vartheta)$.

(c) The function $T(\vartheta):=P(\cos\vartheta)\sin^2\vartheta$ is a
nonnegative trigonometric polynomial of degree at most $m+2$. Each of its
Fourier coefficients is at most the constant one,
$\gamma_0=\frac1\pi\int_{-1}^1P(y)(1-y^2)^{1/2}\,\d y\le\frac1\pi
\int_{-1}^1P$. Hence $T\le(2m+5)\gamma_0$. For $|y|\le1/2$ we have
$\sin^2\vartheta\ge3/4$, so $P(y)\le\frac4{3\pi}(2m+5)\int_{-1}^1P\le
(m+3)\int_{-1}^1P$.
\end{proof}

\begin{lemma}[Localized approximation of powers]\label{lem:x-sc-localized}
Let $p>0$. For every integer $m\ge1$ there is an even polynomial $Q_m$ of
degree at most $m$ with
\[
 ||t|^p-Q_m(t)|\le\frac{C_pm^{-p}}{1+m|t|}\qquad(|t|\le1).
\]
\end{lemma}

\begin{proof}
Define $\varkappa:=\lfloor p\rfloor+2\ge p+1$ and $\ell:=\varkappa$, so that
$2\ell\ge\varkappa+2$ and $2\ell\ge p+3$. Let $\mathbb T:=[-\pi,\pi]$. For
$n\ge1$ let $J_n(s):=\lambda_n^{-1}(\sin(ns/2)/\sin(s/2))^{2\ell}$ with
$\int_{\mathbb T}J_n=1$. As in Step~2 of the proof of
Lemma~\ref{lem:t-tv-approx}, $J_n$ is a nonnegative even trigonometric
polynomial of degree at most $\ell(n-1)$, and Eq.~\eqref{lem:t-eq-jackson},
whose derivation holds for every integer $\ell\ge1$, gives
$J_n(s)\le C_\ell n(1+n|s|)^{-2\ell}$ on $\mathbb T$. Let
$h(\vartheta):=|\cos\vartheta|^p$ and
\[
 \mathcal J_nh(\vartheta):=\int_{\mathbb T}J_n(s)\sum_{i=1}^\varkappa
 (-1)^{i+1}\binom\varkappa i h(\vartheta+is)\,\d s .
\]
The Fourier-coefficient argument of that Step~2 shows that
$\mathcal J_nh$ is a trigonometric polynomial of degree at most
$\ell(n-1)$. It is even and $\pi$-periodic because $h$ is, so it contains
only the frequencies $\cos(2j\vartheta)=T_{2j}(\cos\vartheta)$ and
$\mathcal J_nh(\vartheta)=Q(\cos\vartheta)$ for an even polynomial $Q$ of
degree at most $\ell(n-1)$. As there,
$\mathcal J_nh(\vartheta)-h(\vartheta)=(-1)^{\varkappa+1}\int_{\mathbb T}
J_n(s)\Delta_s^\varkappa h(\vartheta)\,\d s$, where $\Delta_s^\varkappa$ is
the $\varkappa$th forward difference with step $s$.

Let $\delta=\delta(\vartheta)$ be the distance from $\vartheta$ to
$\pi/2+\pi\mathbb Z$, so $|\cos\vartheta|=\sin\delta\in[2\delta/\pi,
\delta]$. Away from $\pi/2+\pi\mathbb Z$, $h=\psi\circ\cos$ with
$|\psi^{(j)}(c)|\le C_p|c|^{p-j}$, and the chain rule gives
$|h^{(\varkappa)}(\vartheta)|\le C_p|\cos\vartheta|^{p-\varkappa}$. If
$\varkappa|s|\le\delta/2$, every point between $\vartheta$ and
$\vartheta+\varkappa s$ has distance at least $\delta/2$ from
$\pi/2+\pi\mathbb Z$, so $|\Delta_s^\varkappa h(\vartheta)|\le
|s|^\varkappa\sup|h^{(\varkappa)}|\le C_p|s|^\varkappa\delta^{p-\varkappa}$.
Otherwise each point $\vartheta+is$, $0\le i\le\varkappa$, has distance at
most $\delta+\varkappa|s|<3\varkappa|s|$ from $\pi/2+\pi\mathbb Z$, so
$|\Delta_s^\varkappa h(\vartheta)|\le2^\varkappa(3\varkappa|s|)^p$. With
$y=n|s|$ and $\Delta:=n\delta$ we obtain
\[
 |\mathcal J_nh(\vartheta)-h(\vartheta)|\le C_p\int_{y\le\Delta/(2\varkappa)}
 \frac{n^{-\varkappa}y^\varkappa\delta^{p-\varkappa}}{(1+y)^{2\ell}}\,\d y
 +C_p\int_{y>\Delta/(2\varkappa)}\frac{n^{-p}y^p}{(1+y)^{2\ell}}\,\d y .
\]
If $\Delta\ge1$, the first term is at most
$C_pn^{-\varkappa}\delta^{p-\varkappa}=C_pn^{-p}\Delta^{p-\varkappa}\le
C_pn^{-p}\Delta^{-1}$. If $\Delta<1$, we use
$y^\varkappa\delta^{p-\varkappa}\le y^pn^{\varkappa-p}$ on the domain of
the first integral, so it is at most $C_pn^{-p}$. The second term is at
most $C_pn^{-p}(1+\Delta/(2\varkappa))^{p+1-2\ell}\le
C_pn^{-p}(1+\Delta)^{-1}$. Altogether
$|\mathcal J_nh-h|(\vartheta)\le C_pn^{-p}(1+n\delta(\vartheta))^{-1}$. For
$t=\cos\vartheta$ we have $|t|\le\delta(\vartheta)$. Finally take
$n:=\lfloor m/\ell\rfloor+1$, so $\ell(n-1)\le m$ and $n\ge m/\ell$, and
$Q_m:=Q$.
\end{proof}

\subsection{Parameters Beyond Exponential Accuracy}

Fix $p\in(0,\infty)\setminus\{2\}$, and let $\delta_0,\delta_1,n_{\rm c}$
be the constants of Lemma~\ref{lem:t-subspace-coloring}. Define
\[
 D_p:=\frac3{\theta_p}+4,\quad C_3:=D_p+6,\quad D_\flat:=20C_3+80,\quad
 c_q:=\frac1{16e\pi^2},\quad c_{\rm dim}:=\frac{\delta_0}{16}.
\]
Let $C_{\rm sp}\ge1$ be a constant, fixed at the end of the proof. For
$r$ and $\epsilon$ define
\[
 R_{\rm 1t}:=(\epsilon^{-2}r^{1/2-p}\Lambda)^\eta,\quad
 N:=\lceil C_{\rm sp}R_{\rm 1t}\rceil,\quad Y:=N^{1/d},\quad
 q:=\lfloor c_qrY\rfloor,\quad\gamma:=\frac{\sqrt r}q,\quad
 Q:=\lceil qr^{3/\theta_p}\rceil,
\]
\[
 \tau:=\frac{r^{-D_p}}q,\quad\mathsf P:=2rL,\quad
 m_\circ:=\frac{\mathsf P^{3r}}N,\quad m_{\rm T}:=\lceil8L\rceil,\quad
 s:=Y^{-1/8},\quad\varsigma:=16s .
\]
In the Lewis setup Eq.~\eqref{lem:t-eq-lewis} let
$F_\natural(x):=F(x)$ if $p<2$ and $F_\natural(x):=r^{-p/2}$ if $p>2$, and
let $E(x):=r^{p/2}F_\natural(x)$ for unit $x$. Thus $E=r^{p/2}F$ below two,
as in Lemma~\ref{lem:t-inverse-energy}, and $E\equiv1$ above two. On
$\mathbb S$ we have $r^{-\bar p}\le F_\natural\le F$, and $E^{-2}\le r^2$.

\begin{lemma}[Consequences of the range of $L$]\label{lem:x-sc-params}
For every $C_{\rm sp}\ge1$ and $B\ge1$ there are $C_{\rm end}\ge1$ and
$r_0$, depending only on $p$, $C_{\rm sp}$ and $B$, such that the following
hold for $r\ge r_0$ and $L\ge C_{\rm end}r\Lambda$.
\begin{enumerate}[label=(\alph*),nosep,leftmargin=*]
\item $L/2\le\log N\le3L$ and $L/(2r)\le\log Y\le4L/r$.
\item $c_qrY/2\le q\le c_qrY$, $q\ge2r$, and $Q+1\le e^{L/2}$.
\item $s\le\mathsf P^{-B}$ and $\tau\le s$.
\item $qs^2\ge8L$ and $qs^2\ge r^2$.
\item $(\pi/s)^d(e+9eL/d)^d\le N^{1/2}$.
\item $r^{p+3}q^{2p+1}\le e^{4L}$.
\end{enumerate}
\end{lemma}

\begin{proof}
Define $x:=L/r\ge C_{\rm end}\Lambda$. We have
$\log R_{\rm 1t}=\eta(2L-2\log2+(\frac12-p)\log r+\log\Lambda)$ and
$\eta\ge1/2$ for $r\ge2p+2$. For large $C_{\rm end}$ the bracket lies
between $L$ and $2L+\log r$, and $\log(2C_{\rm sp})+2L+\log r\le3L$. This
gives (a), using $d\ge3r/4$. For (b), $c_qrY\ge2$ and $Y\ge e^{x/2}$ give the
bounds on $q$, and $\log(Q+1)\le\log(3c_qr^{1+3/\theta_p}Y)\le
C_p\log r+4x\le L/2$. For (c), $s\le e^{-x/16}$, while
$\log\mathsf P=\log(2r^2x)\le2\Lambda+\log x\le2x/C_{\rm end}+\log x$. So
$s\le\mathsf P^{-B}$ once $C_{\rm end}\ge64B$ and $\log x\le x/(32B)$,
which holds for large $C_{\rm end}$ since $x\ge C_{\rm end}\log2$. Also
$\tau\le1/q\le Y^{-1/8}$. For (d), $qs^2\ge(c_qr/2)Y^{3/4}\ge
(c_qr/2)e^{3x/8}$, which exceeds $8rx$ and $r^2$ for large $C_{\rm end}$,
since $e^{3x/8}\ge(2r)^{3C_{\rm end}/8}$. For (e), $(\pi/s)^d=\pi^dN^{1/8}$
and $m_{\rm T}\le9L$, while $N^{3/8}\ge e^{3rx/16}$ by (a). So (e) follows
from $d\log(\pi e(1+18x))\le3rx/16$, true for large $x$. For (f),
$\log(r^{p+3}q^{2p+1})\le(3p+4)\log r+4(2p+1)x\le4rx$ for large $r$.
\end{proof}

We say that the \emph{standing assumptions} hold if $r\ge r_0$ and
$L\ge C_{\rm end}r\Lambda$, where $C_{\rm end}$ and $r_0$ are given by
Lemma~\ref{lem:x-sc-params} for the current $C_{\rm sp}$ and
$B:=2D_\flat+10$, and $r_0$ is enlarged so that finitely many explicit
conditions on $r$ alone, stated where they are used, also hold. Under the
standing assumptions $\mathsf P\ge2r^2$, so part~(c) implies $s\le1/2$,
$32s\le1/(4\sqrt r)$, $s\le r^{-2/\theta_p}/4$ and
$\varsigma\le\mathsf P^{-2D_\flat-6}$. Also $N\ge e^{L/2}$ is as large as
needed, in particular $N\ge\max\{5n_{\rm c},20\}$.

\subsection{The Approximation Criterion}

For a finite set $\mathcal I$, a vector $\mu\in[0,\infty)^{\mathcal I}$
and $f\in\R^{\mathcal I}$, let $\|f\|_{\mu}^2:=\sum_{i}\mu_if(i)^2$ and,
for a subspace $\mathcal W\subseteq\R^{\mathcal I}$,
$\operatorname{dist}_\mu(f,\mathcal W):=\min_{g\in\mathcal W}\|f-g\|_\mu$.

\begin{lemma}[Approximation spaces for query laws]\label{lem:x-sc-criterion}
There is $C_{\rm mean}$, depending only on $p$, such that the following
holds under the standing assumptions, for every $C_{\rm sp}\ge1$. Let
$\mathcal I$ be a finite set, let $u_i\in\mathbb S$ for $i\in\mathcal I$
(repetitions allowed), let $\mu\in[0,\infty)^{\mathcal I}$ with
$\sum_i\mu_i\le1$, and let $\lambda$ be a finitely supported probability
measure on $\mathbb S$. For $x\in\mathbb S$ let
$\kappa_x(i):=|\langle u_i,x\rangle|^p/F_\natural(x)$. Then there is a
subspace $\mathcal W\subseteq\R^{\mathcal I}$ with
$\dim\mathcal W\le c_{\rm dim}N$ and
\begin{equation}\label{lem:x-sc-eq-criterion}
 \int\operatorname{dist}_\mu(\kappa_x,\mathcal W)^2\,\d\lambda(x)
 \le C_{\rm mean}\,r^pq^{-2p}\gamma .
\end{equation}
\end{lemma}

The space $\mathcal W$ depends on $\lambda$. For $\mathcal I''\subseteq
\mathcal I$ and $C\subseteq\mathbb S$ we write $\mathcal I''(C):=
\{i\in\mathcal I'':u_i\in C\}$ and $\mu(\mathcal I''(C)):=
\sum_{i\in\mathcal I''(C)}\mu_i$. Since $\mathcal W$ is a subspace,
$\operatorname{dist}_\mu(\kappa_x,\mathcal W)=\operatorname{dist}_\mu
(|\langle u_\cdot,x\rangle|^p,\mathcal W)/F_\natural(x)$, and
$1/F_\natural\le r^{\bar p}$ on $\mathbb S$.

\begin{proof}
All constants $C_4,C_5,\ldots$ depend only on $p$.

\emph{Step 1: a global polynomial.} Let $P_Q$ be the polynomial of
Lemma~\ref{lem:t-tv-approx} for $(v,\sigma)=(p,0)$ and degree $Q$, with
constant $C_p$. Define $\bar P(t):=2^pP_Q(t/2)$ and $\bar R(t):=|t|^p-\bar P(t)$
on $[-2,2]$. Then $\sup_{[-2,2]}|\bar P|\le2^p(1+C_p)$ and
$\operatorname{TV}_{[-1,1]}(\bar R)\le\operatorname{TV}_{[-2,2]}(\bar R)
\le2^pC_pQ^{-p}$. By Lemma~\ref{lem:t-threshold-form} there is a signed
measure $\bar\beta$ on $(-1,1]$ with $\|\bar\beta\|\le2^pC_pQ^{-p}$ and
\begin{equation}\label{lem:x-sc-eq-split}
 |t|^p=\bar P(t)+\bar R(-1)+\int\mathbf1_{\{t\ge s'\}}\,\d\bar\beta(s')
 \qquad(|t|\le1).
\end{equation}
We claim that for $c\in[-1,1]$ and real $z$ with $|z|\le\tau$,
\begin{equation}\label{lem:x-sc-eq-taylor}
 |\bar P(c+z)-\sum_{b\le m_{\rm T}}\frac{\bar P^{(b)}(c)}{b!}z^b|
 \le e^{-4L}.
\end{equation}
For complex $|z|\le2\tau$ the point $\zeta:=(c+z)/2$ has
$|\operatorname{Re}\zeta|\le(1+2\tau)/2\le5/6$ and
$|\operatorname{Im}\zeta|\le\tau\le1/3$. Since $3\tau Q\le6r^{-4}\le1$,
Lemma~\ref{lem:x-sc-poly}(a) gives $|\bar P(c+z)|\le2^{p+1}e(1+C_p)(Q+1)
=:M_A$. Cauchy's estimate on the disk of radius $2\tau$ bounds the
Taylor coefficients by $M_A(2\tau)^{-b}$, so the left side
of Eq.~\eqref{lem:x-sc-eq-taylor} is at most $M_A2^{-m_{\rm T}}\le
2^{p+1}e(1+C_p)e^{L/2}e^{-5.5L}\le e^{-4L}$, by
Lemma~\ref{lem:x-sc-params}(b).

\emph{Step 2: concentrated rows and queries.} Define
$\mathcal I_0:=\mathcal I$. As long as some $v\in\mathbb S$ satisfies
$\mu(\mathcal I_k(\operatorname{Cap}(v,\tau)))>m_\circ$, choose such a $v$
as $v_{k+1}$, let $\mathcal U_{k+1}$ be the set of indices
$i\in\mathcal I_k$ with $u_i\in\operatorname{Cap}(v_{k+1},\tau)$, and let
$\mathcal I_{k+1}:=\mathcal I_k\setminus\mathcal U_{k+1}$. Each step removes mass more than $m_\circ$, so this stops
after $n_1\le1/m_\circ$ steps. Define $\mathcal U:=\bigcup_k\mathcal U_k$ and
$\mathcal I':=\mathcal I_{n_1}$. Then
$\mu(\mathcal I'(\operatorname{Cap}(v,\tau)))\le m_\circ$ for every
$v\in\mathbb S$. In the same way remove from $\lambda$, one after another,
the remaining parts of caps $\operatorname{Cap}(x^{(k)},\tau)$ of remaining
mass more than $m_\circ$, $k\le n_2\le1/m_\circ$. Let $\lambda''$ be the
removed part and $\lambda':=\lambda-\lambda''$. Every point of
$\supp\lambda''$ lies in some $\operatorname{Cap}(x^{(k)},\tau)$, and
$\lambda'(\operatorname{Cap}(v,\tau))\le m_\circ$ for every $v$.

We use four families of functions on $\mathcal I$. First, the constant
function $\mathbf1$. Second, for $k\le n_1$, the functions
$i\mapsto\mathbf1_{\mathcal U_k}(i)\Phi(u_i)$ with $\Phi\in
\mathcal P_{m_{\rm T}}$. Third, for $k\le n_2$ and multi-indices
$\beta\in\mathbb N^r$ with $|\beta|\le m_{\rm T}$, the functions
$\chi_{k,\beta}(i):=\bar P^{(|\beta|)}(\langle u_i,x^{(k)}\rangle)
u_i^\beta/\beta!$. Since
$\binom{m_{\rm T}+r}r\le(e(1+9L/r))^r\le\mathsf P^r$, these three families
span at most $1+2\mathsf P^r/m_\circ\le1+2N\mathsf P^{-2r}$ dimensions.
Fourth, a residual net. For $x\in\mathbb S$ let
$\vartheta_x(i):=\int\mathbf1_{\{\langle u_i,x\rangle\ge s'\}}\,
\d\bar\beta(s')$. Apply Lemma~\ref{lem:t-integrated-vectors}(b) with
$X:=\mathbb S$, $\ell:=1$, index set $\mathcal I$, weights
$\tau_i:=\mu_i^{1/2}$, measures $\bar\beta$, the test function
$\mathbf1_{\{\langle u_i,x\rangle\ge s'\}}$ and the identity projection.
Its vectors are $(\mu_i^{1/2}\vartheta_x(i))_i$, with
$D\le\|\bar\beta\|$, and their distances are the distances
$\|\vartheta_x-\vartheta_{x'}\|_\mu$. Let $C_{\rm iv}$ be its constant for
$\ell=1$ and $u_{\rm res}:=C_{\rm iv}2^pC_pQ^{-p}(c_{\rm dim}N/20)^{-1/(2d)}$.
We obtain a set $\mathcal N_{\rm res}\subseteq\mathbb S$ of at most
$c_{\rm dim}N/20$ points such that every $\vartheta_x$ is within
$u_{\rm res}$ of some $\vartheta_{x'}$, $x'\in\mathcal N_{\rm res}$, in the
norm $\|\cdot\|_\mu$. Take the functions $\vartheta_{x'}$ and
$\mathbf1_{\mathcal U}\vartheta_{x'}$ for $x'\in\mathcal N_{\rm res}$.

For every $x\in\mathbb S$ and $i\in\mathcal U_k$ we have
$|\langle u_i-v_k,x\rangle|\le\tau$, and Eq.~\eqref{lem:x-sc-eq-taylor} with
$c=\langle v_k,x\rangle$ gives a polynomial $\Phi_{k,x}\in
\mathcal P_{m_{\rm T}}$ with $|\bar P(\langle u_i,x\rangle)-\Phi_{k,x}(u_i)|
\le e^{-4L}$ on $\mathcal U_k$. With Eq.~\eqref{lem:x-sc-eq-split} and the net,
for every $x\in\mathbb S$ some combination $g^{\mathcal U}_x$ of the
second and fourth families satisfies
\begin{equation}\label{lem:x-sc-eq-removed}
 \|\mathbf1_{\mathcal U}(|\langle u_\cdot,x\rangle|^p-g^{\mathcal U}_x)
 \|_\mu\le e^{-4L}+u_{\rm res}.
\end{equation}
If $x\in\operatorname{Cap}(x^{(k)},\tau)$, define $y:=x-x^{(k)}$. By the
multinomial theorem,
$\sum_{|\beta|\le m_{\rm T}}\chi_{k,\beta}(i)y^\beta=\sum_{b\le m_{\rm T}}
\bar P^{(b)}(\langle u_i,x^{(k)}\rangle)\langle u_i,y\rangle^b/b!$. Since
$|\langle u_i,y\rangle|\le\tau$, Eq.~\eqref{lem:x-sc-eq-taylor},
Eq.~\eqref{lem:x-sc-eq-split} and the net give, for every
$x\in\supp\lambda''$,
\begin{equation}\label{lem:x-sc-eq-heavy}
 \operatorname{dist}_\mu(|\langle u_\cdot,x\rangle|^p,\mathcal W)
 \le e^{-4L}+u_{\rm res}
\end{equation}
for any $\mathcal W$ containing the four families.

\emph{Step 3: cells and weights.} Let $\mathcal Z\subseteq\mathbb S$ be a
maximal set with pairwise distances at least $s$. For $y\in\mathcal Z$ let
$C_y$ be the set of points whose nearest point of $\mathcal Z$ is $y$, with
ties broken by a fixed order of $\mathcal Z$. The cells $C_y$ partition
$\mathbb S$ into Borel sets with
$\{u:\|u-y\|_2<s/2\}\subseteq C_y\subseteq\operatorname{Cap}(y,s)$. By
Lemma~\ref{lem:x-sc-caps}, $\omega(C_y)\ge(s/\pi)^d/(2d)$ and
$|\mathcal Z|\le2d(\pi/s)^d$. Define
\[
 f_y:=\frac{\mu(\mathcal I'(C_y))}{\omega(C_y)},\quad
 g_y:=\frac{\lambda'(C_y)}{\omega(C_y)},\quad E_y:=E(y),\quad
 a_y:=(1\vee f_y)^{\frac1{r+2p}},\quad
 b_y:=(1\vee\frac{g_y}{E_y^2})^{\frac1{r+2p}},
\]
and let $a$ and $b$ be the functions on $\mathbb S$ equal to $a_y$ and
$b_y$ on $C_y$. We claim:
\begin{enumerate}[label=(\roman*),nosep,leftmargin=*]
\item $E_y/2\le E(x)\le2E_y$ for $x\in C_y$;
\item $a_y\le\mathsf P^{C_3}$ and $b_y\le\mathsf P^{C_3}$;
\item $\|a^d\|_{r/d}\le2$ and $\|b^d\|_{r/d}\le C_b$, with $C_b$
depending only on $p$.
\end{enumerate}
For (i) there is nothing to prove above two. Let $\|x-x'\|_2\le s$ for unit
$x,x'$, and recall $F\ge1/r$ on $\mathbb S$. Since $p\le2$, Jensen's
inequality and Eq.~\eqref{lem:t-eq-lewis} give
$F(x-x')\le(\sum_i\nu_i\langle u_i,x-x'\rangle^2)^{p/2}\le
(s^2/r)^{p/2}$. If $1\le p<2$, Minkowski's inequality (Lemma~\ref{lem:pre_minkowski}(a)) in $L_p(\nu)$ gives
$F(x)^{1/p}\le F(x')^{1/p}+s/\sqrt r\le F(x')^{1/p}(1+s\sqrt r)$, and
$(1+s\sqrt r)^p\le2$. If $p<1$, then $F(x)\le F(x')+F(x-x')\le
F(x')(1+r^{1-p/2}s^p)\le2F(x')$, since $s\le r^{-2/p}/4$. By symmetry this
proves (i). For (ii), Lemma~\ref{lem:x-sc-caps}(c) covers $C_y$ by at
most $2d(\pi^2s/\tau)^d$ caps of radius $\tau$, each of $\mu$-mass at most
$m_\circ$ on $\mathcal I'$, so
\[
 f_y\le4d^2(\pi^3/\tau)^dm_\circ\le4d^2(\pi^3c_qr^{D_p+1})^d\mathsf P^{3r}
 \le\mathsf P^{(D_p+5)r},
\]
where we used $\tau^{-d}=r^{D_pd}q^d\le(c_qr^{D_p+1})^dN$ and $\pi^3c_q<1$.
The same bound holds for $g_y$, and $E_y^{-2}\le r^2\le\mathsf P$. This
gives (ii). For (iii) define $\alpha:=r/(r+2p)<1$, so $a^r=(1\vee f)^\alpha$
with $f:=\sum_yf_y\mathbf1_{C_y}$. Since $\int f\,\d\omega=
\mu(\mathcal I')\le1$, Jensen's inequality gives
$\int a^r\le1+(\int f)^\alpha\le2$. Let $E_c:=\sum_yE_y\mathbf1_{C_y}$ and
$g:=\sum_yg_y\mathbf1_{C_y}$. By H\"older's inequality with exponents
$1/\alpha$ and $1/(1-\alpha)$, and since $2\alpha/(1-\alpha)=r/p$,
\[
 \int b^r\,\d\omega\le1+(\int g)^\alpha
 (\int E_c^{-r/p}\,\d\omega)^{2p/(r+2p)}
 \le1+(2C_{\rm inv})^{2r/(r+2p)}\le1+(1+2C_{\rm inv})^2 .
\]
Here we used $E_c\ge E/2$ from (i) and $\|E^{-1}\|_{L_{r/p}(\omega)}\le
C_{\rm inv}$, where $C_{\rm inv}$ is the constant of
Lemma~\ref{lem:t-inverse-energy} below two and $C_{\rm inv}:=1$ above two.
Since $\|a^d\|_{r/d}=(\int a^r)^{d/r}$ and similarly for $b$, (iii) follows
with $C_b:=1+(1+2C_{\rm inv})^2$.

\emph{Step 4: features on the remaining rows.} Call a pair
$(y,y')\in\mathcal Z^2$ near if $|\langle y,y'\rangle|\le10s$ and far
otherwise. For $u\in C_y$ and $x\in C_{y'}$,
\begin{equation}\label{lem:x-sc-eq-cellscore}
 |\langle u,x\rangle-\langle y,y'\rangle|\le\|u-y\|_2+\|x-y'\|_2\le2s .
\end{equation}
Let $\mathcal B$ be a partition of $[-s,s]$ into $\lceil2/s\rceil\le3/s$
intervals of length at most $s^2$. For a near pair $(y,y')$ and
$I\in\mathcal B$ let $\mathcal I'_{y,y',I}:=\{i\in\mathcal I'(C_y):
\langle u_i-y,y'\rangle\in I\}$, which covers $\mathcal I'(C_y)$ since
$\|u_i-y\|_2\le s$. Define $m_{yy'}:=\lceil qs^2a_yb_{y'}\rceil$ and
$h_{yy'}:=1/(qa_yb_{y'})$. The far features are the functions
$\mathbf1_{\mathcal I'(C_y)}(i)\Phi(u_i)$, $y\in\mathcal Z$,
$\Phi\in\mathcal P_{m_{\rm T}}$, and the near features are the functions
$\mathbf1_{\mathcal I'_{y,y',I}}(i)\Phi(u_i)$, $(y,y')$ near,
$I\in\mathcal B$, $\Phi\in\mathcal P_{m_{yy'}}$. Let $\mathcal W$ be the
span of all features of Steps~2 and~4.

By Eq.~\eqref{lem:x-sc-eq-dim} and Lemma~\ref{lem:x-sc-params}(e), the far
features span at most $|\mathcal Z|d_{m_{\rm T}}\le4d(\pi/s)^d
(e+9eL/d)^d\le4dN^{1/2}$ dimensions. For the near features note
$m_{yy'}\ge qs^2\ge d$ by Lemma~\ref{lem:x-sc-params}(d), so
$d_{m_{yy'}}\le2(2em_{yy'}/d)^d\le2(4eqs^2a_yb_{y'}/d)^d$. By
Eq.~\eqref{lem:x-sc-eq-cellscore}, near pairs satisfy $|\langle u,x\rangle|\le
12s$ on $C_y\times C_{y'}$, hence
\[
 \sum_{\text{near}}a_y^db_{y'}^d\le\frac{4d^2\pi^{2d}}{s^{2d}}\sum_{\text{near}}
 a_y^db_{y'}^d\omega(C_y)\omega(C_{y'})\le\frac{4d^2\pi^{2d}}{s^{2d}}
 \mathcal M_{a^d,b^d}(12s)\le\frac{4d^2\pi^{2d}}{s^{2d}}\,1536\,C_b\,s\sqrt r,
\]
by Lemma~\ref{lem:x-sc-section}(b) and (iii). The near features therefore
span at most
\[
 \frac3s\cdot2(\frac{4eqs^2}d)^d\sum_{\text{near}}a_y^db_{y'}^d
 \le C_4d^2\sqrt r(\frac{4e\pi^2q}d)^d\le C_4d^2\sqrt r
 (8e\pi^2c_q)^dN=C_4d^2\sqrt r\,2^{-d}N
\]
dimensions, using $q\le c_qrY\le2c_qdY$. With Step~2 and
$|\mathcal N_{\rm res}|\le c_{\rm dim}N/20$,
$\dim\mathcal W\le1+2N\mathsf P^{-2r}+c_{\rm dim}N/10+4dN^{1/2}
+C_4d^2\sqrt r\,2^{-d}N\le c_{\rm dim}N$ for $r\ge r_0$.

\emph{Step 5: pointwise errors on the remaining rows.} Let
$x\in\supp\lambda'$, say $x\in C_{y'}$, and let $i\in\mathcal I'(C_y)$. If
$(y,y')$ is far, define $c:=\langle y,y'\rangle$. By
Eq.~\eqref{lem:x-sc-eq-cellscore}, $z:=(\operatorname{sgn}c)\langle u_i,x\rangle$
lies in $[\frac45|c|,\frac65|c|]$. Let $T_{|c|}$ be the Taylor polynomial of
degree $m_{\rm T}$ of $z\mapsto z^p$ at $|c|$. The residual vanishes at
$|c|$, so Lemma~\ref{lem:t-shifted-tv}(b) with $\rho=1/5$ bounds
$||\langle u_i,x\rangle|^p-T_{|c|}((\operatorname{sgn}c)\langle u_i,x
\rangle)|$ by $C'_p(m_{\rm T}+2)^{\lceil p\rceil}5^{-m_{\rm T}}\le e^{-4L}$, since
$m_{\rm T}\ge8L$ and $L$ is large.
As a function of $u_i$ this approximant lies in $\mathcal P_{m_{\rm T}}$. If
$(y,y')$ is near and $\langle u_i-y,y'\rangle\in I$, let $v_I$ be the
midpoint of $I$ and $\alpha:=\langle y,x\rangle+v_I$, which depends only on
$x,y,y',I$. Then
\[
 \langle u_i,x\rangle-\alpha=(\langle u_i-y,y'\rangle-v_I)+
 \langle u_i-y,x-y'\rangle\in[-\tfrac32s^2,\tfrac32s^2].
\]
Define $m:=m_{yy'}$ and $h:=h_{yy'}$. If $|\alpha|\le8s^2$, let
$w:=12s^2$, so $|\langle u_i,x\rangle|\le w$, and use
$w^pQ_m(\langle u_i,x\rangle/w)$ with $Q_m$ from
Lemma~\ref{lem:x-sc-localized}, whose constant we call $C_{\rm loc}$.
Since $m/w\ge1/(12h)$, its error is at most
\[
 C_{\rm loc}\,w^pm^{-p}(1+\frac{m|\langle u_i,x\rangle|}w)^{-1}
 \le12^{p+1}C_{\rm loc}\,h^p(1+\frac{|\langle u_i,x\rangle|}h
 )^{-1}.
\]
If
$|\alpha|>8s^2$, then $(\operatorname{sgn}\alpha)\langle u_i,x\rangle\in
[\frac34|\alpha|,\frac54|\alpha|]$, and the Taylor polynomial of degree $m$
of $z^p$ at $|\alpha|$ has error at most
$C'_p2^p(m+2)^{\lceil p\rceil}4^{-m}\le e^{-4L}$ by
Lemma~\ref{lem:t-shifted-tv}(b), since $m\ge qs^2\ge8L$ and $L$ is large.
In both cases the
approximant is a polynomial of degree at most $m$ in $u_i$. Combining these
choices on the sets $\mathcal I'(C_y)$ and $\mathcal I'_{y,y',I}$ gives
$g_x\in\mathcal W$ supported on $\mathcal I'$ with
\begin{equation}\label{lem:x-sc-eq-pointwise}
 \mathcal E(x)^2:=\|\mathbf1_{\mathcal I'}(|\langle u_\cdot,x\rangle|^p-g_x)
 \|_\mu^2\le C_5\sum_{y:\,(y,y')\text{ near}}\sum_{i\in\mathcal I'(C_y)}
 \mu_ih_{yy'}^{2p}(1+\frac{|\langle u_i,x\rangle|}{h_{yy'}})^{-2}
 +e^{-8L}.
\end{equation}

\emph{Step 6: a bin estimate.} For integers $k,l\ge0$ let
$\mathcal Z_k:=\{y:2^k\le a_y<2^{k+1}\}$,
$\mathcal Z'_l:=\{y:2^l\le b_y<2^{l+1}\}$ and $h_{kl}:=2^{-k-l}/q$. Let
$\mu_k$ be the measure $\sum_{y\in\mathcal Z_k}\sum_{i\in\mathcal I'(C_y)}
\mu_i\delta_{u_i}$ on $\mathbb S$ and $\lambda_l:=\sum_{y'\in\mathcal Z'_l}
E_{y'}^{-2}\lambda'|_{C_{y'}}$. Their coarse densities are
$f_k:=\sum_{y\in\mathcal Z_k}f_y\mathbf1_{C_y}$ and
$g_l:=\sum_{y'\in\mathcal Z'_l}(g_{y'}/E_{y'}^2)\mathbf1_{C_{y'}}$. We claim
\begin{equation}\label{lem:x-sc-eq-bin}
 I_{kl}:=\iint(1+\frac{|\langle u,x\rangle|}{h_{kl}})^{-2}
 \d\mu_k(u)\,\d\lambda_l(x)\le13\,h_{kl}\sqrt r(1+
 \frac{\mathcal M_{f_k,g_l}(2\varsigma)}{\varsigma\sqrt r}).
\end{equation}
We may assume $\mathcal Z_k,\mathcal Z'_l\ne\emptyset$, so by (ii)
$2^{k+l}\le\mathsf P^{2C_3}$. Write $h:=h_{kl}$ and
$\tau':=\mathsf P^{-2C_3-1}/q$. Then $10\tau'\le h$ and $\tau'\le\tau$.

\emph{Smoothing.} Let $\varkappa_{\tau'}(z,u):=\mathbf1_{\{\|z-u\|_2\le
\tau'\}}/\omega_r(\tau')$, and let $f^*(z):=\int\varkappa_{\tau'}(z,u)\,
\d\mu_k(u)$ and $g^*(z):=\int\varkappa_{\tau'}(z,x)\,\d\lambda_l(x)$. If
$\|z-u\|_2,\|x'-x\|_2\le\tau'$, then
$|\langle z,x'\rangle-\langle u,x\rangle|\le2\tau'\le h/5$, so
$1+|\langle u,x\rangle|/h\ge\frac45(1+|\langle z,x'\rangle|/h)$. Averaging
over $z$ and $x'$ gives $I_{kl}\le\frac{25}{16}I^*$ with
$I^*:=\iint(1+|\langle z,x'\rangle|/h)^{-2}f^*(z)g^*(x')\,\d\omega(z)\,
\d\omega(x')$. The caps of radius $\tau'\le\tau$ have $\mu$-mass at most
$m_\circ$ on $\mathcal I'$ and $\lambda_l$-mass at most $r^2m_\circ$, and
the total masses are at most $1$ and $r^2$. Hence, with
Lemma~\ref{lem:x-sc-caps}(a) and $q^d\le(c_qr)^dN$,
\[
 \Xi:=\|f^*\|_2\|g^*\|_2\le\frac{r^2m_\circ}{\omega_r(\tau')}\le
 2dr^2\mathsf P^{3r}(\frac{\pi c_qr\mathsf P^{2C_3+1}}2)^d
 \le\mathsf P^{(2C_3+6)r}.
\]

\emph{Flatness of the fine pair density.} By
Lemmas~\ref{lem:x-sc-pairs}(c) and~\ref{lem:x-sc-pair-density},
$I^*=\int(1+|t|/h)^{-2}\phi_r(t)\varphi^*(t)\,\d t$ with
$\varphi^*:=\varphi_{f^*,g^*}\ge0$ and $\int\phi_r\varphi^*=\int f^*\int g^*
\le r^2$. Let $k_0:=\lfloor(r-4)/8\rfloor\ge r/10$ and
$\varsigma_*:=\mathsf P^{-D_\flat}\le r^{-1/2}$. The Taylor polynomial of
$\varphi^*$ at $0$ of degree $k_0-1$ differs from $\varphi^*$ on
$[-\varsigma_*,\varsigma_*]$ by at most
$2(8\varsigma_*)^{k_0}\Xi\le2\mathsf P^{-(D_\flat-1)r/10}
\mathsf P^{(2C_3+6)r}\le1$, by Lemma~\ref{lem:x-sc-pair-density}. Adding
one gives a polynomial $U$ of degree at most $k_0-1$ with
$\varphi^*\le U\le\varphi^*+2$ on $[-\varsigma_*,\varsigma_*]$. Since
$\phi_r\ge\sqrt r/8$ there, $\int_{-\varsigma_*}^{\varsigma_*}\varphi^*\le
8r^{3/2}$ and $\int_{-\varsigma_*}^{\varsigma_*}U\le4\varsigma_*+8r^{3/2}$.
Lemma~\ref{lem:x-sc-poly}(c) applied to $U(\varsigma_*y)$ gives
$M_U:=\sup_{|t|\le\varsigma_*/2}U\le(k_0+2)(4+8r^{3/2}\varsigma_*^{-1})
\le\mathsf P^{D_\flat+2}$. Lemma~\ref{lem:x-sc-poly}(b) applied to
$U(\varsigma_*y/2)$ gives $|U'(t)|\le3k_0^2M_U/\varsigma_*\le
\mathsf P^{2D_\flat+5}$ for $|t|\le\varsigma_*/4$. Since
$\varsigma\le\mathsf P^{-2D_\flat-6}\le\varsigma_*/4$, for
$t,t'\in[-\varsigma,\varsigma]$ we get $|\varphi^*(t)-\varphi^*(t')|\le
|U(t)-U(t')|+2\le3$. As $\phi_r\ge\sqrt r/8$ on $[-\varsigma,\varsigma]$,
\[
 \sup_{|t|\le\varsigma}\varphi^*\le3+\frac1{2\varsigma}
 \int_{-\varsigma}^{\varsigma}\varphi^*\le3+
 \frac{4\mathcal M_{f^*,g^*}(\varsigma)}{\varsigma\sqrt r}.
\]

\emph{Conclusion of the bin estimate.} Since $\phi_r\le\sqrt r$ and
$\int(1+|t|/h)^{-2}\d t=2h$, the part of $I^*$ with $|t|\le\varsigma$ is at
most $2h\sqrt r(3+4\mathcal M_{f^*,g^*}(\varsigma)/(\varsigma\sqrt r))$.
The rest is at most $(h/\varsigma)^2r^2\le h\sqrt r$, because $h\le1/q$
and $q\varsigma^2\ge256r^2$ by Lemma~\ref{lem:x-sc-params}(d). Finally,
let $z,x'$ be as in the smoothing and let $z^\circ,x^\circ$ be points of
the cells containing $u$ and $x$. Then
$|\langle z,x'\rangle-\langle z^\circ,x^\circ\rangle|\le2(\tau'+2s)\le
\varsigma$. So $\mathbf1_{\{|\langle z,x'\rangle|\le\varsigma\}}\le
\mathbf1_{\{|\langle z^\circ,x^\circ\rangle|\le2\varsigma\}}$. Integrating
$z,x'$ against $\varkappa_{\tau'}$, $z^\circ,x^\circ$ uniformly over the
cells, and $u,x$ against $\mu_k\otimes\lambda_l$ gives
$\mathcal M_{f^*,g^*}(\varsigma)\le\mathcal M_{f_k,g_l}(2\varsigma)$. This
proves Eq.~\eqref{lem:x-sc-eq-bin}.

\emph{Step 7: summation.} On $C_{y'}$ we have
$F_\natural(x)^{-2}=r^pE(x)^{-2}\le4r^pE_{y'}^{-2}$ by (i). The function
$h\mapsto h^{2p}(1+|t|/h)^{-2}=h^{2p+2}(h+|t|)^{-2}$ is nondecreasing, and
$h_{yy'}\le h_{kl}$ for $y\in\mathcal Z_k$, $y'\in\mathcal Z'_l$. Integrate
Eq.~\eqref{lem:x-sc-eq-pointwise} divided by $F_\natural(x)^2$ against
$\lambda'$ and include all pairs. Since
$\sum_{y'}E_{y'}^{-2}\lambda'(C_{y'})\le r^2$,
\[
 \int\frac{\mathcal E(x)^2}{F_\natural(x)^2}\,\d\lambda'(x)\le
 4C_5r^p\sum_{k,l\ge0}h_{kl}^{2p}I_{kl}+4r^{p+2}e^{-8L}.
\]
By Eq.~\eqref{lem:x-sc-eq-bin}, $\sum_{k,l}h_{kl}^{2p}I_{kl}\le13\sqrt r
\sum_{k,l}h_{kl}^{2p+1}+13\varsigma^{-1}q^{-2p-1}\mathcal M_{F_1,G_1}
(2\varsigma)$, where by bilinearity
$F_1:=\sum_k2^{-k(2p+1)}f_k$ and $G_1:=\sum_l2^{-l(2p+1)}g_l$. The first
sum is at most $4q^{-2p-1}$. On $C_y$ with $y\in\mathcal Z_k$ we have
$f_y\le a_y^{r+2p}=a_y^da_y^{2p+1}$ and $a_y<2^{k+1}$, so
$F_1\le2^{2p+1}a^d$. In the same way $G_1\le2^{2p+1}b^d$. Since
$2\varsigma\le1/(4\sqrt r)$, Lemma~\ref{lem:x-sc-section}(b) and (iii) give
$\mathcal M_{F_1,G_1}(2\varsigma)\le4^{2p+1}\cdot128\,\varsigma\sqrt r\cdot
2C_b$. Hence
$\int\mathcal E^2F_\natural^{-2}\,\d\lambda'\le C_6r^pq^{-2p-1}\sqrt r+
4r^{p+2}e^{-8L}$.

Now let $x\in\supp\lambda'$. The function $g^{\mathcal U}_x+g_x$ lies in
$\mathcal W$, and Eq.~\eqref{lem:x-sc-eq-removed} gives
\[
 \operatorname{dist}_\mu(\kappa_x,\mathcal W)^2\le F_\natural(x)^{-2}
 (2(e^{-4L}+u_{\rm res})^2+2\mathcal E(x)^2).
\]
For
$x\in\supp\lambda''$ we use Eq.~\eqref{lem:x-sc-eq-heavy}. Since
$F_\natural^{-2}\le r^{2\bar p}$ and $2\bar p\le p+2$,
\[
 \int\operatorname{dist}_\mu(\kappa_x,\mathcal W)^2\,\d\lambda\le
 8r^{p+2}(e^{-8L}+u_{\rm res}^2)+2C_6r^pq^{-2p}\gamma+8r^{p+2}e^{-8L}.
\]
By Lemma~\ref{lem:x-sc-params}(f), $r^{p+2}e^{-8L}\le r^pq^{-2p}\gamma$.
Also $1/Y\le c_qr/q$ and $Q^{-2p}\le q^{-2p}r^{-6p/\theta_p}$, so
$r^{p+2}u_{\rm res}^2\le C_7r^{p+3-6p/\theta_p}q^{-2p-1}\le
C_7r^pq^{-2p}\gamma$, because $6p/\theta_p\ge6$. This
proves Eq.~\eqref{lem:x-sc-eq-criterion}.
\end{proof}

\subsection{One Round of Rounding}

\begin{lemma}[One round of rounding]\label{lem:x-sc-round}
There is $C_{\rm rd}$, depending only on $p$, such that the following
holds under the standing assumptions, for every $C_{\rm sp}\ge1$. In the
Lewis setup with row set $I$, let $w\in[0,\infty)^I$ be a probability
vector with $n:=|\supp w|>N$, and let $J:=\{i:0<w_i\le2/n\}$. There are a
subspace $H\subseteq\R^J$ of codimension at most $\delta_0|J|$ and a closed
symmetric convex set $\mathcal K\subseteq H$ with
$\gamma_H(\mathcal K)\ge e^{-\delta_0|J|}$ such that every
$z\in\mathcal K$ satisfies $\sum_{i\in J}w_iz_i=0$ and
\begin{equation}\label{lem:x-sc-eq-round}
 |\sum_{i\in J}w_iz_i|\langle u_i,x\rangle|^p|\le C_{\rm rd}
 (\frac{r^{p+1}\Lambda q^{-2p}\gamma}n)^{1/2}F(x)\qquad(x\in\R^r).
\end{equation}
\end{lemma}

\begin{proof}
By Lemma~\ref{lem:t-light-rows}, $|J|\ge n/2$ and
$\sum_{i\in J}w_i^2\le2/n$. Define $\mu_i:=nw_i^2/4$ for $i\in J$, so
$\sum_i\mu_i\le1/2$, and for unit $x$ let
$v_x:=(w_i\kappa_x(i))_{i\in J}$, with $\kappa_x$ as in
Lemma~\ref{lem:x-sc-criterion}. For a subspace $\mathcal W$ of $\R^J$ let
$\Pi_{\mathcal W}$ be the orthogonal projection onto
$\{(w_if(i))_i:f\in\mathcal W\}$. Then
$\|(\mathrm{Id}-\Pi_{\mathcal W})v_x\|_2^2=\min_{f\in\mathcal W}\sum_iw_i^2
(\kappa_x(i)-f(i))^2=\frac4n\operatorname{dist}_\mu(\kappa_x,\mathcal W)^2$.
Define $\rho^2:=4C_{\rm mean}r^pq^{-2p}\gamma/n$.

\emph{One projection for all queries.} Let $\mathfrak K$ be the compact
convex set of symmetric $J\times J$ matrices $B$ with $0\preceq B\preceq\mathrm{Id}$
and $\tr[B]\le c_{\rm dim}N$, and for unit $x$ let
$\mathfrak K_x:=\{B\in\mathfrak K:v_x^\top(\mathrm{Id}-B)v_x\le\rho^2\}$, a closed
convex set. Let $x_1,\ldots,x_m$ be unit vectors. The set
$\mathcal S:=\{(v_{x_j}^\top(\mathrm{Id}-B)v_{x_j})_{j\le m}:B\in\mathfrak K\}$ is
compact and convex. If it missed the closed convex set
$(-\infty,\rho^2]^m$, a strictly separating vector $y\ne0$ would satisfy
$\langle y,a\rangle>c\ge\langle y,b\rangle$ for $a\in\mathcal S$ and
$b\le\rho^2$ coordinatewise. Letting coordinates of $b$ tend to $-\infty$
gives $y\ge0$, and then $c\ge\rho^2\sum_jy_j$. With
$\lambda:=\sum_jy_j\delta_{x_j}/\sum_jy_j$, every $B\in\mathfrak K$ would
have $\sum_j\lambda(x_j)v_{x_j}^\top(\mathrm{Id}-B)v_{x_j}>\rho^2$. This fails for
$B=\Pi_{\mathcal W}$ with $\mathcal W$ from
Lemma~\ref{lem:x-sc-criterion} for this $\lambda$, applied with
$\mathcal I:=J$. Hence finitely many sets $\mathfrak K_x$ always
intersect, and by compactness of $\mathfrak K$ there is
$B_0\in\bigcap_{x\in\mathbb S}\mathfrak K_x$. Let $\Pi$ be the orthogonal
projection onto the span of the eigenvectors of $B_0$ with eigenvalue at
least $1/2$. Then $\operatorname{rank}\Pi\le2c_{\rm dim}N$ and
$\mathrm{Id}-\Pi\preceq2(\mathrm{Id}-B_0)$, so $\|(\mathrm{Id}-\Pi)v_x\|_2^2\le2\rho^2$ for every unit
$x$.

\emph{The subspace.} Let $P_q$ be the polynomial of
Lemma~\ref{lem:t-tv-approx} for $(v,\sigma)=(p,0)$ and degree $q$, define
$R_q:=|t|^p-P_q$, and by Lemma~\ref{lem:t-threshold-form} write
$R_q(t)=R_q(-1)+\int\mathbf1_{\{t\ge s'\}}\,\d\beta_q(s')$ on $[-1,1]$ with
$\|\beta_q\|\le C_pq^{-p}$. Let
\[
 H:=\{z\in\R^J:\Pi z=0,\ \sum_{i\in J}w_iz_i\Phi(u_i)=0\text{ for all }
 \Phi\in\mathcal P_q\}.
\]
Taking $\Phi=1$ shows $\sum_iw_iz_i=0$ on $H$. By Eq.~\eqref{lem:x-sc-eq-dim},
$d_q\le2(2eq/d)^d\le2(4ec_q)^dN\le c_{\rm dim}N$, so the codimension of $H$
is at most $3c_{\rm dim}N\le\delta_0|J|$. For unit $x$ the function
$u\mapsto P_q(\langle u,x\rangle)+R_q(-1)$ lies in $\mathcal P_q$. Hence for
$z\in H$,
\begin{equation}\label{lem:x-sc-eq-resid}
 \sum_{i\in J}w_iz_i|\langle u_i,x\rangle|^p=\langle z,\tilde v_x\rangle,
 \qquad\tilde v_x:=(w_i\int\mathbf1_{\{\langle u_i,x\rangle\ge s'\}}\,
 \d\beta_q(s'))_{i\in J},
\end{equation}
and therefore $\Pi_Hv_x=\Pi_H\tilde v_x/F_\natural(x)$, where $\Pi_H$ is
the orthogonal projection onto $H$. Since $H\subseteq\ker\Pi$, also
$\Pi_Hv_x=\Pi_H(\mathrm{Id}-\Pi)v_x$, so $\|\Pi_Hv_x\|_2\le\sqrt2\rho$.

\emph{Chaining over energy tiers.} For
$0\le l\le l_{\max}:=\lceil\bar p\log_2r\rceil$ let $t_l:=2^lr^{-\bar p}$
and $\mathcal X_l:=\{x\in\mathbb S:t_l\le F_\natural(x)\le
2t_l\}$. These sets cover $\mathbb S$, and empty sets are skipped. Let
$\mathcal V_l:=\{t_l^{-1}\Pi_H\tilde v_x:x\in\mathcal X_l\}$. Its
points have norm $(F_\natural(x)/t_l)\|\Pi_Hv_x\|_2\le2\sqrt2\rho=:R_0$.
Apply Lemma~\ref{lem:t-integrated-vectors}(b) with $X:=\mathcal X_l
\subseteq\mathbb S$, $\ell:=1$, index set $J$, weights
$\tau_i:=w_i/t_l$, measures $\beta_q$, the test function
$\mathbf1_{\{\langle u_i,x\rangle\ge s'\}}$ and the projection $\Pi_H$. Here
$D\le C_pq^{-p}r^{\bar p}(2/n)^{1/2}=:D_0$, and
$N(\mathcal V_l,u)\le\max\{1,(C_{\rm iv}D_0/u)^{2d}\}$. Since
$2\bar p-p\le2$ and $\gamma^{-1}=q/\sqrt r\le c_q\sqrt r\,Y$,
\[
 (\frac{C_{\rm iv}D_0}{R_0})^{2d}=(\frac{C_{\rm iv}^2C_p^2
 r^{2\bar p-p}}{16C_{\rm mean}\gamma})^d\le(C_8r^{5/2})^dN\le
 e^{C_9r\Lambda}n ,
\]
using $N<n$ in the last step. Hence $N_j:=N(\mathcal V_l,2^{-j/2}R_0)\le e^{C_9r\Lambda}n2^{dj}$. Apply
Lemma~\ref{lem:t-sidak-chaining}(c) with
$\lambda_j:=2((C_9+1)r\Lambda+dj)^{1/2}$. It gives strips in $H$ whose
failure probabilities sum to at most
$\sum_jN_je^{-\lambda_j^2/2}\le ne^{-r\Lambda}\sum_j(2e^{-2})^{dj}\le
ne^{-r\Lambda}$, and on whose intersection
$\sup_{v\in\mathcal V_l}|\langle G,v\rangle|\le2R_0\sum_j\lambda_j2^{-j/2}
\le C_{10}R_0(r\Lambda)^{1/2}$. This bound is the deterministic
inequality Eq.~\eqref{lem:t-eq-chain}, so it holds for every vector of $H$ in
the intersection of the strips. Let $\mathcal K$ be the intersection of
the strips of all tiers. It is closed, symmetric and convex. The total
failure probability is at most $(l_{\max}+1)ne^{-r\Lambda}\le\delta_0n/8$,
and each strip fails with probability at most $1/2$, so
Lemma~\ref{lem:t-sidak-chaining}(d) gives $\gamma_H(\mathcal K)\ge
e^{-\delta_0n/4}\ge e^{-\delta_0|J|}$.

\emph{The error.} Let $z\in\mathcal K$ and $x\in\mathcal X_l$.
By Eq.~\eqref{lem:x-sc-eq-resid},
$\sum_iw_iz_i|\langle u_i,x\rangle|^p=\langle z,\Pi_H\tilde v_x\rangle=
t_l\langle z,t_l^{-1}\Pi_H\tilde v_x\rangle$. As
$t_l\le F_\natural(x)\le F(x)$, its modulus is at most
\[
 t_lC_{10}R_0(r\Lambda)^{1/2}\le F(x)\,C_{10}2\sqrt2\,
 (\frac{4C_{\rm mean}r^{p+1}\Lambda q^{-2p}\gamma}n)^{1/2}.
\]
Both sides of Eq.~\eqref{lem:x-sc-eq-round} are $p$-homogeneous in $x$, which gives the
claim with $C_{\rm rd}:=4\sqrt2\,C_{10}C_{\rm mean}^{1/2}$.
\end{proof}

%% file: 82_proof.tex
\section{Growing Rank Above Two}
\label{sec:joint_above_two}

\subsection{Main Result}

Theorem~\ref{thm:above_two_main} bounds the worst-case support at
subexponential accuracy for $p>2$.
The proof of part~(a) rounds the Lewis masses by iterated partial coloring, as in the
upper bound of Theorem~\ref{thm:below_two_main}(b), but with a single
residual on the whole sphere. Write $|t|^p=P(t)+R(t)$ with an even
polynomial $P$ of degree $q$ of order $L$. Linear constraints taken from a
tensor sketch make the change of the polynomial part negligible on a norming
net, hence on the sphere. The residual $R$ has total variation
$O_p(q^{-p})$, and the entropy of integrated affine tests bounds the
Gaussian width of its contribution. Above two the energy of a unit query is
at least $r^{-p/2}$, which turns the absolute error into the factor
$r^{p+1}$.
Parts~(b) and~(c) are lower bounds for non-even $p$. Part~(b)
amplifies the adaptive-shift input of Lemma~\ref{lem:s-shift} by
Lemma~\ref{lem:lp_lower_tensor_product}, and part~(c) uses overlapping
transverse groups (Lemma~\ref{lem:x-a2-overlap}).

\begin{theorem}[Subexponential accuracy above two]
\label{thm:above_two_main}
Fix $p>2$.
\begin{enumerate}[label=(\alph*),nosep]
\item For $C_p\log(2r)\le L\le r$,
\[
 S_p(r,\epsilon)\lesssim_p \epsilon^{-2}r^{p+1}L^{-2p}.
\]
\item Let $p$ be non-even. For $r\ge C_p L$ define
\[
 G_*=\begin{cases}
      \min\{L^{p/2},r/L\},&2<p<4,\\
      1,&p>4 .
      \end{cases}
\]
Then
\[
 S_p(r,\epsilon)\gtrsim_p
 \frac{r^{p/2}}{\epsilon^2 L^{(3p-1)/2}}
 G_*^{2-p/2}(1+\log\frac{L^{p/2}}{G_*}).
\]
\item Let $p$ be non-even. For sufficiently large $r$
and $C_p\le L\le r/4$,
\[
 S_p(r,\epsilon)\gtrsim_p
     \epsilon^{-2}r^{p/2}L^{-3(p-1)/2}
        \min\{1,r/L^2\}
        (1+\log\max\{1,L/\sqrt r\}).
\]
The lower bound holds even with arbitrary real weights on the
selected original rows.
\end{enumerate}
\end{theorem}
\begin{proof}
\emph{Part (a).}
We take $C_p\ge2$. Then $r\ge L\ge2\log 2>1$, so $r\ge2$. Zero rows do not
matter, so we may use the Lewis setup Eq.~\eqref{lem:t-eq-lewis} of the common
tools, with rows $i\in I$, unit directions $u_i\in\R^r$ and masses
$\nu_i$. Since $p\ge2$, $F(x)\ge r^{-p/2}$ for unit $x$. By homogeneity
and the last sentence of the subsection on the Lewis setup, it suffices to
find $\mu^*\in[0,\infty)^I$ with $|\supp\mu^*|\le N_*$, for the number
$N_*$ chosen in Step~3, such that
\begin{equation}\label{lem:x-a2-eq-goal}
 |\sum_{i\in I}\mu^*_i|\langle u_i,x\rangle|^p-F(x)|\le\epsilon
 F(x)\qquad(x\in S^{r-1}).
\end{equation}

\emph{Parameters.} Let $c_2\in(0,1]$ be the absolute constant fixed in
Step~4, and define
\[
 q:=\lfloor c_2L\rfloor,\qquad \tau_{\rm sk}:=(2r)^{-p-1}.
\]
Step~4 shows $1\le q\le r$. Let $P$ be the even polynomial of degree at
most $q$ given by Lemma~\ref{lem:t-tv-approx} with $(v,\sigma)=(p,0)$, write
$C_{\rm tv}$ for the constant $C_p$ of that lemma, and define
$R(t):=|t|^p-P(t)$. Then
\[
 \sup_{|t|\le1}|R(t)|\le\operatorname{TV}_{[-1,1]}(R)\le C_{\rm tv}q^{-p},
 \qquad \sup_{|t|\le1}|P(t)|\le M_P:=1+C_{\rm tv},
\]
and the absolute monomial coefficients of $P/M_P$ sum to at most
$3^{q+1}$. Let $\mathcal N\subseteq S^{r-1}$ be the norming set of
Lemma~\ref{lem:t-norming-set}(a) for degree $q$, so that
$|\mathcal N|\le(8(q+1)^2)^r$.

\emph{Step 1: the sketch.} Apply Lemma~\ref{lem:t-sketch} with $a=1$,
$r'=r$, the rows $i\in I$, $w_{i,1}=u_i$, $P_i=P/M_P$, $Y=\mathcal N$,
$\mathsf A=3^{q+1}$ and $\tau=\tau_{\rm sk}$. Since $q+1\le2r$, we have
$\log(2|\mathcal N|)\le1+r\log(32r^2)\le6r\log(2r)$, and hence
$\log\log(e|\mathcal N|)\le3\log(2r)$. Also $\log\mathsf A\le3q$ and
$\log(1/\tau_{\rm sk})=(p+1)\log(2r)$. We obtain a matrix $\widetilde W$
indexed by $I\times\mathcal N$ of rank at most $d_{\rm sk}$ such that
\[
 |P(\langle u_i,y\rangle)-M_P\widetilde W_{i,y}|\le M_P\tau_{\rm sk}
 \quad(i\in I,\ y\in\mathcal N),\qquad
 \log(2d_{\rm sk})\le C_{\rm sk}(q+(p+5)\log(2r)),
\]
where $C_{\rm sk}\ge1$ is absolute.

\emph{Step 2: one round.} Let $\mu\in[0,\infty)^I$ have total mass one and
$N:=|\supp\mu|$ positive entries. Let $J:=\{i:0<\mu_i\le2/N\}$ and
$\mu_J:=(\mu_i)_{i\in J}$. By Lemma~\ref{lem:t-light-rows},
$|J|\ge N/2$ and $\|\mu_J\|_2^2\le2/N$. Let $\mathcal B\subseteq\R^J$ be a
basis of the span of the columns $(\widetilde W_{i,y})_{i\in J}$,
$y\in\mathcal N$, so that $|\mathcal B|\le d_{\rm sk}$, and let $H$ be the
subspace of $z\in\R^J$ with
\[
 \sum_{i\in J}\mu_iz_i=0\quad\text{and}\quad
 \sum_{i\in J}\mu_iz_i\xi_i=0\ \text{ for all }\xi\in\mathcal B .
\]
Its codimension is at most $1+d_{\rm sk}$. For $z\in\R^J$ define
\[
 \Psi(z):=\sup_{x\in S^{r-1}}|\sum_{i\in J}\mu_iz_i
 |\langle u_i,x\rangle|^p| ,
\]
a finite seminorm on $\R^J$, hence continuous. We claim that a standard
Gaussian vector $G$ in $H$ satisfies
\begin{equation}\label{lem:x-a2-eq-round}
 \E[\Psi(G)]\le C_{\rm rd}\sqrt r\,q^{-p}N^{-1/2},
\end{equation}
where $C_{\rm rd}$ depends only on $p$. For $x\in S^{r-1}$ we split
\[
 \sum_{i\in J}\mu_iG_i|\langle u_i,x\rangle|^p
 =\mathcal Q(x)+\sum_{i\in J}\mu_iG_iR(\langle u_i,x\rangle),\qquad
 \mathcal Q(x):=\sum_{i\in J}\mu_iG_iP(\langle u_i,x\rangle).
\]

The random function $\mathcal Q$ is a polynomial on $\R^r$ of degree at
most $q$. Fix $y\in\mathcal N$. The column $(\widetilde W_{i,y})_{i\in J}$
lies in the span of $\mathcal B$, so
$\sum_{i\in J}\mu_iG_i\widetilde W_{i,y}=0$ because $G\in H$. Hence
$\mathcal Q(y)=\langle G,\omega^y\rangle$ with
$\omega^y_i:=\mu_i(P(\langle u_i,y\rangle)-M_P\widetilde W_{i,y})$, a
centered Gaussian variable of variance at most
$\|\omega^y\|_2^2\le M_P^2\tau_{\rm sk}^2\|\mu_J\|_2^2\le2M_P^2
\tau_{\rm sk}^2/N$. By Lemma~\ref{lem:t-norming-set}(a), applied to each
realization, and Lemma~\ref{lem:t-sidak-chaining}(a),
\[
 \E[\sup_{x\in S^{r-1}}|\mathcal Q(x)|]
 \le2\,\E[\max_{y\in\mathcal N}|\mathcal Q(y)|]
 \le2M_P\tau_{\rm sk}\sqrt{2/N}\sqrt{2\log(2|\mathcal N|)}
 \le10M_P\tau_{\rm sk}\sqrt{r\log(2r)/N}.
\]
Since $\log(2r)\le2r$ and $q\le r$, we have
$\tau_{\rm sk}\sqrt{\log(2r)}\le(2r)^{-p-1/2}\le q^{-p}$. Thus
$\E[\sup_x|\mathcal Q(x)|]\le10M_P\sqrt r\,q^{-p}N^{-1/2}$.

The residual $R$ is absolutely continuous on $[-1,1]$. By
Lemma~\ref{lem:t-threshold-form}, the signed measure $\beta:=R'(s)\,\d s$
on $(-1,1]$ has $\|\beta\|=\operatorname{TV}_{[-1,1]}(R)\le
C_{\rm tv}q^{-p}$, and $R(z)-R(-1)=\int\mathbf1_{\{z\ge s\}}\,\d\beta(s)$
for $z\in[-1,1]$. The vector $\mu_J$ is orthogonal to $H$, so
$\Pi_H(\mu_iR(-1))_{i\in J}=0$, and therefore, for $x\in S^{r-1}$,
\[
 \sum_{i\in J}\mu_iG_iR(\langle u_i,x\rangle)=\langle G,v_x\rangle,\qquad
 v_x:=\Pi_H(\mu_i\int\mathbf1_{\{\langle u_i,x\rangle\ge s\}}
 \,\d\beta(s))_{i\in J}.
\]
Apply Lemma~\ref{lem:t-integrated-vectors} with $X:=S^{r-1}$, $d:=r-1$,
$\ell:=1$, the index set $J$, $\mathcal Z:=(-1,1]$, $\tau_e:=\mu_e$,
$\beta_e:=\beta$, $V:=C_{\rm tv}q^{-p}$,
$\phi(e,s,x):=\mathbf1_{\{\langle u_e,x\rangle\ge s\}}$, which is the
indicator of one affine test on $\R^r$, and $\Pi:=\Pi_H$. Its quantity
$V\|\tau\|_2$ is at most $\varrho_N:=\sqrt2\,C_{\rm tv}q^{-p}N^{-1/2}$.
Parts~(a) and~(b) give $\|v_x\|_2\le\varrho_N$ for all $x$ and
$N(\mathcal V,u)\le\max\{1,(C_1\varrho_N/u)^{2(r-1)}\}$ for
$\mathcal V:=\{v_x:x\in S^{r-1}\}\subseteq H$ and all $u>0$, where $C_1$ is
the constant $C_\ell$ for $\ell=1$. This entropy bound rests on
Lemma~\ref{lem:t-one-inclusion}. Lemma~\ref{lem:t-sidak-chaining}(b) with
$R:=\varrho_N$, $B_0:=1$, $a_0:=\max\{1,C_1\}$ and $d:=r-1$ gives
\[
 \E[\sup_{x\in S^{r-1}}|\langle G,v_x\rangle|]
 \le14\varrho_N(\sqrt{\log2}+\sqrt{(r-1)\log(ea_0)})
 \le20C_{\rm tv}(1+\sqrt{\log(ea_0)})\sqrt r\,q^{-p}N^{-1/2}.
\]
Adding the two bounds proves Eq.~\eqref{lem:x-a2-eq-round} with
$C_{\rm rd}:=10M_P+20C_{\rm tv}(1+\sqrt{\log(ea_0)})$.

\emph{Step 3: iteration.} Let $\delta_0,\delta_1,n_{\rm c}$ be the
constants of Lemma~\ref{lem:t-subspace-coloring}, and define
\[
 a_*:=2C_{\rm rd}\,r^{(p+1)/2}q^{-p},\qquad
 C_*:=(\frac{2C_{\rm rd}}{1-(1-\delta_1/10)^{1/2}})^2,\qquad
 N_*:=\lceil\frac{C_*r^{p+1}q^{-2p}}{\epsilon^2}\rceil .
\]
Step~4 shows that
\begin{equation}\label{lem:x-a2-eq-budget}
 N_*\ge\max\{5n_{\rm c},\ 20,\ 2/\delta_0,\ 2(1+d_{\rm sk})/\delta_0
 \}.
\end{equation}
We apply Lemma~\ref{lem:t-iterated-coloring} with the index set $I$,
$\mathcal X:=S^{r-1}$, $\Phi_x(\mu):=\sum_i\mu_i|\langle u_i,x\rangle|^p$,
$\omega:=F$, $\mathcal S:=\{\mu\in[0,\infty)^I:\sum_i\mu_i=1\}$,
$b_*:=1/2$, and $a_*,N_*$ as above. Let $\mu\in\mathcal S$ have $N>N_*$
positive entries. The units are the singletons of $J$, so there are
$|J|\ge N/2$ of them. Take $H$ from Step~2 and
$\mathcal K:=\{z\in H:\Psi(z)\le2\E[\Psi(G)]\}$, which is closed, symmetric
and convex. By Eq.~\eqref{lem:x-a2-eq-budget} and $|J|>N_*/2$, the
codimension of $H$ is at most $\delta_0N_*/2<\delta_0|J|$, and
$|J|\ge1/\delta_0$. Lemma~\ref{lem:t-light-rows} then gives
$\gamma_H(\mathcal K)\ge1/2\ge e^{-\delta_0|J|}$. For
$z\in\mathcal K\cap[-1,1]^J$ the masses $\mu^z_i=\mu_i(1+z_i)$ are
nonnegative and have total mass one, because $\sum_{i\in J}\mu_iz_i=0$.
Thus $\mu^z\in\mathcal S$. Since $F(x)\ge r^{-p/2}$,
Eq.~\eqref{lem:x-a2-eq-round} gives, for $x\in S^{r-1}$,
\[
 |\Phi_x(\mu^z)-\Phi_x(\mu)|=|\sum_{i\in J}\mu_iz_i|\langle u_i,x
 \rangle|^p|\le\Psi(z)\le2C_{\rm rd}\sqrt r\,q^{-p}N^{-1/2}
 \le a_*N^{-1/2}F(x).
\]
Lemma~\ref{lem:t-iterated-coloring} with $\mu^0:=\nu$, for which
$\Phi_x(\nu)=F(x)$, gives $\mu^*\in\mathcal S$ with
$|\supp\mu^*|\le N_*$ and
\[
 |\Phi_x(\mu^*)-F(x)|\le\frac{a_*N_*^{-1/2}}{1-(1-\delta_1/10)^{1/2}}F(x)
 \le\epsilon F(x)
\]
by the choice of $C_*$. This is Eq.~\eqref{lem:x-a2-eq-goal}.

\emph{Step 4: the parameters.} Choose $c_2:=1/(2C_{\rm sk})$ and let
$C_p\ge2/c_2$. Then $c_2L\ge c_2C_p\log 4\ge2$, so $q\ge1$ and
$q\ge c_2L/2$, while $q\le L\le r$. By Step~1 and $d_{\rm sk}\ge1$,
\[
 \log(2(1+d_{\rm sk})/\delta_0)\le\log(2/\delta_0)+\log(2d_{\rm sk})
 \le L/2+C_{\rm sk}(p+5)\log(2r)+\log(2/\delta_0).
\]
On the other hand $\epsilon^{-2}=e^{2L}/4$, $r^{p+1}q^{-2p}\ge r^{1-p}$
and $C_*\ge4$, so $\log N_*\ge2L-(p-1)\log r$. For
$L\ge C_p\log(2r)$ with $C_p$ large in terms of $p$, $\delta_0$ and
$n_{\rm c}$, this gives Eq.~\eqref{lem:x-a2-eq-budget}. Finally, for such
$C_p$ we have $\epsilon^{-2}r^{p+1}L^{-2p}\ge e^{2L}r^{1-p}/4\ge1$, and
therefore, with $C_p':=1+C_*(2/c_2)^{2p}$,
\[
 N_*\le1+C_*\epsilon^{-2}r^{p+1}(c_2L/2)^{-2p}
 \le C_p'\,\epsilon^{-2}r^{p+1}L^{-2p}.
\]

\emph{Part (b).}
Note that $G_*\ge1$, since $L\ge\log4>1$ and $r/L\ge C_p\ge1$. Define
$\epsilon_0:=3\epsilon$ and $L_0:=\log(2/\epsilon_0)=L-\log3$. Let $c'_S$
and $d_0$ be the constants of Lemma~\ref{lem:s-shift} for the exponent
$p$ and $\lambda_1:=1/2$, $\lambda_2:=1$, $\beta_0:=1$, and let $c_{\rm t}$
and $d_p$ be the constants $c_p$ and $d_p$ of
Lemma~\ref{lem:lp_lower_tensor_product}. Define
$L_1:=\max\{d_0,\log12,1+2\log3\}$.

\emph{Bounded accuracy.} Let $L<L_1$. Lemma~\ref{lem:lp_noneven_lower}
with $s:=\lceil p/2\rceil\ge2$ gives
\[
 S_p(r,\epsilon)\ge c_{p,s}\min\{r^{p/2}/\epsilon,r^s\}\ge c_{p,s}r^{p/2}.
\]
On the other side
$\epsilon^{-2}=e^{2L}/4<e^{2L_1}$, $L^{-(3p-1)/2}\le1$,
$G_*^{2-p/2}\le L_1$ (for $2<p<4$ because $G_*\le L^{p/2}$ and
$(p/2)(2-p/2)\le1$), and $1+\log(L^{p/2}/G_*)\le1+(p/2)\log L_1$. Hence the
right side of the theorem is at most a constant depending on $p$ times
$r^{p/2}$, and the claim follows after decreasing the implied constant.

\emph{The base.} Let $L\ge L_1$. Define
\[
 d:=\lceil L\rceil,\qquad c_G:=\min\{1,1/(4d_p)\},\qquad
 G:=\max\{1,\lfloor c_GG_*\rfloor\},\qquad N:=Gd,\qquad b:=G^{-1/p}.
\]
We apply Lemma~\ref{lem:s-shift} at accuracy $\epsilon_0$ with $N$ tags,
$H:=G$ transverse copies and transverse dimension $d$. Its hypotheses
hold. First, $d\ge L\ge d_0$. Second, $L_0\le L\le d$, and
$L_0=L-\log3\ge(L+1)/2\ge d/2$ because $L\ge1+2\log3$, so
$\lambda_1d\le L_0\le\lambda_2d$. Third, $0<b\le1$ and $Hb^p=Gb^p=1$.
Fourth, $G\le G_*\le L^{p/2}$, since $c_G\le1$ and $G_*\ge1$ (for $p>4$,
$G=1$), and therefore
\[
 b\sqrt d=(\frac{d^{p/2}}G)^{1/p}\ge
 (\frac{L^{p/2}}{G_*})^{1/p}\ge1=\beta_0 .
\]
The conclusions of Lemma~\ref{lem:s-shift} include
$d/2\le\log\Xi\le3d$ for its dictionary size $\Xi$, and that
Lemma~\ref{lem:s-frame} applies to its width and shift window, so
$2C_0w\le b$ and in particular $w\le b$. Its shift count
$J=\lfloor b/(2C_0w)\rfloor$, with $w=c_0/\sqrt d$ for the constant $c_0$
of Lemma~\ref{lem:s-shift}, satisfies
$J\le b\sqrt d/(2C_0c_0)\le\sqrt d/(2C_0c_0)$. Take its rank to be
$k:=N+Gd=2Gd$, so that the last summand of the input of
Lemma~\ref{lem:s-cost} is absent. This gives a matrix $B$ of rank $k$
such that every reweighting of its rows by real weights, in particular
every nonnegative one, with relative error at most $3\epsilon$ gives
nonzero weight to more than
\[
 T:=c'_S\epsilon_0^{-2}NGd^{-p-1/2}(1+\log_+(b\sqrt d))
\]
rows. By the display above,
$1+\log_+(b\sqrt d)\ge1+\frac1p\log(L^{p/2}/G_*)\ge\frac1p(1+
\log(L^{p/2}/G_*))$. Conversely $d\le2L$ and $G\ge c_GG_*/2$ (shown
below) give $b\sqrt d\le(2^{1+p/2}L^{p/2}/(c_GG_*))^{1/p}$, so this factor,
and $\log J$, are of order $1+\log(L^{p/2}/G_*)$ and not larger.

\emph{Rank.} If $2<p<4$ and $\lfloor c_GG_*\rfloor\ge1$, then
$G\le c_Gr/L$ and $k=2Gd\le4c_Gr\le r/d_p$. Otherwise $G=1$ and
$k=2d\le4L\le r/d_p$ once $C_p\ge4d_p$. Also $\epsilon\le1/6$ since
$L\ge\log12$.

\emph{Amplification.} Lemma~\ref{lem:lp_lower_tensor_product}, applied to
$B$, its rank $k$ and the row count $T$, gives a rank-$r$ matrix that
requires at least $c_{\rm t}(r/k)^{p/2}T$ rows at accuracy $\epsilon$.
With $N=Gd$ and $k=2Gd$,
\[
 c_{\rm t}(\frac r{2Gd})^{p/2}T
 \ge\frac{c_{\rm t}c'_S}{9p\,2^{p/2}}\,\epsilon^{-2}r^{p/2}G^{2-p/2}
 d^{-(3p-1)/2}(1+\log\frac{L^{p/2}}{G_*}).
\]
Here $d^{-(3p-1)/2}\ge(2L)^{-(3p-1)/2}$. If $2<p<4$, then
$G\ge c_GG_*/2$: either $c_GG_*\ge1$ and $\lfloor c_GG_*\rfloor\ge
c_GG_*/2$, or $G=1>c_GG_*$. Hence $G^{2-p/2}\ge(c_G/2)^{2-p/2}
G_*^{2-p/2}$. If $p>4$, then $G=G_*=1$. This proves the bound.

\emph{Part (c).}
Let $c_{\rm ov}$ and $r_{\rm ov}$ be the constants of
Lemma~\ref{lem:x-a2-overlap}, and $c_{\rm t}$ and $d_p$ the constants $c_p$
and $d_p$ of Lemma~\ref{lem:x-a2-signed-tensor}. Define
$C_p:=\log3+5+\sqrt{r_{\rm ov}+1}$ and let $r\ge\max\{r_{\rm ov},
d_p(r_{\rm ov}+26)\}$. All lower bounds below count the rows of nonzero
weight in reweightings by arbitrary real weights with relative error at
most $\epsilon$. Sparsifiers are such reweightings, so these bounds also
bound $S_p(r,\epsilon)$. We distinguish three ranges.

\emph{The range $\sqrt r\le L\le r/4$.} Here $\min\{1,r/L^2\}=r/L^2$ and
$\max\{1,L/\sqrt r\}=L/\sqrt r$, and the right side of the theorem is
$\epsilon^{-2}r^{1+p/2}L^{-(3p+1)/2}(1+\log(L/\sqrt r))$ up to the implied
constant. Lemma~\ref{lem:x-a2-overlap} gives this bound.

\emph{The range $C_p\le L\le\sqrt{r/d_p}$.} Define $\epsilon_0:=3\epsilon$,
$L_0:=\log(2/\epsilon_0)=L-\log3$ and $k:=\lfloor L_0^2\rfloor$. Then
$\epsilon\le1/6$ since $L\ge\log12$, and $L_0\ge5+\sqrt{r_{\rm ov}+1}$.
Hence $k\ge r_{\rm ov}$, $\sqrt k\le L_0$, and $k\ge L_0^2-1\ge4L_0$, that
is, $L_0\le k/4$. Lemma~\ref{lem:x-a2-overlap} at rank $k$ and accuracy
$\epsilon_0$ gives a rank-$k$ matrix $B$ for which every reweighting of its
rows by real weights with relative error at most $3\epsilon$ gives nonzero
weight to more than
\[
 T:=c_{\rm ov}(3\epsilon)^{-2}k^{1+p/2}L_0^{-(3p+1)/2}
\]
rows, where we dropped the factor $1+\log(L_0/\sqrt k)\ge1$. Since
$k\le L^2\le r/d_p$, Lemma~\ref{lem:x-a2-signed-tensor} gives a rank-$r$
matrix that requires at least
\[
 c_{\rm t}(\frac rk)^{p/2}T=\frac{c_{\rm t}c_{\rm ov}}9\,
 \epsilon^{-2}r^{p/2}kL_0^{-(3p+1)/2}
 \ge\frac{c_{\rm t}c_{\rm ov}}{18}\,\epsilon^{-2}r^{p/2}L^{-3(p-1)/2}
\]
rows, using $k\ge L_0^2/2$ and $L_0\le L$. In this range
$\min\{1,r/L^2\}=1$ and $\max\{1,L/\sqrt r\}=1$, so this is the claim.

\emph{The range $\sqrt{r/d_p}<L<\sqrt r$.} Define $k_1:=\lfloor L^2\rfloor$.
Then $L^2>r/d_p\ge r_{\rm ov}+26$, so $k_1\ge r_{\rm ov}$, $L\ge5$ and
$k_1\ge L^2-1\ge4L$. Also $\sqrt{k_1}\le L$ and $k_1<r$.
Lemma~\ref{lem:x-a2-overlap} at rank $k_1$ and accuracy $\epsilon$ gives a
rank-$k_1$ matrix $B_1$. The matrix $\operatorname{diag}(B_1,I_{r-k_1})$
has rank $r$. On queries supported in the first $k_1$ coordinates the
identity rows vanish, so every reweighting of its rows by real weights with
relative error at most $\epsilon$ restricts to such a reweighting of the
rows of $B_1$. Hence more than
\[
 c_{\rm ov}\epsilon^{-2}k_1^{1+p/2}L^{-(3p+1)/2}
 \ge c_{\rm ov}(2d_p)^{-1-p/2}\epsilon^{-2}r^{1+p/2}L^{-(3p+1)/2}
 \ge c_{\rm ov}(2d_p)^{-1-p/2}\epsilon^{-2}r^{p/2}L^{-3(p-1)/2}
\]
rows have nonzero weight, using $k_1\ge L^2/2>r/(2d_p)$ and $r/L^2>1$. In
this range $\min\{1,r/L^2\}=1$ and $\max\{1,L/\sqrt r\}=1$.
\end{proof}

\subsection{Rank Powers at Subexponential Accuracy}

Combining Theorem~\ref{thm:above_two_main}(a) with
Theorem~\ref{thm:lp_above_two}, at
$L=r^\alpha$ with fixed $0<\alpha<1$ the upper rank power
of $\epsilon^2 S_p$ is at most
$\min\{p/2,\,p+1-2p\alpha\}$. In particular the polynomial
approximation improves this power for $\alpha>1/4+1/(2p)$.

For example, when $p=3$ and $r\ge L^{5/2}$, the denominator
in the large-rank lower bound is $L^{13/4}$. When $p>4$, one
transverse group gives a factor of order $1+\log L$ over the
single-shift lower bound. The gap at moderate accuracy remains open.
At $L=r^\alpha$, $0<\alpha<1$ fixed, the lower rank power of
$\epsilon^2S_p$ supplied by Theorem~\ref{thm:above_two_main}(b) is
\[
 \frac p2-\frac{(3p-1)\alpha}{2}
     +(2-p/2)_+\min\{p\alpha/2,\,1-\alpha\}.
\]

By Theorem~\ref{thm:above_two_main}(c), at $L=r^\alpha$ for fixed
$0<\alpha<1$ the lower rank power
of $\epsilon^2S_p$ is at least
\[
 \frac p2-\frac{3(p-1)}2\alpha-(2\alpha-1)_+.
\]
This exceeds the power from Theorem~\ref{thm:above_two_main}(b)
by $\min\{\alpha,1-\alpha\}-(2-p/2)_+\min\{p\alpha/2,1-\alpha\}>0$.

For $p=3$ and $L\asymp r^{2/3}$ this gives
$\Omega(\epsilon^{-2}r^{-5/6}\log r)$, improving the bound of
Theorem~\ref{thm:above_two_main}(b) there by a factor of order
$r^{1/6}$. Theorem~\ref{thm:above_two_main}(a) gives an
upper bound of $O(\epsilon^{-2})$ in the same regime.

\subsection{Overlapping Transverse Groups}

For $\sqrt r\le L\le r/4$ the construction of
Theorem~\ref{thm:above_two_main}(c) uses $N\asymp r$ shared tag
coordinates and $G\asymp_p(r/L)^{p/2}$ overlapping transverse groups. The
transverse part of a row is $v_g\otimes u_t$, where the outer directions
$v_g\in\R^m$ have small inner products $\rho_{hg}$ with bounded $p$th-power
row sums, and the $u_t\in\R^d$, $d:=\lceil L\rceil$, form the dictionary of
Lemma~\ref{lem:s-frame}. A test of group $h$, built from shifted queries as
in Lemma~\ref{lem:s-cost}, also sees every other group $g$, through the
scalar factor $\rho_{hg}$. These responses lie within $e^{-12d}$ of a space
of evaluations of low-degree polynomials, because the test functions are
band limited and $\rho_{hg}$ is small. After projecting the Gaussian input
perturbations off that space, the clipped tests of
Lemma~\ref{lem:s-width} keep their width, and Lemma~\ref{lem:s-pack}
forces order $NG\Xi$ rows. Smaller $L$ is reduced to this range by
Lemma~\ref{lem:lp_lower_tensor_product} or by padding. We use the scalar
test $\chi$ with its density $\varrho$, the function $\psi$ of
Lemma~\ref{lem:s-filter}, the constant $C_0$ of Lemma~\ref{lem:s-scalar},
and the notation $h_j$, $f_j:=w^ph_j$, $\sigma:=w^{p+1/2}$,
$\mathbf f_j(X)$, $\mathbb S_d$, $\mathcal B(\mathcal Q)$ and admissible
coefficients Eq.~\eqref{lem:s-eq-admissible} of the subsection on the
adaptive-shift construction in Appendix~\ref{app:common_tools}.

\begin{lemma}[Nearly orthogonal outer directions]\label{lem:x-a2-outer}
Fix $p>2$ and $0<\rho_*\le1/2$. There are $c_\rho\in(0,1]$ and
$C_\rho\ge1$, depending only on $p$ and $\rho_*$, such that for every
integer $m\ge1$ there are unit vectors $v_1,\ldots,v_G$ spanning $\R^m$,
with $c_\rho m^{p/2}\le G\le2m^{p/2}$, whose inner products
$\rho_{hg}:=\langle v_h,v_g\rangle$ satisfy $|\rho_{hg}|\le\rho_*$ for
$h\ne g$ and $\sum_{g\in[G]}|\rho_{hg}|^p\le C_\rho$ for every $h$.
\end{lemma}

\begin{proof}
Let $m_0\ge4$ be such that $32m^{p+1/2}e^{-\rho_*^2(m-3)/2}\le1/4$ for
$m\ge m_0$. If $m\le m_0$, take an orthonormal basis. Then $G=m$, and
$m\ge m_0^{1-p/2}m^{p/2}$ because $p>2$, so the claim holds with
$c_\rho:=m_0^{1-p/2}$ and any $C_\rho\ge1$.

Let $m>m_0$, define $M:=2\lceil m^{p/2}\rceil\le4m^{p/2}$, and let
$v_1,\ldots,v_M$ be independent uniform points of $S^{m-1}$. For $h\ne g$,
conditioning on $v_h$ and rotating shows that $\rho_{hg}$ has the law of a
coordinate of a uniform unit vector, with density
$c_m(1-t^2)^{(m-3)/2}$ on $[-1,1]$ and $c_m\le\sqrt m$ by
Lemma~\ref{lem:s-gauss}(d). Using $1-t^2\le e^{-t^2}$,
\[
 \Pr[|\rho_{hg}|>\rho_*]\le2\sqrt m\,e^{-\rho_*^2(m-3)/2},\qquad
 \E[|\rho_{hg}|^p]\le\sqrt m\int_\R|t|^pe^{-(m-3)t^2/2}\,\d t
 \le\gamma_pm^{-p/2},
\]
with $\gamma_p:=8^{(p+1)/2}\Gamma((p+1)/2)$, since
$\int_\R|t|^pe^{-at^2}\,\d t=\Gamma((p+1)/2)a^{-(p+1)/2}$ and
$m/(m-3)\le4$. A union bound over at most $M^2\le16m^p$ pairs shows that
$\max_{h\ne g}|\rho_{hg}|\le\rho_*$ outside probability $1/4$. The sum
$S:=\sum_{h\ne g}|\rho_{hg}|^p$ has $\E[S]\le M^2\gamma_pm^{-p/2}\le
4\gamma_pM$, so $S\le16\gamma_pM$ outside probability $1/4$ by Markov's
inequality (Lemma~\ref{lem:pre_markov}). Almost surely any $m$ of the points are linearly independent,
since a proper subspace meets the sphere in a null set. Fix a realization
with all three properties. Let $R_h:=\sum_{g\ne h}|\rho_{hg}|^p$ and keep
the $G:=M/2$ indices with the smallest $R_h$. Each kept index has
$R_h\le2S/M\le32\gamma_p$, for otherwise all $M/2$ discarded indices would
also have $R_h>2S/M$ and $\sum_hR_h>S$. Restricting the family can only
decrease the row sums, so $C_\rho:=1+32\gamma_p$ works. Finally
$m\le m^{p/2}\le G\le2m^{p/2}$, so the kept vectors span $\R^m$.
\end{proof}

\begin{lemma}[Tests on overlapping groups]\label{lem:x-a2-cost}
Fix a non-even $p>0$ and let $C_p$ be the constant of
Lemma~\ref{lem:s-cost}. Let $N,G,m,d,\Xi\ge1$ be integers, let
$u_1,\ldots,u_\Xi\in\R^d$, and let $\mathcal Q$ be a probability law on
$\R^d$ with $\Xi^{-1}\sum_t|\langle u_t,X\rangle|^p\le C_E$ for
$\mathcal Q$-a.e.\ $X$. Let $v_1,\ldots,v_G\in\R^m$ be unit vectors with
$\sum_g|\rho_{hg}|^p\le C_V$ for every $h$, where
$\rho_{hg}:=\langle v_h,v_g\rangle$. Let $0<w\le b\le1$ with $Gb^p\le1$,
and let $s_1,\ldots,s_J\in[0,b]$. For admissible coefficients $c$ and
$h\in[G]$ define $\mathbf B^{(h)}(c)\in\R^{N\times G\times\Xi}$ by
\[
 \mathbf B^{(h)}(c)_{igt}:=\E_X[\sum_{j\in[J]}c_{ij}(X)
 f_j(\rho_{hg}\langle u_t,X\rangle)],\qquad X\sim\mathcal Q .
\]
In $\R^N\oplus(\R^m\otimes\R^d)$ consider the rows
$\alpha_{igt}^{1/p}(e_i,v_g\otimes u_t)$, $(i,g,t)\in[N]\times[G]\times
[\Xi]$, with masses $\alpha_{igt}\in[1/2,3/2]$. Suppose that real weights
on these rows preserve every energy within relative error $\epsilon$, and
let $\widehat\alpha_{igt}$ be the output mass of row $(i,g,t)$, that is,
its weight times $\alpha_{igt}$. Then for every $h\in[G]$ and all
admissible $c$,
\begin{equation}\label{lem:x-a2-eq-cost}
 |\langle\widehat\alpha-\alpha,\mathbf B^{(h)}(c)\rangle|
 \le C_p(1+C_VC_E)\epsilon N\Xi .
\end{equation}
\end{lemma}

\begin{proof}
Let $\mathcal E(x)$ be the output energy minus the input energy at the
query $x$, so $|\mathcal E(x)|\le\epsilon F(x)$. Fix $h$ and $c$, and define
$z:=\widehat\alpha-\alpha$. Fix $X$ with
$\Xi^{-1}\sum_t|\langle u_t,X\rangle|^p\le C_E$ and
$C_X:=\max_i\sum_j|c_{ij}(X)|>0$, since queries with $C_X=0$ contribute
nothing, and write $Z_t:=\langle u_t,X\rangle$. As in the proof of
Lemma~\ref{lem:s-cost}, draw independently at each tag $i$ an index
$J_i\in[J]\cup\{\emptyset\}$ with $\Pr[J_i=j]=(c_{ij}(X))_+/C_X$ and
$\Pr[J_i=\emptyset]=1-\sum_j(c_{ij}(X))_+/C_X$. For $\upsilon\in\R$ let
$y(\upsilon)\in\R^N$ have coordinates $y_i(\upsilon):=-s_{J_i}-w\upsilon$
if $J_i\ne\emptyset$ and $0$ otherwise, and let
$x_\pm(\upsilon):=(y(\upsilon),\pm v_h\otimes X)$. The row $(i,g,t)$ has
inner product $\alpha_{igt}^{1/p}(y_i(\upsilon)\pm\rho_{hg}Z_t)$ with
$x_\pm(\upsilon)$, because $\langle v_g\otimes u_t,v_h\otimes X\rangle=
\rho_{hg}Z_t$.

A tag with $J_i=\emptyset$ contributes
$z_{igt}(|\rho_{hg}Z_t|^p-|-\rho_{hg}Z_t|^p)=0$ to
$\mathcal E(x_+)-\mathcal E(x_-)$. For a tag with $J_i=j$,
Eq.~\eqref{lem:s-eq-filter} and the evenness of $\chi$ give, for all
$\rho,Z\in\R$,
\[
 \frac12\int(|\rho Z-s_j-w\upsilon|^p-|\rho Z+s_j+w\upsilon|^p)
 \psi(\upsilon)\,\d\upsilon=f_j(\rho Z).
\]
Hence
$\frac12\int(\mathcal E(x_+(\upsilon))-\mathcal E(x_-(\upsilon)))
\psi(\upsilon)\,\d\upsilon=\sum_{i:J_i\ne\emptyset}\sum_{g,t}z_{igt}
f_{J_i}(\rho_{hg}Z_t)$, whose expectation over the $J_i$ is
$\sum_{i,g,t}z_{igt}\sum_j((c_{ij}(X))_+/C_X)f_j(\rho_{hg}Z_t)$.

For the energy define $c_p':=\max\{1,2^{p-1}\}$. Every tag coordinate
satisfies $|y_i(\upsilon)|^p\le c_p'(b^p+w^p|\upsilon|^p)$, and
$|y_i\pm\rho_{hg}Z_t|^p\le c_p'(|y_i|^p+|\rho_{hg}|^p|Z_t|^p)$. Since the
masses are at most $3/2$,
\[
 F(x_\pm(\upsilon))\le\frac32c_p'[G\Xi\sum_i|y_i(\upsilon)|^p
 +N\sum_g|\rho_{hg}|^p\sum_t|Z_t|^p]
 \le\frac32c_p'^2N\Xi(1+C_VC_E)(1+|\upsilon|^p),
\]
using $Gb^p\le1$ and $Gw^p\le Gb^p$. Therefore the absolute value of the
last expectation is at most
$\epsilon\frac32c_p'^2N\Xi(1+C_VC_E)\int(1+|\upsilon|^p)|\psi(\upsilon)|
\,\d\upsilon$. The negative parts $(c_{ij}(X))_-$ are treated in the same
way. Multiplying by $C_X$, integrating over $X$ and
using Eq.~\eqref{lem:s-eq-admissible} gives Eq.~\eqref{lem:x-a2-eq-cost} with
$C_p:=3c_p'^2\int(1+|\upsilon|^p)|\psi(\upsilon)|\,\d\upsilon$, which is the
constant of Lemma~\ref{lem:s-cost}.
\end{proof}

\begin{lemma}[Responses off the active group]\label{lem:x-a2-offdiag}
Let $d\ge32$ be an integer, $q_0:=\lfloor d/16\rfloor$, $A_0,c_0>0$,
$w:=c_0/\sqrt d$, and $0\le\rho_*\le c_0e^{-192}/(32eA_0)$. Let
$u_1,\ldots,u_\Xi\in\R^d$ have $\|u_t\|_2\le A_0$, let $E_0\subseteq\R^\Xi$
be the space of vectors $(\Pi(u_t))_{t\in[\Xi]}$ over real polynomials $\Pi$
on $\R^d$ of degree at most $q_0$, and let $P_0$ be the orthogonal projection
onto $E_0^\perp$. Then $\dim E_0\le e^{d/3}$. Moreover, for all $s\in\R$,
all $|\rho|\le\rho_*$ and all $X\in\R^d$ with $\|X\|_2\le\sqrt d$, the
vector
\[
 \mathbf e:=(\frac{w^p}2[\chi(\frac{\rho\langle u_t,X\rangle-s}
 w)-\chi(\frac{\rho\langle u_t,X\rangle+s}w)])_{t\in[\Xi]}
\]
satisfies $\|P_0\mathbf e\|_2\le\sqrt\Xi\,w^pe^{-12d}$.
\end{lemma}

\begin{proof}
The space $E_0$ is spanned by the evaluations of the monomials of degree at
most $q_0$ in $d$ variables, so $\dim E_0\le\binom{d+q_0}{q_0}\le
(e(d+q_0)/q_0)^{q_0}$. For $d\ge32$ we have $q_0\ge d/16-1\ge d/32$, so
$(d+q_0)/q_0\le33$, and $q_0\le d/16$. Hence $\dim E_0\le(33e)^{d/16}\le
e^{d/3}$, since $1+\log33<16/3$.

For real $\theta$ the Taylor polynomial
$T(\theta):=\sum_{n=0}^{q_0}(\mathrm i\theta)^n/n!$ satisfies
$|e^{\mathrm i\theta}-T(\theta)|\le|\theta|^{q_0+1}/(q_0+1)!$. For
$y\in\R$ and $u\in\R^d$ we have
$\chi((\rho\langle u,X\rangle-y)/w)=\int\varrho(\xi)e^{-\mathrm i\xi y/w}
e^{\mathrm i\xi\rho\langle u,X\rangle/w}\,\d\xi$. Let $\Pi_y(u)$ be the real
part of the same integral with the last exponential replaced by
$T(\xi\rho\langle u,X\rangle/w)$. Expanding $\langle u,X\rangle^n$ shows
that $\Pi_y$ is a real polynomial on $\R^d$ of degree at most $q_0$. Since
$\chi$ is real, $|\chi((\rho\langle u,X\rangle-y)/w)-\Pi_y(u)|$ is at most
$\int\varrho(\xi)|e^{\mathrm i\theta}-T(\theta)|\,\d\xi$ with
$\theta:=\xi\rho\langle u,X\rangle/w$. On the support of $\varrho$ we have
$|\xi|\le2$, and $\int\varrho=1$. Hence, for $\|u\|_2\le A_0$, this error is
at most
\[
 \frac{(2\rho_*A_0\sqrt d/w)^{q_0+1}}{(q_0+1)!}
 \le(\frac{2e\rho_*A_0d}{c_0(q_0+1)})^{q_0+1}
 \le(\frac{32e\rho_*A_0}{c_0})^{q_0+1}
 \le e^{-192(q_0+1)}\le e^{-12d},
\]
using $(q_0+1)!\ge((q_0+1)/e)^{q_0+1}$ and $q_0+1\ge d/16$. Applying this
with $y=s$ and $y=-s$ writes $\mathbf e$ as a vector of $E_0$ plus a vector
whose entries are at most $w^pe^{-12d}$ in absolute value. The projection
$P_0$ annihilates the first and contracts the second.
\end{proof}

\begin{lemma}[Overlapping transverse groups]\label{lem:x-a2-overlap}
Fix a non-even $p>2$. There are $c_{\rm ov}>0$ and $r_{\rm ov}$, depending
only on $p$, such that for every integer $r\ge r_{\rm ov}$ and every
$0<\epsilon\le1/2$ with $\sqrt r\le L\le r/4$ some matrix of rank $r$ has
the following property. Every reweighting of its rows by real weights
with relative error at most $\epsilon$ gives nonzero weight to more than
\[
 c_{\rm ov}\epsilon^{-2}r^{1+p/2}L^{-(3p+1)/2}(1+\log(L/\sqrt r))
\]
rows.
\end{lemma}

\begin{proof}
\emph{Constants.} Let $c_{0,F}$, $c_F$ and $C_F$ be the constants of
Lemma~\ref{lem:s-frame} for $p$, $\lambda_1:=1$ and $\lambda_2:=3$, so
that (F1) holds with $A_0=4$. Let $J_0$, $A_1$ and $c_W$ be the constants of
Lemma~\ref{lem:s-width} for $c_F,C_F$. Define
\[
 c_0:=\min\{c_{0,F},\frac1{4C_0J_0}\},\qquad
 \rho_*:=\min\{\frac12,\frac{c_0e^{-192}}{32eA_0}\},
\]
and let $c_\rho,C_\rho$ be the constants of Lemma~\ref{lem:x-a2-outer} for
$p$ and $\rho_*$. Let $C_p$ be the constant of Lemma~\ref{lem:x-a2-cost},
define $C_{\rm cost}:=C_p(1+C_\rho C_F)$ and $\theta_0:=\min\{1,c_W/(8A_1)\}$,
let $c_{\rm pk}$ and $n_0$ be the constants of Lemma~\ref{lem:s-pack} for
$\theta_0$, and define $c_*:=\min\{1,(2c_{\rm pk}A_1/C_{\rm cost})^2\}$.
Finally let $d_1$ be large in terms of all of these, and
$r_{\rm ov}:=d_1^2$.

\emph{Parameters.} Define $d:=\lceil L\rceil$, $m:=\lfloor r/(2d)\rfloor$
and $N:=r-md$. Let $v_1,\ldots,v_G\in\R^m$ be given by
Lemma~\ref{lem:x-a2-outer}, and define
\begin{gather*}
 b:=G^{-1/p},\qquad w:=\frac{c_0}{\sqrt d},\qquad
 J:=\lfloor\frac b{2C_0w}\rfloor,\qquad
 s_j:=\frac b2+(j-1)C_0w\quad(j\in[J]),\\
 \tau:=(\frac{\log J}4)^{1/2},\qquad
 \Theta:=c_*\epsilon^{-2}\sigma^2\tau^2,\qquad \Xi:=\lfloor\Theta\rfloor,
 \qquad \zeta:=w^pe^{-12d}.
\end{gather*}
Since $L\ge\sqrt r\ge d_1$, we have $d\ge d_1$, and we record the facts
used below.
\begin{enumerate}[label=(\alph*),nosep,leftmargin=*]
\item $d\le L+1\le r/2$, so $m\ge1$ and $N\ge r/2$. Also
$m\le r/(2L)\le L/2\le d/2$, as $r\le L^2$, and $m\ge r/(4d)\ge r/(8L)$.
\item $0<b\le1$, $Gb^p=1$, and, since $G\le2m^{p/2}$ and $m\le r/(2d)$,
\[
 y:=b\sqrt d\ge2^{-1/p}(d/m)^{1/2}\ge2^{-1/p}(2d^2/r)^{1/2}
 \ge L/\sqrt r\ge1 .
\]
\item $b/(2C_0w)=y/(2C_0c_0)\ge2J_0y\ge2J_0$. Hence $J\ge J_0$,
$2C_0w\le b$, so $w\le b$, and $J\ge b/(4C_0w)\ge J_0y\ge y$. Also
$J\le\sqrt d/(2C_0c_0)\le d$. The centers $s_j$ lie in $[b/2,b]$ and are
pairwise at distance at least $C_0w$. As $J\ge55$, $1\le\tau^2\le\log d$.
\item $\Theta\le e^{2L}\log d$, so $\log\Xi\le2d+\log\log d\le3d$. Also
$\log\Theta\ge2L-\log4+\log(c_*c_0^{2p+1})-(p+\frac12)\log d\ge d+\log2$,
so $\Xi\ge\Theta/2$ and $d\le\log\Xi\le3d$. In particular
$\Xi\ge e^d\ge J^2$ and $\Xi\ge n_0$. Moreover $d\ge32$ and $d$ is at
least the number $d_F$ of Lemma~\ref{lem:s-frame} for $c_0'=c_0$.
\item $A_1e^{d/6}\le(c_W/2)e^{d/2}$ and
$2c_0^{-1/2}d^{p/2+1/4}e^{-21d/2}\le\min\{A_1,c_W/4\}$.
\end{enumerate}

\emph{The input.} By (c) and (d), Lemma~\ref{lem:s-frame} applies to $d$,
$\Xi$, $w=c_0/\sqrt d$, $b$ and the $s_j$, with the sure event as
$\mathcal E$. Let $u_1,\ldots,u_\Xi\in\R^d$ and $\mathcal Q$ on
$\mathbb S_d$ be given by it, so that (F1)--(F5) hold. The input has the
rows $\alpha_{igt}^{1/p}(e_i,v_g\otimes u_t)\in\R^N\oplus(\R^m\otimes\R^d)
=\R^r$, $(i,g,t)\in[N]\times[G]\times[\Xi]$, with masses
$\alpha\in[1/2,3/2]^{N\times G\times\Xi}$ chosen below, and $n:=NG\Xi$.
It has rank $r$. Indeed, if $\langle y,u_t-u_1\rangle=0$ for all $t$, then
$(-\langle y,u_1\rangle,y)$ is orthogonal to every $(1,u_t)$, so $y=0$ by
(F1). Hence the differences $u_t-u_1$ span $\R^d$, the vectors
$v_g\otimes(u_t-u_1)$ span $\R^m\otimes\R^d$, and together with the rows
$(e_i,v_g\otimes u_1)$ they span $\R^r$.

\emph{The test seminorm.} For an array $z\in\R^{N\times G\times\Xi}$ and
$h\in[G]$ let $D_hz\in\R^{N\times\Xi}$ be the block $(z_{iht})_{i,t}$ and
$O_hz$ the array $z$ with this block replaced by zero. Let
$\mathbf B^{(h)}(c)$ be the arrays of Lemma~\ref{lem:x-a2-cost} for these
data, $X\sim\mathcal Q$ and admissible $c$. Since $\rho_{hh}=1$, the $i$th
row of $D_h\mathbf B^{(h)}(c)$ is $\E_X[\sum_jc_{ij}(X)\mathbf f_j(X)]$. Define
$R:=A_1\sqrt N\sigma\tau$,
\[
 \mathcal B^{(h)}:=\{\mathbf B^{(h)}(c):c\text{ admissible},\
 \|D_h\mathbf B^{(h)}(c)\|_2\le R\},\qquad
 \|z\|_\sharp:=\sum_{h\in[G]}\sup_{\mathbf B\in\mathcal B^{(h)}}
 \langle z,\mathbf B\rangle .
\]
Since $|\chi|\le1$, every entry of $\mathbf B^{(h)}(c)$ has absolute value
at most $w^p$. Each $\mathcal B^{(h)}$ is therefore bounded, it contains
$0$, and it is symmetric under $c\mapsto-c$. Hence $\|\cdot\|_\sharp$ is a
seminorm on $\R^n$.

\emph{Cost.} By (b), (c), (F2) and Lemma~\ref{lem:x-a2-cost} with
$C_E:=C_F$ and $C_V:=C_\rho$, every reweighting of the rows by real
weights with relative error at most $\epsilon$ has an output mass array
$\widehat\alpha$ with
\[
 \|\widehat\alpha-\alpha\|_\sharp\le\epsilon_0:=GC_{\rm cost}\epsilon N\Xi .
\]

\emph{Removing the responses of the other groups.} Let $q_0$, $E_0$ and
$P_0$ be as in Lemma~\ref{lem:x-a2-offdiag} for $d$, $A_0$, $c_0$ and the
$u_t$, and let $\mathcal P$ act on $\R^{N\times G\times\Xi}$ by applying
$P_0$ to each vector $(z_{igt})_{t\in[\Xi]}$. Let $h\in[G]$ and
$\mathbf B=\mathbf B^{(h)}(c)\in\mathcal B^{(h)}$. For $g\ne h$ the vector
$((\mathcal P\mathbf B)_{igt})_t$ equals
$\E_X[\sum_jc_{ij}(X)P_0(f_j(\rho_{hg}\langle u_t,X\rangle))_t]$. We have
$\|u_t\|_2\le A_0$ by (F1), $\|X\|_2=\sqrt d$ on $\mathbb S_d$, which
carries $\mathcal Q$, and $|\rho_{hg}|\le\rho_*$, so
Lemma~\ref{lem:x-a2-offdiag} bounds its norm by
$\E_X[\sum_j|c_{ij}(X)|\sqrt\Xi\,\zeta]\le\sqrt\Xi\,\zeta$. Hence
$\|O_h\mathcal P\mathbf B\|_2\le\sqrt n\,\zeta$, while
$D_h\mathcal P\mathbf B$ is $D_h\mathbf B$ with $P_0$ applied to each row,
so that $\|D_h\mathcal P\mathbf B\|_2\le R$. By (a), (d) and (e),
\begin{equation}\label{lem:x-a2-eq-small}
 G\sqrt\Xi\,\zeta\le2d^{p/2}e^{3d/2}w^pe^{-12d}
 \le2c_0^{-1/2}d^{p/2+1/4}e^{-21d/2}\sigma
 \le\min\{A_1,c_W/4\}\,\sigma\tau ,
\end{equation}
using $G\le2m^{p/2}\le2d^{p/2}$, $\sigma=w^pc_0^{1/2}d^{-1/4}$ and
$\tau\ge1$.

\emph{Upper bound on the range of $\mathcal P$.} Let
$z\in\operatorname{ran}\mathcal P$ and $\mathbf B\in\mathcal B^{(h)}$.
Since $\mathcal P$ is an orthogonal projection,
\[
 \langle z,\mathbf B\rangle=\langle z,\mathcal P\mathbf B\rangle
 =\langle D_hz,D_h\mathcal P\mathbf B\rangle
  +\langle O_hz,O_h\mathcal P\mathbf B\rangle
 \le R\|D_hz\|_2+\sqrt n\,\zeta\|z\|_2 .
\]
Summing over $h$, with $\sum_h\|D_hz\|_2\le\sqrt G\|z\|_2$
and Eq.~\eqref{lem:x-a2-eq-small}, we get
$\|z\|_\sharp\le(\sqrt GR+G\sqrt n\,\zeta)\|z\|_2\le R'\|z\|_2$ for
$R':=2\sqrt GR$, because $G\sqrt n\,\zeta=\sqrt{NG}\,G\sqrt\Xi\,\zeta\le
\sqrt{NG}A_1\sigma\tau=\sqrt GR$.

\emph{Width.} Let $g'$ be a standard Gaussian array in
$\R^{N\times G\times\Xi}$ and $g:=\mathcal Pg'$. Let $\|\cdot\|_*$ be the
seminorm on $\R^{N\times\Xi}$ given by Lemma~\ref{lem:s-width} for
$\mathcal Q$ and the maps $\mathbf f_j$. Its hypotheses are (F3)--(F5),
$J\ge J_0$ and $\Xi\ge J^2$. We claim that
\[
 \sup_{\mathbf B\in\mathcal B^{(h)}}\langle D_hg,D_h\mathbf B\rangle
 \ge\|D_hg\|_* .
\]
Let $\mathbf B'\in\mathcal B(\mathcal Q)$ with $\|\mathbf B'\|_2\le R$. By
the definition of $\mathcal B(\mathcal Q)$ as a closure there are arrays
$\mathbf B'_\nu$ with rows $\E_X[\sum_jc^{(\nu)}_{ij}(X)\mathbf f_j(X)]$ for
admissible $c^{(\nu)}$ that converge to $\mathbf B'$. The set of such
arrays is convex and contains $0$, so
$\min\{1,R/\|\mathbf B'_\nu\|_2\}\mathbf B'_\nu$ is again of this form, has
norm at most $R$ and converges to $\mathbf B'$. Each of these arrays is
$D_h\mathbf B$ for some $\mathbf B\in\mathcal B^{(h)}$, namely
$\mathbf B=\mathbf B^{(h)}$ of the rescaled coefficients, and the claim
follows by taking limits.

Now let $\mathbf B\in\mathcal B^{(h)}$. As in the previous paragraph,
$\langle g,\mathbf B\rangle\ge\langle D_hg,D_h\mathcal P\mathbf B\rangle
-\sqrt n\,\zeta\|g\|_2$, and $\langle D_hg,D_h\mathcal P\mathbf B\rangle
=\langle D_hg,D_h\mathbf B\rangle$ because $D_hg$ is fixed by $P_0$ in each
row. The block $g'_h:=D_hg'$ is a standard Gaussian array, and
$D_hg=g'_h-g''_h$, where $g''_h$ is the projection of $g'_h$ onto the
arrays whose rows lie in $E_0$, a standard Gaussian vector in a space of
dimension $N\dim E_0\le Ne^{d/3}$. Since $\|\cdot\|_*\le R\|\cdot\|_2$,
Lemma~\ref{lem:s-width} and (e) give
\[
 \E[\|D_hg\|_*]\ge\E[\|g'_h\|_*]-R\,\E[\|g''_h\|_2]
 \ge c_WN\sqrt\Xi\,\sigma\tau-A_1N\sigma\tau e^{d/6}
 \ge\frac{c_W}2N\sqrt\Xi\,\sigma\tau ,
\]
as $\sqrt\Xi\ge e^{d/2}$. Summing the claim over $h$, and using
$\E[\|g\|_2]\le\sqrt n$ and Eq.~\eqref{lem:x-a2-eq-small},
\[
 \E[\|g\|_\sharp]\ge G(\frac{c_W}2N\sqrt\Xi\,\sigma\tau-n\zeta)
 =GN\sqrt\Xi(\frac{c_W}2\sigma\tau-G\sqrt\Xi\,\zeta)
 \ge\frac{c_W}4GN\sqrt\Xi\,\sigma\tau\ge\theta_0R'\sqrt n ,
\]
since $R'\sqrt n=2A_1GN\sqrt\Xi\,\sigma\tau$.

\emph{Packing.} The condition $\epsilon_0\le c_{\rm pk}R'\sqrt n$ reads
$C_{\rm cost}\epsilon\sqrt\Xi\le2c_{\rm pk}A_1\sigma\tau$. It holds because
$\Xi\le\Theta=c_*\epsilon^{-2}\sigma^2\tau^2$ and
$c_*\le(2c_{\rm pk}A_1/C_{\rm cost})^2$. Since $n\ge\Xi\ge n_0$,
Lemma~\ref{lem:s-pack}, applied to $\mathcal P$, $\|\cdot\|_\sharp$ and
$R'$, gives $\alpha\in[1/2,3/2]^{N\times G\times\Xi}$ such that every array
within $\epsilon_0$ of $\alpha$ in $\|\cdot\|_\sharp$ has more than
$c_{\rm pk}n$ nonzero entries. We use these masses. By the cost bound,
every reweighting of the rows by real weights with relative error at most
$\epsilon$ has more than $c_{\rm pk}n$ nonzero output masses, and a row
with nonzero output mass has nonzero weight.

\emph{Size.} By (c), $\log J\ge\max\{4,\log y\}\ge(1+\log y)/2$, and by
(b), $\log y\ge\log(L/\sqrt r)\ge0$. Hence
$\tau^2\ge(1+\log(L/\sqrt r))/8$. With (a), $G\ge c_\rho m^{p/2}$,
$\Xi\ge\Theta/2$, $\sigma^2=c_0^{2p+1}d^{-p-1/2}$ and $d\le2L$,
\[
 c_{\rm pk}n\ge c_{\rm pk}\,\frac r2\,c_\rho(\frac r{8L})^{p/2}
 \frac{c_*c_0^{2p+1}}{2\epsilon^2(2L)^{p+1/2}}\cdot
 \frac{1+\log(L/\sqrt r)}8 ,
\]
which is the claimed bound for a suitable $c_{\rm ov}$, since
$r(r/L)^{p/2}L^{-p-1/2}=r^{1+p/2}L^{-(3p+1)/2}$.
\end{proof}

The amplification lemma of Section~\ref{sec:lp_even_lower} is stated for
nonnegative weights. Its proof does not use the sign of the weights.

\begin{lemma}[Tensor amplification with real weights]
\label{lem:x-a2-signed-tensor}
Let $p$, $\epsilon$, $B$, $k$ and $r$ be as in
Lemma~\ref{lem:lp_lower_tensor_product}, with the same constants $c_p$ and
$d_p$. Suppose that every reweighting of the rows of $B$ by real weights
with relative error at most $3\epsilon$ gives nonzero weight to at least
$T$ rows. Then some rank-$r$ matrix has the property that every
reweighting of its rows by real weights with relative error at most
$\epsilon$ gives nonzero weight to at least $c_p(r/k)^{p/2}T$ rows.
\end{lemma}

\begin{proof}
Take the matrix built in the proof of
Lemma~\ref{lem:lp_lower_tensor_product}: the rows $(u_i\otimes b_j)^\top$,
where $i\in[g]$ and $b_j^\top$ runs over the rows of $B$, padded by a direct
sum with an identity matrix if their span has dimension less than $r$. Here
$g\ge c_p(r/k)^{p/2}$, and the matrix
$K_{\ell i}:=|u_\ell^\top u_i|^p$ has $K_{\ell\ell}=1$ and
$\|K-I_g\|_{\infty\to\infty}\le1/2$. Let real weights with relative error
at most $\epsilon$ be given, with $w_{ij}$ on the row $u_i\otimes b_j$, and
define $F(y):=\|By\|_p^p$ and $e_i(y):=\sum_jw_{ij}|b_j^\top y|^p-F(y)$ for $y$
in the coordinates of the row span of $B$. At the query $u_\ell\otimes y$,
with zero padding coordinates, the padding rows vanish and
$|\langle u_i\otimes b_j,u_\ell\otimes y\rangle|^p=K_{\ell i}|b_j^\top y|^p$.
Hence the input energy is $(K\mathbf1)_\ell F(y)$ and the output energy
minus the input energy is $(Ke(y))_\ell$. This identity is linear in the
weights and holds whatever their signs. The relative error gives
$\|Ke(y)\|_\infty\le\epsilon\|K\mathbf1\|_\infty F(y)\le\frac32\epsilon
F(y)$, and $\|Ke(y)\|_\infty\ge\|e(y)\|_\infty-\|(K-I_g)e(y)\|_\infty\ge
\|e(y)\|_\infty/2$. Thus $|e_i(y)|\le3\epsilon F(y)$ for every $i$ and
every $y$, that is, for each $i$ the weights $(w_{ij})_j$ reweight the rows
of $B$ by real weights with relative error at most $3\epsilon$. Each copy
therefore has at least $T$ rows of nonzero weight. Distinct copies consist
of distinct rows, so at least $gT\ge c_p(r/k)^{p/2}T$ rows have nonzero
weight.
\end{proof}

%% file: _3_app.tex
\input{90_proof}

%% file: 90_proof.tex
\section{Common Tools}
\label{app:common_tools}

This appendix proves the lemmas that several results of
Sections~\ref{sec:joint_accuracy}, \ref{sec:high_accuracy}
and~\ref{sec:joint_above_two} share. We first approximate powers of $|t|$ by polynomials
in total variation, with and without a shift. We then state partial
coloring in a subspace at small Gaussian measure and its iteration. Next
we bound the entropy of vectors obtained by integrating affine tests, and
we prove a low-crossing matching and a chaining bound with Gaussian
strips. After that come frames, coverings of the low-energy queries by
caps, and the Gaussian width of the cap residuals. The last tools are
norming sets and tensor sketches for polynomial constraints, and a volume
bound for the unit ball of the Lewis energy. The calibration and
adaptive-shift tools complete the appendix.

\subsection{Conventions}
For $d\ge1$ we write $S^d$ for the unit sphere of $\R^{d+1}$. For an
interval $I$ and a function $R$ on $I$, $\operatorname{TV}_I(R)$ is the total
variation of $R$ on $I$. It equals $\int_I|R'|$ when $R$ is absolutely
continuous. An \emph{affine test} on $\R^n$ is a set
$\{x:\langle a,x\rangle\ge b\}$ or $\{x:\langle a,x\rangle>b\}$ with
$a\in\R^n$ and $b\in\R$. A set \emph{crosses} a pair of points if it
contains exactly one of them. For a subset $\mathcal V$ of a Euclidean
space and $u>0$, $N(\mathcal V,u)$ is the least size of a set
$\mathcal N\subseteq\mathcal V$ such that every point of $\mathcal V$ lies
within distance $u$ of $\mathcal N$. A standard Gaussian vector $G$ in a
subspace $H$ is $\sum_jg_je_j$ for an orthonormal basis $(e_j)$ of $H$
and independent standard Gaussians $g_j$, and $\gamma_H$ is its law. For
every vector $y$ we have $\langle G,y\rangle=\langle G,\Pi_Hy\rangle$,
where $\Pi_H$ is the orthogonal projection onto $H$. For $w\in H$ the
variable $\langle G,w\rangle$ is centered Gaussian with variance
$\|w\|_2^2$.

\subsection{The Lewis Setup}
Let $p>0$ and let $A$ have nonzero rows and rank $r\ge1$. The proof of
Lemma~\ref{lem:lp_above_two_lewis} does not use $p>2$. For every $p>0$
it shows that, after restricting queries to the row span and a linear
change of coordinates, $\|Ax\|_p^p=rF(x)$ with
\begin{equation}\label{lem:t-eq-lewis}
 F(x):=\sum_i\nu_i|\langle u_i,x\rangle|^p,\qquad
 \sum_i\nu_i=1,\qquad \sum_i\nu_iu_iu_i^\top=I_r/r,
\end{equation}
where the $u_i$ are unit vectors of $\R^r$ and $\nu_i>0$. We call this
the \emph{Lewis setup}. If $p\le2$, then $|t|^p\ge t^2$ on $[-1,1]$ gives
$1/r\le F(x)\le1$ for unit $x$. If $p\ge2$, Jensen's inequality gives
$r^{-p/2}\le F(x)\le1$ for unit $x$. If a nonnegative vector $\mu$
satisfies $(1-\epsilon)F(x)\le\sum_i\mu_i|\langle u_i,x\rangle|^p\le
(1+\epsilon)F(x)$ for all $x$, then the weights $\mu_i/\nu_i$ on the
original rows preserve $\|Ax\|_p^p$ within the same factors, with the
same support.

\subsection{Polynomial Approximation in Total Variation}

\begin{lemma}[Approximation of powers in total variation]
\label{lem:t-tv-approx}
Let $v>0$, $\sigma\in\{0,1\}$ and $f(t):=(\operatorname{sgn}t)^\sigma|t|^v$,
with $f(0):=0$. For every integer $q\ge1$ there is a polynomial $P$ of
degree at most $q$ with $P(-t)=(-1)^\sigma P(t)$ and $P(0)=0$ such that
$R:=f-P$ satisfies
\begin{equation}\label{lem:t-eq-tv}
 \sup_{|t|\le1}|R(t)|\le\operatorname{TV}_{[-1,1]}(R)\le C_vq^{-v}.
\end{equation}
For $0<v\le1$ the constant $C_v$ can be taken absolute. Moreover
$\sup_{|t|\le1}|P(t)|\le1+C_vq^{-v}$, and the absolute values of the
monomial coefficients of $P$ sum to at most $3^{q+1}(1+C_vq^{-v})$.
\end{lemma}

For $p>0$ the case $(v,\sigma)=(p,0)$ gives an even polynomial with
$\||t|^p-P\|_{L_\infty[-1,1]}+\operatorname{TV}_{[-1,1]}(|t|^p-P)\le
2C_pq^{-p}$. For $p>1$ the derivative of $|t|^p$ is
$p(\operatorname{sgn}t)|t|^{p-1}$, so the case $(v,\sigma)=(p-1,1)$
approximates it by an odd polynomial with residual variation at most
$pC_{p-1}q^{1-p}$. For $p=1$ the derivative is $\operatorname{sgn}t$, with
a fixed value in $[-1,1]$ at zero, and we take $P:=0$. Its variation on
$[-1,1]$ is $2$, so the same bound holds with constant $2=2q^{1-p}$. This is
the case $v=0$, which the lemma excludes. The constant $C_v$ is not uniform
as $v$ increases to an integer $n\ge2$, because the factor
$1/(\varkappa-v)$ in Step~1 of the proof below, where
$\varkappa:=\lfloor v\rfloor+1$, blows up there. No use in this paper needs
uniformity in $v$.

\begin{proof}
For $t\ne0$ we have $f'(t)=v(\operatorname{sgn}t)^{\sigma+1}|t|^{v-1}$.
This is integrable on $[-1,1]$ since $v>0$, and $f(t)=\int_0^tf'$. Define
$h(\theta):=f'(\cos\theta)$ when $\cos\theta\ne0$. Then $h$ is
$2\pi$-periodic, integrable and even, and
$h(\theta+\pi)=(-1)^{\sigma+1}h(\theta)$ because $f'$ has the parity
opposite to that of $f$. Write $\mathbb T:=[-\pi,\pi]$ with the measure
$\d\theta$.

\emph{Step 1: a modulus of smoothness.}
Define $\varkappa :=\lfloor v\rfloor+1$, so that $\varkappa -1\le v<\varkappa$, and
\[
 \Delta_s^{\varkappa}h(\theta):=\sum_{i=0}^{\varkappa}(-1)^{\varkappa-i}\binom{\varkappa}{i} h(\theta+is).
\]
We show that $\|\Delta_s^{\varkappa}h\|_{L_1(\mathbb T)}\le C_v|s|^v$ for
$|s|\le\pi$. Since $\Delta_{-s}^{\varkappa}h(\theta)=(-1)^{\varkappa}\Delta_s^{\varkappa}h(\theta-\varkappa s)$,
we may assume $0<s\le\pi$. Let $Z:=\pi/2+\pi\mathbb Z$ be the zero set of
the cosine and let $\delta(\theta)\in[0,\pi/2]$ be the distance from
$\theta$ to $Z$. Then
$|\cos\theta|=\sin\delta(\theta)\in[2\delta(\theta)/\pi,\delta(\theta)]$,
and therefore
\begin{equation}\label{lem:t-eq-hbound}
 |h(\theta)|=v|\cos\theta|^{v-1}\le\tfrac\pi2\,v\,\delta(\theta)^{v-1}.
\end{equation}
Off $Z$ we have $h=v\,\psi\circ\cos$ with
$\psi(c):=(\operatorname{sgn}c)^{\sigma+1}|c|^{v-1}$, and
$|\psi^{(l)}(c)|=|(v-1)(v-2)\cdots(v-l)|\,|c|^{v-1-l}$. By the chain rule
$h^{(\varkappa)}$ is a sum of at most $C_\varkappa$ terms, each of the form
$\psi^{(l)}(\cos\theta)$ with $l\le \varkappa$ times a product of derivatives of
the cosine. Since $|\cos\theta|\le1$, every term is at most
$C_v|\cos\theta|^{v-1-\varkappa}$, so
\begin{equation}\label{lem:t-eq-hderiv}
 |h^{(\varkappa)}(\theta)|\le C_v\delta(\theta)^{v-1-\varkappa}\qquad(\theta\notin Z).
\end{equation}
For $0<v<1$ we have $\varkappa =1$ and the explicit bound
$|h'(\theta)|=v(1-v)|\cos\theta|^{v-2}|\sin\theta|
\le(\pi/2)^2v(1-v)\delta(\theta)^{v-2}$. For $v=1$ we have $\varkappa =2$, and $h$
is constant on each interval of $\R\setminus Z$, so $h''=0$ there.

Let $U:=\{\theta\in\mathbb T:\delta(\theta)<2\varkappa s\}$ and
$V:=\mathbb T\setminus U$. For $\theta\in V$, every point $\vartheta$ of
$[\theta,\theta+\varkappa s]$ satisfies
$\delta(\vartheta)\ge\delta(\theta)-\varkappa s\ge\delta(\theta)/2>0$. Thus $h$ is
smooth on this interval and
$\Delta_s^{\varkappa}h(\theta)=\int_{[0,s]^{\varkappa}}h^{(\varkappa)}(\theta+t_1+\cdots+t_\varkappa)\,\d t$.
By Eq.~\eqref{lem:t-eq-hderiv},
$|\Delta_s^{\varkappa}h(\theta)|\le2^{\varkappa+1-v}C_vs^\varkappa\delta(\theta)^{v-1-\varkappa}$. The map
$\delta$ takes each value in $(0,\pi/2)$ at four points of $\mathbb T$ and
has slope $\pm1$, so
\[
 \int_V|\Delta_s^{\varkappa}h|\le 2^{\varkappa+3-v}C_vs^\varkappa\int_{2\varkappa s}^\infty
 \delta^{v-1-\varkappa}\,\d\delta
 =\frac{2^{\varkappa+3-v}C_v(2\varkappa)^{v-\varkappa}}{\varkappa-v}\,s^v .
\]
For $\theta\in U$ and $0\le i\le \varkappa$ the point $\theta+is$ lies within
distance $3\varkappa s$ of $Z$. By Eq.~\eqref{lem:t-eq-hbound},
\[
 \int_U|\Delta_s^{\varkappa}h|\le\sum_{i=0}^{\varkappa}\binom{\varkappa}{i}\int_{\{\delta<3\varkappa s\}}|h|
 \le 2^\varkappa\cdot4\int_0^{3\varkappa s}\tfrac\pi2v\,\delta^{v-1}\,\d\delta
 =2^{\varkappa+1}\pi(3\varkappa s)^v .
\]
This proves the claim. For $0<v<1$ the explicit bound on $h'$ makes the
contribution of $V$ at most $2\pi^2v\,s^v$, and for $v=1$ it vanishes.
Hence for $0<v\le1$ the constant in Step~1 is absolute.

\emph{Step 2: a Jackson-type operator.}
Define $\ell:=\lceil(v+2)/2\rceil$, so $2\ell\ge v+2$, and for an integer
$n\ge1$ let
\[
 J_n(s):=\lambda_n^{-1}(\frac{\sin(ns/2)}{\sin(s/2)})^{2\ell},
 \qquad \int_{\mathbb T}J_n=1 .
\]
Since $(\sin(ns/2)/\sin(s/2))^2=\sum_{|j|<n}(n-|j|)e^{\mathrm ijs}$, the
kernel $J_n$ is a nonnegative even trigonometric polynomial of degree at
most $\ell(n-1)$. On $\mathbb T$ we have $|\sin(s/2)|\ge|s|/\pi$ and
$|\sin(ns/2)/\sin(s/2)|\le n$, so the base is at most $\min\{n,\pi/|s|\}$.
For $|s|\le1/n$ the base is at least $2n/\pi$, since
$\sin(n|s|/2)\ge n|s|/\pi$ and $\sin(|s|/2)\le|s|/2$. Hence
$\lambda_n\ge(2/n)(2n/\pi)^{2\ell}$. Using
$\min\{1,\pi/y\}\le(1+\pi)/(1+y)$ for $y\ge0$, we get
\begin{equation}\label{lem:t-eq-jackson}
 J_n(s)\le C_\ell\,n(1+n|s|)^{-2\ell},\qquad
 \int_{\mathbb T}J_n(s)|s|^v\,\d s
 \le2C_\ell n^{-v}\int_0^\infty(1+y)^{v-2\ell}\,\d y\le2C_\ell n^{-v},
\end{equation}
with $C_\ell:=(\pi(1+\pi)/2)^{2\ell}/2$.

Define
\[
 \mathcal J_nh(\theta):=\int_{\mathbb T}J_n(s)\sum_{i=1}^{\varkappa}(-1)^{i+1}
 \binom{\varkappa}{i} h(\theta+is)\,\d s .
\]
Since $\int J_n=1$ and
$\sum_{i=0}^{\varkappa}(-1)^i\binom{\varkappa}{i}h(\theta+is)=(-1)^{\varkappa}\Delta_s^{\varkappa}h(\theta)$, we
have
\[
 \mathcal J_nh(\theta)-h(\theta)=(-1)^{\varkappa+1}\int_{\mathbb T}J_n(s)
 \Delta_s^{\varkappa}h(\theta)\,\d s .
\]
Minkowski's integral inequality (Lemma~\ref{lem:pre_minkowski}(b)), Step~1 and Eq.~\eqref{lem:t-eq-jackson}
give
\begin{equation}\label{lem:t-eq-l1error}
 \|\mathcal J_nh-h\|_{L_1(\mathbb T)}
 \le\int_{\mathbb T}J_n(s)\|\Delta_s^{\varkappa}h\|_{L_1(\mathbb T)}\,\d s
 \le C'_vn^{-v}.
\end{equation}
For each $i$, the function $\theta\mapsto\int J_n(s)h(\theta+is)\,\d s$
is continuous by continuity of translation in $L_1$. By Fubini's theorem
its $m$th Fourier coefficient is $\hat h(m)\int J_n(s)e^{\mathrm imis}\,\d s$,
which vanishes when $|m|i>\ell(n-1)$, hence when $|m|>\ell(n-1)$. A
continuous function with finitely many nonzero Fourier coefficients is the
corresponding trigonometric polynomial, by Fej\'er's theorem. Thus
$\mathcal J_nh$ is a trigonometric polynomial of degree at most
$D:=\ell(n-1)$. It is even, because $J_n$ and $h$ are even, and
$\mathcal J_nh(\theta+\pi)=(-1)^{\sigma+1}\mathcal J_nh(\theta)$. An even
trigonometric polynomial of degree at most $D$ equals
$\sum_{m\le D}a_m\cos(m\theta)=Q(\cos\theta)$ with
$Q:=\sum_{m\le D}a_mT_m$, where $T_m$ is the Chebyshev polynomial with
$T_m(\cos\theta)=\cos(m\theta)$. Since
$\cos(m(\theta+\pi))=(-1)^m\cos(m\theta)$, only frequencies $m$ of the
parity of $\sigma+1$ occur. Hence $Q(-c)=(-1)^{\sigma+1}Q(c)$, which is the
parity of $f'$.

\emph{Step 3: integration.}
Take $n:=\lfloor(q-1)/\ell\rfloor+1$, so that $D+1\le q$ and
$n\ge q/(2\ell)$, and define $P(t):=\int_0^tQ(c)\,\d c$. Then $P$ has degree
at most $q$, the parity of $f$, and $P(0)=0$. The residual $R=f-P$ is
absolutely continuous with $R'=f'-Q$ almost everywhere. The substitution
$t=\cos\theta$, $0\le\theta\le\pi$, and Eq.~\eqref{lem:t-eq-l1error} give
\[
 \operatorname{TV}_{[-1,1]}(R)=\int_0^\pi|h(\theta)-\mathcal J_nh(\theta)|
 \sin\theta\,\d\theta\le\tfrac12\|h-\mathcal J_nh\|_{L_1(\mathbb T)}
 \le C'_v(2\ell)^vq^{-v}.
\]
Since $R(0)=0$, we also have
$|R(t)|\le\operatorname{TV}_{[-1,1]}(R)$ for $|t|\le1$. For $0<v\le1$ we
have $\ell=2$ and every constant above is absolute.

Finally $|P|\le|f|+|R|\le1+C_vq^{-v}$ on $[-1,1]$. Write
$P=\sum_{j\le q}c_jT_j$. The Chebyshev coefficients
$c_j=(2/\pi)\int_0^\pi P(\cos\theta)\cos(j\theta)\,\d\theta$ for $j\ge1$,
and $c_0$ with $1/\pi$ in place of $2/\pi$, satisfy
$|c_j|\le2\sup_{[-1,1]}|P|$. The recurrence $T_{j+1}(t)=2tT_j(t)-T_{j-1}(t)$
shows that the absolute monomial coefficients of $T_j$ sum to at most
$(1+\sqrt2)^j$. Hence the absolute monomial coefficients of $P$ sum to at
most $2\sup|P|\sum_{j\le q}(1+\sqrt2)^j\le\sqrt2(1+\sqrt2)^{q+1}\sup|P|
\le3^{q+1}\sup|P|$.
\end{proof}

\begin{lemma}[Threshold form of a residual]\label{lem:t-threshold-form}
Let $R$ be absolutely continuous on $[a,b]$ with $R(a)=0$, and let $\beta$
be the signed measure $R'(s)\,\d s$ on $(a,b]$. Then
$\|\beta\|=\operatorname{TV}_{[a,b]}(R)$ and, for every $z\in[a,b]$,
\[
 R(z)=\int\mathbf1_{\{z\ge s\}}\,\d\beta(s)=\int\mathbf1_{\{z>s\}}\,\d\beta(s).
\]
This applies to every residual of $|t|^v$ or $\operatorname{sgn}(t)|t|^v$
with $v>0$ by a polynomial, after subtracting its value at the left
endpoint.
\end{lemma}
\begin{proof}
Since $R$ is absolutely continuous, $R(z)=\int_a^zR'(s)\,\d s$ and
$\operatorname{TV}_{[a,b]}(R)=\int_a^b|R'(s)|\,\d s=\|\beta\|$. The two
forms agree because $\beta$ has no atoms. Such a residual is absolutely
continuous on every bounded interval because $v>0$.
\end{proof}

\begin{lemma}[Variation approximation below one]\label{lem:t-vc-variation}
Let $0<p\le1$ and let $q\ge1$ be an integer. There are an even polynomial
$P$ of degree at most $q$ with $P(0)=0$ and a signed Borel measure $\beta$
on $(0,1]$ with $\|\beta\|\le Cq^{-p}$ such that
\[
 |t|^p-P(t)=\int_{(0,1]}\mathbf1_{\{|t|>s\}}\,\d\beta(s)\qquad(|t|\le1).
\]
Moreover $\operatorname{TV}_{[-1,1]}(|t|^p-P)=2\|\beta\|$. Here $C$ is an
absolute constant.
\end{lemma}
\begin{proof}
Take $P$ from Lemma~\ref{lem:t-tv-approx} with $(v,\sigma)=(p,0)$, whose
constant is absolute since $p\le1$, and define $R:=|t|^p-P$. Then $R$ is even
and absolutely continuous with $R(0)=0$. For $|t|\le1$ we get
$R(t)=R(|t|)=\int_0^{|t|}R'(s)\,\d s=\int_{(0,1]}\mathbf1_{\{|t|>s\}}R'(s)\,\d s$.
Let $\beta:=R'(s)\,\d s$ on $(0,1]$. By evenness
$\|\beta\|=\operatorname{TV}_{[0,1]}(R)=\frac12\operatorname{TV}_{[-1,1]}(R)$,
which is at most $Cq^{-p}/2$.
\end{proof}

\begin{lemma}[Scaled, Taylor and shifted approximations]\label{lem:t-shifted-tv}
Let $p>0$ and let $q\ge1$ be an integer. Let $P$ be the polynomial of
Lemma~\ref{lem:t-tv-approx} for $(v,\sigma)=(p,0)$ and degree $q$, with
constant $C_p$.
\begin{enumerate}[label=(\alph*),nosep,leftmargin=*]
\item For every $\varrho>0$ the polynomial
$z\mapsto(3\varrho)^pP(z/(3\varrho))$ approximates $|z|^p$ on
$[-3\varrho,3\varrho]$ with residual variation at most
$3^pC_p\varrho^pq^{-p}$.
\item For $c>0$ and $0<\rho\le5/7$, the Taylor polynomial $T_c$ of degree
$q$ of $z\mapsto z^p$ at $c$ satisfies
\[
 \operatorname{TV}_{[(1-\rho)c,(1+\rho)c]}(z^p-T_c(z))
 \le C'_pc^p(q+2)^{\lceil p\rceil}\rho^q .
\]
\item Suppose $q\ge p$. For every $\alpha\in\R$ and $\varrho>0$ there is a
polynomial $Q_{\alpha,\varrho}$ of degree at most $q$ such that
$R_{\alpha,\varrho}(s):=|\alpha+\varrho s|^p-Q_{\alpha,\varrho}(s)$
satisfies $R_{\alpha,\varrho}(-1)=0$ and
\begin{equation}\label{lem:t-eq-shifted}
 \sup_{|s|\le1}|R_{\alpha,\varrho}(s)|\le
 \operatorname{TV}_{[-1,1]}(R_{\alpha,\varrho})\le C''_p\varrho^pq^{-p}.
\end{equation}
Moreover $\sup_{|s|\le1}|Q_{\alpha,\varrho}(s)|\le(|\alpha|+\varrho)^p
+C''_p\varrho^pq^{-p}$.
\item (Grouped form.) Let $\varrho\in(0,1]$ and $v\in\R^r$ with
$\|v\|_2\le2$. Split the unit vectors $u$ of $\R^r$ into the near group
$g_0:=\{u:|\langle u,v\rangle|\le2\varrho\}$ and the far groups
$g_{j,\pm}:=\{u:2^j\varrho<\pm\langle u,v\rangle\le2^{j+1}\varrho\}$,
$j\ge1$, of which only those with $2^j\varrho<2$ are nonempty. Each group
is an intersection of at most two affine tests in $u$. There are
polynomials $P_g$ of degree at most $q$, one for each group, and signed
measures $\beta_g$ such that for all unit $u\in g$ and all $x$ with
$\|x-v\|_2\le\varrho$,
\[
 |\langle u,x\rangle|^p=P_g(\langle u,x\rangle)+\int\mathbf1_{\{\langle
 u,x\rangle\ge s\}}\,\d\beta_g(s),\qquad
 \sum_g\|\beta_g\|\le3^pC_p\varrho^pq^{-p}+C^\flat_p(q+2)^{\lceil p\rceil}
 (5/7)^q .
\]
In particular, if $\varrho\ge r^{-b_0}$ for some $b_0>0$ and $q\ge c_0r$ with
$c_0>0$, then $\sum_g\|\beta_g\|\le(3^pC_p+1)\varrho^pq^{-p}$ for
$r\ge r_0(p,b_0,c_0)$.
\end{enumerate}
\end{lemma}
\begin{proof}
(a) The substitution $z=3\varrho t$ maps $[-1,1]$ onto
$[-3\varrho,3\varrho]$ and multiplies the residual by $(3\varrho)^p$.
Total variation is invariant under this change of variable, so the claim
follows from Eq.~\eqref{lem:t-eq-tv}.

(b) Write $z=c(1+w)$ with $|w|\le\rho<1$. The binomial series
$z^p=c^p\sum_{j\ge0}\binom pjw^j$ converges absolutely there, and $T_c$
is its truncation at $j=q$. Define $P_0:=\lceil p\rceil$. For $i<p$ we have
$|p-i|\le P_0\le P_0+i$, and for $i\ge p$ we have $|p-i|\le P_0+i$.
Hence $|\binom pj|\le\binom{P_0+j-1}{j}=\prod_{l=1}^{P_0-1}\frac{j+l}{l}
\le(j+1)^{P_0-1}$. Differentiating the tail termwise and integrating over
the interval of length $2\rho c$, and using
$j|\binom pj|\le(j+1)^{P_0}$, gives
\[
 \operatorname{TV}\le 2c^p\sum_{j>q}j|\binom pj|\rho^j
 \le2c^p\rho^q\sum_{i\ge1}(q+1+i)^{P_0}\rho^i
 \le2c^p\rho^q(q+2)^{P_0}\sum_{i\ge1}(1+i)^{P_0}(5/7)^i ,
\]
where we used $q+1+i\le(q+2)(1+i)$.

(c) If $|\alpha|\le2\varrho$, define
$Q_{\alpha,\varrho}(s):=(3\varrho)^pP((\alpha+\varrho s)/(3\varrho))$.
For $|s|\le1$ the point $(\alpha+\varrho s)/(3\varrho)$ lies in
$[-1,1]$, so the residual variation is at most $3^pC_p\varrho^pq^{-p}$,
as in (a). If $|\alpha|>2\varrho$, then
$|\alpha+\varrho s|=|\alpha|+(\operatorname{sgn}\alpha)\varrho s$ for
$|s|\le1$. Define $Q_{\alpha,\varrho}(s):=T_{|\alpha|}(|\alpha|+
(\operatorname{sgn}\alpha)\varrho s)$ and apply (b) with $c=|\alpha|$ and
$\rho=\varrho/|\alpha|<1/2$. Since $q\ge p$,
\[
 C'_p|\alpha|^p(q+2)^{\lceil p\rceil}(\varrho/|\alpha|)^q
 =C'_p\varrho^p(q+2)^{\lceil p\rceil}(\varrho/|\alpha|)^{q-p}
 \le C'_p2^p\varrho^p(q+2)^{\lceil p\rceil}2^{-q}\le C^\sharp_p\varrho^pq^{-p}.
\]
Hence Eq.~\eqref{lem:t-eq-shifted} holds with $C''_p:=\max\{3^pC_p,C^\sharp_p\}$.
In both cases add a constant to $Q_{\alpha,\varrho}$ so that the residual
vanishes at $-1$. This changes neither the degree nor the variation, and
then $|R_{\alpha,\varrho}(s)|=|R_{\alpha,\varrho}(s)-R_{\alpha,\varrho}(-1)|$
is at most the variation. The bound on $Q_{\alpha,\varrho}$ follows from
$|\alpha+\varrho s|\le|\alpha|+\varrho$.

(d) Let $\|x-v\|_2\le\varrho$ and $z:=\langle u,x\rangle$, so
$|z-\langle u,v\rangle|\le\varrho$. For $u\in g_0$ we have
$z\in[-3\varrho,3\varrho]$, and we take $P_{g_0}(z):=(3\varrho)^p
P(z/(3\varrho))$ with residual variation bounded by (a). For
$u\in g_{j,\pm}$ we have $\pm z\in[2^{j-1}\varrho,3\cdot2^j\varrho]
=[(1-\frac57)c_j,(1+\frac57)c_j]$ with $c_j:=7\cdot2^{j-2}\varrho$, and we
take $P_{g_{j,\pm}}(z):=T_{c_j}(\pm z)$ with residual variation at most
$C'_pc_j^p(q+2)^{\lceil p\rceil}(5/7)^q$ by (b). Since $c_j<7/2$ for the
nonempty far groups and the $c_j$ grow geometrically, the far groups
contribute at most $C^\flat_p(q+2)^{\lceil p\rceil}(5/7)^q$ in total. In each
group, add a constant to $P_g$ so that the residual vanishes at the left
endpoint of the interval of $z$, and apply
Lemma~\ref{lem:t-threshold-form}. For the last claim note that
$(q+2)^{\lceil p\rceil}(5/7)^q\varrho^{-p}q^p\le(q+2)^{\lceil p\rceil+p}
(5/7)^qr^{pb_0}$ tends to zero as $r\to\infty$ when $q\ge c_0r$.
\end{proof}

\subsection{Partial Coloring in a Subspace and Its Iteration}

\begin{lemma}[Partial coloring in a subspace with explicit constants]
\label{lem:t-subspace-coloring}
There are absolute constants $\delta_0,\delta_1\in(0,1)$ and an integer
$n_{\rm c}$ with the following property. Let $n\ge n_{\rm c}$, let
$H\subseteq\R^n$ be a subspace of codimension at most $\delta_0n$, and let
$\mathcal K\subseteq H$ be a closed symmetric convex set with
$\gamma_H(\mathcal K)\ge e^{-\delta_0n}$. Then some
$z\in\mathcal K\cap[-1,1]^n$ has at least $\delta_1n$ coordinates with
$|z_i|=1$.
\end{lemma}
\begin{proof}
We use Lemma~\ref{lem:pre_rothvoss} with $\varepsilon:=1/60000$ and its
$\delta$ and $n_{\rm c}$, and take $\delta_0:=\delta/2$ and
$\delta_1:=\varepsilon$.
For $\rho>0$ the set $\mathcal K\cap\rho B_2^n$ is a compact symmetric
convex subset of $H$, and its measure increases to
$\gamma_H(\mathcal K)\ge e^{-\delta n/2}>e^{-\delta n}$ as
$\rho\to\infty$. Hence $\mathcal K':=\mathcal K\cap\rho B_2^n$ satisfies
the hypotheses for large $\rho$. It has nonempty interior in $H$ because
its measure is positive.
\end{proof}

When $\gamma_H(\mathcal K)\ge1/2$ this is
Lemma~\ref{lem:subspace_partial_coloring}. The rounding steps below update
masses through \emph{units}. Let $\mu\in[0,\infty)^I$ for a finite set $I$.
A unit of $\mu$ is either a singleton $\{i\}\subseteq\supp\mu$ or an
ordered pair $(i,j)$ of distinct indices of $\supp\mu$. For a set
$\mathcal P$ of pairwise disjoint units and $z\in[-1,1]^{\mathcal P}$ we
define $\mu^z$ as follows. For a singleton unit $e=\{i\}$ define
$\mu^z_i:=\mu_i(1+z_e)$. For a pair unit $e=(i,j)$ define
$\tau_e:=\min\{\mu_i,\mu_j\}$, $\mu^z_i:=\mu_i+\tau_ez_e$ and
$\mu^z_j:=\mu_j-\tau_ez_e$. Indices in no unit keep their mass. Then
$\mu^z\ge0$. Pair units preserve the total mass, and singleton units
preserve it when $\sum_{e=\{i\}}\mu_iz_e=0$.

\begin{lemma}[Iterated partial coloring]\label{lem:t-iterated-coloring}
Let $\delta_0,\delta_1,n_{\rm c}$ be as in Lemma~\ref{lem:t-subspace-coloring}.
Let $I$ and $\mathcal X$ be sets with $I$ finite, and for $x\in\mathcal X$
let $\Phi_x:\R^I\to\R$ and $\omega(x)>0$. Let
$\mathcal S\subseteq[0,\infty)^I$, $N_*\ge\max\{5n_{\rm c},20\}$,
$a_*>0$ and $b_*>0$. Suppose that for every $\mu\in\mathcal S$ with
$N:=|\supp\mu|>N_*$ there are a set $\mathcal P$ of pairwise disjoint units of $\mu$ with
$|\mathcal P|\ge N/5$, a subspace $H\subseteq\R^{\mathcal P}$
of codimension at most $\delta_0|\mathcal P|$ and a closed symmetric convex
set $\mathcal K\subseteq H$ with $\gamma_H(\mathcal K)\ge
e^{-\delta_0|\mathcal P|}$, such that every
$z\in\mathcal K\cap[-1,1]^{\mathcal P}$ satisfies $\mu^z\in\mathcal S$ and
\[
 |\Phi_x(\mu^z)-\Phi_x(\mu)|\le a_*N^{-b_*}\omega(x)\qquad(x\in\mathcal X).
\]
Then for every $\mu^0\in\mathcal S$ there is $\mu^*\in\mathcal S$ with
$\supp\mu^*\subseteq\supp\mu^0$, $|\supp\mu^*|\le N_*$, and
\[
 |\Phi_x(\mu^*)-\Phi_x(\mu^0)|\le\frac{a_*N_*^{-b_*}}{1-(1-\delta_1/10)^{b_*}}
 \,\omega(x)\qquad(x\in\mathcal X).
\]
\end{lemma}
\begin{proof}
Let $\mu\in\mathcal S$ with $N>N_*$, and take $\mathcal P,H,\mathcal K$ as
in the hypothesis. Then $|\mathcal P|\ge N/5>n_{\rm c}$, so
Lemma~\ref{lem:t-subspace-coloring} gives
$z\in\mathcal K\cap[-1,1]^{\mathcal P}$ with at least
$\delta_1|\mathcal P|$ saturated
units, that is, units $e$ with $|z_e|=1$. Call $z_e=-1$ the deleting sign
of a singleton unit. For a pair unit $e=(i,j)$, the sign $z_e=-1$ deletes
$i$ if $\mu_i\le\mu_j$, and the sign $z_e=+1$ deletes $j$ if
$\mu_j\le\mu_i$. Every saturated unit is deleting for $z$ or for $-z$, and
$-z\in\mathcal K\cap[-1,1]^{\mathcal P}$ by symmetry. Replacing $z$ by
$-z$ if necessary, at least $\delta_1|\mathcal P|/2\ge\delta_1N/10$
indices of
$\supp\mu$ get mass zero in $\mu^z$. No index outside $\supp\mu$ gains
mass. Hence $\mu':=\mu^z\in\mathcal S$ has
$|\supp\mu'|\le(1-\delta_1/10)N$, $\supp\mu'\subseteq\supp\mu$, and
$|\Phi_x(\mu')-\Phi_x(\mu)|\le a_*N^{-b_*}\omega(x)$.

Starting from $\mu^0$, repeat this step while the support exceeds $N_*$.
The supports $N_0>N_1>\cdots$ decrease, so the process stops at some
$\mu^{n_{\rm r}}=:\mu^*$ with $N_{n_{\rm r}}\le N_*$, after $n_{\rm r}$
steps. If $n_{\rm r}=0$ there is nothing to prove. Otherwise
$N_{n_{\rm r}-1}>N_*$, and reading the steps backwards gives
$N_l\ge(1-\delta_1/10)^{-(n_{\rm r}-1-l)}N_*$ for $l<n_{\rm r}$. The
triangle inequality gives
\[
 |\Phi_x(\mu^*)-\Phi_x(\mu^0)|\le\sum_{l<n_{\rm r}}a_*N_l^{-b_*}\omega(x)
 \le a_*N_*^{-b_*}\omega(x)\sum_{j\ge0}(1-\delta_1/10)^{jb_*}.\qedhere
\]
\end{proof}

\begin{lemma}[Light rows]\label{lem:t-light-rows}
Let $\delta_0$ be as in Lemma~\ref{lem:t-subspace-coloring}. Let
$\mu\ge0$ have $N$ positive entries and total mass one, and let
$J:=\{i:0<\mu_i\le2/N\}$ be the set of light rows. Then $|J|\ge N/2$ and
$\sum_{i\in J}\mu_i^2\le2/N$. With the
singletons of $J$ as units we get $|\mathcal P|\ge N/2$. With the pairs of
a maximal matching of $J$ we get $|\mathcal P|\ge\lfloor N/4\rfloor\ge N/5$
when $N\ge20$, and the
pair weights satisfy $\tau_e\le2/N$. If
$\mathcal K=\{z\in H:\Psi(z)\le2\E[\Psi(G)]\}$ for a seminorm $\Psi$ and a
standard Gaussian vector $G$ in $H$, then
$\gamma_H(\mathcal K)\ge1/2\ge e^{-\delta_0|\mathcal P|}$ for
$|\mathcal P|\ge(\log2)/\delta_0$.
\end{lemma}
\begin{proof}
Fewer than $N/2$ entries exceed $2/N$ because the total mass is one, so
$|J|\ge N/2$, and
$\sum_{i\in J}\mu_i^2\le(2/N)\sum_{i\in J}\mu_i\le2/N$. A maximal matching
of $J$ has $\lfloor|J|/2\rfloor\ge\lfloor N/4\rfloor$ pairs, and
$\lfloor N/4\rfloor\ge N/5$ for $N\ge20$. A pair weight is the smaller of
the masses of two light rows, so $\tau_e\le2/N$. The last claim is
Markov's inequality (Lemma~\ref{lem:pre_markov}).
\end{proof}

\subsection{Entropy of Integrated Affine Tests}

\begin{lemma}[Sign patterns of affine tests]\label{lem:t-sign-patterns}
Let $h_1,\ldots,h_m$ be affine tests on $\R^n$ and, for $X\subseteq\R^n$,
let $\Pi(X)$ be the number of distinct vectors
$(\mathbf1_{\{x\in h_1\}},\ldots,\mathbf1_{\{x\in h_m\}})$ with $x\in X$.
\begin{enumerate}[label=(\alph*),nosep,leftmargin=*]
\item If $m\ge n\ge1$, then $\Pi(\R^n)\le(em/n)^n$.
\item If $n=d+1$ with $d\ge1$ and $m\ge d$, then $\Pi(S^d)\le(8em/d)^d$.
\end{enumerate}
\end{lemma}
\begin{proof}
Write $h_l:=\{x:\langle a_l,x\rangle\ge b_l\}$ or
$h_l:=\{x:\langle a_l,x\rangle>b_l\}$. Pick one witness point in $X$ for
each realized vector. Replace $b_l$ by $b_l-\epsilon_0$ for a weak test
and by $b_l+\epsilon_0$ for a strict test, where $\epsilon_0>0$ is smaller
than every nonzero value of $|\langle a_l,x\rangle-b_l|$ at the witnesses.
The affine functions $g_l(x):=\langle a_l,x\rangle-b'_l$ obtained in this
way do not vanish at any witness, and $g_l>0$ at a witness exactly when the
witness lies in $h_l$. Distinct witnesses therefore have distinct strict
sign vectors of $(g_1,\ldots,g_m)$.

(a) A function $g_l$ with $a_l=0$ is a nonzero constant. The others vanish
on hyperplanes, and points with distinct strict sign vectors lie in
distinct open regions of the arrangement of these at most $m$ hyperplanes.
Lemma~\ref{lem:small_p_regions} bounds the number of regions by
$\sum_{j\le n}\binom mj\le(em/n)^n$.

(b) Choose a unit vector $e$ of $\R^{d+1}$ different from all witnesses and
orthonormal coordinates with $e=e_{d+1}$. Stereographic projection
$\varsigma(w):=(2w,|w|^2-1)/(|w|^2+1)$ is a bijection from $\R^d$ onto
$S^d\setminus\{e\}$. The functions
$G_l(w):=(|w|^2+1)\,g_l(\varsigma(w))$ are polynomials of degree at most
two in $w\in\R^d$. They have the signs of $g_l$, and none vanishes
identically because $g_l\ne0$ at the witnesses. Points with distinct
strict sign vectors of $(G_1,\ldots,G_m)$ lie in distinct connected
components of $\R^d\setminus\bigcup_l\{G_l=0\}$. Lemma~\ref{lem:pre_warren} with degree two bounds their number by $(8em/d)^d$ for
$m\ge d$.
\end{proof}

\begin{lemma}[Packing under a trace bound]\label{lem:t-one-inclusion}
Let $d\ge1$ and $a\ge1$. Let $(\Omega,P)$ be a probability space, let $Y$
be a set with at most three elements, and let $\mathcal F$ be a family of
measurable maps $\Omega\to Y$ such that, for every $b\ge d$ and all
$\omega_1,\ldots,\omega_b\in\Omega$,
\[
 |\{(f(\omega_1),\ldots,f(\omega_b)):f\in\mathcal F\}|\le(ab/d)^d .
\]
If $0<\delta\le1$ and $f_1,\ldots,f_H\in\mathcal F$ satisfy
$\Pr[f_i\ne f_j]\ge\delta$ for all $i\ne j$, then
$H\le2(30a\log(ea)/\delta)^d$.
\end{lemma}
\begin{proof}
\emph{The one-inclusion bound.} For a finite set $W\subseteq\{0,1\}^m$ let
$V(W)$ be its VC dimension, the largest size of a set of coordinates on
which $W$ realizes all patterns. The one-inclusion graph of $W$ joins two
elements that differ in exactly one coordinate. We use that this graph has at
most $V(W)|W|$ edges. We argue by induction on $m$, the case $m=0$ being
trivial. Let $W'$ be the projection of $W$ to the first $m-1$ coordinates
and $W'':=\{y:(y,0),(y,1)\in W\}$. Then $|W|=|W'|+|W''|$ and
$V(W')\le V(W)$. If $W''\ne\emptyset$, then $V(W'')\le V(W)-1$, since a set
shattered by $W''$ together with the last coordinate is shattered by $W$.
There are $|W''|$ edges in the last coordinate. Every other edge projects
to an edge of $W'$, and two of them have the same projection only if they
are $\{(y,0),(y',0)\}$ and $\{(y,1),(y',1)\}$, in which case $\{y,y'\}$ is
an edge of $W''$. By induction the number of edges is at most
$|W''|+V(W)|W'|+(V(W)-1)|W''|=V(W)|W|$. Now let $\varpi:W\to\{1,2,\ldots\}$
be multiplicities and give the edge $\{y,y'\}$ the weight
$\min\{\varpi(y),\varpi(y')\}$. Apply the bound to the level sets
$W_l:=\{y:\varpi(y)\ge l\}$, which have $V(W_l)\le V(W)$. Since
$\min\{\varpi(y),\varpi(y')\}=\sum_{l\ge1}\mathbf1_{\{y,y'\in W_l\}}$, the
total weight is at most $V(W)\sum_y\varpi(y)$.

\emph{Two values.} Suppose $Y\subseteq\{0,1\}$. Let $V$ be the largest $b$
such that some $\omega_1,\ldots,\omega_b$ are shattered by $\mathcal F$. If
$V\ge d$, the hypothesis gives $2^V\le(aV/d)^d$. With $y:=V/d\ge1$ this
reads $y\log2\le\log a+\log y\le\log a+y/e$, so $y\le3.1\log a$. Hence
$V\le V_0:=4d\log(ea)$. Define $m:=\lceil4V_0/\delta\rceil+d$ and let
$\omega_1,\ldots,\omega_m$ be independent with law $P$. Let $W$ be the set
of traces $(f_i(\omega_1),\ldots,f_i(\omega_m))$, $1\le i\le H$, with
multiplicities $\varpi$. Then $V(W)\le V$, and the total edge weight is
at most $VH$.

We bound the expected weight from below. Fix $l\le m$ and condition on
$(\omega_j)_{j\ne l}$. The indices $i$ form clusters of equal traces off
$l$. A cluster of size $t$ splits on $\omega_l$ into $b$ indices with value
$0$ and $c$ with value $1$. This gives an edge in coordinate $l$ of weight
$\min\{b,c\}\ge bc/t$ when $b,c>0$. Since $\omega_l$ is independent of the
others, $\E[bc]=\sum \Pr[f_i\ne f_j]\ge\delta t(t-1)/2$, the sum being over
pairs in the cluster. Hence the expected weight in coordinate $l$ is at
least $\delta(H-\E[T_l])/2$, where $T_l$ is the number of distinct traces
off $l$. By exchangeability $\E[T_l]=\E[T]$, where $T$ counts the traces on
$m-1$ independent samples. Summing over $l$ gives
$m\delta(H-\E[T])/2\le VH$. Since $m\ge4V/\delta$ and $m-1\ge d$, we obtain
$H\le2\E[T]\le2(a(m-1)/d)^d$, and $a(m-1)/d\le a(16\log(ea)/\delta+1)\le
17a\log(ea)/\delta$.

\emph{Three values.} Let $Y$ carry the uniform probability and define
$g_i(\omega,y):=\mathbf1_{\{f_i(\omega)=y\}}$ on $\Omega\times Y$. If
$f_i(\omega)\ne f_j(\omega)$, then $g_i(\omega,y)\ne g_j(\omega,y)$ exactly
for $y\in\{f_i(\omega),f_j(\omega)\}$, so the disagreement probability of
$g_i,g_j$ is $\frac23\Pr[f_i\ne f_j]\ge2\delta/3$. The traces of the $g_i$
on $(\omega_1,y_1),\ldots,(\omega_b,y_b)$ are functions of the traces of
the $f_i$ on $\omega_1,\ldots,\omega_b$, so the trace hypothesis holds with
the same $a$. The two-valued case gives
$H\le2(17a\log(ea)\cdot\frac3{2\delta})^d\le2(30a\log(ea)/\delta)^d$.
\end{proof}

\begin{lemma}[A low-crossing matching on the sphere]\label{lem:t-low-crossing}
Let $d\ge2$ and let $u_1,\ldots,u_N\in S^d$, with repetitions allowed.
There are $\lfloor N/2\rfloor$ pairwise disjoint pairs of indices such that
every affine test on $\R^{d+1}$ crosses at most
\[
 C(N^{1-1/d}+d\log(eN))
\]
of the pairs $\{u_i,u_j\}$. Consequently the intersection of $\ell$ affine
tests crosses at most $\ell$ times this number of pairs. Here $C$ is
absolute.
\end{lemma}
\begin{proof}
Let $\mathcal T$ be the set of traces $\{i:u_i\in h\}$ of affine tests $h$.
On a finite set every strict test has the trace of a weak one, obtained by
raising its threshold slightly. The trace of
$\{u:\langle a,u\rangle\ge b\}$ is the pattern at the point
$(a,b)\in\R^{d+2}$ of the $N$ affine tests $\{(a,b):\langle a,u_i\rangle
-b\ge0\}$. Lemma~\ref{lem:t-sign-patterns}(a), or the bound $2^N$ when
$N<d+2$, gives $\log|\mathcal T|\le(d+2)\log(eN)$. For each
$\tau\in\mathcal T$ fix an affine test $h_\tau$ with trace $\tau$. For
$u\in S^d$ define $f_u(\tau):=\mathbf1_{\{u\in h_\tau\}}$ on $\mathcal T$.
By Lemma~\ref{lem:t-sign-patterns}(b) the family $\{f_u:u\in S^d\}$
satisfies the trace hypothesis of Lemma~\ref{lem:t-one-inclusion} with
$a=8e$. Hence, for every probability on $\mathcal T$, a family of these
functions with pairwise disagreement at least $\delta$ has at most
$2(C_1/\delta)^d$ members, where $C_1:=240e\log(8e^2)$.

We pair greedily with multiplicative weights. Give every trace weight
$\Omega(\tau):=1$ and let all indices be live. While at least two indices
are live, let $n$ be their number and $\delta_n:=\min\{1,2C_1n^{-1/d}\}$,
and view $\Omega/\Omega(\mathcal T)$ as a probability on $\mathcal T$. If
$\delta_n<1$, the live functions $f_{u_i}$ cannot be pairwise
$\delta_n$-separated, since $2(C_1/\delta_n)^d=2^{1-d}n<n$. In either case
there are live $i\ne j$ with disagreement at most $\delta_n$, which is the
normalized weight of the traces crossing $\{u_i,u_j\}$. Pair $i$ with $j$,
make them dead, and double $\Omega(\tau)$ for each trace $\tau$ crossing the
pair. This multiplies $\Omega(\mathcal T)$ by at most
$1+\delta_n\le e^{\delta_n}$.

Since $\sum_{n\le N}n^{-1/d}\le N^{1-1/d}/(1-1/d)\le2N^{1-1/d}$, the final
total weight satisfies $\log\Omega(\mathcal T)\le(d+2)\log(eN)
+4C_1N^{1-1/d}$. A trace crossing $c$ of the pairs ends with weight $2^c$,
so $c\log2$ is at most this bound. An affine test crosses exactly the pairs
crossed by its trace. Finally, if $h_1\cap\cdots\cap h_\ell$ crosses a
pair, one endpoint lies in every $h_m$ and the other misses some $h_m$,
which then crosses the pair.
\end{proof}

\begin{lemma}[Integrated coefficient vectors]\label{lem:t-integrated-vectors}
Let either $X\subseteq\R^n$ and $d:=n$, or $X\subseteq S^d$ with $d\ge1$.
Let $\ell\ge1$ be an integer, $\mathcal P$ a finite index set and
$(\mathcal Z,\mathfrak Z)$ a measurable space. For $e\in\mathcal P$ let
$\tau_e\ge0$ and let $\beta_e$ be a finite signed measure on $\mathcal Z$
with $\|\beta_e\|\le V$. Let $\phi:\mathcal P\times\mathcal Z\times X\to
[-1,1]$ take at most three values, let $(e,z)\mapsto\phi(e,z,x)$ be
measurable for each $x$, and suppose that for each $(e,z)$ there are
affine tests $h_1,\ldots,h_\ell$ on the ambient space of $X$ such that
$\phi(e,z,x)$ is a function of
$(\mathbf1_{\{x\in h_1\}},\ldots,\mathbf1_{\{x\in h_\ell\}})$. Let $\Pi$ be
an orthogonal projection of $\R^{\mathcal P}$, define
\[
 v_x:=\Pi(\tau_e\int\phi(e,z,x)\,\d\beta_e(z))_{e\in\mathcal P},
 \qquad\mathcal V:=\{v_x:x\in X\},\qquad D:=V\|\tau\|_2 .
\]
\begin{enumerate}[label=(\alph*),nosep,leftmargin=*]
\item There is a probability $P$ on $\mathcal P\times\mathcal Z$ such that
$\|v_x\|_2\le D$ and $\|v_x-v_y\|_2^2\le4D^2\Pr[\phi(\cdot,x)\ne\phi(\cdot,y)]$
for all $x,y\in X$.
\item $N(\mathcal V,u)\le\max\{1,(C_\ell D/u)^{2d}\}$ for every $u>0$,
where $C_\ell$ depends only on $\ell$.
\item If $\beta_e=\beta$ for all $e$, and for every $z\in\mathcal Z$ and
$x\in X$ at most $\chi$ indices $e$ have $\phi(e,z,x)\ne0$, then
$\|v_x\|_2\le\|\beta\|\max_e\tau_e\sqrt\chi$.
\item Let $\mathcal N\subseteq\mathcal V$ be a set such that every point of
$\mathcal V$ is within $u$ of $\mathcal N$, and let $\Pi'$ be the
orthogonal projection onto the subspace of the range of $\Pi$ orthogonal to
$\mathcal N$. Then $\|\Pi'v_x\|_2\le u$ for every $x$, and
$N(\Pi'\mathcal V,u')\le N(\mathcal V,u')$ for every $u'>0$.
\end{enumerate}
\end{lemma}
\begin{proof}
(a) Define $\bar v_x:=(\tau_e\int\phi(e,z,x)\,\d\beta_e(z))_e$, so that
$v_x=\Pi\bar v_x$. Since $|\phi|\le1$ and $\Pi$ contracts,
$\|v_x\|_2^2\le\sum_e\tau_e^2\|\beta_e\|^2\le D^2$. Let
$\Xi:=\sum_e\tau_e^2\|\beta_e\|\le V\|\tau\|_2^2$. If $\Xi=0$, all $v_x$
vanish. Otherwise let $P(\{e\}\times\d z):=\tau_e^2|\beta_e|(\d z)/\Xi$. By
the Cauchy--Schwarz inequality,
$(\int g\,\d\beta_e)^2\le\|\beta_e\|\int g^2\,\d|\beta_e|$, and
$|\phi(e,z,x)-\phi(e,z,y)|\le2$. Hence
\[
 \|v_x-v_y\|_2^2\le\sum_e\tau_e^2\|\beta_e\|\int(\phi(e,z,x)-\phi(e,z,y))^2
 \,\d|\beta_e|(z)\le4V\Xi\,\Pr[\phi(\cdot,x)\ne\phi(\cdot,y)].
\]

(b) The maps $\phi(\cdot,x)$, $x\in X$, are measurable functions on
$(\mathcal P\times\mathcal Z,P)$ with at most three values. On observations
$(e_1,z_1),\ldots,(e_b,z_b)$ their traces are determined by the patterns of
$\ell b$ affine tests at $x$. By Lemma~\ref{lem:t-sign-patterns} there are
at most $(e\ell b/n)^n$ of them in the first case and $(8e\ell b/d)^d$ in
the second, for $b\ge d$. Thus Lemma~\ref{lem:t-one-inclusion} applies
with $a:=8e\ell$. Let $0<u\le2D$ and $\delta:=u^2/(4D^2)$. Take $X'\subseteq X$
of maximal size such that the $\phi(\cdot,x)$, $x\in X'$, are pairwise
$\delta$-separated. It has at most $2(30a\log(ea)/\delta)^d$ elements, and
by maximality every $x\in X$ has some $x'\in X'$ with disagreement below
$\delta$, hence $\|v_x-v_{x'}\|_2<u$ by (a). This gives
$N(\mathcal V,u)\le(C_\ell D/u)^{2d}$ with
$C_\ell:=(240\,a\log(ea))^{1/2}$. For $u>2D$ a single point suffices
because the diameter of $\mathcal V$ is at most $2D$.

(c) Minkowski's integral inequality (Lemma~\ref{lem:pre_minkowski}(b)) gives
$\|\bar v_x\|_2\le\int\|(\tau_e\phi(e,z,x))_e\|_2\,\d|\beta|(z)
\le\|\beta\|\max_e\tau_e\sqrt\chi$, and $\Pi$ contracts.

(d) Choose $w\in\mathcal N$ with $\|v_x-w\|_2\le u$. Since $\Pi'w=0$ and
$\Pi'$ contracts, $\|\Pi'v_x\|_2=\|\Pi'(v_x-w)\|_2\le u$. A set of centers
for $\mathcal V$ at scale $u'$ maps under $\Pi'$ to a set of centers for
$\Pi'\mathcal V$ at the same scale.
\end{proof}

\begin{lemma}[Chaining with Gaussian strips]\label{lem:t-sidak-chaining}
Let $H\subseteq\R^n$ be a subspace and $G$ a standard Gaussian vector in
$H$.
\begin{enumerate}[label=(\alph*),nosep,leftmargin=*]
\item If $X_1,\ldots,X_M$ are centered jointly Gaussian with
$\E[X_l^2]\le\sigma^2$, then
\[
 \E[\max_l|X_l|]\le\sigma\sqrt{2\log(2M)}.
\]
\end{enumerate}
Let $\mathcal V\subseteq H$ be nonempty with $\sup_{v\in\mathcal V}
\|v\|_2\le R$, and define $N_j:=N(\mathcal V,2^{-j/2}R)$, assumed finite for
$j\ge1$.
\begin{enumerate}[label=(\alph*),nosep,leftmargin=*,resume]
\item We have
\[
 \E[\sup_{v\in\mathcal V}|\langle G,v\rangle|]\le2R\sum_{j\ge1}
 2^{-j/2}\sqrt{\log(2N_j)} .
\]
In particular, if
$N(\mathcal V,u)\le B_0\max\{1,(a_0R/u)^{2d}\}$ for all $u>0$, with
$B_0,a_0\ge1$ and $d\ge1$, then
\[
 \E[\sup_{v\in\mathcal V}|\langle G,v\rangle|]\le14R(\sqrt{\log(2B_0)}
 +\sqrt{d\log(ea_0)}).
\]
\item For every sequence $\lambda_j\ge1$ there are countably many
$w_l\in H$ and $c_l\ge0$ with
$\sum_lP(|\langle G,w_l\rangle|>c_l)\le\sum_jN_je^{-\lambda_j^2/2}$ such that
$\sup_{v\in\mathcal V}|\langle G,v\rangle|\le2R\sum_j\lambda_j2^{-j/2}$ on
the event $\bigcap_l\{|\langle G,w_l\rangle|\le c_l\}$. If
$N_j\le\kappa_0n_02^{dj}$ for all $j$, with $\kappa_0>0$, $n_0\ge1$ and
$d\ge1$, then $\lambda_j:=2\sqrt{dj}$ gives total failure probability at
most $\kappa_0n_0/2$ and $\sup_v|\langle G,v\rangle|\le17\sqrt d\,R$ on
the event. The strips of level $j$ have the form
$\{|\langle G,w_l\rangle|\le\lambda_j\|w_l\|_2\}$, so each fails with
probability at most $e^{-\lambda_j^2/2}$. The bound on the supremum is
deterministic: it holds with $G$ replaced by any $y\in H$ satisfying
$|\langle y,w_l\rangle|\le c_l$ for all $l$.
\item (\v{S}id\'ak's inequality for strips.) For
countably many $w_l\in H$ and $c_l\ge0$ with
$p_l:=\Pr[|\langle G,w_l\rangle|>c_l]\le1/2$,
\[
 \Pr[\bigcap_l\{|\langle G,w_l\rangle|\le c_l\}]\ge\prod_l(1-p_l)
 \ge\exp(-2\sum_lp_l).
\]
\end{enumerate}
\end{lemma}
\begin{proof}
(a) For $\lambda>0$, Jensen's inequality gives
\[
 \E[\max_l|X_l|]\le\lambda^{-1}\log\E[\sum_l(e^{\lambda X_l}
 +e^{-\lambda X_l})]\le\lambda^{-1}\log(2M)+\lambda\sigma^2/2 .
\]
Take $\lambda:=\sqrt{2\log(2M)}/\sigma$.

(b) and (c). Let $\mathcal N_0:=\{0\}$ and, for $j\ge1$, let
$\mathcal N_j\subseteq\mathcal V$ have $N_j$ points with every point of
$\mathcal V$ within $2^{-j/2}R$ of $\mathcal N_j$. Give each
$w\in\mathcal N_j$ a parent $\pi(w)\in\mathcal N_{j-1}$ with
$\|w-\pi(w)\|_2\le2^{-(j-1)/2}R$. This is possible since $w\in\mathcal V$,
and for $j=1$ since $\|w\|_2\le R$. Let $\mathcal W_j$ be the set of the at most
$N_j$ links $w-\pi(w)$, $w\in\mathcal N_j$. They lie in $H$ and have length
at most $\sqrt2\,2^{-j/2}R$. For $v\in\mathcal V$ and $J\ge1$ pick
$w\in\mathcal N_J$ with $\|v-w\|_2\le2^{-J/2}R$. Following parents from $w$
to $0$ writes $\langle G,w\rangle$ as a sum of one link per level. Letting
$J\to\infty$ gives, for every realization,
\begin{equation}\label{lem:t-eq-chain}
 \sup_{v\in\mathcal V}|\langle G,v\rangle|
 \le\sum_{j\ge1}\max_{w\in \mathcal W_j}|\langle G,w\rangle| .
\end{equation}
The left side is lower semicontinuous in $G$, hence measurable. Part (a)
bounds the expectation of the $j$th term by
$\sqrt2\,2^{-j/2}R\sqrt{2\log(2N_j)}$, which proves the first claim of (b).
Under the covering hypothesis, $N_j\le B_0a_0^{2d}2^{jd}$, so
$\sqrt{\log(2N_j)}\le\sqrt{\log(2B_0)}+\sqrt{2d\log a_0}+\sqrt{jd\log2}$.
With $\sum_j2^{-j/2}\le2.42$ and $\sum_j\sqrt j\,2^{-j/2}\le4.15$ this gives
the second claim. For (c), take the strips
$\{|\langle G,w\rangle|\le\lambda_j\|w\|_2\}$ for $w\in \mathcal W_j$. Each fails
with probability at most $\Pr[|g|>\lambda_j]\le e^{-\lambda_j^2/2}$ for a
standard Gaussian $g$, since $\lambda_j\ge1$. On their intersection
Eq.~\eqref{lem:t-eq-chain} gives the bound
$\sum_j\lambda_j\sqrt2\,2^{-j/2}R$. With $\lambda_j:=2\sqrt{dj}$ we get
$\sum_jN_je^{-2dj}\le\kappa_0n_0\sum_j(2e^{-2})^j\le\kappa_0n_0/2$ and
$2R\sum_j2\sqrt{dj}\,2^{-j/2}\le17\sqrt d\,R$.

(d) The strips $\{y\in H:|\langle y,w_l\rangle|\le c_l\}$ are symmetric and
convex. For finitely many strips the first inequality is
Lemma~\ref{lem:pre_sidak} applied to the jointly Gaussian variables
$\langle G,w_l\rangle$. For countably many strips it follows by taking a
decreasing limit. The second inequality uses $\log(1-x)\ge-2x$ for
$0\le x\le1/2$.
\end{proof}

\subsection{Frames and Caps}

\begin{lemma}[A common frame for two energies]\label{lem:t-two-energy-frame}
Let $0<p\le2$, let $u_1,\ldots,u_n\in\R^r$ be nonzero, and let
$\nu_1,\ldots,\nu_n>0$. There are a set $S\subseteq[n]$ with
$|S|\le C_{\rm fr}r/\min\{p,1\}$ and weights $w_i>0$, $i\in S$, such that
for every $x\in\R^r$
\begin{align*}
 \tfrac12\sum_i\nu_i|\langle u_i,x\rangle|^p&\le\sum_{i\in S}w_i
 |\langle u_i,x\rangle|^p\le2\sum_i\nu_i|\langle u_i,x\rangle|^p,\\
 \tfrac12\sum_i\nu_iu_iu_i^\top&\preceq\sum_{i\in S}w_iu_iu_i^\top
 \preceq2\sum_i\nu_iu_iu_i^\top .
\end{align*}
Here $C_{\rm fr}$ is absolute. In the Lewis setup the weights satisfy
$\sum_{i\in S}w_i|\langle u_i,x\rangle|^p\le2F(x)$ and
$I_r/(2r)\preceq\sum_{i\in S}w_iu_iu_i^\top\preceq2I_r/r$, and taking
traces gives $\sum_{i\in S}w_i\le2$.
\end{lemma}
\begin{proof}
We delete rows in rounds and keep two energies under control at once. In a
round the current weights are $\mu_i>0$ on a set of $N$ indices. Define
$f_i(x):=\mu_i|\langle u_i,x\rangle|^p$, $g_i(x):=\mu_i\langle u_i,x\rangle^2$,
$\mathsf F:=\sum_if_i$ and $\mathsf Q:=\sum_ig_i$. Let
\[
 s:=80r(1+\frac1{\min\{p,1\}}),\qquad\vartheta:=2\sqrt{s/N},
\]
and assume $N\ge4s$. We construct $z\in[-1,1]^N$ with at least $N/4$
entries equal to $-1$ such that, for every $x$,
\begin{equation}\label{lem:t-eq-frame-round}
 |\sum_iz_if_i(x)|\le C_1\vartheta\,\mathsf F(x),\qquad
 |\sum_iz_ig_i(x)|\le C_1\vartheta\,\mathsf Q(x),
\end{equation}
with $C_1$ absolute.

\emph{The two cones.} Let $\mathcal D_f:=\{y\in\R^N:\sum_iy_if_i(x)\le0
\text{ for all }x\}$, define $\mathcal D_g$ in the same way, and define
$\mathcal D:=\mathcal D_f\cap\mathcal D_g$, a closed convex cone. Let $G$ be
standard Gaussian in $\R^N$ and $\mathsf D$ a positive diagonal matrix. By
the symmetry $G\mapsto-G$,
$\Pr[\mathsf DG\in\mathcal D_f]=\Pr[\sum_iG_i|\langle(\mathsf D_{ii}\mu_i)^{1/p}
u_i,x\rangle|^p\ge0\text{ for all }x]$. The rows
$(\mathsf D_{ii}\mu_i)^{1/p}u_i$ have rank at most $r$, so
Lemma~\ref{lem:general_lp_cone} bounds this probability below by
$\exp(-80r/\min\{p,1\})$. In the same way, with exponent two and rows
$(\mathsf D_{ii}\mu_i)^{1/2}u_i$, we get
$\Pr[\mathsf DG\in\mathcal D_g]\ge e^{-80r}$. If $y\in\mathcal D_f$ and
$y'\le y$ coordinatewise, then $y'\in\mathcal D_f$ because $f_i\ge0$, and
similarly for $\mathcal D_g$. Hence both events are nonincreasing functions
of $G$ in the coordinatewise order. By Harris's inequality,
proved below for completeness,
\begin{equation}\label{lem:t-eq-harris}
 \Pr[\mathsf DG\in\mathcal D]\ge \Pr[\mathsf DG\in\mathcal D_f]\,
 \Pr[\mathsf DG\in\mathcal D_g]\ge e^{-s}.
\end{equation}
For Harris's inequality, let $\varphi,\psi$ be bounded, measurable and
nonincreasing on $\R^N$, and let the coordinates of $X$ be independent.
For $N=1$ and an independent copy $X'$, the product
$(\varphi(X)-\varphi(X'))(\psi(X)-\psi(X'))$ is nonnegative, and its
expectation is $2\E[\varphi\psi]-2\E[\varphi]\,\E[\psi]$. For $N>1$, condition
on the first $N-1$ coordinates. The case $N=1$ applies to the sections,
and the conditional expectations of $\varphi$ and $\psi$ are nonincreasing
functions of the first $N-1$ coordinates, so induction applies.

\emph{Every halfspace containing a small cube captures volume.} Let
$J:=\mathcal D\cap3\sqrt NB_2^N$. Fix $a\in\R^N\setminus\{0\}$ and let
$J_a:=J\cap\{y:\langle a,y\rangle\le\vartheta\|a\|_1\}$. Take
$\mathsf D_{ii}:=\min\{1,\|a\|_1/(N|a_i|)\}$ if $a_i\ne0$ and
$\mathsf D_{ii}:=1$ otherwise. Then $\|\mathsf Da\|_\infty\le\|a\|_1/N$ and
$\|\mathsf Da\|_1\le\|a\|_1$, so $\|\mathsf Da\|_2^2\le\|a\|_1^2/N$. Since
$\log y\le y$, we get $-\log\det\mathsf D\le\sum_iN|a_i|/\|a\|_1=N$. The
Gaussian tail gives $\Pr[\langle a,\mathsf DG\rangle>\vartheta\|a\|_1]\le
\exp(-\vartheta^2N/2)=e^{-2s}$. Since $\mathsf D_{ii}\le1$,
$\Pr[\|\mathsf DG\|_2>3\sqrt N]\le e^{-9N/4}\E[e^{\|G\|_2^2/4}]
=2^{N/2}e^{-9N/4}\le e^{-N}$. With Eq.~\eqref{lem:t-eq-harris}, $N\ge4s$ and
$s\ge160$ we get $\Pr[\mathsf DG\in J_a]\ge e^{-s}-e^{-2s}-e^{-N}\ge e^{-s}/2$.
The density of $\mathsf DG$ is at most $(2\pi)^{-N/2}e^N$, so
$\operatorname{vol}(J_a)\ge\frac12e^{-s-N}(2\pi)^{N/2}\ge e^{-2N}$.

\emph{Symmetrization.} For a convex body $\mathcal K$ define
$v_{\mathcal K}(b):=\operatorname{vol}(\mathcal K\cap(2b-\mathcal K))$.
Fubini's theorem gives $\int v_{\mathcal K}(b)\,\d b=2^{-N}
\operatorname{vol}(\mathcal K)^2$, and $v_{\mathcal K}$ vanishes outside
$\mathcal K$ by convexity. Applied to $\mathcal K=J_a$, this yields
$b\in J_a$ with $v_J(b)\ge v_{J_a}(b)\ge2^{-N}\operatorname{vol}(J_a)\ge
2^{-N}e^{-2N}$. Let $\mathcal U:=\{b:v_J(b)\ge2^{-N}e^{-2N}\}$. It is
closed since $v_J$ is continuous, bounded since $\mathcal U\subseteq J$,
and convex by the Brunn--Minkowski inequality
(Lemma~\ref{lem:pre_brunn_minkowski}), because
$J\cap(2b-J)$ contains $\lambda(J\cap(2b_1-J))+(1-\lambda)(J\cap(2b_2-J))$
for $b=\lambda b_1+(1-\lambda)b_2$. We just showed that $\mathcal U$ meets
every halfspace $\{\langle a,y\rangle\le\vartheta\|a\|_1\}$. If
$\mathcal U$ missed the cube $\vartheta[-1,1]^N$, a strictly separating
$a\ne0$ would give $\langle a,b\rangle>\max_{y\in\vartheta[-1,1]^N}
\langle a,y\rangle=\vartheta\|a\|_1$ for all $b\in\mathcal U$, which is
impossible. Choose $b\in\mathcal U\cap\vartheta[-1,1]^N$ and define
$T:=(J-b)\cap(b-J)$. This is a symmetric convex body of volume
$v_J(b)\ge2^{-N}e^{-2N}$. For $y\in T$ we have $b\pm y\in J$, hence
$\|y\|_2\le3\sqrt N$ and $\pm\sum_iy_if_i(x)\le-\sum_ib_if_i(x)\le
\vartheta\mathsf F(x)$, and the same with $g_i$ and $\mathsf Q$.

\emph{Partial coloring.} Let $Q_N:=[-1,1]^N$. The tiles $2\mathbf j+Q_N$,
$\mathbf j\in\mathbb Z^N$, that meet $T$ lie in $5\sqrt NB_2^N$ and have disjoint
interiors, so there are at most
$\operatorname{vol}(5\sqrt NB_2^N)/2^N\le(5\sqrt{2\pi e}/2)^N\le11^N$ of
them. Here we used $\operatorname{vol}(B_2^N)=\pi^{N/2}/\Gamma(1+N/2)$ and
$\Gamma(1+y)\ge(y/e)^y$. For $A_{\mathbf j}:=T\cap(2\mathbf j+Q_N)$, the symmetry of $T$ and
$Q_N$ gives $(A_{\mathbf j}-A_{\mathbf j})/2\subseteq T\cap Q_N$.
The Brunn--Minkowski inequality (Lemma~\ref{lem:pre_brunn_minkowski}) then gives
$\operatorname{vol}(A_{\mathbf j})\le\operatorname{vol}(T\cap Q_N)$. Hence
$\operatorname{vol}(T\cap Q_N)\ge11^{-N}\operatorname{vol}(T)\ge
(22e^2)^{-N}$. Lemma~\ref{lem:volume_partial_coloring} with $c=1/(22e^2)$
gives an absolute $C_1$ and $z\in C_1(T\cap Q_N)\cap Q_N$ with at least
$N/2$ entries of modulus one. Replacing $z$ by $-z$ if necessary, at least
$N/4$ entries equal $-1$. Since $z/C_1\in T$,
Eq.~\eqref{lem:t-eq-frame-round} holds. Multiplying $\mu_i$ by $1+z_i$ keeps
the weights nonnegative, deletes at least $N/4$ rows, and multiplies both
energies by factors in $[1-C_1\vartheta,1+C_1\vartheta]$ at every $x$.

\emph{Iteration.} Start from $\mu=\nu$ and repeat the round while
$N>N_\circ:=C's$, where $C'\ge4$ is a large absolute constant. If $n\le
N_\circ$, take $S:=[n]$ and $w:=\nu$. Otherwise let $N_0>N_1>\cdots>N_{n_{\rm r}-1}>
N_\circ\ge N_{n_{\rm r}}$ be the supports before the rounds and after the
last one. Reading backwards, $N_l\ge(4/3)^{n_{\rm r}-1-l}N_{n_{\rm r}-1}$,
so the errors $\vartheta_l:=2\sqrt{s/N_l}$ satisfy
\[
 \sum_{l<n_{\rm r}}C_1\vartheta_l\le2C_1\sqrt{s/N_\circ}\sum_{j\ge0}(3/4)^{j/2}
 \le15C_1/\sqrt{C'}\le1/4
\]
for $C'$ large. Both energies are then multiplied by factors in
$[1-\sum_lC_1\vartheta_l,\exp(\sum_lC_1\vartheta_l)]\subseteq[1/2,2]$. The
final support has at most $N_\circ\le160C'r/\min\{p,1\}$ indices, and we
let $w$ be the final weights. The statements for the Lewis setup follow
from Eq.~\eqref{lem:t-eq-lewis} and $\tr[2I_r/r]=2$.
\end{proof}

\begin{lemma}[Grid nets for $\ell_p$ balls]\label{lem:t-pball-net}
Let $0<p\le2$, let $1\le k\le n$ be integers and let $R>0$. For
$y\in\R^n$ with $\sum_i|y_i|^p\le R^p$, let $S\subseteq[n]$ index $k$
entries of largest modulus, and let $y^0$ be the vector that agrees on $S$
with $y$ rounded toward zero to the grid $Rk^{-1/p}\mathbb Z$ and vanishes
off $S$. Then $e:=y-y^0$ satisfies $|e_i|\le|y_i|$ for every $i$ and
\[
 \|e\|_\infty\le Rk^{-1/p},\qquad\|e\|_p\le R,\qquad
 \|e\|_2\le Rk^{1/2-1/p},
\]
and $\|e\|_1\le Rk^{1-1/p}$ if $p\le1$. The set $\mathcal N_{n,k}(R)$ of all
such centers $y^0$ satisfies $\log|\mathcal N_{n,k}(R)|\le k\log(en/k)+C_pk$.
\end{lemma}
\begin{proof}
On $S$ the rounding moves each entry toward zero by less than
$Rk^{-1/p}$. Hence $|e_i|\le|y_i|$ and $|e_i|\le Rk^{-1/p}$ there. Off $S$ we have $e_i=y_i$, and
$|y_i|$ is at most the $(k+1)$st largest modulus $y^*_{k+1}$. Since
$(k+1)(y^*_{k+1})^p\le\sum_i|y_i|^p\le R^p$, this gives
$|e_i|\le Rk^{-1/p}$. Hence $|e_i|\le|y_i|$ for all $i$,
$\|e\|_\infty\le Rk^{-1/p}$ and $\|e\|_p\le\|y\|_p\le R$. Since $p\le2$,
\[
 \|e\|_2^2\le\|e\|_\infty^{2-p}\|e\|_p^p\le R^2k^{1-2/p},
\]
and for $p\le1$ similarly
$\|e\|_1\le\|e\|_\infty^{1-p}\|e\|_p^p\le Rk^{1-1/p}$.

A center is determined by $S$ and by the integers $m_i$ with
$y^0_i:=Rk^{-1/p}m_i$, $i\in S$. Since $|y^0_i|\le|y_i|$, they satisfy
$\sum_{i\in S}|m_i|^p\le k\sum_i|y_i|^p/R^p\le k$. With
$\Sigma_p:=\sum_{m\in\mathbb Z}e^{-|m|^p}<\infty$, the number of
$(m_i)\in\mathbb Z^S$ with $\sum|m_i|^p\le k$ is at most
$e^k\sum_{(m_i)}e^{-\sum|m_i|^p}=(e\Sigma_p)^k$. There are
$\binom nk\le(en/k)^k$ choices of $S$.
\end{proof}

\begin{lemma}[An ensemble of unweighted frames]\label{lem:t-frame-ensemble}
In the Lewis setup let $0<p\le2$. There are an absolute constant $C_0\ge2$
and a family $\mathcal L$ of at most $C_pr\log(2r)$ lists, each a list
$\mathbf u=(u_{i_1},\ldots,u_{i_{|\mathbf u|}})$ of input directions with
repetitions allowed, such that every list satisfies
\[
 r\le|\mathbf u|\le2C_0r,\qquad\tfrac{C_0}2I_r\preceq\sum_{l}
 u_{i_l}u_{i_l}^\top\preceq2C_0I_r,
\]
and every unit $x$ has a list $\mathbf u\in\mathcal L$ with
$\sum_l|\langle u_{i_l},x\rangle|^p\le C'rF(x)$. Here $C'$ is absolute.
\end{lemma}
\begin{proof}
Let $C_0\ge2$ be a large absolute constant and $m_0:=\lceil C_0r\rceil$.
Choose an integer $b\ge0$ such that $n:=2^bm_0\ge100/\min_i\nu_i$, and
integers $\varpi_i$ with $|\varpi_i-n\nu_i|\le1\le n\nu_i/100$. Let $\mathcal U$ be the
multiset containing $u_i$ with multiplicity $\varpi_i$. For a submultiset $W$
write $A_W:=\sum_{u\in W}uu^\top$, with multiplicity. Then
$0.99(n/r)I_r\preceq A_{\mathcal U}\preceq1.01(n/r)I_r$ by
Eq.~\eqref{lem:t-eq-lewis}.

We split $\mathcal U$ repeatedly into two parts, down to depth $b$. Define
$n_j:=n/2^j$. Suppose a node $W$ at depth $j<b$ satisfies
$A_W\succeq n_jI_r/(2r)$. The vectors $A_W^{-1/2}u$, $u\in W$, have
$\sum(A_W^{-1/2}u)(A_W^{-1/2}u)^\top=I_r$ and squared lengths at most
$\delta_j:=2r/n_j$. By Lemma~\ref{lem:pre_mss} with two parts, $W$ splits into
$W_1,W_2$ with $\|\sum_{u\in W_l}A_W^{-1/2}uu^\top A_W^{-1/2}\|\le
(1/\sqrt2+\sqrt{\delta_j})^2$. Since $A_{W_1}=A_W-A_{W_2}$, both children
satisfy
\[
 (1-\varepsilon_j)\tfrac12A_W\preceq A_{W_l}\preceq(1+\varepsilon_j)
 \tfrac12A_W,\qquad\varepsilon_j:=2\sqrt{2\delta_j}+2\delta_j .
\]
Since $\delta_j=2^{j-b+1}r/m_0$, we have
$\sum_{j<b}\varepsilon_j\le10\sqrt{r/m_0}+4r/m_0\le14/\sqrt{C_0}\le0.1$
for $C_0$ large. By induction on $j$, every node at depth $j$ satisfies
$0.8\,(n_j/r)I_r\preceq A_W\preceq1.2\,(n_j/r)I_r$, which in particular
gives the hypothesis $A_W\succeq n_jI_r/(2r)$ used at each split. The
$2^b$ leaves have $n_b=m_0\in[C_0r,C_0r+1]$, so their covariance sums lie
between $(C_0/2)I_r$ and $2C_0I_r$, and their sizes $\tr[A_W]$ lie between
$r$ and $2C_0r$.

The leaves partition $\mathcal U$. A uniformly random leaf contains each
element of $\mathcal U$ with probability $2^{-b}$, so for unit $x$ its
expected $p$-energy sum is $2^{-b}\sum_i\varpi_i|\langle u_i,x\rangle|^p\le
1.01m_0F(x)$. By Markov's inequality (Lemma~\ref{lem:pre_markov}) the leaf sum is at most
$2.02m_0F(x)$ with probability at least $1/2$. Let $\mathcal N$ be a
Euclidean $r^{-2/p}$-net of $S^{r-1}$ with
$\log|\mathcal N|\le r\log(1+2r^{2/p})\le C_pr\log(2r)$. Sample
$\lceil\log_2|\mathcal N|\rceil+1$ independent uniform leaves. By a union
bound, with positive probability every $y\in\mathcal N$ has a sampled leaf
with sum at most $2.02m_0F(y)$. Fix such a sample as $\mathcal L$. For
$p\le2$ and reals $a,b$ we have $|a+b|^p\le2(|a|^p+|b|^p)$. If $x$ is a
unit vector and $y\in\mathcal N$ has $\|x-y\|_2\le r^{-2/p}$, the good leaf
$\mathbf u$ of $y$ therefore satisfies
\[
 \sum_l|\langle u_{i_l},x\rangle|^p\le2\sum_l|\langle u_{i_l},y\rangle|^p
 +2|\mathbf u|r^{-2}\le4.04m_0(2F(x)+2r^{-2})+4C_0r^{-1}.
\]
Since $F(x)\ge1/r$ and $m_0\le2C_0r$, the right side is at most
$C'rF(x)$.
\end{proof}

\begin{lemma}[Covers of the low-energy queries]\label{lem:t-cap-cover}
In the Lewis setup let $0<p<2$, let $t>0$, let $1\le k\le r$ be an integer,
and define $X_t:=\{x\in S^{r-1}:F(x)\le t/r\}$.
\begin{enumerate}[label=(\alph*),nosep,leftmargin=*]
\item (Maurey cover.) If $1\le p<2$, then $X_t$ is covered by at most
$\binom{2C_{\rm fr}r+k}{k}$ Euclidean balls of radius
$4r^{1-1/p}t^{1/p}k^{-1/2}$.
\item (Thresholded cover.) If $1\le p<2$, then $X_t$ is covered by at most
$\binom{2C_{\rm fr}r+k}{k}$ Euclidean balls of radius $6t^{1/p}k^{1/2-1/p}$.
\item (Ensemble cover.) For every $0<p<2$, $X_t$ is covered by at most
$\exp(C_p(\log(2r)+k\log(er/k)))$ Euclidean balls of radius
$C_pt^{1/p}k^{1/2-1/p}$.
\end{enumerate}
\end{lemma}
\begin{proof}
(a) Let $S,w$ be given by Lemma~\ref{lem:t-two-energy-frame}, define
$\Sigma:=\sum_{j\in S}w_ju_ju_j^\top\succeq I_r/(2r)$, and
$b_j:=\Sigma^{-1}u_j/(2r)$, so that $\|b_j\|_2\le1$. For $x\in X_t$ write
$z_j:=\langle u_j,x\rangle$. Then $x=\Sigma^{-1}\Sigma x=2r\sum_jw_jz_jb_j$.
By Hölder's inequality and $\sum_jw_j\le2$,
$2r\sum_jw_j|z_j|\le2r\,2^{1-1/p}(2t/r)^{1/p}=4r^{1-1/p}t^{1/p}=:B_t$. Thus
$x/B_t$ is a combination of the points $\pm b_j$ as in Lemma~\ref{lem:pre_maurey}. The
centers are $B_t$ times the means of $k$ points from the set
$\{0\}\cup\{\pm b_j\}$ of size at most $2|S|+1$, and there are at most
$\binom{2|S|+k}k$ such means.

(b) With the notation of (a), let $\xi:=(t/k)^{1/p}$ and split
$x=x_{\rm s}+x_{\rm l}$, where
\[
 x_{\rm s}:=\Sigma^{-1}\sum_{|z_j|\le\xi}w_jz_ju_j,\qquad
 x_{\rm l}:=2r\sum_{|z_j|>\xi}w_jz_jb_j .
\]
The matrix with
columns $\sqrt{w_j}\Sigma^{-1/2}u_j$ has operator norm one, and
$\|\Sigma^{-1/2}\|\le\sqrt{2r}$. Since $z_j^2\le\xi^{2-p}|z_j|^p$ for
$|z_j|\le\xi$ and $\sum_jw_j|z_j|^p\le2F(x)\le2t/r$,
\[
 \|x_{\rm s}\|_2\le\sqrt{2r}(\sum_{|z_j|\le\xi}w_jz_j^2)^{1/2}
 \le\sqrt{2r}\,\xi^{1-p/2}(2t/r)^{1/2}=2t^{1/p}k^{1/2-1/p}.
\]
Since $p\ge1$, we have $|z_j|\le\xi^{1-p}|z_j|^p$ for $|z_j|>\xi$, so the
coefficients of $x_{\rm l}$ in the atoms $b_j$ have absolute sum at most
$2r\xi^{1-p}(2t/r)=4t\xi^{1-p}$. Lemma~\ref{lem:pre_maurey} with this multiplier
approximates $x_{\rm l}$ within $4t\xi^{1-p}k^{-1/2}=4t^{1/p}k^{1/2-1/p}$.
The centers are $4t\xi^{1-p}$ times means of $k$ points of
$\{0\}\cup\{\pm b_j\}$, which gives the count and the radius.

(c) Let $\mathcal L$ be the family of Lemma~\ref{lem:t-frame-ensemble}.
For $x\in X_t$ take a list $\mathbf u$ with $\sum_l|\langle u_{i_l},x
\rangle|^p\le C'rF(x)\le C't$, and let $y\in\R^{|\mathbf u|}$ be its
analysis vector $y_l:=\langle u_{i_l},x\rangle$. Let $U$ be the matrix
with columns $u_{i_l}$, so $UU^\top\succeq(C_0/2)I_r$. The left inverse
$U^\dagger:=(UU^\top)^{-1}U$ satisfies $U^\dagger y=x$ and
$\|U^\dagger\|=\lambda_{\min}(UU^\top)^{-1/2}\le1$. By
Lemma~\ref{lem:t-pball-net} with $n=|\mathbf u|\ge r\ge k$ and
$R=(C't)^{1/p}$, some $y^0\in\mathcal N_{|\mathbf u|,k}(R)$ has
$\|y-y^0\|_2\le(C't)^{1/p}k^{1/2-1/p}$, so
$\|x-U^\dagger y^0\|_2\le(C't)^{1/p}k^{1/2-1/p}$. Thus the points
$U^\dagger y^0$, $y^0\in\mathcal N_{|\mathbf u|,k}(R)$, cover the queries
witnessed by $\mathbf u$ at radius $C_pt^{1/p}k^{1/2-1/p}$. The
union over $\mathbf u\in\mathcal L$ has logarithmic size at most
$\log(C_pr\log(2r))+k\log(2eC_0r/k)+C_pk$.
\end{proof}

A \emph{cap} is a triple $c=(v_c,\varrho_c,t_c)$ with $v_c\in\R^r$,
$\varrho_c\in(0,1]$ and $t_c>0$. In the Lewis setup with $p<2$, its query
set is
\[
 X_c:=\{x\in S^{r-1}:\|x-v_c\|_2\le\varrho_c,\ t_c/(2r)\le F(x)\le t_c/r\}.
\]

\begin{lemma}[Cap families and the radius bound]\label{lem:t-cap-family}
In the Lewis setup let $0<p<2$ and let $1\le k\le r$ be an integer. There
is a family $\mathcal C$ of caps such that every $x\in S^{r-1}$ lies in
$X_c$ for some $c\in\mathcal C$, and every $c\in\mathcal C$ has
$X_c\ne\emptyset$, $\|v_c\|_2\le2$, $t_c\in\{2^j:0\le j\le\lceil\log_2r
\rceil\}$ and
\begin{equation}\label{lem:t-eq-radius}
 \varrho_c=\min\{1,C_pt_c^{1/p}k^{1/2-1/p}\},\qquad
 \frac{r\varrho_c^p}{t_c}\le C_p^p\,rk^{p/2-1},\qquad
 \varrho_c\ge r^{1/2-1/p},
\end{equation}
with some $C_p\ge1$. Moreover:
\begin{enumerate}[label=(\roman*),nosep,leftmargin=*]
\item $\log|\mathcal C|\le C'_p(\log(2r)+k\log(er/k))$. For $1\le p<2$
one may take $C_p:=6$ and $|\mathcal C|\le(2+\log_2r)
\binom{2C_{\rm fr}r+k}{k}$.
\item For every $a_0>0$ there are $c_{p,a_0}>0$ and $r_{p,a_0}$ such that
$k:=\max\{1,\lfloor c_{p,a_0}r\rfloor\}$ gives $\log|\mathcal C|\le a_0r$
for $r\ge r_{p,a_0}$. Then $\varrho_c\le C_{p,a_0}t_c^{1/p}r^{1/2-1/p}$
and $r\varrho_c^p/t_c\le C_{p,a_0}r^{p/2}$.
\end{enumerate}
\end{lemma}
\begin{proof}
For unit $x$ we have $1/r\le F(x)\le1$, so the tiers
$t/(2r)\le F(x)\le t/r$ with $t=2^j$, $0\le j\le\lceil\log_2r\rceil$, cover
$S^{r-1}$. Fix a tier $t$ and define
$\varrho_t:=\min\{1,C_pt^{1/p}k^{1/2-1/p}\}$, where $C_p\ge1$ is at
least the radius constant of Lemma~\ref{lem:t-cap-cover}(c), or $C_p=6$
with part (b) of that lemma when $p\ge1$. If $\varrho_t=1$, use the single cap $(0,1,t)$, whose ball
contains $S^{r-1}$, and drop it if its query set is empty. Otherwise take the cover of $X_t$ from
Lemma~\ref{lem:t-cap-cover}, enlarge its radii to $\varrho_t$, and keep the
caps with $X_c\ne\emptyset$. A kept cap contains a unit vector within
$\varrho_t\le1$ of its center, so $\|v_c\|_2\le2$. In both cases
$\varrho_c\le C_pt^{1/p}k^{1/2-1/p}$ gives
$r\varrho_c^p/t_c\le C_p^prk^{p/2-1}$. Since $t\ge1$, $C_p\ge1$, $k\le r$
and $1/2-1/p<0$, we also have $\varrho_c\ge\min\{1,r^{1/2-1/p}\}
=r^{1/2-1/p}$. Adding over at most $2+\log_2r$ tiers gives (i). For (ii)
use (i) with $k\log(er/k)\le2c_{p,a_0}\log(e/c_{p,a_0})\,r$ and take
$c_{p,a_0}$ small and then $r$ large. The radius bounds follow from
$k\ge c_{p,a_0}r/2$.
\end{proof}

\begin{lemma}[Gaussian width of cap residuals]\label{lem:t-cap-width}
In the Lewis setup let $0<p<2$, let $q\ge2$ be an integer, and let
$\mathcal C$ be a family as in Lemma~\ref{lem:t-cap-family} for some
$1\le k\le r$. For $c\in\mathcal C$ and each row $i$ let
$Q_{c,i}:=Q_{\alpha,\varrho}$ and $R_{c,i}:=R_{\alpha,\varrho}$ be given by
Lemma~\ref{lem:t-shifted-tv}(c) with $\alpha:=\langle u_i,v_c\rangle$ and
$\varrho:=\varrho_c$. For $x\in X_c$ define
$s_{c,i}(x):=\langle u_i,x-v_c\rangle/\varrho_c\in[-1,1]$, so that
$|\langle u_i,x\rangle|^p=Q_{c,i}(s_{c,i}(x))+R_{c,i}(s_{c,i}(x))$. Let $J$
be a set of rows, $\mu\in[0,\infty)^J$, $H\subseteq\R^J$ a subspace and $G$
a standard Gaussian vector in $H$. Then
\begin{equation}\label{lem:t-eq-round}
 \E[\sup_{c\in\mathcal C}\sup_{x\in X_c}\frac1{F(x)}
 |\sum_{i\in J}\mu_iG_iR_{c,i}(s_{c,i}(x))|]
 \le C_p\,rk^{p/2-1}q^{-p}\|\mu\|_2(\sqrt r+\sqrt{\log(2|\mathcal C|)}).
\end{equation}
If moreover $\sum_{i\in J}\mu_iy_iQ_{c,i}(s_{c,i}(x))=0$ for all
$y\in H$, $c\in\mathcal C$ and $x\in X_c$, then the left side
of Eq.~\eqref{lem:t-eq-round} bounds
$\E[\sup_{x\in S^{r-1}}F(x)^{-1}|\sum_{i\in J}\mu_iG_i|\langle u_i,x
\rangle|^p|]$.
\end{lemma}
\begin{proof}
By Lemma~\ref{lem:t-threshold-form} and Eq.~\eqref{lem:t-eq-shifted}, for
$|s|\le1$ we have $R_{c,i}(s)=\int\mathbf1_{\{s\ge z\}}\,\d\beta_{c,i}(z)$
with $\|\beta_{c,i}\|\le C''_p\varrho_c^pq^{-p}$. For fixed $(i,z)$ the set
$\{x:s_{c,i}(x)\ge z\}=\{x:\langle u_i,x\rangle\ge\langle u_i,v_c\rangle
+\varrho_cz\}$ is an affine test on $\R^r$. Fix $c$ and apply
Lemma~\ref{lem:t-integrated-vectors} with $X:=X_c\subseteq\R^r$, $d:=r$,
$\ell:=1$, $\mathcal P:=J$, $\tau_i:=(2r/t_c)\mu_i$, $\beta_i:=\beta_{c,i}$,
$\phi(i,z,x):=\mathbf1_{\{s_{c,i}(x)\ge z\}}$ and $\Pi:=\Pi_H$. The
resulting vectors $v^c_x$ satisfy $\langle G,v^c_x\rangle=(2r/t_c)
\sum_{i\in J}\mu_iG_iR_{c,i}(s_{c,i}(x))$. By Eq.~\eqref{lem:t-eq-radius},
\[
 D_c:=C''_p\varrho_c^pq^{-p}\frac{2r}{t_c}\|\mu\|_2\le
 D:=2C''_pC_p^p\,rk^{p/2-1}q^{-p}\|\mu\|_2 .
\]
Let $\mathcal V$ be the union over $c\in\mathcal C$ of the families
$\{v^c_x:x\in X_c\}$. By Lemma~\ref{lem:t-integrated-vectors}(b) it has
radius at most $D$ and
$N(\mathcal V,u)\le|\mathcal C|\max\{1,(C_1D/u)^{2r}\}$ for all $u>0$,
where $C_1$ is absolute. Then
Lemma~\ref{lem:t-sidak-chaining}(b) gives
\[
 \E[\sup_{v\in\mathcal V}|\langle G,v\rangle|]\le14D(\sqrt{\log(2|
 \mathcal C|)}+\sqrt{r\log(eC_1)}).
\]
Since $F(x)\ge t_c/(2r)$ on $X_c$, this
proves Eq.~\eqref{lem:t-eq-round}. Under the additional hypothesis, every
$x\in S^{r-1}$ lies in some $X_c$, and there
$\sum_i\mu_iG_i|\langle u_i,x\rangle|^p=\sum_i\mu_iG_iR_{c,i}(s_{c,i}(x))$
because $G\in H$.
\end{proof}

\subsection{Norming Sets and Tensor Sketches}

\begin{lemma}[Norming sets]\label{lem:t-norming-set}
Let $q,r\ge1$ be integers and $\rho_q:=1/(2q(q+1))$.
\begin{enumerate}[label=(\alph*),nosep,leftmargin=*]
\item Let $\mathcal N\subseteq S^{r-1}$ be such that every point of
$S^{r-1}$ is within geodesic distance $\rho_q$ of $\mathcal N$. Then
$\sup_{S^{r-1}}|P|\le2\max_{\mathcal N}|P|$ for every polynomial $P$ on
$\R^r$ of degree at most $q$. Such sets exist with
$|\mathcal N|\le(8(q+1)^2)^r$.
\item There is a set $Y\subseteq B_2^r$ with $|Y|\le(8(q+1)^2)^{r+1}$ such
that $\sup_{B_2^r}|P|\le2\max_Y|P|$ for every polynomial $P$ on $\R^r$ of
degree at most $q$.
\end{enumerate}
\end{lemma}
\begin{proof}
(a) Let $M:=\sup_{S^{r-1}}|P|$, attained at $x_*$, and take $y\in\mathcal N$
at geodesic distance $\gamma\le\rho_q$ from $x_*$. If $\gamma>0$, write
$y=\cos\gamma\,x_*+\sin\gamma\,w$ with a unit vector $w\perp x_*$. The
function $T(\theta):=P(\cos\theta\,x_*+\sin\theta\,w)$ is a trigonometric
polynomial $\sum_{|m|\le q}c_me^{\mathrm im\theta}$ with $|T|\le M$, so
$|c_m|\le M$ and $|T'|\le\sum_{|m|\le q}|m||c_m|\le q(q+1)M$. Hence
$|P(y)|\ge|P(x_*)|-\gamma q(q+1)M\ge M/2$. For existence take a maximal
set of points of $S^{r-1}$ with pairwise geodesic distances at least
$\rho_q$. Their Euclidean distances are at least $2\rho_q/\pi$, so the
balls of radius $\rho_q/\pi$ around them are disjoint and lie in
$(1+\rho_q/\pi)B_2^r$. The volume ratio gives
$|\mathcal N|\le(1+\pi/\rho_q)^r\le(8(q+1)^2)^r$.

(b) The coordinate projection $\pi(y_1,\ldots,y_{r+1}):=(y_1,\ldots,y_r)$
maps $S^r$ onto $B_2^r$. For $P$ on $\R^r$ of degree at most $q$, the
polynomial $P\circ\pi$ on $\R^{r+1}$ has the same degree bound and
$\sup_{S^r}|P\circ\pi|=\sup_{B_2^r}|P|$. Apply (a) on $S^r$ and let $Y$ be
the image under $\pi$ of the set obtained there.
\end{proof}

\begin{lemma}[Tensor sketches of polynomial constraints]\label{lem:t-sketch}
Let $a,q,r'\ge1$ be integers, $0<\tau<1/2$ and $\mathsf A\ge1$. Let
$Y\subseteq B_2^{r'}$ be finite. For $i=1,\ldots,M$ let
$w_{i,1},\ldots,w_{i,a}\in B_2^{r'}$, and let $P_i$ be a polynomial in $a$
variables of degree at most $q$ in each variable whose monomial
coefficients have absolute values summing to at most $\mathsf A$. Then the
$M\times|Y|$ matrix
$(P_i(\langle w_{i,1},y\rangle,\ldots,\langle w_{i,a},y\rangle))_{i,y}$
is within $\tau$ entrywise of a matrix of rank at most $d_{\rm sk}$, where
\[
 \log(2d_{\rm sk})\le C_a(q+\log(2r')+\log\mathsf A+\log(1/\tau)
 +\log\log(e|Y|)).
\]
More precisely, one may take $d_{\rm sk}:=1+(q+1)^ad_0$ with
\[
 \log d_0\le C+2\log(8\mathsf A/\tau)+\log(2ar')
 +\log\log(e+16aq\mathsf A/\tau)+\log\log(e|Y|),
\]
where $C$ is an absolute constant.
\end{lemma}

For $a=1$ and $\sup_{[-1,1]}|P_i|\le1$, the last paragraph of the proof of
Lemma~\ref{lem:t-tv-approx} gives $\mathsf A=3^{q+1}$. The lemma is used
with $Y$ the set of Lemma~\ref{lem:t-norming-set}(b), with the net of
Lemma~\ref{lem:t-norming-set}(a), which lies in $S^{r'-1}\subseteq
B_2^{r'}$, or with any finite net of queries in the unit ball. For a true
polynomial $Q$ of degree at most $q$ on $\R^{r'}$ and $Y$ as in
Lemma~\ref{lem:t-norming-set}(b), control on $Y$ also gives control on the
sphere, since $\sup_{S^{r'-1}}|Q|\le\sup_{B_2^{r'}}|Q|\le2\max_Y|Q|$.

\begin{proof}
Write $P_i(t)=\sum_{\mathbf j}a_{i,\mathbf j}t_1^{j_1}\cdots t_a^{j_a}$ over
multi-indices $\mathbf j\in\{0,\ldots,q\}^a$, and define
$|\mathbf j|:=j_1+\cdots+j_a$. For vectors $w_1,\ldots,w_a,y$,
\[
 \prod_{l\le a}\langle w_l,y\rangle^{j_l}
 =\langle w_1^{\otimes j_1}\otimes\cdots\otimes w_a^{\otimes j_a},\,
 y^{\otimes|\mathbf j|}\rangle
\]
in the Euclidean tensor power $(\R^{r'})^{\otimes|\mathbf j|}$, and both
tensors have norm at most one when all vectors lie in $B_2^{r'}$. Define
$b:=\tau/(8\mathsf A)$. Let $\mathcal U\subseteq B_2^{r'}$ be a set such
that every point of $B_2^{r'}$ is within $b/(aq)$ of $\mathcal U$, with
$|\mathcal U|\le(1+2aq/b)^{r'}$. We need $\mathcal U\subseteq B_2^{r'}$ so
that $|\langle w,y\rangle|\le1$ for $w\in\mathcal U$ and $y\in Y$. For each
$i$ and $l$ choose $w'_{i,l}\in\mathcal U$ with
$\|w_{i,l}-w'_{i,l}\|_2\le b/(aq)$. All inner products involved have
modulus at most one, so telescoping over the $|\mathbf j|\le aq$ factors
gives
\[
 |\prod_l\langle w_{i,l},y\rangle^{j_l}
 -\prod_l\langle w'_{i,l},y\rangle^{j_l}|\le|\mathbf j|\,\frac b{aq}
 \le b\qquad(y\in Y).
\]

For each $\mathbf j\ne0$ apply
Lemma~\ref{lem:pre_jl} with distortion $\varepsilon:=b/3$ to the
finite set consisting of $0$, the tensors
$w_1^{\otimes j_1}\otimes\cdots\otimes w_a^{\otimes j_a}$ with
$w_l\in\mathcal U$, and the tensors $y^{\otimes|\mathbf j|}$ with $y\in Y$.
It has at most $1+|\mathcal U|^a+|Y|$ points of norm at most one. This
gives a map $f_{\mathbf j}$ into $\R^{d_0}$, with
$d_0:=\lceil24\varepsilon^{-2}\log(2+|\mathcal U|^a+|Y|)\rceil$, that
preserves all squared distances in the set up to factors $1\pm\varepsilon$.
Define $g_{\mathbf j}:=f_{\mathbf j}-f_{\mathbf j}(0)$. For points $\xi,\eta$
of the set, polarization and $\|\xi-\eta\|_2\le2$ give
\[
 |\langle g_{\mathbf j}(\xi),g_{\mathbf j}(\eta)\rangle-\langle\xi,\eta\rangle|
 \le\tfrac\varepsilon2(\|\xi\|_2^2+\|\eta\|_2^2+\|\xi-\eta\|_2^2)
 \le3\varepsilon=b .
\]
Let $\xi_{i,\mathbf j}:=w'^{\otimes j_1}_{i,1}\otimes\cdots\otimes
w'^{\otimes j_a}_{i,a}$ and define
\[
 \widetilde M_{i,y}:=a_{i,0}+\sum_{\mathbf j\ne0}a_{i,\mathbf j}\langle
 g_{\mathbf j}(\xi_{i,\mathbf j}),g_{\mathbf j}(y^{\otimes|\mathbf j|})\rangle .
\]
Each term with $\mathbf j\ne0$ is a matrix of rank at most $d_0$, so
$\widetilde M$ has rank at most $1+((q+1)^a-1)d_0$. By the two displays,
$|\widetilde M_{i,y}-P_i(\langle w_{i,1},y\rangle,\ldots,\langle w_{i,a},y
\rangle)|\le2b\sum_{\mathbf j}|a_{i,\mathbf j}|\le2b\mathsf A=\tau/4$.

Finally define $d_{\rm sk}:=1+(q+1)^ad_0$. We have $\log(1/b)=\log(8\mathsf A/\tau)$
and $\log(2+|\mathcal U|^a+|Y|)\le\log3+ar'\log(1+2aq/b)+\log(e|Y|)$.
Using $\log(x+y)\le\log(1+x)+\log(1+y)$ for $x,y\ge0$, we get
$\log d_0\le C+2\log(1/b)+\log(2ar')+\log\log(e+2aq/b)+\log\log(e|Y|)+1$,
which is the explicit form, and hence the first bound with a constant $C_a$
depending only on $a$.
\end{proof}

\subsection{The Volume of the Lewis Unit Ball}

\begin{lemma}[Volume bound from the geometric Brascamp--Lieb inequality]
\label{lem:t-brascamp-lieb}
Let $p>0$, let $b_1,\ldots,b_m$ be unit vectors of $\R^s$, and let $c_i>0$
satisfy $\sum_ic_ib_ib_i^\top=I_s$. Define $\Upsilon(y):=\sum_ic_i|\langle b_i,y
\rangle|^p$. Then
\[
 \operatorname{vol}_s\{y\in\R^s:\Upsilon(y)\le1\}\le\operatorname{vol}_s(B_p^s),
 \qquad B_p^s:=\{y:\sum_{j\le s}|y_j|^p\le1\}.
\]
\end{lemma}
\begin{proof}
Lemma~\ref{lem:pre_brascamp_lieb}
states that
\[
 \int_{\R^s}\prod_if_i(\langle b_i,y\rangle)^{c_i}\,\d y\le
 \prod_i(\int_\R f_i)^{c_i}
\]
for nonnegative measurable functions $f_i$ on $\R$. With
$f_i(t):=e^{-|t|^p}$ and $\sum_ic_i=\tr[I_s]=s$ it gives
$\int e^{-\Upsilon}\le(2\Gamma(1+1/p))^s$. The vectors $b_i$ span $\R^s$, so $\Upsilon$
is continuous, positive off the origin and $p$-homogeneous. Writing
$e^{-\Upsilon(y)}=\int_{\Upsilon(y)}^\infty e^{-\lambda}\,\d\lambda$ and using Tonelli's
theorem and homogeneity,
\[
 \int_{\R^s}e^{-\Upsilon}=\int_0^\infty e^{-\lambda}\operatorname{vol}_s\{\Upsilon\le
 \lambda\}\,\d\lambda=\Gamma(1+s/p)\operatorname{vol}_s\{\Upsilon\le1\}.
\]
The same computation for $\Upsilon_0(y):=\sum_j|y_j|^p$ gives
$(2\Gamma(1+1/p))^s=\Gamma(1+s/p)\operatorname{vol}_s(B_p^s)$. Comparing the
two identities proves the claim.
\end{proof}

\begin{lemma}[Inverse moment of the Lewis energy]\label{lem:t-inverse-energy}
In the Lewis setup let $p>0$, and let $\omega$ be the normalized surface
measure on $S^{r-1}$. Then
$\operatorname{vol}_r\{x:rF(x)\le1\}\le\operatorname{vol}_r(B_p^r)$ and
\[
 \|F^{-1}\|_{L_{r/p}(\omega)}=(\int_{S^{r-1}}F^{-r/p}\,\d\omega)
 ^{p/r}\le C_pr^{p/2},\qquad
 C_p:=pe\,(3/2)^{p/2}(2\Gamma(1+1/p)/\sqrt\pi)^p .
\]
Equivalently $E:=r^{p/2}F$ satisfies $\|E^{-1}\|_{L_{r/p}(\omega)}\le
C_p$.
\end{lemma}
\begin{proof}
By Eq.~\eqref{lem:t-eq-lewis}, $rF(x)=\sum_i(r\nu_i)|\langle u_i,x\rangle|^p$
with $\sum_i(r\nu_i)u_iu_i^\top=I_r$, so
Lemma~\ref{lem:t-brascamp-lieb} with $s=r$ gives the volume bound. The set
$\{rF\le1\}$ is star-shaped with radial function $(rF(\theta))^{-1/p}$ on
the sphere, and $F>0$ there. Integration in polar coordinates gives
\[
 \operatorname{vol}_r\{rF\le1\}=\operatorname{vol}_r(B_2^r)\int_{S^{r-1}}
 (rF)^{-r/p}\,\d\omega ,
\]
hence $\int F^{-r/p}\,\d\omega\le r^{r/p}\operatorname{vol}(B_p^r)/
\operatorname{vol}(B_2^r)$. For $y>0$ we have
$\Gamma(1+y)\ge\int_y^\infty t^ye^{-t}\,\d t\ge(y/e)^y$, and
$\Gamma(1+y)=\E[\mathcal E^y]\le(\E[\mathcal E^{\lceil y\rceil}])^{y/\lceil
y\rceil}\le\lceil y\rceil^y\le(1+y)^y$ for an exponential variable
$\mathcal E$, by Lyapunov's inequality and $n!\le n^n$. The volume formulas
$\operatorname{vol}(B_p^r)=(2\Gamma(1+1/p))^r/\Gamma(1+r/p)$ and
$\operatorname{vol}(B_2^r)=\pi^{r/2}/\Gamma(1+r/2)$ then give
\[
 (\frac{\operatorname{vol}(B_p^r)}{\operatorname{vol}(B_2^r)})^{p/r}
 =(\frac{2\Gamma(1+1/p)}{\sqrt\pi})^p
 \frac{\Gamma(1+r/2)^{p/r}}{\Gamma(1+r/p)^{p/r}}
 \le(\frac{2\Gamma(1+1/p)}{\sqrt\pi})^p(\frac{3r}2)^{p/2}
 \frac{pe}{r}.
\]
Multiplying by $r$ gives the bound on $\|F^{-1}\|_{L_{r/p}(\omega)}$, and
dividing by $r^{p/2}$ gives the bound for $E$.
\end{proof}

The remaining common tools are a sampling lemma, which turns a common
approximation space for the row functions into a sparsifier of
original rows, and three lower-bound constructions built from shifted
scalar tests.

\subsection{Calibrated Sampling}
Let $A$ have nonzero rows $a_i^\top$, $i\in I$, where repeated rows
are allowed, and rank $r$. Lemma~\ref{lem:lp_above_two_lewis}, which
holds for every $p>0$, gives a linear map $x\mapsto y$ from $\R^n$ onto
$\R^r$, probabilities $\nu_i$ and unit vectors $u_i\in\R^r$ such that,
with $c_i:=r\nu_i>0$,
\begin{equation}\label{lem:s-eq-lewis}
 a_i^\top x=c_i^{1/p}u_i^\top y,\qquad
 \sum_{i\in I}c_iu_iu_i^\top=I_r,\qquad
 \sum_{i\in I}c_i=r .
\end{equation}
We call the $c_i$ Lewis weights. If several rows are copies of one
original row, a reweighting of the copies is combined into a
reweighting of that row.

\begin{lemma}[Calibrated sampling]\label{lem:s-cal}
Fix $0<p<2$, $\beta\ge1$ and $B_0>0$, and define
$\theta_p:=\min\{p,1\}$. There are an absolute constant
$c_{\rm cal}>0$ and a constant $C\ge1$ depending only on $p$, $\beta$
and $B_0$ with the following property. Let $A$ be as
in Eq.~\eqref{lem:s-eq-lewis}, and let $\pi$ be a probability vector on $I$
with
\begin{equation}\label{lem:s-eq-importance}
 \pi_i\ge\frac{c_i}{\beta r}\qquad(i\in I).
\end{equation}
Define $\kappa_x(i):=|a_i^\top x|^p/\pi_i$, so that
$F(x):=\|Ax\|_p^p=\E_\pi[\kappa_x]$, and $\mathbb B:=\{x:F(x)\le1\}$. Let
$\mathcal V$ be a space of real functions on $I$ with
$1\le\dim\mathcal V\le e^{L/8}$, and let $\delta\ge r^{-B_0}$ satisfy
\begin{equation}\label{lem:s-eq-cal-hyp}
 \operatorname{dist}_{L_2(\pi)}(\kappa_x,\mathcal V)\le\delta
 \qquad(x\in \mathbb B).
\end{equation}
Suppose that $L\ge C\log(2r)$.
\begin{enumerate}[label=(\alph*),nosep,leftmargin=*]
\item There are a set $S\subseteq I$ with
$|S|\le C\epsilon^{-2}r\delta^2\log(2r)$ and positive weights
$\widehat w_i$, $i\in S$, such that
$|\sum_{i\in S}\widehat w_i|a_i^\top x|^p-F(x)|\le\epsilon F(x)$ for
every $x$, and $\sum_{i\in S}\widehat w_i\pi_if(i)=\E_\pi[f]$ for every
$f\in\mathcal V$.
\item There are a finite set $\widetilde I$, a map
$\iota:\widetilde I\to I$, a probability vector $\widetilde\pi$ on
$\widetilde I$ and a function $\omega:\widetilde I\to(0,2]$,
constructed from $\pi$ and a fixed basis of $\mathcal V$ and not from
the sample, with
$\sum_{\iota(\tilde\imath)=i}\omega(\tilde\imath)
\widetilde\pi(\tilde\imath)=\pi_i$ for every $i$, and with the
following property. Let $M\ge e^{1.8L}$ be an integer with
$\delta\le c_{\rm cal}\epsilon\sqrt M$, and let $\widehat\pi$ be the
empirical law of $M$ independent draws from $\widetilde\pi$. For
queries $x,x'$ define
$\Delta_{xx'}:=\omega\cdot(\kappa_x-\kappa_{x'})\circ\iota$ and
$\widehat d(x,x'):=\|\Delta_{xx'}\|_{L_2(\widehat\pi)}
+\|\Delta_{xx'}\|_{L_2(\widetilde\pi)}$. Suppose that, outside an
event of probability at most $e^{-3L/2}$,
\begin{equation}\label{lem:s-eq-entropy}
 \frac1{\sqrt M}\int_0^{8\delta}
   \sqrt{\log N(\mathbb B,\widehat d,t)}\,\d t\le c_{\rm cal}\epsilon ,
\end{equation}
where $N(\mathbb B,\widehat d,t)$ is the least number of $\widehat d$-balls of
radius $t$ with centers in $\mathbb B$ that cover $\mathbb B$. Then the conclusion
of~(a) holds with $|S|\le M$ in place of the stated bound.
\end{enumerate}
\end{lemma}

\begin{proof}
Throughout, $C_1,C_2,\ldots$ denote constants depending only on $p$,
$\beta$ and $B_0$, and the constant $C$ in the lower limit for $L$
is taken large in terms of them. We prove~(b) first.

\emph{Pointwise bounds.} Write $y,y'$ for the images of the queries
$x,x'$ under the map in Eq.~\eqref{lem:s-eq-lewis}. Since
$|u_i^\top y|\le\|y\|_2$ and $p<2$,
\[
 \|y\|_2^2=\sum_ic_i|u_i^\top y|^2
 \le\|y\|_2^{2-p}\sum_ic_i|u_i^\top y|^p=\|y\|_2^{2-p}F(x),
\]
so $\|y\|_2\le1$ on $\mathbb B$. By Eq.~\eqref{lem:s-eq-importance},
$\kappa_x(i)=(c_i/\pi_i)|u_i^\top y|^p\le\beta r|u_i^\top y|^p$. For
$|s|,|s'|\le1$ we have $||s|^p-|s'|^p|\le2|s-s'|^{\theta_p}$, by the
mean value theorem if $1\le p<2$ and by subadditivity of $t\mapsto t^p$
if $p<1$. Hence for $x,x'\in \mathbb B$ and every $i$,
\begin{equation}\label{lem:s-eq-cal-pointwise}
 0\le\kappa_x(i)\le\beta r,\qquad
 |\kappa_x(i)-\kappa_{x'}(i)|\le2\beta r\|y-y'\|_2^{\theta_p}.
\end{equation}
In particular $\|\kappa_x\|_{L_2(\pi)}^2\le\beta r\,\E_\pi[\kappa_x]
\le\beta r$ on $\mathbb B$. As $0\in\mathcal V$, we may assume
$\delta\le\sqrt{\beta r}$.

\emph{A lift with bounded importance.} The construction below uses
only $\pi$ and a fixed basis of $\mathcal V$. Let $k:=\dim\mathcal V$, fix a
basis of $\mathcal V$, and write $\Psi(i)\in\R^k$ for the vector of its
values at $i$. Let $o_1,\ldots,o_{2k}$ be the vectors
$\pm e_1,\ldots,\pm e_k$ of $\R^k$. Define $\widetilde I:=I\times[2k]$
and $\iota(i,l):=i$, give $(i,l)$ the mass $\mu(i,l):=\pi_i/(2k)$, and
define $v(i,l):=(\Psi(i),o_l)\in\R^{2k}$. By Eq.~\eqref{lem:s-eq-importance}
every $\pi_i$ is positive, and the basis is linearly independent, so
$\E_\mu[\Psi\Psi^\top]$ is positive definite. The averages over $l$ of
$o_l$ and of $\Psi(i)o_l^\top$ vanish, and the average of $o_lo_l^\top$
is $I_k/k$. Hence $\E_\mu[vv^\top]$ is positive definite and
$c_\mu:=\min_{\|e\|_2=1}\E_\mu[|e^\top v|]>0$.

For symmetric positive definite $Q$ define
$\Upsilon(Q):=\E_\mu[\|Qv\|_2]-\log\det Q$. If $e$ is a unit eigenvector
for the largest eigenvalue of $Q$, then
$\|Qv\|_2\ge|e^\top Qv|=\|Q\|_{\rm op}|e^\top v|$, so
$\Upsilon(Q)\ge c_\mu\|Q\|_{\rm op}-2k\log\|Q\|_{\rm op}$. Also
$\Upsilon(Q)\to\infty$ when $Q$ approaches a singular matrix within a
bounded set. Thus $\Upsilon$ attains its minimum. It is strictly
convex, since $Q\mapsto\|Qv\|_2$ is convex and $-\log\det$ is strictly
convex, so the minimizer $Q$ is unique. Every $v(i,l)$ is nonzero, so
$\Upsilon$ is differentiable. For symmetric $H$ the derivative of
$\Upsilon(Q+tH)$ at $t=0$ is $\tr[H\Sigma Q]-\tr[HQ^{-1}]$, where
$\Sigma:=\E_\mu[vv^\top]/\|Qv\|_2$. Since $\tr[H\Sigma Q]=\tr[H(\Sigma Q+Q\Sigma)]/2$ for
symmetric $H$, the minimizer satisfies $\Sigma Q+Q\Sigma=2Q^{-1}$. In an
orthonormal eigenbasis of $Q$ with eigenvalues $q_a$ this reads
$\Sigma_{ab}(q_a+q_b)=2\delta_{ab}/q_a$. Hence $\Sigma=Q^{-2}$ and
\begin{equation}\label{lem:s-eq-l1-lewis}
 \E_\mu[\frac{(Qv)(Qv)^\top}{\|Qv\|_2}]=Q\Sigma Q=I_{2k}.
\end{equation}
For a signed permutation matrix $\varsigma$ of $\R^k$, the bijection
$(i,l)\mapsto(i,l')$ with $o_{l'}=\varsigma o_l$ preserves $\mu$ and
maps $v$ to $Uv$, where $U:=I_k\oplus\varsigma$ is orthogonal. Hence
$\Upsilon(U^\top QU)=\Upsilon(Q)$, and uniqueness gives $U^\top QU=Q$.
Taking $\varsigma=-I_k$ shows that the off-diagonal blocks of $Q$
vanish. Taking all signed permutations shows that the lower diagonal
block is $q_0I_k$ for some $q_0>0$.

Define $\lambda_Q:=\|Qv\|_2$. Then $\lambda_Q\ge q_0\|o_l\|_2=q_0$
everywhere. The trace of the lower block of Eq.~\eqref{lem:s-eq-l1-lewis}
gives $q_0^2\E_\mu[1/\lambda_Q]=k$, so $k\le q_0$. The full trace gives
$\E_\mu[\lambda_Q]=2k$. Define
\[
 \widetilde\pi:=\frac{\lambda_Q\mu}{2k},\qquad
 \omega:=\frac{2k}{\lambda_Q}\le\frac{2k}{q_0}\le2,\qquad
 \Phi:=\frac{Qv}{\lambda_Q}.
\]
Then $\widetilde\pi$ is a probability vector and
$\omega\widetilde\pi=\mu$. Hence
$\E_{\widetilde\pi}[\omega\cdot f\circ\iota]=\E_\pi[f]$ for every
function $f$ on $I$, and in particular
$\sum_{\iota(\tilde\imath)=i}\omega(\tilde\imath)
\widetilde\pi(\tilde\imath)=\pi_i$. Moreover $\Phi$ is a unit
vector, $\Phi=\omega Qv/(2k)$, and Eq.~\eqref{lem:s-eq-l1-lewis} gives
$\E_{\widetilde\pi}[\Phi\Phi^\top]=I_{2k}/(2k)$.

Let $\mathcal G$ be the span of the coordinates of $\Phi$, a space
of functions on $\widetilde I$. As $Q$ is invertible,
$\mathcal G=\{\omega\cdot a^\top v:a\in\R^{2k}\}$, so
$\omega\cdot f\circ\iota\in\mathcal G$ for every $f\in\mathcal V$. Let
$P$ be the orthogonal projection of $L_2(\widetilde\pi)$ onto
$\mathcal G$. Then $Pf=2k(\E_{\widetilde\pi}[f\Phi])^\top\Phi$.
Since $\E_{\widetilde\pi}[(e^\top\Phi)^2]=1/(2k)$ for every unit
vector $e$, Cauchy--Schwarz gives the pointwise bound
\begin{equation}\label{lem:s-eq-proj}
 |Pf|\le2k\sup_{\|e\|_2=1}\E_{\widetilde\pi}[f\,e^\top\Phi]
     \le\sqrt{2k}\,\|f\|_{L_2(\widetilde\pi)} .
\end{equation}

\emph{Sampling and correction.} Let $\xi_1,\ldots,\xi_M$ be
independent draws from $\widetilde\pi$ with empirical law
$\widehat\pi$, and define
$\widehat\Gamma:=\E_{\widehat\pi}[\Phi\Phi^\top]$. When
$\widehat\Gamma$ is invertible, define
$\zeta:=\widehat\Gamma^{-1}(\E_{\widetilde\pi}[\Phi]
-\E_{\widehat\pi}[\Phi])$ and
\[
 \widehat w_i:=\frac1{M\pi_i}\sum_{l:\,\iota(\xi_l)=i}
   \omega(\xi_l)(1+\zeta^\top\Phi(\xi_l)).
\]
For a function $f$ on $I$ this gives
$\sum_i\widehat w_i\pi_if(i)=\E_{\widehat\pi}
[(1+\zeta^\top\Phi)\,\omega\cdot f\circ\iota]$. For
$a\in\R^{2k}$,
\[
 \E_{\widehat\pi}[(1+\zeta^\top\Phi)a^\top\Phi]
 =a^\top\E_{\widehat\pi}[\Phi]+a^\top\widehat\Gamma\zeta
 =a^\top\E_{\widetilde\pi}[\Phi] .
\]
Thus the reweighted sample integrates every function of $\mathcal G$
exactly. Taking $f\in\mathcal V$ gives the exactness assertion
of~(a), since
$\E_{\widetilde\pi}[\omega\cdot f\circ\iota]=\E_\pi[f]$. For a query
$x$ let $R_x:=(I-P)(\omega\cdot\kappa_x\circ\iota)$. Since
$F(x)=\E_{\widetilde\pi}[\omega\cdot\kappa_x\circ\iota]$ and
$\omega\cdot\kappa_x\circ\iota-R_x\in\mathcal G$, the error is
\begin{equation}\label{lem:s-eq-cal-error}
 \sum_i\widehat w_i|a_i^\top x|^p-F(x)
 =(\E_{\widehat\pi}[R_x]-\E_{\widetilde\pi}[R_x])
   +\zeta^\top\E_{\widehat\pi}[\Phi R_x].
\end{equation}
Orthogonality gives $\E_{\widetilde\pi}[\Phi R_x]=0$, so the second
term is a product of two sampling errors.

For $x\in \mathbb B$ and $f\in\mathcal V$, the function $R_x$ is the error of
the best approximation of $\omega\cdot\kappa_x\circ\iota$ from
$\mathcal G$, which contains $\omega\cdot f\circ\iota$. As
$\omega\le2$,
\[
 \|R_x\|_{L_2(\widetilde\pi)}^2
 \le\E_{\widetilde\pi}[\omega^2\cdot(\kappa_x-f)^2\circ\iota]
 \le2\E_{\widetilde\pi}[\omega\cdot(\kappa_x-f)^2\circ\iota]
 =2\|\kappa_x-f\|_{L_2(\pi)}^2 ,
\]
so Eq.~\eqref{lem:s-eq-cal-hyp} gives $\|R_x\|_{L_2(\widetilde\pi)}\le
\sqrt2\delta$. By Eq.~\eqref{lem:s-eq-cal-pointwise} and Eq.~\eqref{lem:s-eq-proj},
with $b_R:=8\beta r\sqrt{2k}$,
\begin{equation}\label{lem:s-eq-cal-res}
 |R_x|\le b_R,\qquad
 |R_x-R_{x'}|\le b_R\|y-y'\|_2^{\theta_p}\qquad(x,x'\in \mathbb B).
\end{equation}

\emph{Four events.} We use $k\le e^{L/8}$ and $M\ge e^{1.8L}$
throughout. Let $\Omega_1$ be the event
$\|2k\widehat\Gamma-I\|_{\rm op}\le1/2$. The summands
$2k\Phi(\xi_l)\Phi(\xi_l)^\top-I$ are independent and centered,
with squared Hilbert--Schmidt norm $4k^2-2k$. Hence
$\E[\|2k\widehat\Gamma-I\|_{\rm HS}^2]\le4k^2/M$, and Chebyshev gives
$\Pr[\Omega_1^c]\le16k^2/M\le16e^{-1.55L}$. On $\Omega_1$ we have
$\|\widehat\Gamma^{-1}\|_{\rm op}\le4k$ and
$\E_{\widehat\pi}[f^2]\le\frac32\E_{\widetilde\pi}[f^2]$ for $f\in\mathcal G$.
Let $\Omega_2$ be the event
$\|\E_{\widehat\pi}[\Phi]-\E_{\widetilde\pi}[\Phi]\|_2\le10M^{-1/2}$.
As $\E[\|\E_{\widehat\pi}[\Phi]-\E_{\widetilde\pi}[\Phi]\|_2^2]\le1/M$,
we have $\Pr[\Omega_2^c]\le1/100$. On $\Omega_1\cap\Omega_2$,
$\|\zeta\|_2\le40kM^{-1/2}\le40e^{-0.775L}<1/2$, so every
$\widehat w_i$ is positive.

Let $\mathcal N\subseteq \mathbb B$ be maximal such that the images $y$ of its
points are $e^{-4L/\theta_p}$-separated. These images lie in the unit
ball of $\R^r$, so Lemma~\ref{lem:pre_nets}
gives $|\mathcal N|\le(3e^{4L/\theta_p})^r$, and every
$x\in \mathbb B$ has some $x'\in\mathcal N$ with
$\|y-y'\|_2\le e^{-4L/\theta_p}$. Let $\Omega_3$ be the event
$\E_{\widehat\pi}[R_x^2]\le3\delta^2$ for every $x\in\mathcal N$. The
summands $R_x(\xi_l)^2$ lie in $[0,b_R^2]$ and have mean at most
$2\delta^2$, so Hoeffding's inequality (Lemma~\ref{lem:pre_hoeffding})
and a union bound give
$\Pr[\Omega_3^c]\le|\mathcal N|\exp(-2M\delta^4/b_R^4)$. Let
$\Omega_4$ be the event $|\E_{\widehat\pi}[\Phi_aR_x]|\le e^{-0.8L}$
for every $x\in\mathcal N$ and every coordinate $a$ of $\Phi$. These
summands lie in $[-b_R,b_R]$ and have mean zero, so
$\Pr[\Omega_4^c]\le4k|\mathcal N|\exp(-Me^{-1.6L}/(2b_R^2))$. Now
$b_R^2\le128\beta^2r^2e^{L/8}$, $\delta\ge r^{-B_0}$ and
$\log(4k|\mathcal N|)\le L+r(2+4L/\theta_p)$. Therefore
\[
 \frac{2M\delta^4}{b_R^4}\ge
   \frac{e^{1.55L}}{8192\beta^4r^{4+4B_0}},\qquad
 \frac{Me^{-1.6L}}{2b_R^2}\ge\frac{e^{0.075L}}{256\beta^2r^2},
\]
and both failure probabilities are at most $e^{-2L}$ once $C$ is large
in terms of $p$, $\beta$ and $B_0$. No upper limit on $L$ is used here
or below.

On $\Omega_3$, every $x\in \mathbb B$ has $x'\in\mathcal N$ as above,
and Eq.~\eqref{lem:s-eq-cal-res} gives
$\|R_x\|_{L_2(\widehat\pi)}\le\sqrt3\delta+b_Re^{-4L}\le2\delta$,
because $b_Re^{-4L}\le16\beta re^{-3.9L}\le0.2r^{-B_0}$. On
$\Omega_1\cap\Omega_2\cap\Omega_4$, every $x\in \mathbb B$ satisfies
\[
 |\zeta^\top\E_{\widehat\pi}[\Phi R_x]|
 \le\|\zeta\|_2(\sqrt{2k}\,e^{-0.8L}+b_Re^{-4L})
 \le\frac{80k\sqrt{2k}}{\sqrt M}e^{-0.8L}
 \le120e^{-1.51L}\le\frac\epsilon4 .
\]

\emph{The main term.} Define
$Z:=\sup_{x\in \mathbb B}|(\E_{\widehat\pi}[R_x]-\E_{\widetilde\pi}[R_x])|$. The map
$x\mapsto R_x$ is continuous, so the supremum may be taken over an
increasing sequence of finite sets $\mathbb B'\subseteq \mathbb B$. By monotone
convergence it suffices to bound the expectation of the supremum over
$\mathbb B'$ uniformly in $\mathbb B'$. Let $\sigma_1,\ldots,\sigma_M$ be independent
uniform signs, independent of the draws. Symmetrization
(Lemma~\ref{lem:pre_symmetrization}) in $\ell_\infty(\mathbb B')$, applied to the
centered vectors $(R_x(\xi_l)-\E_{\widetilde\pi}[R_x])_{x\in \mathbb B'}$,
together with $|\E_{\widetilde\pi}[R_x]|\le\sqrt2\delta$ and
$\E[|\sum_l\sigma_l|]\le\sqrt M$, gives
\[
 \E[\sup_{x\in \mathbb B'}|(\E_{\widehat\pi}[R_x]-\E_{\widetilde\pi}[R_x])|]
 \le\frac2{\sqrt M}\E[\sup_{x\in \mathbb B'}|X_x|]+\frac{2\sqrt2\delta}{\sqrt M},
 \qquad X_x:=\frac1{\sqrt M}\sum_l\sigma_lR_x(\xi_l).
\]
Conditionally on the draws, $(X_x)$ is a centered process with
sub-gaussian increments,
$\|X_x-X_{x'}\|_{\psi_2}\le C_H\|R_x-R_{x'}\|_{L_2(\widehat\pi)}$ with
$C_H$ absolute (Lemma~\ref{lem:pre_subgaussian}(b)). Let
$\Omega_5$ be the event in Eq.~\eqref{lem:s-eq-entropy}, and define
$\Omega_*:=\Omega_1\cap\Omega_3\cap\Omega_5$, so that
$\Pr[\Omega_*^c]\le16e^{-1.55L}+e^{-2L}+e^{-1.5L}\le3e^{-1.5L}$. Since
$R_x-R_{x'}=(I-P)\Delta_{xx'}$, on $\Omega_1$
\[
 \|R_x-R_{x'}\|_{L_2(\widehat\pi)}
 \le\|\Delta_{xx'}\|_{L_2(\widehat\pi)}
   +\sqrt{3/2}\,\|P\Delta_{xx'}\|_{L_2(\widetilde\pi)}
 \le\sqrt{3/2}\,\widehat d(x,x').
\]
On $\Omega_3$ the class $\{R_x:x\in \mathbb B\}$ has diameter at most $4\delta$
in $L_2(\widehat\pi)$. Covering numbers of $\mathbb B'$ with centers in $\mathbb B'$
are at most those of $\mathbb B$ at half the radius. Points of $\mathbb B'$ at
$L_2(\widehat\pi)$-distance zero have identical $X_x$, so after
identifying them the pseudometric becomes a metric. Dudley's inequality
(Lemma~\ref{lem:pre_subgaussian}(c)), with an absolute constant $C_D$, therefore
gives on $\Omega_*$
\[
 \E_\sigma[\sup_{x\in \mathbb B'}X_x]
 \le C_DC_H\int_0^{4\delta}\sqrt{\log N(\mathbb B,\widehat d,t/\sqrt6)}\,\d t
 \le\sqrt6C_DC_H\sqrt M\,c_{\rm cal}\epsilon .
\]
The process is symmetric and $\E_\sigma[X_{x_0}]=0$ for a fixed
$x_0\in \mathbb B'$. Hence $\E_\sigma[\sup_{\mathbb B'}|X_x|]\le\E_\sigma[|X_{x_0}|]
+2\E_\sigma[\sup_{\mathbb B'}X_x]\le2\delta+2\sqrt6C_DC_H\sqrt Mc_{\rm cal}
\epsilon$ on $\Omega_*$. Off $\Omega_*$ we use
$M^{-1/2}\E_\sigma[\sup|X_x|]\le b_R$. Altogether, using
$\delta\le c_{\rm cal}\epsilon\sqrt M$,
\[
 \E[Z]\le(7+10C_DC_H)\,c_{\rm cal}\epsilon+6b_Re^{-1.5L}.
\]
Choose $c_{\rm cal}:=(40(7+10C_DC_H))^{-1}$. Since
$6b_Re^{-1.5L}\le48\sqrt2\beta re^{-1.43L}\le\epsilon/40$ for large
$C$, we get $\E[Z]\le\epsilon/20$, and Markov's inequality (Lemma~\ref{lem:pre_markov}) gives
$\Pr[Z>\epsilon/2]\le1/10$.

\emph{Conclusion of~(b).} The events $\Omega_1,\ldots,\Omega_4$ and
$\{Z\le\epsilon/2\}$ fail with total probability less than one. On
their intersection the weights are positive and,
by Eq.~\eqref{lem:s-eq-cal-error}, every $x\in \mathbb B$ satisfies
$|\sum_i\widehat w_i|a_i^\top x|^p-F(x)|\le\epsilon/2+\epsilon/4$.
Both sides are positively homogeneous of degree $p$ in $x$, and both
vanish when $F(x)=0$. Hence the relative error is at most $\epsilon$
for every query. The support has at most $M$ points.

\emph{Proof of~(a).} By Eq.~\eqref{lem:s-eq-cal-pointwise} and $\omega\le2$,
$\widehat d(x,x')\le8\beta r\|y-y'\|_2^{\theta_p}$ on $\mathbb B$. A maximal
subset of $\mathbb B$ whose images are $\rho$-separated, $0<\rho\le1$, has at
most $(3/\rho)^r$ points (Lemma~\ref{lem:pre_nets}). Hence, with probability one,
$\log N(\mathbb B,\widehat d,t)\le r\log3+(r/\theta_p)\log(8\beta r/t)$ for
$0<t\le8\beta r$. Define $t_1:=8e^2\beta r$. For $0<a\le t_1/e$ we have
$\int_0^a\sqrt{\log(t_1/t)}\,\d t\le2a\sqrt{\log(t_1/a)}$, since the
derivative of the right side is at least $\sqrt{\log(t_1/a)}$. With
$a:=8\delta\le8\sqrt{\beta r}\le t_1/e$ and $\delta\ge r^{-B_0}$,
\[
 \int_0^{8\delta}\sqrt{\log N(\mathbb B,\widehat d,t)}\,\d t
 \le8\delta\sqrt{r\log3}
   +16\delta\sqrt{\frac r{\theta_p}\log\frac{e^2\beta r}\delta}
 \le C_1\delta\sqrt{r\log(2r)} .
\]
Take $M:=\lceil C_2\epsilon^{-2}r\delta^2\log(2r)\rceil$ with
$C_2:=(C_1^2+1)/c_{\rm cal}^2$. Then Eq.~\eqref{lem:s-eq-entropy} holds with
probability one and $\delta\le c_{\rm cal}\epsilon\sqrt M$. Moreover
$M\ge C_2e^{2L}r^{1-2B_0}/8\ge e^{1.8L}$ for large $C$. Part~(b)
gives~(a) with $|S|\le M\le2C_2\epsilon^{-2}r\delta^2\log(2r)$.
\end{proof}

\begin{lemma}[Sharper entropy factors]\label{lem:s-cal-entropy}
In the setting of Lemma~\ref{lem:s-cal}(b), the following hold. First,
$\widehat d(x,x')\le\sqrt2\|\kappa_x-\kappa_{x'}\|_{L_2(\xi)}$ for the
positive measure
$\xi(i):=\sum_{\iota(\tilde\imath)=i}\omega(\tilde\imath)^2
(\widehat\pi(\tilde\imath)+\widetilde\pi(\tilde\imath))$ on $I$, whose
total mass is at most $8$. Second,
$\E_{\widetilde\pi}[\omega^2\cdot\kappa_x\circ\iota]\le2F(x)$, and
on $\mathbb B$ the summands $\omega^2\cdot\kappa_x\circ\iota$ are at most
$4\beta r$ and satisfy the H\"older bound Eq.~\eqref{lem:s-eq-cal-pointwise}
with an extra factor $4$. Third, the draws are independent with law
$\widetilde\pi$ and $M\ge e^{1.8L}$. Fourth, outside an event of
probability at most $(3e^{4L/\theta_p})^r\exp(-M/(32\beta^2r^2))\le
e^{-2L}$,
\[
 \sup_{x\in\mathbb B}\E_{\widehat\pi}[\omega^2\cdot\kappa_x\circ\iota]\le3,
 \qquad\text{hence}\qquad
 \sup_{x\in\mathbb B}\sum_i\xi(i)\kappa_x(i)\le5 .
\]
\end{lemma}
\begin{proof}
Write $\Delta:=\Delta_{xx'}$. For the first fact,
$\widehat d(x,x')\le\sqrt2(\|\Delta\|_{L_2(\widehat\pi)}^2
+\|\Delta\|_{L_2(\widetilde\pi)}^2)^{1/2}
=\sqrt2\|\kappa_x-\kappa_{x'}\|_{L_2(\xi)}$, and the mass of $\xi$ is at
most $8$ because $\omega\le2$. For the second, $\omega\le2$ gives
$\E_{\widetilde\pi}[\omega^2\cdot\kappa_x\circ\iota]\le
2\E_{\widetilde\pi}[\omega\cdot\kappa_x\circ\iota]=2F(x)$, and the other
two bounds follow from $\omega\le2$
and Eq.~\eqref{lem:s-eq-cal-pointwise}. The third fact is part of the
setting. For the fourth, at a point of the net $\mathcal N$ of the proof
of Lemma~\ref{lem:s-cal} the summands
$\omega^2\cdot\kappa_x\circ\iota(\xi_l)$ lie in $[0,4\beta r]$ and have
mean at most $2$. Hoeffding's inequality (Lemma~\ref{lem:pre_hoeffding})
bounds the probability that
their empirical mean exceeds $5/2$ by $\exp(-M/(32\beta^2r^2))$. For
$x\in\mathbb B$ and its net point $x'$, the H\"older bound gives
$|\omega^2(\kappa_x-\kappa_{x'})\circ\iota|\le8\beta re^{-4L}\le1/2$.
The union bound over the net and $M\ge e^{1.8L}$ give the failure
probability $e^{-2L}$ once $C$ is large.
\end{proof}

\subsection{Scalar Tests and Gaussian Facts}
Fix once and for all a nonnegative even smooth function $\varrho$ on
$\R$ with support in $\{1\le|\xi|\le2\}$ and $\int\varrho=1$, and define
\[
 \chi(y):=\int_\R\varrho(\xi)e^{\mathrm i\xi y}\,\d\xi .
\]
Constants depending on $\varrho$ are treated as absolute. Throughout,
$G$, $G_1$, $G_2$, $Z$ and $W$ denote standard Gaussian variables,
$\phi(t):=e^{-t^2/2}/\sqrt{2\pi}$, and $\gamma_1$ is the standard
Gaussian measure on $\R$. We call $p>0$ non-even if it is not an even
integer.

\begin{lemma}[A band-limited scalar test]\label{lem:s-filter}
The function $\chi$ is real, even and Schwartz, with $\chi(0)=1$ and
$|\chi|\le1$, and its Fourier transform is supported in
$\{1\le|\xi|\le2\}$. For every non-even $p>0$ there is a real even
Schwartz function $\psi$ with $\int t^j\psi(t)\,\d t=0$ for every
integer $j\ge0$ such that, for all $w>0$ and $s,z\in\R$,
\begin{equation}\label{lem:s-eq-filter}
 w^p\chi(\frac{z-s}w)=\int_\R|z-s-wt|^p\psi(t)\,\d t .
\end{equation}
\end{lemma}

\begin{proof}
The first assertions hold because $\chi$ is the Fourier transform of
an even smooth compactly supported probability density. Let
$n:=\lfloor p/2\rfloor$, let $T_n(u):=\sum_{j=0}^n(-1)^ju^{2j}/(2j)!$
be the Taylor polynomial of $\cos$ of degree $2n$, and define
$R_n:=\cos-T_n$. Then $R_0=\cos-1\le0$, with equality only on
$2\pi\mathbb Z$. Since $R_n''=-R_{n-1}$ and $R_n(0)=R_n'(0)=0$, we have
$R_n(u)=-\int_0^u(u-v)R_{n-1}(v)\,\d v$, and induction shows that
$(-1)^{n+1}R_n\ge0$ on $\R$, with equality only on $2\pi\mathbb Z$ if
$n=0$ and only at $0$ if $n\ge1$. In particular $R_n$ is not
identically zero on any interval. Also $|R_n(u)|\le u^{2n+2}/(2n+2)!$
and $|R_n(u)|\le(n+2)(1+u^{2n})$. As $2n<p<2n+2$, the integral
$I_p:=\int_0^\infty R_n(u)u^{-p-1}\,\d u$ converges absolutely and is
nonzero. Substituting $u=\xi|z|$ and using that $R_n$ is even,
\begin{equation}\label{lem:s-eq-cosine}
 |z|^p=\frac1{I_p}\int_0^\infty R_n(\xi z)\,\xi^{-p-1}\,\d\xi
 \qquad(z\in\R).
\end{equation}
Let $\psi$ be the even function with Fourier transform
$\widehat\psi(\xi)=2I_p|\xi|^{p+1}\varrho(\xi)$, that is,
\[
 \psi(t)=\frac{2I_p}\pi\int_0^\infty\xi^{p+1}\varrho(\xi)\cos(\xi t)\,\d\xi .
\]
It is real, even and Schwartz, and all its moments vanish since
$\widehat\psi$ vanishes near the origin.

Fix $y\in\R$. Since $R_n$ has constant sign,
$\int_0^\infty|R_n(\xi(y-t))|\xi^{-p-1}\,\d\xi=|I_p||y-t|^p$, which is
integrable against $|\psi(t)|\,\d t$. By Eq.~\eqref{lem:s-eq-cosine} and
Fubini,
\[
 \int|y-t|^p\psi(t)\,\d t
 =\frac1{I_p}\int_0^\infty\xi^{-p-1}
   \int\psi(t)R_n(\xi(y-t))\,\d t\,\d\xi .
\]
As $T_n(\xi(y-t))$ is a polynomial in $t$ and the moments of $\psi$
vanish, the inner integral equals
$\int\psi(t)\cos(\xi(y-t))\,\d t$, which is
$\cos(\xi y)\widehat\psi(\xi)$ because $\psi$ is even. Hence
\[
 \int|y-t|^p\psi(t)\,\d t=2\int_0^\infty\varrho(\xi)\cos(\xi y)\,\d\xi
 =\chi(y).
\]
Replacing $y$ by $(z-s)/w$ and multiplying by $w^p$
gives Eq.~\eqref{lem:s-eq-filter}.
\end{proof}

\begin{lemma}[Gaussian facts]\label{lem:s-gauss}
Let $h_n:=He_n/\sqrt{n!}$ be the normalized Hermite polynomials,
which form an orthonormal basis of $L_2(\gamma_1)$, and for
$f\in L_2(\gamma_1)$ write $\widehat f_n:=\E[f(G)h_n(G)]$.
\begin{enumerate}[label=(\alph*),nosep,leftmargin=*]
\item Let $(G_1,G_2)$ be a centered Gaussian pair with unit variances
and correlation $u$. For all $f,g\in L_2(\gamma_1)$,
\[
 \E[f(G_1)g(G_2)]=\sum_{n\ge0}\widehat f_n\widehat g_nu^n .
\]
\item Let $E_1,E_2\subseteq\R^k$ be subspaces with orthogonal
projections $P_1,P_2$, define $\lambda:=\|P_1P_2\|_{\rm op}$, and let $g$
be standard Gaussian in $\R^k$. Let $W_a:\R^k\to E_a$, $a=1,2$, be
odd measurable maps with $\E[\|W_a(g)\|_2^2]<\infty$ such that $W_a(g)$
depends only on $P_ag$. Then
\[
 |\E[\langle W_1(g),W_2(g)\rangle]|
 \le\lambda^2(\E[\|W_1(g)\|_2^2]\,\E[\|W_2(g)\|_2^2])^{1/2}.
\]
\item Let $\tau\ge1$, define
$\vartheta_\tau(t):=\operatorname{sgn}(t)\min\{1,\tau(|t|-\tau)_+\}$,
$\varpi:=\E[|\vartheta_\tau(G)|]$ and $u_0:=1/(10(\tau+1)^2)$, and let
$c_\phi:=e^{-3/2}/2$. Then
$c_\phi\phi(\tau)/\tau\le\varpi\le2\phi(\tau)/\tau$,
$\tau\varpi\le\E[G\vartheta_\tau(G)]\le(2/c_\phi)\tau\varpi$ and
$\E[\vartheta_\tau(G)^2]\le\varpi$. For $(G_1,G_2)$ as in~(a),
$\E[|\vartheta_\tau(G_1)\vartheta_\tau(G_2)|]\le\varpi$, and if
$|u|\le u_0$, then, with an absolute constant $C_\phi$,
\[
 |\E[\vartheta_\tau(G_1)\vartheta_\tau(G_2)]|\le C_\phi|u|\tau^2\varpi^2,
 \qquad
 \E[|\vartheta_\tau(G_1)\vartheta_\tau(G_2)|]\le C_\phi\varpi^2 .
\]
\item Let $d\ge3$ and let $\theta$ be the angle between two independent
uniform points of a centered sphere in $\R^d$. Then $\theta$ has density
$c_d\sin^{d-2}\theta$ on $[0,\pi]$ with $c_d\le\sqrt d$, so a coordinate
of a uniform unit vector has density $c_d(1-t^2)^{(d-3)/2}$ on
$[-1,1]$. Moreover $\Pr[\sin\theta\le u]\le4u^{d-1}$ for $0<u\le1/2$,
$\E[1/\sin\theta]\le6$, and
$\E[e^{-\sin^2\theta/w^2}]\le e^{-1/(4w^2)}+4w^{d-1}\Gamma((d+1)/2)$
for every $w>0$.
\end{enumerate}
\end{lemma}

\begin{proof}
(a) This is Mehler's formula, and we include the short
proof. For $s,t\in\R$,
$\E[\exp(sG_1-s^2/2+tG_2-t^2/2)]=e^{stu}$. The generating function
$e^{sx-s^2/2}=\sum_ns^nHe_n(x)/n!$ converges in $L_2(\gamma_1)$, its
terms having norms $|s|^n/\sqrt{n!}$. The product of the two series
therefore converges in $L_1$, and the double series
$\sum_{n,l}s^nt^l\E[He_n(G_1)He_l(G_2)]/(n!\,l!)$ converges absolutely
to $e^{stu}$. Comparing coefficients gives
$\E[h_n(G_1)h_l(G_2)]=\delta_{nl}u^n$. Since
$|\E[f(G_1)g(G_2)]|\le\|f\|_{L_2(\gamma_1)}\|g\|_{L_2(\gamma_1)}$ by
Cauchy--Schwarz, the bilinear form is continuous, and expanding $f$
and $g$ in the basis gives~(a).

(b) Let $k_a:=\dim E_a$. Choose orthonormal bases
$e_1,\ldots,e_{k_1}$ of $E_1$ and $f_1,\ldots,f_{k_2}$ of $E_2$ with
$\langle e_a,f_b\rangle=\delta_{ab}\lambda_a$, where
$\lambda_1\ge\lambda_2\ge\cdots\ge0$ are the singular values of the
matrix $(\langle e'_a,f'_b\rangle)$ for arbitrary orthonormal bases.
Rotating both bases by the factors of a singular value decomposition
gives such bases. The matrix represents $P_1$ restricted to $E_2$, so
every $\lambda_a\le\lambda$. Define $\alpha_a:=\langle g,e_a\rangle$ and
$\beta_b:=\langle g,f_b\rangle$. The vector $(\alpha,\beta)$ is
Gaussian with $\alpha$ and $\beta$ standard and
$\operatorname{Cov}(\alpha_a,\beta_b)=\delta_{ab}\lambda_a$, so the
pairs $(\alpha_a,\beta_a)$, $a\le\min\{k_1,k_2\}$, and the remaining
coordinates are independent. Write $W_1(g)=\sum_a\Phi_a(\alpha)e_a$
and $W_2(g)=\sum_b\Psi_b(\beta)f_b$. The functions $\Phi_a,\Psi_b$ are
odd, hence centered, and
$\E[\langle W_1,W_2\rangle]=\sum_a\lambda_a\E[\Phi_a(\alpha)\Psi_a(\beta)]$.

We claim that $|\E[\Phi(\alpha)\Psi(\beta)]|\le\lambda\|\Phi\|_2\|\Psi\|_2$
whenever $\E[\Phi(\alpha)]=0$. The products
$H_{\mathbf n}(\alpha):=\prod_ah_{n_a}(\alpha_a)$ form an orthonormal
basis for functions of $\alpha$, and similarly for $\beta$. By
independence and~(a),
$\E[H_{\mathbf n}(\alpha)H_{\mathbf n'}(\beta)]=\prod_a\lambda_a^{n_a}$
if $\mathbf n=\mathbf n'$ and both vanish outside the first
$\min\{k_1,k_2\}$ coordinates, and it is zero otherwise. Expanding
$\Phi$ and $\Psi$, using continuity of the bilinear form as in~(a) and
the vanishing of the constant coefficient of $\Phi$, we get
$|\E[\Phi\Psi]|\le\sum_{\mathbf n\ne0}\lambda^{|\mathbf n|}
|\widehat\Phi_{\mathbf n}\widehat\Psi_{\mathbf n}|
\le\lambda\|\Phi\|_2\|\Psi\|_2$. Applying the claim to each $a$ and
Cauchy--Schwarz over $a$,
\[
 |\E[\langle W_1,W_2\rangle]|
 \le\lambda^2\sum_a\|\Phi_a\|_2\|\Psi_a\|_2
 \le\lambda^2(\sum_a\|\Phi_a\|_2^2)^{1/2}
           (\sum_a\|\Psi_a\|_2^2)^{1/2},
\]
which is~(b).

(c) Since $\vartheta_\tau$ vanishes on $[-\tau,\tau]$ and equals
$\pm1$ for $|t|\ge\tau':=\tau+1/\tau$,
$\Pr[|G|\ge\tau']\le\varpi\le\Pr[|G|\ge\tau]$. The Gaussian tail
bounds $(1/t-1/t^3)\phi(t)\le\Pr[G>t]\le\phi(t)/t$ for $t>0$
(Lemma~\ref{lem:pre_gaussian_tail}), together with $2\le\tau'\le2\tau$ and
$\phi(\tau')\ge e^{-3/2}\phi(\tau)$, give the first two inequalities:
$\Pr[|G|\ge\tau']\ge\frac32\phi(\tau')/\tau'\ge c_\phi\phi(\tau)/\tau$. As $t\vartheta_\tau(t)\ge\tau|\vartheta_\tau(t)|$,
$\E[G\vartheta_\tau(G)]\ge\tau\varpi$, and
$\E[G\vartheta_\tau(G)]\le\E[|G|\mathbf 1_{|G|\ge\tau}]=2\phi(\tau)
\le(2/c_\phi)\tau\varpi$. Also $\vartheta_\tau^2\le|\vartheta_\tau|$, and
Cauchy--Schwarz gives $\E[|\vartheta_\tau(G_1)\vartheta_\tau(G_2)|]
\le\E[\vartheta_\tau(G)^2]\le\varpi$.

Now let $|u|\le u_0$, and write $\phi_u$ for the density of
$(G_1,G_2)$. As $x^2-2uxy+y^2\ge(1-|u|)(x^2+y^2)$,
$\phi_u(x,y)\le(2\pi)^{-1}(1-u^2)^{-1/2}
\exp(-(x^2+y^2)/(2(1+|u|)))$. Substituting
$x=\sqrt{1+|u|}\,x'$ and $y=\sqrt{1+|u|}\,y'$,
\[
 \E[|\vartheta_\tau(G_1)\vartheta_\tau(G_2)|]
 \le\Pr[|G_1|\ge\tau,|G_2|\ge\tau]
 \le2\Pr[|G|\ge\tau(1+|u|)^{-1/2}]^2
 \le\frac{8\cdot2\phi(\tau)^2}{\tau^2}e^{\tau^2|u|},
\]
and $\tau^2|u|\le1/10$. This is at most $C\varpi^2$. For the first
estimate define $S(u):=\E[\vartheta_\tau(G_1)\vartheta_\tau(G_2)]
=\sum_{n\ {\rm odd}}\widehat{(\vartheta_\tau)}_n^2u^n$ by~(a), since
$\vartheta_\tau$ is odd. The function $\vartheta_\tau$ is Lipschitz
with a.e.\ derivative $\vartheta_\tau'=\tau\mathbf 1_{\tau<|t|<\tau'}$.
Gaussian integration by parts, $\E[f'(G)q(G)]=\E[f(G)(Gq(G)-q'(G))]$
for polynomials $q$, and the recurrence
$He_n=xHe_{n-1}-He_{n-1}'$ give
$\widehat{(\vartheta'_\tau)}_{n-1}=\sqrt n\,\widehat{(\vartheta_\tau)}_n$.
Hence~(a) applied to $\vartheta_\tau'$ gives, for $|u|<1$,
\[
 S'(u)=\sum_nn\widehat{(\vartheta_\tau)}_n^2u^{n-1}
 =\E[\vartheta_\tau'(G_1)\vartheta_\tau'(G_2)]
 =\tau^2\Pr_u[\tau<|G_1|<\tau',\ \tau<|G_2|<\tau'].
\]
On the square $\{|x|,|y|\le\tau+1\}$ the ratio
$\phi_u(x,y)/(\phi(x)\phi(y))$ is at most
$(1-u^2)^{-1/2}\exp(4|u|(\tau+1)^2/(2(1-u^2)))\le2$ for $|u|\le u_0$.
Thus $0\le S'(u')\le2\tau^2(2\phi(\tau)/\tau)^2\le(8/c_\phi^2)\tau^2\varpi^2$
for $|u'|\le u_0$. Since $S(0)=0$, the mean value theorem gives the
first estimate.

(d) The angle has density $c_d\sin^{d-2}\theta$ on $[0,\pi]$, where
$c_d^{-1}:=\int_0^\pi\sin^{d-2}\theta\,\d\theta\ge
\int_{|\theta-\pi/2|\le d^{-1/2}}(1-(\theta-\pi/2)^2/2)^{d-2}\,\d\theta
\ge2d^{-1/2}(1-1/(2d))^{d-2}\ge d^{-1/2}$. For $0<u\le1/2$,
substituting $s=\sin\theta$ on both halves of $[0,\pi]$,
\[
 \Pr[\sin\theta\le u]\le2c_d\int_0^u\frac{s^{d-2}\,\d s}{\sqrt{1-s^2}}
 \le\frac{4\sqrt d}{\sqrt3(d-1)}u^{d-1}\le4u^{d-1}.
\]
Then $\E[1/\sin\theta]\le2+\int_2^\infty\Pr[\sin\theta<1/s]\,\d s
\le2+4\int_2^\infty s^{1-d}\,\d s\le6$. Finally, for
$f(u):=e^{-u^2/w^2}$,
$\E[f(\sin\theta)]=\int_0^\infty|f'(u)|\Pr[\sin\theta<u]\,\d u$.
The part $u>1/2$ contributes at most $f(1/2)$, and the part
$u\le1/2$ at most $\int_0^\infty4u^{d-1}(2u/w^2)e^{-u^2/w^2}\,\d u
=4w^{d-1}\Gamma((d+1)/2)$.
\end{proof}

\begin{lemma}[Empirical halfspace probabilities]\label{lem:s-halfspace}
For every $a'>0$ there is $d_H$ such that the following holds for
$d\ge d_H$ and $\Xi\ge e^{a'd}$. Let $Z_1,\ldots,Z_\Xi$ be independent
standard Gaussian vectors in $\R^d$. With probability at least
$1-2(\sqrt2/e)^d$,
\begin{enumerate}[label=(\roman*),nosep,leftmargin=*]
\item $\|Z_t\|_2\le2(d+\log\Xi)^{1/2}$ for every $t$, and
\item for every unit $v\in\R^d$ and every interval $J'\subseteq\R$,
open, closed or half-open and possibly unbounded,
$|\Xi^{-1}\#\{t:\langle Z_t,v\rangle\in J'\}-\Pr[G\in J']|
\le2\Xi^{-1/5}$.
\end{enumerate}
\end{lemma}

\begin{proof}
Define $H_0:=2(d+\log\Xi)^{1/2}$. Since $\E[e^{\|Z\|_2^2/4}]=2^{d/2}$,
Markov's inequality (Lemma~\ref{lem:pre_markov}) gives
$\Pr[\|Z_t\|_2>H_0]\le2^{d/2}e^{-d-\log\Xi}$, and a union bound over
$t$ gives~(i) outside probability $(\sqrt2/e)^d$.

Define $\rho_1:=\Xi^{-1/4}$. Let $\mathcal N_1$ be a $\rho_1$-net of the
unit sphere with $|\mathcal N_1|\le(3/\rho_1)^d$
(Lemma~\ref{lem:pre_nets}), and let
$\mathcal S_1:=\rho_1\mathbb Z\cap[-2H_0,2H_0]$. For each
$(v',s')\in\mathcal N_1\times\mathcal S_1$, Hoeffding's inequality
(Lemma~\ref{lem:pre_hoeffding}) bounds the probability that
$|\Xi^{-1}\#\{t:\langle Z_t,v'\rangle\le s'\}-\Pr[G\le s']|>\rho_1$ by
$2e^{-2\Xi\rho_1^2}=2e^{-2\sqrt\Xi}$. The logarithm of the number of
pairs is at most $d\log(3\Xi^{1/4})+\log(5H_0\Xi^{1/4})$, which is
$o(\sqrt\Xi)$ uniformly for $\Xi\ge e^{a'd}$ as $d\to\infty$. Hence
for $d\ge d_H$ all these estimates hold outside probability
$e^{-\sqrt\Xi}\le(\sqrt2/e)^d$.

On both events, fix a unit $v$ and $s\in\R$, and choose
$v'\in\mathcal N_1$ with $\|v-v'\|_2\le\rho_1$. If $s>H_0$ the
empirical frequency of $\{\langle Z,v\rangle\le s\}$ is one and
$\Pr[G\le s]\ge1-e^{-H_0^2/2}$ by the Chernoff bound
$\Pr[G>t]\le e^{-t^2/2}$ of Lemma~\ref{lem:pre_gaussian_chernoff}, and similarly if $s<-H_0$. Otherwise
choose $s_\pm\in\mathcal S_1$ with
$s+H_0\rho_1\le s_+\le s+(H_0+1)\rho_1$ and
$s-(H_0+1)\rho_1\le s_-\le s-H_0\rho_1$. As
$|\langle Z_t,v-v'\rangle|\le H_0\rho_1$,
\[
 \mathbf 1_{\langle Z_t,v'\rangle\le s_-}
 \le\mathbf 1_{\langle Z_t,v\rangle\le s}
 \le\mathbf 1_{\langle Z_t,v'\rangle\le s_+},
\]
and the Gaussian density is at most one, so the empirical frequency of
$\{\langle Z,v\rangle\le s\}$ differs from $\Pr[G\le s]$ by at most
$\rho_1+(H_0+1)\rho_1\le\Xi^{-1/5}$ for $d\ge d_H$. The same bound for
strict inequalities follows by letting $s$ increase to its limit. An
interval is the difference of two such halfline events, which
gives~(ii).
\end{proof}

\begin{lemma}[Scalar estimates for shifted tests]\label{lem:s-scalar}
There are absolute constants $C_0\ge4$, $c_3>0$ and
$C_3,C_4,C_5,C_6$ with the following properties. Let $0<w\le1/8$, let
$s_1,\ldots,s_J\in(0,1]$, define
$h_j(z):=\frac12(\chi((z-s_j)/w)-\chi((z+s_j)/w))$, and let
$Y:=\{\pm s_j\}$. Let $(Z,W)$ be a centered Gaussian pair with unit
variances and correlation $\rho\in(-1,1)$.
\begin{enumerate}[label=(\alph*),nosep,leftmargin=*]
\item $\E[h_j(Z)^2]\le C_3w$, and if $2s_j\ge C_0w$, then
$\E[h_j(Z)^2]\ge c_3w$.
\item If the $2J$ points of $Y$ are pairwise at distance at least $w$,
then $\sum_j|h_j(z)|\le C_4$ for every $z$.
\item If the points of $Y$ are pairwise at distance at least $C_0w$,
then for every $j$
\[
 \E[|h_j(Z)|\sum_{l\ne j}|h_l(Z)|]\le\frac{C_5w}{C_0^2}\le\frac{c_3w}{16}.
\]
\item $|\E[h_j(Z)h_l(W)]|\le C_6w\exp(-(1-\rho^2)/(2w^2))$. If all
$s_j$ lie in one translate of $w\mathbb Z$, then for every $l$,
$\sum_j|\E[h_j(Z)h_l(W)]|^2\le C_6^2w^2\exp(-(1-\rho^2)/w^2)$.
\item $\E[h_j(Z)^2h_l(W)^2]\le C_3^2w^2(1-\rho^2)^{-1/2}$.
\end{enumerate}
\end{lemma}

\begin{proof}
Define $\chi_y(z):=\chi((z-y)/w)$ and $C_2:=\sup_y(1+|y|)^2|\chi(y)|$.
(a) Since $h_j^2\le(\chi_{s_j}^2+\chi_{-s_j}^2)/2$,
$\int h_j^2\le w\|\chi\|_2^2$, and the Gaussian density is at most
$(2\pi)^{-1/2}$. Define $C_3:=\|\chi\|_2^2/\sqrt{2\pi}$. For the lower
bound, $\chi(y)=\int\varrho(\xi)\cos(\xi y)\,\d\xi\ge1-2y^2\ge7/8$ for
$|y|\le1/4$. For $|z-s_j|\le w/4$ we have
$(z+s_j)/w\ge2s_j/w-1/4\ge C_0-1/4$, so
$|\chi_{-s_j}(z)|\le C_2C_0^{-2}\le1/4$ once $C_0\ge2\sqrt{C_2}+1$. Thus
$h_j\ge1/4$ on an interval of length $w/2$ inside $[-1/32,33/32]$,
where the Gaussian density is at least $\phi(2)$. Take
$c_3:=\phi(2)/32$.

(b) A half-open interval of length $w$ contains at most one point of
$Y$, so
\[
 \sum_j|h_j(z)|\le\frac12\sum_{y\in Y}\frac{C_2}{(1+|z-y|/w)^2}
 \le C_2\sum_{n\ge0}\frac1{(1+n)^2}\le2C_2=:C_4 .
\]

(c) For $y\ne y'$ in $Y$ and $\Delta:=(y'-y)/w$,
\[
 \E[|\chi_y(Z)\chi_{y'}(Z)|]\le\frac{w}{\sqrt{2\pi}}
 \int|\chi(t)\chi(t-\Delta)|\,\d t
 \le\frac{4C_2^2w}{\sqrt{2\pi}}(1+|\Delta|/2)^{-2},
\]
since $\max\{|t|,|t-\Delta|\}\ge|\Delta|/2$ and
$\int(1+|t|)^{-2}\,\d t=2$. For fixed $y$, the $n$-th nearest other
point of $Y$ on either side is at distance at least $nC_0w$. Summing
over the two points $y=\pm s_j$ and all $y'\in Y\setminus\{\pm s_j\}$
gives the bound with $C_5:=56C_2^2/\sqrt{2\pi}$. Enlarging $C_0$ we
may assume $C_5/C_0^2\le c_3/16$.

(d) For centers $y,y'$, Fubini and the characteristic function of
$(Z,W)$ give
\[
 \E[\chi_y(Z)\chi_{y'}(W)]=\int\!\!\int\varrho(\xi)\varrho(\zeta)
 e^{-\mathrm i(\xi y+\zeta y')/w}
 \exp(-\frac{(\zeta+\rho\xi)^2+(1-\rho^2)\xi^2}{2w^2})
 \,\d\zeta\,\d\xi .
\]
The inner integral is at most $\|\varrho\|_\infty\sqrt{2\pi}\,w$ in
absolute value, and $|\xi|\ge1$ on the support of $\varrho$. Writing
$h_j(Z)h_l(W)$ as a quarter of a signed sum of four such products
gives the first bound for every $C_6\ge\sqrt{2\pi}\|\varrho\|_\infty$. We take
$C_6:=4\pi(1+\|\varrho\|_\infty)^2$. For the second, fix $y'$ and let
$\Phi(\xi)$ be $\varrho(\xi)$ times the inner integral, a function
supported in $[-2,2]\subseteq(-\pi,\pi)$ with
$|\Phi|\le\sqrt{2\pi}\|\varrho\|_\infty^2we^{-(1-\rho^2)/(2w^2)}$.
For $y=y_0+nw$, $n\in\mathbb Z$, the expectation above is
$\int_{-\pi}^\pi\Phi(\xi)e^{-\mathrm i\xi y_0/w}e^{-\mathrm in\xi}\,\d\xi$,
a Fourier coefficient. Lemma~\ref{lem:pre_parseval} gives
$\sum_n|\E[\chi_{y_0+nw}(Z)\chi_{y'}(W)]|^2\le2\pi\int|\Phi|^2
\le16\pi^2\|\varrho\|_\infty^4w^2e^{-(1-\rho^2)/w^2}$. The points
$s_j$ and $-s_j$ lie in two translates of $w\mathbb Z$, so summing the
four products and using $(\sum_{1}^4a_m)^2\le4\sum a_m^2$ gives the
second bound.

(e) The density of $(Z,W)$ is at most $(2\pi)^{-1}(1-\rho^2)^{-1/2}$,
and $\int h_j^2\int h_l^2\le w^2\|\chi\|_2^4$.
\end{proof}

\subsection{The Adaptive-Shift Construction}
We now fix the dimension $d$ of the transverse queries, a dictionary
size $\Xi$, a width $w$ and centers $s_1,\ldots,s_J$. The numbers $N$
of tags, $H$ of transverse copies, $J$ of centers and $\Xi$ of
dictionary rows are local to the lemmas below. In particular $H$ is
unrelated to the $H$ of Section~\ref{sec:lp_above_two}. Define
\[
 h_j(z):=\tfrac12(\chi((z-s_j)/w)-\chi((z+s_j)/w)),\qquad
 f_j:=w^ph_j,\qquad \sigma:=w^{p+1/2},\qquad\bar\sigma:=\sqrt\Xi\,\sigma .
\]
For rows $u_1,\ldots,u_\Xi\in\R^d$ and $X\in\R^d$ write
$\mathbf f_j(X):=(f_j(\langle u_t,X\rangle))_{t\in[\Xi]}\in\R^\Xi$ for
the test vectors, and $v_j(X):=\mathbf f_j(X)/\|\mathbf f_j(X)\|_2$
when $\mathbf f_j(X)\ne0$. Let $\mathbb S_d:=\sqrt d\,S^{d-1}$. For a
probability law $\mathcal Q$ on $\R^d$ and an integer $N\ge1$, let
$\mathcal B=\mathcal B(\mathcal Q)\subseteq\R^{N\times\Xi}$ be the
closure of the set of arrays $\mathbf B$ whose $i$-th row is
$\E_X[\sum_jc_{ij}(X)\mathbf f_j(X)]$, $X\sim\mathcal Q$, over measurable
coefficients with
\begin{equation}\label{lem:s-eq-admissible}
 \E_X[\max_{i\in[N]}\sum_{j\in[J]}|c_{ij}(X)|]\le1 .
\end{equation}
This set is convex and symmetric. When the $\mathbf f_j$ are bounded
it is compact.

\begin{lemma}[A dictionary of shifted tests]\label{lem:s-frame}
Fix $p>0$ and $0<\lambda_1<\lambda_2$, and let $C_0$ be the constant
of Lemma~\ref{lem:s-scalar}. There are $c_0>0$, depending only on
$\lambda_2$, and $c_F\in(0,1]$ and $C_F\ge1$, depending only on $p$,
$\lambda_1$ and $\lambda_2$, with the following property. Let
$0<c_0'\le c_0$. Then there is $d_F$, depending only on $p$,
$\lambda_1$, $\lambda_2$ and $c_0'$, such that the following holds
for integers $d\ge d_F$ and $\Xi$ with
$\lambda_1d\le\log\Xi\le\lambda_2d$. Let $w:=c_0'/\sqrt d$, let
$2C_0w\le b\le1$, and let $s_1,\ldots,s_J\in[b/2,b]$ be pairwise at
distance at least $C_0w$. Let $\mathcal E$ be any event for
independent $N(0,I_d/d)$ vectors $u_1,\ldots,u_\Xi$ with
$\Pr[\mathcal E]\ge9/10$. Then there are $u_1,\ldots,u_\Xi$ in
$\mathcal E$ and a probability law $\mathcal Q$ on $\mathbb S_d$ with
the following properties.
\begin{enumerate}[label=(F\arabic*),nosep,leftmargin=*]
\item The $u_t$ are distinct, the vectors $(1,u_t)$ span $\R^{d+1}$, and
$\|u_t\|_2\le A_0:=2(1+\lambda_2)^{1/2}$ for every $t$.
\item For $\mathcal Q$-a.e.\ $X$,
$\Xi^{-1}\sum_t|\langle u_t,X\rangle|^p\le C_F$.
\item For $\mathcal Q$-a.e.\ $X$ and every $j$,
$c_F\bar\sigma\le\|\mathbf f_j(X)\|_2\le C_F\bar\sigma$.
\item For $\mathcal Q$-a.e.\ $X$, the Gram matrix of
$v_1(X),\ldots,v_J(X)$ has all its eigenvalues in $[1/2,3/2]$.
\item Let $X,Y\sim\mathcal Q$ and $j,l$ uniform in $[J]$ be
independent. Then $U:=\langle v_j(X),v_l(Y)\rangle$ satisfies
$\E[U^2]\le C_F/\Xi$ and $\E[U^4]\le C_F/\Xi^2$.
\end{enumerate}
\end{lemma}

\begin{proof}
Define $u_*:=\min\{1/2,e^{-5\lambda_2}\}$, $\rho_0:=(1-u_*^2)^{1/2}$ and
$c_0:=\min\{u_*(10\lambda_2)^{-1/2},1/8\}$, so $w\le1/8$. The points
of $Y:=\{\pm s_j\}$ are pairwise at distance at least $C_0w$, since
$s_j+s_l\ge b\ge2C_0w$. Also $J\le1+b/(2C_0w)\le\sqrt d/c_0'$. Let
$u_1,\ldots,u_\Xi$ be independent $N(0,I_d/d)$ vectors, and let
$\sigma_d$ be the uniform probability on $\mathbb S_d$. For fixed
$X\in\mathbb S_d$ the scores $Z_t:=\langle u_t,X\rangle$ are
independent standard Gaussians. Define $m_p:=\E[|G|^p]$.

\emph{Good queries.} Call $X\in\mathbb S_d$ good for the dictionary if
for every $j$
\[
 \frac1\Xi\sum_th_j(Z_t)^2\in[\frac{c_3w}2,2C_3w],\qquad
 \frac1\Xi\sum_t|h_j(Z_t)|\sum_{l\ne j}|h_l(Z_t)|\le\frac{2C_5w}{C_0^2},
\]
and $\Xi^{-1}\sum_t|Z_t|^p\le2m_p$. Fix $X$. The summands in the first
two conditions lie in $[0,1]$ and $[0,C_4]$, and by
Lemma~\ref{lem:s-scalar}(a)--(c) their means lie in $[c_3w,C_3w]$ and
are at most $C_5w/C_0^2$, while their variances are at most $C_3w$
and $C_4C_5w/C_0^2$. Bernstein's inequality (Lemma~\ref{lem:pre_bernstein})
bounds the failure probability of each condition by
$2\exp(-c_7\Xi w)$, with $c_7>0$ absolute. Chebyshev's inequality (Lemma~\ref{lem:pre_markov})
bounds that of the last condition by $m_{2p}/(m_p^2\Xi)$. Since
$\Xi w\ge c_0'e^{\lambda_1d}/\sqrt d$ and $J\le\sqrt d/c_0'$, the
probability $q_0$ that a fixed $X$ is not good tends to zero as
$d\to\infty$. By rotation invariance it does not depend on $X$.

On good $X$, $\|\mathbf f_j(X)\|_2^2=w^{2p}\sum_th_j(Z_t)^2$ lies
between $c_3\bar\sigma^2/2$ and $2C_3\bar\sigma^2$. The off-diagonal
row sums of the Gram matrix of the $v_j(X)$ are at most
$\sum_t|h_j(Z_t)|\sum_{l\ne j}|h_l(Z_t)|/\min_l\sum_th_l(Z_t)^2
\le4C_5/(c_3C_0^2)\le1/4$, by Lemma~\ref{lem:s-scalar}(c). Gershgorin's
theorem (Lemma~\ref{lem:pre_gershgorin}) places its eigenvalues in
$[3/4,5/4]$.

\emph{Correlations.} Let $X,Y$ be independent with law $\sigma_d$,
independent of the dictionary, let $\theta$ be their angle and
$\rho:=\cos\theta$, and let $j,l$ be independent and uniform. Define
$S:=\sum_tB_t$ with $B_t:=h_j(\langle u_t,X\rangle)h_l(\langle
u_t,Y\rangle)$. Given $X,Y,j,l$, the $B_t$ are independent copies of
$h_j(Z)h_l(W)$, where $(Z,W)$ has correlation $\rho$. By
Lemma~\ref{lem:s-gauss}(d), $\Pr[|\rho|>\rho_0]=\Pr[\sin\theta<u_*]
\le4u_*^{d-1}\le\Xi^{-4}$ for large $d$, as $\Xi\le e^{\lambda_2d}$.
If $|\rho|\le\rho_0$, Lemma~\ref{lem:s-scalar}(d),(e) give
$|\E[B_t]|\le C_6w\exp(-u_*^2d/(2c_0'^2))\le\Xi^{-5}$, by the choice of
$c_0$ and for large $d$, and
$\E[B_t^2]+\E[B_t^4]\le2C_3^2w^2/u_*$, since $|h_j|\le1$. For the
centered sum $S-\Xi\E[B_t]$ we get second moment at most
$2C_3^2\Xi w^2/u_*$ and fourth moment at most
$C\Xi w^2+C\Xi^2w^4\le2C\Xi^2w^4$, with $C$ depending on $\lambda_2$,
because $\Xi w^2\ge1$ for large $d$. Hence $\E[S^2]\le C'\Xi w^2$ and
$\E[S^4]\le C'\Xi^2w^4$ on $\{|\rho|\le\rho_0\}$. If $X$ and $Y$ are
good, then $|\langle v_j(X),v_l(Y)\rangle|\le2|S|/(c_3\Xi w)$. Let
$\mathbf 1_{XY}$ indicate that both are good. Using $|U|\le1$ on
$\{|\rho|>\rho_0\}$, the expectations over the dictionary, the queries
and the indices satisfy
\[
 \E[U^2\mathbf 1_{XY}]\le\Xi^{-4}+\frac{4C'}{c_3^2\Xi}\le\frac{C''}\Xi,
 \qquad
 \E[U^4\mathbf 1_{XY}]\le\Xi^{-4}+\frac{16C'}{c_3^4\Xi^2}
   \le\frac{C''}{\Xi^2}.
\]

\emph{Selection.} Let $q_1$ be the $\sigma_d$-measure of the queries
that are not good, and let $q_2,q_4$ be the two expectations above for
a fixed dictionary. Then $\E[q_1]=q_0$, $\E[q_2]\le C''/\Xi$ and
$\E[q_4]\le C''/\Xi^2$. By Markov's inequality (Lemma~\ref{lem:pre_markov}),
$\Pr[q_1>1/2]+\Pr[q_2>40C''/\Xi]+\Pr[q_4>40C''/\Xi^2]\le2q_0+1/20$.
By Lemma~\ref{lem:s-halfspace}(i) applied to $\sqrt du_t$,
$\|u_t\|_2\le2(1+\log\Xi/d)^{1/2}\le A_0$ for every $t$ outside
probability $2(\sqrt2/e)^d$. Almost surely the $u_t$ are distinct and,
as $\Xi\ge d+1$, the $(1,u_t)$ span $\R^{d+1}$. For large $d$ these
events and $\mathcal E$ hold simultaneously with positive probability.
Fix such a dictionary and let $\mathcal Q$ be $\sigma_d$ conditioned on
the good set, which has measure at least $1/2$. Conditioning
multiplies $q_2$ and $q_4$ by at most $4$. This gives (F1)--(F5) with
$c_F:=(c_3/2)^{1/2}$ and $C_F:=\max\{1,2m_p,(2C_3)^{1/2},160C''\}$.
\end{proof}

\begin{lemma}[Tests against original rows]\label{lem:s-cost}
Fix a non-even $p>0$. There is $C_p$, depending only on $p$, with the
following property. Let $N,H,d,\Xi\ge1$ be integers, let
$u_1,\ldots,u_\Xi\in\R^d$, and let $\mathcal Q$ be a probability law
on $\R^d$ with $\Xi^{-1}\sum_t|\langle u_t,X\rangle|^p\le C_E$ for
$\mathcal Q$-a.e.\ $X$. Let $0<w\le b\le1$ with $Hb^p\le1$, and let
$s_1,\ldots,s_J\in[0,b]$. Let $r\ge N+Hd$ and write
$\R^r=\R^N\oplus(\R^d)^H\oplus\R^{r-N-Hd}$, with $\iota_\ell$ the
inclusion of the $\ell$-th copy of $\R^d$. Consider the rows
$(e_i,\iota_\ell u_t,0)$ with masses $\alpha_{it}\in[1/2,3/2]$, for
$(\ell,i,t)\in[H]\times[N]\times[\Xi]$, together with the unit coordinate
rows of the last summand. Suppose that real weights on these rows
preserve every energy within relative error $\epsilon$, and let
$\widehat\alpha_{\ell it}$ be the output mass of row $(\ell,i,t)$, that is,
its weight times $\alpha_{it}$. Then for every $\ell\in[H]$ and every
$\mathbf B\in\mathcal B(\mathcal Q)$,
\begin{equation}\label{lem:s-eq-cost}
 \sum_{i,t}(\widehat\alpha_{\ell it}-\alpha_{it})\mathbf B_{it}
 \le C_p(1+C_E)\epsilon N\Xi .
\end{equation}
\end{lemma}

\begin{proof}
Let $\mathcal E(x)$ be the output energy minus the input energy at the
query $x$, so $|\mathcal E(x)|\le\epsilon F(x)$. Fix $\ell$, write
$z_{it}:=\widehat\alpha_{\ell it}-\alpha_{it}$, and let $\mathbf B$ have
coefficients $c_{ij}$ satisfying Eq.~\eqref{lem:s-eq-admissible}. The bound
then passes to the closure. Fix $X$ with
$C_X:=\max_i\sum_j|c_{ij}(X)|>0$, since queries with $C_X=0$
contribute nothing. At each tag $i$ draw independently an index
$J_i\in[J]\cup\{\emptyset\}$ with $\Pr[J_i=j]=(c_{ij}(X))_+/C_X$ and
$\Pr[J_i=\emptyset]=1-\sum_j(c_{ij}(X))_+/C_X\ge0$. For $\upsilon\in\R$
let $y(\upsilon)\in\R^N$ have coordinates $-s_{J_i}-w\upsilon$ if
$J_i\ne\emptyset$ and $0$ otherwise. Let $x_\pm(\upsilon)$ be the query
with tag coordinates $y(\upsilon)$, transverse coordinates $\pm X$ in the
$\ell$-th copy, and zero elsewhere. One variable $\upsilon$ is shared by all
tags. This matters because $\int\psi=0$, so independent variables at
different tags would average every term to zero.

In $\mathcal E(x_+)-\mathcal E(x_-)$ the rows of the other copies and
the coordinate rows cancel, because their inner products with
$x_\pm$ agree. A tag with $J_i=\emptyset$ contributes
$z_{it}(|\langle u_t,X\rangle|^p-|-\langle u_t,X\rangle|^p)=0$. For a tag
with $J_i=j$, Eq.~\eqref{lem:s-eq-filter} and the evenness of $\chi$ give
\[
 \frac12\int(|Z-s_j-w\upsilon|^p-|Z+s_j+w\upsilon|^p)\psi(\upsilon)\,\d\upsilon
 =f_j(Z)\qquad(Z\in\R).
\]
Hence
$\frac12\int(\mathcal E(x_+(\upsilon))-\mathcal E(x_-(\upsilon)))\psi(\upsilon)\,\d\upsilon
=\sum_{i:J_i\ne\emptyset}\sum_tz_{it}f_{J_i}(\langle u_t,X\rangle)$,
whose expectation over the $J_i$ is
$\sum_{i,t}z_{it}\sum_j((c_{ij}(X))_+/C_X)f_j(\langle u_t,X\rangle)$.

For the energy, every tag coordinate satisfies
$|y_i(\upsilon)|^p\le c_p'(b^p+w^p|\upsilon|^p)$ with
$c_p':=\max\{1,2^{p-1}\}$. Hence
\[
 F(x_\pm(\upsilon))\le\frac32[c_p'\sum_{i,t}
   (|y_i(\upsilon)|^p+|\langle u_t,X\rangle|^p)
   +(H-1)\Xi\sum_i|y_i(\upsilon)|^p]
 \le\frac32c_p'^2N\Xi\,(C_E+1+|\upsilon|^p),
\]
using $Hb^p\le1$ and $Hw^p\le Hb^p$. Therefore the absolute value of
the last expectation is at most
$\epsilon\frac32c_p'^2N\Xi\int(C_E+1+|\upsilon|^p)|\psi(\upsilon)|\,\d\upsilon$.
The negative parts are treated in the same way with
$(c_{ij}(X))_-$. Multiplying by $C_X$, integrating over $X$ and
using Eq.~\eqref{lem:s-eq-admissible} gives Eq.~\eqref{lem:s-eq-cost} with
$C_p:=3c_p'^2\int(1+|\upsilon|^p)|\psi(\upsilon)|\,\d\upsilon$.
\end{proof}

\begin{lemma}[Packing bounded masses]\label{lem:s-pack}
For every $\theta_0\in(0,1]$ there are $c_{\rm pk}\in(0,1/2]$ and
$n_0$, depending only on $\theta_0$, with the following property. Let
$n\ge n_0$, let $\mathcal P$ be an orthogonal projection of $\R^n$, and
define $g:=\mathcal Pg'$ for a standard Gaussian vector $g'$ in $\R^n$. Let
$\|\cdot\|_*$ be a seminorm on $\R^n$ with $\|z\|_*\le R\|z\|_2$ for
$z\in\operatorname{ran}\mathcal P$, and suppose
$\E[\|g\|_*]\ge\theta_0R\sqrt n$. If $0\le\epsilon_0\le c_{\rm pk}R\sqrt n$,
then there is $\alpha\in[1/2,3/2]^n$ such that every
$\widehat\alpha\in\R^n$ with $\|\widehat\alpha-\alpha\|_*\le\epsilon_0$
has more than $c_{\rm pk}n$ nonzero coordinates. For $\mathcal P=I$ the
vector $g$ is standard Gaussian.
\end{lemma}

\begin{proof}
By Gaussian concentration (Lemma~\ref{lem:pre_gaussian_concentration}) there is an
absolute $c_G>0$ such that every $R$-Lipschitz function $f$ on $\R^n$
satisfies $\Pr[|f(g')-\E[f(g')]|\ge t]\le2\exp(-c_Gt^2/R^2)$ for $t>0$. Choose $\varkappa\ge1$
with $1-2e^{-\varkappa^2/2}\ge\exp(-c_G\theta_0^2/16)$, then
$c_{\rm pk}\le\min\{1/2,\theta_0/(8\varkappa)\}$ with
$c_{\rm pk}(\log(e/c_{\rm pk})+\log(5+16\varkappa/\theta_0))
\le c_G\theta_0^2/16$, and $n_0\ge8/(c_G\theta_0^2)$. Define
$\mathcal U:=\theta_0R\sqrt n$ and $\ell:=\lfloor c_{\rm pk}n\rfloor$.
Suppose, for a contradiction, that every $\alpha\in[1/2,3/2]^n$ admits
$\widehat\alpha$ with at most $\ell$ nonzero coordinates and
$\|\widehat\alpha-\alpha\|_*\le\epsilon_0$.

For $T\subseteq[n]$ with $|T|=\ell$ let
$\mathcal A_T:=-2\varkappa\mathbf 1+\operatorname{span}\{e_t:t\in T\}$ and
$\mathcal Y_T:=\{y\in\mathcal A_T:\|y\|_*\le\varkappa R\sqrt n
+\mathcal U/4\}$. Let $\mathcal N_T\subseteq\mathcal Y_T$ be maximal
with pairwise $\|\cdot\|_*$-distances larger than $\mathcal U/4$. In
the quotient of $\operatorname{span}\{e_t:t\in T\}$ by the kernel of
$\|\cdot\|_*$, a normed space of dimension at most $\ell$, the balls of
radius $\mathcal U/8$ about the points of $\mathcal N_T$ are disjoint
and lie in a ball of radius $2\varkappa R\sqrt n+\mathcal U/2
+\mathcal U/8$. By Lemma~\ref{lem:pre_nets},
$|\mathcal N_T|\le(5+16\varkappa/\theta_0)^\ell$. Define
$\mathcal N:=\bigcup_T\mathcal N_T$, so that
$\log|\mathcal N|\le\ell\log(en/\ell)+\ell\log(5+16\varkappa/\theta_0)
\le c_G\theta_0^2n/16$.

Let $g\in\operatorname{ran}\mathcal P\cap Q_\varkappa$, where
$Q_\varkappa:=[-\varkappa,\varkappa]^n$, and
$\alpha:=\mathbf 1+g/(2\varkappa)\in[1/2,3/2]^n$. Take $\widehat\alpha$
as assumed, $T\supseteq\operatorname{supp}\widehat\alpha$ with
$|T|=\ell$, and $y:=2\varkappa(\widehat\alpha-\mathbf 1)\in\mathcal A_T$.
Then $\|g-y\|_*=2\varkappa\|\alpha-\widehat\alpha\|_*\le\mathcal U/4$ and
$\|y\|_*\le R\|g\|_2+\mathcal U/4\le\varkappa R\sqrt n+\mathcal U/4$, so
$y\in\mathcal Y_T$, and some $y'\in\mathcal N_T$ has
$\|g-y'\|_*\le\mathcal U/2$. Here the bound $\|g\|_*\le R\|g\|_2$ is used
only for $g\in\operatorname{ran}\mathcal P$, while $y$, $y'$ and
$\mathcal A_T$ need not lie in $\operatorname{ran}\mathcal P$. For fixed
$y'$, convexity and the symmetry of $g$ give
$\E[\|g-y'\|_*]\ge\E[\|g\|_*]\ge\mathcal U$. The function
$g'\mapsto\|\mathcal Pg'-y'\|_*$ is $R$-Lipschitz on $\R^n$, since
$\|\mathcal P(g'-g'')\|_*\le R\|\mathcal P(g'-g'')\|_2\le R\|g'-g''\|_2$.
Hence
\[
 \Pr[g\in Q_\varkappa]\le\sum_{y'\in\mathcal N}
   \Pr[\|g-y'\|_*\le\mathcal U/2]
 \le2|\mathcal N|e^{-c_G\theta_0^2n/4}
 <e^{-c_G\theta_0^2n/16},
\]
where the last step uses $n\ge n_0$. The coordinates
$g_i=\langle\mathcal Pe_i,g'\rangle$ form a centered Gaussian vector with
$\operatorname{Var}g_i=\|\mathcal Pe_i\|_2^2\le1$. \v{S}id\'ak's
inequality (Lemma~\ref{lem:pre_sidak}) and the Chernoff bound (Lemma~\ref{lem:pre_gaussian_chernoff}) give
$\Pr[g\in Q_\varkappa]\ge\prod_i\Pr[|g_i|\le\varkappa]
\ge(1-2e^{-\varkappa^2/2})^n\ge e^{-c_G\theta_0^2n/16}$, a contradiction.
\end{proof}

\begin{lemma}[Clipped adaptive tests]\label{lem:s-width}
For all $0<c_F\le1\le C_F$ there are $J_0\ge55$, depending only on
$C_F$, and $A_1,c_W>0$, depending only on $c_F$ and $C_F$, with the
following property. Let $N,\Xi,J$ be positive integers with $J\ge J_0$
and $\Xi\ge J^2$, let $\sigma>0$ and $\bar\sigma:=\sqrt\Xi\,\sigma$.
Let $\mathcal Q$ be a probability law on a measurable space
$\mathcal X$, and let $\mathbf f_1,\ldots,\mathbf f_J:\mathcal X\to\R^\Xi$
be measurable maps such that, for $\mathcal Q$-a.e.\ $X$,
$c_F\bar\sigma\le\|\mathbf f_j(X)\|_2\le C_F\bar\sigma$ for every $j$
and the Gram matrix of the unit vectors $v_j(X)$ has all its
eigenvalues in $[1/2,3/2]$. Suppose that, for independent
$X,Y\sim\mathcal Q$ and independent uniform $j,l\in[J]$, the correlation
$U:=\langle v_j(X),v_l(Y)\rangle$ satisfies $\E[U^2]\le C_F/\Xi$ and
$\E[U^4]\le C_F/\Xi^2$. Define
$\mathcal B$ by Eq.~\eqref{lem:s-eq-admissible}, define $\tau:=(\log(J)/4)^{1/2}$
and $R:=A_1\sqrt N\sigma\tau$, and let
\[
 \|z\|_*:=\sup\{\langle z,\mathbf B\rangle:\mathbf B\in\mathcal B,\
   \|\mathbf B\|_2\le R\}\qquad(z\in\R^{N\times\Xi}).
\]
Then $\|z\|_*\le R\|z\|_2$, and $\E[\|g\|_*]\ge c_WN\bar\sigma\tau$ for a
standard Gaussian array $g$.
\end{lemma}

\begin{proof}
The first assertion is Cauchy--Schwarz. Let $V_X$ be the
$\Xi\times J$ matrix with columns $v_j(X)$, $\Gamma_X:=V_X^\top V_X$,
and $P_X:=V_X\Gamma_X^{-1}V_X^\top$ the orthogonal projection onto the
column span. Let $\vartheta_\tau$, $\varpi$ and $c_\phi$ be as in
Lemma~\ref{lem:s-gauss}(c), and note $\tau\ge1$ as $J\ge55$. Let $c_G$ be the
Gaussian concentration constant from the proof of
Lemma~\ref{lem:s-pack}, and define
$u_J:=2\exp(-2c_GJ\varpi^2/(3\tau^2))$. By Lemma~\ref{lem:s-gauss}(c),
$\varpi\ge c_\phi\phi(\tau)/\tau\ge c_\phi'J^{-1/8}/\tau$ with
$c_\phi':=c_\phi/\sqrt{2\pi}$, so
$J\varpi^2/\tau^2\ge c_\phi'^2J^{3/4}/\tau^4$. Hence
$u_J/\varpi\to0$, $Ju_J/\varpi^2\to0$ and $\varpi J^2\to\infty$ as
$J\to\infty$. Choose $J_0\ge55$ so that $u_J\le\varpi/2$,
$6C_FJu_J/\varpi^2\le1$ and $\varpi J^2\ge1$ for $J\ge J_0$. Then
$\varpi\Xi\ge\varpi J^2\ge1$. Below, $v$ denotes $v_j(X)$ for
$X\sim\mathcal Q$ and $j$ uniform in $[J]$, and $v'$ an independent
copy.

\emph{The unclipped average.} For $g\in\R^\Xi$ and good $X$ define
\[
 U_X(g):=\frac{\bar\sigma}{J\varpi}\sum_j
   \vartheta_\tau(\langle g,v_j(X)\rangle)v_j(X),\qquad
 \mathbf A(g):=\E_X[U_X(g)]=\frac{\bar\sigma}\varpi
   \E_v[\vartheta_\tau(\langle g,v\rangle)v] .
\]
For unit $v,v'$ the pair $(\langle g,v\rangle,\langle g,v'\rangle)$ is
standard Gaussian with correlation $\langle v,v'\rangle$. By Fubini and
Lemma~\ref{lem:s-gauss}(a), with the odd coefficients
$\widehat\vartheta_n$ of $\vartheta_\tau$,
\[
 \E_g[\|\mathbf A(g)\|_2^2]
 =\frac{\bar\sigma^2}{\varpi^2}\E_{v,v'}[\langle v,v'\rangle
   \sum_{n\ {\rm odd}}\widehat\vartheta_n^2\langle v,v'\rangle^n]
 \le\frac{\bar\sigma^2}{\varpi^2}(\widehat\vartheta_1^2\E[U^2]
   +\textstyle\sum_n\widehat\vartheta_n^2\,\E[U^4]),
\]
as $|U|\le1$. Here $\widehat\vartheta_1=\E[G\vartheta_\tau(G)]\le
(2/c_\phi)\tau\varpi$ and
$\sum_n\widehat\vartheta_n^2=\E[\vartheta_\tau(G)^2]\le\varpi$. Hence
$\E[\|\mathbf A(g)\|_2^2]\le(4/c_\phi^2)C_F\sigma^2\tau^2
+C_F\sigma^2/(\varpi\Xi)\le C_7\sigma^2\tau^2$ with
$C_7:=5C_F/c_\phi^2$.

\emph{Clipping at each query.} Define
$\mathcal M_X(g):=\sum_j|\vartheta_\tau(\langle g,v_j(X)\rangle)|$,
whose mean is $J\varpi$, and
\[
 U'_X(g):=\min\{1,\frac{2J\varpi}{\mathcal M_X(g)}\}U_X(g),
\]
with multiplier one when $\mathcal M_X(g)=0$. Then
$U'_X(g)=\sum_jc_j(X,g)\mathbf f_j(X)$ with
\[
 \sum_j|c_j(X,g)|
 \le\min\{1,\frac{2J\varpi}{\mathcal M_X(g)}\}
   \frac{\bar\sigma\,\mathcal M_X(g)}{J\varpi\,c_F\bar\sigma}
 \le\frac2{c_F}
\]
for every $X$ and $g$. Since $|\vartheta_\tau'|\le\tau$ and
$\|V_X\|_{\rm op}\le(3/2)^{1/2}$, the map $g\mapsto\mathcal M_X(g)$ is
$\tau(3J/2)^{1/2}$-Lipschitz, so
$\Pr[\mathcal M_X>2J\varpi]\le u_J$ by Gaussian concentration. As
$t\vartheta_\tau(t)\ge\tau|\vartheta_\tau(t)|\ge0$,
\[
 \E_g[\langle g,U'_X(g)\rangle]
 \ge\frac{\bar\sigma\tau}{J\varpi}\E[\mathcal M_X
   \mathbf 1_{\mathcal M_X\le2J\varpi}]
 \ge\frac{\bar\sigma\tau}{J\varpi}(J\varpi-Ju_J)
 \ge\frac{\bar\sigma\tau}2,
\]
using $\mathcal M_X\le J$.

\emph{The clipped part is short on average.} Define
$W_X:=U'_X-U_X$ and $\mathbf W(g):=\E_X[W_X(g)]$. The map $W_X$ is odd,
takes values in the column span of $V_X$, depends on $g$ only through
$P_Xg$, and vanishes unless $\mathcal M_X>2J\varpi$. Since
$\|U_X(g)\|_2\le(\bar\sigma/(J\varpi))\|V_X\|_{\rm op}\sqrt J$, we have
$\E_g[\|W_X(g)\|_2^2]\le3\bar\sigma^2u_J/(2J\varpi^2)$. For good $X,Y$
let $\lambda_{XY}:=\|P_XP_Y\|_{\rm op}$. As $V_X\Gamma_X^{-1/2}$ has
orthonormal columns and $\|\Gamma_X^{-1/2}\|_{\rm op}\le\sqrt2$,
\[
 \lambda_{XY}=\|V_X\Gamma_X^{-1}(V_X^\top V_Y)\Gamma_Y^{-1}V_Y^\top\|_{\rm op}
 \le2\|V_X^\top V_Y\|_{\rm op},\qquad
 \lambda_{XY}^2\le4\sum_{j,l}\langle v_j(X),v_l(Y)\rangle^2 ,
\]
so $\E_{X,Y}[\lambda_{XY}^2]\le4J^2\E[U^2]\le4C_FJ^2/\Xi$.
Lemma~\ref{lem:s-gauss}(b) and Fubini give
\[
 \E_g[\|\mathbf W(g)\|_2^2]
 =\E_{X,Y}[\E_g[\langle W_X(g),W_Y(g)\rangle]]
 \le\E_{X,Y}[\lambda_{XY}^2]\cdot\frac{3\bar\sigma^2u_J}{2J\varpi^2}
 \le\frac{6C_FJu_J}{\varpi^2}\sigma^2\le\sigma^2 .
\]
Define $\mathbf A'(g):=\E_X[U'_X(g)]=\mathbf A(g)+\mathbf W(g)$. Then
$\E[\|\mathbf A'(g)\|_2^2]\le2(C_7+1)\sigma^2\tau^2$.

\emph{Tags.} Let $g_1,\ldots,g_N$ be the rows of $g$, and define
$\mathsf C_i:=\langle g_i,\mathbf A'(g_i)\rangle
=\E_X[\langle g_i,U'_X(g_i)\rangle]\ge0$. Then
$\E[\mathsf C_i]\ge\bar\sigma\tau/2$. Also
\[
 \langle g,U'_X(g)\rangle
 \le\min\{1,\frac{2J\varpi}{\mathcal M_X}\}
   \frac{\bar\sigma\mathcal M_X}{J\varpi}\max_j|\langle g,v_j(X)\rangle|
 \le2\bar\sigma\max_j|\langle g,v_j(X)\rangle| .
\]
By Jensen's inequality
and $\E[\max_jG_j^2]\le4\log\E[\sum_je^{G_j^2/4}]=4\log(\sqrt2J)$ for
standard Gaussians $G_j$,
$\E[\mathsf C_i^2]\le4\bar\sigma^2\E_X[\E_g[\max_j\langle g,v_j(X)\rangle^2]]
\le128\bar\sigma^2\tau^2$. Let $A_2:=4(256(C_7+1))^{1/2}$ and
$\widetilde{\mathbf A}(g_i):=\mathbf A'(g_i)
\mathbf 1_{\|\mathbf A'(g_i)\|_2\le A_2\sigma\tau}$. By Cauchy--Schwarz
and Markov's inequality (Lemma~\ref{lem:pre_markov}),
\[
 \E[\langle g_i,\widetilde{\mathbf A}(g_i)\rangle]
 \ge\frac{\bar\sigma\tau}2-(\E[\mathsf C_i^2])^{1/2}
   \Pr[\|\mathbf A'(g_i)\|_2>A_2\sigma\tau]^{1/2}
 \ge\frac{\bar\sigma\tau}2-\frac{(256(C_7+1))^{1/2}}{A_2}\bar\sigma\tau
 =\frac{\bar\sigma\tau}4 .
\]
The array with rows $\widetilde{\mathbf A}(g_i)$ has coefficients
$c_{ij}(X):=c_j(X,g_i)$ or zero, with $\sum_j|c_{ij}(X)|\le2/c_F$ for
every $X$ and $i$. Hence, divided by $2/c_F$, it satisfies
Eq.~\eqref{lem:s-eq-admissible} and lies in $\mathcal B$, and its Euclidean
norm is at most $(c_F/2)\sqrt NA_2\sigma\tau$ for every $g$. Define
$A_1:=c_FA_2/2$ and $c_W:=c_F/8$. Since $\|g\|_*\ge0$ and
$\|g\|_*\ge(c_F/2)\sum_i\langle g_i,\widetilde{\mathbf A}(g_i)\rangle$,
we get $\E[\|g\|_*]\ge c_WN\bar\sigma\tau$.
\end{proof}

\begin{lemma}[The adaptive-shift lower bound]\label{lem:s-shift}
Fix a non-even $p>0$, $0<\lambda_1<\lambda_2$ and $\beta_0>0$. There
are $c_0,c_*,c_S,c_S'>0$ and $d_0$, depending only on $p$,
$\lambda_1$, $\lambda_2$ and $\beta_0$, with the following property.
Let $N,H\ge1$ and $d\ge d_0$ be integers, let $\lambda_1d\le L\le
\lambda_2d$, and let $0<b\le1$ with $Hb^p\le1$ and $b\sqrt d\ge\beta_0$.
Define
\begin{gather*}
 w:=\frac{c_0}{\sqrt d},\qquad J:=\lfloor\frac b{2C_0w}\rfloor,
 \qquad s_j:=\frac b2+(j-1)C_0w\quad(j\in[J]),\\
 \tau:=(\frac{\log J}4)^{1/2},\qquad
 \Xi:=\lfloor c_*\epsilon^{-2}w^{2p+1}\tau^2\rfloor .
\end{gather*}
Then $\lambda_1d\le\log\Xi\le3\lambda_2d$, and
Lemma~\ref{lem:s-frame} applies with $\lambda_1$, $3\lambda_2$ and
$c_0'=c_0$ to these $w$, $b$ and $s_j$. Let $u_1,\ldots,u_\Xi$ and
$\mathcal Q$ be given by Lemma~\ref{lem:s-frame} for any admissible
event $\mathcal E$. Then there are masses
$\alpha\in[1/2,3/2]^{N\times\Xi}$ such that for every $r\ge N+Hd$ the
rank-$r$ input of Lemma~\ref{lem:s-cost} built from these data has the
following property: every reweighting of its rows by real weights
with relative error at most $\epsilon$ gives nonzero weight to more
than $c_SNH\Xi$ of the rows $(e_i,\iota_\ell u_t,0)$. Moreover
\[
 c_SNH\Xi\ge c_S'\epsilon^{-2}NHd^{-p-1/2}(1+\log_+(b\sqrt d)),
\]
so $S_p(r,\epsilon)$ exceeds the right side.
\end{lemma}

\begin{proof}
Let $c_0(3\lambda_2)$, $c_F$ and $C_F$ be the constants of
Lemma~\ref{lem:s-frame} for $p$, $\lambda_1$ and $3\lambda_2$, and let
$J_0$, $A_1$ and $c_W$ be those of Lemma~\ref{lem:s-width} for
$c_F,C_F$. Define $c_0:=\min\{c_0(3\lambda_2),\beta_0/(4C_0J_0)\}$,
$\theta_0:=c_W/A_1$, let $c_{\rm pk},n_0$ be as in
Lemma~\ref{lem:s-pack} for $\theta_0$, let $C_p$ be as in
Lemma~\ref{lem:s-cost}, and define
$c_*:=\min\{1,(c_{\rm pk}A_1/(C_p(1+C_F)))^2\}$. Finally $d_0$ is
large in terms of all of these.

\emph{Parameters.} Since $b/(2C_0w)=b\sqrt d/(2C_0c_0)\ge2J_0$, we
have $J\ge J_0\ge2$ and $J\ge b/(4C_0w)$. The centers lie in
$[b/2,b]$, are $C_0w$-separated, and $b\ge2C_0w$, so $w\le b$. Also
$J\le\sqrt d/(2C_0c_0)$. As $c_*\le1$ and $\tau^2\le\log d$, we have
$\log\Xi\le2L+\log\log d\le3\lambda_2d$. As $\epsilon^{-2}=e^{2L}/4$,
$\tau\ge1$ and
$w^{2p+1}=c_0^{2p+1}d^{-p-1/2}$, we have
$\log\Xi\ge2L-(p+\tfrac12)\log d+\log(c_*c_0^{2p+1}/8)\ge\lambda_1d$
for $d\ge d_0$. In particular $\Xi\ge J^2$ and $N\Xi\ge\Xi\ge n_0$.

\emph{The seminorm.} By (F3)--(F5), Lemma~\ref{lem:s-width} applies to
$\mathcal Q$ and the test vectors. It gives a seminorm $\|\cdot\|_*$ on
$\R^{N\times\Xi}$ with $\|z\|_*\le R\|z\|_2$, where
$R:=A_1\sqrt N\sigma\tau$, and
$\E[\|g\|_*]\ge c_WN\bar\sigma\tau=\theta_0R\sqrt{N\Xi}$. Its set
$\mathcal B$ is the set $\mathcal B(\mathcal Q)$ of
Lemma~\ref{lem:s-cost}.

\emph{Cost and packing.} Consider any valid reweighting, and for
$\ell\in[H]$ let $\widehat\alpha^{(\ell)}:=(\widehat\alpha_{\ell it})_{i,t}$. By
Lemma~\ref{lem:s-cost} with (F2),
$\|\widehat\alpha^{(\ell)}-\alpha\|_*\le\epsilon_0:=C_p(1+C_F)\epsilon N\Xi$.
The condition $\epsilon_0\le c_{\rm pk}R\sqrt{N\Xi}$ reads
$\Xi\le(c_{\rm pk}A_1/(C_p(1+C_F)))^2\epsilon^{-2}\sigma^2\tau^2$, which
holds by the choice of $c_*$ and $\Xi$. Lemma~\ref{lem:s-pack} gives
$\alpha\in[1/2,3/2]^{N\times\Xi}$ such that every array within
$\epsilon_0$ of $\alpha$ in $\|\cdot\|_*$ has more than $c_{\rm pk}N\Xi$
nonzero entries. This $\alpha$ is used in every copy. Each of the $H$
blocks $\widehat\alpha^{(\ell)}$ therefore has more than $c_{\rm pk}N\Xi$
nonzero entries, and a row with nonzero output mass has nonzero weight.
Take $c_S:=c_{\rm pk}$.

\emph{Size.} The quantity $c_*\epsilon^{-2}w^{2p+1}\tau^2$ is at least
$2$, so $\Xi$ is at least half of it. From $J\ge J_0\ge2$ and
$J\ge b\sqrt d/(4C_0c_0)$ we get
$\log J\ge c(1+\log_+(b\sqrt d))$ with $c>0$ depending on $C_0c_0$.
This gives the displayed bound. The input has rank $r$ by (F1), so
$S_p(r,\epsilon)$ exceeds the bound.
\end{proof}

\begin{lemma}[The adaptive-shift input has regular components]
\label{lem:s-regular}
In Lemma~\ref{lem:s-shift} take for $\mathcal E$ the event of
Lemma~\ref{lem:s-halfspace} for the vectors $\sqrt du_t$, with
$a'=\lambda_1$. For large $d$ it has probability at least $9/10$.
There are $B_1\ge3$ and $c_V>0$, depending only on $\lambda_1$ and
$\lambda_2$, such that for every choice of masses in $[1/2,3/2]$ the
resulting input has the following form. Its energy is
\[
 F(x)=\sum_gc_g\,\E[\theta_g(Z_g)|\langle m_g+Z_g,x\rangle|^p],
\]
with one term for each pair $(\ell,i)\in[H]\times[N]$ and one for each
coordinate row, so $NH+(r-N-Hd)\le2r^2$ terms. For a pair, $c_g:=\Xi$, $Z_g$ is
uniform on the $\Xi$ points $\iota_\ell(u_t-\bar u)$, where
$\bar u:=\Xi^{-1}\sum_tu_t$, $m_g:=(e_i,\iota_\ell\bar u,0)$ and
$\theta_g(\iota_\ell(u_t-\bar u)):=\alpha_{it}$. For a coordinate row
$e_k$, $c_g:=1$, $Z_g:=0$, $m_g:=e_k$ and $\theta_g:=1$. Each $Z_g$ is
centered with covariance $Q_g$, satisfies
$\E[e^{\langle v,Z_g\rangle}]\le e^{B_1^2v^\top Q_gv/2}$ for every $v$,
and satisfies
\[
 \Pr[0<|a+\langle v,Z_g\rangle|\le z\sqrt{v^\top Q_gv}]\le3z
\]
for every $v$ with $v^\top Q_gv>0$, every $a\in\R$ and every
$e^{-c_Vd}\le z\le1$. In particular, if $1\le p<2$, $B\ge B_1$,
$L^{-1/2}\ge e^{-c_Vd}$ and $2r^2\le e^{cL}$ for the constant $c$ of
Lemma~\ref{lem:joint-regular-components}, then the input belongs to
the class of that theorem.
\end{lemma}

\begin{proof}
The energy identity is a regrouping of the rows, and $\theta_g$ is well
defined because the $u_t$ are distinct (F1). Coordinate rows give
degenerate terms, for which both conditions are empty or trivial. Fix
a pair and write $\zeta:=\langle u,v\rangle$ for unit $v\in\R^d$, with
$u$ uniform on the dictionary. All $u_t$ lie in the ball of radius
$A_0':=2(1+3\lambda_2)^{1/2}$. Under $N(0,I_d/d)$, $\zeta$ has law
$N(0,1/d)$. By Lemma~\ref{lem:s-halfspace}(ii), applied to
$\sqrt d\zeta$, every interval has empirical and Gaussian
probabilities differing by at most $2\Xi^{-1/5}\le2e^{-\lambda_1d/5}$.

\emph{Mean and covariance.} By symmetry of the Gaussian law,
\[
 |\langle\bar u,v\rangle|=|\int_0^{A_0'}(\Pr_{\rm emp}[\zeta>s]
   -\Pr_{\rm emp}[\zeta<-s])\d s|\le4A_0'\Xi^{-1/5}.
\]
Similarly,
$\E_{\rm emp}[\zeta^2]=\int_0^{A_0'}2s\Pr_{\rm emp}[|\zeta|>s]\,\d s$ and
$1/d=\int_0^\infty2s\Pr_G[|\zeta|>s]\,\d s$ differ by at most
$2A_0'^2\Xi^{-1/5}+(4/d)e^{-dA_0'^2/2}\le1/(10d)$ for large $d$, by
the Chernoff bound $\Pr[|G|>t]\le2e^{-t^2/2}$ of
Lemma~\ref{lem:pre_gaussian_chernoff}. Hence
the empirical covariance $Q$ of $u$ satisfies
$0.8\|v\|_2^2/d\le v^\top Qv\le1.2\|v\|_2^2/d$, and $v^\top Q_gv
=v_\ell^\top Qv_\ell$, where $v_\ell$ is the $\ell$-th transverse block of $v$.

\emph{Interval condition.} Let $v_\ell\ne0$, define $\hat v:=v_\ell/\|v_\ell\|_2$,
and let $e^{-c_Vd}\le z\le1$ with $c_V:=\lambda_1/5$. The event in
question is contained in $\{\langle u,\hat v\rangle\in J'\}$ for an
interval $J'$ of length at most $2z(1.2/d)^{1/2}$. Its Gaussian
probability is at most $(d/(2\pi))^{1/2}\cdot2z(1.2/d)^{1/2}\le0.9z$,
so its empirical probability is at most $0.9z+2e^{-c_Vd}\le3z$.

\emph{Moment generating function.} Let $\zeta_c:=\zeta-\langle\bar u,v
\rangle$, so $|\zeta_c|\le2A_0'$. For $0\le s\le2A_0'\sqrt d$ and
$s_0:=4A_0'\Xi^{-1/5}\sqrt d\le1$,
\[
 \Pr_{\rm emp}[\sqrt d|\zeta_c|>s]\le\Pr[|G|>s-s_0]+2\Xi^{-1/5}
 \le2e^{1/2-s^2/4}+2e^{-\lambda_1s^2/(20A_0'^2)},
\]
using the same Chernoff bound, $(s-s_0)^2\ge s^2/2-s_0^2$ and
$s^2\le4A_0'^2d$. For larger $s$
the probability vanishes. Thus $\sqrt d\zeta_c$ is a centered
sub-gaussian variable with $\psi_2$-norm bounded in terms of
$\lambda_1,\lambda_2$, and Lemma~\ref{lem:pre_subgaussian}(a) gives
$\E[e^{\lambda\sqrt d\zeta_c}]\le e^{C_{\rm sg}^2\lambda^2}$ for all real
$\lambda$, with $C_{\rm sg}$ depending on $\lambda_1,\lambda_2$. For
$v\in\R^r$ with $v_\ell\ne0$ this yields
$\E[e^{\langle v,Z_g\rangle}]\le e^{C_{\rm sg}^2\|v_\ell\|_2^2/d}\le
e^{B_1^2v^\top Q_gv/2}$ for $B_1^2\ge2C_{\rm sg}^2/0.8$, and the case $v_\ell=0$
is trivial. The last assertion of the lemma collects these facts:
the masses lie in $[1/2,3/2]\subseteq[B^{-1},B]$, the reference laws are
finitely supported, and the number of terms is $NH+(r-N-Hd)\le2r^2$,
as $N\le r$ and $H\le r$.
\end{proof}

\subsection{One Tag at High Accuracy}
At accuracy $L\ge ar$ we use a single tag, a much smaller width and a
dictionary whose size may be far larger than $e^{Cd}$. The Gaussian
threshold is then too high for the clipping of
Lemma~\ref{lem:s-width}. Instead a single test integrates the
coefficient mass over the queries, and a tail bound for the
correlations replaces their fourth moment.

\begin{lemma}[A one-tag dictionary]\label{lem:s-onetag-frame}
Fix $p>0$ and $a'>0$. There are absolute constants $c_0,c_1>0$ and
constants $c_{\rm ot}>0$, $C_{\rm ot}\ge1$ and $d_{\rm ot}$, depending
only on $p$ and $a'$, with the following property. Let
$d\ge d_{\rm ot}$ and $\Xi\ge e^{a'd}$, and define
\[
 w:=c_0d^{-1/2}\Xi^{-1/d},\quad J:=\lfloor1/(8w)\rfloor+1,\quad
 s_j:=\tfrac12+4(j-1)w\ (j\in[J]),\quad \tau:=(c_1d\log(2d))^{1/2}.
\]
Let $\vartheta_\tau$, $\varpi$ and $u_0$ be as in
Lemma~\ref{lem:s-gauss}(c), and let $\mathcal E$ be an event for
independent $N(0,I_d/d)$ vectors $u_1,\ldots,u_\Xi$ with
$\Pr[\mathcal E]\ge9/10$. Then there are $u_1,\ldots,u_\Xi$ in
$\mathcal E$ and a probability law $\mathcal Q$ on $\mathbb S_d$ such
that the $(1,u_t)$ span $\R^{d+1}$, such that for $\mathcal Q$-a.e.\
$X$ and every $j$
\[
 \Xi^{-1}\sum_t|\langle u_t,X\rangle|^p\le C_{\rm ot},\qquad
 c_{\rm ot}\bar\sigma\le\|\mathbf f_j(X)\|_2\le C_{\rm ot}\bar\sigma,
\]
and such that, for independent $X,Y\sim\mathcal Q$ and independent
uniform $j,l\in[J]$, $U:=\langle v_j(X),v_l(Y)\rangle$ satisfies
\begin{equation}\label{lem:s-eq-onetag-tail}
 \E[U^2]\le\frac{C_{\rm ot}}\Xi,\qquad
 \Pr[|U|>u_0]\le\frac{C_{\rm ot}\varpi}\Xi .
\end{equation}
\end{lemma}

\begin{proof}
Let $c_0:=\min\{e^{-3}/2,1/C_0\}/16$ and $c_1:=1/2$. Then $w\le1/8$,
$s_J\le1$, $2s_j\ge1\ge C_0w$, $1/J\le8w$, and the $s_j$ lie in
$\frac12+w\mathbb Z$. Let $u_1,\ldots,u_\Xi$ be independent
$N(0,I_d/d)$ and let $\sigma_d$ be uniform on $\mathbb S_d$. Call $X$
good if $\Xi^{-1}\sum_th_j(\langle u_t,X\rangle)^2\in[c_3w/2,2C_3w]$
for every $j$ and $\Xi^{-1}\sum_t|\langle u_t,X\rangle|^p\le2m_p$, where
$m_p:=\E[|G|^p]$. As in the proof of Lemma~\ref{lem:s-frame}, a fixed
$X$ fails to be good with probability at most
$2Je^{-c_7\Xi w}+m_{2p}/(m_p^2\Xi)$. Since
$\Xi w\ge c_0\Xi^{1/2}/\sqrt d$ and $J\le\sqrt d\,\Xi^{1/d}/c_0$, this
tends to zero uniformly in $\Xi\ge e^{a'd}$. On good queries the
lengths satisfy the stated bounds with $c_{\rm ot}:=(c_3/2)^{1/2}$.

Let $X,Y\sim\sigma_d$ be independent, independent of the dictionary,
with angle $\theta$, $\rho:=\cos\theta$ and $z:=\sin\theta$, and let
$j,l$ be independent and uniform. Define
$\mu_{jl}:=\E[h_j(Z)h_l(W)]$ for $(Z,W)$ of correlation $\rho$, and
$B_t:=h_j(\langle u_t,X\rangle)h_l(\langle u_t,Y\rangle)$. If $X$ and
$Y$ are good, then
$|U|\le2(\Xi|\mu_{jl}|+|\sum_t(B_t-\mu_{jl})|)/(c_3\Xi w)$.

\emph{Second moment.} By Lemma~\ref{lem:s-scalar}(d) and $1/J\le8w$,
for each $l$ we have
$J^{-1}\sum_j(\mu_{jl}/w)^2\le8C_6^2we^{-z^2/w^2}$. By
Lemma~\ref{lem:s-gauss}(d),
\[
 8C_6^2w\,\E[e^{-z^2/w^2}]
 \le8C_6^2w(e^{-1/(4w^2)}+4w^{d-1}\Gamma((d+1)/2))
 \le\frac{e^{-2d}}\Xi
\]
for large $d$, because $\Gamma((d+1)/2)\le((d+1)/2)^{d/2}$,
$w^d((d+1)/2)^{d/2}\le(2c_0)^d\Xi^{-1}$ and
$1/(4w^2)=d\Xi^{2/d}/(4c_0^2)\ge3d+\log\Xi$. By
Lemma~\ref{lem:s-scalar}(e) and $\E[1/z]\le6$, the fluctuation
satisfies
$\E[|\sum_t(B_t-\mu_{jl})|^2]\le\Xi\E[B_t^2]\le C_3^2\Xi w^2\E[1/z]
\le6C_3^2\Xi w^2$. Hence, in expectation over the dictionary, the
queries and the indices,
$\E[U^2\mathbf 1_{XY}]\le(8/c_3^2)(e^{-2d}+6C_3^2)/\Xi$, where
$\mathbf 1_{XY}$ indicates that both queries are good.

\emph{Tail.} Define $c':=c_3/4$. If both queries are good,
$|\mu_{jl}|\le c'u_0w$ and $|\sum_t(B_t-\mu_{jl})|\le c'u_0\Xi w$, then
$|U|\le u_0$. By Lemma~\ref{lem:s-scalar}(d),
$|\mu_{jl}|>c'u_0w$ forces
$z<w\ell_0$ with $\ell_0:=(2\log(C_6/(c'u_0)))^{1/2}$, and, for fixed
$\theta$ and $l$, at most $C_6^2/(c'u_0)^2$ indices $j$ satisfy it.
Their fraction is at most $8C_6^2w/(c'u_0)^2$. By
Lemma~\ref{lem:s-gauss}(d), as $w\ell_0\le1/2$,
\[
 \Pr[|\mu_{jl}|>c'u_0w]\le\frac{8C_6^2w}{c'^2u_0^2}\cdot4(w\ell_0)^{d-1}
 =\frac{32C_6^2\ell_0^{d-1}}{c'^2u_0^2}\cdot\frac{c_0^d}{d^{d/2}\Xi}
 \le\frac{e^{-d\log(2d)/3}}\Xi
\]
for large $d$, since $u_0^{-1}=O(d\log(2d))$ and $\ell_0^2=O(\log d)$.
Given the queries and indices, the $B_t-\mu_{jl}$ are independent,
centered, bounded by $2$, with variance at most $\E[h_j(Z)^2]\le C_3w$.
Bernstein's inequality (Lemma~\ref{lem:pre_bernstein}) gives
$\Pr[|\sum_t(B_t-\mu_{jl})|>c'u_0\Xi w]\le2\exp(-c''u_0^2\Xi w)$ with
$c''>0$ absolute. Here $u_0^2\Xi w\ge c_0\Xi^{1/2}/(100\sqrt d(\tau+1)^4)
\ge\Xi^{1/3}$ for large $d$. By Lemma~\ref{lem:s-gauss}(c),
$\varpi\ge c_\phi\phi(\tau)/\tau\ge\exp(-d\log(2d)/4-\log\tau-4)$.
Hence both probabilities are at most $\varpi/(2\Xi)$ for large $d$,
and $\E[\Pr[|U|>u_0,\ X,Y\ {\rm good}]]\le\varpi/\Xi$.

\emph{Selection.} As in the proof of Lemma~\ref{lem:s-frame}, Markov's
inequality (Lemma~\ref{lem:pre_markov}) and $\Pr[\mathcal E]\ge9/10$ give a dictionary in
$\mathcal E$ with distinct rows, with the $(1,u_t)$ spanning, with
good-query measure at least $1/2$, and with the two expectations above
at most $40$ times their means. Conditioning $\sigma_d$ on the good
set multiplies them by at most $4$. This
gives Eq.~\eqref{lem:s-eq-onetag-tail}.
\end{proof}

\begin{lemma}[One tag at high accuracy]\label{lem:s-onetag}
Fix a non-even $p>0$ and $a>0$. There are $c_{\rm 1t}>0$ and $r_0$,
depending only on $p$ and $a$, such that the following holds for
$r\ge r_0$ and $L\ge ar$. Define $d:=r-1$, $\eta:=(r-1)/(r+2p)$ and
$R_{\rm 1t}:=(\epsilon^{-2}r^{1/2-p}\log(2r))^\eta$.
\begin{enumerate}[label=(\alph*),nosep,leftmargin=*]
\item There is a rank-$r$ input with rows $(1,u_t)\in\R\times\R^d$,
$t\in[\Xi]$, and masses in $[1/2,3/2]$ such that every reweighting of
its rows by real weights with relative error at most $\epsilon$ gives
nonzero weight to more than $c_{\rm 1t}R_{\rm 1t}$ rows. In
particular $S_p(r,\epsilon)>c_{\rm 1t}R_{\rm 1t}$.
\item Let $Q_0\ge1$ and $1\le q\le Q_0rR_{\rm 1t}^{1/(r-1)}$, and let
$r\ge r_0'$ with $r_0'$ depending on $p$, $a$ and $Q_0$. The input
in~(a) can be chosen with the following additional property. Define
$\ell_t:=\|(1,\sqrt du_t)\|_2$ and $\tilde u_t:=(1,\sqrt du_t)/\ell_t$, and
for masses $\alpha_t\in[1/2,3/2]$ let $\nu$ be the probability law
giving $\tilde u_t$ mass proportional to $\alpha_t\ell_t^p$. Then for every
choice of masses and every unit $x\in\R^r$,
\[
 \E_\nu[|\tilde u^\top x|^p]\ge c_{\rm fl}r^{-p/2},\qquad
 \E_\nu[(1+q|\tilde u^\top x|)^{-2}]\le A_{\rm sl}\sqrt r/q,
\]
where $c_{\rm fl}>0$ and $A_{\rm sl}$ depend only on $p$. For example,
if $q=\lfloor c(r-1)M^{1/(r-1)}\rfloor$ with $M=\lceil CR_{\rm 1t}\rceil$,
$c\le1\le C$, as in Lemma~\ref{lem:regular-high-accuracy}, then
$R_{\rm 1t}\ge1$ gives
$q\le c(r-1)(CR_{\rm 1t}+1)^{1/(r-1)}\le2rR_{\rm 1t}^{1/(r-1)}$ once
$r-1\ge\log_2(2C)$, so $Q_0=2$ applies with $r_0'$ depending on $C$.
The law $\nu$
is the input of~(a) after the invertible change of coordinates
$(y_0,y)\mapsto(y_0,\sqrt d\,y)$ of the rows, with row lengths
absorbed into the masses and a common normalization. Hence it
requires the same support.
\end{enumerate}
\end{lemma}

\begin{proof}
Let $c_0,c_1$ be the absolute constants and
$c_{\rm ot},C_{\rm ot}$ the constants of
Lemma~\ref{lem:s-onetag-frame} for $p$ and $a'=a$.

\emph{A single adaptive test.} Let $\mathcal Q$ and the test vectors be
as in Lemma~\ref{lem:s-onetag-frame}, and write $v$ for $v_j(X)$ with
$X\sim\mathcal Q$ and $j$ uniform. For $g\in\R^\Xi$ define
\[
 \mathbf A(g):=\frac{\bar\sigma}\varpi\E_v[\vartheta_\tau(\langle g,v\rangle)v],
 \qquad
 P_\tau(g):=\frac1\varpi\E_v[|\vartheta_\tau(\langle g,v\rangle)|] .
\]
Then $\E_g[P_\tau]=1$. Splitting according to $|U|\le u_0$ and using
Lemma~\ref{lem:s-gauss}(c) and Eq.~\eqref{lem:s-eq-onetag-tail},
\[
 \E_g[\|\mathbf A(g)\|_2^2]
 =\frac{\bar\sigma^2}{\varpi^2}\E_{v,v'}[U\,
   \E_g[\vartheta_\tau(\langle g,v\rangle)\vartheta_\tau(\langle g,v'\rangle)]]
 \le\bar\sigma^2(C_\phi\tau^2\E[U^2]+\frac{\Pr[|U|>u_0]}\varpi)
 \le C_8\sigma^2\tau^2,
\]
and in the same way
$\E[P_\tau^2]\le\varpi^{-2}(C_\phi\varpi^2+\varpi\Pr[|U|>u_0])\le C_8$, with
$C_8:=2C_{\rm ot}C_\phi$. By the Paley--Zygmund inequality (Lemma~\ref{lem:pre_paley_zygmund}),
$\Pr[P_\tau\ge1/2]\ge1/(4C_8)$. By Markov's inequality (Lemma~\ref{lem:pre_markov}), the event
$\mathcal E_A:=\{P_\tau\ge1/2,\ P_\tau\le C_9,\
\|\mathbf A\|_2\le C_9\sigma\tau\}$ has probability at least
$1/(8C_8)$ for $C_9:=16C_8$. On $\mathcal E_A$, as
$t\vartheta_\tau(t)\ge\tau|\vartheta_\tau(t)|$,
$\langle g,\mathbf A(g)\rangle\ge\bar\sigma\tau P_\tau(g)\ge
\bar\sigma\tau/2$.

\emph{Admissibility.} We have
$\mathbf A(g)=\E_X[\sum_jc_j(X)\mathbf f_j(X)]$ with
\[
 c_j(X):=\frac{\bar\sigma\,\vartheta_\tau(\langle g,v_j(X)\rangle)}
   {J\varpi\|\mathbf f_j(X)\|_2},\qquad
 \E_X[\sum_j|c_j(X)|]\le\frac{P_\tau(g)}{c_{\rm ot}} .
\]
This step charges the coefficient mass integrated over the queries, so it needs no
concentration over the centers at individual queries. With $N=1$,
the array
$\widetilde{\mathbf A}(g):=(c_{\rm ot}/C_9)\mathbf A(g)\mathbf 1_{\mathcal E_A}$
lies in $\mathcal B(\mathcal Q)$ and has norm at most
$c_{\rm ot}\sigma\tau$. For $R:=c_{\rm ot}\sigma\tau$ define
$\|\cdot\|_*$ as in Lemma~\ref{lem:s-width}. Then
$\|z\|_*\le R\|z\|_2$ and
\[
 \E[\|g\|_*]\ge\E[\langle g,\widetilde{\mathbf A}(g)\rangle]
 \ge\frac{c_{\rm ot}}{C_9}\cdot\frac{\bar\sigma\tau}2\cdot\frac1{8C_8}
 =\theta_0R\sqrt\Xi,\qquad \theta_0:=\frac1{16C_8C_9}.
\]

\emph{Parameters and packing.} Let $c_{\rm pk},n_0$ be as in
Lemma~\ref{lem:s-pack} for $\theta_0$, let $C_p$ be as in
Lemma~\ref{lem:s-cost}, and define
$c_2:=\min\{1,(c_{\rm pk}c_{\rm ot}/(C_p(1+C_{\rm ot})))^2
c_0^{2p+1}c_1\}$ and
\[
 \Xi:=\lfloor c_2(\epsilon^{-2}d^{1/2-p}\log(2d))^\eta\rfloor .
\]
Since $\eta=d/(d+2p+1)$ and $L\ge a(d+1)$, we have
$\log\Xi\ge\eta(2L-C_p'\log d)-C\ge ad$ for $d$ large, so
Lemma~\ref{lem:s-onetag-frame} applies. Apply Lemma~\ref{lem:s-cost}
with $N=H=1$, $b=1$, $C_E=C_{\rm ot}$ and the centers $s_j\in[1/2,1]$.
The rows $(1,u_t)$ have rank $r=d+1$. For a valid reweighting,
$\|\widehat\alpha-\alpha\|_*\le\epsilon_0:=C_p(1+C_{\rm ot})\epsilon\Xi$.
The condition $\epsilon_0\le c_{\rm pk}R\sqrt\Xi$ reads
$\Xi\le(c_{\rm pk}c_{\rm ot}/(C_p(1+C_{\rm ot})))^2\epsilon^{-2}
\sigma^2\tau^2$. As $\sigma^2\tau^2=c_0^{2p+1}c_1d^{1/2-p}\log(2d)
\Xi^{-(2p+1)/d}$, it is equivalent to
$\Xi^{1/\eta}\le c_2'\epsilon^{-2}d^{1/2-p}\log(2d)$ with $c_2'\ge c_2$,
which holds because $c_2\le1$ and $\eta\le1$ give
$c_2X^\eta\le(c_2X)^\eta$. Lemma~\ref{lem:s-pack} with $n=\Xi\ge n_0$
gives masses for which every valid reweighting has more than
$c_{\rm pk}\Xi$ nonzero output masses. Finally
$\Xi\ge\frac12c_2(\epsilon^{-2}d^{1/2-p}\log(2d))^\eta\ge cR_{\rm 1t}$,
since the bases differ by a factor bounded in terms of $p$ and
$\eta\le1$. This proves~(a), with constants uniform in $L\ge ar$.

\emph{The energy floor in~(b).} Write $Z_t:=\sqrt du_t$, standard
Gaussian in $\R^d$, and a unit query as $x=(x_0,v)$. Then
$\tilde u_t^\top x=(x_0+\langle v,Z_t\rangle)/\ell_t$. For Gaussian $Z$ the
variable $S:=(x_0+\langle v,Z\rangle)^2$ has $\E[S]=1$ and $\E[S^2]\le3$,
so the Paley--Zygmund inequality (Lemma~\ref{lem:pre_paley_zygmund}) gives
$\Pr[S\ge1/2]\ge1/12$. The set
$\{|x_0+\langle v,z\rangle|\ge2^{-1/2}\}$ is the complement of an
interval event in direction $v$ (or everything if $v=0$). We include in
$\mathcal E$ the event of Lemma~\ref{lem:s-halfspace} for the $Z_t$ and
the event $\Xi^{-1}\sum_t\ell_t^p\in[\E[\ell^p]/2,2\E[\ell^p]]$, which holds
with probability tending to one by Chebyshev. Then at least $\Xi/24$
indices have $|x_0+\langle v,Z_t\rangle|\ge2^{-1/2}$, and
$\E[\ell^p]\le C_pd^{p/2}$. Hence
\[
 \E_\nu[|\tilde u^\top x|^p]=\frac{\sum_t\alpha_t|x_0+\langle v,Z_t\rangle|^p}
   {\sum_t\alpha_t\ell_t^p}
 \ge\frac{2^{-p/2}\Xi/48}{3\Xi\E[\ell^p]}\ge c_{\rm fl}r^{-p/2}.
\]

\emph{A population slab estimate.} Let $Z$ be standard Gaussian,
$R_0:=\|Z\|_2$, $\ell:=(1+R_0^2)^{1/2}$ and $\tilde u:=(1,Z)/\ell$. Write
$Z=R_0\Theta$ with $\Theta$ uniform on $S^{d-1}$ and independent of
$R_0$. The coordinate density of $\Theta$ is
$c_d(1-t^2)^{(d-3)/2}\le\sqrt d\,e^{-(d-3)t^2/2}$, by
Lemma~\ref{lem:s-gauss}(d). Fix $\delta>0$ and a unit $x=(x_0,v)$. The
event $|\tilde u^\top x|\le\delta$ is
$|x_0+R_0\langle v,\Theta\rangle|\le \delta\ell$. Conditionally on $R_0$:
if $|x_0|<1/2$, then $\|v\|_2\ge\sqrt3/2$ and
$\langle v/\|v\|_2,\Theta\rangle$ lies in an interval of length
$2\delta\ell/(R_0\|v\|_2)$, so the probability is at most
$3\delta\ell\sqrt d/R_0$. If $|x_0|\ge1/2$ and $\delta\ell\ge1/4$, the probability
is at most $1\le4\delta\ell$. If $|x_0|\ge1/2$ and $\delta\ell<1/4$, then on the
event $|R_0\langle v,\Theta\rangle|\ge|x_0|-\delta\ell\ge|x_0|/2$, so the event
is empty if $v=0$. Otherwise $\langle v/\|v\|_2,\Theta\rangle$ lies in an interval of
length $2\delta\ell/(R_0\|v\|_2)$ at distance at least
$\beta_1:=|x_0|/(2R_0\|v\|_2)$ from zero. As
$1/(R_0\|v\|_2)\le4\beta_1$, the probability is at most
$8\delta\ell\beta_1\sqrt de^{-(d-3)\beta_1^2/2}\le8\delta\ell$ for $d\ge5$. In
all cases the conditional probability is at most
$8\delta\ell(1+\sqrt d/R_0)$. Since $\E[\ell^{p+1}]\le C_pd^{(p+1)/2}$ and
$\E[\ell^{p+1}/R_0]\le2^{(p+1)/2}(\E[1/R_0]+\E[R_0^p])\le C_pd^{p/2}$,
\[
 \E[\ell^p\mathbf 1_{|\tilde u^\top x|\le\delta}]
 \le8\delta\,\E[\ell^{p+1}(1+\sqrt d/R_0)]\le C_{10}d^{(p+1)/2}\delta .
\]
For $s\in\R$,
\[
 (1+q|s|)^{-2}\le\mathbf 1_{|s|\le1/q}
   +\sum_{k\ge1}4^{1-k}\mathbf 1_{|s|\le2^k/q}.
\]
Hence $\E[\ell^p(1+q|\tilde u^\top x|)^{-2}]\le5C_{10}d^{(p+1)/2}/q=:\mu_0$.

\emph{Uniform empirical comparison.} Define $H_1:=3(d+\log\Xi)^{1/2}$,
an upper bound for every $\ell_t$ on the event of
Lemma~\ref{lem:s-halfspace}. For unit $x$ let
$k_x(Z):=\ell^p(1+q|\tilde u^\top x|)^{-2}\mathbf 1_{\ell\le H_1}$, with
values in $[0,H_1^p]$, mean at most $\mu_0$ and variance at most
$H_1^p\mu_0$. Bernstein's inequality (Lemma~\ref{lem:pre_bernstein}) gives
$\Pr[\Xi^{-1}\sum_tk_x(Z_t)>2\mu_0]\le2\exp(-3\Xi\mu_0/(8H_1^p))$.
For unit $x,x'$, $|k_x-k_{x'}|\le2qH_1^p\|x-x'\|_2$. Take a
net of the unit sphere of $\R^{d+1}$ with mesh
$\rho_q:=\min\{1,\mu_0/(2qH_1^p)\}$ and at most $(3/\rho_q)^{d+1}$
points (Lemma~\ref{lem:pre_nets}). As
$\Xi^{1/d}\ge(R_{\rm 1t}/C)^{1/d}$ for a constant $C$ depending on
$p$, we have $q\le4Q_0d\,\Xi^{1/d}$ for large $d$. Then
$\Xi\mu_0/H_1^p\ge\Xi^{1/3}$, while the logarithm of the net size is
$O_{p,Q_0}(d\log d+\log\Xi+d\log\log\Xi)$, which is $o(\Xi^{1/3})$
uniformly in $\Xi\ge e^{ad}$. Hence,
outside probability tending to zero, every unit $x$ satisfies
$\Xi^{-1}\sum_t\ell_t^p(1+q|\tilde u_t^\top x|)^{-2}\le3\mu_0$. We include
this event in $\mathcal E$. Dividing by
$\Xi^{-1}\sum_t\alpha_t\ell_t^p\ge\E[\ell^p]/4\ge c_pd^{p/2}$ and using
$\alpha_t\le3/2$ gives
$\E_\nu[(1+q|\tilde u^\top x|)^{-2}]\le A_{\rm sl}\sqrt r/q$. All events included
in $\mathcal E$ have probability tending to one, so
$\Pr[\mathcal E]\ge9/10$ for large $r$. The change of coordinates is
invertible and acts on queries contragrediently, and multiplying all
masses by a common factor does not change relative errors. This
proves~(b).
\end{proof}